\documentclass[12pt,a4paper, titlepage, headinclude, withmarginpar]{book}
\usepackage[utf8]{inputenc}
\usepackage[british]{babel}
\usepackage{xargs}                      
\usepackage{cite}
\usepackage[pdftex,dvipsnames]{xcolor}  

\usepackage{amsmath, amsthm, amsfonts, amssymb}
\usepackage[alphabetic]{amsrefs}   
\usepackage{mdframed}
\usepackage{bm}   
\usepackage{slashed}     
\usepackage{xcolor}    
\usepackage[pdftex]{graphicx}
\usepackage[all]{xy}  
\usepackage{ifthen}
\usepackage{tikz}
\usepackage{sseq} 

\usepackage[left=3cm,right=3cm,top=2.5cm,bottom=2cm,includeheadfoot]{geometry}

\usepackage{index}
\usepackage{epigraph} 
\usepackage[colorinlistoftodos,prependcaption,textsize=tiny,disable]{todonotes} 
\usepackage[intoc, english, refpage]{nomencl}                    
\makenomenclature

\usetikzlibrary{matrix,arrows}

\usepackage{wrapfig}   

\usepackage{hyperref}  
\hypersetup{colorlinks=true, linkcolor=blue, citecolor=blue, urlcolor=blue,}

\theoremstyle{plain} 
\newtheorem{proposition}{Proposition}[section]
\newtheorem{theorem}[proposition]{Theorem}
\newtheorem{lemma}[proposition]{Lemma}
\newtheorem{cor}[proposition]{Corollary}

\newtheorem*{ConjectureNon}{Conjecture}

\newtheorem{thmx}{Theorem} 

\theoremstyle{definition}
\newtheorem{definition}[proposition]{Definition}

\newtheorem{setup}[proposition]{Setup}

\theoremstyle{remark} 
\newtheorem{rem}[proposition]{Remark}
\newtheorem{example}[proposition]{Example}

\newcommandx{\unsure}[2][1=]{\todo[linecolor=red,backgroundcolor=red!25,bordercolor=red,#1]{#2}}
\newcommandx{\change}[2][1=]{\todo[linecolor=blue,backgroundcolor=blue!25,bordercolor=blue,#1]{#2}}
\newcommandx{\info}[2][1=]{\todo[linecolor=OliveGreen,backgroundcolor=OliveGreen!25,bordercolor=OliveGreen,#1]{#2}}
\newcommandx{\improvement}[2][1=]{\todo[linecolor=Plum,backgroundcolor=Plum!25,bordercolor=Plum,#1]{#2}}
\newcommandx{\thiswillnotshow}[2][1=]{\todo[disable,#1]{#2}}

\newcommandx{\colim}[1][1=]{\mathop{\underrightarrow{\mathrm{colim}}_{#1}}}
\newcommandx{\SingMetNon}[2][1=M, 2=\bullet, usedefault]{\mathcal{R}_{#2}(#1)}
\newcommandx{\SingMet}[2][1=M, 2=\bullet, usedefault]{\mathcal{R}^+_{#2}(#1)} 
\newcommandx{\AuxSingMet}[1][1=\bullet]{{A_{#1}}}              
\newcommandx{\ConcSet}[2][2=M, 1=\bullet, usedefault]{\widetilde{\mathcal{R}^{+}_{#1}}(#2)} 
\newcommandx{\BlockMetrics}[2][2=M,1=\bullet]{\widetilde{\mathcal{R}_{#1}}(#2)}
\newcommandx{\BordPscSet}[2][2=d, 1=\bullet, usedefault]{\mathbf{Bord}^{+}_{#1}(#2)}
\newcommandx{\BordSet}[2][2=d, 1=\bullet, usedefault]{\mathbf{Bord}_{#1}(#2)}
\newcommandx{\genConcSet}[2][2=d, 1 = \bullet, usedefault]{\mathrm{\mathbf{Conc}}_{#1}(#2)}
\newcommandx{\ConcBndl}[2][2=d, 1 = \bullet, usedefault]{\mathrm{\mathbf{ConcBndl}}_{#1}(#2)}
\newcommandx{\pscgenConcSet}[2][2=d, 1 = \bullet, usedefault]{\mathrm{\mathbf{Conc}}^{+}_{#1}(#2)}
\newcommandx{\pscPseudoIso}[2][1=M, usedefault]{{C}^{\mathrm{psc}}_{#2}(#1)}
\newcommandx{\hofib}[1][1=\bullet,usedefault]{\mathrm{hofib}_{#1}}
\newcommand{\Diff}{\mathrm{Diff}}
\newcommandx{\BlockDiff}[1][1=\bullet]{\widetilde{\mathrm{Diff}_{#1}}}
\newcommand{\hAut}{\mathrm{hAut}}
\newcommandx{\BlockhAut}[1][1=\bullet]{\widetilde{\mathrm{hAut}_{#1}}}
\newcommand{\BoxCat}{\Box}   
\newcommand{\cSet}{\mathbf{cSet}}  
\newcommandx{\face}[2][1=\varepsilon, usedefault]{\partial^{#1}_{#2}}
\newcommandx{\degen}[1]{\sigma_{#1}}
\newcommandx{\CubeIncl}[2][1=\varepsilon, usedefault]{{\delta^{#1}_{#2}}}
\newcommandx{\CIncl}[2][1=\varepsilon, usedefault]{\mathrm{incl}_{#2, #1 R}}
\newcommandx{\CubeProj}[1]{p_{#1}}
\newcommandx{\CubeBox}[3][1=i, 2=\varepsilon, 3=n]{\sqcap^{#3}_{(#1,#2)}}
\newcommandx{\StandCube}[2][1=\bullet, usedefault]{\Box [#2]_{#1}} 
\newcommandx{\quader}[3][1=i, 2=\epsilon, 3=R, usedefault]{{\mathrm{Cub}_{#1,#2}^{#3}}}
\newcommandx{\subplane}[2][1=n]{\mathbb{R}^{#1}\langle #2 \rangle}
\newcommandx{\auxmet}[2][1=g]{{}_{#2}#1}
\newcommandx{\aux}[1][1=]{\mathrm{aux}_{#1}}
\newcommandx{\Aux}[1][1=]{\mathrm{Aux}_{#1}}
\newcommandx{\grad}[2][1=]{\mathrm{grad}(#2)}
\newcommandx{\hemi}[2][1=n, 2=+, usedefault]{D^{#1}_{#2}}
\newcommandx{\dicehemi}[2][1=n,2=+,usedefault]{\mathord{\text{\textit{\DJ}}}^{#1}_{#2}}
\newcommandx{\Torp}[1][1=n]{\mathrm{Torp}^{#1}}
\newcommandx{\capsubspace}[1]{\mathord{E_{#1}^{\varepsilon_{#1}}}}
\newcommandx{\normdirect}[1][1=p]{\mathord{I(#1)}}
\newcommandx{\Sing}[3][1=\bullet, 2=, usedefault]{S_{#1}^{#2}\left(#3\right)}
\newcommandx{\SingSmooth}[3][1=\bullet, 2=\infty, usedefault]{\mathcal{S}_{#1}^{#2}\left(#3\right)}
\newcommandx{\inddif}[3][1=\mathrm{H}, 2=\cdot, 3=g_0, usedefault]{{\mathrm{inddif}^{#1}(#2,#3)}}
\newcommandx{\KO}[1]{{\mathrm{KO}^{#1}}}
\newcommandx{\Kom}[1][1=]{{\mathbf{Kom}_{#1}}}

\newcommandx{\Pseud}[3][1=m,2=E,3=F, usedefault]{\Psi\mathrm{DO}^{#1}(#2,#3)}
\newcommand{\PseudOp}{\Psi\mathrm{DO}}
\newcommandx{\PseudCl}[3][1=m,2=E,3=F, usedefault]{\overline{\Psi\mathrm{DO}}^{#1}(#2,#3)}
\newcommand{\PseudOpCl}{\overline{\Psi\mathrm{DO}}}
\newcommandx{\BlockPseud}[4][1=m,2=E,3=F,4=\bullet, usedefault]{\Psi\mathrm{DO}^{#1}_{#4}(#2,#3)}
\newcommandx{\PseudDir}[2][1=, 2=M , usedefault]{\Psi\mathrm{Dir}_{#1}(#2)}
\newcommandx{\InvPseudDir}[2][1= , 2=M , usedefault]{\Psi\mathrm{Dir}^{\times}_{#1}(#2)}
\newcommandx{\PathDirac}[1][1=\bullet]{{\mathcal{X}_{#1}}}
\newcommandx{\AuxInvDir}[1][1=\bullet]{{B_{#1}}}
\newcommandx{\BlockDirac}[2][1=\bullet, 2=M, usedefault]{\widetilde{\Psi\mathrm{Dir}}_{#1}(#2)}
\newcommandx{\CliffRight}{\mathbf{r}}  
\newcommandx{\InvBlockDirac}[2][1=\bullet, 2=M, usedefault]{\widetilde{\Psi\mathrm{Dir}}^{\times}_{#1}(#2)}
\newcommand{\Dirac}{{\slashed{\mathfrak{D}}}}
\newcommand{\DiracSt}{{\slashed{\mathit{D}}}}  
\newcommand{\Err}{\mathfrak{Err}}
\newcommand{\lowBnd}{\mathrm{low}} 
\newcommand{\Shift}{\mathrm{Sh}} 

\newcommandx{\Sobolev}[3][1= ]{H_{#1}^{#2}(#3)}
\newcommandx{\normSo}[2][1=s]{||#2||_{#1}}
\newcommandx{\ind}[1][1=]{\mathrm{index}_{#1}}
\newcommandx{\smoothabs}[1]{\langle #1\rangle}   
\newcommand{\vol}{\mathrm{vol}}

\newcommand{\symb}{\mathrm{symb}}   
\newcommand{\Symb}{\mathrm{Symb}}   
\newcommand{\Gram}{\mathrm{Gram}}
\newcommand{\preGauge}{p\mathcal{G}}
\newcommand{\fullGauge}{\mathcal{G}}

\newcommand{\DecompMap}{\Phi}   

\newcommand{\euclmetric}{\langle \cdot, \cdot\rangle}   
\newcommand{\pr}{\mathord{\mathrm{pr}}}                 
\newcommand{\restrict}{\mathord{\upharpoonleft}}        
\DeclareMathOperator*{\im}{\mathrm{im}}
\newcommand{\N}{\mathord{\mathbb{N}}}                   
\newcommand{\Z}{\mathord{\mathbb{Z}}}                   
\newcommand{\Q}{\mathord{\mathbb{Q}}}                   
\newcommand{\R}{\mathord{\mathbb{R}}}                   
\newcommand{\C}{\mathord{\mathbb{C}}}                   
\newcommand{\I}{\mathord{\mathrm{i}}}                   

\newcommand{\GL}{\mathord{\mathrm{GL}}}

\newcommand{\Or}{\mathrm{O}}
\newcommand{\SOr}{\mathrm{SO}}
\newcommand{\Spin}{\mathord{\mathrm{Spin}}}
\newcommandx{\Cliff}[2][1=n, 2=0, usedefault]{Cl_{#1,#2}}  
\newcommand{\Clifford}{Cl}
\newcommand{\Spinor}{\mathfrak{S}}
\newcommand{\Hom}{\mathrm{Hom}}

\newcommand{\scal}{\mathord{\mathrm{scal}}}

\newcommand{\eps}{\varepsilon}
\newcommand{\fateps}{{\bm{\varepsilon}}}
\newcommand{\fateta}{{\bm{\eta}}}
\newcommand{\fatalpha}{{\bm{\alpha}}}
\newcommand{\incl}{\mathrm{incl}}
\newcommand{\diff}{\mathrm{d}}
\newcommand{\diffj}{\text{\dj}}

\newcommand{\dice}{\mathrm{dice}}
\newcommand{\id}{\mathrm{id}}

\newcommand{\Riem}{\mathrm{Riem}}
\newcommand{\const}{\mathrm{const}}
\newcommand{\supp}{\mathrm{supp}\,}
\newcommand{\susp}{\mathrm{susp}}

\newcommandx{\End}[3][1=]{\mathrm{End}(#2,#3)}
\newcommand{\evenodd}{\mathrm{e/o}}
\newcommand{\extension}{\mathrm{ex}}
\newcommand{\dom}{{\mathrm{dom}\, }}
\newcommand{\per}{\mathrm{per}}  
\newcommand{\Cliffmult}{\mathbf{c}} 
\newcommand{\cycl}{\mathrm{cycl}}
\newcommand{\Cycl}{\mathrm{Cycl}}
\newcommand{\trace}{\mathrm{tr}}
\newcommand{\op}{\mathrm{op}}
\newcommand{\placeholder}{\text{--}}
\newcommand{\Pos}{\mathrm{Pos}}

\newcommand{\deff}{\mathrel{\vcenter{\offinterlineskip\hbox{.}\vskip-.80ex\hbox{.}}}\joinrel \hskip 1pt =} 

\usepackage[T1]{fontenc}

\makenomenclature
\makeindex

\begin{document}
 
 \pagenumbering{gobble}  

 \newgeometry{left=3cm,right=2.5cm,top=2.5cm,bottom=2cm}
\begin{center}
\huge
\vspace*{0,45in}
\textbf{Concordances in Positive Scalar Curvature and Index Theory}
\\
\vspace*{0.9in}
\normalsize
\textsc{Dissertation}\\
\textsc{zur Erlangung des mathematisch-naturwissenschaftlichen Doktorgrades}\\
\textsc{"Doctor rerum naturalium"}\\
\textsc{der Georg-August Universit\"at G\"ottingen}\\
\vspace*{0.9in}
\textsc{im Promotionsstudiengang Mathematical Sciences}\\
\textsc{der Georg-August University School of Science (GAUSS)}\\
\vspace*{0.9in}
\textsc{vorgelegt von}\\
\textsc{Thorsten Hertl}\\
\textsc{aus G\"ottingen}\\
\normalsize
\vspace*{0.9in}
\textsc{G\"ottingen, 2022}
\end{center}
\restoregeometry 

\newpage

\newgeometry{left=3cm,right=2.5cm,top=2.5cm,bottom=2cm}

\textbf{Betreuungsausschuss}
\begin{itemize}
    \item[ ] \textbf{Prof. Dr. Thomas Schick} 
    \begin{itemize}
        \item[ ] Mathematisches Institut, Georg-August-Universität Göttingen
    \end{itemize}
    \item[ ] \textbf{Prof. Dr. Wolfgang Steimle} 
    \begin{itemize}
        \item[ ] Institut für Mathematik, Universität Augsburg
    \end{itemize}
\end{itemize}

\vspace{1cm}

\textbf{Mitglieder der Prüfungskommission}
\begin{itemize}
    \item[ ] Referent: \textbf{Prof. Dr. Thomas Schick} 
    \begin{itemize}
        \item[ ] Mathematisches Institut, Georg-August-Universität Göttingen
    \end{itemize}
    \item[ ] Koreferent: \textbf{Prof. Dr. Wolfgang Steimle} 
    \begin{itemize}
        \item[ ] Institut für Mathematik, Universität Augsburg
    \end{itemize}
\end{itemize}

\vspace{1cm}

\textbf{Weitere Mitglieder der Prüfungskommission}
\begin{itemize}
    \item[ ] \textbf{Prof. Dr. Victor Pidstrygach}
    \begin{itemize}
        \item[ ] Mathematisches Institut, Georg-August Universität Göttingen
    \end{itemize}
    \item[ ] \textbf{Prof. Dr. Frank Gounelas}
    \begin{itemize}
        \item[ ] Mathematisches Institut, Georg-August Universität Göttingen
    \end{itemize}
    \item[ ] \textbf{Prof. Dr. Dorothea Bahns}
    \begin{itemize}
        \item[ ] Mathematisches Institut, Georg-August Universität Göttingen
    \end{itemize}
    \item[ ] \textbf{Prof. Dr. Stephan Huckemann}
    \begin{itemize}
        \item[ ] Institut für mathematische Stochastik, Georg-August Universität Göttingen
    \end{itemize}
\end{itemize}

\vspace{1cm}

\textbf{Tag der mündlichen Prüfung:} 02.09.2022
\restoregeometry
\newpage

\Huge\noindent\textbf{Danksagung:}
\vspace{1cm}
\normalsize

\noindent Diese Dissertation ist das Endprodukt einer wunderbaren Zeit voller Schaffensperioden und Durststrecken gleichermaßen.
Sie hätte niemals das Licht der Welt erblickt, hätte ich nicht meine beiden Betreuer Prof. Dr. Thomas Schick und Prof. Dr. Wolfgang Steimle um mich gehabt, die mir das Vertrauen geschenkt haben, an dem Schwerpunktsprojekt \glqq Diffeomorphisms and the topology of positive scalar curvature\glqq{} mitzuarbeiten.
In der Zeit ließen sie mir den Freiraum, meine eigenen Ideen auszuprobieren, standen mir aber auch mit Rat und Tat zur Seite, wenn es nicht wie gewünscht lief.
Die regelmäßigen Gespräche mit Wolfgang gerade zu Beginn meiner Promotion haben mich davor bewahrt den roten Faden zu verlieren und von der Menge des neu zu lernenden Materials erschlagen zu werden. 
Thomas' Expertise und Intuition, die sich in den 12 Jahren, in denen ich ihn kenne, nie als falsch heraus gestellt hat, hat mich nie an dem Erfolg des indextheoretischen Projektes dieser Dissertation zweifeln lassen, obwohl die Vorarbeit mehr als 1.5 Jahre veranschlagt hat.
Ich bin beiden zu größtem Dank verpflichtet; vielen, vielen Dank!

Neben meinen Betreuern konnte ich mich auf die Hilfe vieler Korrekturleser verlassen.
Vielen Dank geht an Jonas Baltes, Tim Höpfner, Collin Mark Joseph, Alexei Kudryashov, Simon Naarmann, Daniel Räde, Hemanth Saratchandran und Julian Seipel. 
Besonderen Dank gilt Georg Frenck, der mehr Fehler gefunden hat als es mir lieb gewesen wäre und für seine Idee, die von Hand gezeichneten Grafiken zu verbessern, sowie Moritz Meisel, der einen Löwenanteil der gesamten Diss gelesen hat.

In meinen 2.5 Jahren in Augsburg hätte ich mich nicht so wohl gefühlt, wäre ich nicht so herzlich am Lehrstuhl für Differentialgeometrie aufgenommen worden.
Ich werde die Seerunden und den Kuchen vor dem Oberseminar nie vergessen.
Ganz besonders bedanken möchte ich mich hier bei Mauricio Bustamante, Benedikt Hunger, Alexei Kudryashov, Moritz Meisel, Hemanth Saratchandran, Eric Schlaarmann, sowie Alexandra Linde mit denen ich viele angeregte mathematische und nicht mathematische Diskussionen geführt habe, sowie außerhalb der Arbeitsgruppe, Pavel Hajek und Nick Wittig mit seiner Crew.

Wenn ich der Mathematik überdrüssig wurde, war ich froh einen sicheren Rückzugsort zu haben. Familie Gerstenacker und mein Mitbewohner Franz Korn haben dafür gesorgt, dass ich in Gersthofen eine neue Heimat gefunden habe.

In meinen 12 Jahren in Göttingen habe ich viele gute Freunde getroffen, die mich auf die ein oder andere Art und Weise unterstützt haben und deren Bekanntschaft ich nicht missen will. Hier seien namentlich erwähnt: Simone Cecchini, Jakob Dittmer, Anne Prepeneit, Simon Naarmann, Mehran Seyed Hosseini, Jonas Breustedt, Joel Meier und Tim Höpfner, sowie Christian Döring und Elias Most, die mich von Beginn an begleitet haben.

Ohne die Unterstützung meiner Familie mich zu meinen Bedingungen studieren zu lassen und meinen Interessen nachzugehen, wäre diese Arbeit nie entstanden, mein Dank gebührt euch.  

\newpage

\section*{Abstract}

We apply the strategy to study of diffeomorphisms via block diffeomorphisms to the world of positive scalar curvature (psc) metrics. 
For each closed psc manifold $M$, we construct the cubical set $\ConcSet$ of all psc block metrics, which only encodes concordance information of psc metrics within its homotopy type.
We show that $\ConcSet$ is a cubical Kan set, give a geometric description for the group structure of the combinatorial homotopy groups, and construct a comparison map from the cubical model $\SingMet$ of the space of psc metrics on $M$ to $\ConcSet$.
Next, we build a concordance-themed model for real $K$-theory based on the notion of invertible block Dirac operators and use it to factor the index difference through $\ConcSet$.
In the final part of this thesis, we construct the psc Hatcher spectral sequence, which is a non-index-theoretic tool to get information about the difference of $\SingMet$ and $\ConcSet$.

 \pagenumbering{arabic}
 \frontmatter
 \tableofcontents
 
 \mainmatter

 
 \chapter{Introduction}

The overall goal of this thesis is to systematically study the space of positive scalar curvature (psc) metrics on closed manifolds from a concordance view point.

\subsubsection*{Background}

The scalar curvature $\scal(g)$ of a Riemannian metric $g$ is a geometric quantity of the Riemannian manifold $(M,g)$ that measures how far this Riemannian metric deviates from being flat.
Among the commonly considered curvatures (sectional, Ricci, and scalar) in Riemannian geometry, scalar curvature is the weakest one.
It is a scalar-valued function $\scal(g) \colon M \rightarrow \R$ that measures on small scales how fast the volume of a geodesic ball grows compared to a ball in euclidean space.
More precisely, the scalar curvature at $p \in M$ appears in the Taylor series expansion \cite{gray1974volume}*{Theorem 3.1}
\begin{equation*}
    \frac{\vol(B_\eps(p) \subseteq M^d )}{\vol\bigl(B_\eps(0) \subseteq \R^d \bigr)} = 1 - \frac{\scal(g)(p)}{6(d+2)}\eps^2 + \mathcal{O}(\eps^4).
\end{equation*}
In particular, if the scalar curvature is positive, then the volume grows slower on those scales where the second order term dominates the higher order terms.

For surfaces, i.e., $\dim \, M = 2$, the scalar curvature is twice the Gauss curvature, so assuming positivity is a strict condition.
Indeed, the Gauss-Bonnet theorem 
\begin{equation*}
    \int_M \scal(g) \diff \, \vol_g = 4\pi \chi(M) = 4\pi\bigl(2-2\mathrm{genus}(M)\bigr)
\end{equation*}
together with the classification theorem for closed surfaces imply that only the sphere $S^2$ and the real projective plane $\R P^2$ can carry a positive scalar curvature metric, which they do!

In higher dimensions, the curvature notions differ, and it is a bit surprising that there are obstructions against the existence of positive scalar curvature metrics at all.
Lichnerowicz \cite{lichnerowicz1963spineurs} observed that on spin manifolds, i.e., orientable Riemannian manifolds with vanishing second Stiefel-Whitney class, the (Clifford-linear) Dirac operator $\Dirac_g$ is closely related to the scalar curvature of its underlying Riemannian metric:
\begin{equation*}
    \Dirac_g^2 = \nabla^\ast \nabla + \frac{\scal(g)}{4}.
\end{equation*}

This formula implies that if the scalar curvature is positive, then the Dirac operator is invertible.
On closed manifolds, the Dirac operator is an (unbounded) Fredholm operator and its kernel generates a class in real $K$-theory.
By Hitchin's extension of the Atiyah-Singer index theorem, this class is a purely topological quantity of the $d$-dimensional manifold $M$, the so called $\alpha$\emph{-invariant} or \emph{topological index} $\alpha(M) \in KO^{-d}(\mathrm{pt})$.
In particular, this class is independent of the metric that was used to define the Dirac operator and one infers that the non-vanishing of $\alpha(M)$ is an obstruction against positive scalar curvature.

On the constructive side, the surgery result of Gromov-Lawson \cite{gromov1980classification} is an important tool to produce examples of manifolds that carry a positive scalar curvature metric.
It says, in its simplest form, that if $M_2$ is obtained from $M_1$ by a surgery of codimension $\geq 3$ and if $M_1$ carries a psc metric, then $M_2$ also carries a psc metric.
The surgery result implies that a spin manifold of dimension $\geq 5$ carries a psc metric if and only if its singular bordism class $[M,f] \in \Omega^{\Spin}(B\pi_1(M))$ has a representative that carries a psc metric.
Here, $f \colon M \rightarrow B\pi_1(M)$ is a continuous map to the classifying space of the fundamental group of $M$ that induces the identity between the fundamental groups.

Stephan Stolz used in \cite{stolz1992simply} the surgery result to prove that, for simply connected spin manifolds of dimension $\geq 5$, the vanishing of the $\alpha$-invariant is not only necessary, but also sufficient by showing that all elements in the kernel of the $\alpha$-invariant
\begin{equation*}
    \alpha \colon \Omega_\ast^{\Spin}(\mathrm{pt}) \rightarrow KO^{-\ast}(\mathrm{pt})
\end{equation*}
are represented by psc-manifolds.
For other fundamental groups, the analogous question whether the (higher) $\alpha$-invariant, which we are not going to define here, is also a sufficient condition for existence became known as the Gromov-Lawson-Rosenberg conjecture and was proven to be false in general by the counterexample of Thomas Schick \cite{schick1998counterexample}. 

The stable Gromov-Lawson-Rosenberg conjecture claims that if the higher $\alpha$-invariant vanishes, then $M \times \mathrm{Bott}^N$ carries a positive scalar curvature metric for all sufficiently large $N \in \N$, see \cite{rosenberg2006manifolds}*{p.19 ff} for details.
Here, $\mathrm{Bott}$ is an eight-dimensional, simply-connected, closed manifold with $\alpha(\mathrm{Bott})=1 \in KO^{-8}(\mathrm{pt}) \cong \Z$, in particular it cannot carry a positive scalar curvature metric.
This manifold is referred to as \emph{Bott-manifold} and can be chosen to be Ricci-flat.
The stable Gromov-Lawson-Rosenberg conjecture, formulated in the 90s, is a consequence of the Novikov conjecture and is state of the art for the exisitence question.

\vspace{0.5cm}

Although not completely solved, the question whether $M$ carries a psc metric or not is reasonably well understood: it is governed by the stable Gromov-Lawson-Rosenberg conjecture. 
Thus, the focus has turned toward the question: 
\begin{center}
    How many essentially different psc metrics can a manifold carry if it carries at least one psc metric?
\end{center}
More precisely, on a closed manifold $M$, we consider the set $\Riem^+(M)$
of all Riemannian metrics on $M$ with positive scalar curvature.
This is a subset of sections on $M$ with values in bilinear forms on tangent spaces, so we equip it with the smooth Whitney topology - the usual topology on the space of sections.
We are interested in the homotopy type of $\Riem^+(M)$.

If $\dim M = 2,3$, then the space $\Riem^+(M)$ is either empty or contractible, see \cite{rosenberg149metrics}*{Theorem 3.4} and \cite{bamler2019ricci}*{Theorem 1.1}.
If $\dim M \geq 4$ and $M$ is spin, the situation changes dramatically.
Using a secondary version of the $\alpha$-invariant, the so called \emph{index difference}, many non-trivial elements in the homotopy groups of $\Riem^+(M)$ were discovered, see for example \cite{hitchin1974harmonic}, \cite{carr1988construction}, \cite{hanke2014space},  \cite{botvinnik2017infinite}, and \cite{frenck2021SphericalKappa}.
A priori, there are two versions of the index difference, one is due to Hitchin \cite{hitchin1974harmonic}, the other one is due to Gromov-Lawson \cite{gromov1983positive}.
Let us now describe these in turn.

Recall that $\mathrm{Fred}^{d,0}(H)$, the space of $\Cliff[d][0]$-linear Fredholm operators on a separable Hilbert space with a $\Cliff[d][0]$-module structure, is a classifying space for real $K$-theory.
This means that there is a natural isomorphism
\begin{equation*}
    [X; \mathrm{Fred}^{d,0}(H)] \cong KO^{-d}(X)
\end{equation*}
for all finite CW-complexes $X$.
Here, $\Cliff[d][0]$ denotes the real Clifford algebra, the universal algebra generated by $e_1,\dots, e_d$ under the relations $e_ie_j + e_je_i = -2\delta_{ij}$.

Fix a base point $g_0 \in \Riem^+(M)$.
For each psc metric $g \in \Riem^+(M)$, convex combination gives a path of Riemannian metrics $g_t$ (possibly without positive scalar curvature).
This path gives rise to a path of Dirac operators $\Dirac_{g_t}$ that is invertible near the end points.
Its bounded transform is a path of Fredholm operators with invertible ends\footnote{We ignore for this moment that the Hilbert space itself depend on the metric, too.}. 
The space of invertible operators is contractible by Kuiper's theorem, see \cite{ebert2017indexdiff}*{Lemma 2.8} for the Clifford-linear version, so we can close this path to a loop.
This procedure yields a continuous map
\begin{equation*}
    \mathrm{inddif}_H(\placeholder,g_0) \colon \Riem^+(M) \rightarrow \Omega\mathrm{Fred}^{d,0}(H).
\end{equation*}

The Gromov-Lawson version is constructed differently.
For a psc metric $g \in \Riem^+(M)$, we consider the metric $G$ on $M \times \R$ that is given by
\begin{equation*}
    G = \begin{cases}
       g_0 + \diff t^2, &\text{on } M \times \R_{\geq -1}, \\
       \frac{1-t}{2}g_0 + \frac{(1+t)}{2}g + \diff t^2, & \text{on } M \times [-1,1], \\
       g + \diff t^2, & \text{on } M \times \R_{\geq 1}.
    \end{cases}
\end{equation*}
Its Dirac operator is still an unbounded Fredholm operator, so its bounded transform is an element in $\mathrm{Fred}^{d+1,0}(H)$.
This gives a continous map
\begin{equation*}
    \mathrm{inddif}_{GL}(\placeholder,g_0) \colon \Riem^+(M) \rightarrow \mathrm{Fred}^{d+1,0}(H).
\end{equation*}
The main result of \cite{ebert2017indexdiff} shows that the two versions essentially agree, that means, composed with the weak homotopy equivalence $\mathrm{Fred}^{d+1,0}(H) \rightarrow \Omega \mathrm{Fred}^{d,0}(H)$, the Gromov-Lawson version of the index difference map is homotopic to the Hitchin version.

The important feature of the Gromov-Lawson version is its invariance under \emph{concordance} and not just \emph{isotopy}.
Let us recall these concepts.
An isotopy between two psc metrics is nothing but a smooth path between them, so the term isotopy is just a synonym for smooth homotopy.
A concordance between two psc metrics on the other hand is a psc metric $G$ on $M \times [-1,1]$ that agrees near the end points $M \times \{j\}$ with the metrics $g_j + \diff t^2$.
We call two psc metrics \emph{isotopic/concordant} if there is an {isotopy/concordance} between them.
Being isotopic implies being concordant.
It is an open conjecture whether the converse is also true.
\begin{ConjectureNon}[Concordance-Implies-Isotopy]
   If $g_0,g_1 \in \Riem^+(M)$ are concordant, then they are also isotopic.
\end{ConjectureNon}
There are counterexamples to this conjecture in dimension 4, see \cite{ruberman2002positive}, but this conjecture might still be true in higher dimension.

We observe a huge discrepancy:
Although we are interested in the homotopy type of the space $\Riem^+(M)$ (its isotopy type), our invariants at hand are invariant under concordance.
Thus, it makes sense to have a systematic look at $\Riem^+(M)$ from a concordance perspective.
However, results in this direction are rather scarce.
Stolz \cite{stolz1998concordance} considered the set of concordance classes of $\Riem^+(M)$ and showed that it only depends on the first Stiefel-Whitney class $w_1(M)$ and the group extension 
\begin{equation*}
    \xymatrix{1 \ar[r] & \Z_2 \ar[r] & \hat{\pi} \ar[r] & \pi_1(M) \ar[r] & 1}
\end{equation*}
that is classified by the second Stiefel-Whitney class $w_2(M)$, see \cite{stolz1998concordance} for details.

Motivated by this discrepancy and the results of Stolz, we systematically study the space of psc metrics from a concordance view point in this thesis.
Our ansatz is inspired from the study of diffeomorphism groups of high-dimensional manifolds started in the late 60's of the last century, which we are going to recall.


\subsubsection*{The theory of block diffeomorphisms as blue print}

The key idea to study the space of diffeomorphism groups on a compact manifold $M$ is to ``approximate'' the diffeomorphism group $\Diff(M)$ by the larger group of \emph{block diffeomorphisms} $\BlockDiff[ ](M)$.
This group of block diffeomorphism cannot be written down in a closed form but needs to be constructed with tools from combinatorial topology, like simplicial or cubical sets.
We use cubical sets here because we will also use them in the main body of this thesis. 
The group $\BlockDiff[ ](M)$ is thus the geometric realisation of the cubical group whose $n$-cubes are given by diffeomorphisms on $M \times [-1,1]^n$ that map each face $M\times \{x_i=\eps\}$, where $\eps = \pm 1$, to itself.
It contains $\Diff(M)$, or rather its combinatorial counterpart, as the cubical subgroup whose $n$-cubes are diffeorphisms $\varphi$ on $M \times [-1,1]^n$ that make the following diagram commute:
\begin{equation*}
    \xymatrix{M \times [-1,1]^n \ar[rr]^\varphi \ar[rd]_{\pr_2} && M\times [-1,1]^n \ar[ld]^{\pr_2}\\
    &[-1,1]^n.& }
\end{equation*}

The advantage of $\BlockDiff[ ](M)$ over $\Diff(M)$ is that it can be better handled with tools from algebraic topology like homotopy theory or surgery theory.
To explain this, we first need to introduce the topological monoid of \emph{homotopy automorphisms} $\hAut(M)$, which consists of all continuous self-maps that are homotopy equivalences.
There is also a block version $\BlockhAut[ ](M)$ that arises out of $\hAut(M)$ in the same way as $\BlockDiff[ ](M)$  arises out of $\Diff(M)$; but by \cite{dold1963partitions}*{Theorem 6.1} the canonical inclusion $\hAut(M) \rightarrow \BlockhAut[ ](M)$ is a homotopy equivalence.
The inclusion $\BlockDiff[ ](M) \rightarrow \BlockhAut[ ](M) \simeq \hAut(M)$ gives rise to a homotopy fibration
\begin{equation*}
    \xymatrix{\hAut(M)/\BlockDiff[ ](M) \ar[rr]&& B\BlockDiff[ ](M) \ar[rr] && B\hAut(M).}
\end{equation*}
The classifying space $B\hAut(M)$ is a purely homotopy-theoretic object and can be approached with tools from algebraic topology, like obstruction theory \cite{rutter1997spaces}.
The fibre $\hAut(M)/\BlockDiff[ ](M)$ has a close connection to surgery theory as it is weakly homotopy equivalent to the fibre in Quinn's surgery fibration
\begin{equation*}
    \xymatrix{\hAut(M)/\BlockDiff[ ](M) \ar[r]^-\simeq &\mathcal{S}(M) \ar[r] &\mathrm{Maps}(M,G/O) \ar[r] & \mathbb{L}(\Z [\pi_1M]), }
\end{equation*}
see \cite{berglund2012homological}*{Section 3} for a clear and concise discussion on the Surgery fibration and the weak equivalence statement.
The long exact sequence of homotopy groups of this fibration yields that $\pi_k(\hAut(M)/\BlockDiff[ ](M))$ appear in the long exact surgery sequence
\begin{equation*}
    \xymatrix@C-1em{L_{k+d+1}(\Z [\pi_1M]) \ar[r] & S_\partial(M \times D^k) \ar[r] & \mathcal{N}_\partial(M \times D^k) \ar[r] & L_{k+d}(\Z [\pi_1M])  \\
    & \ar[u]_\cong \pi_k\bigl(\hAut(M)/\BlockDiff[ ](M)\bigr). && }
\end{equation*}
We refer to \cite{Lueck2020surgery} for details on the surgery sequence.

The $L$-groups $L_\ast(\Z [\pi_1(M)])$ have a completely algebraic description, while the \emph{normal bordism groups} $\mathcal{N}_\partial(M \times D^k)$ can be studied via (stable) homotopy theory.
The upshot is that, despite its differential-topological origin of $\BlockDiff[ ](M)$, one can extract a lot information out of $\BlockDiff[ ](M)$ using only algebraic and homotopy-theoretic tools, which reduces the difficulty of the problem in practice.

To get information about the difference between $\Diff(M)$, the object of interest, and $\BlockDiff[ ](M)$, the object that is easier to study, one uses the fibration
\begin{equation*}
    \xymatrix{\BlockDiff[ ](M)/\Diff(M) \ar[r] & B\Diff(M) \ar[r] & \Diff(M).}
\end{equation*}
In \cite{hatcher1978concordance}, Hatcher examines the fibre $\BlockDiff[ ](M)/\Diff(M)$ with the help of the group of \emph{pseudo isotopies}
\begin{equation*}
    C(M) \deff \{\varphi \in \Diff(M \times [0,1]) \, : \, \varphi = \id \text{ near } M \times \{0\} \cup \partial M \times [0,1]\}.
\end{equation*}
More precisely, he constructs a first quadrant spectral sequence that starts with $E^1_{p,q} = \pi_q(C(M \times [0,1]^p))$ and converges to $\pi_{p+q+1}(\BlockDiff[ ](M)/\Diff(M))$.

If $q \ll \dim M$, then one can actually compute the differentials $d^1_{p,q}$ in the Hatcher spectral sequence and relate the results to Waldhausen's algebraic $K$-theory.
Farell and Hsiang \cite{farrell1978rational} used these tools to compute the rational homotopy groups $\pi_k(B\Diff(S^n)) \otimes \Q$ for $k \ll n$.
For completeness, the precise results is the following:
\begin{equation*}
   \pi_k(B\Diff(S^n)) \otimes \Q = \begin{cases}
     \Q, & \text{if } k \equiv 0 \mod 4, n \text{ even}, \\
     \Q \oplus \Q, & \text{if } k \equiv 0 \mod 4, n \text{ odd}, \\
     0, & \text{otherwise}
   \end{cases}
\end{equation*}
in the range $0 < k < n/6 - 7$.

\subsubsection*{Statement of the results}

The overall aim of this thesis is to develop the theoretical foundations to carry over the strategy of the study of diffeomorphism groups described in the last section to the world of psc metrics.
We will introduce the space $\ConcSet[ ]$ of \emph{psc block metrics}, which should be thought of as an analog of $\BlockDiff[ ](M)$ in the psc realm.

As before, the space $\ConcSet[ ]$ cannot be written down in a closed form but needs to be constructed. 
Our construction tool is cubical set theory because we believe that cubes are better suited to study curvature related questions than simplices due to the canonical identification $I^n \times I^m = I^{n+m}$.
In fact, we will mostly work in the category of cubical sets, instead of the category of topological spaces.
This is unproblematic because, in this thesis, we are only interested in homotopy-theoretic questions and the two categories have the same homotopy category.

To keep this thesis self-contained, we recall the foundations of cubical set theory we are going to use.
We first give a general introduction to cubical set theory and survey the existing literature before we discuss examples and variations of general concepts that are specifically adapted to this thesis. 

We then model the approximation space $\ConcSet[ ]$ by the cubical set $\ConcSet$ whose $n$-cubes can be thought of as psc metrics on $M \times [-1,1]^n$ that agree near each face $M \times \{x_i = \eps\}$ with the product metric $g \restrict_{M \times \{x_i = \eps\}} \oplus \diff x_i^2$, where $\eps = \pm 1$.

Usually, cubical set theory only works well for so called \emph{Kan-sets}, which should be thought of as cubical sets with ``enough cubes''.
The first main result of this thesis can be considered as a door-opener to the combinatorial world:

\begin{thmx}\label{KanSet - MainThm}
 If $M$ is a closed psc manifold, then $\ConcSet$ is a Kan set.
\end{thmx}

An important consequence of this result is that the homotopy groups of $\ConcSet[ ]$ have a purely combinatorial description.
For example, the set of path components $\pi_0(\ConcSet[ ])$ is precisely the set of concordance classes of $\Riem^+(M)$ and the fundamental group $\pi_1(\ConcSet[ ],g_0)$ is the set of concordance classes of self-concordances classes of $g_0$.
We know from Stolz' results in \cite{stolz1998concordance} that $\pi_0(\ConcSet)$ only depends on the fundamental group of $M$ if $M$ is spin, and we expect a similar result for the higher homotopy groups $\pi_n(\ConcSet)$.
 
Cubical set theory also provides an abstract group structure on the combinatorial versions of the homotopy groups.
However, this group structure is implicit and, in practice, not really useful for calculations.
We therefore present a more geometric and explicit description.
In the case of the fundamental group, the geometric product of two self-concordances $G_1, G_2$ of $g_0$ on $M \times [-1,1]$ can be thought as the union $G_1 \cup G_2$ on $M \times [-1,3]$, which we re-identify with $M \times [-1,1]$.

Of course, we would like to relate the original space $\Riem^+(M)$ to $\ConcSet[ ]$.
A cubical model for $\Riem^+(M)$ is the cubical set  $\SingMet$ whose $n$-cubes are smooth maps $[-1,1]^n \rightarrow \Riem^+(M)$ that near each hyperface $\{x_i = \eps\}$ are independent of $x_i$ .
The cubical map $\SingMet \rightarrow \ConcSet$ that assigns to each map of psc metrics its adjoint metric
\begin{equation*}
    \bigl(g \colon [-1,1]^n \rightarrow \Riem^+(M)\bigr) \mapsto g(\placeholder) + \diff x_1^2 + \dots + \diff x_n^2 \in \Riem^+(M \times [-1,1]^n)
\end{equation*}
is not well defined because the adjoint metric can have non-positive scalar curvature.
However, if we restrict this map to a homotopy equivalent cubical subset of $\SingMet$ whose cubes are sufficiently slowly parametrised, then the adjoint metric has positive scalar curvature as well so that we end up with a well defined map on a weakly equivalent subset $\susp_\bullet \colon \SingMet \dashrightarrow \ConcSet$. \vspace{0.5cm}

The next task is of course to find out whether the index difference (or rather its cubical model) factors through $\ConcSet$.
If this is true, then $\ConcSet$ would have many non-trivial homotopy groups whenever $\dim M \geq 6$ because the index difference induces non-trivial homomorphisms
\begin{equation*}
    \pi_n(\mathrm{inddif}(\placeholder,g_0)) \colon \pi_k(\Riem^+(M^d),g_0) \rightarrow KO^{-(d+k+1)}(\mathrm{pt})
\end{equation*}
whenever its target group is non-trivial; see \cite{botvinnik2017infinite}. 

The starting point is the observation made in \cite{ebert2017indexdiff} that the space of \emph{invertible pseudo Dirac operators} $\InvPseudDir[ ]$ on a closed, spin, psc manifold $M$ of dimension $d$ is a classifying space for $KO^{-d}$.
We prefer this model over $\mathrm{Fred}^{d,0}(H)$, because, in this picture, the index difference is the continuous map $\Riem^+(M) \rightarrow \InvPseudDir[ ]$ that assigns to a psc metric $g$ its Dirac operator $\Dirac_g$.
As for the space of psc metrics, a cubical model for $\InvPseudDir$ is given by the cubical set whose $n$-cubes are smooth maps $[-1,1]^n \rightarrow \InvPseudDir$ that near each hyperface $\{x_i = \eps\}$ are independent of $x_i$ and the continuous map $g \mapsto \Dirac_g$ induces a cubical map between the cubical models so that we end up with two cubical maps 
\begin{equation*}
    \xymatrix{ \SingMet \ar[r]^-{\Dirac_\bullet} \ar@{-->}[d]_{\susp_\bullet} & \InvPseudDir[\bullet] \\ \ConcSet & }
\end{equation*}
that we would like to complete with a cubical map $\ConcSet \rightarrow \InvPseudDir[\bullet]$ to a homotopy commutative diagram.

If such a map exists, then it must have a complicated description because an element in $\ConcSet[n]$ is a metric on $M \times [-1,1]^n$ while an element in $\InvPseudDir[n]$ is a map $[-1,1]^n \rightarrow \InvPseudDir[ ]$.
To overcome this problem, we apply the same procedure we used to construct the cubical set $\ConcSet$ out of the space $\Riem^+(M)$ to the space $\InvPseudDir$.
This yields the cubical set $\InvBlockDirac$ of \emph{invertible block Dirac operators}, whose $n$-cubes can be thought of as pseudo Dirac operators on $M \times [-1,1]^n$ that decompose near the boundary in the same manner as Dirac operators of block metrics would do.

In complete analogy to the map $\susp_\bullet \colon \SingMet \dashrightarrow \ConcSet$ there is a cubical map $\susp_\bullet \colon \InvPseudDir[\bullet] \dashrightarrow \InvBlockDirac$, well-defined only on a homotopy equivalent subset of $\InvPseudDir[\bullet]$, that assigns to a map of operators $[-1,1]^n \rightarrow \InvPseudDir$ its ``adjoint operator'' on $M \times [-1,1]^n$.

We can now formulate the next main theorem of this thesis.
\begin{thmx}\label{BlockModelKTheory - MainThm}
 The map $\susp_\bullet \colon \InvPseudDir[\bullet] \dashrightarrow \InvBlockDirac$ is a weak homotopy equivalence.
\end{thmx}
The important consquence of this theorem is that $\InvBlockDirac$ is another model for the classifying space of real $K$-theory that is big enough to host the Dirac operator of a general psc block metric.
We therefore get a cubical map $\Dirac_\bullet \colon \ConcSet \rightarrow \InvBlockDirac$ that assigns to a psc block metric its Dirac operator.

\begin{thmx}\label{Factorisation - MainThm}
 The diagram 
 \begin{equation*}
     \xymatrix{ \SingMet \ar[rr]^{\Dirac_\bullet} \ar@{-->}[d]_{\susp_\bullet} && \InvPseudDir[\bullet] \ar@{-->}[d]^{\susp_\bullet} \\ \ConcSet \ar[rr]_{\Dirac_\bullet} && \InvBlockDirac}
 \end{equation*}
 commutes up to homotopy.
\end{thmx}
Due to the nature of the construction, we interpret $\Dirac_\bullet \colon \SingMet \rightarrow \InvPseudDir[\bullet]$ as the cubical model of the Hitchin index difference and $\Dirac_\bullet \colon \ConcSet \rightarrow \InvBlockDirac$ as cubical model for the Gromov-Lawson index difference.

An immediate consequence of Theorem \ref{Factorisation - MainThm} is the main result of \cite{ebert2017indexdiff} that the Hitchin and the Gromov-Lawson version of the index difference agree.
However, Theorem \ref{Factorisation - MainThm} shows even more, because the index difference factors through $\pi_n(\ConcSet)$, which we believe to depend only on the fundamental group of the spin manifold $M$.
The image of the index difference in the real $K$-theory therefore only depends on the fundamental group.

\vspace{0.5cm}

Of course, as in the case of the diffeomorphism groups, we would like to know how much $\Riem^+(M)$ and $\ConcSet[ ]$ or, equivalently, their cubical models, $\SingMet$ and $\ConcSet$, differ from each other.
The concordance-implies-isotopy conjecture translates to the statement that $\susp_\bullet \colon \SingMet \dashrightarrow \ConcSet$ induces an injection (and hence a bijection) on the set of path components.
More daringly, we could conjecture:
\begin{ConjectureNon}[Strong Concordance-Implies-Isotopy]
   The cubical map 
   \begin{equation*}
      \susp_\bullet \colon \SingMet \dashrightarrow \ConcSet  
   \end{equation*}
   is a weak homotopy equivalence.
\end{ConjectureNon}
Although the concordance-implies-isotopy conjecture is true if $\dim M \leq 3$ (because $\Riem^+(M)$ is contractible), the stronger version is false. 
We will present counterexamples in Chapter \ref{Factorisation Indexdiff - Chapter}.
In fact, we show even more:
\begin{thmx}\label{ConterExampleStrongConcVsIsotopy - MainThm}
  If $M$ is a psc manifold of dimension $2$ or $3$, then the suspension map $\susp_\bullet \colon \SingMet \dashrightarrow \ConcSet$ is not a weak homotopy equivalence.
\end{thmx}
To get information about the difference between $\SingMet$ and $\ConcSet$, in particular for high-dimensional manifolds, for which we have no counterexamples, we mimick the approach from the study of diffeomorphism groups and consider the homotopy fibre $\mathrm{hofib}_\bullet$ of the comparision map $\susp_\bullet \colon \SingMet \dashrightarrow \ConcSet$.
Inspired by the result of Hatcher \cite{hatcher1978concordance}, we construct a psc-analog of the Hatcher spectral sequence, that approximates the homotopy groups of the homotopy fibre by the homotopy groups of the space of \emph{psc pseudo isotopies} $C^{\mathrm{psc}}(M \times I^p)$. 
This set consists of all psc block metrics on $M \times [-1,1]^{p+1}$ that agree near each hyperface $M \times \{x_i = \eps\}$, where $\eps = \pm 1$, with the product metric $g_0 \oplus \euclmetric_{\R^p}$ of the base point $g_0 \in \Riem^+(M)$ and the euclidean metric, except near the last front face $M \times \{x_{p+1}=1\}$. 

\begin{thmx}\label{HSS - MainThm}
 There is a first quadrant spectral sequence $E^r_{p,q}$ starting with the homotopy groups $E^1_{pq} = \pi_q\bigl(C^{\mathrm{psc}}(M \times I^p)\bigr)$ and converging to $\pi_{p+q}(\mathrm{hofib}_\bullet)$ if $p+q\geq 2$.
\end{thmx}

\subsubsection*{Outline of this thesis}

A short outline of this thesis is the following one (for a more detailed description, see the introduction of each individual chapter).
Chapter \ref{Chapter - Cubical Sets} discusses the required material on cubical set theory. 
Section \ref{Section - Foundations of Cubical Sets} summarises the required general theory and surveys the existing literature, while Section \ref{Section - Examples of Cubical Sets} presents suitable modifications of the general theory and examples needed for this thesis.
Chapter \ref{Chapter - Cubical Versions of psc} focuses on the construction of $\ConcSet$, its Kan property, i.e. Theorem \ref{KanSet - MainThm}, and the geometric addition.
Chapter \ref{The Operator Concordance Set - Chapter} is the operator-theoretic analog of Chapter \ref{Chapter - Cubical Versions of psc}. There we construct $\InvBlockDirac$, show that it is a Kan set, and that it is a model for real $K$-theory. 
In particular, Theorem \ref{BlockModelKTheory - MainThm} is the main result of Section \ref{Section - Comparison to K-theory}.
In Chapter \ref{Factorisation Indexdiff - Chapter}, we prove that the index difference factors through $\ConcSet$, i.e. Theorem \ref{Factorisation - MainThm}. 
As an application, we have a more in-depth discussion of how this thesis and the cubical set $\ConcSet$ relate to previous works in the field of positive scalar curvature metrics.
In particular, we provide with the proof of Theorem \ref{ConterExampleStrongConcVsIsotopy - MainThm} counterexamples to the strong concordance-implies-isotopy conjecture.
Finally, in Chapter \ref{The PSC Hatcher Spectral Sequence - Chapter}, we construct the psc Hatcher spectral sequence, which proves Theorem \ref{HSS - MainThm}.
To keep this thesis self-contained, we give a thorough discussion of Sobolev spaces and pseudo differential operators in the Appendix.
At the end of this thesis we present a list of symbols that are used in more than one section.

 \chapter{Cubical Set Theory}\label{Chapter - Cubical Sets}

Cubical set theory is the foundational tool from combinatorial topology we use to construct the main objects of this thesis, the concordance set $\ConcSet$ in Chapter \ref{Chapter - Cubical Versions of psc} and the operator concordance set $\InvBlockDirac$ in Chapter \ref{The Operator Concordance Set - Chapter}.
These ``approximation spaces'' of $\Riem^+(M)$ and $KO^{-(\dim M+1)}$ cannot be written in a closed form but need to be constructed; and cubical sets are the tool of our choice.
Although simplicial sets are more popular, cubical sets fit our geometric setup better and this is the reason why we work with them.

We will lay out in Section \ref{Section - Foundations of Cubical Sets} the foundation of cubical sets needed in this thesis.
The main result is Theorem \ref{Equivalence of Homotopy Categories - Thm}.
Roughly speaking, it says that the category of topological spaces and the category of cubical set carry the same information that are preserved under homotopy. 
In other words, if one is interested in homotopy-theoretical questions, then one can work, in theory, equally well in both categories\footnote{In practice, which category is more convenient depends on the question. In this thesis, it will be more convenient to work in the category of cubical sets.}.

In Section \ref{Section - Examples of Cubical Sets}, we will present conceptual examples that will be applied later in this thesis.

\section{Foundations of Cubical Sets}\label{Section - Foundations of Cubical Sets}

The presentation of this section mostly follows the articles of Rosa Antolini \cite{antolini_cubical_2000} and John Jardine \cite{jardine2002cubical}, \cite{jardine2006categorical}.
Different authors use different conventions for cubical sets.
The cubical sets we will use in this thesis have \emph{no connections}.
We made the effort to collect everything we need, because the literature on cubical sets is quite scattered and, in contrast to simplicial sets, there is no standard reference.
Although the original article of Kan \cite{KanAbstractHomotopyI} contains most of the results presented here, most theorems are not proven.
We do not claim any originality for this section; our contribution is merely to bring everything in consistent notation, and present proof-sketches when needed or hint at analogous statements in simplicial set theory from which the proof can be modified straightforwardly.

Throughout this section, we use the convention $I = [-1,1]$ and $\Z_2 = \{\pm 1\}$.

\begin{definition}\label{Box Category - Def}
 The \emph{box category} $\BoxCat$ is the category whose objects are $I^n$ 
 \nomenclature[Cub]{$I$}{unit interval $[-1,1]$}
 for all $n \in \N_0$ and whose morphisms are generated by the following ones:
 \begin{itemize}
  \item The morphism $\CubeIncl{i} \colon I^{n-1} \hookrightarrow I^n$ includes $\eps$ in the $i$-th position of the tuple. 
  \nomenclature[Cub]{$\CubeIncl{i}$}{Injective map $I^{n-1}\hookrightarrow I^n$ that includes $\eps$ into the $i$-th entry}
  Here, $(i,\eps)$ ranges over $\{1,\dots,n\} \times \Z_2$.
  \item The morphism $\CubeProj{i} \colon I^n \rightarrow I^{n-1}$ projects the $i$-th component away. 
  \nomenclature[Cub]{$\CubeProj{i}$}{Linear map $I^n \rightarrow I^{n-1}$ that projects the $i$-th component away.}
  Here, $i$ ranges over $\{1,\dots,n\}$ if $n>1$.
  If $n = 1$, then $\CubeProj{1} \colon I^1 \rightarrow I^0 = \{0\}$ is the zero map.
 \end{itemize}
\end{definition}

The geometric idea is that $\CubeIncl{i}$ identifies $I^{n-1}$ with the face $\{x_i = \eps\}$ of $I^n$ and that $\CubeProj{i}$ projects the $n$-dimensional cube to the $(n-1)$-dimensional cube along the $x_i$-axis. 
All $\CubeIncl{i}$ are the restrictions of affine maps $\CubeIncl{i} \colon \R^{n-1} \rightarrow \R^n$ and all $\CubeProj{i}$ are restrictions of linear maps. 
We will also denote restriction of these maps to other subsets by $\CubeIncl{i}$ and $\CubeProj{i}$. 
More generally, we will denote the insertion of $s \in \R$ into the $i$-th position by $\CubeIncl[s]{i}$.

Combinatorial calculations imply the next lemma.

\begin{lemma}\label{Cocubical Identities - Lemma}  
 The maps $\CubeIncl{i}$ and $\CubeProj{i}$ satisfy the \emph{co-cubical identities}
 \begin{align*}
  \CubeIncl[\omega]{j}\CubeIncl{i} &= \CubeIncl{i}\CubeIncl[\omega]{j-1}, \quad \  \text{ if } i<j, \\
  \CubeProj{j}\CubeProj{i} &= \CubeProj{i}\CubeProj{j+1}, \quad \ \text{ if } i \leq j,  \\
  \CubeProj{j} \CubeIncl{i} &= \begin{cases}
     \CubeIncl{i}\CubeProj{j-1}, &\text{if } i < j, \\
     \id, &\text{if } i = j, \\
     \CubeIncl{i-1}\CubeProj{j}, &\text{if } i > j. 
    \end{cases}
 \end{align*}
\end{lemma}

An important consequence of these identities is the following characterisation of morphisms in the box category.

\begin{lemma}\label{Box Morphism Ordering Lemma}
 For every morphism $\varphi \in \Hom_\BoxCat(I^n,I^m)$ exists an $r \in \N_0$, natural numbers $n \geq j_1 \geq\dots \geq j_{n-r}$, and pairs $(i_1,\eps_{i_1}), \dots, (i_{m-r},\eps_{m-r})$ with $i_1 \leq i_2 \dots \leq i_{m-r}$ such that 
 \begin{equation*}
  \varphi = \CubeIncl[\eps_{i_{m-r}}]{i_{m-r}} \circ \dots \circ \CubeIncl[\eps_{i_1}]{i_1} \circ \CubeProj{j_1} \circ \dots \circ \CubeProj{j_{n-r}}.
 \end{equation*}
\end{lemma}

\begin{definition}\label{Cubical Sets categorial - Def}
 A \emph{cubical set} $X_\bullet$ is a contravariant functor $X \colon \BoxCat \rightarrow \mathbf{Set}$. 
 A \emph{cubical map} $f_\bullet \colon X_\bullet \rightarrow Y_\bullet$ is a natural transformation between the contravariant functors $X_\bullet$ and $Y_\bullet$.
\end{definition}

Unpacked, Definition \ref{Cubical Sets categorial - Def} reads as follows.

\begin{definition}\label{Cubical Sets elementary - Def}
 A cubical set $X_\bullet$ is a sequence of sets $(X_n)_{n \in \N_0}$ together with \emph{connecting maps} $\face{i} \deff X(\CubeIncl{i}) \colon X_{n} \rightarrow X_{n-1}$ 
 \nomenclature[Cub]{$\face{i}$}{face map of a cubical set}
 for all $(i,\eps) \in \{1,\dots,n\} \times \Z_2$ with $n\geq 1$, and $\degen{i} \deff X(\CubeProj{i}) \colon X_n \rightarrow X_{n+1}$ 
 \nomenclature[Cub]{$\degen{i}$}{degeneracy map of a cubical set}
 for all $i \in \{1,\dots,n+1\}$ with $n\geq 0$ that satisfies the \emph{cubical identities}
 \begin{align*}
  \face{i} \face[\omega]{j} &= \face[\omega]{j-1}\face{i}, \quad \ \text{ if } i<j, \\
  \degen{i}\degen{j} &= \degen{j+1} \degen{i}, \quad \ \text{ if } i \leq j,  \\
  \face{i} \degen{j} &= \begin{cases}
     \degen{j-1}\face{i}, &\text{if } i < j, \\
     \id, &\text{if } i = j, \\
     \degen{j} \face{i-1}, &\text{if } i > j. 
    \end{cases}
 \end{align*}
 A \emph{cubical map} $f_\bullet \colon X_\bullet \rightarrow Y_\bullet$ is a sequence of maps $f_n \colon X_n \rightarrow Y_n$ for all $n \in \N_0$ that commute with all connecting maps.  
\end{definition}

We refer to an element of $X_n$ as an $n$-cube of $X_\bullet$ and sometimes call $X_n$ the set of all $n$-cubes (of $X_\bullet$).
All cubical sets together with cubical maps form a category, which we denote by $\cSet$.

\begin{example}\label{The Singular Set - Example}
 Let $X$ be a topological space. The \emph{singular set} $S_\bullet(X)$ 
 \nomenclature[Cub]{$S_\bullet(X)$}{The singular set of a topological space $X$}
 is the cubical set consisting of 
 \begin{equation*}
   S_n(X) \deff \{s \colon I^n \rightarrow X \text{ continuous }\}
 \end{equation*}
 and whose connecting maps are
 \begin{equation*}
  \face{i} s = s \circ \CubeIncl{i} \qquad \text{ and } \qquad \degen{i}s = s \circ \CubeProj{i}.
\end{equation*}  
 In other words, the connecting maps are given by pre-composing a singular cube with a morphism of the box category. 
 
 A continuous map $f \colon X \rightarrow Y$ induces a map of cubical sets by post-composing it to singular cubes.
 In this way, the singular set is a functor $S_\bullet \colon \mathbf{Top} \rightarrow \cSet$.   
\end{example}

\begin{example}\label{Cubical Subset - Example}
 Let $X_\bullet$ be a cubical set and $(A_n)_{n\in\N_0}$ a sequence of subsets of $(X_n)_{n\in \N_0}$ to which the connecting maps restrict. 
 These restrictions turn the sequence of subsets into a cubical set $A_\bullet$ such that the level-wise inclusion is a cubical map. 
 We call $A_\bullet$ a cubical subset of $X_\bullet$.
 For example, if $Y \subseteq X$ is a topological subspace, then $S_\bullet(Y)$ is a cubical subset of $S_\bullet(X)$.
 
 Let $f_\bullet \colon X_\bullet \rightarrow Y_\bullet$ be a cubical map and let $A_\bullet \subseteq X_\bullet$ and $B_\bullet \subseteq Y_\bullet$ be cubical subsets. The preimage $f_\bullet^{-1}(B_\bullet)$ and the image $f_\bullet(A_\bullet)$ are cubical subsets of $X_\bullet$ and $Y_\bullet$, respectively.
 Furthermore, intersections of cubical subsets are cubical subsets and unions of cubical subsets are cubical subsets. 
\end{example}

In this thesis, we use cubical sets as a blueprint to construct topological spaces by gluing cubes together. 
This interpretation begs the question how to model the standard cubes in this manner.

\begin{example}\label{Standard Cube - Example}
 The \emph{(cubical)} $n$\emph{-cube} $\StandCube{n}$ 
 \nomenclature[Cub]{$\StandCube{n}$}{The (cubical) $n$-cube}
 is the cubical set given by the functor
 that sends an object $I^m$ to $\mathrm{Hom}_\BoxCat(I^m,I^n)$ and a morphism $\varphi \colon I^k \rightarrow I^m$ to its pre-composition $\varphi^\ast \colon \mathrm{Hom}_\BoxCat(I^m,I^n) \rightarrow \mathrm{Hom}_\BoxCat(I^k,I^n)$.
 We will often just write $\ast$ for $\StandCube{0}$.
 A morphism $\psi \colon I^n \rightarrow I^r$ in $\Box$ induces a cubical map $\psi_\ast \colon \StandCube{n} \rightarrow \StandCube{r}$ by post-composition. 
 In this way, $\StandCube{\placeholder} \colon \Box \rightarrow \cSet$ is a covariant functor.
\end{example}

The cubical cubes can be considered as the building blocks of cubical sets. 
The next lemma, which is an immediate application of the Yoneda lemma, makes this statement more precise.

\begin{lemma}\label{Building Block Lemma}
 The evaluation at $I^n$ is a natural bijection 
 \begin{equation*}
  \mathrm{ev}_{I^n} \colon \mathrm{Hom}_\cSet(\StandCube{n},X_\bullet) \rightarrow X_n
 \end{equation*}  
 for every cubical set $X_\bullet$.
\end{lemma}

\begin{example}\label{Basepoint - Example}
 A base point $x_0$ of a cubical set $X_0$ is a cubical map $x_0 \colon \StandCube{0} \rightarrow X_\bullet$. 
 Using Lemma \ref{Building Block Lemma}, we may identify the cubical map with its image $x_0 \in X_0$. 
 The image of the map is the cubical subset $\{x_0\}_\bullet$ given by $\{x_0\}_n = \{\degen{n-1}\dots \degen{1} x_0\}$.
\end{example}

For the sake of readability, we will denote $\degen{n-1}\dots \degen{0}x_0$ again by $x_0$ if $x_0$ serves as a base point.  

\begin{example}\label{Cubical n-horn - Example}
 For $(i,\eps) \in \{1,\dots,n\} \times \Z_2$, the \emph{cubical $n$-horn}
 \nomenclature[Cub]{$\CubeBox$}{Cubical $n$-horn open towards $(i,\eps)$.}
 $\CubeBox$ is the cubical subset of $\StandCube{n}$ given by
 \begin{equation*}
  \CubeBox \deff \bigcup \Bigl\{\CubeIncl[\omega]{j}(\StandCube{n-1}) \, : \, (j,\omega) \in \{1,\dots,n\}\times \Z_2 \setminus \{(i,\eps)\} \Bigr\}.
\end{equation*}  
 It is the combinatorial model of the boundary of $I^n$ with the open face $\CubeIncl{i}(I^{n-1})^\circ$ removed.
 In view of Lemma \ref{Building Block Lemma}, every cubical map $f_\bullet \colon \CubeBox \rightarrow X_\bullet$ uniquely determines a set of elements
 \begin{equation*}
  \{x_{(j,\omega)} \in X_{n-1} \, : \, \face[\omega]{j}x_{(k,\eta)} = \face[\eta]{k-1}x_{(j,\omega)} \text{ for } j<k, \, (j,\omega), (k,\eta) \neq (i,\eps)\},
\end{equation*}  
where $(j,\omega)$ ranges over $\{1,\dots,n\} \times \Z_2 \setminus \{(i,\eps)\}$.
Conversely, each collection of this type uniquely determines a cubical map $f_\bullet \colon \CubeBox \rightarrow X_\bullet$.
\end{example}

\begin{example}\label{The combinatorial sphere}
 The \emph{cubical sphere} $\partial \StandCube{n}$ is the cubical set given by
 \begin{equation*}
  \partial \StandCube{n} \deff \bigcup \Bigl\{\CubeIncl[\omega]{j}(\StandCube{n-1}) \, : \, (j,\omega) \in \{1,\dots,n\} \times \Z_2 \Bigr\}.
 \end{equation*}
 \nomenclature[Cub]{$\partial \StandCube{n}$}{Cubical sphere}
 It is a combinatorial model for the boundary $\partial I^n = S^{n-1}$.
 In view of Lemma \ref{Building Block Lemma}, every cubical map $f_\bullet \colon \partial \StandCube{n} \rightarrow X_\bullet$ uniquely determines a set of elements 
 \begin{equation*}
  \{x_{(j,\omega)} \in X_{n-1} \, : \, \face[\omega]{j}x_{(k,\eta)} = \face[\eta]{k-1}x_{(j,\omega)} \text{ for } j<k\},
\end{equation*}  
where $(j,\omega)$ ranges over $\{1,\dots,n\} \times \Z_2$.
Conversely, each collection of this type uniquely determines a cubical map $f_\bullet \colon \partial \StandCube{n} \rightarrow X_\bullet$.
\end{example}

If we think of a cubical set as a blueprint to construct a topological space, then we eventually want to construct this topological space from it.
This is what the geometric realisation functor does. 

\begin{definition}\label{Geometric Realisation - Def}
 Let $X_\bullet$ be a cubical set. Equip each set $X_n$ with the discrete topology.
 Its \emph{geometric realisation} is the topological space obtained by endowing the set
 \begin{equation*}
   |X_\bullet| = \left( \bigsqcup_{n=0} I^n \times X_n \right)\biggl\slash \sim 
 \end{equation*}
 \nomenclature[Cub]{|X_\bullet|}{Geometric Realisation of a cubical set}
 with the quotient topology.
 Here, $\sim$ is the equivalence relation that is generated by $(\varphi(v),x) \sim (v,X(\varphi)(x))$ for all $v \in I^m$, $x \in X_n$, and $\varphi \in \mathrm{Hom}_\BoxCat(I^m,I^n)$.
 
 The geometric realisation of a cubical map $f_\bullet \colon X_\bullet \rightarrow Y_\bullet$ is the continuous map $|f_\bullet| \colon |X_\bullet| \rightarrow |Y_\bullet|$ that is induced by the sequence of continuous maps $\id \times f_n \colon I^n\times X_n \rightarrow I^n \times Y_n$.  
\end{definition} 

The proof of the following lemma can be found in \cite{jardine2006categorical}*{p.89}.

\begin{lemma}\label{Adjointness Geo/Sing - Lemma} 
 The geometric realisation is a functor $| \placeholder | \colon \cSet \rightarrow \mathbf{Top}$.
 It is the left adjoint to the singular set functor
 \begin{equation*}
   |\placeholder | \colon \cSet \rightleftarrows \mathbf{Top} :\! S_\bullet.
 \end{equation*}
\end{lemma}

\begin{example}
 As expected, we have $|\StandCube{n}| = I^n$ and $\partial |\!\CubeBox\!| = I^n\setminus \CubeIncl{i}(I^{n-1})^\circ$.
 However, only few continuous maps $I^n \rightarrow I^m$ are geometric realisations of cubical maps $\StandCube{n} \rightarrow \StandCube{m}$, namely, the morphisms of the box category. 
\end{example}

Let us now discuss products of cubical sets. 
The (categorical) product of two cubical set $X_\bullet$ and $Y_\bullet$ is the cubical set $(X \times Y)_\bullet$ whose set of $n$-cubes is given by
\begin{equation*}
    (X \times Y)_n = X_n \times Y_n
\end{equation*}
\nomenclature[Cub]{$X \otimes Y$}{Reduced product of two cubical sets}
and whose connecting maps are applied component-wise.
In contrast to simplicial sets, the categorical product does not commute with the geometric realisation.
Indeed, Jardine shows in \cite{jardine2002cubical} that $|\StandCube{1} \times \StandCube{1}| \simeq S^1$ instead of $I^1 \times I^1$.
This is the reason why Kan introduced in \cite{KanAbstractHomotopyI} a different product that fits the geometric expectations.

\begin{definition}\label{Reduced Product}
 Let $X_\bullet$ and $Y_\bullet$ be cubical sets. 
 Define the \emph{reduced product} $(X \otimes Y)_\bullet$ to be the cubical set whose set of $n$-cubes is the set of equivalence classes
 \begin{equation*}
  (X \otimes Y)_\bullet = \bigcup_{k=0}^n X_k \times Y_{n-k} \biggl / \sim, 
 \end{equation*}  
  of the equivalence relation generated by $(\degen{k+1}x,y) \sim (x,\degen{1}y)$ for all $x \in X_{k}$, $y \in Y_{n-k}$, and $0 \leq k \leq n$.
  The connecting maps are given by
  \begin{align*}
      \face{i}(x,y) = \begin{cases}
       (\face{i}x,y), & \text{if } i\leq k, \\
       (x,\face{i-k}y), & \text{if } i>k,
      \end{cases} \quad \text{ and } \quad 
      \degen{i}(x,y) = \begin{cases}
        (\degen{i}x,y), & \text{if } i \leq k, \\
        (x,\degen{i-k}y), & \text{if } i > k.
      \end{cases}
  \end{align*}
\end{definition} 

\begin{example}
 As expected, the cubical sets $(\StandCube[ ]{m} \otimes \StandCube[ ]{n})_\bullet$ and $\StandCube{m+n}$ are isomorphic.
 An isomorphism is given by the unique cubical map that sends the element $\id_{I^{m+n}} \in \StandCube{m+n}$ to the equivalence class of $(\id_{I^m},\id_{I^n})$.
 Its inverse is the cubical map induced by
 \begin{align*}
     \StandCube[k]{m} \times \StandCube[l]{n} &\rightarrow \StandCube[k+l]{m+n}, \\
     (\varphi,\psi) &\mapsto (\varphi \circ \CubeProj{k+1} \circ \dots \circ \CubeProj{k+l},\psi \circ \CubeProj{1}^{\circ k}). 
 \end{align*}
 
 To check that both maps are indeed inverse to each other, it is enough to check that they send $\id_{I^{m+n}}$ and $(\id_{I^m},\id_{I^n})$ to each other, as all other elements in those cubical sets can be obtained either from $\id_{I^{m+n}}$ or $(\id_{I^{m}},\id_{I^n})$ using a composition of connecting maps.
\end{example}

\begin{lemma}
 For all cubical sets $X_\bullet$, $Y_\bullet$, the geometric realisation $|(X \otimes Y)_\bullet|$ is homeomorphic to $|X_\bullet| \times |Y_\bullet|$ if the latter is equipped with the CW-topology.
\end{lemma}
\begin{proof}
 For the proof we refer to Jardine \cite{jardine2006categorical}*{p.92},
 who introduced $\otimes$ in a different but equivalent manner.
\end{proof}

We can use the reduced product to model homotopies combinatorially.

\begin{definition}\label{Combinatorial Homotopy - Def}
 Let $f_\bullet, g_\bullet \colon X_\bullet \rightarrow Y_\bullet$ be two cubical maps and let $A_\bullet \subseteq X_\bullet$ be a cubical subset. 
 A \emph{homotopy} between $f_\bullet$ and $g_\bullet$ that is \emph{stationary} on $A_\bullet$ is a cubical map $H_\bullet \colon (\StandCube[ ]{1} \otimes X)_\bullet \rightarrow Y_\bullet$ such that the following diagram commutes:
 \begin{equation*}
   \xymatrix{ (\StandCube[ ]{0} \otimes X)_\bullet \ar@{=}[r] \ar@<-0.5ex>[d]_{\CubeIncl[1]{0}} \ar@<0.5ex>[d]^{\CubeIncl[-1]{0}} 
   & X_\bullet  \ar@<0.5ex>[rr]^{f_\bullet} \ar@<-0.5ex>[rr]_{g_\bullet} &
   & Y_\bullet \, \\
   (\StandCube[ ]{1} \otimes X)_\bullet \ar@/_1pc/[rrru]_{H_\bullet}  & &  \\
   (\StandCube[ ]{1} \otimes A)_\bullet \ar[u]_{\mathrm{incl}} \ar[r]^{{\CubeProj{0}}_\ast} & (\StandCube[ ]{0} \otimes A)_\bullet \ar@{=}[rr] && \:A_\bullet \ar[uu]_{H_\bullet|_{A_\bullet}}.}
 \end{equation*}
 The maps $f_\bullet$ and $g_\bullet$ are called \emph{homotopic relative to} $A_\bullet$, if there is a homotopy between them that is stationary on $A_\bullet$.
\end{definition}
The diagram in the previous definitons are, in fact, two diagrams in one. 
The upper triangle describes a homotopy between two given cubical maps, the lower square describes what it means to be stationary. 

In light of Lemma \ref{Building Block Lemma}, we may say that two elements $x_{-1}, x_1 \in X_n$ are homotopic if and only if their representing maps $f_{x_i} \colon \StandCube{n} \rightarrow X_\bullet$ agree on $\partial \StandCube{n}$ and are homotopic relative to $\partial \StandCube{n}$.

For concrete calculations, the following equivalent description in terms of elements turns out to be more convenient.

\begin{definition}\label{Element Homotopy - Def}
 Two elements $x_{-1}, x_1 \in X_n$ are \emph{homotopic} if $\face{i}x_{-1} = \face{i}x_1$ for all $(i,\eps) \in \{1,\dots,n\} \times \Z_2$ and if there is an $h \in X_{n+1}$ that satisfies
 \begin{equation*}
  \face[\pm 1]{1}h = x_{\pm 1} \qquad \text{ and } \qquad \face{i}h = \degen{1}\face{i-1}x_{\pm 1}.
 \end{equation*}
 We call such elements $h$ a \emph{homotopy} between $x_{-1}$ and $x_1$.
\end{definition}

Being homotopic is not an equivalence relation on $\mathrm{Hom}_\cSet(X_\bullet,Y_\bullet)$ in general. 
However, it is an equivalence relation if the target satisfy an extension condition that was introduced by Kan in \cite{KanAbstractHomotopyI}.

\begin{definition}\label{Kan condition and Kan fibration - Def}
 A cubical set $K_\bullet$ is a \emph{Kan set} or satisfies the \emph{Kan extension property} if every commutative outer square has a dotted filler
 \begin{equation*}
  \xymatrix@C+1em{\CubeBox \ar[d] \ar[r] & K_\bullet \ar[d] \\ \StandCube{n} \ar[r] \ar@{.>}[ru] & \:\ast. }
 \end{equation*}
 
 More generally, a cubical map $p_\bullet \colon E_\bullet \rightarrow B_\bullet$ is a \emph{Kan fibration} if every commutative outer square has a dotted filler
 \begin{equation*}
  \xymatrix@C+1em{\CubeBox \ar[d] \ar[r] & E_\bullet \ar[d]^{p_\bullet} \\ \StandCube{n} \ar[r] \ar@{.>}[ru] & \:B_\bullet .  }
 \end{equation*}
 
 A cubical set $X_\bullet$ is \emph{combinatorially contractible} if every commutative outer square has a dotted filler
\begin{equation*}
  \xymatrix@C+1em{\partial \StandCube{n}  \ar[d] \ar[r] & X_\bullet \ar[d] \\ \StandCube{n} \ar[r] \ar@{.>}[ru] & \:\ast.}
 \end{equation*}
\end{definition}

The  name 'combinatorially contractible' is motivated by Proposition \ref{CombContract - Prop}. 
In light of the Yoneda lemma (see also Example \ref{Cubical n-horn - Example} and Example \ref{Standard Cube - Example}) we can give an element-wise description of Kan sets and Kan fibrations.

\begin{lemma}\label{Elementwise Kan - Lemma}
 A cubical set $K_\bullet$ is a Kan set if and only if, for all $n \in \N$ and for each collection
 \begin{equation*}
  \bigl\{x_{(j,\omega)} \in K_{n-1} \, : \, \face[\omega]{j}x_{(k,\eta)} = \face[\eta]{k-1}x_{(j,\omega)} \text{ for } j<k; (j,\omega),(k,\eta) \neq (i,\eps) \bigr\},
\end{equation*}  
there is an element $x \in K_n$ such that $\face[\omega]{j} x = x_{(j,\omega)}$.
Here, $(j,\omega)$ ranges over $\{1,\dots,n\} \times \Z_2 \setminus \{(i,\eps)\}$.

More generally, a map $p_\bullet \colon E_\bullet \rightarrow B_\bullet$ is a \emph{Kan fibration} if and only if, for each collection
\begin{equation*}
  \bigl\{e_{(j,\omega)} \in E_{n-1} \, : \, \face[\omega]{j}e_{(k,\eta)} = \face[\eta]{k-1}e_{(j,\omega)} \text{ for } j<k; (j,\omega),(k,\eta) \neq (i,\eps) \bigr\}
\end{equation*}
and each $b \in B_n$ with $p_{n-1}(e_{(j,\omega)}) = \face[\omega]{j}b$, there is an $e \in E_n$  that satisfies $\face[\omega]{j}e = e_{(j,\omega)}$ and $p_n(e) = b$.
Again, $(j,\omega)$ ranges over $\{1,\dots,n\} \times \Z_2 \setminus \{(i,\eps)\}$.

A cubical set $X_\bullet$ is combinatorially contractible if and only if, for all $n \in \N$ and each collection 
 \begin{equation*}
  \bigl\{x_{(j,\omega)} \in X_{n-1} \, : \, \face[\omega]{j}x_{(k,\eta)} = \face[\eta]{k-1}x_{(j,\omega)} \text{ for } j<k \bigr\},
 \end{equation*}
 there is an element $x \in X_{n}$ such that $\face[\omega]{j}x = x_{(j,\omega)}$. 
 Here, $(j,\omega)$ ranges over $\{1,\dots,n\} \times \Z_2$.
\end{lemma}

There are further formal consequences that can be proven analogously to their simplicial counterparts.

\begin{lemma}\label{KanFormalImplication - Lemma}
 Let $p_\bullet \colon E_\bullet \rightarrow B_\bullet$ be a Kan fibration and let $b_0$ be a base point of $B_\bullet$. 
  \begin{enumerate}
   \item[(i)] If the base $B_\bullet$ is a Kan set, then the total set $E_\bullet$ is a Kan set.
   \item[(ii)] The fibre $F_\bullet \deff p^{-1}_\bullet(\{b_0\})$ is a Kan set.
  \end{enumerate}
\end{lemma}
\begin{proof}
 To prove the first claim, we consider the commutative diagram
 \begin{equation*} 
  \xymatrix@C+2em{ & \CubeBox \ar@<0.5ex>[rrr] \ar@{->}'[d][dd] &&& E_\bullet \ar[ld]_{p_\bullet} \ar[dd] \\
  \CubeBox \ar[rrr] \ar@{=}[ru] \ar[dd]&&& B_\bullet \ar[dd] &\\
  &\StandCube{n} \ar@<-0.5ex>@/^2pc/@{.>}[rrruu] \ar@{-}'[r][rr] && \ar[r] & b_0\\
  \StandCube{n} \ar@{=}[ru] \ar@/_1.7pc/@{-->}[rrruu] \ar[rrr] &&& \:b_0. \ar@{=}[ru] &}
 \end{equation*}
 The dashed arrow exists because $B_\bullet$ is Kan. Since $p_\bullet$ is a Kan fibration, we can lift the dashed arrow, which gives the dotted arrow.
 
 The proof of the second claim is captured in the diagram
 \begin{equation*}
  \xymatrix{ \CubeBox \ar[r] \ar[d] & F_\bullet \ar@{^{(}->}[r] \ar[d] & E_\bullet \ar[d]^{p_\bullet} \\ \StandCube{n} \ar[r] \ar@{.>}[ru] \ar@{-->}[rru] & b_0 \ar@{^{(}->}[r] & \:B_\bullet.}
 \end{equation*}
 The dashed arrow exists because $p_\bullet$ is a Kan fibration. 
 By construction, it takes values in the fibre, so the dotted arrow exists.
\end{proof}

\begin{example}
 The singular set $S_\bullet(X)$ of a topological space $X$ is a Kan set.
 Indeed, each singular cubical $n$-horn $f_\bullet \colon \CubeBox \rightarrow S_\bullet(X)$ yields a continuous map $f \colon |\!\CubeBox\!| \rightarrow X$. Clearly, $|\!\CubeBox\!| = \partial I^n \setminus (\CubeIncl{i}(I^{n-1}))^\circ$ is a deformation retract of $I^n$. If $r \colon I^n \rightarrow |\!\CubeBox\!|$ is a retract, then $f \circ r \in S_n(X)$ is an extension of the given singular $n$-horn $f_\bullet$. 
 
 With a similar argument, one can show that $S_\bullet(p) \colon S_\bullet(E) \rightarrow S_\bullet(B)$ is a Kan fibration if $p \colon E \rightarrow B$ is a Serre fibration.  
\end{example}

\begin{lemma}
 For all Kan sets $K_\bullet$ and all cubical sets $X_\bullet$, being homotopic is an equivalence relations on $\mathrm{Hom}_\cSet(X_\bullet, K_\bullet)$.
\end{lemma}
\begin{proof}
 The proof is analogous to the proof for the simplicial counterpart sketched in \cite{curtis1971simplicial}*{Lemma 1.15ff}. 
\end{proof}

We define the $n$-th homotopy group of a cubical set to be the $n$-th homotopy group of its geometric realisation.
However, if the cubical set is a Kan set, then there is a purely combinatorial description of the homotopy groups. 
The following notation is taken from \cite{antolini_cubical_2000}.

\begin{definition}\label{combitarial loop space - Def}
 Let $k_0$ be a base point of the Kan set $K_\bullet$.
 Recall the convention that we abbreviate $\degen{n-1} \dots \degen{0}k_0 \in K_n$ to $k_0$.
 For $n \geq 1$, define the set of $n$-fold loops as
 \begin{equation*}
  \Omega_{k_0}^n K_\bullet \deff \bigl\{k \in K_n \, : \, \face{i} k = k_0 \in K_{n-1} \text{ for all } (i,\eps) \in \{1,\dots,n\}\times \Z_2\bigr\}.
 \end{equation*}
 \nomenclature[Cub]{$\Omega_{k_0}^n K_\bullet$}{Set of $n$-fold loops in $K_\bullet$}
\end{definition}

There is also a version for pairs of Kan sets, see \cite{antolini_cubical_2000}.

\begin{definition}\label{Homotopy groups - Def}
 For every cubical set $X_\bullet$, define $\pi_0(X_\bullet)$ to be the set of all (combinatorial) path components
 \begin{equation*}
  \pi_0(X_\bullet) \deff X_0/\sim_{\text{homotopy}}\!.
 \end{equation*}
 If $K_\bullet$ is a Kan set with base point $k_0$, then define the higher (combinatorial) homotopy groups $\pi_n(K_\bullet, k_0)$ via
 \begin{equation*}
  \pi_n(K_\bullet,k_0) \deff \Omega^n_{k_0}(K_\bullet)/\sim_{\text{homotopy}},
 \end{equation*} 
 \nomenclature[Cub]{$\pi_j(K_\bullet)$}{(Combinatorial) Homotopy group of the Kan set $K_\bullet$}
 where $\sim_{\text{homotopy}}$ denotes the element-wise homotopy relation from Definition \label{Element Homotopy - Def}.
\end{definition}

We can use the Kan condition to put a purely combinatorial group structure on $\pi_n(K_\bullet)$ for all $n > 0$. 
For two given elements $[a], [b] \in \pi_n(K_\bullet,k_0)$, we pick representatives $a, b \in \Omega_{k_0}^nK_\bullet$, respectively.
These two elements form an $\CubeBox[1][1][n+1]$-horn with $k_{(1,-1)}=a$, $k_{(2,1)}=b$, and $k_{(j,\omega)} = k_0$ for all other $(j,\omega) \in \{1,\dots,n+1\} \times \Z_2 \setminus (1,1)$.
By the Kan extension property, this cubical horn has a filler $x \in K_{n+1}$, and the cubical identities imply that $\face[1]{1}x \in \Omega^n_{k_0}K_\bullet$.
We define $[a]+[b] \deff [\face[1]{1}x]$.
\begin{lemma}\cite{KanAbstractHomotopyI}*{Theorem 7}\label{CombinatorialAddition - Lemma}
 This operation is independent of all choices and turns $\pi_n(K_\bullet,k_0)$ into a group for all $n \geq 1$.
\end{lemma}
There are also relative homotopy groups and connecting maps $\partial_n \colon \pi_n(K_\bullet, L_\bullet) \rightarrow \pi_{n-1}(L_\bullet)$ that fit into a long exact sequence, see  \cite{antolini_cubical_2000} for details. 
Furthermore, different choices of base points within the same path component yield isomorphic groups, see \cite{KanAbstractHomotopyI}*{Theorem 8} (the isomorphisms are not canonical, though). 

Adapting the straightforward (but cumbersome) proof of \cite{curtis1971simplicial}*{Theorem 2.8}
to the cubical setup, we get the following fibration result.

\begin{lemma}\label{CubicalLESHmptyGroups - Lemma}
 Let $E_\bullet \rightarrow B_\bullet$ be a Kan fibration, let $b_0$ be a base point of the Kan set $B_\bullet$, and let $F_\bullet = p_\bullet^{-1}(b_0)$ be the fibre.
 Then, for all $n \geq 1$, there are maps of pointed sets $\partial_n \colon \pi_n(E_\bullet) \rightarrow \pi_{n-1}(F_\bullet)$ that are group homomorphisms if $n\geq 2$, and that give rise to a long exact sequence of homotopy groups
 \begin{equation*}
  \xymatrix{\dotsb \ar[r] & \pi_n(F_\bullet) \ar[r] & \pi_n(E_\bullet) \ar[r] & \pi_n(B_\bullet) \ar[r] & \pi_{n-1}(F_\bullet) \ar[r] & \dotsb .}
\end{equation*}  
\end{lemma}

Lastly, we would like to know that cubical sets carry the same homotopy information within them as topological spaces. 

\begin{theorem}\label{Unit Weak Equivalence - Thm}\cite{antolini_cubical_2000}*{Theorem 2.9}
 Let $K_\bullet$ be a Kan set. The  unit of the adjunction in Lemma \ref{Adjointness Geo/Sing - Lemma}, 
 \begin{equation*}
    \psi_\bullet \colon K_\bullet \rightarrow S_\bullet(|K_\bullet|) \quad \text{given by} \quad k \mapsto \bigl( x \mapsto [x,k] \bigr),   
 \end{equation*}
 induces isomorphisms on all homotopy groups.
\end{theorem} 

Although it is easy to prove that the homotopy groups $\pi_n(S_\bullet(X))$ and $\pi_n(X)$ are isomorphic, the theorem itself is quite deep. 
For example, it implies that every cubical map $\StandCube{n} \rightarrow S_\bullet(|K_\bullet|)$ that is constant on the boundary is homotopic relative to its boundary to a cubical map $\StandCube{n} \rightarrow K_\bullet$.
In \cite{antolini_cubical_2000}, the theorem is not proven.
Instead, the source refers to the book of May \cite{may1992simplicial} by analogy. 
We will also not prove this theorem in detail, but guide through the existent literature that can be used as a blue-print to prove this theorem following the strategy of Milnor \cite{milnor1957geometric}.

\begin{proof}
 We first show that $K_\bullet$ and $S_\bullet(|K_\bullet|)$ have the same homology groups.
 The proof of \cite{milnor1957geometric}*{Lemma 5} carries over without change
 if we use homology groups for cubical sets whose definition is along the lines in \cite{dubrovin1990modern}. 
 
 The next step is to prove that $K_\bullet$ and $S_\bullet(|K_\bullet|)$ have isomorphic fundamental groups.
 Following \cite{milnor1957geometric}, first observe that each cubical Kan set $K_\bullet$ has a minimal subcomplex $M_\bullet$, which is, by definition, a cubical subset that satisfies the Kan extension property, that is a (combinatorial) strong deformation retract of $K_\bullet$, and in that every element $x \in M_n$ is the unique representative of its homotopy class\footnote{homotopy class in the sense of Definition \ref{Element Homotopy - Def}} for all $n \in \N$.
 The existence proof can be carried out, for example, along the line of \cite{curtis1971simplicial}*{Prop 1.23}.
 By minimality, $M_\bullet$ consists of only one vertex per path component, so the geometric realisation $|M_\bullet|$ is a CW complex with a single zero cell per path component.
 The fundamental group of such a CW complexes has one generator per edge and one relation per attaching map.
 But the combinatorial group structure of $\pi_1(M_\bullet)$ is a just a reformulation of these relations, so the unit induces an isomorphism $\pi_1(M_\bullet) \rightarrow \pi_1(S_\bullet(|M_\bullet|)) = \pi_1(|M_\bullet|)$ because it sends each element of the combinatorial fundamental group to an edge.
 
 Since the unit $\psi_\bullet$ is natural with respect to cubical maps and since the inclusion $M_\bullet \hookrightarrow K_\bullet$ is a homotopy equivalence, the unit for $K_\bullet$ induces an isomorphism on the fundamental group because of the commutative diagram
 \begin{equation*}
     \xymatrix{\pi_1(K_\bullet) \ar[rr]^-{\pi_1(\psi_\bullet)} && \pi_1(S_\bullet(|K_\bullet|)) \\
     \pi_1(M_\bullet) \ar[u]^{\pi_1(\mathrm{incl}_\bullet)}_-\cong \ar[rr]^\cong_-{\pi_1(\psi_\bullet)} && \:\pi_1(S_\bullet(|M_\bullet|)) \ar[u]_{\pi_1(S_\bullet(|\mathrm{incl}_\bullet|))}^{\cong}. }
 \end{equation*}
 
 We prove that the unit induces an isomorphism on higher homotopy groups by applying a relative Hurewicz theorem in the category of cubical sets to the pair of universal coverings $(S_\bullet(|\tilde{K}_\bullet|), \tilde{K}_\bullet)$ to show that all relative homotopy groups vanish.
 Note that the homology of this pair of cubical sets vanish by an argument analogous to the one that simplicial and singular homology of a simplicial set\footnote{Hatcher speaks of $\Delta$-sets or semisimplicial complexes. However, these are nothing but the geometric realisation of a simplicial set.} agree, see \cite{hatcher2001algebraic}*{Chapter 2.1}.
 
 The relative Hurewicz theorem for cubical sets can be proven as follows. 
 First, one establishes a cubical version of the homotopy addition theorem as in \cite{curtis1971simplicial}*{Theorem 2.4}
 and uses it to prove an absolute Hurewicz theorem for cubical sets as in \cite{curtis1971simplicial}*{Theorem 3.12}. 
 Then the relative Hurewicz theorem can be derived by modifying the proof of the relative Hurewicz theorem in \cite{stocker2013algebraische}.
\end{proof}  

The main theorem of this chapter is the cubical extension theorem, which was proven in \cite{antolini_cubical_2000}.
The key ingredient of its proof is Theorem \ref{Unit Weak Equivalence - Thm}. 

\begin{theorem}[Cubical Approximation Theorem] \cite{antolini_cubical_2000}*{Theorem 2.11}

\noindent Let $(X_\bullet, A_\bullet)$ be a pair of cubical sets and let $K_\bullet$ be a Kan set.
 For every continuous map $g \colon |X_\bullet| \rightarrow |K_\bullet|$ whose restriction to $|A_\bullet|$ is the geometric realisation of a cubical map $f_\bullet \colon A_\bullet \rightarrow K_\bullet$, there is a cubical map $F_\bullet \colon X_\bullet \rightarrow K_\bullet$ such that $F_\bullet |_{A_\bullet} = f_\bullet$ and $|F_\bullet|$ is homotopic to $g$ relative to $|A_\bullet|$.
\end{theorem}

A formal consequence of this (relative) homotopy extension theorem is that the (naive) homotopy categories of cubical sets and CW complexes are isomorphic. 
Informally speaking, this means that if we are only interested in homotopy-theoretic question, then it does not matter in which category we work.

\begin{theorem}\label{Equivalence of Homotopy Categories - Thm} {$ $}
 \begin{enumerate}
  \item For all CW complexes $X$, the co-unit of the adjunction between topological spaces and cubical sets
  \begin{equation*}
     \phi \colon |S_\bullet X| \rightarrow X, \quad \text{given by } \quad [t,\sigma] \mapsto \sigma(t),
  \end{equation*}
  is a homotopy equivalence.
  \item If $K_\bullet$ is a Kan set and $X_\bullet$ is an arbitrary cubical set, then the geometric realisation induces a bijection  between the homotopy classes of maps
  \begin{equation*}
    |\placeholder | \colon [X_\bullet; K_\bullet] \rightarrow [|X_\bullet|; |K_\bullet|].
  \end{equation*}
  \item The singular set functor induces a bijection between the homotopy classes of maps
  \begin{equation*}
    S_\bullet \colon [X,Y] \rightarrow [S_\bullet(X),S_\bullet(Y)]
  \end{equation*}
  for all CW complexes.
 \end{enumerate}
\end{theorem}
\begin{proof}
 To prove the first claim, observe that, by definition, the co-unit $\phi$ is a weak homotopy equivalence if and only if $S_\bullet(\phi) \colon S_\bullet(|S_\bullet X|) \rightarrow S_\bullet X$ induces an isomorphism on all (combinatorial) homotopy groups.
 From the decomoposition
 \begin{equation*}
     \xymatrix{ S_\bullet X \ar[rr]^-{\psi_\bullet} \ar@/_1.5pc/[rrrr]_{\id} && S_\bullet(|S_\bullet X|) \ar[rr]^-{S_\bullet \phi} && S_\bullet X,}
 \end{equation*}
 we deduce that $S_\bullet(\phi)$ induces an isomorphism on all (combinatorial) homotopy groups because $\psi_\bullet$ and the identity do it.
 By Whitehead's theorem, $\phi$ is an homotopy equivalence.
 
 The second claim follows immediately from the Cubical Approximation Theorem: 
 To prove surjectivity, we apply it to a given continuous map $g \colon |X_\bullet| \rightarrow |Y_\bullet|$.
 To prove injectivity, we apply the Cubical Approximation Theorem to a homotopy $H$ relating $|f_\bullet|$ and $|g_\bullet|$ so that we deform $H$ relative to $|X_\bullet| \times \partial I = |(X \otimes \partial \StandCube[ ]{1})_\bullet|$ to a cubical homotopy $|H_\bullet|$.
 
 The third claim follows from the first claim by adjunction
 \begin{equation*}
     \xymatrix{[X,Y] \ar[rr]^-{S_\bullet} \ar@/_1.5pc/[rrrr]_{\placeholder \circ \phi} && [S_\bullet X, S_\bullet Y] \ar[rr]^-{\mathrm{Ad}}_-{\cong} && [|S_\bullet X|, Y]. }
 \end{equation*}
 As $\phi$ is a homotopy equivalence, it induces a bijection between the homotopy classes and so does $S_\bullet$.
\end{proof}

As an application, we show that the geometric realisation of a combinatorially contractible cubical set is contractible.

\begin{proposition}\label{CombContract - Prop}
 Let $X_\bullet$ be combinatorially contractible. Then $X_\bullet$ is Kan and its geometric realisation is contractible.
\end{proposition}
\begin{proof}
 For a given a cubical $n$-horn
 \begin{equation*}
  \xymatrix@C+1em{\CubeBox  \ar[d] \ar[r]^{f} & X_\bullet \ar[d] \\ \StandCube{n} \ar[r] & \:\ast, }
 \end{equation*}
 the cubical map $f$ can be extended to $\partial \StandCube{n}$ as follows:
 The cubical subset $\CubeBox \cap \CubeIncl[\eps]{i}\left(\StandCube{n-1}\right)$ can be identified with $\partial \StandCube{n-1}$.
 Since $X_\bullet$ is combinatorially contractible, the restriction of $f$ to $\CubeBox \cap \CubeIncl[\eps]{i}\left(\StandCube{n-1}\right)$ can be extended to a map $f_{(i,\eps)} \colon \CubeIncl[\eps]{i}\left(\StandCube{n-1}\right) \rightarrow X_\bullet$.
 Using again that $X_\bullet$ is combinatorially contractible, we get the extension
 \begin{equation*}
     \xymatrix@C+1em{\CubeBox \ar@{^{(}->}[r] \ar[d] & \partial \StandCube{n} \ar[r]^-{f \cup f_{(i,\eps)}} \ar[d] & X_\bullet \ar[d] \\
     \StandCube{n} \ar@{=}[r] & \StandCube{n} \ar[r] \ar@{.>}[ru] & \:\ast,}
 \end{equation*}
 which shows that $X_\bullet$ is Kan.

 To show that $X_\bullet$ is contractible, we show that all homotopy groups are trivial.
 This is clear for $\pi_0$ because the definition of combinatrial contractibility for $n=1$ precisely states that every two points can be joined by a path.
 For the higher homotopy groups, let $x \in \Omega^n_{x_0}(X_\bullet)$ be given.
 It can be represented by a map $f_x \colon \StandCube{n} \rightarrow X_\bullet$ that restricts on $\partial \StandCube{n}$ to the constant map with value $k_0$.
 If we identify $\StandCube{n}$ with $\CubeIncl[-1]{1}(\StandCube{n}) \subseteq \partial \StandCube{n+1}$, then we can extend $f_x$ by the constant map with value $x_0$ to the cubical set $\partial \StandCube{n+1}$.
 As $X_\bullet$ is combinatorially contractible, the extension has a filler $h \colon \StandCube{n+1} \rightarrow X_\bullet$.
 Under the identification $\StandCube{n+1} \cong \left(\StandCube[ ]{1}\otimes \StandCube[ ]{n}\right)_\bullet$, the filler corresponds to a homotopy between the constant map $x_0$ and $f$.
 Thus, $x$ represents the zero element and $\pi_n(X,x_0)$ is therefore trivial.
 
 Theorem \ref{Unit Weak Equivalence - Thm} and the remark below it implies that $|X_\bullet|$ is a CW complex with vanishing homotopy groups, so it is contractible by Whitehead's theorem.
\end{proof}

Theorem \ref{Equivalence of Homotopy Categories - Thm} can actually be improved. 
Cisinski proved in \cite{cisinski2006prefaisceaux}, 
see also Jardine \cite{jardine2006categorical}, that the category of cubical sets carries a model structure.
The weak equivalences are cubical maps whose geometric realisations induce isomorphisms on all (topological) homotopy groups.
The fibrations are given by Kan fibrations. 
The cofibrations are just cubical maps that are level-wise injective.
The geometric realisation and the singular set functor form a Quillen equivalence between $\cSet$, equipped with the just mentioned model structure, and the category $\mathbf{Top}$ of topological space endowed with the Quillen model structure.\footnote{Oversimplifying: The Quillen model structure favours the CW approximation of a topological space over the space itself.} 

\section{Examples of Cubical Sets}\label{Section - Examples of Cubical Sets}

In this subsection we discuss some conceptual examples that we will need later in this thesis.
We begin with the construction of a different model for the singular set of topological spaces.

\subsection{Block Versions of the Singular Set}

\begin{definition}
 Let $X$ be a topological space. A continuous map $f \colon \R^n \rightarrow X$ is a \emph{block map} if there is a $\rho > 0$ such that, for all $1 \leq i \leq n$, the map $f$ satisfies
 \begin{equation*}
  f(x_1,\dots, x_i,\dots, x_n) = f(x_1, \dots, \mathrm{sgn}(x_i)\rho , \dots , x_n)
 \end{equation*}  
 if $|x_i| \geq \rho$.
\end{definition} 

We think of block maps as singular cubes with unknown side length $2\rho$.

\begin{definition}
 For a topological space $X$, let $S_\bullet^{bl}(X)$ be the cubical set whose set of $n$-cubes is given by
 \begin{equation*}
  S_n^{bl}(X) \deff \{f \colon \R^n \rightarrow X \, : \, f \text{ block map}\}
 \end{equation*}
 and whose connecting maps are given by
 \begin{align*}
  \face{i}f = \lim_{R\to\infty} f \circ \CubeIncl[R\eps]{i} \quad \text{ and } \quad \degen{i}f = f \circ \CubeProj{i}.
 \end{align*}
\end{definition}

There is a canonical inclusion $\iota_\bullet \colon S_\bullet(X) \hookrightarrow S_\bullet^{bl}(X)$ given by elongation.

\begin{lemma}\label{BlockInclusionWeakEqui - Lemma}
 The cubical set $S_\bullet^{bl}(X)$ is a Kan set and the canonical inclusion $\iota_\bullet$ is a weak homotopy equivalence.
\end{lemma}
\begin{proof}
 It is straightforward to verify that the inclusion is a cubical map.
 
 For a given singular cubical $n$-horn $f_\bullet \colon \CubeBox \rightarrow S_\bullet(X)$ that is represented by the finite collection of block maps $\{f_{(j,\omega)}\}_{(j,\omega) \neq (i,\eps)}$, we find a uniform $\rho$ such that each $f_{(j,\omega)}$ is the elongation of its restriction to $\rho I^{n-1}$.
 The compatibility $\face[\omega]{j}f_{(k,\eta)} = \face[\eta]{k-1}f_{(j,\omega)}$ implies that we can unify these maps to a continuous map $f\colon \rho |\!\CubeBox\!| \rightarrow X$. 
 Pick a retract $r \colon \rho I^n \rightarrow \rho |\!\CubeBox\!|$. 
 The constant extension of $f \circ r$ in normal direction to each hypersurface $\{\eps x_i = \rho\}$ yields a filler for the given cubical $n$-horn. 
 Thus, $S_\bullet^{bl}(X)$ is a Kan set.
 
 We will use the combinatorial homotopy groups to show that the inclusion $\iota_\bullet$ is a weak homotopy equivalence.
 Clearly, the inclusion induces a surjective map on the set of path components for $S_0(X) = S_0^{bl}(X) = X$.
 It is injective because for each block path $y \in S_1^{bl}(X)$ joining $x_{-1}$ and $x_{1}$ there is a $\rho > 0$ such that $y$ is the elongation of $y|_{\rho I}$. Then $y(\rho \cdot) \colon I \rightarrow X$ is an element in $S_1(X)$ that joins $x_{-1}$ with $x_1$. 
 
 For the higher degree case, the proof is similar. 
 We first prove surjectivity.
 Every element in $\pi_n(X_\bullet,x_0)$ is represented by a block map $f \colon \R^n \rightarrow X$ that maps each point outside some $\rho I^n$ to $x_0$.
 Abbreviate the image of $f|_{\rho I^n}(\rho \cdot \placeholder) \colon I^n \rightarrow X$ under $\iota_n$ to $\hat{f}$. 
 A joining homotopy $h \in S_{n+1}^{bl}(X)$ can be constructed as follows. 
 Pick a smooth function $\chi \colon \R \rightarrow \R$ that is identically $1$ on $\{t \leq -1\}$ and identically $\rho$ on $\{t \geq 1\}$. 
 Then $h(t,x) \deff f(\chi(t)x)$ is a block map that satisfies $\face[-1]{1} h = f$, $\face[1]{1}(h) = \hat{f}$, and $\face{j}h = x_0$ if $j>1$. 
 
 We prove injectivity by applying the previous argument to homotopies.
 Let $h \in S_{n+1}^{bl}(X)$ be a homotopy that relates the extended singular cubes $\iota_n(f_1)$ and $\iota_n(f_{-1})$.
 There is a $\rho > 0$ such that $h$ is the elongation of $h|_{\rho I^n}$.
 Thus, $h(\rho \cdot \placeholder) \colon I^{n+1} \rightarrow X$ is a homotopy between $\iota_n(f_1)|_{\rho I^n}(\rho \cdot \placeholder)$ and $\iota_n(f_{-1})|_{\rho I^n}(\rho \cdot \placeholder)$.
 But $\iota(f)_{\rho I^n}(\rho \cdot \placeholder)$ is homotopic to $f$ relative $\partial I^n$ for all continuous maps $f\in S_n(X)$.
 An example homotopy is provided by $\iota(f)|_{\rho_t I^n}(\rho_t \cdot \placeholder)$, where $\rho_t = t\rho + (1-t)\cdot 1$.
 Thus, $f_1$ and $f_{-1}$ already represent the same element in $\pi_n(S_\bullet(X))$ and injectivity is proven.
\end{proof}

If the reference space $X$ is a smooth object, which means in our context that $X$ is a finite dimensional smooth manifold or an open subset of an affine Fréchet space\footnote{All affine Fréchet spaces we encounter are closed affine subspaces of an ambient Fréchet space, so the ready may safely only think of those.}, then we can also speak about smooth block maps.
A map $\R^n \rightarrow X$ into a Fréchet space $X$ is smooth if all partial derivatives of it exists.
We denote the cubical subset of $S_\bullet^{bl}(X)$ consisting of all smooth block maps by $\SingSmooth{X}$.

\begin{lemma}\label{SmoothInclusionSingSpace - Lemma}
 The cubical set $\SingSmooth{X}$ is Kan and the inclusion $\SingSmooth{X} \hookrightarrow S^{bl}_\bullet(X)$ is a weak homotopy equivalence.
\end{lemma}
The proof of this lemma requires to generalise Whitney's Approximation Theorem to Fréchet space valued functions.
We follow the formulation of \cite{lee2013smooth}*{Theorem 6.21 and Theorem 6.26}. 
\begin{proposition}[Whitney's Approximation Theorem]\label{Whitney's Approximation - Prop}
 Let $M$ be a smooth manifold with corners, $X$ be an open subset of an affine Fréchet space with Fréchet metric $d$, and $f \colon M \rightarrow X$ a continuous map. Then the following two assertions hold true:
 \begin{enumerate}
     \item For each continuous function $\delta \colon M \rightarrow \R_{>0}$, there is a smooth map $g\colon M \rightarrow X$ satisfying $d(f,g) < \delta$. 
     If $f$ was already smooth near a closed subset $A\subset M$, then $g$ can be chosen such that $f|_A = g|_A$.
     \item The function $f$ is homotopic to a smooth function $g$. If $f$ is already smooth near $A$, then $g$ can be chosen to be equal to $f$ on $A$ and the homotopy can be chosen to be stationary on $A$.
 \end{enumerate}
\end{proposition}
\begin{proof} 
  The proof of Theorem 6.21 in \cite{lee2013smooth} carries over without too much effort. 
  However, we need to be more careful because we have to deal with an unbounded number of different semi-norms, which is why we present the argument in full detail.
  
  Recall that countably many semi-norms of a Fréchet space give rise to a metric $d$ via the following construction:
  \begin{equation*}
      d(x,y) \deff \sum_{k=1} \frac{1}{2^k} \cdot \frac{||x-y||_k}{1 + ||x-y||_k}. 
  \end{equation*}
  It follows that for all $n \in \N$ there is a $N_n \in \N$ such that if
  \begin{equation*}
      ||x-y||_k < \frac{1}{2(n+1)} \text{ for all } k \leq N_n  \Longrightarrow d(x,y) < \frac{1}{2n}.
  \end{equation*}
  Note that $d(x,y) < 1$ and that we may assume without restriction that $\delta < 1$, because we may replace $\delta$ with $\delta \cdot (1+\delta)^{-1}$ otherwise.
  For all $n \in \N$, we set $M_n \deff \delta^{-1}([1/(n+1),1/n))$ so that $M$ is a disjoint union of all $M_n$'s.
  
  Let $V$ be an open neighbourhood of $A$ on which $f$ is already smooth and set $V_n \deff \{x \in V \, : \, 1/n > \delta(x) > 1/(2n)\}$.
  For a given $x \in M\setminus A$, we find a small open neighbourhood $U_x$ in $M\setminus A$ such that all $y \in U_x$ satisfy
  \begin{equation*}
      \frac{1}{n} \delta(y) > \frac{1}{2n} \quad \text{ and } \quad ||f(y) - f(x)||_k < \frac{1}{2(n+1)} \text{ for all }k \leq N_n.
  \end{equation*}
  Pick a locally finite, countable subcover $\{U_i\}_{i \in \N}$ of $\{U_x \, : \, x \in M\setminus A\}$ and extend it with $\{V_n\}_{n\in \N}$ to a locally finite, countable cover of $M$.
  Rename $V_i$ into $U_{-(i-1)}$ so that the open cover is $\{U_i\}_{i \in \Z}$. 
  
  Pick a partition of unity $(\varphi_i)_{i \in \Z}$ that is subordinated to that cover and define
  \begin{equation*}
      g(y) \deff \sum_{i>0} \varphi_i(y) f(x_i) + \sum_{i \leq 0} \varphi_i(y) f(y).
  \end{equation*}
  Of course, we have $f|_A = g|_A$ by construction.
  
  If $y \in M\setminus A$, then we observe that
  \begin{equation*}
      d\left(g(y),f(y)\right) = d\left(\sum_{i>0}\varphi_i(y)f(x_i), \sum_{i>0} \varphi_i(y)f(y)\right).
  \end{equation*}
  By assumption, $\varphi_i(y) \neq 0$ for only finitely many $i \in \N$.
  Let $\Bar{i}$ be the minimal number with that property.
  This implies that $1/ \Bar{i} > \delta(y) > 1/(2\Bar{i})$.
  Thus, it suffices to show that 
  \begin{equation*}
      \left|\left| \sum_{i>0} \varphi_i(y)f(x_i) - \sum_{i>0} \varphi_i(y)f(y)\right|\right|_k < \frac{1}{2(\Bar{i}+1)} 
  \end{equation*}
  for all $k \leq N_{\Bar{i}}$.
  But this follows straightforwardly from 
  \begin{align*}
      \left|\left| \sum_{i>0} \varphi_i(y) f(x_i) - \sum_{i>0} \varphi_i(y) f(y) \right|\right|_k &= \left|\left|  \sum_{i>0} \varphi_i(y)(f(x_i) - f(y))\right|\right|_k \\
      &= \sum_{i>0} \varphi_i(y) ||f(x_i) - f(y)||_k \\
      &< \sum_{i>0}\varphi_i(y) \frac{1}{2(i+1)} \\
      &< \sum_{i>0}\varphi_i(y) \frac{1}{2(\Bar{i}+1)} =\frac{1}{2(\Bar{i}+1)}
  \end{align*}
  This now implies 
  \begin{equation*}
      d(g(y),f(y)) = d\left(\sum_{i>0}\varphi_i(y)f(x_i), \sum_{i>0} \varphi_i(y)f(y)\right) < \frac{1}{2\bar{i}} < \delta(y).
  \end{equation*}
  
  The second statement follows immediately from the first. 
  Since $X$ is an open subset of an affine Fréchet space there is a positive function $\delta \colon M \rightarrow \R$ such that all functions $\tilde{f}$ with values in this affine Fréchet space with $d(\tilde{f},f) < \delta$ still take values in $X$.
  If $g$ is a function as in the first statement, then a homotopy between $f$ and $g$ is given by convex combination:
  \begin{equation*}
      H \colon M \times [0,1] \rightarrow X, \qquad (m,t) \mapsto (1-t)f(m) + tg(m). 
  \end{equation*}
  (The proof of the first statement can be easily adapted to show that 
  \begin{equation*}
    d(H(m,s),f(m)) < \delta(m),  
  \end{equation*}
  so that $H$ indeed takes values in $X$).
  If $f$ and $g$ agree on $A$, then $H$ is stationary on $A$.
\end{proof}

\begin{proof}[Proof of Lemma \ref{SmoothInclusionSingSpace - Lemma}]
 It is straightforward to verify that the inclusion is a cubical map.
 
 For a given a singular cubical $n$-horn $f_\bullet \colon \CubeBox \rightarrow S_\bullet(X)$ that is represented by the finite collection of block maps $\{f_{(j,\omega)}\}_{(j,\omega) \neq (i,\eps)}$, we find a $\rho$ such that each $f_{(j,\omega)}$ is the elongation of its restriction to $\rho I^{n-1}$.
 Let $\quader[n][i,\eps][\rho] \deff \{x \in \R^n \, : \, \eps x_i \geq -\rho, x_j \in \rho I^n \text{ if } j \neq i\}$ be the infinite cuboid that is open towards $\eps \infty$ in the $i$-th coordinate.
 The compatibility $\face[\omega]{j}f_{(k,\eta)} = \face[\eta]{k-1}f_{(j,\omega)}$ implies that we can unify these maps to a continuous map $f\colon \left(\quader[n][i,\eps][\rho]\right)^c \rightarrow X$.
 Pick an open subset $U \subset \left(\quader[n][i,\eps][\rho]\right)^c$ that contains $\left(\quader[][][R]\right)^c$ for some $R > \rho$ and that is diffeomorphic to $\R^n$ relative $\left(\quader[][][R]\right)^c$.
 If we denote the diffeomorphism by $\phi$, then $f \circ \phi^{-1} \colon \R^n \rightarrow X$ yields a filler for the given cubical $n$-horn.
 
 We will use the combinatorial homotopy groups to show that the inclusion is a weak equivalence.
 
  For surjectivity, we use Proposition \ref{Whitney's Approximation - Prop} to find a smooth block map that is sufficiently $\mathcal{C}^0$-close to any chosen representative of a given element of $\pi_n(X)$.
  A homotopy between the given representative and its approximation can be given by convex combination.
  
  For injectivity, we apply the previous argument to homotopies.
\end{proof}

The following criterion is useful for constructing cubical subsets of singular sets.
Recall that $p_j \colon \R^n \rightarrow \R^{n-1}$ are the linear maps that project the $j$-th component away.

\begin{definition}
 Let $X$ be a smooth manifold or an open subset of an affine Fréchet space. 
 For $0 \leq k \leq \infty$, let $f_n \colon \mathcal{C}^k(\R^n,X) \rightarrow \mathcal{C}^0(\R^n,\R)$ be a sequence of maps.
 This sequence is called \emph{stable} if $f_n(\psi \circ \CubeProj{j}) = f_{n-1}(\psi)\circ \CubeProj{j}$ for all $1 \leq j \leq n$ 
 This sequence is called \emph{local} if, for any open $U \subset \R^n$, the equation $f_n(\psi|_U) = f_{n}(\psi)|_U$ holds.
\end{definition}
We formulate the next Lemma only in the smooth case for definiteness. 
Of course, it also holds for any other degree of differentialbility.
\begin{lemma}\label{local stable criterium - Lemma}
 Let $X$ be an open subset of an affine Fréchet space.
 Let $U$ be an open subset of $\R$ and $(f_n)_{n\in \N}$ be a local and stable sequence of maps $\mathcal{C}^\infty(\R^n,X) \rightarrow \mathcal{C}^0(\R^n,\R)$.
 Then the preimage $f_n^{-1}(\mathcal{C}^0(\R^n,U)) \cap \SingSmooth[n]{X}$ forms a cubical subset of $\SingSmooth{X}$.
\end{lemma}
\begin{proof}
 Let $\varphi \colon \R^n \rightarrow X$ be a block map such that $f_n(\varphi) \colon \R^n \rightarrow U$. 
 We need to show that the images of $f_{n-1}(\face{i}\varphi)$ and $f_{n+1}(\degen{i}\varphi)$ also lie in $U$. 
 From stability we get
 \begin{equation*}
  f_{n+1}(\degen{i} \varphi)(y) = f_{n+1}(\varphi \circ \CubeProj{i})(y) = f_{n}(\varphi)(\CubeProj{i}(y)) \in U
\end{equation*}  
for all $y \in \R^{n+1}$.
Since $\varphi$ is a block map, there is a sufficiently large $R > 0$ such that $\varphi|_{\{\eps x_i > R\}} = \face{i}\varphi \circ \CubeProj{i} |_{\{\eps x_i > R\}}$. 
Because $f_n$ is local and stable, we have
\begin{align*}
 f_{n-1}(\face{i}\varphi)(y) &= f_{n-1}(\face{i}\varphi)(\CubeProj{i}\CubeIncl[2R\eps]{i} y) \overset{\text{stab}}{=} f_n((\face{i}\varphi)\circ \CubeProj{i})(\CubeIncl[2R\eps]{i}y) \\
 &= f_n(\varphi|_{\{\eps x_i > R\}})(\CubeIncl[2R\eps]{i}y) \overset{\text{loc.}}{=} f_n(\varphi)(\CubeIncl[2R\eps]{i}(y)) \in U 
\end{align*}
for all $y \in \R^{n-1}$.
\end{proof} 

\subsection{The Topological Space of Block Maps}\label{Subsection - The Topological Space of Block Maps}

When we construct the psc Hatcher spectral sequence, we need to consider the set of block maps as a topological space (with a topology distinct from the discrete one).
The main result in this subsection is that the inclusion of the space of smooth block loops into the loop space (modeled by continuous block maps) is a weak homotopy equivalence.
In this subsection, $X$ denotes either a smooth manifold or an open subset of an affine Fréchet space.

\begin{definition}
 For each $R> 0$ set
 \begin{equation*}
     {}_R\Omega^q X \deff \{f \colon RI^q \rightarrow X \, : \,  f = x_0 \text{ near } \partial RI^q\} 
 \end{equation*}
 and 
 \begin{equation*}
     {}_R\Omega^{q,\infty} X \deff \{f \colon RI^q \rightarrow X \, : \, f \text{ smooth, }f = x_0 \text{ near } \partial RI^q\}
 \end{equation*}
 endowed with the smooth Whitney topology.     
\end{definition}
For all $R < S$, we have embeddings ${}_R\Omega^{p,(\infty)} X \hookrightarrow {}_S\Omega^{p,(\infty)} X$ given by elongation.
We use these spaces to define a block version of (smooth) loop spaces.
\begin{definition}
 Let 
 \begin{equation*}
     \Omega^q X \deff \{f \colon \R^q \rightarrow X \, : \, f \text{ block map}, f = x_0 \text{ near } \infty\} 
 \end{equation*}
 and
 \begin{equation*}
     \Omega^{q,\infty} X \deff \{f \colon \R^q \rightarrow X \, : \, f \text{ smooth block map}, f = x_0 \text{ near } \infty\}
 \end{equation*}
 be equipped with the colimit topology induced by the embeddings ${}_R\Omega^{p,(\infty)} X \hookrightarrow {}_S\Omega^{p,(\infty)} X$.
\end{definition}

We wish to show that the canonical inclusion $\Omega^{q,\infty} X \hookrightarrow \Omega^q X$ is a weak homotopy equivalence.
The proof is based on the following strengthening of the exponential law.

\begin{lemma}\label{Exponential Adjunction Diffeo - Lemma}
 Let $X$ be an open subset of an affine Fréchet space or a finite dimensional manifold and let $M$, $N$ be smooth, compact manifolds, possibly with faces. 
 Then the adjoint map 
 \begin{equation*}
     \mathcal{C}^\infty(M, \mathcal{C}^\infty(N,X)) \rightarrow \mathcal{C}^\infty(M \times N; X)
 \end{equation*}
 is a diffeomorphism between the smooth spaces.
\end{lemma}
\begin{proof}
 We first consider the case, where $X$ is an open subset of an affine Fréchet space for this is the case of interest. 
 We may assume that $X$ is an Fréchet space because smoothness is a local property and preserved by translation with a fixed element.
 Every Fréchet space is convenient in the sense of Kriegl-Michor \cite{kriegl1997convenient}*{Theorem 2.14(7)}.
 It thus becomes a Frölicher space whose smooth curves are precisely those that are infinitely often differentiable \cite{kriegl1997convenient}*{Section 23.1}.
 If $M$ is a compact manifold, then \cite{kriegl1997convenient}*{Lemma 42.5} yields that the set-theoretic identity between the set of smooth maps with a priory different smooth structures
 \begin{equation*}
     \mathfrak{C}^\infty(M,X) \xrightarrow{\id} \mathcal{C}_{\text{Frö}}^\infty(M,X)
 \end{equation*}
 is an isomorphism in the category of Frölicher spaces.
 Furthermore, the set-theoretic identity
 \begin{equation*}
     \mathfrak{C}^\infty(M,X) \xrightarrow{\id} \mathcal{C}^\infty_{\text{Fréchet}}(M,X),
 \end{equation*}
 where the right-hand side carries the smooth weak topology,
 is an isomorphism of locally convex vector space by the remark above Lemma 42.5 in \cite{kriegl1997convenient}.
 The claim follows now from \cite{kriegl1997convenient}*{Theorem 23.2}, which states that the exponential map
 \begin{equation*}
     \mathcal{C}_{\text{Frö}}^\infty(Y_1 \times Y_2, Z) \cong \mathcal{C}_{\text{Frö}}^\infty(Y_1,\mathcal{C}_{\text{Frö}}^\infty(Y_2,Z))
 \end{equation*}
 is an isomorphism of Frölicher spaces for all Frölicher spaces $X,Y,$ and $Z$.
 Thus, we have 
 \begin{equation*}
     \mathcal{C}_{\text{Fréchet}}^\infty(M,\mathcal{C}_{\text{Fréchet}}^\infty(N,X)) \cong \mathcal{C}_{\text{Fréchet}}^\infty(M\times N,X).\qedhere
 \end{equation*}
\end{proof}
For completeness, if $X$ is a finite dimensional smooth manifold, then the corresponding statement of Lemma \ref{Exponential Adjunction Diffeo - Lemma} still holds, provided $\mathcal{C}^\infty(N,X)$ is given an appropriate smooth structure. 
We refer to \cite{kriegl1997convenient}*{§42} for more details, in particular \cite{kriegl1997convenient}*{Theorem 42.14}.

\begin{lemma}\label{SmoothLoopInclusionWHE - Lemma}
 If $X$ is a finite-dimensional smooth manifold or an open subset of an affine Fréchet space, then 
 the canonical inclusion $\Omega^{q,\infty} X \hookrightarrow \Omega^q X$ is a weak homotopy equivalence.
\end{lemma}

For the proof, we face the following nuisance of the category of topoligcal spaces: It is in general not true that a continuous map $f \colon K \rightarrow \mathrm{colim}_{i \in \mathcal{I}} X_i$ from a compact space into a colimit factors through an element $X_i$ of the diagram.
 However, it is true in most cases of interest by the following result of \cite{dugger2004topological}*{Lemma A3}.
 \begin{lemma}\label{Relative T1}
   Let $(\varphi_{ij}\colon X_i \rightarrow X_j)$ be a sequential diagram of inclusions that are \emph{relatively $T_1$}, $i.e.$, for all $U \subseteq X_i$ open and all $y \in X_j \setminus \varphi_{ij}(U)$ there is an open $V \subseteq X_j$ such that $\varphi_{ij}^{-1}(V) \subseteq U$ and $y \notin V$. 
   
   Then, for each continuous map $f \colon K \rightarrow \colim X_i$ from a compact domain $K$, there is an $X_i$ such that $f$ factors through it.
 \end{lemma}
 The lemma can be applied in our setup because of the following elementary observation.
 \begin{lemma}\label{EmbeddingsReloativeT1 - Lemma}
  Let $X_1 \hookrightarrow X_2 \hookrightarrow X_3 \hookrightarrow \dotsb$ be a sequence of embeddings of $T_1$-spaces.
  Then each embedding is relative $T_1$.
 \end{lemma}
 \begin{proof}
  Fix a composition $\varphi_{ij} \colon X_i \hookrightarrow X_j$. 
  It is again an embedding.
  Let $U \subseteq X_i$ be open and $y \in X_j \setminus \varphi_{ij}(U)$ be given.
  The set $\Tilde{V} \deff X_j\setminus\{y\}$ is open as $X_j$ is $T_1$.
  Since $\varphi_{ij}$ is an embedding, we can find an open subset $\Tilde{U} \subseteq X_j$ such that $\varphi_{ij}^{-1}(\Tilde{U}) = U$.
  Thus, $V \deff \Tilde{U} \cap \Tilde{V}$ is still open, satisfies $y \notin V$, and
  \begin{equation*}
      \varphi_{ij}^{-1}(V) = \varphi_{ij}^{-1}(\Tilde{U}) \cap \varphi_{ij}^{-1}(\Tilde{V}) = U \cap \varphi_{ij}^{-1}(\tilde{V}) \subseteq U. \qedhere
  \end{equation*}
 \end{proof}

\begin{proof}[Proof of Lemma \ref{SmoothLoopInclusionWHE - Lemma}]
 It suffices to check that each continuous map $\varphi \colon M \rightarrow \Omega^q X$, where $M$ is a compact, finite-dimensional manifold, can be deformed into a map with values in $\Omega^{q,\infty}X$ relative to each closed subset on which $\varphi$ already takes values in $\Omega^{q,\infty}X$.
 
 Since $\Omega^q X$ is a sequential colimit of embeddings between Hausdorff spaces, there is an $R>0$ such that $\varphi$ factors through ${}_R\Omega^q X$ by Lemma \ref{EmbeddingsReloativeT1 - Lemma} and \ref{Relative T1}.
 Under the exponential law, it corresponds to a continuous map $\mathrm{Ad}(\varphi) \colon M \times RI^n \rightarrow X$ that is constant near the (topological) boundary $\partial RI^n$.
 Approximation theory, see Proposition \ref{Whitney's Approximation - Prop}, yields a homotopy $H \colon M \times RI^n \times [0,1] \rightarrow X$ between $\mathrm{Ad}(\varphi)$ and a smooth map $\psi \colon M \times RI^n \rightarrow X$ that is still constant near the (topological) boundary.
 If $\psi$ is smooth on $A$, then $H$ can be chosen to be stationary there.
 Now, $\mathrm{Ad}^{-1}(H)$ is a homotopy between $\varphi$ and $\mathrm{Ad}^{-1}(\psi)$ by Lemma \ref{Exponential Adjunction Diffeo - Lemma}.
 It is relative to the given closed subspace on which $\varphi$ already takes values in $\Omega^{q,\infty}X$ by the choice of $H$.
\end{proof}

\subsection{Morphism Set and Path Sets}

For two cubical sets $X_\bullet$ and $Y_\bullet$, we will construct a cubical set $\mathbf{Hom}_\bullet(X,Y)$ that has $\mathrm{Hom}_\cSet(X_\bullet,Y_\bullet)$ as $0$-cubes.
Define $\mathbf{Hom}_\bullet(X,Y)$ to the cubical set with
\begin{equation*}
 \mathbf{Hom}_n(X,Y) \deff \mathrm{Hom}_\cSet((\StandCube[ ]{n} \otimes X)_\bullet , Y_\bullet)
\end{equation*}
and whose connecting maps are given by 
\begin{equation*}
 \face{i} \deff \placeholder \circ ({\CubeIncl{i}}_\ast \otimes \id) \qquad \text{and} \qquad \degen{i} \deff \placeholder \circ ({\CubeProj{i}}_\ast \otimes \id).
\end{equation*}
It is the right adjoint functor to the reduced product
\begin{equation*}
 \mathrm{Hom}_\cSet((X \otimes Z)_\bullet, Y_\bullet) \overset{1:1}{\rightleftarrows} \mathrm{Hom}_\cSet(X_\bullet, \mathbf{Hom}_\bullet(Z,Y)), 
\end{equation*}
see \cite{jardine2006categorical}*{p.96}, which is different compared to simplicial sets, where the $\mathbf{Hom}_\bullet$ is the right adjoint to the (categorical) product.

Our primary interest lies in the \emph{cubical path set}
\begin{equation*}
 P_\bullet X \deff \mathbf{Hom}_\bullet(\StandCube[ ]{1}, X).
\end{equation*}
It comes with two cubical maps
\begin{equation*}
 \mathrm{ev}_\eps \deff \placeholder \circ {\CubeIncl{1}}_\ast \colon P_\bullet X = \mathbf{Hom}_\bullet( 
 \StandCube[ ]{1},X) \rightarrow \mathbf{Hom}_\bullet(\StandCube[ ]{0},X) \cong X_\bullet.
\end{equation*}
The latter identification comes from the Yoneda Lemma.
These two cubical maps have a section
\begin{equation*}
 s \deff \placeholder \circ {\CubeProj{1}}_\ast \colon X_\bullet \cong \mathbf{Hom}_\bullet(\StandCube[ ]{0},X) \rightarrow P_\bullet X,
\end{equation*}
which can be thought of as mapping an element to the constant path with that element as its single value.

\begin{lemma}\label{EvaluationFib-Lemma}
 If $X_\bullet$ is a Kan set, then 
 \begin{equation*}
  (\mathrm{ev}_{1},\mathrm{ev}_{-1}) \colon P_\bullet X \rightarrow (X \times X)_\bullet
 \end{equation*}
 is a Kan fibration.
\end{lemma}
\begin{proof}
 Let a commutative diagram 
 \begin{equation*}
  \xymatrix{\CubeBox \ar[r]^y \ar@{^{(}->}[d] & P_\bullet X \ar[d] \\ \StandCube{n} \ar[r]_-{(f,g)} & (X \times X)_\bullet }
 \end{equation*}
 be given.
 Via the adjunction between $\otimes$ and $\mathbf{Hom}_\bullet$, the map $y$ corresponds to a cubical map $\mathrm{Ad}(y) \colon (\CubeBox \otimes \StandCube[ ]{1})_\bullet \rightarrow X_\bullet$.
 We use the decomposition 
 \begin{equation*}
     \CubeBox[i][\eps][n+1] = \CubeBox \otimes \StandCube{1} \cup \StandCube{n} \otimes \CubeIncl[1]{1}(\StandCube{0}) \cup \StandCube{n} \otimes \CubeIncl[-1]{1}(\StandCube{0})
 \end{equation*}
  to extend $\mathrm{Ad}(y)$ to $\CubeBox[i][\eps][n+1]$ as follows:
 \begin{equation*}
  \mathrm{Ad}(y) \cup f \cup g \colon \CubeBox[i][\eps][n+1] \rightarrow X_\bullet.
 \end{equation*}
 As $X_\bullet$ is Kan, this map extends to a filler $F \colon \StandCube{n+1} \rightarrow X_\bullet$. 
 It is not hard to verify that under the isomorphism $\StandCube{n+1} \cong (\StandCube[ ]{n} \otimes \StandCube[ ]{1})_\bullet$ and the adjunction with $\mathbf{Hom}_\bullet$ the filler $F$ yields a filler $\mathrm{Ad}^{-1}(F) \colon \StandCube{n} \rightarrow P_\bullet X$ for the initial diagram. 
\end{proof}

As for topological spaces, we will use the cubical path set to turn every cubical map between Kan sets into a Kan fibration by replacing its domain with a homotopy equivalent cubical set.

For each cubical map $f_\bullet \colon X_\bullet \rightarrow K_\bullet$ between Kan sets, we define the \emph{mapping path space} $P_\bullet(f)$ as the following pullback
\begin{equation*}
\xymatrix{ P_\bullet(f) \ar[rr] \ar[d] && P_\bullet K \ar[d]^{(\mathrm{ev}_{1},\mathrm{ev}_{-1})} \\ X_\bullet \times K_\bullet \ar[rr]_{f_\bullet \times \id} && K_\bullet \times K_\bullet. }
\end{equation*}

The universal property of pullbacks implies that $P_\bullet(f) \rightarrow  X_\bullet \times K_\bullet$ is a Kan fibration.
Composing it with the projection to $X_\bullet$ and $K_\bullet$ yields the maps $p_\bullet$ and $q_\bullet$, respectively. 
It is easy to see that the projection from a product of Kan sets to one of its factors is a Kan fibration and that fibrations are closed under compositions. 
Thus, the maps $p_\bullet$ and $q_\bullet$ are also Kan fibrations. 
The map $p_\bullet$ has a section $s_\bullet \colon X_\bullet \rightarrow P_\bullet(f)$, which is induced by $(\id,f_\bullet, s \circ f_\bullet) \colon X_\bullet \rightarrow X_\bullet \times K_\bullet \times P_\bullet K$.

In view of the Yoneda Lemma \ref{Building Block Lemma}, we have the identification of cubical sets
\begin{equation*}
    P_\bullet X \cong X_{\bullet + 1},
\end{equation*}
where the connecting maps of the latter cubical sets agree with the one of $X_\bullet$ with index shifted by $1$, for example, ${}_{X_{\bullet + 1}}\face{i} = {}_{X_{\bullet}}\face{i+1}$.
Under this identification and $\mathbf{Hom}_\bullet(\StandCube[ ]{0},X) \cong X_\bullet$, the fibration $\mathrm{ev}_\varepsilon$ corresponds to $\face{1}$.
Thus, the mapping path space can be identified with
\begin{equation*}
 P_n(f) \cong \{(x,k) \in X_n \times K_{n+1} \, : \, f_n(x) = \face[1]{1}(k)\}.
\end{equation*}
\begin{lemma} 
 The map $s_\bullet \circ p_\bullet \colon P_\bullet(f) \rightarrow P_\bullet(f)$ is homotopic to the identity.
\end{lemma}
\begin{proof}
 Let $F_\bullet = p_\bullet^{-1}(\{x_0\})$ be the fibre of $p_\bullet$. 
 It is explicitly given by 
 \begin{equation*}
  F_n \deff \{(x_0,k) \in X_n \times K_{n+1} \, : \, \face[1]{1}k = f_n(x_0)\},
 \end{equation*}
 so it agrees with the pointed path space $P_{f(x_0)}K_\bullet$.
 The latter is combinatorially contractible by an argument similar to the proof Lemma \ref{EvaluationFib-Lemma} that $(\mathrm{ev}_{1},\mathrm{ev}_{-1})$ is a Kan fibration.
 The long exact sequence for Kan fibrations between Kan sets implies that $p_\bullet$ must be a weak homotopy equivalence.
 It follows from Whitehead's Theorem that $|p_\bullet| \colon |P_\bullet(f)| \rightarrow |X_\bullet|$ is a homotopy equivalence. 
 The cubical approximation theorem implies that $p_\bullet$ must be a homotopy equivalence, too.
 This, together with the equation $p_\bullet \circ s_\bullet = \id$, implies that $s_\bullet$ is the homotopy inverse of $p_\bullet$. 
\end{proof}

\begin{definition}
 Let $f_\bullet \colon X_\bullet \rightarrow K_\bullet$ be a cubical map between Kan sets. 
 The \emph{homotopy fibre} of $f_\bullet$ at the point $k_0 \in K_0$ is the Kan set
 \begin{equation*}
    \hofib(f_\bullet;k_0) \deff q_\bullet^{-1}(\{k_0\}) = \{(x,k) \in X_n \times K_{n+1} \, : \, \face[1]{1}k= f_n(x), \, \face[-1]{1} k = k_0\}.
 \end{equation*}
 \nomenclature[Cub]{$\hofib(f_\bullet)$}{Homotopy fibre of the cubical map $f_\bullet$}
\end{definition}

\begin{rem}
   In contrast to simplicial sets, it is unknown to the author, whether the geometric realisation of a Kan fibration $ F_\bullet \xrightarrow{\iota} E_\bullet \xrightarrow{p} B_\bullet$ is a Serre fibration.
   However, Lemma \ref{CubicalLESHmptyGroups - Lemma} and Theorem \ref{Equivalence of Homotopy Categories - Thm} imply that it is a homotopy fibration, which means that the composition $|p| \circ |\iota|$ is homotopic to a constant map $\const_{b_0}$ and, for some null-homotopy $h$, the induced fibre comparison map
   \begin{equation*}
       |F_\bullet| \rightarrow \mathrm{hofib}(|p_\bullet|) \quad \text{given by} \quad f \mapsto \bigl(h(b_0,\placeholder), f, b_0\bigr)
   \end{equation*}
   is a weak homotopy equivalence.
   In particular, the canonical map
   \begin{equation*}
       |\hofib(f_\bullet;k_0)| \rightarrow \mathrm{hofib}(|f_\bullet|;|k_0|)
   \end{equation*}
   is a weak homotopy equivalence.
   Informally speaking, the cubical homotopy fibre models the correct thing!
\end{rem}


 \chapter{Cubical Versions of Positive Scalar Curvature}\label{Chapter - Cubical Versions of psc}

We will apply the theory of cubical sets developed in Chapter \ref{Chapter - Cubical Sets} to positive scalar curvature.
Cubical set theory allows us to construct the central object of this thesis, the concordance set $\ConcSet$ of positive scalar curvature metrics, which cannot be written in a closed form.
This cubical set comes with a comparison map $\susp_\bullet\colon \SingMet \dashrightarrow \ConcSet$ defined on a weakly equivalent cubical subset of $\SingMet$, the cubical set analog of the space of psc metrics $\Riem^+(M)$.
Although we already describe the mapping rule of $\susp_\bullet$ in this section, the construction of its domain, the weakly equivalent subset, and the proof that $\susp_\bullet$ takes values in $\ConcSet$ will be carried out in Chapter \ref{Factorisation Indexdiff - Chapter} after we have developed the index theoretical counterpart (because the index theoretical counter also influences the choices made in the construction of the weakly equivalent subset). 

The main part of this chapter is devoted to show that the concordance set satisfy the Kan extension property.
This is an important property, as it allows us to apply the tools from cubical set theory to the concordance set.
For example, the Kan property gives us a combinatorial description of the homotopy groups of $\ConcSet$, which is an inevitable tool for the construction of the psc Hatcher spectral sequence in Chapter \ref{The PSC Hatcher Spectral Sequence - Chapter}.

In Section \ref{Section - Foundations on the Concordance Set} we introduce the concept of \emph{block metrics} on $M\times \R^n$ and present some of their elementary properties. 
We describe the cubical analog $\SingMet$ of $\Riem^+(M)$ and construct the concordance set $\ConcSet$ as the cubical set of all block metrics with positive scalar curvature.
Section \ref{Section - Construction of Special Submanifolds} serves as a preparatory section and is quite technical. 
There, we will construct, for each block metric, a certain folation of hyperplanes in $\R^n$ that bend well adapted to the chosen block metric.
This foliation will then be used in Section \ref{Section - Kan Property of the Concordance Set} to modify a block metric with positive scalar curvature so that the resulting block metric still has positive scalar curvature.
This modification is the key ingredient to prove the Kan property, what we will do at the end of Section \ref{Section - Kan Property of the Concordance Set}.
We finish this chapter by giving a different and geometric motivated description of the groups structures on $\pi_n(\ConcSet)$ in Section \ref{Geometric Addition - Section}.

In the entire section, let $M$ be a smooth closed, manifold of dimension $d$.

\section{Foundations}\label{Section - Foundations on the Concordance Set}

We first translate the space of positive scalar curvature metrics into the combinatorial world by considering the cubical set of block maps into the space of Riemannian metrics $\Riem(M)$
\nomenclature[Conc]{$\Riem(M)$}{Space of Riemannian metrics}
and the space of psc metrics $\Riem^+(M)$ 
\nomenclature[Conc]{$\Riem^+(M)$}{Space of psc metrics}
instead of the spaces themselves.

\begin{definition}
 Let $\SingMetNon$ 
 \nomenclature[Conc]{$\SingMetNon$}{Cubical Set of smooth block maps into $\Riem(M)$}
 denote the cubical set whose $n$-cubes consists of smooth block maps $g \colon \R^n \rightarrow \Riem(M)$.
 The connecting maps are given by 
 \begin{equation*}
  \face{i}g = \lim_{R \to \infty} g \circ \CubeIncl[R\eps]{i}  \qquad \text{ and } \qquad \degen{i} g = g \circ \CubeProj{i}.
 \end{equation*}
 Let $\SingMet$ 
 \nomenclature[Conc]{$\SingMet$}{Cubical Set of smooth block maps into $\Riem^+(M)$}
 denote the cubical subset consisting of all smooth block maps with values in $\Riem^+(M)$.
\end{definition}
 These cubical sets are indeed the correct models because by Lemma \ref{BlockInclusionWeakEqui - Lemma} and \ref{SmoothInclusionSingSpace - Lemma} and Theorem \ref{Equivalence of Homotopy Categories - Thm}(1) give the following zig-zag of weak equivalences
 \begin{equation*}
     \xymatrix{ |\mathcal{R}_\bullet^{(+)}(M)| \ar[r]^-\simeq & \bigl|S_\bullet^{bl}\bigl(\Riem^{(+)}(M)\bigr)\bigr| & \bigl|S_\bullet\bigl(\Riem^{(+)}(M)\bigr)\bigr| \ar[l]_-\simeq \ar[r]^-\simeq & \Riem^{(+)}(M), }
 \end{equation*}
 where the first map is the inclusion of smooth block maps into continuous block maps, the second map is the inclusion of singular cubes into continuous block maps described before Lemma \ref{BlockInclusionWeakEqui - Lemma}, and the last map is the co-unit of the adjunction between $|\placeholder|$ and $S_\bullet(\placeholder)$  described in Theorem \ref{Unit Weak Equivalence - Thm}.
 
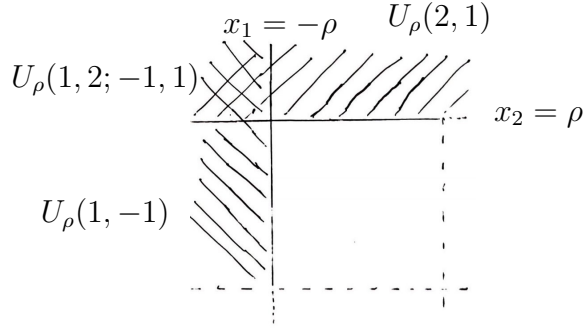
\begin{figure}
    \centering
    \begin{tikzpicture}
		\node at (0,0) {\includegraphics[width=0.25\textwidth]{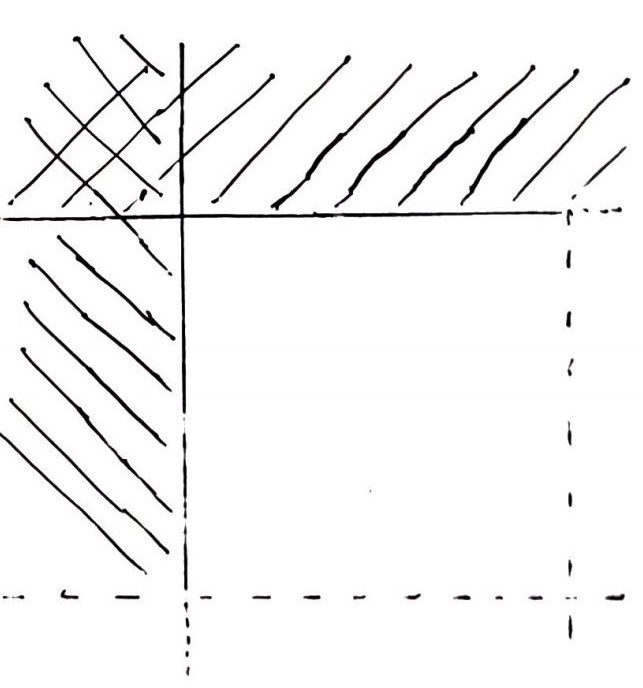}};
		\node at (1.4,2.1) {$U_\rho(2,1)$};
		\node at (-3,-.5) {$U_\rho(1,-1)$};
		\node at (-3,1.3) {$U_{\rho}(1,2;-1,1)$};
		\node at (2.7,.8) {$x_2 = \rho$};
		\node at (-.7,2) {$x_1 = -\rho$};
	\end{tikzpicture}
    \caption{An example for case $\fateps = 1 \mapsto -1$, denoted by $(1,-1)$, for the case $\fateps = 2 \mapsto 1$, denoted by $(2,1)$, and the union of these two maps, denoted by $(1,2;1,-1)$,}
    \label{fig:UNeighbourhood}
\end{figure}

 Now, we are going to construct the combinatorial reference space.
 First we have to make some organisational definitions.
 For a visualisation, we refer to figure \ref{fig:UNeighbourhood}.
 \begin{definition}\label{Block Domains - Def}
 Set $\mathbf{n} \deff \{1,\dots,n\}$.
 We define for $\fateps \colon \dom \, \fateps \subseteq \mathbf{n} \rightarrow \Z_2$, 
 \nomenclature[Conc]{$\fateps$}{A map $\dom \fateps \subseteq \{1,\dots,n\} \rightarrow \Z_2 = \{\pm 1\}$}
 the sets
  \begin{align*}
   U_\rho(\fateps) &\deff M \times \{x \in \R^n \, : \, \fateps(i)x_i > \rho \text{ for } i \in \dom \fateps\}, \nomenclature[Conc]{$U_\rho(\fateps)$}{Open subset of $M \times \R^n$ determined by $\rho$ and $\fateps$} \\
   \R^n(\fateps) &\deff \mathrm{Map}(\mathbf{n}\setminus \dom \fateps,\R) \cong \{x \in \R^n \, : \, \fateps(i)x_i = \rho\}, 
   \nomenclature[Conc]{\(\R^n(\fateps)\)}{Affine hyperplane of $\R^n$ or $\R^{n - |\dom \fateps|}$}\\
   \R^\fateps_\rho &\deff \{x \in \R^{|\dom\fateps|} \, : \, \fateps(i_j)x_j > \rho; \, j \mapsto i_j \text{ strictly monotone}\}.
   \nomenclature[Conc]{$\R^\fateps_\rho$}{Open subset of $\R^{|\dom \fateps|}$}
 \end{align*} 
 \end{definition}
 
Note that the isomorphism is canonical once $\rho$ is fixed. 
Using permutation of coordinates we can identify $U_\rho(\fateps)$ with $M \times \R^n(\fateps) \times \R^\fateps_\rho$.
To be precise, all of these identification are of the form $\id_M \times A$, where $A \in GL_n(\R)$ is a permutation matrix.

For each $\fateps$, and all $\rho$ we set $\CubeIncl[\rho\fateps]{ } \colon \R^{n - |\dom \fateps|} \rightarrow \R^n$ that includes $\rho\fateps(i)$ into the $i$-th position of the tuple for all $i \in \dom \fateps$.
The image is the affine hyperplane $\{x \in \R^n \, : \, \fateps(i)x_i=\rho\}$.

\begin{definition}\label{block metric - def}
  A Riemannian metric $g$ on $M \times \R^n$ is a \emph{block metric} if there is a $\rho > 0$ such that $g$ restricts on all $U_\rho(\fateps)$ to the following product metric  
  \begin{equation*}
    g|_{U_\rho(\fateps)} = \CubeIncl[ \rho \fateps]{ }^\ast g \oplus \euclmetric =: g\restrict_{M \times \R^n(\fateps) }\oplus \euclmetric. 
    \nomenclature[Conc]{$\restrict$}{Restriction to a codimension > 0 submanifold}
    \nomenclature[Conc]{$\euclmetric$}{Standard Euclidean Metric}
  \end{equation*}
  We say that $g$ \emph{decomposes outside of} $M \times \rho I^n$.
\end{definition}

In fact, we can weaken this definition. 

\begin{lemma}\label{Criteria block metric - lemma}
 Let $g$ be a Riemannian metric on $M \times \R^n$. 
 If there is a $\rho > 0$ such that
 \begin{equation*}
  g|_{\{\eps x_i > \rho\}} = \CubeIncl[\rho\eps]{i}^\ast g \oplus \diff x_i^2 = g\restrict_{\{x_i = \eps\}} \oplus  \diff x_i^2
\end{equation*}  
for all $(i, \eps) \in \{1,\dots,n\}\times \Z_2$, then $g$ is a block metric.
\end{lemma}

Informally, we think of a block metric as a metric on $M \times \rho I^n$ that is a product metric near each hyperface $\{x_i = \rho \eps\}$.

Since $M\times \R^n$ is an open manifolds, we have to make sure that block metrics are not covered by Gromov's h-principle. 
From this perspective, the next proposition tells us that block metrics are indeed interesting.

\begin{proposition}\label{PropertiesBlockMetric - Prop}
 Every block metric $g$ on $M \times \R^n$ has the following properties:
 \begin{enumerate}
  \item The metric $g$ is complete.
  \item If $\scal(g) > 0$, then there is a constant $c>0$ such that $\scal(g) \geq c$.
 \end{enumerate}  
\end{proposition}
\begin{proof}
 For the first statement, we will prove that $M \times \R^n$ is complete with respect to the path-length metric induced by $g$.
Consider a Cauchy sequence $\{(m_j,p_j)\}$ of points in $M \times \R^n$. Assume that $g$ decomposes outside the cube $M \times \rho I^n$.
Then, either infinitely many points lie in one of the sets $M \times Q^\eps_k$, where $Q^\eps_k \deff \{x \in \R^n \, | \, \eps \cdot x_k \geq 1.2\rho \}$, 
or infinitely many points lie in the cube $M \times 2\rho I^n$. 
In the latter case, there exists a limit due to compactness. 
In the first case, we use that $g$ is a product metric of a block metric on $M \times \CubeIncl[2\rho\eps]{k}(\R^{n-1})$ and the euclidean metric on $\R$. 
Since a product metric of two complete metrics is again complete, we conclude the existence of a limit point. 
Thus, the given block metric is complete.

 For the second statement, we make an induction over the dimension.
 For a block metric on $M\times \R^0 = M$ the statement is clear because $M$ is compact.
 Now let $g$ be a block metric on $M\times \R^n$ such that it decomposes outside of $M \times \rho I^n$ and has positive scalar curvature.
 Due to compactness, $\scal(g) \geq \Tilde{c}$ on $M \times 2\rho I^n$ for some $\Tilde{c} > 0$.
 By the induction hypothesis, the same is true for $\scal(g)$ on $\{\eps x_j > \rho\}$. 
 If $c$ the minimum of the infima of $\scal(g)$ over all these domains, then $c$ is positive and satisfies $\scal(g) \geq c$. 
\end{proof}

\begin{definition}
 Let $\BlockMetrics$ 
 \nomenclature[Conc]{$\BlockMetrics$}{Cubical Set of Block Metric on $M$}
 be the cubical set whose set of $n$-cubes is given by
 \begin{equation*}
  \BlockMetrics[n] \deff \{g \in \Riem(M \times \R^n) \, : \, g \text{ block metric }\}
 \end{equation*}
 and whose connecting maps are
 \begin{equation*}
  \face{i}g = \lim_{R\to \infty} \CubeIncl[R\eps]{i}^\ast g \qquad \text{ and } \qquad \degen{i}g = \CubeProj{j}^\ast g + \diff x_i^2.
 \end{equation*}
 The concordance set $\ConcSet$ 
 \nomenclature[Conc]{$\ConcSet$}{Cubical Set of psc block metrics on $M$}
 is the sequence of subset of $\BlockMetrics$ whose $n$-cubes have \emph{positive scalar curvature}. 
\end{definition}

\begin{lemma}
 The connecting maps of $\BlockMetrics$ are well defined and satisfy the cubical identities. 
 The connecting maps restrict to $\ConcSet$, so $\ConcSet$ is a cubical subset of $\BlockMetrics$.
\end{lemma}
\begin{proof}
 If $g$ is a block metric that decomposes outside of $M \times \rho I^n$, then so do all of its faces $\face{i}g$.
 Thus,
 \begin{align*}
     \face{i}\face[\omega]{j}g = \left(\CubeIncl[\rho\omega]{j}\CubeIncl[\rho\eps]{i}\right)^\ast g = \left(\CubeIncl[\rho\omega]{i}\CubeIncl[\rho\eps]{j-1}\right)^\ast g  = \face[\omega]{j-1}\face{i} g
 \end{align*}
 for all $i < j$.
 The identity
 \begin{align*}
     {\CubeIncl[R\varepsilon]{i}}^\ast \diff x_j = \diff x_j \circ D \CubeIncl[R \eps]{i} = \begin{cases}
       \diff x_{j-1}, & \text{if } i < j \\
       0, & \text{if } i = j, \\
       \diff x_{j}, & \text{if } i> j
     \end{cases}
 \end{align*}
implies 
 \begin{align*}
    \face{i}\degen{j} g &= \lim_{R \to \infty} {\CubeIncl[R \eps]{i}}^\ast \left(\CubeProj{j}^\ast g + \diff x_j^2\right) = \lim_{R \to \infty} {\CubeIncl[R \eps]{i}}^\ast \CubeProj{j}^\ast g + {\CubeIncl[R \eps]{i}}^\ast\diff x_j^2 \\
    &= \begin{cases}
      \lim_{R \to \infty} \CubeProj{j-1}^\ast {\CubeIncl[R\eps]{i}}^\ast g + \diff x_{j-1}^2, & \text{if } i<j, \\
      g, & \text{if } i = j,\\
      \lim_{R \to \infty} \CubeProj{j}{\CubeIncl[R \eps]{i-1}}^\ast g + \diff x_j^2 & \text{if } i > j.
    \end{cases}\\
    &= \begin{cases}
      \degen{j-1}\face{i} g, & \text{if } i < j, \\
      g, & \text{if } i =j, \\
      \degen{j}\face{i-1}g, & \text{if } i>j.
    \end{cases}
 \end{align*}
 Lastly, the equation
 \begin{align*}
      \CubeProj{i}^\ast \diff x_j = \diff x_{j} \circ \CubeProj{i} = \begin{cases}
        \diff x_{j+1}, & \text{if } i \leq j, \\
        \diff x_j, & \text{if } i > j,
      \end{cases}
 \end{align*}
 gives, for $i \leq j$, the equation
 \begin{align*}
     \degen{i}\degen{j}g &= \CubeProj{i}^\ast( \CubeProj{j}^\ast g + \diff x_j^2) + \diff x_i^2  \\
     &= \CubeProj{i}^\ast \CubeProj{j}^\ast g + \CubeProj{i}^\ast \diff x_j^2 + \diff x_i^2 \\
     &= \CubeProj{j+1}^\ast (\CubeProj{i}^\ast g + \diff x_i^2) + \diff x_{j+1}^2 \\
     &= \degen{j+1}\degen{i} g,
 \end{align*}
 so the first claim is proven.
 
 The connecting maps send positive scalar curvature metrics to positive scalar curvature metrics because scalar curvature is additive under taking products, see Lemma \ref{Scalar Curvature Formuala - Prop} below.
 In more detail,
 \begin{align*}
     \scal(\degen{i}g) &= \scal(\CubeProj{i}^\ast g + \diff x_i^2) = \scal(g \oplus \euclmetric) = \scal(g) > 0,
 \end{align*}
 and, since $g|_{\{\eps x_i > \rho\}}$ is the product metric of $\face{i}g$ and the euclidean metric, we have
 \begin{align*}
     \scal(\face{i}g) &= \scal(\CubeProj{i}^\ast \face{i}g + \diff x_i^2)|_{\{\eps x_i> \rho\}} = \scal(g)|_{\{\eps x_i > \rho\}} > 0,
 \end{align*}
 which implies the second claim.
\end{proof}
We would like to know whether $\BlockMetrics$ and $\ConcSet$ are Kan sets.
For $\BlockMetrics$ even more is true.

\begin{proposition}\label{BlockMetric Combi Contrac - Prop}
 The cubical set $\BlockMetrics$ is combinatorially contractible.
\end{proposition}

The proposition follows immediately from the following two elementary lemmas.

\begin{lemma}\label{Metric Patching - Lemma}
 Let $U, V$ be an open cover of a not necessarily compact manifold $N$.
 Assume there is a submanifold $H$ such that $U \cap V$ is diffeomorphic to $H \times (-1,1)$.
 If $g_U$ and $g_V$ are Riemannian metrics on $U$ and $V$, respectively, then there is a Riemannian metric $g$ on $N$ such that 
 \begin{equation*}
  g|_{U \setminus \overline{V}} = g_U|_{U \setminus \overline{V}} \quad \text{ and } \quad g|_{V \setminus \overline{U}} = g_V|_{V \setminus \overline{U}}. 
\end{equation*}  
\end{lemma}
\begin{proof}
 Pick a smooth function $\chi \colon (-1,1) \rightarrow [0,1]$ that is identically zero near $-1$ and identically $1$ near $1$.
 Let $\phi \colon U \cap V \rightarrow H \times (-1,1)$ be a diffeomorphism.
 Define on $U \cap V$ the Riemannian metric 
 \begin{equation*}
     g|_{U \cap V} \deff \left(\chi \circ \pr_2 \circ \phi\right)  \cdot g_U + \left(1-\chi\circ \pr_2 \circ \phi\right) \cdot g_V
 \end{equation*}
  
 It can be extended via $g_U$ and $g_V$ to the rest of $M$ and the result $g$ is the desired metric. 
\end{proof}

\begin{lemma}\label{Metric Suspension Union - Lemma}
 For $(j,\omega)$ and $(k,\eta) \in \{1,\dots,n\} \times \Z_2$, pick two block metrics $g_{(j,\omega)}$ and $g_{(k,\eta)}$ on $M \times \R^{n}$ that decompose outside of $M \times \rho I^n$. 
 Assume that $j \leq k$ and that $\face[\omega]{j}g_{(k,\eta)} = \face[\eta]{k-1}g_{(j,\omega)}$.
 Then, the metric $\degen{j}g_{(j,\omega)}$ and $\degen{k}g_{(k,\eta)}$ agree on the set $\{\omega x_j > \rho, \eta x_k > \rho\}$. 
\end{lemma}
\begin{proof}
 Abbreviate $\{\omega x_j > \rho, \eta x_k > \rho\}$ to $U$.
 If $j = k$ there is nothing to prove for either $U$ is empty or $(j,\omega) = (k,\eta)$. 
 For $j<k$, the calculation
 \begin{align*}
  \degen{j}g_{(j,\omega)}|_U &= \degen{k} \face[\eta]{k} \degen{j} g_{(j,\omega)} |_U \\
  &= \degen{j} \degen{k-1} \face[\eta]{k-1} g_{(j,\omega)} |_U\\
  &= \degen{j} \degen{k-1} \face[\omega]{j}g_{(k,\eta)} |_U \\
  &= \degen{j} \face[\omega]{j} \degen{k} g_{(k,\eta)} |_U = \degen{k} g_{(k,\eta)}|_U
 \end{align*}
 shows the claim. 
 Here, we used the block form of $g_{(j,\omega)}$ in the first line, cubical identities in the second line, and the assumed compatibility in the third line.
\end{proof}

\begin{proof}[Proof of Proposition \ref{BlockMetric Combi Contrac - Prop}] 
 Let $\partial \StandCube{n} \rightarrow \BlockMetrics$ be given by the family of block metrics
 \begin{equation*}
  \{g_{(j,\omega)} \in \BlockMetrics[n-1]\, : \, \face[\omega]{j}g_{(k,\eta)} = \face[\eta]{k-1} g_{(j,\omega)} \text{ for } j < k\},
 \end{equation*}
 where $(j,\omega), (k,\eta) \in \{1,\dots,n\} \times \Z_2$.
 Since we are given only finitely many block metrics, we can pick a sufficient large $\rho > 0$ such that all $g_{(j,\omega)}$ decompose outside of $M \times \rho I^{n-1}$.
 The compatibility assumption yields that the metrics $\degen{j}g_{(j,\omega)}|_{\{\omega x_j > \rho\}}$ agree on the intersection of their domain by Lemma \ref{Metric Suspension Union - Lemma}.
 Thus, these metrics form a Riemannian metric $g$ on $M \times \R^n \setminus \rho I^n$.
 Pick some Riemannian metric $g_0$ on $M \times (\rho + 2) I^n$ and apply Lemma \ref{Metric Patching - Lemma}.
 The result $g_{fill}$ is obviously a block metric because it agrees with $g$ outside a compact subset and therefore satisfies $\face[\omega]{j}g_{fill} = g_{(j,\omega)}$.
 Thus, it is a desired filler for the initially given cubical sphere, so $\BlockMetrics$ is combinatorially contractible. 
\end{proof}

To prove that $\ConcSet$ is a Kan set is more complicated. 
In the next section, we will do some preparations and finish off the proof in Section \ref{Section - Kan Property of the Concordance Set}.

We close this section by introducing the comparison map between $\SingMetNon$ and $\BlockMetrics$.

\begin{definition}\label{Suspension Map - Definition}
 For each smooth block map $g \in \SingMetNon[][n]$, we define $\susp_n(g)$ as the adjoint metric on $M \times \R^n$.
 More precisely, we set
 \begin{equation*}
  \susp_n(g)_{(m,x)} = g(x)_m + \diff x_1^2 + \dots + \diff x_n^2
 \end{equation*}
 \nomenclature[Conc]{$\susp_n(g)$}{Suspension of the block map $\R^n \rightarrow \Riem^{(+)}(M)$}
 and call it the \emph{suspension} of $g$.
\end{definition}

\begin{lemma}
 The suspension of every $g \in \SingMetNon[][n]$ is a block metric and all suspension maps assemble to a cubical map $\susp_\bullet \colon \SingMetNon \rightarrow \BlockMetrics$. 
\end{lemma}
\begin{proof}
 Since $g$ is a block map, there is a $\rho > 0$ such $g|_{\{\pm x_i > \rho\}}$ does not depend on $x_i$ anymore.
 Thus, on $\{\eps x_i > \rho\}$, the suspension $\susp_n(g)$ is the product of $\CubeIncl[\rho\eps]{i}\susp_n(g)$ and $\diff x_i^2$. 
 Lemma \ref{Criteria block metric - lemma}. 
 implies that $\susp_n g$ is a block metric.
 
 The second statement follows from the calculations
 \begin{align*}
     \face{i} \susp_n(g)_{(m,x)} &= \lim_{R \to \infty} {\CubeIncl[R\eps]{i}}^\ast \left(g(\placeholder)_{(\placeholder)} + \diff x_1^2 + \dots + \diff x_n^2\right)_{(m,x)} \\
     &= \lim_{R \to \infty} g(\CubeIncl[R\eps]{i}(x))_m + {\CubeIncl[R\eps]{i}}^\ast \diff x_1^2 + \dots +  {\CubeIncl[R\eps]{i}}^\ast \diff x_n^2 \\
     &= \lim_{R \to \infty} (\face{i}g)_m(x) + \diff x_1^2 + \dots + \diff x_{n-1}^2 \\
     &= \susp_{n-1}(\face{i}g)_{(m,x)}
 \end{align*}
 and
 \begin{align*}
     \degen{i}\susp_n(g)_{(m,x)} &= \CubeProj{i}^\ast \left( g(\placeholder)_{(\placeholder)} + \diff x_1^2 + \dots + \diff x_n^2\right)_{(m,x)} + \diff x_i^2 \\
     &= g(p_i(x))_m + \CubeProj{i}^\ast \diff x_1^2 + \dots + \CubeProj{i}^\ast \diff x_n^2 + \diff x_i^2 \\
     &= \susp_{n+1}(\degen{i}g)_{(m,x)}.
 \end{align*}
\end{proof}

Unfortunately, it is not true that the suspension restricts to a cubical map $\susp_\bullet \colon \SingMet \rightarrow \ConcSet$.
If the parameter of the block map $g$ 'runs too fast', the scalar curvature of $\susp_n(g)$ might be non-positive at some points.
This can be already seen in the following example: The euclidean metric on $2D^n \setminus \overset{\circ}{D^n}$ expressed in polar coordinates is the suspension of the curve $[1,2]\rightarrow \mathrm{Riem}^+(S^{n-1})$ that is given by $r \mapsto r^2\cdot g_{S^{n-1}}$,
where $g_{S^{n-1}}$ is the standard round metric.

Later we will construct a weakly equivalent Kan subset $A_\bullet$ of $\SingMet$ so that $\susp_\bullet|_{A_\bullet}$ takes values in $\ConcSet$.
But we would like this cubical subset $\AuxSingMet$ to be adapted to the index theory we will develop Chapter \ref{The Operator Concordance Set - Chapter}. 
Therefore, we postpone the construction of $\AuxSingMet$ to Chapter \ref{Factorisation Indexdiff - Chapter}.

 \section{Construction of Special Submanifolds}\label{Section - Construction of Special Submanifolds}


The main result of this section, Theorem \ref{Embedding Theorem}, is an embedding result.
It allows us to find, for each given block metric, a family of neatly embedded hyperplanes in the sense that the planes bend only in those direction in which the block metric is euclidean.
More precisely, the unit normal vector field of these hyperplanes always lies in a subspace on which the given block metric restricts to the euclidean metric.

The application of this theorem in the proof of the Kan property are two-fold.
The existence of these hyperplanes provide an open, cylindrical neighbourhood in which we can deform the block metric in a way such that the result has product structure on one end of the cylinder.
Secondly, the theorem provides a diffeomorphism that we will use to ``push'' the modified metric on (an open neighbourhood) the given initial $n$-horn onto the entire cube. 


We start with definitions of certain subsets of $\R^n$ that we will need in the proof of the embedding theorem and throughout the thesis overall.

\begin{figure}[htbp]
    \centering
    \begin{tikzpicture}
      \node at (0,0) {\includegraphics[width=0.5\textwidth]{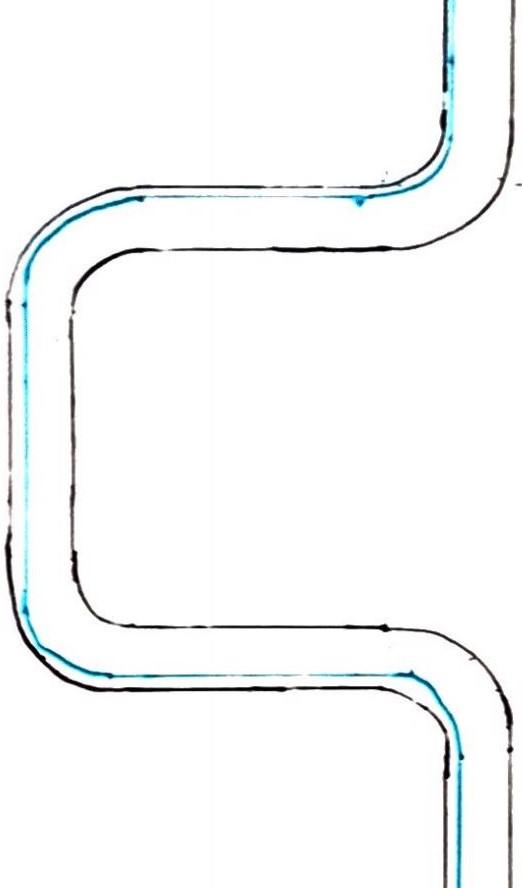}};
      \draw[dashed] (-5,1) -- (5,1);
      \draw[dashed] (-5,-1) -- (5,-1);
      \draw[dashed] (-1,-6.3) -- (-1,6.3);
      \draw[dashed] (1,-6.3) -- (1,6.3);
      
      \draw[thick] (-1,1) -- (6.3,1);
      \draw[thick] (-1,1) -- (-1,-1);
      \draw[thick] (-1,-1) -- (6.3,-1);
      
      \node at (0,5) {$V_{\rho-4}(0,1)$};
      \node at (0,-5) {$V_{\rho-4}(0,1)$};
      
      \node at (4,1.8) {$V_{\rho-4}(1,1)$};
      \node at (4,-1.8) {$V_{\rho-4}(1,-1)$};
      
      \node at (-5,0) {$V_{\rho-4}(-1,0)$};
      \node at (-4,5) {$V_{\rho-4}(-1,1)$};
      
      \node at (-4,-5) {$V_{\rho-4}(-1,-1)$};
      \node at (0,0) {$V_\rho(0,0)$};
      \node at (4,0) {$V_{\rho}(1,0)$};
      
      \node at (1.3,.6) {$\quader[n][(1,1)][\rho-4]$};
      
      \node at (0,3.2) {$\Phi(H\times [1,2])$};
      \node at (6,-3) {$H = H\times\{1.9\}$};
      \draw[->] (6,-3.2) parabola[bend at end] (2.89,-5);
      
      \node at (7,1.2) {$x_2 = \rho-4$};
      \node at (7,-1.2) {$x_2 = -(\rho-4)$};
      \node at (1,7) {$x_1 = \rho-4$};
      \node at (-1,6.5) {$x_1 = -(\rho - 4)$};
      
      \node at (-4.4,-3.3) {$U_{left}$};

    \end{tikzpicture}
    \caption{The shape of $\Phi(H\times[1,2])$ in $\R^n$, here for $n=2$. The set $H$ is canonically identified with $\Phi(H \times \{1.9\})$.}
    \label{fig:EmbeddingTheorem}
\end{figure}
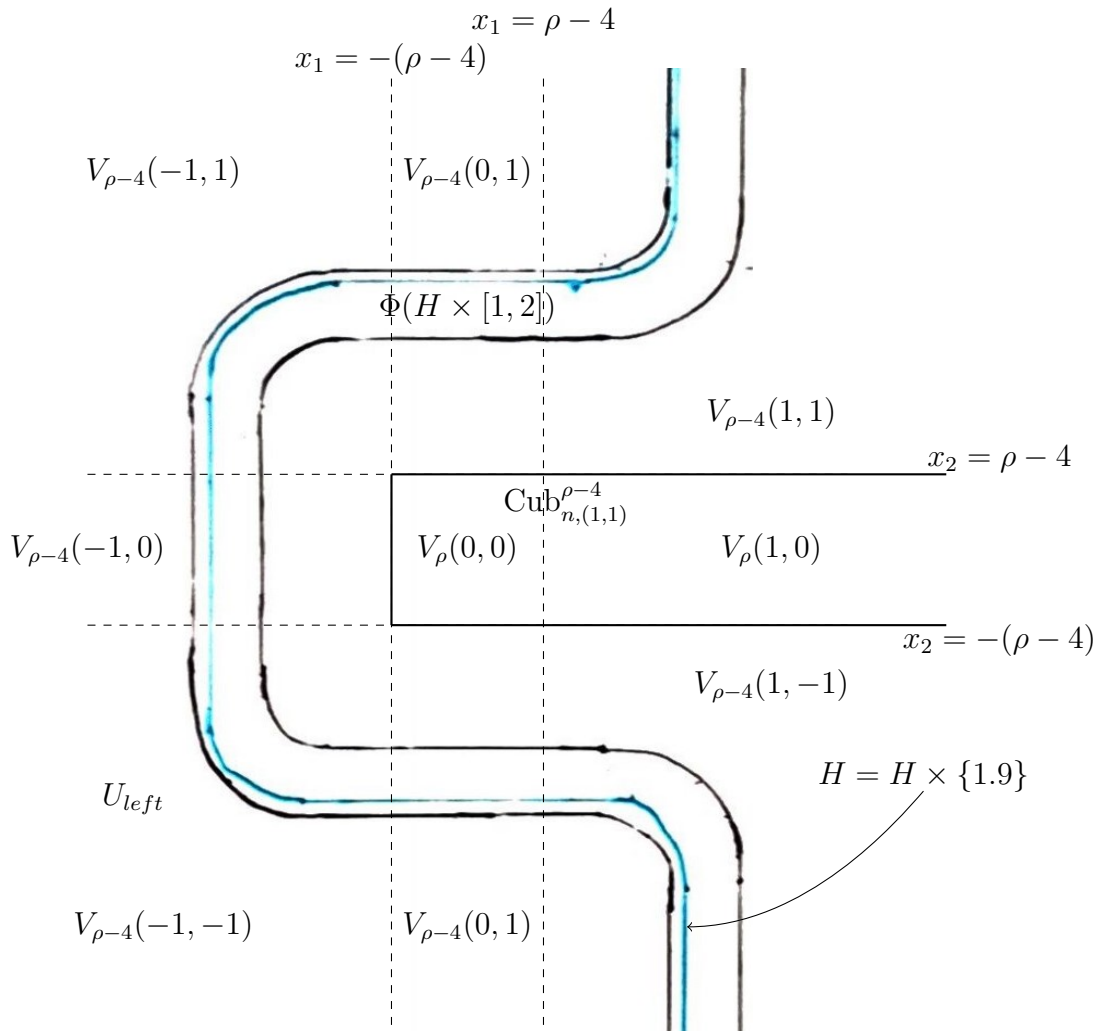

\begin{definition}
 We define, for each $\fatalpha \colon \{1,\dots,n\} \rightarrow \Z_3 = \{0,\pm1\}$
 \nomenclature[Conc]{$\fatalpha$}{Element of $\Z_3^n = \{0,\pm1\}^n$}
 and each $\rho > 0$, the open subsets
 \begin{equation*}
  V_\rho(\fatalpha) \deff \{x \in \R^n \, : \, \fatalpha(i)x_i > \rho-1 \text{ if } i \in \supp \fatalpha, x_i \in (-\rho,\rho) \text{ if } \fatalpha(i) = 0\}.
 \end{equation*}
 \nomenclature[Conc]{$V_\rho(\fatalpha)$}{Open subset of $\R^n$}
 
 For $\rho > 0$ and $(i,\eps) \in \{1,\dots,n\} \times \Z_2$, define the \emph{infinite half cuboid open towards} $(i,\eps)$ as
 \begin{equation*}
     \quader[n][i,\eps][\rho] \deff \{x \in \R^n \, : \, \eps x_i \geq -\rho,\  x_j \in [-\rho,\rho] \text{ if } j \neq i\}.
 \end{equation*}
 \nomenclature[Conc]{$\quader[n][i,\eps][\rho]$}{Infinite half cuboid open towards $(i,\eps)$}
\end{definition}
The open subsets $V_\rho(\fatalpha)$ should be thought of as an open partition of $\R^n$ with only one relative compact set $V_\rho(\mathbf{0})$, see figure \ref{fig:EmbeddingTheorem} for a visualisation in which we chose $\rho$ so large that the difference between $\rho$ and $\rho - 1$ is not visible.

For each $\fatalpha \colon \{1,\dots,n\} \rightarrow \Z_3$, let $\hat{\fatalpha} = \fatalpha|_{\supp \fatalpha}$\nomenclature[Conc]{$\hat{\fatalpha}$}{Restriction of $\fatalpha$ to its support}. 
The tupel $\hat{\fatalpha}$ parameterises the coordinates in those directions $V_\rho(\fatalpha)$ is unbounded and, more importantly for our purpose, in which coordinate-direction a block metric decomposes.
Note also that $U_{\rho-1}(\fatalpha)$ is the union of all $V_\rho(\fatalpha')$ with $\supp \fatalpha' \supseteq \supp \fatalpha$.

\begin{theorem}[Embedding Theorem] \label{Embedding Theorem}
 For all $\rho > 5$, there is a smooth hypersurface $H \subset \R^n \setminus \quader[n][1,1][\rho-4]$ and an embedding $\Phi \colon H \times [1,2] \rightarrow \R^n \setminus \quader[n][1,1][\rho-4]$ with the following properties for all $\fatalpha$:
 \begin{enumerate}
  \item \begin{itemize}
          \item[(i)] 
          The permutation that identifies $V_{\rho-4}(\fatalpha)$ with $(-\rho+4,\rho-4)^{\mathrm{Null}(\fatalpha)} \times \R^{\hat{\fatalpha}}_{\rho-5}$ identifies $H \cap V_{\rho-4}(\fatalpha)$ with $(-\rho + 4,\rho-4)^{|\mathrm{Null}(\fatalpha)|} \times H^{\hat{\fatalpha}}$ for some hypersurface $H^{\hat{\fatalpha}} \subseteq \R^{\hat{\fatalpha}}_{\rho-5}$.
          \item[(ii)] Under this identification, the embedding $\Phi$ corresponds to $\id \times \Phi^{\hat{\fatalpha}}$.
          \item[(iii)] Outside a compact set, $\Phi_t(x) = x - (t-1.9)e_1$. 
        \end{itemize}
  \item The hypersurface $H$ separates $\R^n$ into two unbounded open subsets.
        If $U_{left}$ denotes the component that contains $\Phi(H \times \{2\})$, then there is a diffeomorphism 
        \begin{equation*}
          \Psi_n \colon U_{left} \cup \Phi(H\times [1,2]) \rightarrow \{x_1 \leq \rho + 1\}
        \end{equation*}
        that is the identity on $\Phi(H \times [1.8,2])$ and outside a compact subset of $\Phi(H \times [1,2])$.
  \item There are smooth maps $h \colon [1,2] \rightarrow \Riem(H)$ and $f \colon [1,2] \rightarrow \R_{>0}$ that are constant outside a compact set $K$ with values $\euclmetric_{\R^{n-1}}$ and $1$, respectively, such that
  \begin{equation*}
   \left(\Phi^\ast \euclmetric\right)_{(p,t)} = h(t)_p + f^2(p,t) \diff t^2. 
  \end{equation*}
  In fact, on $\Phi^{-1}(H \cap V_{\rho-4}(\fatalpha))$, the curve of metrics $h$ decomposes into
  $t \mapsto \euclmetric_{\R^n(\hat{\fatalpha})} \oplus h^{\hat{\fatalpha}}(t)$ and the map $f$ factors through $\R^{\hat{\fatalpha}}_{\rho-5}$.
 \end{enumerate}
\end{theorem}
Unfortunately, the proof of this theorem is quite technical and the entire section is devoted to it.
The theorem is therefore formulated in a manner so that the reader can safely take it for granted during his or her first reading. 
We will not need the ingredients of its proof later, with the following definition as the only exception.
However, it might be helpful to have a look at the pictures presented here.

\begin{definition}
 For $\rho > 2$, let $\aux[\rho] \colon \R \rightarrow \R$ be a smooth, convex function that satisfies 
 \begin{equation*}
   \aux[\rho]^{-1}(\{0\}) = \R_{\leq \rho - 2} \qquad \text{and} \qquad \aux[\rho]|_{\R_{\geq \rho}} = \id - (\rho -1).
 \end{equation*}
 Define the \emph{dice function }$\dice_{n,\rho} \colon \R^n \rightarrow \R$ via 
 \begin{equation*}
  \dice_{n,\rho}(x) \deff \sum_{j=1}^n \aux[\rho](|x_j|).
 \end{equation*}
\end{definition}

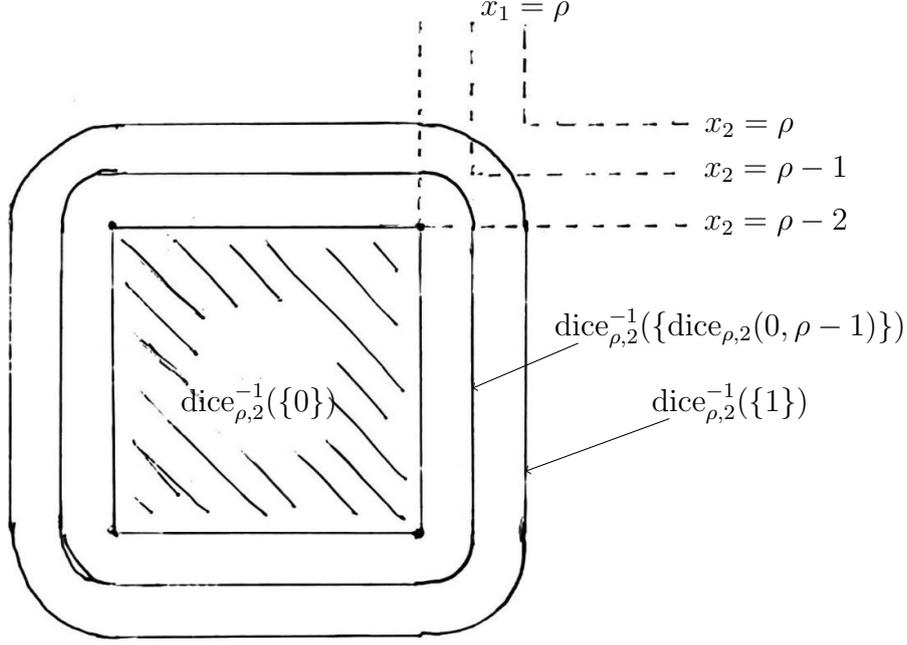
\begin{figure}[htbp]
    \centering
    \begin{tikzpicture}
      \node at (0,0) {\includegraphics[width=0.6\textwidth]{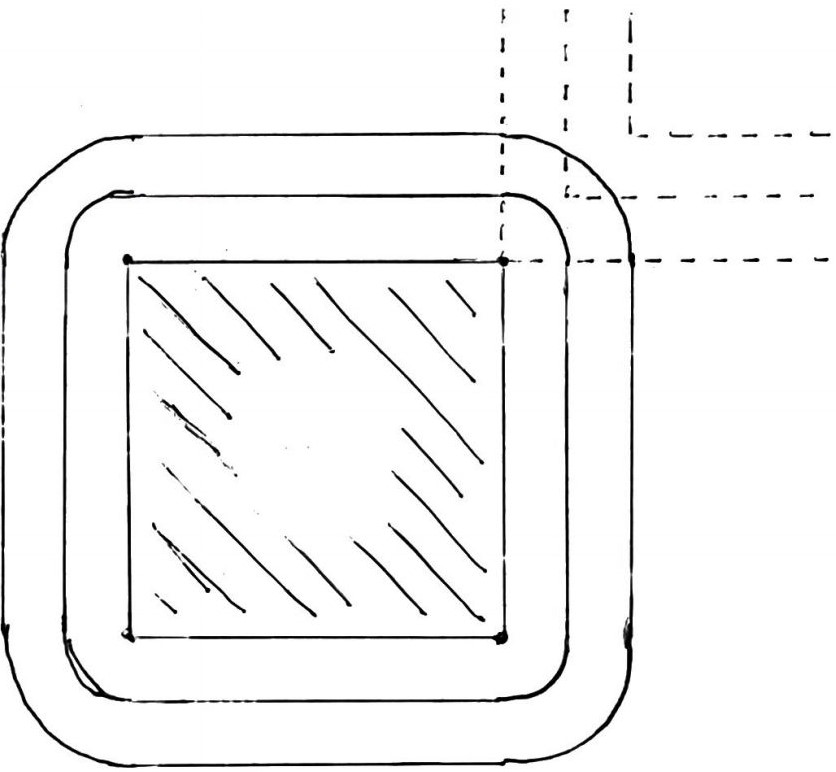}};
      \node at (-1.2,-1) {$\dice^{-1}_{\rho,2}(\{0\})$};
      \node at (5.6,1.4) {$x_2=\rho-2$};
      \node at (5.6,2.1) {$x_2=\rho-1$};
      \node at (5.25,2.65) {$x_2=\rho$};
      \node at (2.3,4.2) {$x_1=\rho$};
      \node at (5,0) {$\dice_{\rho,2}^{-1}(\{\dice_{\rho,2}(0,\rho-1)\})$};
      \node at (5,-1) {$\dice^{-1}_{\rho,2}(\{1\})$};
      \draw[->] (3,-.2) -- (1.6,-0.8);
      \draw[->] (4.2,-1.2) -- (2.3,-1.9);
    \end{tikzpicture}
    \caption{The level sets of a dice function. Recall that $\dice_{\rho,2}(0,\rho) = 1$ by construction.}
    \label{fig:LevelSetDiceFunction}
\end{figure}

The level set of $\dice_{n,\rho}$ of values around $1$ look like a cube with smoothed corners, see Figure \ref{fig:LevelSetDiceFunction}, hence the name.
In the following, we will abbreviate $\dice_{2,\rho}$ to $\dice$.
Fix $0 < \eps < (\rho-2)/2$ and set 
\begin{equation*}
 \dice^l \deff \dice|_{\{x_1 < \eps\}}.
\end{equation*} 
Furthermore, we define 
\begin{align*}
 \dice^o(x_1,x_2) &\deff -\dice(x_1,x_2 - (2\rho + 1)) + 3  \quad \text{ on } (-\eps,\rho+1]\times [\rho,\rho+4], \\
 \dice^u(x_1,x_2) &\deff -\dice(x_1,x_2 + (2\rho+1)) + 3  \quad \text{ on } (-\eps,\rho+1]\times [-\rho,-(\rho+4)] \\
 &\, = \dice^o(x_1,-x_2).
\end{align*}

\begin{lemma}
 The functions $\dice^l$ and $\dice^o$ agree on $(-\eps,\eps) \times [\rho,\rho + 1]$, while $\dice^l$ and $\dice^u$ agree on $(-\eps,\eps) \times [-(\rho+1),-\rho]$. 
\end{lemma}
\begin{proof}
 For $(x_1,x_2) \in (-\eps,\eps) \times [\rho,\rho+1]$, we have
 \begin{align*}
  \dice^o(x_1,x_2) &= -\dice(x_1,x_2-(2\rho+1)) + 3 = -\aux[\rho](|x_2 - (2\rho + 1)|) + 3 \\
  &= -\left(|x_2-(2\rho+1)| - (\rho-1)\right) + 3 \\
  &= - \left(2\rho + 1 - x_2 - (\rho-1)\right) + 3 \\
  &= x_2 - (\rho - 1) = \dice^l(x_1,x_2), 
 \end{align*}
 where we used $|x_2-(2\rho+1)| \in [\rho, \rho+1]$ to go from the first to the second line.
 The second statement follows from the first because $\dice^l(x_1,x_2) = \dice^l(x_1,-x_2)$ and $\dice^u(x_1,x_2) = \dice^o(x_1,-x_2)$.
\end{proof}

On $[\rho,\rho+1] \times [\rho+3,\rho+4]$, we have 
\begin{align*}
 \dice^o(x_1,x_2) &= -\left(\aux[\rho](|x_1|) + \aux[\rho](|x_2 - (2\rho+1)|) \right) + 3 \\
 &= - \left(x_1 - (\rho-1)\right) + 0 + 3 \\
 &= - (x_1 - \rho) + 2   
\end{align*}
and, similarly, on $[\rho,\rho+1] \times [-(\rho+4), -(\rho+3)]$ we have
\begin{align*}
 \dice^u(x_1,x_2) = - (x_1 - \rho) + 2. 
\end{align*}
This allows us to merge these functions.
For 
\begin{align*}
    Q_1 &\deff \bigl([-\eps, \rho] \times [\rho,\rho+1]\bigr)  \cup \bigl([\rho -2, \rho +1] \times [\rho,\rho+4]\bigr),\\
    Q_2 &\deff \bigl([-\eps, \rho] \times [-(\rho+1),-\rho]\bigr) \cup  \bigl([\rho -2, \rho +1] \times [-(\rho+4),\rho]\bigr),
\end{align*}
we define
\begin{align*}
    \mathtt{d}_\rho(x) \deff \begin{cases}
      \dice^o(x), & \text{if } x \in Q_1, \\
      \dice^u(x), & \text{if } x \in Q_2, \\
      \dice^l, & \text{if } x \in \{x_1 < \eps\}, \\
      -(\pr_1 - \rho) + 2, & \text{if } x \in [\rho, \rho+1] \times [-(\rho+3),\rho+3]^c.
    \end{cases}
\end{align*}

\begin{lemma}\label{Regular Value two dim function - Lemma}
 The set of regular values for $\dice_{n,\rho}$ is $\R_{>0}$.
 The interval $[1,2]$ is a set of regular values for $\mathtt{d}_\rho$.
\end{lemma}
\begin{proof}
 Since $\aux[\rho]$ is a smooth and convex function, its derivative is monotonically increasing. 
 As its null-set is $\R_{\rho-2}$, the derivative $\mathrm{aux}'_\rho$ must be positive on $\R_{>\rho-2}$. 
 The first claim now follows from
 \begin{equation*}
     \diff \, \dice_{n,\rho} = \sum_{j=1}^n \mathrm{sgn}(x_j) \mathrm{aux}'_\rho(|x_j|) \diff x_j.
 \end{equation*}
 
 It suffices to prove the second claim for each summand.
 For 
 \begin{equation*}
  -(\pr_1 - \rho) +2 \colon [\rho,\rho+1] \times [-(\rho +3),\rho+3]^c \rightarrow \R,
 \end{equation*}
 there is nothing to prove.
 The first claim implies that $\dice^o$ and $\dice^u$ have only one singular value, namely $3$.
 Since $\dice^l$ is the restriction of $\dice$ to $\{x_1 < \eps\}$, it only has one singular value, namely $0$.
 Thus, $[1,2]$ is indeed a set of regular values for $\mathtt{d}_\rho$.
\end{proof}

Similar arguments also yield the following lemma.

\begin{lemma}\label{dice vanishing at corners - Lemma}
 If $\partial_j \mathtt{d}_{\rho}(x) \neq 0$, then $|x_j| > \rho - 2$. 
\end{lemma}
\begin{proof}
 It suffices to prove the statement for each summand.
 For $-(\pr_1 - \rho) +2$ there is nothing to prove as all points in its domain satisfy the claimed inequalities.
 
 All $x \in \mathrm{dom}(\dice^o)$ satisfy $|x_2| > \rho - 2$, so we have to check the statement only for $\partial_1 \dice^o$.
 Since $\aux[\rho]'$ is monotonically increasing and $\aux[\rho]^{-1}(\{0\}) = \R_{\rho-2}$ the condition
 \begin{equation*}
  (\partial_1 \dice^o)(x) = \partial_1 \aux(|x_1|) = \mathrm{sgn}(x_1) \aux[\rho]'(|x_1|) \neq 0
 \end{equation*}
 is equivalent to $|x_1| > \rho-2$.
 
 Since $\dice^u(x_1,x_2) = \dice^o(x_1,-x_2)$ the statement holds true also for $\dice^u$.
 
 For $\dice^l$, we deduce from
 \begin{equation*}
  \partial_j \dice^l(x) = \mathrm{sgn}(x_j) \aux[\rho]'(|x_j|) \neq 0, 
 \end{equation*}  
 that $|x_j| > \rho-2$. 
\end{proof}

We would like to generalise this construction to arbitrary dimension in a 'dice-rotational' manner.
To this end, we need an abstract extension result.

\begin{lemma}\label{Abstract Dice Extension - Lemma}
 For all $\rho > 5$, there is a function $\widetilde{\dice}_{n-1,\rho} \colon \R^{n-1} \rightarrow \R_{\geq 0}$ that satisfies:
 \begin{itemize}
  \item[(i)] The preimages of $\R_{\geq \rho-2}$ under $\widetilde{\dice}_{n-1,\rho}$ and $\dice_{n-1,\rho-2} + \rho -3$ agree. 
  Furthermore, the two functions agree on $\widetilde{\dice}_{n-1,\rho}^{-1}(\R_{\geq \rho - 2})$.
  \item[(ii)] $\diff \, \widetilde{\dice}_{n-1,\rho} \neq 0$ on $\R^n\setminus \{0\}$.
  \item[(iii)] The origin is the unique local minimum of $\widetilde{\dice}_{n-1,\rho}$ with value $0$.
 \end{itemize}
\end{lemma} 
\begin{proof}
  Pick a smooth function $\chi \colon \R \rightarrow [0,1]$ that is monotonically increasing, satisfies $\chi' \leq 1.2$, and that has the following preimages: $\chi^{-1}(\{0\}) = \R_{\leq \rho - 3}$ and $\chi^{-1}(\{1\}) = \R_{\geq \rho - 2}$.
  Pick further $\mathtt{C}^{-1} > (n-1) \cdot (\rho - 2)^2 / (\rho - 3)$ and abbreviate $\dice_{n-1,\rho-2} + (\rho - 3)$ to $\varphi$.
  Define
  \begin{equation*}
      \widetilde{\dice}_{n-1,\rho}(x) \deff \varphi(x) \cdot \chi(\varphi(x)) + \mathtt{C} (1 - \chi(\varphi(x))) \cdot ||x||^2.
  \end{equation*}
  The two functions $\widetilde{\dice}_{n-1,\rho}$ and $\varphi$ agree on $\varphi^{-1}(\R_{\geq \rho - 2})$.
  To prove $(i)$ it is enough to show that $\widetilde{\dice}_{n-1,\rho} < \rho - 2$ on the complement $\varphi^{-1}(\R_{<\rho-2})$ so that the preimage of $\R_{\geq \rho-2}$ under the two functions agree. 
  By the choice of $\mathtt{C}$, we have $\mathtt{C}||x||^2 < \rho - 3$ on $\varphi^{-1}(\R_{\leq \rho -2})$ so that the second summand satisfies 
  \begin{equation*}
      \mathtt{C} (1 - \chi(\varphi(x))) \cdot ||x||^2 < \rho - 3.
  \end{equation*}
  Since $\varphi \geq \rho - 3$, we conclude $\widetilde{\dice}_{n-1,\rho} < \varphi$, whenever $\chi \circ \varphi \neq 0$ and $(i)$ follows.
  
  Clearly, $\widetilde{\dice}_{n-1,\rho} \geq 0$ and the origin gets mapped to zero. Thus, $(iii)$ follows from $(ii)$.
  
  To prove $(ii)$, we calculate
  \begin{align*}
      \begin{split}
        \qquad \diff_x \widetilde{\dice}_{n-1,\rho}  &= \bigl(\chi( \varphi(x)) + (\varphi(x) - \mathtt{C}||x||^2 )\chi'(\varphi(x))\bigr) \diff_x\varphi \\
        &\qquad + 2\mathtt{C} (1-\chi(\varphi(x))) \cdot \sum_{j=1}^{n-1} x_j \cdot \diff x_j
      \end{split}\\ 
      \begin{split}
          &= \sum_{j=1}^{n-1} \Bigl( \bigl(\chi( \varphi(x)) + (\varphi(x) - \mathtt{C}||x||^2) \chi'( \varphi(x)) \bigr) \aux[\rho-2]'(|x_j|) \mathrm{sgn}(x_j)  \\
          & \qquad \qquad + 2\mathtt{C} x_j(1 - \chi(\varphi(x)))  \Bigr) \diff x_j.
      \end{split}
  \end{align*}
  Since $(\varphi - \mathtt{C}||x||^2) \geq 0$ on $\varphi^{-1}(\R_{\leq \rho-2})$, in particular on $\supp(\chi' \circ \varphi)$, the two summands have the same parity.
  The second summand $\mathtt{C}x_j(1-\chi(\varphi(x))$ vanishes only if $x_j = 0$ or $\varphi(x) \geq \rho - 2$.
  In the latter case, the first summand does not vanish.
  Thus, $\diff_x \widetilde{\dice}_{n-1,\rho}$ only vanishes at the origin.
\end{proof}

We use this function to define

\begin{equation*}
 \mathfrak{d}_{n,\rho}(x) \deff \mathtt{d}_{\rho}(x_1,\widetilde{\dice}_{n-1,\rho}(x_2,\dots,x_n) )
\end{equation*} 
as real valued map with domain 
\begin{equation*}
 \mathrm{dom}(\mathfrak{d}_{n,\rho}) = \bigl\{x \in \R^n \, : \, (x_1,\widetilde{\dice}_{n-1,\rho}(x_2,\dots,x_n)) \in \mathtt{d}_{\rho}^{-1}([1,2])\bigr\}.
\end{equation*}
Observe that $\mathfrak{d}_{2,\rho} = \mathtt{d}_\rho$ because $\dice_{1,\rho-2}(x_2) + \rho - 3= |x_2|$ on $\dice_{1,\rho-2}^{-1}(\R_{\geq \rho - 2})$ because $\mathtt{d}_\rho(x_1,x_2) = \mathtt{d}_\rho(x_1,|x_2|)$ and $\mathtt{d}_\rho$ is independent of $x_2$ on $\mathrm{dom}\, \mathtt{d}_\rho \cap \{|x_2| < \rho - 2\}$.

\begin{lemma}\label{Shape Final Dice - Lemma}
 Restricted to $\mathfrak{d}_{n,\rho}^{-1}([1,2])$, the map $\mathfrak{d}_{n,\rho}$ satisfies:
 \begin{itemize}
  \item[(i)] If $\partial_j \mathfrak{d}_{n,\rho}(x) \neq 0$, then $|x_j|>\rho-4$.
  \item[(ii)] Every point in $\mathfrak{d}_{n,\rho}^{-1}([1,2])$ is a regular point.
\end{itemize}  
\end{lemma} 
\begin{proof} 
 We start to prove the first statement.
 For $j=1$, the assumption
 \begin{equation*}
   (\partial_1\mathfrak{d}_{n,\rho})(x) = (\partial_1 \mathtt{d}_\rho)(x_1,\widetilde{\dice}_{n-1,\rho}(x_2,\dots , x_n) ) \neq 0,
 \end{equation*}
 implies $|x_1| > \rho - 2$ by Lemma \ref{dice vanishing at corners - Lemma}.
 
 For $j \geq 2$, the assumption
 \begin{equation*}
  \partial_j\mathfrak{d}_{n,\rho}(x) = (\partial_2\mathtt{d}_\rho)(x_1, \widetilde{\dice}_{n-1,\rho}(\dots)) \cdot \partial_j \widetilde{\dice}_{n-1,\rho}(x_2,\dots,x_n) \neq 0
\end{equation*}  
 implies that the first factors do not vanish, which concludes  $\widetilde{\dice}_{n-1,\rho}(x_2,\dots,x_n) > \rho -2$.
 Since $\widetilde{\dice}_{n-1,\rho} = \dice_{n-1,\rho-2} + (\rho - 3)$ on this domain, we have
 \begin{align*}
  \partial_j \widetilde{\dice}_{n-1,\rho}(x_2,\dots,x_n) &= \partial_j \dice_{n-1,\rho-2}(x_2,\dots,x_n) \\
  &= \mathrm{sgn}(x_j)\aux[\rho-2]'(|x_j|) \neq 0,
 \end{align*}
 which implies $|x_j| > \rho - 4$.
 
 
 
 
 Assume that $x \in \mathfrak{d}_{n,\rho}^{-1}([1,2])$ were a singular point.
 We know from Lemma \ref{Regular Value two dim function - Lemma} that $(x_1,\widetilde{\dice}_{n-1,\rho}(x_2,\dots,x_n))$ is a regular point of $\mathtt{d}_\rho$.
 If $\partial_1 \mathtt{d}$ does not vanish at this point, then $\mathfrak{d}_{n,\rho}(x) \neq 0$ and $x$ would be a regular point for $\mathfrak{d}_{n,\rho}$.
 Thus, $\partial_j \widetilde{\dice}_{n-1,\rho}(x) = 0$ for all $2 \leq j \leq n$, which implies that $(x_2,\dots,x_n)$ is the origin of $\R^{n-1}$.
 From the equation
 \begin{align*}
     0 &= \partial_1 \mathfrak{d}_{n,\rho}(x) = (\partial_1 \mathtt{d}_\rho)(x_1,\widetilde{\dice}_{n-1,\rho}(x_2,\dots,x_n))=(\partial_1 \mathtt{d}_\rho)(x_1,0)
 \end{align*}
 follows $\mathfrak{d}_{n,\rho}(x) = \mathtt{d}_\rho(x_1,0) \notin [1,2]$, which contradicts our domain assumption.
\end{proof}

\begin{lemma}\label{Flow Target Control - Lemma}
 The set $\mathfrak{d}_{n,\rho}^{-1}([1,2])$ does not intersect the closure of the cuboid
 $\quader[n][1,1][\rho-4] = \{x \in \R^{n} \, : \, x_1 \geq -(\rho-4), \, x_j \in [-(\rho-4),\rho-4]\}$.
\end{lemma}
\begin{proof}
 Each point $x \in \quader[n][1,1][\rho-4]$ satisfies $\dice_{n-1,\rho-2}(x_2,\dots,x_n) = 0$.
 Condition $(i)$ of $\widetilde{\dice}_{n-1,\rho}$ yields $\widetilde{\dice}_{n-1,\rho}(x) < \rho - 2$.
 Thus, $x$ is either not in the domain of $\mathfrak{d}_{n,\rho}$ or
 \begin{align*}
    \mathfrak{d}_{n,\rho}(x) &= \dice^l(x_1,\widetilde{\dice}_{n-1,\rho}(x_2,\dots,x_n)) \\
    &= \aux[\rho](|x_1|) + \aux[\rho](|\widetilde{\dice}_{n-1,\rho}(x_2,\dots,x_n)|)\\
    &= 0. 
 \end{align*}
\end{proof}

\begin{definition}
 Let $H \deff \mathfrak{d}_{n,\rho}^{-1}(\{1.9\})$ and let $\Phi \colon H \times [1,2] \rightarrow \mathfrak{d}_{n,\rho}^{-1}([1,2])$ be the restriction of the (shifted) gradient flow of $\mathfrak{d}_{n,\rho}$ so that $\Phi_{1.9} = \id_H$.
\end{definition}

Recall that the (shifted) gradient flow is the unique smooth map 
\begin{equation*}
 \Phi \colon \mathrm{dom} \, \Phi \subset \R^n \times \R \rightarrow \R^n
\end{equation*}
that solves the initial value problem
\begin{equation*}
 (\partial_t\Phi)_{t_0}(x) = \frac{ \grad{ \mathfrak{d}_{n,\rho}}}{||\mathfrak{d}_{n,\rho}||^2}(\Phi_{t_0}(x)) \qquad \text{ and }\qquad \Phi_{1.9}(x) = x. 
\end{equation*}

\begin{proof}[Proof of the Embedding Theorem]
Since $\mathfrak{d}_{n,\rho}$ is non degenerate on $\mathfrak{d}_{n,\rho}^{-1}([1,2])$ there is a diffeomorphism 
\begin{equation*}
 \xymatrix{\Phi \colon H \times [1,2] \ar[rr]^\cong \ar[dr]_{\pr_2} & & \mathfrak{d}_{n,\rho}^{-1}([1,2]) \ar[dl]^{\mathfrak{d}_{n,\rho}} \\
 & [1,2], & }
\end{equation*}
see for example \cite{hirsch1997differential}*{p.153ff} for a proof.
By Lemma \ref{Flow Target Control - Lemma}, the image of $\Phi$ lies in  $\R^n \setminus \quader[n][1,1][\rho-4]$. 
In particular, $H$ is a hypersurface of this open manifold. 

We can describe $\Phi$ even in more details:
Recall from Definition \ref{Block Domains - Def} that $\R^n(\fateps) = \mathrm{Abb}(\mathbf{n}\setminus \mathrm{dom}\, \fateps,\R)$, and that
$\R^\fateps_{\rho-5} = \{x\in \R^{|\mathrm{dom} \fateps|} \, : \, \fateps(i_j) x_j > \rho - 5\}$.
By permuting coordinates, we identify $V_{\rho-4}(\fatalpha)$ with $(-(\rho-4),\rho-4)^{|\mathrm{Null}(\fatalpha)|} \times \R^{\hat{\fatalpha}}_{\rho-5}$.
Lemma \ref{Shape Final Dice - Lemma} implies that $\mathfrak{d}_{n,\rho}$ does not depend on the points parametrised by $\mathrm{Null} (\fatalpha)$, so under this identification, $H \cap V_{\rho-4}(\fatalpha)$ corresponds to $(-\rho+4,\rho-4)^{|\mathrm{Null}(\fatalpha)|} \times H^{\hat{\fatalpha}}$, where $H^{\hat{\fatalpha}} = \mathfrak{d}_{n,\rho}\restrict_{\R^{\hat{\fatalpha}}_{\rho-5}}(\{1.9\})$.
Furthermore, under this identification, the gradient flow $\Phi$ corresponds to $\id \times \Phi^{\hat{\fatalpha}}$.
Note also that $\mathfrak{d}_{n,\rho}(x) = 2-(x_1 - \rho)$ if $\widetilde{\dice}_{n-1,\rho}(x_2,\dots,x_n) > \rho + 3$, so $\grad{ \mathfrak{d}_{n,\rho} } = -e_1$ outside the compact subset $K \deff \widetilde{\dice}_{n-1,\rho}^{-1}(\R_{\leq \rho+3}) \cap \mathfrak{d}_{n,\rho}^{-1}([1,2])$.
Thus, $\Phi_t(x) = x + (t-1.9)e_1$ outside this compact subset. 
Statement 1. is therefore proven. 

Next we prove statement 3.
The gradient of a function is always perpendicular to its level sets, so
\begin{equation*}
 \Phi^\ast \euclmetric_{(x,t)} = h(t)_x + f^2(x,t) \diff t^2,
\end{equation*}
where $h(t)$ is the induced metric $\euclmetric \restrict_{\Phi_t(H)}$ and $f^2 = ||\grad{ \mathfrak{d}_{n,\rho}}||^{-2}$.
On $H \cap V_{\rho-4}(\fatalpha)$, we can use the decomposition of the gradient flow to get the finer decomposition
\begin{equation*}
 \Phi^\ast\euclmetric = \euclmetric\restrict_{\R^n(\hat{\fatalpha})} + h^{\hat{\fatalpha}} + (f^{\hat{\fatalpha}})^2 \diff t^2,
\end{equation*}
where $h^{\hat{\fatalpha}}(t)$ is the induced metric $\euclmetric \restrict_{\Phi^{\hat{\fatalpha}}_t(H^{\hat{\fatalpha}})}$ and $(f^{\hat{\fatalpha}})^2 = ||\grad{ \mathfrak{d}_{n,\rho}}||^{-2}$, which depends only on those coordinates that are parametrised by $\supp \fatalpha$.
Since $\Phi$ is a translation outside of $K$, the maps $h$ and $f$ are constant with values $\euclmetric\restrict_{\R^{n-1}}$ and $1$, respectively.
Thus, statement 3. is proven.

To prove statement 2., we will study the shape of $\mathfrak{d}_{n,\rho}([1,2])$ and how it lies in $\R^n$.
We will consider the case $n=2$ first.

Since $\mathfrak{d}_{2,\rho}^{-1}([1,2]) = \mathtt{d}_\rho^{-1}([1,2])$ it is easy to see that this set separates $\R^2$ into two unbounded regions $\hat{U}_{left}$ and  $\hat{U}_{right}$, see Figure \ref{fig:SepartionOfPlane}.
We set
\begin{equation*}
    U_{left} \deff \hat{U}_{left} \cup \mathfrak{d}_{2,\rho}^{-1}(\R_{>1.9}) \qquad \text{and} \qquad U_{right} \deff \hat{U}_{right} \cup \mathfrak{d}_{2,\rho}^{-1}(\R_{<1.9})
\end{equation*}
so that $H = \mathfrak{d}_{2,\rho}^{-1}(\{1.9\})$.
\begin{figure}[htbp]
    \centering
    \begin{tikzpicture}
      \node at (0,0) {\includegraphics[width=.6\textwidth]{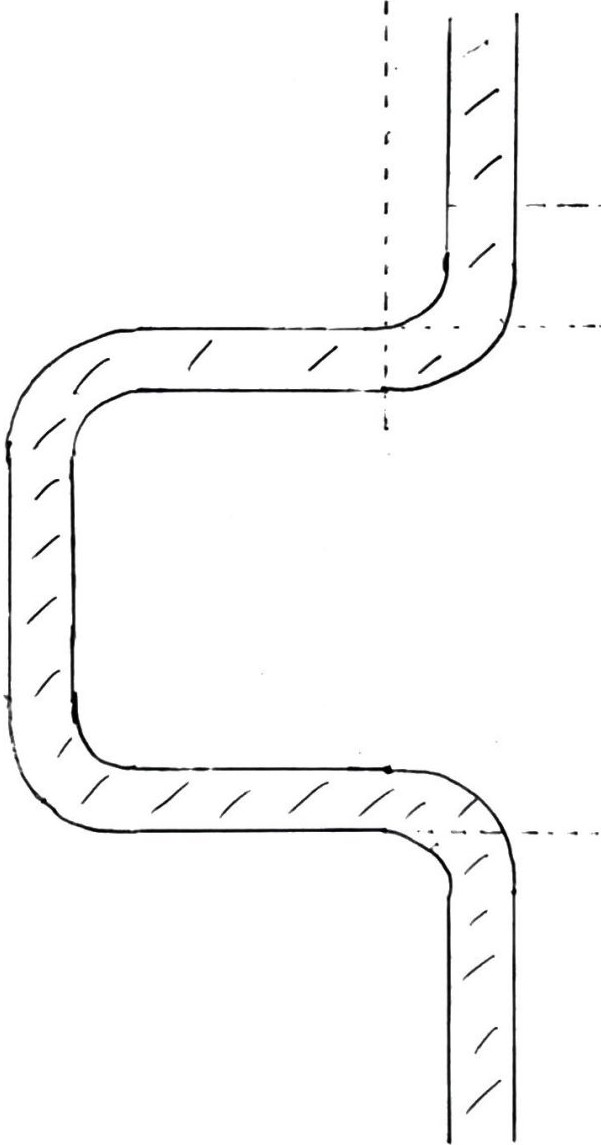}};
      \node at (1.4,2) {$x_1 = \rho-1$};
      \node at (2.2,9.1) {$x_1 = \rho$};
      \node at (3.4,8.6) {$x_1 = \rho + 1$};
      \node at (5.6,3.7) {$x_2 = \rho + 1$};
      \node at (5.6,5.5) {$x_2 = \rho + 3$};
      \node at (5.8,-4.1) {$x_2 = -(\rho + 1)$};
      \node at (-1.4,3.2) {$\dice_{\rho,2}^{-1}([1,2])$};
      \node at (-1.2, -7) {$U_{left}$};
      \node at (3,0) {$U_{right}$};
    \end{tikzpicture}
    \caption{The separation of $\R^2$ by $\mathfrak{d}_{2,\rho}^{-1}([1,2])$.}
    \label{fig:SepartionOfPlane}
\end{figure}
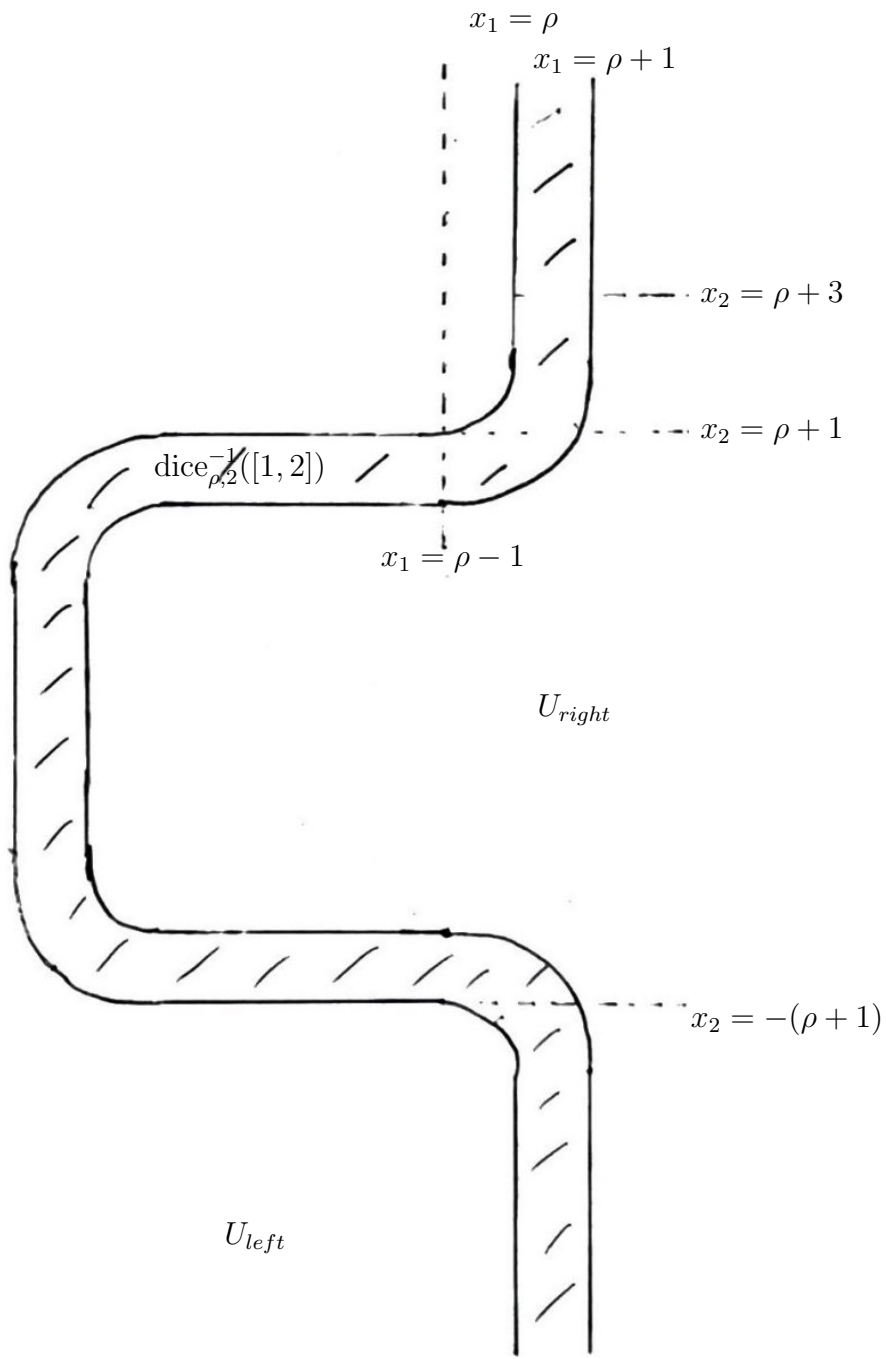
Visual reasons imply the existence of a diffeomorphism $\Psi \colon \mathtt{d}_\rho^{-1}([1,2]) \rightarrow \Psi(\mathtt{d}_\rho^{-1}([1,2]))$ that is the identity near the $\mathtt{d}_\rho^{-1}([1.8,2])$ and outside a compact set, that is of the form $\Psi(x_1,x_2) = (\Psi^{(1)}(x_1),x_2)$ near the line $\R e_1$, and those extension to $U_{left}$ via the identity maps $U_{left} \cup \mathfrak{d}_{n,\rho}^{-1}([1,2])$ to $\{x_1 \leq \rho + 1\}$, see Figure \ref{fig:GraphicalDiffeo}.
\begin{figure}[htbp]
    \centering
    \begin{tikzpicture}
      \node at (0,0) {\includegraphics[width=.9\textwidth]{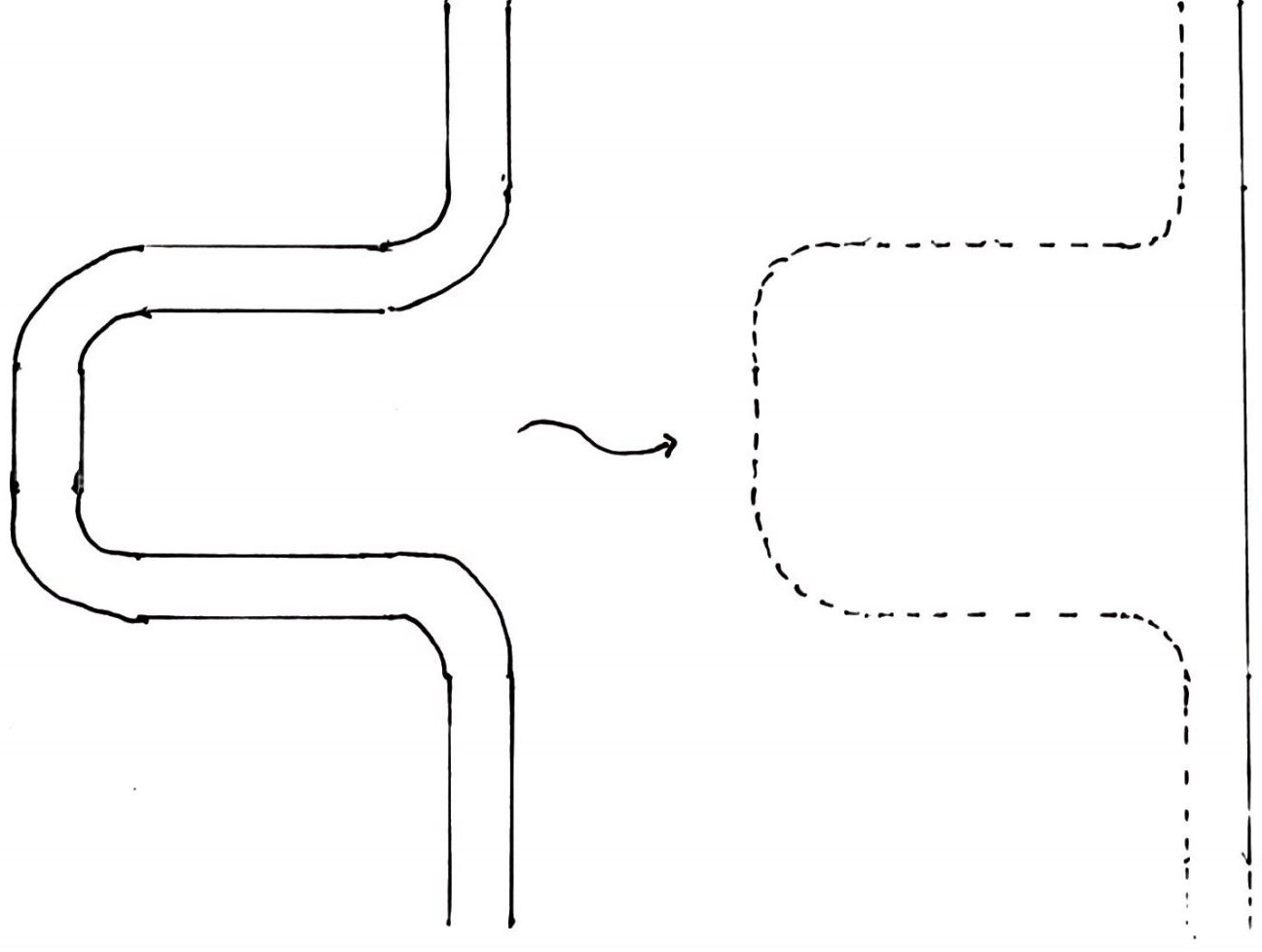}};
      \node at (6.7,5.7) {${x_1 = \rho+1}$};
      \node at (-1.1,5.7) {${x_1 = \rho+1}$};
      \node at (-2,5.2) {$x_1 = \rho$};
      \node at (-4,2.1) {$\mathtt{d}_\rho^{-1}([1,2])$};
      \node at (-4,-4) {$U_{left}$};
      \node at (-3.5,-.5) {$\mathtt{d}_\rho^{-1}(\{1\})$};
      \node at (-3.5,-1.9) {$\mathtt{d}_\rho^{-1}(\{2\})$};
      \node at (3.5,-2.1)  {$\mathtt{d}_\rho^{-1}(\{2\})$};
      \node at (0,0.7) {$\Psi$};
      \node at (-0.2,-.4) {pull it straight!};
    \end{tikzpicture}
    \caption{The diffeomorphism $\Psi$ ``stretches'' $\mathtt{d}_\rho^{-1}([1,2])$ such that $\mathtt{d}_\rho^{-1}(\{1\})$ becomes a straight line and $\mathtt{d}_\rho^{-1}(\{2\})$ remains unchanged.}
    \label{fig:GraphicalDiffeo}
\end{figure}

For higher dimensions, we extend the diffeomorphism in a 'dice-symmetrical' manner.
Let $\Xi \colon \R^{n-1} \setminus \{0\} \times \R_{>0} \rightarrow \R^{n-1} \setminus \{0\}$
be the map that sends $(x,t)$ to $\Xi_t(x)$, the image of $x$ under the (shifted) gradient flow that lies in the levelset $\widetilde{\dice}_{n-1,\rho}^{-1}(\{t\})$.
Define
\begin{align*}
 \begin{split}
     \Psi_n(x) \deff \Bigl( &\Psi^{(1)}(x_1,\widetilde{\dice}_{n-1,\rho}(x_2,\dots,x_n)) ,  \\
      & \ \Xi_{\Psi^{(2)}(x_1,\widetilde{\dice}_{n-1,\rho}(x_2,\dots,x_n))}(x_2,\dots,x_n) \Bigr),
 \end{split}
\end{align*} 
where $\Psi^{(j)}$ denotes $\pr_j \circ \Psi$.
Its inverse is given by 
\begin{align*}
 \begin{split}
     \Psi_n^{-1}(y) \deff \Bigl( &(\Psi^{-1})^{(1)}(y_1,\widetilde{\dice}_{n-1,\rho}(y_2,\dots,y_n), \\
     & \ \Xi_{(\Psi^{-1})^{(2)}(y_1,\widetilde{\dice}_{n-1,\rho}(y_2,\dots,y_n))}(y_2,\dots,y_n)\Bigr).
 \end{split}
\end{align*}
The two maps are smooth and inverse to each other.
As 
\begin{equation*}
    \Xi_{\Psi^{(2)}(x_1,\widetilde{\dice}_{n-1,\rho}(x_2,\dots,x_n))}(x_2,\dots,x_n) = (x_2,\dots,x_n)
\end{equation*}
in a neighbourhood of $\R e_1$,
this assignment extends to a smooth diffeomorphism
\begin{equation*}
 \Psi_n \colon \mathrm{dom}( \mathfrak{d}_{n,\rho}) \rightarrow \{x \in \R^n \, : \, (x_1,\widetilde{\dice}_{n-1,\rho}(x_2,\dots,x_n)) \in \mathrm{im} \Psi(\mathrm{dom}(\mathtt{d}_{\rho}))\}.
\end{equation*}
Since $\Psi$ is the identity on $\mathtt{d}_\rho^{-1}([1.8,2])$ and outside a compact set, the diffeomorphism $\Psi_n$ is the identity on $\mathfrak{d}_{n,\rho}^{-1}([1.8,2])$ and outside a compact set.
We extend it via the identity 'to the left' to $U_{left} \deff \{x \in \R^n \, : \, (x_1,\widetilde{\dice}_{n-1,\rho}(x_2,\dots,x_n)) \in U_{left}\}$. 
This gives a diffeomorphism
\begin{equation*}
 \Psi_n \colon U_{left} \cup \mathfrak{d}_{n,\rho}^{-1}([1,2]) \rightarrow \{x_1 \leq \rho + 1\}
\end{equation*}
with the same properties as $\Psi = \Psi_2$.
Note that this particularly implies that $H$ is diffeomorphic to $\R^{n-1}$.
\end{proof}

 \section{Kan Property of the Concordance Set}\label{Section - Kan Property of the Concordance Set}

The aim of this section is to prove Theorem \ref{KanSet - MainThm}, which can be thought of as a can opener to the combinatorial world. 
The precise formulation is the following one:

\begin{theorem}\label{Concordance Set is Kan - Theorem}
 The concordance set $\ConcSet$ is a Kan set for every closed psc manifold $M$.
\end{theorem}

The proof requires the following technical ingredient.

\begin{lemma}[Metric Modification]\label{Metric Modification - Lemma}
 Let $g$ be a block metric on $M \times \R^n$ that decomposes outside of $M \times \rho/2I^n$ for $\rho > 20$. Assume further that $\face[\omega]{j}g$ is a psc metric for all $(j,\omega) \neq (1,1)$.
 Let $H$ and $\Phi$ be as in the Embedding Theorem \ref{Embedding Theorem}. 
 and let ${}_R\Phi \deff R \cdot \Phi \circ (\id \times R^{-1}) \colon H \times [R,2R] \rightarrow M \times \R^n$.
  
 Then, for all $R\geq 1$, there is a metric on $G_R$ on $M \times R\Phi(H \times [1,2])$ that satisfies:
 \begin{itemize}
  \item[(o)] There is an $R_0 > 0$ such that $\scal(G_R) > 0$ if $R>R_0$.
  \item[(i)] The pull back $({}_R\Phi)^\ast G_R$ and $({}_R\Phi)^\ast g$ agree on $M \times H \times [1.8R,2R]$ and on the complement of a compact subset of $M \times H \times [R,2R]$.
  \item[(ii)] The pull back $({}_R\Phi)^\ast G_R$ is a product metric on $M \times H \times [R,1.2R]$. More precisely,
  \begin{equation*}
   ({}_R\Phi)^\ast G_R = ({}_R\Phi)^\ast G_R\restrict_{M \times H \times \{R\}} \oplus \diff t^2.
  \end{equation*}
 \end{itemize}
\end{lemma}
\begin{proof}
 Recall the open subsets
 \begin{equation*}
  V_\rho(\fatalpha) \deff \{x \in \R^n \, : \, \fatalpha(i)x_i > \rho - 1 \text{ if } i \in \supp \fatalpha; x_i \in (-\rho,\rho) \text{ otherwise}\},
 \end{equation*}
 for all $\fatalpha \colon \{1,\dots, n\} \rightarrow \Z_3$. 
 Also recall the notation ${\hat{\fatalpha}} =\fatalpha |_{\supp \fatalpha}$.
 We  have the following finite decomposition
 \begin{equation*}
  \R^n\setminus \quader[n][1,1][\rho] = \bigcup \{V_\rho(\fatalpha) \, : \, \fatalpha \neq \mathbf{0}, (1,0,\dots,0)\}.
 \end{equation*}
 Let $H$ and $\Phi$ be as in the Embedding Theorem. 
 Then the subset $R\Phi(H \times [1,2])$ is still in the complement of $\quader[n][1,1][\rho]$ for all $R\geq 1$.
 The block metric $g$ decomposes on $M \times V_{R\rho}({\fatalpha})$ into $\face[{\hat{\fatalpha}}]{ }g \oplus \euclmetric$. 
 Together with Embedding Theorem \ref{Embedding Theorem} 
 this implies
 \begin{align*}
  {}_R\Phi^\ast(g|_{V_{R\rho}(\fatalpha})) &= ( ( \id \times (\Phi^{\hat{\fatalpha}} \circ (\id \times R^{-1})) ) )^\ast R^\ast(\face[\hat{\fatalpha}]{ }g \oplus \euclmetric ) \\
   &= R^\ast \face[\hat{\fatalpha}]{ }g \oplus R^2 (\id \times R^{-1})^\ast {\Phi^{\hat{\fatalpha}}}^\ast \euclmetric \\
  &= R^\ast \face[\hat{\fatalpha}]{{ }} g \oplus R^2 h^{\hat{\fatalpha}}(\placeholder,R^{-1}\cdot) + (f^{\hat{\fatalpha}})^2(\placeholder, R^{-1} \cdot ) \diff t^2  \\
  &= R^\ast \face[\hat{\fatalpha}]{ } g \oplus R^2{}_R h^{\hat{\fatalpha}} + {}_R(f^{\hat{\fatalpha}})^2 \diff t^2\\
  &=: h_R^{\hat{\fatalpha}} + f_R^2 \diff t^2
 \end{align*}   
  These decompositions assemble to the (global) decomposition
 \begin{equation*}
  ({}_R\Phi)^\ast g = h_R + f^2_R \diff t^2.
 \end{equation*}
 
 Define the auxiliary metric $k_{\mathrm{aux},R} \in \mathrm{Riem}(M \times H \times [R,2R])$ via 
 \begin{equation*}
  k_{\mathrm{aux},R}(m,x,t) \deff h_R(m,x,1.1\cdot R) + \diff t^2
 \end{equation*}
 and pick a smooth function $\chi \colon [1,2] \rightarrow [0,1]$ that vanishes identically on $[1,1.2]$ and that is identically $1$ on $[1.8,2]$.
 Set further ${}_R\chi(t) = \chi(R^{-1}t)$.
 
 We define $k_{R} \in \Riem(M \times H \times [R,2R])$ via
 \begin{equation*}
  k_{R} \deff (1 - {}_R\chi) \left( k_{\mathrm{aux},R} +  \diff t^2\right) + {}_R\chi \left( h_{R} + f^2_R \diff t^2\right)
 \end{equation*}
 and define 
 \begin{equation*}
     G_R \deff ({}_R\Phi)_\ast(k_R).
 \end{equation*}
  
  In the following, we drop the lower $R$ if $R=1$.
  To prove statement $(i)$, we claim that metric agrees with ${}_R\Phi^\ast(g)$ on $M \times H \times [1.8R,2R]$ and on the complement of some compact set so that $G_R$ agrees with $g$ on $M \times {}_R\Phi(H \times [R,1.2\cdot2])$ and away from a compact subset.
  The first claim follows immediately from $\chi \equiv 1$ on $[1.8,2]$. 
  The second claim follows from the existence of a compact set $C$ such that $H \setminus C$ is an open submanifold of $\{x_1 = \rho + 1\}$ and such that $\Phi$ restricts to $(x,t) \mapsto x - (t-1.9)e_1$ on $H \setminus C \times [1,2]$, so that, on $H \setminus C \times [R,2R]$, the pull back ${}_R\Phi^\ast g$ decomposes into
  \begin{equation*}
    ({}_R\Phi)^\ast g = (R\Phi \circ \id \times R^{-1})^\ast(\face[1]{1}g \oplus \diff x_1^2) = R^\ast(\face[1]{1}g) +  \diff t^2.
  \end{equation*}
  Thus, on $H\setminus C \times [R,2R]$, the curve $h_R \equiv R^\ast(\face[1]{1}g)$ and the map $f_R \equiv 1$ are constant so that $k_R = ({}_R\Phi)^\ast g$ there.
  
  To prove $(o)$, we claim that $k_R$ (and hence $G_R$) is a psc metric, provided that $R$ is sufficiently large.
  Indeed, on $V_{R\rho}({\fatalpha})$ with $\fatalpha \notin \{\mathbf{0}, (1,0,\dots,0)\}$ the metric $k_R$ decomposes into
 \begin{align*}
  k_R &= R^\ast \face[{\hat{\fatalpha}}]{ }g \oplus \left( (1 -{}_R\chi)(h_R(\cdot,1.1R) + \diff t^2) + {}_R\chi ( h_R + f_R^2 \diff t^2) \right)|_{V_{R\rho}(\fatalpha) \times [R,2R]} \\
  &= R^\ast \face[{\hat{\fatalpha}}]{ }g \oplus \left( (1 - {}_R\chi)(R^2{}_Rh^{\hat{\fatalpha}}(\cdot,1.1R) + \diff t^2) + {}_R\chi (R^2{}_Rh^{\hat{\fatalpha}} + ({}_Rf^{\hat{\fatalpha}})^2 \diff t^2) \right) \\
  &= R^\ast \face[{\hat{\fatalpha}}]{ }g \oplus R^2{}_R\left( (1 - \chi) h^{\hat{\fatalpha}}(\cdot,1.1R) + \chi  h^{\hat{\fatalpha}} \right) + {}_R\left((1-\chi) + \chi (f^{\hat{\fatalpha}})^2\right)\diff t^2
 \end{align*}
 Since $\Phi_t(x) = x - (t-1.9)e_1$ outside a compact set $C$, the metrics
  $k_1$ and $\Phi^\ast g$ agree outside of $C$ and $h_R$ is a psc metric there.
  Thus, the metrics $k_{\mathrm{aux},R}$ and $h_R^{\hat{\fatalpha}}$ agree outside $C$ for all $R\geq 1$ and have positive scalar curvature there. 
 Furthermore, for $R=1$, we find  positive constants $\mathrm{Const}({\fatalpha})$ such that
 \begin{equation*}
  \left|\scal \left( (1 -\chi)(h^{\hat{\fatalpha}}(\cdot,1.1) + \diff t^2) + \chi (h^{\hat{\fatalpha}} + (f^{\hat{\fatalpha}})^2 \diff t^2) \right)\right| \leq \mathrm{Const}(\fatalpha)
 \end{equation*}
 because the domain where the two metrics disagree is compact.
 
 Corollary \ref{Rescaling Behaviour - Cor} implies that 
 \begin{equation*}
    \scal\left(R^2{}_R\left( (1 - \chi) h^{\hat{\fatalpha}}(\cdot,1.1) + \chi  h^{\hat{\fatalpha}} \right) + {}_R\left((1-\chi) + (f^{\hat{\fatalpha}})^2\right)\diff t^2\right) \xrightarrow{R \to \infty} 0. 
 \end{equation*}
 Furthermore, $\scal(R^\ast \face[{\hat{\fatalpha}}]{ } g)= \scal(\face[{\hat{\fatalpha}}]{ } g) \circ R$ is bounded from below by a positive constant that is independent of $R$, for $\face[{\hat{\fatalpha}}]{ } g$ is a block metric with positive scalar curvature. 
 Putting these facts together, we can choose $R$ sufficiently large such that $\scal(k_R)$ is positive on $(H \cap V_{\rho-4}(\fatalpha)) \times [R,2R]$.
Since all $H \cap V_{\rho-4}(\fatalpha)$ form a finite cover of $H$, we find an $R_0>0$ such such that $\scal(k_R) > 0$ for all $R > R_0$.   

A choice for the desired metric is given by the unique metric $G_R$ on $M \times {}_R\Phi(H \times [R,2R])$ that satisfies $({}_R\Phi)^\ast G_R = k_R$. 
\end{proof}

We can now prove the Kan property.

\begin{proof}[Proof of Theorem \ref{Concordance Set is Kan - Theorem}]
 Assume we are given a cubical $n$-horn $\CubeBox$ that is represented by the following set of block metrics 
 \begin{equation*}
  \{g_{(j,\omega)} \in \ConcSet[n-1] \, : \, \face[\omega]{j}g_{(k,\eta)} = \face[\eta]{k-1} g_{(j,\omega)} \text{ for } \, j<k; (j,\omega), (k,\eta) \neq (i,\eps) \}.
 \end{equation*}
 The problem is symmetric in $(i,\eps)$ so we may assume that $(i,\eps) = (1,1)$.
 We pick a sufficient large $\rho > 20$ such that all given block metrics $g_{(j,\omega)}$ decompose outside of $M \times \rho/2I^{n-1}$. 
 By assumption, the block metrics $g_{(j,\omega)}$ have matching faces, so the Riemannian metrics $\degen{j} g_{(j,\omega)}|_{\{\omega x_j > \rho - 4\}}$ agree on the intersection on their domains, see Lemma \ref{Metric Patching - Lemma}.
 Thus, the union of all $\degen{j} g_{(j,\omega)}|_{\{\omega x_j > \rho - 4\}}$ with $(j,\omega) \neq (1,1)$ forms a smooth psc metric $g$ on $M \times \R^n\setminus \quader[n][1,1][\rho-4]$.
 
 By the proof of Proposition \ref{BlockMetric Combi Contrac - Prop}, we can find a block metric on $M \times \R^n$ that agrees with $g$ on $M \times \R^n \setminus \quader[n][1,1][\rho-2]$ as follows:
 First, the proof shows that we can extend $g\restrict_{\{x_1 = c\}}$ to a block metric on $\{x_1 = c\}$ that agrees with $g$ on $M \times \{c\}\times \R^{n-1} \setminus (\rho-2)I^{n-1}$ for all $c > \rho -4$. 
 As $g\restrict_{\{x_1 = c\}}$ does not depend on $c$, the extension can be chosen to be independent of $c$.
 These extensions assemble to a smooth metric $g_{\mathrm{aux}}$ on $M \times \R^n \setminus (\rho - 4)I^n$.
 Applying the proof of Proposition \ref{BlockMetric Combi Contrac - Prop} one more time, yields a block metric that agrees with $g_\mathrm{aux}$ outside of $(\rho-2)I^n$, in particular, it agrees with $g$ on $M \times \R^n \setminus \quader[n][1,1][\rho-2]$.
 We denote this extension again with $g$.
 Of course, we do not claim that the scalar curvature of this extension of $g$ is positive on the whole manifold $M \times \R^n$.
 
 Let $H$, $\Phi$, and $U_{left} \subset \R^n\setminus \quader[n][1,1][\rho-4]$ be as in the Embedding Theorem \ref{Embedding Theorem}.
 By the Metric Modification Lemma \ref{Metric Modification - Lemma},
 we find, for all sufficient large $R$, a psc metric $G_R$ on $M \times R\Phi(H \times [1,2])$ that agrees with $g$ outside a compact subset of $M \times R\Phi(H \times [1,2])$ and that can be extended on $M \times \left(R \cdot U_{left} \setminus R\Phi(H \times [1,2])\right)$ to a smooth psc metric via $g$.
 This extension is again denoted by $G_R$.
 The set $U_{left} \cup \Phi(H \times [1,2])$ is diffeomorphic to $\{x_1 \leq \rho + 1\}$ via $\Psi_n$, by Embedding Theorem \ref{Embedding Theorem}, 2.,
  the set $R\cdot U_{left} \cup R \Phi(H \times [1,2])$ is diffeomorphic to $\{x_1 \leq (\rho+1) R\}$ via ${}_R\Psi_n \deff R \Psi_n(R^{-1} \cdot)$.
 Since ${}_R\Phi^\ast G_R$ is a product metric on $H \times [R,1.2R]$, the metric $({}_R\Psi_n^{-1})^\ast G_R$ is also collared with collar length $0.2R$.
 This means that there is a (necessarily unique) collar map ${}_R\kappa \colon M \times \{x_1 = (\rho+1)R\} \times (-0.2R ,0] \rightarrow M \times \{x_1 \leq (\rho+1)R \}$ such that the ${}_R\kappa^\ast ({}_R\Psi_n^{-1})^\ast G_R$ is a product metric.
 
 The collar ${}_R\kappa$ might not be the standard one $(m,x,t) \mapsto (m,x + te_1)$ but they agree outside a compact set. 
 This follows from the fact that ${}_R\Psi_n$ is the identity away from a compact set, so $({}_R\Psi_n^{-1})^\ast G_R$ agrees with $g$ away from a compact set.
 Since the set of collars that agree outside a fixed compact set is contractible 
 we find a diffeomorphism ${}_R\Theta \colon M \times \{x_1 \leq (\rho + 1)R\} \rightarrow M \times \{x_1 \leq (\rho + 1)R\}$ that is the identity outside a compact set and that maps the standard collar to ${}_R\kappa$.
 In fact, the proof of this statement allows us to assume that ${}_R\Theta \circ {}_R\kappa$ agrees with the standard collar on $M \times \{x_1 = (\rho +1)R\} \times (-R/10,0]$.
 
 We can now extend the metric ${}_R\Theta^\ast ({}_R\Psi_n^{-1})^\ast G_R$ from $\{x_1 \leq (\rho + 1)R\}$ to $M \times \R^n$ with the product metric $({}_R\Psi_n^{-1})^\ast G_R \restrict_{\{x_1 = (\rho + 1)R\}} \oplus \diff x_1^2$.
 The resulting metric $g_{\mathrm{fill}}$ is a block metric that has positive scalar curvature.
 It further satisfies $\face[\omega]{j} g_{\mathrm{fill}} = g_{(j,\omega)}$ for all $(j,\omega)$ because it agrees with $g$ on the complement of a sufficient large cuboid $\quader[n][1,1][?\cdot R]$.
 
 Thus, $g_{\mathrm{fill}}$ is a filler for given $n$-horn, so $\ConcSet$ is a Kan set. 
\end{proof}

 \section{Geometric Addition}\label{Geometric Addition - Section}

Since $\ConcSet$ is a Kan set, the (combinatorial) homotopy groups are well defined and carry a (combinatorical) group structure.
These group structures are formally defined through the Kan condition and, as we do not keep track of the fillers, they are only defined on homotopy classes.
This raises the question, whether there is a geometrical meaningful representative for it.

Theorem \ref{Geometric Addition - Thm} gives answers this question.
Informally, the theorem says that the addition on $\pi_n(\ConcSet)$ is given by gluing two block metrics together along one of their matching faces.

\subsubsection*{Angel Rotation}

We start with describing a construction that extends a metric on $N \times \R_{\geq 0}$ to $N \times \R_{\geq 0}^2$ by ``rotating it around the origin''. 
In the following, let $N$ be a smooth, not necessarily closed manifold.
\begin{definition}
  Denote by $\Shift_{j,R}$ 
  \nomenclature{$\Shift_{j,R}$}{Shift operator: Push-foward with diffeomorphism $x \mapsto x + Re_j$}
  the shift diffeomorphism on $\R^n$ given by $x \mapsto x + Re_j$.
  If $g$ is a Riemannian metric on $N \times \R^n$, then we define
  \begin{equation*}
      \Shift_{j,R}(g) \deff {\Shift_{j,R}}_\ast(g) = \Shift_{j,-R}^\ast(g).
  \end{equation*}
\end{definition}

\begin{proposition}\label{Construction-AngleRotation - Prop}
  Let $g$ be a Riemannian metric on $N \times \R_{\geq 0}$ that has a product structure away from $N \times [r_1,r_2]$. 
  We can construct a metric $g_R^\angle$ 
  \nomenclature{$g_R^\angle$}{metric obtained from $g$ by ``angle-rotation'' with rotation speed $1/R$}
  on $N \times \R_{\geq 0}^2$ for all $R \geq 1$ with the following properties:
  \begin{itemize}
      \item[(o)] The metric $g^\angle_R$ is a product metric away from some compact subset.
      \item[(i)] If $g$ is a psc metric, then there is a $R(g) > 0$ such that $g^\angle_R$ is a psc metric for all $R > R(g)$.
      \item[(ii)] $g^\angle_R$ restricts to $\Shift_{R(g)}$ on $N \times \{0\} \times \R_{\geq 0}$ and on $N \times \R_{\geq 0} \{0\}$ under the canonical identification with $N \times \R_{\geq 0}$.
      \item[(iii)] $g^\angle_R$ is a product metric on $\{x_j \leq R / 2\}$ and the product structures are compatible on the intersection.
      \item[(v)] If $g$ is the suspension of $h \colon \R_{\geq 0} \rightarrow \Riem^{(+)}(N)$, then $g^\angle_R$ is the suspension of a function $\mathfrak{h}_R \colon \R_{\geq 0}^2 \rightarrow \Riem^{(+)}(N)$.
      Furthermore, if $h$ is constant so is $\mathfrak{h}_R$.
  \end{itemize}
\end{proposition}

\begin{proof}
Let $K$ be a compact, convex, point symmetric body inside $I^2$ with smooth boundary $\partial K$. 
We further assume $K$ agrees with $I^2$ inside $\{|x_j|\leq  1/2\}$.
Let $\gamma \colon [0,l] \rightarrow \partial K \cap (\R_{\geq 0})^2$ be the surjective curve that is parametrised by arclength and satisfies $\gamma(0) = e_1$. 
This particularly implies that 
\begin{align*}
     \gamma(\varphi) = \begin{cases}
      e_1 + \varphi e_2, & \text{on } [0,1/2], \\
      e_2 + (l-\varphi)e_2, & \text{on } [l- 1/2, l].
     \end{cases}
\end{align*}
Let $v \colon [0,l] \rightarrow \R^2$ be the unique normalised vector field that is perpendicular to $\gamma'$ and satisfies $\det(v,\gamma')>0$.
In particular, $v(0) = e_1$.
Furthermore, let $\kappa$ be the curvature of $\gamma$.
Recall that the curvature\footnote{The sign appears because we chose $v$ such that $(v,\gamma')$ is positively oriented.} of $\gamma$ is defined by $\gamma''(\varphi) = - \kappa(\varphi) v(\varphi)$.
The curvature is non-negative as $\gamma$ parametrises the boundary of a convex body.

If we replace $K$ by $R \cdot K$, then the corresponding objects $\gamma_R$, $v_R$, and $\kappa_R$ have $[0,R\cdot l]$ as domain and satisfy
\begin{equation*}
    \gamma_R(\varphi) = R\gamma(R^{-1} \varphi), \quad v_R(\varphi) = v(R^{-1}\varphi), \quad \kappa_R(\varphi) = R^{-1} \kappa(R^{-1}\varphi).
\end{equation*}

For all $R \geq 1$, this datum defines a diffeomorphism
\begin{align*}
    \Xi \colon (0,\infty) \times [0,R\cdot l] &\rightarrow (\R_{\geq 0})^2 \setminus R\cdot K \\
    (r,\varphi) &\mapsto \gamma_R(\varphi) + r \cdot v_R(\varphi),
\end{align*}
which can be thought as a variant of polar coordinates on the quadrant $\R^2_{\geq 0}$.
Pulling back the euclidean metric yields
\begin{equation*}
    \bigl(\Xi_R^\ast \euclmetric\bigr)_{(r,\varphi)} = \diff r^2 + (1+r\kappa_R(\varphi))^2 \diff \varphi^2.
\end{equation*}

This motivates the following construction, which is visualised in Figure \ref{fig:AngleRotation}.
\begin{figure}[htbp]
    \centering
    \includegraphics[width=\textwidth]{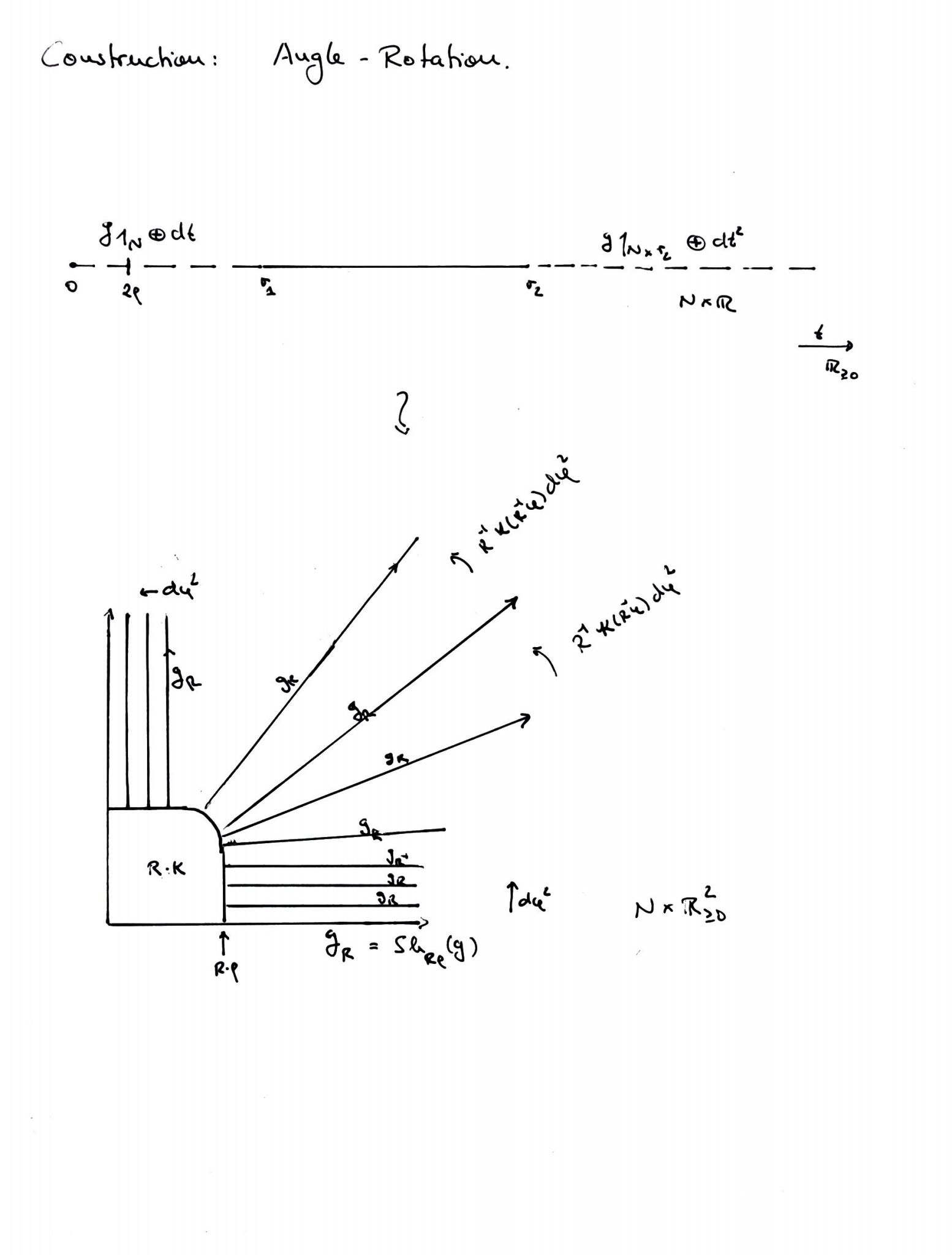}
    \caption{The metric $g$ gets rotates around the smooth, convex body $R\cdot K$.}
    \label{fig:AngleRotation}
\end{figure}
Given $g \in \Riem(N \times \R_{\geq 0})$ that is a product metric on $N \times \lbrack 0, \eps)$, for some $\epsilon > 0$.
We extend $g$ via the product structure to $N \times \R$.
Define $g_R^\angle$ to be unique Riemannian metric on $N \times (\R_{\geq 0})^2 \setminus R\cdot K^\circ$ that satisfies
\begin{align*}
    \bigl((\id \times \Xi_{R})^\ast g_R^\angle\bigr)_{(n,r,\varphi)} = g_{(n,r)} + (1+r\kappa_R(\varphi))^2 \diff \varphi^2.
\end{align*}
Since $g = g\restrict_N \oplus \diff r^2$ away from $M \times [r_1,r_2]$, we conclude that
\begin{equation}\label{Eq: Defining rotated metric}
  (\id \times \Xi_R)^\ast g_R^\angle = g\restrict_N \oplus (\diff r^2 + (1+r\kappa_R(\varphi))^2 \diff \varphi^2) = g\restrict_N \oplus \Xi_R^\ast\euclmetric
\end{equation}
away from $M \times [r_1,r_2] \times [0,R\cdot l]$.

As a consequence, $(\id \times \Xi_R)^\ast g_R^\angle$ has positive scalar curvature away from $M \times [r_1, r_2]\times[0,R\cdot l])$, if $g$ has positive scalar curvature.
 Since the $r$-coordinate is bounded on $M \times [r_1, r_2]\times[0,R\cdot l]$, we can choose $R$ sufficiently large such that the 2-jets of $(\id \times \Xi_R)^\ast g_R^\angle$ and the product metric $g \oplus \diff \varphi$ are arbitrarily close on $M \times [r_1, r_2]\times[0,R\cdot l]$
 This implies that $g_R^\angle$ is a psc metric if $R$ is sufficient large.

The metric $g_R^\angle$ is a product metric on $\{x_j \leq R/2\}$.
Indeed, on $\Xi_R^{-1}(\{x_1 \leq R/2\})$ we have $\Xi_R(r,\varphi) = (R +r)e_2 -  \varphi e_1$, so that $\CubeProj{1} \circ (\id \times \Xi_R) = \Shift_R \circ \CubeProj{1}$, where $\CubeProj{1}$ projects the first component away.  
Thus, the identity
\begin{equation*}
    g_R^\angle |_{\{x_1 \leq R/2\}} = \Shift_{R\, \ast}(g) \oplus \diff x_1^2 = \Shift_{R\, \ast}(\CubeProj{1}^\ast g)  + \diff x_1^2
\end{equation*}
follows from the calculation
\begin{align*}
    (\id \times \Xi_R)^\ast(\Shift_{R\, \ast}(g) \oplus \diff x_1^2) &= (\id \times \Xi_R)^\ast(\CubeProj{1}^\ast \Shift_{-R}^\ast(g) + \diff x_1^2) \\
    &= \CubeProj{1}^\ast g + (\id \times \Xi_R)^\ast \diff x_1^2 \\
    &= g \oplus \diff \varphi^2 = (\id \times \Xi_R)^\ast g_R^\angle.
\end{align*}

We use Equation (\ref{Eq: Defining rotated metric}) to extend $g_R^\angle$ from $N\times \R_{\geq 0}^2 \setminus RK^\circ$ to $N \times (\R_{\geq 0})^2$ via $g\restrict_N + \diff x_1^2 + \diff x_2^2$ and denote the result again with $g_R^{\angle}$. 
The extension satisfies $(o)$, $(i)$, $(ii)$, and $(iii)$ by construction.

Lastly, if $g = \susp(h)$ is the suspension of a curve of metrics, then $g_R^\angle$ is the suspension of 
\begin{align*}
    \mathfrak{h}_R \colon (x_1,x_2) \mapsto \begin{cases}
     h(\pr_1(\Xi_R^{-1}(x_1,x_2))), & \text{if } (x_1,x_2) \in (\R_{\geq 0})^2 \setminus R\cdot K, \\
     h(0), & \text{if } (x_1,x_2) \in (\R_{\geq 0})^2 \cap R\cdot K.
    \end{cases}
\end{align*}
The map is well defined and smooth, as $h$ is constant near zero and thus $\Shift_R(h) = h(\cdot - R)$ can be constantly extended to $[0,R]$.
As before, one checks that $(\id \times \Xi_R)^\ast \susp(\mathfrak{h}_R) = (\id \times \Xi_R)^\ast g_R^\angle = h + \diff r^2 + (1 + r\kappa_R(\varphi))^2 \diff \varphi^2$ and claim $(iv)$ follows.
\end{proof}

\subsubsection*{Geometric Addition on the Concordance Set}

We will use the previous extension construction to give a geometrical description of the group structures.

\begin{definition}
  For $1 \leq j \leq n$, let $g, h \in \BlockMetrics[n]$ be two block metrics with matching faces, $\mathrm{i.e.}$, $\face[1]{j}g = \face[-1]{j}h$.
  If $g$ and $h$ decompose away from $M \times RI^n$, then we set
  \begin{equation*}
      g +_{j,R} h = \Shift_{j,-2R}(g)|_{\{x_j \leq 0\}} \cup \Shift_{j,2R}(h)|_{\{x_j\geq 0\}}.
      \nomenclature{$ g +_{j,R} h$}{Geometric addition in the $j$-th coordinate}
  \end{equation*}
\end{definition}
It is easy to see that $g +_{j,R}h$  is a block metric whose faces are given by
\begin{equation*}
    \face{i}(g +_{j,R}h) = \begin{cases}
      \face[-1]{j} g, & \text{if } i=j, \eps = -1, \\
      \face[1]{j} h, & \text{if } i=j, \eps = 1, \\
      \face{i}g +_{j-1,R} \face{i}h, & \text{if } i < j, \\
      \face{i}g +_{j,R} \face{i}h, & \text{if } i > j.
    \end{cases}
\end{equation*}
Furthermore, $g+_{j,R}h \in \ConcSet[n]$ has positive scalar curvature, if $g$ and $h$ have positive scalar curvature.
Note further, that if $g = \susp(\Bar{g})$ and $h = \susp(\Bar{h})$ are suspension of block maps, then 
\begin{equation*}
   g+_{j,R}h = \susp\bigl(\Bar{g}(\cdot + 2Re_j)|_{\{x_j \leq 0\}} \cup \Bar{h}(\cdot - 2Re_j)|_{\{x_j \geq 0\}}  \bigr)    
\end{equation*}
is again a suspension of a block map.

These operations represent the addition on the combinatorial homotopy groups. 
\begin{theorem}\label{Geometric Addition - Thm}
  There are $n$-many group structures $+_j$ on $\pi_n(\ConcSet,g_0)$ that are defined as follows:
  If $g, h$ are representatives of $[g], [h] \in \pi_n(\ConcSet)$ that decompose away from $M \times R I^n$, then 
  \begin{equation*}
      [g] +_j [h] \deff [g +_{j,R} h].
  \end{equation*}
  If $r_j$ denotes the reflection at the hyperplane $\{x_j = 0\}$, then $[r_j^\ast g]$ is the inverse element of $[g]$ with respect to $+_j$.
  The structures are Eckmann-Hilton related to each other.
  Furthermore, $+_1$ agrees with the group structure provided by cubical set theory.
\end{theorem}
\begin{proof}
 We start with the proof that the operations are well defined.
 Two different choices $R_0$, $R_1$ yield isotopic metrics.
 Indeed, an isotopy between $g +_{j,R_0} h$ and $g +_{j,R_1} h$ is given by $s \mapsto g +_{j,R_s} h$, where $R_s \deff (1-2)R_0 + sR_1$.
 Since isotopy implies concordance (after a reparametrisation), the two block metrics represent the same element in $\pi_n(\ConcSet)$.
 
 Let $G \in \ConcSet[n+1]$ be a homotopy between $g_{-1}$ and $g_1$ and let $H \in \ConcSet[n+1]$ be a homotopy between $h_{-1}$ and $h_1$.
 Assume that $G$ and $H$ decompose outside of $M \times RI^n$.
 Then $G +_{j+1,R} H \in \ConcSet[n+1]$ is a homotopy bewteen $g_{-1} +_{j,R} g_1 $ and $h_{-1} +_{j,R} h_1$.
 
 The two previous observations imply that $+_j$ is a well defined operation on $\pi_n(\ConcSet,g_0)$.
 
 The proof of associativity agrees with the proof that the group structure of homotopy groups is associative.
 The neutral element is given by the base point $g_0$.
 
 If $g$ decomposes outside of $M \times \rho I^n$, so does $r_j^\ast g$.
 The homotopy between $r_j^\ast g +_{j,\rho} g$ is indicated in Figure \ref{fig:ReflectionInverse}.
 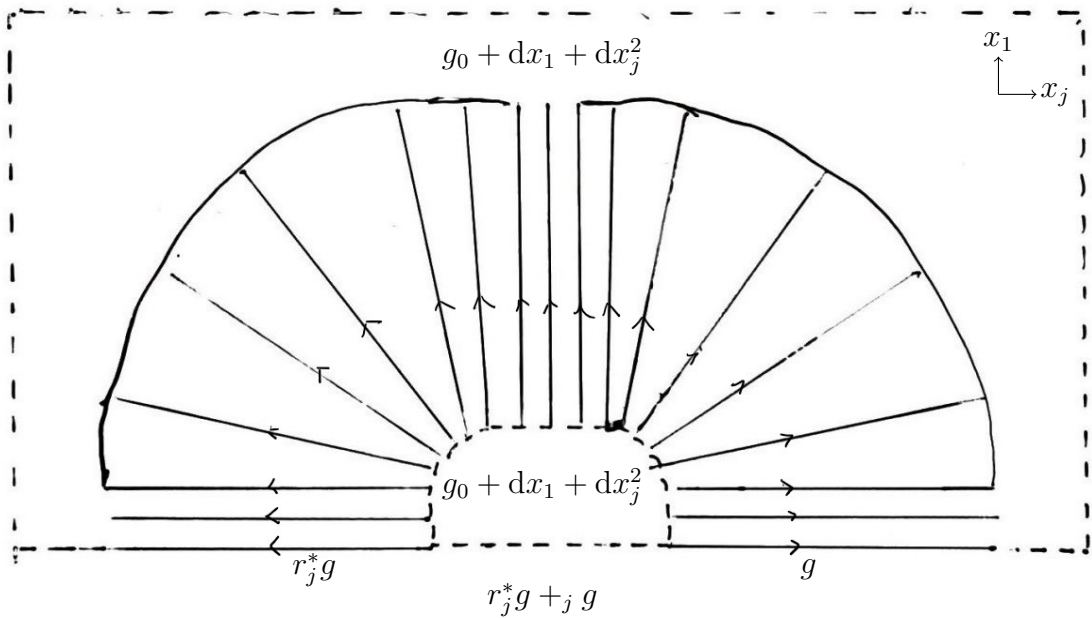
\begin{figure}[htpb]
     \centering
     \begin{tikzpicture}
      \node at (0,0) {\includegraphics[width=\textwidth]{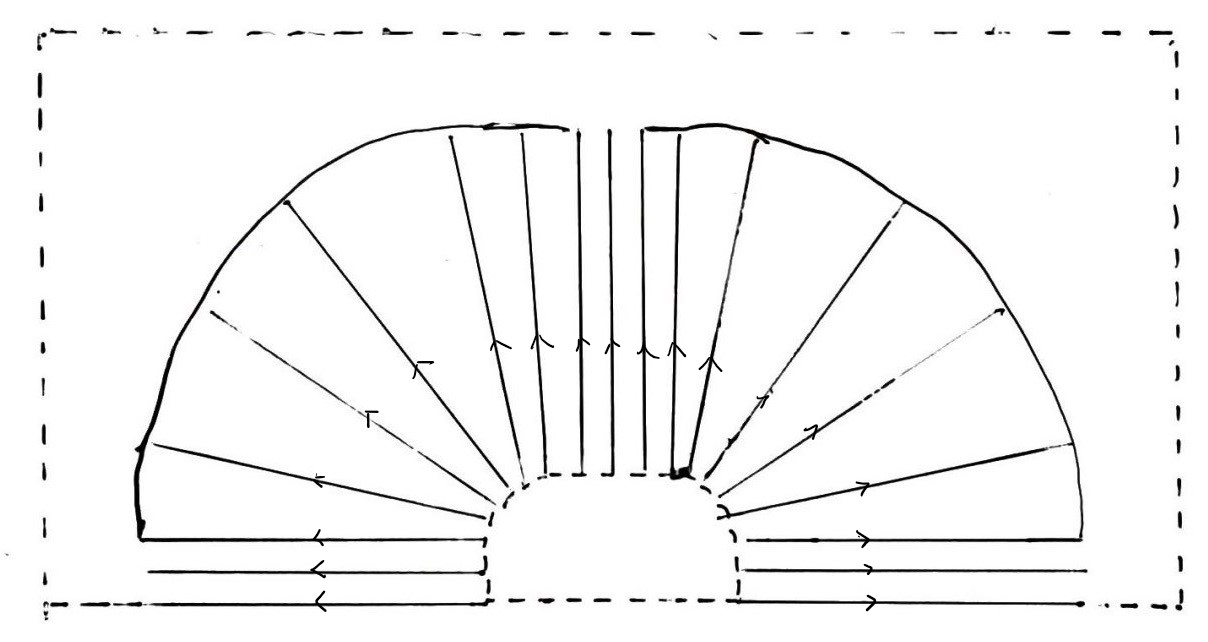}};
      \draw[->]  (6,2.5) -- (6.5,2.5);
      \draw[->]  (6,2.5) -- (6,3);
      \node at (6.75,2.5) {$x_j$};
      \node at (6,3.2) {$x_1$};
      \node at (3.5,-3.8) {$g$};
      \node at (-3,-3.8) {$r_j^\ast g$};
      \node at (0,3) {$g_0 + \diff x_1 + \diff x_j^2$};
      \node at (0, -2.7) {$g_0 + \diff x_1 + \diff x_j^2$};
      \node at (0,-4.2) {$r_j^\ast g +_j g$};
     \end{tikzpicture}
     \caption{A homotopy between $r_j^\ast g +_j g$ and the base point $g_0$. The lower quadrants are not shown. Each line segment can be identified with the lower right one and under this identification, the metric on the line segment corresponds to $g$.}
     \label{fig:ReflectionInverse}
 \end{figure}
 In formulas, if $g_R^{\angle, j}$ is the metric obtained from Proposition \ref{Construction-AngleRotation - Prop} applied to the coordiantes $x_1,x_j$, then, provided $R$ is sufficient large, the metric
 \begin{align*}
    \begin{split}
        G \deff & g_R^{\angle,j}|_{\{x_1 \geq 0, \, x_j \geq 0\}} \, \cup  \, (g + \diff x_1^2) |_{\{x_1 \leq 0, \, x_j \geq 0\}} \, \cup \\ 
        & r_j^\ast g_R^{\angle,j} |_{\{x_1 \geq 0, \, x_j \leq 0\}} \, \cup \, (r_j^\ast g + \diff x_1^2) |_{\{x_1 \leq 0, \, x_j \leq 0\}}
    \end{split}
 \end{align*}
 is a homotopy in $\ConcSet[n+1]$ between $g_0$ and $r_j^\ast g +_{j,R} g$ for all $R \geq \rho$ sufficiently large.
 
 The proof that two of these group structures are Eckmann-Hilton related is the same as for homotopy groups.
 
 The proof that $+_1$ agrees with $+_{top}$, the group structure coming from cubical set theory, is given by a picture.
 Recall that we need to find an element $H \in \ConcSet[n+1]$ that satisfies $[\face[-1]{1}H] = a$, $[\face[1]{1}H] = a +_1 b$, and $[\face[1]{2}H] = b$, while $[\face{i}H] = [g_0]$ for all other faces.
 Such a filler is given indicated in Figure \ref{fig:GeoAddVsTopAdd}.
 \begin{figure}[htpb]
     \centering
     \includegraphics[width=0.8\textwidth]{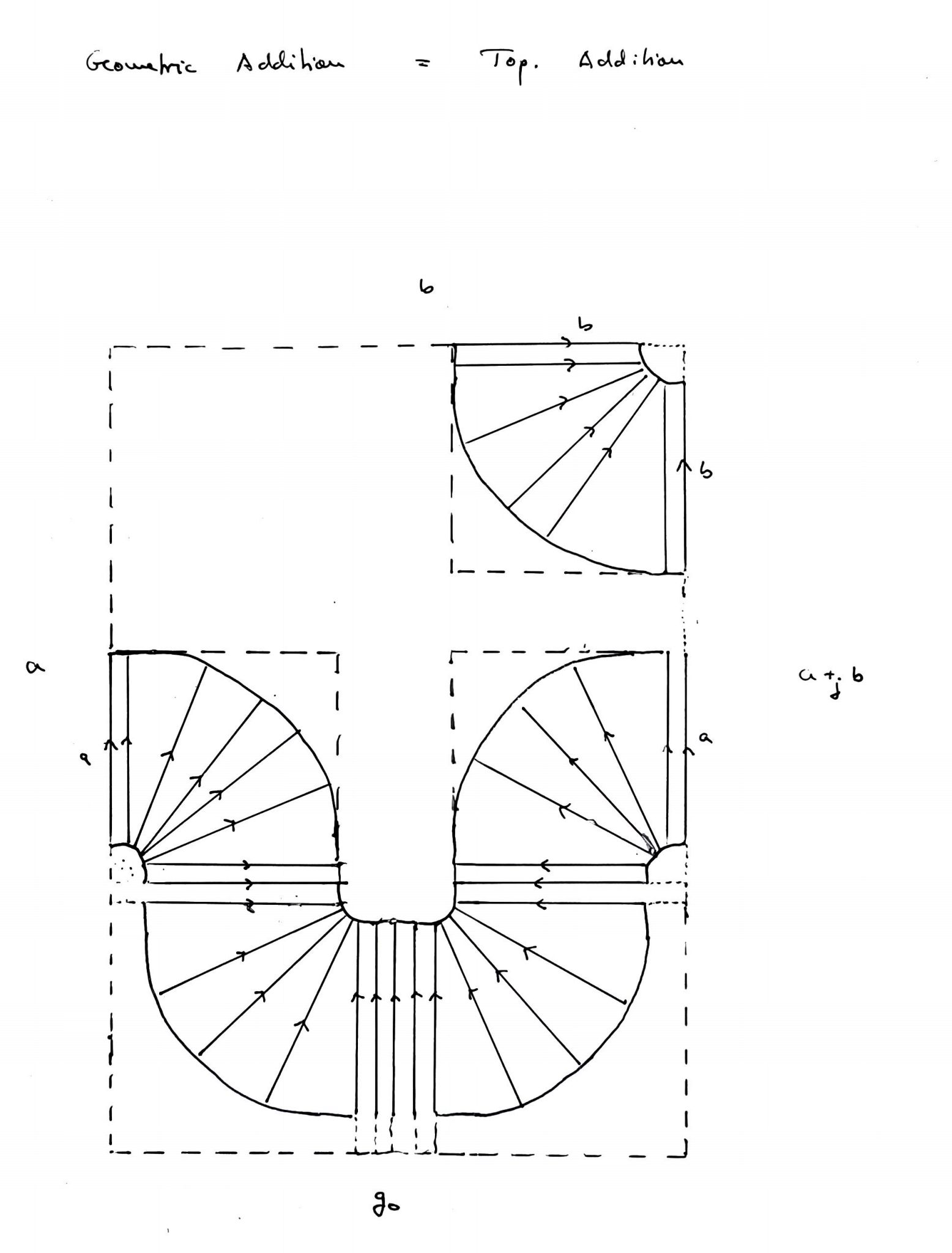}
     \caption{A homotopy between $+_1$ and $+_{top}$. In the blank space, the metric $H$ agrees with $g_0$.}
     \label{fig:GeoAddVsTopAdd}
 \end{figure}
\end{proof}



 \chapter{The Operator Concordance Set}\label{The Operator Concordance Set - Chapter}

Recall that one of our main goals is to factor the index difference through $\ConcSet$.
The naive approach would be to translate the index difference $\Riem^+(M) \rightarrow KO^{-d}$ into the category of cubical sets via the singular set functor (or a weakly equivalent version, like the cubical set of (smooth) block maps) and then try to factor the corresponding cubical map $\SingMet \rightarrow KO_\bullet^{-d}$ through $\ConcSet$.
However, the resulting map $\ConcSet \rightarrow KO^{-d}_\bullet$ would be highly non-canonical and difficult to describe because, roughly speaking, $KO_\bullet^{-d}$ is modelled using the notion of \emph{isotopy} while $\ConcSet$ is build on the notion of \emph{concordances}.

Instead, we construct a new model for real $K$-theory in this section.
It should be the geometric realisation of a cubical set that is the target of a 'nice' cubical map from $\ConcSet$ so that the proof of the factorisation result carried out in Chapter \ref{Factorisation Indexdiff - Chapter} is of mild difficulty.

The starting point is the observation made in \cite{ebert2017indexdiff} that the space $\InvPseudDir[ ]$ of invertible \emph{pseudo Dirac operators} on a closed, $d$-dimensional, spin manifold $M$, provided it is non-empty, is a classifying space for real $K$-theory.
In this model, the index difference is the map $\Riem^+(M) \rightarrow \InvPseudDir[ ]$ that assigns to a psc metric $g$ its Dirac operator $\Dirac_g$.
Our new model for real $K$-theory emerges in the same manner as $\ConcSet$ emerges out of $\Riem^+(M)$.
The cubical set $\InvBlockDirac$ of all invertible \emph{block Dirac operators} comes with a comparison map $\susp_\bullet \colon \InvPseudDir[\bullet] \dashrightarrow \InvBlockDirac$ defined on a weakly equivalent subset of $\InvPseudDir[\bullet]$.
The main result of this chapter is that this operator suspension map $\susp_\bullet$ is a weak homotopy equivalence.

In Section \ref{Foundations on Spin Geometry - Section}, we recall the foundations of spin geometry like the construction of the Clifford-linear spinor bundle and the Clifford linear Dirac operator and study how these geometric objects depend on the underlying metric. 
We also translate the results of \cite{ebert2017indexdiff} into our language and establish the space of invertible block Dirac operators as a classifying space for real $K$-theory and give a description of the (Hitchin) index difference.

Section \ref{Blockoperators - Section} presents the foundations of block operators on spinor bundles. In particular, we lay the focus on how local these pseudo differential operators are and how much they increase the support of sections.

In Section \ref{Analytical Properties of Block Operators - Section} we study the analytic properties of these operators and discuss several cut-and-paste arguments.

Section \ref{Section - Foundations of the Operator Concordance Set} is the operator theoretic analog of Section \ref{Section - Foundations on the Concordance Set}.
There we construct the operator concordance set $\InvBlockDirac$ and define the operator suspension map $\susp_\bullet$.
We will construct, in Section \ref{Section - The Comparison Map}, a weakly homotopy equivalent subset $B_\bullet \hookrightarrow \InvPseudDir[\bullet]$ so that the restriction $\susp_\bullet|_{B_\bullet}$ takes values in $\InvBlockDirac$.

In Section \ref{Section - Kan Property of the Operator Concordance Set} and \ref{SymbolMapKanFib - Section}, we show that $\InvBlockDirac$ is a Kan set followed by the even stronger statement that the cubical map $\InvBlockDirac \rightarrow \BlockMetrics$ that assigns to a block Dirac operator its underlying Riemannian metric is a cubical Kan fibration.

Section \ref{Operator Theoretic Addition - Section} is the operator theoretic analog of Section \ref{Geometric Addition - Section}.
It introduces an operator theoretic addition which agrees with the group structure on $\pi_n(\InvBlockDirac)$ coming from cubical set theory.
As an application, we show that the Gromov-Lawson index difference is a group homomorphism $\pi_n(\InvBlockDirac) \rightarrow KO^{-(d+n+1)}(\mathrm{pt})$.

Finally, in Section \ref{Section - Comparison to K-theory}, we use all previous acquired results to show that the operator suspension $\susp_\bullet \colon \InvPseudDir[\bullet] \rightarrow \InvBlockDirac$ is a weak equivalence.
We show that $\pi_n(\InvPseudDir[\bullet])$ and $\pi_n(\InvBlockDirac)$ can be identified with certain $KK$-theory groups and that, under this identification, the homomorphism $\pi_n(\susp_\bullet)$ corresponds to the Kasparov-product with the fundamental class $\alpha_{n+1}$, which is known to be an isomorphism.

\section{Foundations on Spin Geometry}\label{Foundations on Spin Geometry - Section}

We recall and summarise the necessary material on which we will base our results.
First, we start with the construction of the Clifford-right linear spinor bundle and examine the role of the underlying Riemannian metrics.
In particular, our focus lies on the identification of the two spinor bundles with different underlying metrics.
This allows us to define the universal spinor bundle over $M \times \Riem(M)$ that considers all Riemannian metric at once.

The main task of this section is the construction of the space of invertible pseudo Dirac operators $\InvPseudDir$ of a closed spin manifold $M$ of dimension $d>0$. 
Results of Johannes Ebert \cite{ebert2017indexdiff} imply that if $\InvPseudDir$ is not empty, then it is a classifying space for $KO^{- (\dim M + 1)}$.

\subsection*{The Spin Package}

We recall the foundations of classical spin geometry and set up notational conventions. 
Everything that is covered here can be found in the book \cite{LawsonMichelsonSpin} except for Lemma \ref{OrientationChangeSpinor - Lemma} and \ref{PinorFunctoriality - Lemma}.
Although these results are probably known to experts, the author has not found them in the literature. 

\begin{definition}\label{Clifford Algebra - Def}
 Let $(V,q)$ be a vector space with a quadratic form. 
 The \emph{Clifford algebra} $Cl(V,q)$
 is the unital algebra that contains $V$ and satisfies the following universal property: 
 Every injective linear map $f \colon V \hookrightarrow \mathcal{A}$ into a unital algebra $\mathcal{A}$ that satisfies $f(v)^2 = - q(v) \cdot 1$ 
 extends uniquely to a unital algebra homomorphism $Cl(f) \colon Cl(V,q) \rightarrow \mathcal{A}$.  
\end{definition}
\begin{definition}
 Every Clifford algebra has a canonical grading, which is given by the parity operator $Cl(-\id)$. 
 We call this grading the \emph{even/odd}-grading and it is denoted by $\evenodd$.
\end{definition}

If $V$ is a real vector space, we can consider it as Real vector space $(V \otimes \C, \id \otimes \overline{\cdot})$. 
Then the Clifford algebra is Real as well. 
The Real structure $\overline{\cdot}$ commutes with the parity operator.

\begin{rem}
 Every Real vector space or Real vector bundle that we consider in this thesis comes from a real vector space or bundle via the above construction. 
 We will not notationally distinguish them. 
 The reason why we actually need to work with Real vector bundles is that we consider pseudo differential operators on these bundles, which are locally defined with the help of Fourier transformation, see Appendix \ref{Appendix Pseudos - Chapter} for details, which is only defined on complex vector spaces and that the principal symbol of pseudo differential operator does not in general preserve the real bundle.
\end{rem}

\begin{example}
 We denote the vector space $V = \R^p \times \R^q$ endowed with the quadratic form $q(x,y) \deff ||x||^2_{\mathrm{eucl}} - ||y||^2_{\mathrm{eucl}}$ with $\R^{p,q}$.
 The standard basis is denoted by $e_1\dots,e_p,\eps_1,\dots,\eps_q$.
 Its Clifford algebra is multiplicatively generated by the standard basis of $\R^{p,q}$. The generators satisfy the following relations
 \begin{align*}
  e_ie_j + e_je_i = -2 \delta_{ij}, \\
  \eps_i\eps_j + \eps_j \eps_i = 2\delta_{ij},
 \end{align*}
 and $Cl(\R^{p,q})$ is the universal algebra with these generators and relations. 
 
 Let $I$ and $I'$ be strictly monotonically increasing ordered tupels with entries between $1$ and $p$, and $1$ and $q$, respectively.
 If $e_I$ denotes $e_{i_1}\cdot \dots \cdot e_{i_{|I|}}$ (with the convention $e_{\emptyset} = 1$), then the set of all $e_I \cdot \eps_{I'}$ form a basis of $Cl(\R^{p,q})$. 
 There is a unique quadratic form $||\cdot ||^2$ that turns that basis into a (signed) orthonormal basis. 
 The Clifford algebra endowed with all this data, is denoted by $\Cliff[p][q]$.
 \nomenclature{$\Cliff[p][q]$}{Clifford algebra of $\R^{p,q}$}
 
 The Clifford algebra carries a structure of a $C^\ast$-algebra. 
 Indeed, the left multiplication map $l \colon \Cliff[p][q] \rightarrow \mathrm{End}(\Cliff[p][q])$ is an injective algebra homomorphism whose image is closed under taking adjoints.
 It suffices to show this for the generators $e_1,\dots, e_p, \eps_1,\dots, \eps_q$.
 The linear maps $l(e_i)$ and $l(\eps_j)$ are isometries as they map the canonical orthonormal basis of $\Cliff[p][q]$ to itself up to sign.
 Together with the Clifford relations this implies
 \begin{equation*}
     l(e_i)^\ast = l(-e_i)\qquad \text{and} \qquad l(\eps_j)^\ast = l(\eps_j).
 \end{equation*}
 
 In the following we only write $e_i\cdot$ and $\eps_j\cdot$ for $l(e_i)$ and $l(\eps_j)$ unless we want to emphasise the representation $l$.
 \nomenclature{$l$}{Representation of $\Cliff[p][q]$ on itself by left-multiplication}
 Usually we will denote the $C^\ast$-norm on $\Cliff[p][q]$ with $||\placeholder||_{\mathrm{op}}$. 
 Note that $||v||_{\mathrm{op}} = ||v||$ for vectors $v \in \R^{p,0}$ but $||v||_{\mathrm{op}} \neq ||v||$ for general elements $v \in \Cliff[p][0]$.
 \todo[inline]{Mündliche Prüfung: Find out the relation between these two norms.}
\end{example}
The most important case for us is $q=0$. We will restrict our discussions of spin geometry to this case, although it holds in larger generality.
\begin{definition}
 We denote by $\Spin(d)$
 \nomenclature{$\Spin(d)$}{Spin group}
 the subgroup of $\Cliff[d][0]$ that is multiplicative generated by $\{v\cdot w \, : \, v,w \in \R^d \subseteq \Cliff[d][0], \, ||v||= ||w||=1\}$ and we denote by $\mathrm{Pin}^-(d)$ 
 \nomenclature{$\mathrm{Pin}^-$}{Pin group}
 the subgroup of $\Cliff[d][0]$ that is multiplicatively generated by $\{v \, : \, v \in \R^d \subseteq \Cliff[d][0], \, ||v||=1\}$.
\end{definition}

 Note that left multiplication with elements in $\mathrm{Pin}^-(d)$ is an isometry on $\Cliff[d][0]$.
 The conjugation with elements of $\mathrm{Pin}^-(d)$ is another isometric action on $\Cliff[d][0]$.
 It preserves the subspace $\R^d \subset \Cliff[d][0]$ and thus yields a homomorphism $\lambda \colon \mathrm{Pin}^-(d) \rightarrow \Or(d)$.
 One can show that it is a two-sheeted covering, see for example \cite{LawsonMichelsonSpin}*{Theorem I.2.10}.
 Restricted to $\Spin(d)$, this yields a two-sheeted cover $\lambda \colon \Spin(d) \rightarrow \SOr(d)$.

To generalise these concepts to manifolds, the best approach seems to be via principal frame bundles and representations.
A basic reference for principal bundles, their associated bundles and their differential geometry is \cite{baum2009eichfeldtheorie}.
In the following, let $M^d$ be a smooth manifold of dimension $d$ (it may have corners and does not need to be compact).

\begin{definition}
  For each manifold $M^d$, the \emph{bundle of frames} is given by
  \begin{equation*}
      \GL(M) \deff  \{\R^d \xrightarrow{\cong} TM_x \, : \, \text{isomorphisms}\} \rightarrow M.
  \end{equation*}
  It is a principal $\GL(d)$-bundle whose $\GL(d)$-right action is given by pre-composition.
  If $M^d$ is additionally orientable, we can define, for each orientation, the subbundle
 \begin{equation*}
     \GL^+(M) \deff \{\R^d \xrightarrow{\cong} T_xM \, : \, \text{orientation preserving isomorphisms} \} \rightarrow M.
 \end{equation*}    
 It is a principal $\GL^+(d)$-bundle.
\end{definition}
Conversely, a choice of a $\GL^+(d)$-reduction of $\GL(M)$ is a choice of an orientation.

\begin{definition}
 For a Riemannian metric $g$, we define \emph{bundle of orthonormal frames} to be 
 \begin{equation*}
     \Or(M,g) \deff \{(\R^d,\euclmetric) \xrightarrow{\cong} (T_xM,g_x) \, : \, \text{isometries}\} \rightarrow M.
 \end{equation*}
 If $M^d$ is additionally orientable, we can define, for each orientation, the subbundle
 \begin{equation*}
     \SOr(M,g) \deff \Or(M,g) \cap \GL^+(M) \nomenclature{$\SOr(M,g)$}{Principal bundle of oriented, orthonormal frames}.
 \end{equation*}
\end{definition}
A choice of a Riemannian metric on $M^d$ is equivalent to the choice of an $\Or(d)$-reduction of $\GL(M)$.
Conversely, a choice of an $\Or(d)$-reduction of $\GL(M)$ gives a Riemannian metric because an inner product is uniquely determined by one of its orthonormal bases. 

A \emph{spin structure} on $M$ is the choice of a $\Spin(d)$-reduction of $\SOr(M,g)$, that is a choice of a principal $\Spin(d)$ bundle $P_{\Spin,g}$ and a $\lambda$-equivariant map $\xi \colon P_{\Spin,g} \rightarrow \SOr(M,g)$ over $M$; in other words, the following diagram, in which the horizontal arrows denote the right action, commutes:
\begin{equation*}
    \xymatrix{P_{\Spin,g} \times \Spin(d)  \ar[r] \ar[d]_{\xi \times}^{\lambda} & P_{\Spin,g}  \ar[d]^\xi \\
    \SOr(M,g) \times \SOr(d) \ar[r] & \SOr(M,g).}
\end{equation*}

Not every orientable manifold carries a spin structure. 
The existence of a spin structure is obstructed by the first and second Stiefel-Whitney classes of $M$ \cite{LawsonMichelsonSpin}*{Theorem II.1.7}.
Note that the choice of a spin structure involves the choice of a metric. 
Since manifolds can carry several inequivalent spin structures\footnote{Indeed, the set of (equivalence classes of) spin-structure is parameterised by $H^1(M,\Z_2)$, see \cite{LawsonMichelsonSpin}*{Theorem II.1.7}}, we have to make sure to pick spin structures consistently when varying the metric. 
This can be done as follows:
Let $\widetilde{\GL(d)}^+ \rightarrow \GL(d)^+$ the unique two-sheeted cover that extends the canonical group homomorphism $\Spin(d) \rightarrow \SOr(d)$. The vanishing of the first and second Stiefel-Whitney classes implies the existence of a $\widetilde{\GL(d)}^+$-reduction $\rho \colon P_{\widetilde{\GL(d)}^+} \rightarrow \GL^+(M)$. 
We fix such a structure. 
For each metric and each choice of $\widetilde{\GL(d)^+}$-structure, the $\Spin(d)$-reduction of our choice is given by the following pull-back\footnote{It can be shown that any spin structure can be obtained in that manner}
\begin{equation*}
 \xymatrix{P_{\Spin,g} \ar@{^{(}->}[r] \ar[d] & P_{\widetilde{\GL^+}} \ar[d]^\rho \\ \SOr(M,g) \ar@{^{(}->}[r] & \GL^+(M).} 
\end{equation*} 

With the help of principal bundles we can generalise the Clifford algebra construction to vector bundles over manifolds. 
Note that this construction does not require $M$ to be spin.
\begin{definition}
 Let $M^d$ be a manifold with Riemannian metric $g$.
 The \emph{Clifford bundle} is given by 
\begin{equation*}
 Cl(TM,g) \deff \SOr(M,g) \times_{Cl(\mathrm{taut})} \Cliff[d][0] = P_{\Spin,g} \times_{\mathrm{conj}} \Cliff[d][0],
\end{equation*}
\nomenclature{$Cl(TM,g)$}{Clifford bundle of $TM$, sometimes only $Cl(M,g)$}
where $\mathrm{taut}$ is the tautological representation of $\SOr(d)$ on $\R^d$ and $\mathrm{conj}$ is the representation that is given by conjugations of elements of $\Spin(n)$ in $\Cliff[d][0]$.
\end{definition}
 The two representations, $Cl(\mathrm{taut})$ and $l$, are isometric. 
 The bundle $Cl(TM,g)$ is a bundle of algebras with fibres $Cl(T_xM,g_x)$. 
 It canonically includes $TM$. 
 The Riemannian metric $g$ on $TM$ canonically extends to a metric on $Cl(TM,g)$ so that all products $e_I$ form an orthonormal basis. 
 Since $Cl(\mathrm{taut})$ preserves the grading, $Cl(TM,g)$ is a graded algebra bundle. 
 The Levi-Cevita connection $\nabla^{L.C.}$ on $(TM,g)$ extends uniquely to a metric connection $\nabla^{L.C.}$ on $Cl(TM,g)$ that is a derivation with respect to the algebra structure \cite{LawsonMichelsonSpin}*{Proposition II.4.8}, meaning that, for all smooth sections $\sigma, \tau$, we have
\begin{equation*}
 \nabla^{L.C.}(\sigma \tau) = (\nabla^{L.C.}\sigma)\tau + \sigma \nabla^{L.C.}(\tau).
\end{equation*} 
The grading and the Real structure are parallel with respect to that connection.

If we use the left-multiplication as representation, we obtain the spinor bundle.
In contrast to the conjugation, this representation does not factor through $\SOr(d)$, which is the reason why this construction only works over spin manifolds.
\begin{definition}
 Let $M^d$ be a spin manifold with Riemannian metric $g$.
 The associated \emph{spinor bundle} is given by 
\begin{equation*}
 \Spinor_g \deff P_{\Spin, g} \times_{\Spin(d),l} \Cliff[d][0].
\end{equation*} 
\nomenclature{$\Spinor_g$}{Spinor bundle of the metric $g$}
It is a $Cl(TM,g)$-left module bundle and the module structure is denoted with $\Cliffmult$.
\nomenclature{$\Cliffmult$}{$Cl(TM,g)$-multiplication or left module action on $\Spinor_g$}
\end{definition}
The spinor bundle also has a canonical $\Cliff[d][0]$-right action $\mathbf{r}$ 
\nomenclature{$\mathbf{r}$}{Right-action of $\Cliff[d][0]$ on $\Spinor_g$}
induced by right multiplication. 
It is well defined because left and right multiplication commute.
The left-multiplication $l$ acts by isometries, so $\Spinor_g$ carries a canonical Riemannian metric induced by the inner product on $\Cliff[d][0]$.
Multiplication with unit vectors induces isometries with respect to this metric.
Since $\Spin(d)$ consists only of even elements in $\Cliff[d][0]$, the left multiplication preserves the grading of $\Cliff[d][0]$, so the grading on $\Cliff[d][0]$ carries over to $\Spinor_g$\footnote{This is the very reason why we restrict to spin structures. 
The construction still makes sense, if we consider only Pin structures, but we would then loose the grading, on which index theory relies.}. 
We call this grading again $\evenodd$. 
\nomenclature{$\evenodd$}{Even/Odd-grading on spinor bundle}
Right multiplication with fixed odd elements are odd linear maps.
The Levi-Cevita connection extends canonically to a metric connection $\nabla^{\Spinor}$ on $\Spinor_g$. 
\nomenclature{$\nabla^{\Spinor_g}$}{Spinor connection}
The resulting metric connection is uniquely characterised by the following properties:
\begin{itemize}
 \item[(i)] $\nabla^{\Spinor}$ is a $Cl(TM,g)$ module derivation,
 \item[(ii)] $\nabla^{\Spinor}$ is even,
 \item[(iii)] $\nabla^{\Spinor}$ is linear with respect to the $\Cliff[d][0]$-right action.
\end{itemize} 
Thus, $\Spinor$ is a Dirac-bundle in the sense of \cite{LawsonMichelsonSpin}*{Def. II.5.2}.

Dirac bundles are precisely those bundles that allow the definition of a Dirac operator.
We give definition in the case of the spinor bundle.
\begin{definition}\label{Def - Dirac operator}
 The \emph{Dirac operator} 
 \nomenclature{$\Dirac_g$}{Dirac operator of the Riemannian metric $g$}
 of a Riemannian metric $g$ on $M$ is defined by the composition
 \begin{equation*}\label{Eq: Def Dirac Operator}
  \xymatrix{\Dirac_g \colon \Gamma_c(\Spinor_g) \ar[r]^-{\nabla^{\mathrm{\Spinor}}} & \Gamma_c(T^\vee M \otimes \Spinor_g) \ar[r]^-{g_\sharp \otimes \id} & \Gamma_c(TM \otimes \Spinor_g) \ar[r]^-{\Cliffmult(\cdot)} & \Gamma_c(\Spinor_g), } 
 \end{equation*}
 where $g_\sharp \colon T^\vee M \rightarrow TM$ is defined by $g(g_\sharp(\xi),\placeholder) = \xi$.
\end{definition}
The inner product on $\Spinor_g$, denoted by $h$, and the volume form $\vol_g$
\nomenclature{$\vol_g$}{Volume form on $M$ induced by $g$}
on $M$ induced by the Riemannian metric $g \in \Riem(M)$ allow us to equip $\Gamma_c(\Spinor_g)$ with an inner product, see Appendix \ref{Sobolev Spaces - Chapter} for the general construction.
It is given by
\begin{equation*}
    (\sigma_1,\sigma_2) \deff \int_M h_x(\sigma_1(x),\sigma_2(x)) \diff \vol_g(x).
\end{equation*}
The completion of $\Gamma_c(\Spinor_g)$ is denoted by $L^2(\Spinor_g)$ and is called the space of \emph{square integrable spinors} $\Spinor_g$.
The $Cl_{d,0}$-right action $\mathbf{r}$ and the $\Z_2$-grading $\evenodd$ on $\Spinor_g$ carry over to a $\Cliff[d][0]$-right action and a $\Z_2$-grading on $L^2(\Spinor_g)$, which we denote with the same symbols.
The Dirac operator acts as an unbounded operator on this space.
More precisely, we have the following theorem, which summarises results from \cite{LawsonMichelsonSpin}*{Chapter II}.
\begin{theorem}
 The Dirac operator $\Dirac_g$ is an elliptic, odd, symmetric,  
 first-order, $Cl_{d,0}$-right linear differential operator with principal symbol $i\Cliffmult(g_\sharp(\cdot))$.
 If $(M^d,g)$ is complete, then $\Dirac_g$ is self-adjoint. 
 If $M$ is closed, then $\Dirac_g$ is an (unbounded) Fredholm operator, that is, it has a finite dimensional kernel. 
\end{theorem}

The relation between the Dirac operator and the scalar curvature of a Riemannian metric is established by the Lichnerowicz-formula, which is of outstanding importance.

\begin{theorem}[Lichnerowicz-Formula]
  For every Riemannian metric $g$, the formula  
  \begin{equation}
      \Dirac_g^2 \sigma = \left({\nabla^{\Spinor_{g}}}^\ast \circ \nabla^{\Spinor_g} + \frac{\scal(g)}{4}\right) \sigma, 
  \end{equation}
 holds true for all $\sigma \in \Gamma_c(\Spinor_g)$.
 Here ${\nabla^{\Spinor_{g}}}^\ast \colon \Gamma_c(T^\vee M\otimes \Spinor_g) \rightarrow \Gamma_c(\Spinor_g)$ is the formal adjoint of the spinor connection.
\end{theorem}
Since $\nabla^\ast \nabla$ is positive semi-definite, we immediately deduce that a psc-metric contradicts the existence of smooth harmonic spinors that can be approximated by compactly supported spinors. 
\begin{cor}
 Let $g \in \Riem^+(M)$, then the kernel of $\Dirac_g$ in the minimal domain is trivial.
\end{cor}

The Lichnerowicz-formula combined with the Atiyah-Singer index theorem gives a powerful topological obstruction to positive scalar curvature metrics on closed spin manifolds. 
If the underlying manifold is closed, then the kernel of the Dirac operator is a finite dimensional graded $\Cliff[d][0]$-module and hence represents an element in real $K$-theory $KO^{-d}(\mathrm{pt})$, see \cite{LawsonMichelsonSpin}*{Theorem I.9.21}, which we call the \emph{analytical index}.
The Clifford-linear Atiyah-Singer index theorem, proved in \cite{LawsonMichelsonSpin}*{Chapter III.16}, relates this $KO$-theory class with the $\alpha$-invariant, or \emph{topological index}, which is homomorphism $\alpha \colon \Omega_\ast^{\Spin}(\mathrm{pt}) \rightarrow KO^{-\ast}(\mathrm{pt})$, see \cite{LawsonMichelsonSpin}*{Theorem II.7.14}.
\begin{theorem}
 For a closed spin manifold $(M^d,g)$, the analytic index and the topological index agree:
 \begin{equation*}
     [\ker \Dirac_g] = \alpha(M) \in KO^{-d}(\mathrm{pt}).
 \end{equation*}
\end{theorem}

The Atiyah-Singer index formula implies that if we change the orientation on $M^d$, then the analytical index of $\Dirac_g$ needs to change the sign.
The effect of the chosen orientation is less transparent in the Clifford-linear setup in contrast to the setup that uses the irreducible spinor bundle, where a change of orientation results in the change of the grading.

\begin{lemma}\label{OrientationChangeSpinor - Lemma}
 Let $M$ be a connected spin manifold of dimension $d$. 
 If $M^{\mathrm{op}}$ is the same manifold equipped with the opposite orientation, then there is an \emph{odd} isomorphism of $Cl(M,g)$-left module bundles
 \begin{equation*}
     \mathrm{Sw} \colon \Spinor(M) \xrightarrow{\cong} \Spinor(M^{\mathrm{op}}).
 \end{equation*}
 The map $\mathrm{Sw}$ is not necessarily right-linear but $\mathrm{Sw}^\ast(\mathbf{r})$ and $\mathbf{r}$ are homotopic through $\Cliff[d][0]$-right actions on $\Spinor(M)$. 
\end{lemma}
\begin{proof}
If a Riemannian manifold $(M^d,g)$ is orientable, then the bundle of orthonormal frames has an $SO(d)$-reduction 
\begin{equation*}
    \Or(M,g) \cong \SOr(M,g) \times_{\SOr(d)} \Or(d),
\end{equation*}
and a choice of an orientation agrees with a choice of a connected component.
Every linear map $A \in \Or(d) \setminus \SOr(d)$ allows us to change the orientation by right multiplication\footnote{Warning! This map is not $\Or(d)$-equivariant.}
\begin{equation*}
    \xymatrix{\Or(M,g) \ar@{=}[r] & \SOr(M,g) \times_{\SOr(d)} \Or(d)  \ar[r]^{\id \times \cdot A} & \SOr(M,g) \times_{\SOr(d)} \Or(d) \ar@{=}[r] & \Or(M,g) }
\end{equation*}
Since $\mathrm{Pin}^-(d) \rightarrow \Or(d)$ is a two-sheeted cover, there is an element $\phi \in \mathrm{Pin}^-(d)\setminus \Spin(d)$, unique up to sign, such that $\lambda(\phi) = A$.

If $\xi \colon P_{\mathrm{Pin}^-,g} \rightarrow \Or(M,g)$ is the $\mathrm{Pin}^-$ structure that restricts to the chosen spin structure, then the chosen element $\phi$ induces the change of orientation
\begin{equation*}
    \xymatrix{P_{\mathrm{Pin}^-,g} \ar@{=}[r] \ar[dr]_\xi & P_{\Spin,g}\times_{\Spin(d)} \mathrm{Pin}^-(d) \ar[d] \ar[rr]^{\id \times \cdot \phi} && P_{\Spin,g}\times_{\Spin(d)} \mathrm{Pin}^-(d) \ar[d] & \ar@{=}[l] \ar[dl]^\xi P_{\mathrm{Pin}^-,g} \\
    & \Or(M,g) \ar[rr]_{\id \times \cdot A} && \Or(M,g). & }
\end{equation*}
The identity
\begin{equation*}
    \Spinor(M) = P_{\Spin,g} \times_{\Spin(d),l} \Cliff[d][0] \cong P_{\Spin(d)} \times_{\Spin(d)} \mathrm{Pin}^-(d) \times_{\mathrm{Pin}^-(d),l} \Cliff[d][0]
\end{equation*}
shows that the choice of the embedding of $P_{\Spin(d),g}$ into $P_{\mathrm{Pin}^-(d),g}=P_{\Spin(d),g} \times_{\Spin(d)} \mathrm{Pin}^-(d)$ does not affect the isomorphism type as we can extend the right multiplication with $\phi$ from $\mathrm{Pin}^-(d)$ to $\Cliff[d][0]$. 
Thus, we deduce that $\Spinor(M) \cong \Spinor(M^{\mathrm{op}})$ and that the isomorphism $\mathrm{Sw}$ between them is given by right multiplication with $\phi$, which is odd as $\phi \notin \Spin(d)$.

If $d$ is odd, then $-\id \notin \SOr(d)$ and we can choose $\phi = e_1 \cdot \hdots \cdot e_d$, which is in the centre of the Clifford algebra.
In this case, $\mathrm{Sw}$ is $\Cliff[d][0]$-right linear.

If $d$ is even, then we need to argue differently.
Since $\mathrm{Sw}^\ast(\mathbf{r}) = \mathbf{r}(\phi^{-1}(\placeholder)\phi)$, it suffices to show that the conjugation with $\phi^{-1}$ is homotopic through algebra endomorphisms to $\id_{\Cliff[d][0]}$.
By a theorem of Cartan-Dieudonne, see \cite{LawsonMichelsonSpin}*{Theorem I.2.7}, every orthogonal matrix can be decomposed into a sequence of reflections along hyperplanes, so it suffices to assume that $\phi = v \in \R^d \subseteq \Cliff[d][0]$ is a single vector with length $1$.
By the universal property of $\Cliff[d][0]$, algebra homomorphisms correspond to injective linear maps $A \colon \R^d \rightarrow \Cliff[d][0]$ such that $A(x)^2 = - ||x||^2$. 
It thus suffices to find a homotopy between the reflection at the hyperplane perpendicular to $v$ and the identity through such maps.
If $w$ is another vector with length $1$ that is perpendicular to $v$, we can define such a homotopy $A_s$ as the linear extension of the following map
\begin{equation*}
    A_s(v) = \cos(s)v + \sin(s) v \cdot w,  \qquad A_s(x) = x, \text{ if } x \in (\R v)^\perp. 
\end{equation*}
Anti-commutativity guarantees that $A_s(x)^2 = -||x||^2$ for all $x \in \R^n$. 
Clearly, $A_0 = \id$ and $A_{\pi}$ is the reflection at the hyperplane perpendicular to $v$, which finishes the proof.  
\end{proof}

The observation from the previous proof offers the opportunity to lift isometries, which are not necessarily orientation preserving, to homomorphisms between spinor bundles.
However, we cannot expect that every isometry $\varphi \colon (M,g_1) \rightarrow (M,{g_2})$ lifts to an isometry of the associated spinor bundles.
The following criterion is sufficient.
\begin{definition}
 Let $\xi_j \colon P_{\mathrm{Pin}^-,g_j} \rightarrow \Or(M,g_j)$ be $\mathrm{Pin}^-(d)$-reductions of $\Or(M,g_j)$. 
 An isometry $F\colon (TM,g_1) \rightarrow (TM,g_2)$ is \emph{Pin-structure preserving} if the induced map between the principal bundles of frames lifts to an equivariant map between the Pin-structures.  
 The isometry $F$ is \emph{spin structure preserving}, if it is orientation preserving and Pin-structure preserving. 
\end{definition}
\begin{lemma}\label{PinorFunctoriality - Lemma}
 A Pin-structure preserving isometry $F \colon (TM,g_1) \rightarrow (TM,g_2)$ induces an isometry $\Spinor(F) \colon \Spinor_{g_1} \rightarrow \Spinor_{g_2}$ that is $\Cliff[d][0]$-right linear and $Cl(M,g_1)$-$Cl(M,g_2)$ equivariant.
 If $F$ additionally preserves the orientation, then $\Spinor(F)$ is even.
\end{lemma}
\begin{proof}
 A lift $\tilde{F} \colon P_{\mathrm{Pin^-},g_1} \rightarrow P_{\mathrm{Pin^-},g_2}$ of $F$ to the Pin structures induces a map on spinor bundles via
 \begin{equation*}
     \xymatrix{ \Spinor(F) \colon \Spinor_{g_1} \ar@{=}[r] & P_{\mathrm{Pin^-,g_1}} \times_{\mathrm{Pin}^-(d)} {\Cliff[d][0]} \ar[rr]^{\tilde{F}\times \id} && P_{\mathrm{Pin^-,g_2}} \times_{\mathrm{Pin}^-(d)} {\Cliff[d][0]} \ar@{=}[r] & \Spinor_{g_2}}. 
 \end{equation*}
 The desired properties follow immediately.
 
 If $F$ is orientation preserving, then $\tilde{F}$ restricts to an equivariant map on $P_{\Spin,g_1} \rightarrow P_{\Spin,g_2}$, so that $\Spinor(F)$ preserves the grading.
\end{proof}

We close this section by discussing the spinor bundle and the Dirac operator of a product manifold.

\begin{example}\label{Spinor bundle product manifold - Example}
  Recall from \cite{LawsonMichelsonSpin} that, for two quadratic spaces $(V_i,q_i)$, the Clifford algebra decomposes as
  \begin{equation*}
      \Clifford(V_1 \oplus V_2,q_1 \oplus q_2) \cong \Clifford(V_1,q_1) \hat{\otimes} \Clifford(V_2,q_2).
  \end{equation*}
  The symbol $\hat{\otimes}$ indicates that we have a graded tensor product and a graded algebra structure.
  That means that the algebra structure on the right hand side is given by the linear extension of
  \begin{equation*}
      (v_1 {\otimes} v_2) \cdot (w_1 {\otimes} w_2) = (-1)^{|v_2|\cdot |w_1|} (v_1 \cdot w_1) {\otimes} (v_2\cdot w_2),
  \end{equation*}
  where $v_1,w_1 \in \Clifford(V_1,q_1)$ and $v_2,w_2 \in \Clifford(V_2,q_2)$ are elements of pure degree, i.e., eigenvectors of the grading. 
  The symbol $|v_j| \in \{0,1\}$ denotes the degree of $v_j$; it is zero if and only if $v_j$ is even.
  A similar formula holds for the left representation $l \colon \Clifford(V_1 \oplus V_2,q_1 \oplus q_2) \rightarrow \mathrm{End}(\Clifford(V_1 \oplus V_2,q_1 \oplus q_2))$.
  
  In the case that $(M,g) = (M_1 \times M_2, g_1 \oplus g_2)$ is a Riemannian product, we have the reductions
  \begin{equation*}
     SO(M_1,g_1) \times SO(M_2,g_2) \subseteq SO(M,g) \qquad \text{and} \qquad P_{\Spin,g_1} \times_{\Z_2} P_{Spin,g_2} \subset P_{Spin,g}
  \end{equation*}
  so that the decomposition of Clifford algebras generalises to bundles
  \begin{equation*}
      \Clifford(M_1 \times M_2, g_1 \oplus g_2) \cong \Clifford(M_1,g_1) \hat{\boxtimes} \Clifford(M_2,g_2) 
  \end{equation*}
  and
  \begin{equation*}
      \Spinor(M_1 \times M_2, g_1 \oplus g_2) \cong \Spinor(M_1,g_1) \hat{\boxtimes} \Spinor(M_2,g_2),
  \end{equation*}
  where $\hat{\boxtimes}$ denotes the graded exterior tensor product of two graded vector bundles.
  Recall that the exterior tensor product 
  \nomenclature{$\boxtimes$}{Exterior Tensor Product}
  of two vector bundles $E \rightarrow X$ and $F \rightarrow Y$ is the vector bundle on $E \boxtimes F \rightarrow X \times Y$ whose fibre at $(x,y)$ is $E_x \otimes F_y$.
  If the bundles are graded with gradings $\iota_E$ and $\iota_F$, then the grading on $E \hat{\boxtimes} F$ is given by $\iota_{E \hat{\boxtimes} F} = \iota_E \otimes \iota_F$.
  
  These isomorphisms are isomorphisms of algebras and left-modules, respectively. 
  They also preserve the grading, are isometric, and, in case of spinor bundles, respect also the right action from the Clifford algebra under the identification $\Cliff[d_1+d_2] \cong \Cliff[d_1][0] \hat{\otimes} \Cliff[d_2][0]$.
  
  The $\Clifford(M_1 \times M_2)$-left action restricted to $TM_1 \oplus TM_2$ corresponds to
  \begin{equation*}
      \Cliffmult(v_1 \oplus v_2)(\sigma_1 \otimes \sigma_2) = \Cliffmult_1(v_1)(\sigma_1) \otimes \sigma_2 + \evenodd(\sigma_1) \otimes \Cliffmult_2(v_2)(\sigma_2).
  \end{equation*}
  In the product case, the musical isomorphism $(g_1 \oplus g_2)_\sharp \colon T^\vee M_1 \oplus T^\vee M_2 \rightarrow TM \oplus TM$ is given by
  \begin{equation*}
      (g_1 \oplus g_2)_\sharp(\xi_1,\xi_2) = (g_1)_\sharp(\xi_1) + (g_2)_\sharp(\xi_2).
  \end{equation*}
  By uniqueness, the spinor-connection on $\Spinor(M_1 \times M_2)$ corresponds to
  \begin{equation*}
     \nabla^{\Spinor(M_1)} \otimes \id + \id \otimes \nabla^{\Spinor(M_2)} 
  \end{equation*}
  under the decomposition. 
  In conclusion, the Dirac operator of a product metric is given by
  \begin{equation*}
      \Dirac_{g_1 \oplus g_2} = \Dirac_{g_1} \otimes \id + \evenodd \otimes \Dirac_{g_2}.
  \end{equation*}
  However, we often write $\Dirac_{g_1} \boxtimes \id + \evenodd \boxtimes \Dirac_{g_2}$ instead, to emphasise the product structure of the underlying metric.
 \end{example}

\subsubsection{The Universal Spinor Bundle}

We would like to consider all spinor bundles at once, $\mathrm{i.e.}$, we would like to construct a bundle $\Spinor \rightarrow \Riem(M) \times M$ with fibre $\Spinor_{g_m}$. In order to do that we have to address how the change of the underlying metric affects the constructions.
The ansatz is taken from \cite{ebert2017indexdiff}, but we will give a less abstract presentation and provide more details.

First, we recall the \emph{Gram endomorphism}. 
A bilinear form $b$ on a finite dimensional vector space $V$ defines a canonical map $b^\flat \colon V \rightarrow V^\vee$ via $v \mapsto b(v,\placeholder)$.
The bilinear form $b$ is non-degenerate if and only if $b^\flat$ is an isomorphism.
If $g$ is an inner product on $V$, we denote the inverse of $g^\flat$ with $g_\sharp$ and define the \emph{Gram endomorphism} of $b$ (with respect to $g$) as $\mathrm{Gram}_g(b) \deff g_\sharp \circ b^\flat$.
\nomenclature{$\mathrm{Gram}_g(b)$}{Gram-endomorphism of $b$ w.r.t $g$}
\begin{lemma}
 If $g$ and $g_0$ are inner products on $V$, then $\Gram_{g}(g_0)$ is self-adjoint with respect to $g$ and $g_0$.
 The Gram endomorphism is furthermore positive definite and satisfies $g_0(\Gram_{g_0}(g)v,w) = g(v,w)$ for all $v,w \in V$.
 
 For each invertible linear map $A \colon W \rightarrow V$, the Gram endomorphism satisfies 
 \begin{equation*}
     \Gram_{A^\ast g_0}(A^\ast g) = A^{-1} \circ \Gram_{g_0}(g) \circ A. 
 \end{equation*}
\end{lemma}
\begin{proof}
 We first prove that $\Gram_{g_0}(g)$ is self-adjoint with respect to $g_0$.
 This follows from the calculation
 \begin{align*}
     g_0(\Gram_{g_0}(g)v_1,v_2) &= g_0({g_0}_\sharp \circ g^\flat v_1,v_2) = g(v_1,v_2) = g(v_2,v_1) \\
     &= \dots = g_0(\Gram_{g_0}(g)v_2,v_1) \\
     &= g_0(v_1,\Gram_{g_0}(g)v_2).
 \end{align*}
 This calculation shows further that
 \begin{equation*}
     g_0(\Gram_{g_0}(g)v_1,v_2) = g(v_1,v_2).
 \end{equation*}
 In fact, since $g_0$ is non-degenerate, this equation uniquely determines $\Gram_{g_0}(g)$.
 
 Self-adjointness with respect to $g$ is now easily verified. 
 Using the previous identity in the first equation and self-adjointness with respect to $g_0$ in the second equation, we derive
 \begin{align*}
     g(\Gram_{g_0}(g)v_1,v_2) &= g_0(\Gram^2_{g_0}(g)v_1,v_2) = g_0(v_1,\Gram_{g_0}(g)^2v_2) \\
     &=  g_0(\Gram^2_{g_0}(g)v_2,v_1) =  g(\Gram_{g_0}(g)v_2,v_1).
 \end{align*}
 
 To prove that $\Gram_{g_0}(g)$ is positive definite, it suffices to prove that all eigenvalues are positive.
 Let $v$ be an eigenvector corresponding to the eigenvalue $\lambda$.
 Then 
 \begin{align*}
     \lambda g_0(v,v) = g_0(\Gram_{g_0}(g)v,v) = g(v,v) > 0
 \end{align*}
 implies that $\lambda$ must be positive.
 
 The final identity follows from
 \begin{align*}
     A^\ast g_0(A^{-1}\circ \Gram_{g_0}(g) \circ A w_1,w_2) &= g_0(\Gram_{g_0}(g)Aw_1,Aw_2) \\
     &= g(Aw_1,Aw_2) \\
     &= A^\ast g(w_1,w_2)
 \end{align*}
 and the fact that $A^\ast g_0$ is non-degenerate.
\end{proof}

\begin{definition}\label{TangentPreGauge - Definition}
 Let $g$ and $g_0$ be two inner products on $V$. 
 We define the \emph{pre-gauge} map $\tau_{g,g_0} \in \mathrm{End}(V)$
 \nomenclature{$\tau_{g,g_0}$}{Pre-Gauge map on a single vector space} to be the positive definite square root of $\Gram_{g_0}(g)$.
\end{definition}
 The pre-gauge map inherits many properties of the Gram endomorphism.
\begin{lemma}\label{Pregauge Properties - Lemma}
 The pre-gauge map $\tau_{g,g_0}$ is self-adjoint with respect to $g$ and $g_0$. 
 It is an isometry $\tau_{g,g_0} \colon (V,g) \rightarrow (V,g_0)$ and positive definite.
 Furthermore, $\tau_{g,g_0}$ is uniquely determined by these properties. 
 In particular, $\tau_{g_0,g} = \tau_{g,g_0}^{-1}$ and $\tau_{A^\ast g,A^\ast g_0} = A^{-1} \circ \tau_{g,g_0} \circ A$.
\end{lemma}
\begin{proof}
 The root of a self-adjoint endomorphism is self-adjoint. 
 The pre-gauge map $\tau_{g,g_0}$ is positive definite by definition.
 It is an isometry because
 \begin{equation*}
     g_0(\tau_{g,g_0}v,\tau_{g,g_0}w) = g_0(\tau^2_{g,g_0}v,w) = g_0(\mathrm{Gram}_{g_0}(g)v,w) = g(v,w).
 \end{equation*}
 
 An endomorphism $A$ that is self-adjoint with respect to $g_0$, positive definite, and an isometry $(V,g) \rightarrow (V,g_0)$ necessarily squares to the Gram endomorphism $\mathrm{Gram}_{g_0}(g)$. 
 Since the square root is unique on the set of positive definite self-adjoint operators, the operators $A$ and $\tau_{g,g_0}$ agree.
 
 The equations
 \begin{equation*}
     \tau_{g_0,g} = \tau_{g,g_0}^{-1} \qquad \text{and} \qquad \tau_{A^\ast g, A^\ast g_0} = A^{-1} \circ \tau_{g,g_0} \circ A
 \end{equation*}
 hold because the operators on the right hand sides satisfy the properties of the operators on the left hand side.
\end{proof}
In the most important example $(V,g) = (\R^n,\euclmetric)$, the Gram endomorphism $\Gram_{\euclmetric}(g)$ represented in the standard basis of $\R^n$ gives the more commonly known \emph{Gram matrix} of the metric $g$.
This coordinate expression is useful for the proof that the pre-gauge map depends smoothly on the inner products.
\begin{proposition}\label{Smooth Square Root - Prop}
 The map $A \mapsto \sqrt{A}$ that assigns a positive definite, self-adjoint matrix its positive definite square root is smooth.
\end{proposition}
\begin{proof}
The exponential map 
\begin{equation*}
   \exp \colon \{A \in \R^{n\times n}\, : \, A^\ast = A\} \rightarrow \{A \in \R^{n\times n}\, : \, A^\ast = A, \, A > 0\}    
\end{equation*}
 is surjective because for each self-adjoint, positive definite $A$ there is a $U \in \Or(n)$ and positive numbers $\lambda_1, \dots, \lambda_n$ such that $U^\ast A U = \mathrm{dia}(\lambda_1, \dots, \lambda_n) =: T$. 
If $B = U \log(T) U^\ast$, then $\exp(B) = A$.

We claim that, for all self-adjoint, positive definite $A$, the differential 
\begin{equation*}
 D_A\exp \colon \{A \in \R^{n\times n}\, : \, A^\ast = A\} \rightarrow \{A \in \R^{n\times n}\, : \, A^\ast = A\}
\end{equation*}
is surjective and therefore an isomorphism. 
Indeed, from
\begin{equation*}
 (D_A\exp)(H) = \sum_{k=1}^\infty \frac{1}{k!} \sum_{j=1}^k A^{j-1}HA^{k-j}
\end{equation*}
we conclude 
\begin{equation*}
 U^\ast D_A\exp(H) U = D_{U^\ast A U}(U^\ast H U).
\end{equation*}
Since $H$ is arbitrary, we may replace $H$ by $U H U^\ast$ and assume that $A = T$ is diagonal with entries $\lambda_1,\dots,\lambda_n$.
If $H=(h_{ij})$, then $(T^k H T^{l})_{ij} = \lambda_i^k \lambda_j^l h_{ij}$. 
Thus,
\begin{align*}
 (D_T\exp)(H) &= \left( \left(\sum_{k=1}^\infty \frac{1}{k!} \sum_{\alpha= 1}^k \lambda^{\alpha-1}_i \lambda_j^{k-\alpha}\right)h_{ij}\right) \\
              &=: (r_{ij}h_{ij}),
\end{align*}
where $r_{ij} > 0$. 
If $B = (b_{ij})$ is given, then, for $H= (h_{ij}) \deff (r_{ij}^{-1}b_{ij})$, we have $(D_T\exp)(H) = B$.

By the open mapping theorem, $\exp$ is locally invertible around each positive definite, self-adjoint matrix $A$.
Call the local inverse $\mathrm{log}$. Then $\sqrt{A} = \exp(1/2\mathrm{log}(A))$ is smooth.
\end{proof}

\begin{cor}\label{PreGauge Smooth - Cor}
 For a fixed but arbitrary inner product $g_0$, the maps 
 \begin{gather*}
     \Pos(V) \rightarrow \mathrm{Aut}(V), \\
     g \mapsto \tau_{g,g_0} \quad \text{ and } \quad g \mapsto \tau_{g_0,g}
 \end{gather*}
 are smooth.
\end{cor}
\begin{proof}
 From $\tau_{g,g_0} = \tau_{g_0,g}^{-1}$ follows that one map is smooth if and only if the other one is.
 
 Fix an isometry $A \colon (\R^n,\euclmetric) \rightarrow (V,g_0)$.
 Clearly, $g \mapsto \tau_{g,g_0}$ is smooth if and only if $g \mapsto A^{-1} \circ \tau_{g_0,g} \circ A$ is smooth.
 By Lemma \ref{Pregauge Properties - Lemma}, the latter map agrees with $g \mapsto \tau_{A^\ast g, \euclmetric}$, which is smooth by Proposition \ref{Smooth Square Root - Prop}.
\end{proof}

We wish to extend the pre-gauge map to bundles. 
This is done, most transparently, in terms of frame bundles and their associated bundles.
Let $\mathrm{Pos}(\R^n)$ be the set of all positive definite, symmetric bilinear forms on $\R^n$.
It is a $\GL(n)$ space with left action $A.b \deff (A^{-1})^\ast b$.
Let $\mathrm{Iso}^+(\R^n)$ be the set of all orientation preserving isomorphisms on $\R^n$. 
It is a $\GL(n)$ space with left action $A.\varphi \deff A^{-1}\varphi A$.

Let $\mathrm{Pos}(M)$ 
\nomenclature{$\mathrm{Pos}(M)$}{Bundle of symmetric, positive definite bilinear forms on $M$}
be the bundle of all positive definite, symmetric bilinear forms on $TM$ and $\mathrm{Iso}^+(TM)$ be the bundle of all bundle orientation preserving automorphisms of $TM$. 
These bundles can be obtained from $\GL(M)$ by using the described group actions in the Borel construction.
Lemma \ref{Pregauge Properties - Lemma} implies that the pre-gauge is a $\GL(n)$ equivariant map between
$\mathrm{Pos}(\R^n)^2$ and $\mathrm{Iso}^+(\R^n)$. 
It gives rise to bundle maps
\begin{equation*}
 \xymatrix@C+2em{\mathrm{Pos}(M)\times_M \mathrm{Pos}(M) \ar[r]^-{\tau_{\placeholder \, , \, \placeholder}} 
 & \mathrm{Iso}^+(M) \\
 \GL(M) \tilde{\times} \mathrm{Pos}(\R^n)^2 \ar[u]^\cong \ar[r]^-{\id \times \tau_{\placeholder \, , \, \placeholder}}
 & \GL(M) \tilde{\times} \mathrm{Iso}^+(\R^n). \ar[u]_\cong}
\end{equation*} 

Precomposing the pre-gauge map yields a map
\begin{equation*}
    \xymatrix@C+1.9em{\Pos(M) \times_M \Pos(M) \times_M \GL(M) \ar[r]^-{\tau_{\placeholder,\placeholder}} & \mathrm{Iso}^+(M) \times_M \GL(M) \ar[r]^-\circ & \GL(M)}.
\end{equation*}
Convex combination yields a homotopy equivalence $\Pos(M)\times_M \Pos(M) \simeq M \times \{g_0,g_0\}$, where $g_0$ is a section of $\Pos(M)\rightarrow M$, in other words, a Riemannian metric on $M$. 
Covering theory now provides a unique lift
\begin{equation*}
    \xymatrix@C+2em{ \Pos(M) \times_M \Pos(M) \times_M \widetilde{\GL}^+(M) \ar[d]_{\id \times }^{\rho} \ar@{-->}[rr]^-{\tilde{\tau}_{\placeholder,\placeholder} \circ \, \placeholder } && \widetilde{\GL}^+(M) \ar[d]^\rho \\ \Pos(M) \times_M \Pos(M) \times_M \GL^+(M) \ar[rr]^-{\tau_{\placeholder,\placeholder} \circ \, \placeholder } && \GL^+(M)}
\end{equation*}
determined by $\tilde{\tau}_{g_0,g_0} = \id$.

Let $\SOr(M,g)$ and $\SOr(M,g_0)$ be two orthonormal frame bundles corresponding to the Riemannian metrics $g, g_0$. 
The pre-gauge map $\tau_{g,g_0}$ maps $\SOr(M,g)$ to $\SOr(M,g_0)$ and the lift $\tilde{\tau}_{g,g_0}$ therefore maps $P_{\Spin,g}$ to $P_{\Spin,g_0}$.
The lifted pre-gauge map yields a graded algebra homomorphism
\nomenclature{$\Phi_{g,g_0}$}{Pre-Gauge map between Clifford algebra bundles}
\nomenclature{$\Phi_{g}$}{$\Phi_{g,g_0}$ but base point $g_0$ omitted}
\begin{equation*}
 \xymatrix{Cl(TM,g) \ar@{=}[d] \ar[rr]^{\Phi_{g,g_0}} &&  Cl(TM,g_0) \ar@{=}[d] \\
 P_{\Spin,g} \times_{\mathrm{conj}}Cl_{d,0} \ar[rr]^-{\tilde{\tau}_{g,g_0}\times \id}&& P_{\Spin,g_0} \times_{\mathrm{conj}}Cl_{d,0} }
\end{equation*} 
and a $Cl_{d,0}$-right linear map
\nomenclature{$\preGauge_{g,g_0}$}{Pre-Gauge map between spinor bundles}
\nomenclature{$\preGauge_{g}$}{$\preGauge_{g,g_0}$but base point $g_0$ omitted}
\begin{equation*}
 \xymatrix{ \Spinor_g \ar@{=}[d] \ar[rr]^{\preGauge_{g,g_0}} && \Spinor_{g_0} \ar@{=}[d] \\
 P_{\Spin ,g} \times_l Cl_{d,0} \ar[rr]^-{\tilde{\tau}_{g,g_0} \times \id} & & P_{\Spin, g_0} \times_l Cl_{d,0}.}
\end{equation*}
Of course, these identifications are isometric, grading preserving maps. 
Furthermore, the identifications maps between the spinor bundles are $Cl(TM,g)$-$Cl(TM,g_0)$ left-module maps.
That means the diagram 
\begin{equation*}
 \xymatrix{Cl(TM,g) \times \Spinor_g \ar[r]^-\Cliffmult \ar[d]_{\Phi_{\placeholder,g_0} \times}^{\preGauge_{g,g_0}} & \Spinor_g \ar[d]^{\preGauge_{g,g_0}} \\
 Cl(TM,g_0) \times \Spinor_{g_0} \ar[r]^-\Cliffmult & \Spinor_{g_0}}
\end{equation*}
commutes. Consequently, using these identifications, we get (set-theoretical) bijections
\begin{equation}\label{eq: UniversalClifford}
 \xymatrix{Cl(M) \deff \bigsqcup_{g \in \Riem(M)} Cl(TM,g) \ar[rr]^-{\Phi_{\placeholder,g_0}} & & Cl(TM,g_0)\times \Riem(M)}
\end{equation}
\begin{equation}\label{eq: UniversalSpinor}
    \xymatrix{\Spinor \deff \bigsqcup_{g \in \Riem(M)} \Spinor_g \ar[rr]^{\preGauge_{\placeholder,g_0}} & & \Spinor_{g_0} \times \Riem(M).}
\end{equation}
We call $\preGauge$ the \emph{pre-gauge} map.
If $g_0$ is fixed, then we usually write $\Phi_{g}$ and $\preGauge_{g}$ for $\Phi_{g,g_0}$ and $\preGauge_{g,g_0}$, respectively.

While the maps $\Phi$ and $\preGauge$ are fibre-wise isometries, the pre-gauge map does not induce isometries between the associated Hilbert spaces of square integrable spinors. 
The reason is that the volume form in the definition of the inner product on the sections are also different.
In order to fix that we have to make a conformal change.
\begin{definition}
 For $g, g_0 \in \Riem(M)$ let $\alpha^2_{g,g_0} \colon M \rightarrow M$ be the density of $\mathrm{vol}_g$ with respect to $\mathrm{vol}_{g_0}$.
 We call $\fullGauge_{g,g_0} \deff \alpha_{g , g_0} \cdot \preGauge_{g,g_0}$ the \emph{Gauge map}.
 \nomenclature{$\fullGauge_{g,g_0}$}{Gauge map between spinor bundles}
 \nomenclature{$\fullGauge_{g}$}{$\fullGauge_{g,g_0}$ but base point $g_0$ fixed}
 We simply write $\fullGauge_g$ if $g_0$ is fixed. 
\end{definition}
The following lemma is easily verified.
\begin{lemma}
 The induced map of $\fullGauge_g$ yields an isometry of Hilbert spaces
 \begin{equation*}
  \fullGauge_g \colon L^2(M,g;\Spinor_g) \rightarrow L^2(M,g_0;\Spinor_{g_0}).
 \end{equation*}
\end{lemma}
\begin{rem}
 Note that if $M$ is not compact, then $\preGauge_{g,g_0} \colon \Gamma_c(\Spinor_g) \rightarrow \Gamma_c(\Spinor_{g_0})$ may not extend to a bounded linear map between $L^2(\Spinor_g) \rightarrow L^2(\Spinor_{g_0})$. 
 It extends to a bounded operator if and only if the map $\alpha_{g,g_0}$ is bounded.
\end{rem}


\begin{definition}\label{Universal Bundles - Def}
 Let 
 \begin{equation*}
     Cl(M) \rightarrow \Riem(M) \times M \qquad \text{and} \qquad \Spinor \rightarrow \Riem(M)\times M
 \end{equation*}
 be the bundles whose topologies are the unique ones such that $\Phi_{\placeholder, g_0}$ and $\preGauge_{\placeholder,g_0}$ in equation (\ref{eq: UniversalClifford}) and (\ref{eq: UniversalSpinor}) become homeomorphisms if the right hand side is equipped with the product topology.
\end{definition}

Informally, we think of these as bundles of bundles over $\Riem(M)$.
From now on, we assume that $M$ has no boundary.

\begin{definition}
 Let $\PseudOp^k(\Spinor) \rightarrow \Riem(M)$ be the bundle that is trivialised by 
 \begin{equation}
  \fullGauge \circ - \circ \fullGauge^{-1} \colon \PseudOp^k(\Spinor) \rightarrow \Riem(M) \times \PseudOp^k(\Spinor_{g_0});
 \end{equation}
 here, $\Riem(M)$ carries the smooth Fréchet topology (the smooth weak topology) and $\PseudOp^k(\Spinor_{g_0})$ carries the amplitude topology introduced in Appendix \ref{Appendix Pseudos - Chapter}.
 The union $\PseudOp^\bullet(\Spinor)$ carries the topology of a (trivial) Fréchet space bundle, even a 
 Fréchet algebra bundle if $M$ is compact. 
 The same construction applies to the Atiyah-Singer closure$\PseudOpCl^k(\Spinor)$ defined in Appendix \ref{Appendix Pseudos - Chapter}.
\end{definition}
\begin{rem}
 If $M$ is compact, we can trivialise $\PseudOp^\bullet(\Spinor)$ also with the conjugation of $\preGauge$.
\end{rem}
We are specifically interested in a certain subspace of $\PseudOp^1(\Spinor)$.
\begin{definition}\label{PseudoDiracOperator - Def}
 An element $P \in \PseudOp^1(\Spinor)$ is called a \emph{pseudo Dirac operator} if it is a
 symmetric, odd, $\Cliff[d][0]$-linear operator, whose restriction to each fibre $P_g$ has the same principal symbol as the Dirac operator $\Dirac_g$ that is $\symb_1(P_g) = i\Cliffmult(g_\sharp(\cdot) )$,
 \nomenclature{$\symb_1(P)$}{Principal Symbol of an order $1$ operator} where $\Cliffmult = \Cliffmult_{g}$ denotes the left-module action of $Cl(TM,g)$ on $\Spinor_g$.
 The set of all pseudo Dirac operators is denoted by $\PseudDir$. 
 \nomenclature{$\PseudDir$}{Set of pseudo Dirac operators on $M$}
\end{definition} 
Note that conjugation with $\fullGauge$ does not map $\PseudDir$ to $\PseudDir[ ][\Spinor_{g_0}] \times \Riem(M)$. 
Indeed, since $\fullGauge_g$ is a bundle endomorphism, symbol calculus yields
\begin{align*}
 \symb_1(\fullGauge_g \circ P \circ \fullGauge_g^{-1}) &= \fullGauge(g) \circ \symb_1(P)\circ \fullGauge_g^{-1} \\
 &= \fullGauge_g \circ  i\Cliffmult_g({g}_\sharp(\placeholder)) \circ \fullGauge_g^{-1} \\
 &= \preGauge_g \circ  i\Cliffmult_g({g}_\sharp(\placeholder)) \circ \preGauge_g^{-1} \\
 &= i\Cliffmult_{g_0}(\Phi_g \circ g_\sharp(\placeholder) ),
\end{align*}
which shows that $\fullGauge_g \circ P \circ \fullGauge_g^{-1}$ has not the correct principal symbol, which would be $i\Cliffmult_{g_0}({g_0}_\sharp(\placeholder))$.
Nevertheless, $\PseudDir$ is a fibre bundle over $\Riem(M)$ with fibre $\PseudDir_g = \Psi\mathrm{Dir}(\Spinor_{g})$ the affine space of all pseudo Dirac operators on $\Spinor_{g}$. 
\nomenclature{$\PseudDir_g$}{pseudo Dirac operators on $\Spinor_{g}$}
It is an affine bundle over $\Pseud[0][\Spinor][\Spinor]$ because the principal symbol is specified. 
It is even a trivial bundle because it has a globally non-vanishing smooth section, which sends a Riemannian metric to its Dirac operator.

\begin{definition}
 A map $P \colon \R^n \rightarrow \PseudOp^k(\Spinor,\Spinor)$ is smooth, if 
 \begin{align*}
     \fullGauge_\ast (P) \colon \R^n &\rightarrow \Riem(M) \times \PseudOp^k(\Spinor_{g_0},\Spinor_{g_0}), \\
     t &\mapsto \fullGauge_{g_{P(t)},g_0} \circ P(t) \circ \fullGauge^{-1}_{g_P(t),g_0} 
 \end{align*}
 is a smooth map into a Fréchet space. 
 Here $g_{P(t)}$ denotes the underlying metric of $P(t)$.
\end{definition}
Note that this definition does not depend on the choice of $g_0$, because different choices result in conjugation with a fixed, smooth vector bundle isomorphism.

\begin{theorem}\label{DiracOperatorSmooth - Theorem}
 The Dirac operator defines a smooth section
 \begin{equation*}
  \Dirac \colon \Riem(M) \rightarrow \PseudDir
 \end{equation*}
 in the sense that it maps smooth maps $\R^n \rightarrow \Riem(M)$ to smooth maps $\R^n \rightarrow \PseudOp^1(\Spinor)$.
\end{theorem}
\begin{proof}
 We need to show that $g \mapsto \fullGauge_{g,g_0} \circ \Dirac_g \circ \fullGauge_{g,g_0}^{-1} =: {\fullGauge_{g,g_0}}_\ast \Dirac_g$ is a smooth map
 \begin{equation*}
     \Riem(M) \rightarrow \Pseud[1][\Spinor_{g_0}][\Spinor_{g_0}].
 \end{equation*}
 The topologies of $\Riem(M)$ and $\PseudOp(\Spinor)$ are generated by a family of semi-norms defined in  in terms of local coordinates.
 Thus, it suffices to show in local coordinates that smooth Riemannian metrics are send to smooth amplitudes. 
 The calculation
 \begin{align*}
      {\fullGauge_{g,g_0}}_\ast \Dirac_g &= \alpha_{g,g_0} \preGauge_{g,g_0} \circ \Dirac_g \circ \alpha_{g,g_0}^{-1} \preGauge^{-1}_{g,g_0} \\
      &= {\preGauge_{g,g_0}}_\ast \Dirac_g + \preGauge_{g,g_0} \circ \Cliffmult_g(g_\sharp(-\alpha_{g,g_0}^{-1}\diff \alpha_{g,g_0}) \otimes \id ) \circ \preGauge_{g,g_0}^{-1} \\
      &= {\preGauge_{g,g_0}}_\ast \Dirac_g - \preGauge_{g,g_0} \circ \Cliffmult_g(\alpha_{g,g_0}^{-1}\mathrm{grad}_g(\alpha_{g,g_0}) \otimes \id ) \circ \preGauge_{g,g_0}^{-1} \\
      &= {\preGauge_{g,g_0}}_\ast \Dirac_g - \alpha_{g,g_0}^{-1}\Cliffmult_{g_0}(\Phi_{g,g_0}(\mathrm{grad}_g(\alpha_{g,g_0})) \otimes \id)
 \end{align*}
 allows us to conjugate with $\preGauge$ instead because adding a smooth endomorphism does not affect the smoothness of the map.
 
 Consider the following diagram
 \begin{equation*}
     \xymatrix@C+1.4em{ \Gamma(\Spinor_g) \ar[r]^-{\nabla^{\Spinor_g}} \ar[d]_{\preGauge_{g,g_0}} & \Gamma(T^\vee M \otimes \Spinor_g) \ar[r]^{g_\sharp \otimes \id} \ar[d]_{\preGauge_{g,g_0}} & \Gamma(TM \otimes \Spinor_{g}) \ar[rr]^{\Cliffmult_g} \ar[d]_{\preGauge_{g,g_0}} && \Gamma(\Spinor_g) \ar[d]_{\preGauge_{g,g_0}} \\
     \Gamma(\Spinor_{g_0}) \ar[r]^-{{\preGauge_{g,g_0}}_\ast \nabla^{\Spinor_{g}}} & \Gamma(T^\vee M \otimes \Spinor_{g_0}) \ar[r]^{g_\sharp \otimes \id} & \Gamma(TM \otimes \Spinor_{g_0}) \ar[rr]^{\Cliffmult_{g_0} \circ (\Phi_{g,g_0} \otimes \id)} &&  \Gamma(\Spinor_{g_0} ).}
 \end{equation*}
 It is clear that $g \mapsto g_\sharp \otimes \id$ and $g \mapsto \Cliffmult_{g_0} \circ (\Phi_{g,g_0} \otimes \id)$ are induced by smooth maps of vector bundles (because $\tau_{\placeholder, g_0}$ depends smoothly on the first entry).
 Since the Fréchet space of smooth bundle endomorphism embeds into the Fréchet space of pseudo differential operators of order zero (equipped with the amplitude topology), see Example \ref{ClassicalPseudos - Example} together with Theorem \ref{PseudosFrechetNorm - Theorem},
 the assignments $g \mapsto g_\sharp \otimes \id$ and $g \mapsto \Cliffmult_{g_0} \circ (\Phi_{g,g_0} \otimes \id)$ yield smooth maps 
 \begin{equation*}
     \Riem(M) \rightarrow \Pseud[0][T^\vee M \otimes \Spinor_{g_0}][TM \otimes \Spinor_{g_0}]  \bigr)  
 \end{equation*}
and
\begin{equation*}
    \Riem(M) \rightarrow \Pseud[0][TM \otimes \Spinor_{g_0}][\Spinor_{g_0}].
\end{equation*}
 
 The difference of two connections is an endomorphism-valued $1$-form.
 Thus it suffices to show that the map
 \begin{align*}
     \Riem(M) &\rightarrow \Omega^1(M;\mathrm{End} \,\Spinor_{g_0}), \\
     g &\mapsto {\preGauge_{g,g_0}}_\ast \nabla^{\Spinor_g} - \nabla^{\Spinor_{g_0}} = \preGauge_{g,g_0} \circ \nabla^{\Spinor_g} \circ \preGauge_{g,g_0}^{-1} - \nabla^{\Spinor_{g_0}}
 \end{align*}
 is smooth.
 
 First note that the analogous map for the Levi-Cevita connections 
 \begin{align*}
     \Riem(M) &\rightarrow \Omega^1(M;\mathrm{End}\, TM ) \\
     g &\mapsto {\preGauge_{g,g_0}}_\ast \nabla^{L.C.,g} - \nabla^{L.C.,g_0}
 \end{align*}
 is smooth due to the following calculation in local coordinates.
 
 Let $\partial_1, \dots, \partial_d$ be the standard basis in the chart domain $U$.
 In these coordinates we can express $g$ and $g_0$ as matrices $(g_{ij})$ and $(g_{0 \, ij})$, so that ${g_{(0)}}_{ij} = g_{(0)}(\partial_i,\partial_j)$.
 The Levi-Cevita connection is given by
 \begin{align*}
     \nabla_i^{L.C.,g} &= \nabla_{\partial_i}^{L.C.,g} = \partial_i + {}^g\Gamma_i,
 \end{align*}
 where the Christoffel symbols $\Gamma_i = (\Gamma_{ij}^k)$ are an algebraic expression of the $1$-jet of $g$. 
 More precisely, we have the well-known formula
 \begin{equation*}
     2\Gamma_{ij}^k = \sum_{\alpha = 1}^d g^{k\alpha} \bigl(\partial_j g_{i\alpha} + \partial_i g_{\alpha j} - \partial_\alpha g_{ij}\bigr),
 \end{equation*}
 where $g^{k\alpha}$ denotes the $k\alpha$-entry of $(g_{ij})^{-1}$.
 Furthermore, in this basis, the Gram endomorphism is given by
 \begin{equation*}
     \mathrm{Gram}_{g_0}(g) = ({g_0}_{ij})^{-1} \cdot (g_{ij}),
 \end{equation*}
 which depends smoothly on $g$.
 By Corollary \ref{PreGauge Smooth - Cor} also $\tau_{g,g_0}$ depends smoothly on $g$.
 A straightforward calculation now yields
 \begin{align*}
     {\preGauge_{g,g_0}}_\ast \nabla^{L.C.,g} &= \preGauge_{g,g_0} \circ \nabla^{L.C.,g} \circ \preGauge_{g,g_0}^{-1} \\
     &= \tau_{g,g_0} \cdot (\partial_i + {}^g\Gamma_i) \cdot \tau_{g,g_0}^{-1} \\
     &= \partial_i + \tau_{g,g_0} \cdot {}^g\Gamma_i \cdot \tau_{g,g_0}^{-1} - \tau_{g,g_0} \cdot \tau_{g,g_0}^{-1} \cdot \partial_i \tau_{g,g_0} \cdot \tau_{g,g_0}^{-1} \\
     &= \partial_i + \tau_{g,g_0} \cdot {}^g\Gamma_i \cdot \tau_{g,g_0}^{-1} - \partial_i \tau_{g,g_0} \cdot \tau_{g,g_0}^{-1}.
 \end{align*}
 As all summands depend smoothly on $g$, the map $g \mapsto {\preGauge_{g,g_0}}_\ast \nabla^{L.C.,g} - \nabla^{L.C.,g_0}$ is smooth.
 
 Let $e_1,\dots, e_d$ be a $g_0$-orthonormal frame on $U$ and set $f_j = \tau_{g_0,g}e_j$.
 These orthonormal frames give canonical (up to a sign) orthonormal frames on $\Spinor_{g_0}$ and $\Spinor_g$, which we denote with $(e_I)$ and $(f_I)$, respectively.
 Here, $I$ runs over all ordered tupels $I = (i_1 < \dots < i_r)$ with $1 \leq i_j \leq d$.  
 Using the formula for the Christoffel symbols of the spinor connection, which can be found in \cite{LawsonMichelsonSpin}*{Theorem II.4.14}, we get
 \begin{align*}
     2 \cdot \Bigl({\preGauge_{g,g_0}}_\ast \nabla^{\Spinor_g}\Bigr)_{e_i}(e_I) &= 2 \cdot \preGauge_{g,g_0} \circ \nabla^{\Spinor_g}_{e_i} \circ \preGauge_{g,g_0}^{-1}(e_I) \\
     &= 2 \cdot \preGauge_{g,g_0} \Bigl( \nabla^{\Spinor_g}_{e_i} f_I\Bigr) \\
     &= \preGauge_{g,g_0} \Bigl( \bigl(\sum_{\alpha < \beta } g(\nabla_{e_i}^{L.C.,g} f_\alpha, f_\beta ) \cdot f_\alpha \cdot f_\beta \bigr) \cdot_{g} f_I\Bigr) \\
     &= \Bigl(\sum_{\alpha < \beta } g(\nabla_{e_i}^{L.C.,g} \tau_{g_0,g}e_\alpha, \tau_{g_0,g}e_\beta ) \cdot e_\alpha \cdot e_\beta \Bigr) \cdot e_I.
 \end{align*}
 We already know that $g(\nabla_{e_i}^{L.C.,g} e_\alpha, e_\beta) = {}^g\Gamma_{i\alpha}^\beta$ depends smoothly on $g$.
 Thus, the assignment
 \begin{equation*}
     g \mapsto {\preGauge_{g,g_0}}_\ast \nabla^{\Spinor_g} - \nabla^{\Spinor_{g_0}}
 \end{equation*}
 is a smooth map $\Riem(M) \rightarrow \Gamma\bigl(M;\mathrm{Hom}(\Spinor_{g_0},T^\vee M \otimes \Spinor_{g_0})\bigr)$
 and hence it is also a smooth map
 \begin{equation*}
     \Riem(M) \rightarrow \Pseud[0][\Spinor_g][T^\vee M \otimes \Spinor_g].
 \end{equation*}
 
 This implies that $g \mapsto \Dirac_g$ is smooth as it is a composition of smooth maps.
\end{proof}

The next result is classical and follows immediately from the theory developed in \cite{LawsonMichelsonSpin}*{Chapter III.5}.

\begin{proposition}
 Let $M$ be a closed manifold and let $P \in \PseudDir$ be a pseudo Dirac operator with underlying Riemannian metric $g \in \Riem(M)$.
 Then $P$ is a self-adjoint, unbounded Fredholm operator on $L^2(\Spinor_g)$ with domain $H^1(\Spinor_g)$.
\end{proposition}
\begin{proof}
 By definition, $P$ is symmetric, so $P - \Dirac_g$ is a symmetric pseudo differential operator of order zero on $\Spinor_g$.
 Thus the difference extends to a bounded, symmetric, and hence self-adjoint, operator on $L^2(\Spinor_g)$.
 Since the sum of a self-adjoint operator with a self-adjoint bounded operator is again self-adjoint \cite{higson2000analytic}*{Exercise 1.9.20}, $P = \Dirac_g + (P - \Dirac_g)$ is self-adjoint.
 
 G\aa rdings inequality, see \cite{LawsonMichelsonSpin}*{Thm III.5.2}, implies that the minimal domain of $P$ agrees with $H^1(\Spinor_g)$.
 It is classical that elliptic pseudo differential operators on vector bundles over closed manifolds are unbounded Fredholm operators, see \cite{LawsonMichelsonSpin}*{Thm III.5.2}.
\end{proof}

\subsubsection{Relations to \textit{K}-theory}

Our interests in pseudo Dirac operators arise from their connection to $KO$-theory, which we will outline in this section.
All presented results have already appeared in \cite{ebert2017indexdiff} in a different language and we do not claim originality for the results presented here; rather we translate his results into our language, fill some gaps, and give more details.
Needless to say, the presentation here closely follows \cite{ebert2017indexdiff}.
Throughout this section, we assume that $M^d$ is a \emph{closed} spin manifold of positive dimension.

We start with recalling the connection between $KO$-theory and Clifford-linear Fredholm operators.
Let $H$ be a separable real Hilbert space with $\Z_2$-grading $\iota$ and graded $\Cliff[d][0]$-right action $\mathbf{r}$, which we assume to be a $\ast$-homomorphism.
For brevity, we simply call it a $\Cliff[d][0]$-Hilbert space.
Call $\omega_{d,0} \deff \iota \mathbf{r}(e_1 \cdots e_d)$ the \emph{Chirality element}.
Recall that $\Cliff[d][0]$ has, up to isomorphism, exactly one irreducible representation if $d \not\equiv 0 \mod 4$, and exactly two irreducible representations if $d \equiv -1\mod 4$.
In the latter case, the irreducible representations are distinguished as different eigenspaces of the Chirality element.

The triple $(H,\iota,\mathbf{r})$ is called \emph{ample} if $H$ contains all irreducible representations with infinite multiplicity.

\begin{definition}
 Let $(H,\iota,\mathbf{r})$ be an ample $\Cliff[d][0]$-Hilbert space.
 We call an operator $F$ on $H$ an $\Cliff[d][0]$\emph{-Fredholm operator} if it is Clifford-linear, odd, and self-adjoint.
 If $d \equiv -1 \mod 4$, we further require that $\omega_{p,0}F\iota$ is neither essentially positive nor essentially negative, that means, there exist infinite dimensional, closed subspaces that are orthogonal to each other and on those $F$ restricts to a positive or negative operator, respectively.
 
 Let $\mathrm{Fred}^{d,0}(H)$ 
 \nomenclature{$\mathrm{Fred}^{d,0}(H)$}{$\Cliff[d][0]$-Fredholm operators}
 be the space of all $\Cliff[d][0]$-Fredholm operators on $H$ equipped with the norm topology.
\end{definition}

By a classical result of Atiyah and Singer, these spaces represent real $K$-theory.

\begin{theorem}[Theorem \cite{atiyah1969indexskew}]
 If $H$ is ample, then $\mathrm{Fred}^{d,0}(H)$ is a classifying space for $KO^{-d}$. That means there is a natural bijection
 \begin{equation*}
     \mathrm{ind} \colon [X,\mathrm{Fred}^{d,0}(H)] \rightarrow KO^{-d}(X)
 \end{equation*}
 for every compact $CW$-complex $X$.
\end{theorem}
A detailed construction of this map was carried out in the master's thesis of Jonathan Glöckle \cite{gloeckle2019master}*{Theorem 2.17}.
For more details on this space like Bott-periodicity or Morita equivalences, we refer the reader to \cite{ebert2017indexdiff}.
To identify the space of invertible pseudo Dirac operators as a classifying space for $KO$-theory, we first need to introduce auxiliary spaces.

\begin{definition}
 Let $\InvPseudDir$ be the subspace of $\PseudDir$ consisting of all invertible pseudo Dirac operators and denote $[-1,1]$ by $I$. 
  Provided $\InvPseudDir_{g_0}$ is non-empty for a fixed metric $g_0$, we define, for a fixed choice of base point $B \in \InvPseudDir$, the spaces
 \begin{equation*}
  \PathDirac[ ] \deff \{P \colon I \rightarrow \PseudDir \, : \, P \text{ smooth, } P(-1) = B, \ P(1) \in \InvPseudDir\}
 \end{equation*}
 and
 \begin{equation*}
  \PathDirac[ ]_{g_0} \deff \{P \colon I \rightarrow \PseudDir_{g_0} \, : \, P \text{ smooth, } P(-1) = B, \ P(1) \in \InvPseudDir_{g_0}\}.
 \end{equation*}
\end{definition}

\begin{lemma}
 If $\InvPseudDir_{g_0}$ is non-empty, then the evaluation map
 \begin{equation*}
  \mathrm{ev}_1 \colon \PathDirac[ ] \rightarrow \InvPseudDir
 \end{equation*}
 is a weak homotopy equivalence that also restricts to a weak homotopy equivalence between $\PathDirac[ ]_{g_0}$ and $\InvPseudDir_{g_0}$.
\end{lemma}
\begin{proof}
 We only prove the first statement; for the second statement, one simply needs to add the subscript ${}_{g_0}$ to $\InvPseudDir[ ]$ and the space $\mathcal{Y}$ defined below.

 Let $\mathcal{Y}$ be the set of all \emph{continuous} paths $P \colon I \rightarrow \PseudDir$ having the same boundary conditions as elements in $\PathDirac[ ]$. 
 Smoothing theory yields that the inclusion of $\PathDirac[ ]$ into $\mathcal{Y}$ is a weak homotopy equivalence. 
 
 The map $\mathrm{ev}_1 \colon \mathcal{Y} \rightarrow \InvPseudDir$ is a fibration. 
 Indeed, given a commutative square
 \begin{equation*}
  \xymatrix{ A \times 0\ar[rr]^F \ar[d] && \mathcal{Y} \ar[d]^{\mathrm{ev}_1} \\ A \times I \ar[rr]^{h} && \InvPseudDir }
 \end{equation*}
 we define the lift $\mathcal{H} \colon A \times I \rightarrow \mathcal{Y}$ by 
 \begin{equation*}
   \mathcal{H}(a,s)(t) \deff \frac{(1-t)}{2}\cdot F(a) + \frac{(1+t)}{2}\cdot h(a,s).
 \end{equation*}  
 This is a continuous map because the affine structure on $\PseudDir$ is (jointly) continuous with respect to the Fréchet structure. 
 It is also a lift for $\mathrm{ev}_1(\mathcal{H}) = h$, showing that $\mathrm{ev}_1$ is a fibration.
 
 The homotopy fibre of $\mathrm{ev}_1$ is therefore the preimage of $B$, which is the set of all closed curves in $\PseudDir$ that start (and end) at $B$. 
 This subspace is clearly convex and therefore contractible, implying that $\mathcal{Y}$ and $\InvPseudDir$ are weak homotopy equivalent. 
\end{proof}
The auxiliary path spaces $\PathDirac[{}]$ and $\PathDirac[{}]_{g_0}$ help to identify $\InvPseudDir$ as
classifying space for $KO$-theory.

We remark that, for a closed manifold $M^d$ of positive dimension $d$, the $\Cliff[d][0]$-Hilbert space $H = L^2(M,\Spinor_{g_0})$ is always ample.
Indeed, the right action of the Chirality element decomposes $\Spinor_{g_0}$ into non-zero subbundles $\Spinor_{g_0}^\pm$ on which the Chirality element acts as $\pm \id$, hence 
\begin{equation*}
    L^2(M,\Spinor_{g_0}) \cong L^2(M,\Spinor_{g_0}^+) \oplus L^2(M,\Spinor_{g_0}^-)
\end{equation*} 
is a decomposition into infinite dimensional eigenspaces of the eigenvalues $\pm1$ of the Chirality element.

The following theorem is a collection of results proved by Johannes Ebert in \cite{ebert2017indexdiff}*{Section 4.1}:
\begin{theorem}\label{Summary Eberts Results - Thm}
 Assume that $\InvPseudDir_{g_0}$ is non-empty. 
 Then the diagram of continuous maps
 \begin{equation*}
  \xymatrix{\PathDirac[{}] \ar[r] & \Omega \mathrm{Fred}^{d,0}\bigl(L^2(M,\Spinor_{g_0})\bigr) \ar@{=}[d] \\ \PathDirac[{}]_{g_0} \ar@{^{(}->}[u] \ar[r]^-\simeq & \Omega \mathrm{Fred}^{d,0}\bigl(L^2(M,\Spinor_{g_0})\bigr) }
\end{equation*}
 commutes. 
 The horizontal maps are given by
 \begin{equation*}
  P_t \mapsto \frac{\fullGauge_{g_t,g_0} P_t \fullGauge_{g_t,g_0}^{-1}}{(1 + (\fullGauge_{g_t,g_0}P_t\fullGauge_{g_t,g_0}^{-1})^2)^{1/2}}
\end{equation*}  
and the lower map is a weak homotopy equivalence. 
Here $\Omega \mathrm{Fred}^{d,0}(L^2(\Spinor_{g_0}))$ denotes the space of paths into $\mathrm{Fred}^{d,0}(L^2(\Spinor_{g_0}))$ whose start and end-points are invertible operators.
\end{theorem}
\begin{proof}
 For fixed $t$, these maps certainly take values in the space of self-adjoint, odd, Clifford-linear Fredholm operators.
 From the discussion in \cite{ebert2017indexdiff}*{Section 4.1} it is not clear why the horizontal map takes values in the subspace with the additional condition that $\omega F \iota$ is neither essentially positive nor essentially negative if $d \equiv -1 \mod 4$.
 
 In order to show this, first recall from \cite{LawsonMichelsonSpin}*{Theorem III.5.8} that the spectrum of every self-adjoint elliptic pseudo differential operator $P$ of positive order on a bundle over a closed manifold consists of countably many eigenvalues only, that the eigenspaces are finite dimensional, and that the eigenvalues go to infinity.
 The operator $P/\sqrt{1+P^2}$ is the unique operator that acts as $\lambda/\sqrt{1+\lambda^2} \cdot \id$ on the $\lambda$-eigenspace of $P$.
 
 The Chirality element satisfies
 \begin{equation*}
     \omega_{d,0}^2 = (-1)^{\frac{d(d+1)}{2}}  \quad \text{ and } \quad \omega_{d,0}^\ast = (-1)^{\frac{d(d+1)}{2}} \omega_{d,0}.
 \end{equation*}
 Thus, it is a self-adjoint, odd involution if $d \equiv -1 \mod 4$.
 Since $P$, and therefore $P/\sqrt{1 + P^2}$ too, commute with $\omega_{d,0}$, we can simultaneously diagonalise the two operators.  
 
 The grading $\evenodd$ anti-commutes with $P$, so the grading sends the $\lambda$-eigenspace of $P$ to the $-\lambda$-eigenspaces of $P$.
 The corresponding statement holds true for the $\pm 1$-eigenspaces of $\omega_{d,0}$ for the same reason.
 
 Abbreviate $P/\sqrt{1+P^2}$ to $F$.
 Since the grading $\evenodd$ commutes with $\omega_{d,0}F$, the operator $\omega_{d,0} F \evenodd$ is positive definite on
 \begin{equation*}
     H^+ \deff \bigoplus_{\lambda > 0} \mathrm{Eig}(P;\lambda) \cap \mathrm{Eig}(\omega_{d,0};1) \oplus \mathrm{Eig}(P;-\lambda) \cap \mathrm{Eig}(\omega_{d,0};-1) 
 \end{equation*}
 and negative definite on 
 \begin{equation*}
     H^- \deff \bigoplus_{\lambda > 0} \mathrm{Eig}(P;\lambda) \cap \mathrm{Eig}(\omega_{d,0};-1) \oplus \mathrm{Eig}(P;-\lambda) \cap \mathrm{Eig}(\omega_{d,0};1). 
 \end{equation*}
 By the discussion above, the two subspaces are infinite dimensional. 
 Thus, $\omega_{d,0} F \evenodd$ cannot be essentially positive or essentially negative.
 
 The upper horizontal map is continuous because conjugation with invertible functions is a continuous operation as $\PseudOp^{\bullet}(\Spinor_{g_0})$ is a Fréchet algebra, and bounded functional calculus is continuous by \cite{ebert2017indexdiff}*{Proposition 3.7}. 
 The weak equivalence statement is \cite{ebert2017indexdiff}*{Proposition 4.3}. 
\end{proof}

The space $\InvPseudDir_{g_0}$ is unsuitable for our purposes.
Because we want to compare different metrics, we need the bigger space $\InvPseudDir$. 
Luckily, the canonical inclusion is a weak homotopy equivalence.
This is the main result of \cite{ebert2017indexdiff}*{Section 4.2} formulated in our language. 

\begin{theorem}\label{Inclusion PseudDir fixed WHE - Thm}
 The inclusion $\InvPseudDir_{g_0} \hookrightarrow \InvPseudDir$ is a weak homotopy equivalence.
\end{theorem}

The identification between $\InvPseudDir$ and $\Omega \mathrm{Fred}^{d,0}(L^2(\Spinor_{g_0}))$ is the space-level version of the index difference of Hitchin \cite{hitchin1974harmonic}.

\begin{definition}\label{inddifHitchin - Definition}
 Fix base points $g_0 \in \Riem(M)$ and $B \in \InvPseudDir$.
 For any choice of section $s \colon \InvPseudDir \rightarrow \mathcal{X}$ of $\mathrm{ev}_1$, we call the composition
 \begin{equation*}
     \xymatrix{\mathrm{inddif}(\placeholder,B) \colon \InvPseudDir \ar@<0.5ex>[rr]^-s && \mathcal{X} \ar[rr]^-\simeq \ar@<0.5ex>[ll]^-{\mathrm{ev}_1} && \Omega \mathrm{Fred}^{d,0}(L^2(\Spinor_{g_0})) }
 \end{equation*}
 the \emph{(operator) index difference}.
 
 If $g_0$ is a psc-metric then we call the pre-composition with the Dirac operator
 \begin{equation*}
     \xymatrix{\alpha(\placeholder,g_0) \colon \Riem^+(M) \ar[rr]^-{\Dirac_{(\placeholder)}} && \InvPseudDir \ar[rr]^-{\mathrm{inddif}(\placeholder, \Dirac_{g_0})} && \Omega \mathrm{Fred}^{d,0}(L^2(\Spinor_{g_0}))  }
 \end{equation*}
 the \emph{index difference}.
\end{definition}

We summarise our achievement of this subsection in the following theorem.

\begin{theorem}\label{InvPseudDir Classifying Space - Thm}
  Let $M^d$ be a closed spin manifold of positive dimension.
  If the space of all invertible pseudo Dirac operators $\InvPseudDir$ is non-empty, then it is a classifying space for $KO^{-(d+1)}(\mathrm{pt})$. 
  The natural transformation is given by the \emph{Hitchin index difference}
  \begin{equation*}
      \xymatrix{ [X;\InvPseudDir ] \ar@{-->}[d] \ar[rrr]^-{\mathrm{inddif}(\placeholder,B)_\ast} &&& [X;\Omega \mathrm{Fred}^{d,0}(L^2(\Spinor_{g_0}))] \ar[d]^{\mathrm{adjunction}}_\cong \\ 
      KO^{-(d+1)}(X) & KO^{-d}(\Sigma X) \ar[l]_-\cong && [\Sigma X;\mathrm{Fred}^{d,0}(L^2(\Spinor_{g_0}))] \ar[ll]^-{\mathrm{ind}}_-\cong }
  \end{equation*}
\end{theorem}

Our primary goal is to show that the index difference factors through the concordance set.
To this end, we will construct another model for real $K$-theory that is built out of $\InvPseudDir$ in the same way as $\ConcSet$ is built out of $\Riem^+(M)$.
We will pursue this construction in Section \ref{Section - Foundations of the Operator Concordance Set} by ``blockifying'' the objects of this section.
In the next section, we lay the necessary analytical groundwork.

\section{Block Operators}\label{Blockoperators - Section}
In this section, we discuss the foundations of \emph{block pseudo differential operators.}
These are operators on $M \times \R^n$ that can be formally thought of as pseudo differentials operators with the property that
if we think of $\R^n$ as a cube of infinite length, then they decompose into a product operator near each face.
We will discuss the required combinatorial language to deal with those operators, present a construction of a partition of unity that respects the decomposition of a block pseudo differential operator near infinity, and discuss the application of these partitions to the support of these block operators.
Finally, we define the operator suspension of a family of pseudo Dirac operator and show that it is a block operator.

We assume that the reader is familiar with the concept of pseudo differential operators and their extension of Atiyah and Singer, otherwise he or she may consult Appendix \ref{Appendix Pseudos - Chapter} and the references therein.

\subsection*{Block Bundles and Basics of Block Operators}\label{BlockBdl and Basics of BlockOp - Subsection}

The seemingly harmless space $\R^n$ can be decomposed into products $\R^p \times \R^{n-p}$ in various ways.
It will be important to keep track of these decompositions, which requires notation that seems quite heavy (and maybe unnecessary) at the beginning.
For $\dom \fateps \subseteq \mathbf{n}= \{1,\dots, n\}$ and $\fateps \colon \dom \fateps \rightarrow \Z_2$ we define, for a given $M$ the sets
 \begin{align*}
  U_R(\fateps) \deff M \times \{x \in \R^n \, : \, \fateps(i)x_i > R \text{ for } i \in \dom \fateps\}, \\
  \R^n(\fateps) \deff \mathrm{Map}(\mathbf{n} \setminus \dom \fateps,\R) \cong \{x\in \R^n \, : \, \fateps(i)x_i = R\}\\
  \R^{\fateps}_R \deff \{x\in \R^{|\dom \fateps|} \, : \, \fateps(i_j)x_j > R, \, i_j \in \dom \fateps, \, i_j < i_{j+1}\}.
 \end{align*}
 We would like to emphasise that the identification between $\mathrm{Map}(\mathbf{n} \setminus \dom \fateps,\R)$ and $\{x\in \R^n \, : \, \fateps(i)x_i = R\}$ is only canonical once $R$ is fixed.

 For such an $\fateps$ let
 \begin{equation*}
  c(\fateps)^{-1} \colon \mathbf{n} \rightarrow  \mathbf{n} \setminus \dom \fateps \sqcup \dom \fateps 
 \end{equation*}
 be the unique map that partitions $\mathbf{n}$ into $\mathbf{n} \setminus \dom \fateps$ and $\dom \fateps$ in an order preserving manner. 
 We write $\fateta \geq \fateps$ if $\dom \fateta \supseteq \dom \fateps$ and $\fateta|_{\dom \fateps} = \fateps$.
 For all $\fateta \geq \fateps$ let $c(\fateta,\fateps)^{-1}$ be the unique map that partitions $\dom \fateta^c$ into $\dom \fateta \setminus \dom \fateps$ and $\dom \fateps^c$ in an order preserving manner and let $\mathrm{per}(\fateta, \fateps)^{-1}$ be the unique map that partitions $\dom \fateta$ into $\dom \fateta - \fateps$ and $\dom \fateps$ in an order preserving manner. 
 Here, $\fateta - \fateps$ denotes $\fateta$ restricted to $\dom \fateta \setminus \dom \fateps$.
 
 For each $\fateta \geq \fateps$ we have a commutative diagram
 \begin{equation*}
  \xymatrix{ \mathbf{n} \ar[r]^-{c(\fateta)^{-1}} \ar@{=}[d] 
  & \dom \fateta^c \sqcup \dom \fateta \ar[rr]^-{\id \sqcup \per(\fateta,\fateps)^{-1}} 
  && \dom \fateta^c \sqcup \dom \fateta - \fateps \sqcup \dom \fateps 
  \\ \mathbf{n} \ar[rrr]^-{c(\fateps)^{-1}} 
  & {}
  && \dom \fateps^c \sqcup \dom \fateps \ar[u]_{c(\fateta, \fateps)^{-1} \sqcup \id}.}
 \end{equation*}
 These maps induce linear maps on Euclidean spaces by pulling back their inverse.
 For example, 
 \begin{align*}
     c(\fateps) \colon \R^n &\rightarrow \R^n(\fateps) \times \R^{\fateps} \\
     (x_1, \dots, x_n) &\mapsto (x_{c(\fateps)^{-1}(1)},\dots,x_{c(\fateps)^{-1}(|\dom \fateps^c|)};x_{c(\fateps)^{-1}(|\dom \fateps^c|+1)},\dots,x_{c(\fateps)^{-1}(n)}).
 \end{align*}
 
 We make the following notational convention:
 For a block metric $g$ on $M \times \R^n$, we denote by $g(\fateps)$ the restriction to $M \times \R^n(\fateps)$ if we interpret $\R^n(\fateps)$ as $\{x\in \R^n \, : \, \fateps(i)x_i = R\}$ and by $\face[\fateps]{ }g$ the restriction to $M \times \R^n(\fateps)$ if we interpret $\R^n(\fateps)$ as $\R^{|\dom(\fateps)^c|}$.
 The identification between these two is the pull-back with $\CubeIncl[R\fateps]{ }$.
 
 The pull-back of $c(\fateps)$ coordinates yields an isometry
 \begin{equation*}
  c(\fateps) \deff \id_M \times c(\fateps)\colon  \left(U_R(\fateps),g|_{U(\fateps)}\right) \rightarrow \bigl(M\times \R^n(\fateps) \times \R_R^\fateps , \underbrace{g(\fateps) \oplus \euclmetric}_{= c(\fateps)_\ast g|_{U_R(\fateps)}}\bigr). 
 \end{equation*}
 This fits into the following diagram of isometries:
 \begin{equation*}
  \xymatrix{ U_R(\fateta) \ar@{^{(}->}[d]\ar[r]^-{c(\fateta)} 
  & (M\times \R^n)(\fateta) \times \R_R^{\fateta} \ar[rr]^-{\id \times \per(\fateta,\fateps)}  
  && (M \times \R^n)(\fateta) \times \R_R^{\fateta - \fateps} \times \R_R^\fateps \ar[d]^{c(\fateta, \fateps)} 
  \\  U_R(\fateps) \ar[rrr]^{c(\fateps)}
  &
  && (M \times \R^n)(\fateps) \times \R_R^\fateps.}
 \end{equation*}
 
 Since all of these maps are induced by permutations on Euclidean spaces, they are Pin-structure preserving isometries because each Euclidean space has, up to equivalence, only one $\mathrm{Pin}^-$-structure.  
 
 Functoriality of the spinor bundle construction, see Lemma \ref{PinorFunctoriality - Lemma}, yields a splitting on the associated spinor bundles:
 \begin{equation*}
  \DecompMap_\fateps = \Spinor(c(\fateps)) \colon \Spinor_{U_R(\fateps)} \rightarrow \Spinor_{g(\fateps)} \boxtimes \Cliff[\fateps] = \Spinor_{g(\fateps) \oplus \euclmetric},
 \end{equation*}
 where $\Cliff[\fateps]$ denotes $\Cliff[|\fateps|]$.
 This gives rise to morphisms between spinor bundles:
 \begin{equation*}
  \xymatrix@C+1em{\Spinor|_{M \times U_R(\fateta)} \ar@{^{(}->}[d] \ar[r]^-{\DecompMap_{\fateta}} 
  & \Spinor_{g(\fateta)} \boxtimes \Cliff[\fateta] \ar[rr]^-{\id \times \Spinor(\per(\fateps,\fateta))}
  && \Spinor_{g(\fateta)} \boxtimes \Cliff[\fateta - \fateps] \boxtimes \Cliff[\fateps] \ar[d]^{\mathcal{C}_{\fateta,\fateps}} \\ 
  \Spinor|_{M \times U_R(\fateps)} \ar[rrr]^-{\DecompMap_{\fateps}}    
  &
  && \Spinor_{g(\fateps)} \boxtimes \Cliff[\fateps].}
 \end{equation*}
 
 \begin{definition}\label{BlockPseudo - Def}
  Let $g$ be a block metric on $M \times \R^n$ and let $A \in \{\id, \evenodd\}$. 
  A pseudo differential operator $P \in \PseudOpCl^m(\Spinor_g)$ is of \emph{block form} (\emph{of type} $A$) if there are $R>r>0$ such that: 
  \begin{itemize}
   \item[(1)] The block metric $g$ decomposes outside of $M \times rI^n$.
   \item[(2)] If $\sigma \in \Gamma_c(\Spinor_g)$ is supported within $M \times (RI^n)^\circ$, then
   $P\sigma$ is supported within $M \times (RI^n)^\circ$.
   \item[(3)] The operator $P$ restricts\footnote{This means that if a section is supported within $U_r(\fateps)$, then $P\sigma$ is also supported within this set.} to $M \times U_r(\fateps)$ and there is a $P(\fateps) \in \PseudOpCl^m(\Spinor_{g(\fateps)})$ such that
   \begin{equation*}
     (\DecompMap_\fateps)_\ast P|_{M \times U_r(\fateps)} = P(\fateps) \boxtimes \id + A \boxtimes D(\fateps), 
   \end{equation*}     
   where $P(\fateps) \in \PseudOpCl^m(\Spinor_{g(\fateps)})$ 
   \nomenclature{$P(\fateps)$}{Factor of a block operator at the corner point $\fateps$}
   and $D(\fateps)$ is a differential operator of order $m$ on the trivial bundle $\Cliff[\fateps] \rightarrow \R^\fateps_r$ with constant coefficients.
   
   We call $M \times RI^n$ \emph{the core} of $P$.
  \end{itemize}
 \end{definition}
 
 The constants $R$ and $r$ are not unique. 
 We will see in Lemma \ref{CoreLowerBound - Lemma} below that property (2) and (3) hold for all larger constants $R'>R$ and $r' > r$ as well.

 \begin{lemma}\label{UniquenessSubOperators - Lemma}
  Let $P \in \PseudOpCl^m(\Spinor_g)$ be a block operator of type $A$.
  Then the operators $P(\fateps)$ and $D(\fateps)$ are uniquely determined by $P$.
 \end{lemma}
 \begin{proof}
  We prove the statement by induction over $|\dom \fateps|$.
  On $\dom \fateps = \emptyset$, there is nothing to prove.
  For the induction step, it suffices to prove the following statement:
  Every operator of the form $P = Q \boxtimes \id + A \boxtimes D \in \PseudOpCl^m(\Spinor_{g_M} \boxtimes \Cliff)$ on a \emph{not necessarily closed} spin manifolds $M$ with product metrics $g = g_M \oplus \euclmetric$ uniquely deterimes $Q$ and $D$.
  We will prove this statement by expressing $Q$ and $D$ through $P$.
  
  Let $x^1, \dots, x^n$ be the standard coordinates of $\R^n$ and $\partial_j$ the derivation into the $j$-th direction.
  The operator $D$ can be written as 
  \begin{equation*}
   D = \sum_{|\fatalpha| \leq m}  a_\fatalpha \partial^\fatalpha,
  \end{equation*}
  for some constants $a_\fatalpha \in \mathrm{End}(\Cliff)$. 
  For a point $p \in \R^n$, choose smooth functions $f_1,\dots,f_n$ with compact support satisfying $f_j(x) = (x_j - p_j)$ in some neighbourhood $V$ of $p$. 
  Choose further $u \in \Gamma_c(M,\Spinor_{g_M})$ and $v \in \Gamma_c(V,\Cliff)$.
  Finally, for each multiindex $\fatalpha$ with $|\fatalpha| \leq m$, define the multiplication operator $\mu_\fatalpha$ by
  \begin{align*}
    \mu_\fatalpha \colon \Gamma_c(M \times \R^n, \Spinor_{g_M} \boxtimes \Cliff) &\rightarrow \Gamma_c(M \times \R^n, \Spinor_{g_M} \boxtimes \Cliff) \\
    \sigma &\mapsto \prod_{j = 1}^n f_j^{\alpha_j} \cdot \sigma.
  \end{align*}   
  The Leibniz rule implies
  \begin{equation*}
    \lbrack P, \mu_\fatalpha\rbrack(u \otimes v)(m,p) = A(u)(m) \otimes a_\fatalpha v(p),
  \end{equation*}
  so $D$ is uniquely determined by $P$.
  
  It remains to determine $Q$ from $P$ and $D$. 
  Since $P - A \boxtimes D$ is $\mathcal{C}^\infty(\R^n)$-linear on $\Gamma_c(M \times \R^n, \Spinor_{g_M} \boxtimes \Cliff)$ we can define a linear map
  \begin{equation*}
   \tilde{Q} \colon \Gamma_c(M , \Spinor_{g_M}) \rightarrow \Gamma(M, \Spinor_{g_M})
  \end{equation*}   
 as follows: Fix a point $p \in \R^n$ and a function $f \in \mathcal{C}^\infty(\R^n,\R)$ that is identically $1$ near $x$ and compactly supported. Set then
 \begin{equation*}
   \tilde{Q}(\sigma)(m) \deff (P - A \boxtimes D)(\sigma \otimes f \cdot 1)(m,p).
 \end{equation*}
 Plugging in the decomposition gives
 \begin{align*}
  \tilde{Q}(\sigma)(m) &= (P - A \boxtimes D)(\sigma)(m,p)\\ 
                       &= Q\boxtimes \id (\sigma \otimes f \cdot 1)(m,p)\\
                       &= Q(\sigma)(m) \otimes f(p) = Q(\sigma)(m) \otimes 1, 
 \end{align*}
 which can be identified with $Q(\sigma)$ under the inclusion $\Spinor_{g_M} \subseteq \Spinor_{g_M} \boxtimes \Cliff$. 
 \end{proof}
 
 \begin{cor}\label{BlockDescend - Cor}
  If $P$ is a block operator on $\Spinor_g$, then $P(\fateps)$ is a block operator on $\Spinor_{g(\fateps)}$. 
  Equivalently, if $\CubeIncl[2\rho\fateps]{ } \colon \R^n(\fateps) \xrightarrow{\cong} \{x\in \R^n \, : \, \fateps(i)x_i = 2\rho\}$ denotes the isometric inclusion that plugs $2\fateps(j)\rho$ into the $j$-th coordinate for all $j \in \dom \fateps$, then $\CubeIncl[2\rho\fateps]{ }^\ast P(\fateps) = (\CubeIncl[2\rho\fateps]{ })^{-1}_\ast P(\fateps)$ is a block operator on $\Spinor_{\face[\fateps]{ }g}$ for all $\rho>R$, where $R$ denotes the size of the core of $P$. 
  These block operators satisfy conditions (1)-(3) for the same constants as $P$.
 \end{cor}
 \begin{proof}
  Since $g$ decomposes outside of $rI^n$, it also decomposes outside of $\rho I^n$.
  If we identify $\R^n(\fateps)$ with the affine subspace $\{x \in \R^n \, : \, \fateps(j)x_j = \rho\}$, then the restricted metric $g(\fateps) = g\restrict_{\R^n(\fateps)}$ is again a block metric that decomposes outside of $rI^n \cap \R^n\fateps) \cong r I^{n - |\dom \fateps|}$.

 Set $m \deff n-|\dom \fateps|$.
 We verify axiom $(2)$ for ${\CubeIncl[2r\fateps]{ }}^\ast P(\fateps)$ by contra-position.
 Assume there is a section $\sigma \in \Gamma_c(\Spinor_{\face[]{\fateps}g(\fateps)})$ with support in $M \times (RI^m)^\circ$ such that ${\CubeIncl[2r\fateps]{ }}^\ast P(\fateps)(\sigma)$ is not supported within $M \times (RI^m)^\circ$.
 Pick a smooth function $f \colon \R \rightarrow [0,1]$ that is supported within $(r,R)$ and define on $R_r^\fateps$ the real valued function 
 \begin{equation*}
     f^\fateps(x) = \prod_{j=1}^{|\dom \fateps|} f(x_j).
 \end{equation*}
 Then the section $\sigma \otimes f^\fateps$ is supported within the interior of $M \times RI^m \times RI^{|\dom \fateps|}$ but, by assumption, the support of $\bigl(P(\fateps) \boxtimes \id + A \boxtimes D\bigr)(\sigma \otimes f^\fateps)$ is not supported within $M \times RI^m \times RI^{|\dom \fateps|}$.
 Thus $\tau \deff {\Phi_\fateps}^{-1}_\ast(\sigma \otimes f^\fateps)$ is a section of $\Spinor_g$ that is supported within the interior of $M \times RI^n$.
 But the support of $P(\tau)$ does not lie in $M \times RI^n$, contradicting axiom $(2)$ of $P$.
 
  We first give an informal proof for (3).
  We know that
  \begin{equation*}
      \Phi_{\CubeIncl[\fateps]{ }(\fateta)\, \ast} P|_{U_r(\CubeIncl[\fateps]{ }(\fateta))} = P\bigl(\CubeIncl[\fateps]{ }(\fateta))) \boxtimes \id + \evenodd \boxtimes D(\CubeIncl[\fateps]{ }(\fateta))\bigr).
  \end{equation*}
  If we permute the (fixed) coordinates parametrised by $\dom \, \fateps$ ``to the right'', using $\mathrm{per}(\CubeIncl[\fateps]{ }(\fateta)),\fateps)$ we can identify $\R^{\CubeIncl[\fateps]{ }(\fateta))}_r$ with $\R^{\CubeIncl[\fateps]{ }(\fateta)) - \fateps}_r \times \R^{\fateps}_r$.
  The uniqueness statement of Lemma \ref{UniquenessSubOperators - Lemma} implies (with abuse of notation)
  \begin{equation*}
      P(\fateps) = P\bigl(\CubeIncl[\fateps]{ }(\fateta)\bigr) \otimes \id + \mathrm{per}(\CubeIncl[\fateps]{ }(\fateta),\fateps)_\ast D\bigl(\CubeIncl[\fateps]{ }(\fateta)\bigr).
  \end{equation*}
  Now the formal proof of (3): The inclusion $\CubeIncl[2r\fateps]{ }$ maps $U^{m}_r(\fateta)$ into $U^n_r(\CubeIncl[\fateps]{ }(\fateta))$. 
  Furthermore, the following diagram commutes:
  \begin{equation*}
   \xymatrix{ M \times U^{m}_r(\fateta) \ar[rr]^-{c^m(\fateta)} \ar@{^{(}->}[d]^{\CubeIncl[2r\fateps]{ }} && M \times \R^{m}(\fateta) \times \R^{\fateta}_r \ar@/^2pc/[ddl]^{ \id \times \CubeIncl[\fateps]{ } } \\
   M \times U^n_r(\CubeIncl[\fateps]{ }(\fateta)) \ar[r]_-{c(\CubeIncl[\fateps]{ }(\fateta))} &
   M \times \R^n(\CubeIncl[\fateps]{ }(\fateta)) \times \R_r^{\CubeIncl[\fateps]{ }(\fateta)} \ar[d]^{\mathrm{per}(\CubeIncl[\fateps]{ }(\fateta),\fateps)} &  \\
   & M \times \R^{n}(\CubeIncl[\fateps]{ }(\fateta)) \times \R^{\CubeIncl[\fateps]{ }(\fateta) - \fateps}_r \times \R_r^{\fateps} &} 
  \end{equation*}    
  This implies  
  \begin{equation*}
   (\Phi_\fateta)_\ast ({\CubeIncl[2r\fateps]{ }}^\ast P(\fateps)) ) = (\id \times \CubeIncl[\fateps]{ })^\ast \bigl(P(\CubeIncl[\fateps]{ }(\fateta)) \otimes \id + A \otimes \per(\CubeIncl[\fateps]{ }(\fateta), \fateps)_\ast D(\CubeIncl[\fateps]{ }(\fateta))\bigr),   
  \end{equation*}
  so ${\CubeIncl[2r\fateps]{ }}^\ast P(\fateps)$ is indeed of block form.
 \end{proof}
 
 The converse of this lemma is also true.

 \begin{lemma}\label{BlockSuspends - Lemma}
  Let $m>0$ and $P \in \PseudOpCl^m(\Spinor_g)$ be of block form of type $A$.
  For each differential operator $D$ of order $m$ with constant coefficients on $\mathcal{C}_c^\infty(\R,\Cliff[1][0])$, the operator
  \begin{equation*}
      Q \deff P \boxtimes \id + A \boxtimes D \in \PseudOpCl^m(\Spinor_{g} \boxtimes \Cliff[1][0]) = \PseudOpCl^m(\Spinor_{g \oplus \diff x_{n+1}^2})
  \end{equation*}
  is also block operator of type $A$.
 \end{lemma}
 \begin{proof}
  If $m >0$, then $Q$ is an element of $\PseudOpCl^m(\Spinor_{g \oplus \diff x_{n+1}^2})$ by Theorem \ref{ASTensor - Theorem}.
  
  Axiom (1) is trivially satisfied.
  Axiom (2) also holds because $A \boxtimes D$ acts as a differential operator in the additional dimension.
  Axiom (3) requires some definitions. 
  Let $\fateps \colon \mathbf{n+1} \rightarrow \Z_2$ be given.
  
  If $n+1 \in \dom \fateps$, then $U_r(\fateps) = U_r(\CubeProj{n+1}(\fateps)) \times \{\fateps(n+1) x_{n+1} > r\}$ and $\Phi_{\fateps} = \Phi_{\CubeProj{n+1}(\fateps)} \times \id$.
  In this case, we have
  \begin{equation*}
      P(\fateps) = Q(\CubeProj{n+1}(\fateps)) \boxtimes \id_{\R} \text{ and } D^Q(\fateps) = D^P(\CubeProj{n+1}(\fateps)) \boxtimes \id.
  \end{equation*}
  
  If $n+1 \notin \dom \fateps$, then $U_r(\fateps) = U_r^n(\fateps) \times \R$, where $U_r^n(\fateps) \subseteq \R^n$, and $\Phi_\fateps$ is given by the decomposition
  \begin{equation*}
      \xymatrix{ \Spinor|_{M \times U_r(\fateps)} = \Spinor\restrict_{M \times U_r^n(\fateps)}\boxtimes \Cliff[1][0] \ar[rr]^-{\Phi_\fateps \boxtimes \id} && \Spinor_{g(\fateps)} \boxtimes \Cliff[\fateps][0] \boxtimes \Cliff[1][0] \ar[d]^{\mathrm{swap}}_{\id \boxtimes} \\
      && \Spinor_{g(\fateps)}\boxtimes \Cliff[1][0] \boxtimes \Cliff[\fateps][0]. }
  \end{equation*}
  In this case, we have
  \begin{equation*}
      Q(\fateps) = P(\fateps) \boxtimes \id + A \boxtimes D \quad \text{ and } \quad D^Q(\fateps) = D^P(\fateps).
  \end{equation*}
  \todo[inline]{Figure out if $\Cycl(1:1+|\dom \fateps|$ or $\Cycl(1:1+|\dom \fateps|^{-1}$. Define the correct map in an ad-hoc manner as the conceptual notation will be used later.}
  
  With these definitions the equation 
  \begin{equation*}
      {\Phi_\fateps}_\ast Q|_{M \times U_r(\fateps)} = Q(\fateps) \boxtimes \id + A \boxtimes D^Q(\fateps)
  \end{equation*}
  is easily verified so that $Q$ is a block operator.
 \end{proof}
 
 Pseudo differential operators are not local by nature. 
 Block operators on the other hand are local in certain directions on special classes of open sets.
 In order to use local-to-global principles for block operators, we need partition of unities that are adapted to the specific block form of a block operator.
 
 \begin{definition}\label{BlockPartOfUnity - Def}
  Let $\chi \colon \R \rightarrow [0,1]$ be a compactly supported, smooth, and even function that is identically $1$ near the origin.
  We denote the extensions by zero of 
  \begin{equation*}
      \chi^{-1} \deff (1 - \chi)|_{\R_{\leq 0}}, \quad \chi^0 \deff \chi, \quad \text{ and } \quad \chi^1 \deff (1 - \chi)|_{\R_{\geq 0}} 
  \end{equation*}
  with the same symbols.
  More generally, we define on $\R^n$ or $M \times \R^n$ the partition of unity
  \begin{equation*}
      \biggl\{\chi^\fatalpha(m,x) \deff \prod_{j=1}^n \chi^{\fatalpha(j)}(x_j)\biggr\}.
  \end{equation*}
  parameterised ${\alpha \colon \mathbf{n}\rightarrow \Z_3}$.
  \nomenclature{$\chi^\fatalpha$}{Partition of unity adapted to the block form}
 \end{definition}
 
 As an application of this partition, we prove that the constants in the definition  of a block operator are actually lower bounds.
 
 \begin{lemma}\label{CoreLowerBound - Lemma}
  Let $P \in \PseudOpCl^m(\Spinor_g)$ be a block operator and let $R > r$ be constants such that the conditions (1)-(3) of $P$ hold.
  Then these conditions also hold for all $R' > R$ and $r'>r$ with $R' > r'$, 
 \end{lemma}
 \begin{proof}
  Condition (1) is clearly satisfied for all $r'>r$.
  
  We continue with the proof of condition (3). 
  Let $\sigma \in \Gamma_c(\Spinor_g)$ be a section with $\supp \, \sigma \subseteq U_{r'}(\fateps)$.
  Since the support is compact, we find a slightly bigger $r'' > r'$ such that $\supp \, \sigma \subseteq U_{r''}(\fateps)$.
  Let $\chi \colon \R \rightarrow [0,1]$ be a smooth function that is identically $1$ near $\R_{\geq r''}$ and identically zero near $\R_{\leq r'}$.
  
  For $\fateps \colon \dom \fateps \subseteq \mathbf{n} \rightarrow \Z_2$ define
  \begin{equation*}
      \chi^\fateps(x) \deff \prod_{j \in \dom \fateps} \chi(\fateps(j)x_j).
  \end{equation*}
  This function is supported within $U_{r'}(\fateps)$.
  If we do not notationally distinguish $\chi^\fateps$ and $(\Phi_\fateps)_\ast \chi^\fateps = \chi^\fateps \circ \Phi_\fateps^{-1}$, we have
  \begin{align*}
      (\Phi_\fateps)_\ast(P\chi^\fateps \sigma) &= P(\fateps)\boxtimes \id(\chi^\fateps\sigma) + A \boxtimes D(\chi^\fateps \sigma) \\
      &= \chi^\fateps \cdot (P(\fateps)\boxtimes \id)(\sigma) + D(\chi^\fateps) \cdot (A \boxtimes \id)(\sigma) + \chi^\fateps \cdot (A \boxtimes D)(\sigma),
  \end{align*}
  which is supported within $U_{r'}(\fateps)$, too. 
  Thus, $P$ also restricts to $U_{r'}(\fateps)$ and necessarily decomposes there as in (3).
  
  Condition (2) uses the construction of Definition \ref{BlockPartOfUnity - Def}.
  We will prove condition (2) inductively over $n$, the dimension of the Euclidean factor of $M \times \R^n$.
  For $n=0$ there is nothing to prove, and we can continuue with $n \geq 1$.
  Let $\chi \colon \R \rightarrow [0,1]$ be a smooth and even function that is identically $1$ near $[-r,r]$ and that vanishes near $\R\setminus(-R,R)$.
  Let $\{\chi^\fatalpha\}_{\fatalpha \in \Z_3^n}$ be the corresponding partition of unity on $M \times \R^n$ from Definition \ref{BlockPartOfUnity - Def}.
  We clearly have
  \begin{equation*}
      P(\sigma) = P(\chi^\mathbf{0} \sigma) + \sum_{\fatalpha \in \Z_3^n \setminus \mathbf{0}} P(\chi^\fatalpha \sigma).
  \end{equation*}
  The summand $\chi^{\mathbf{0}}\sigma$ is supported within $M \times (RI^n)^\circ$, so, by condition (2) of $P$, also $P(\chi^\mathbf{0} \sigma)$ is supported within $M \times (RI^n)^\circ$.
  
  Recall the notation $\hat{\fatalpha} \deff \fatalpha|_{\supp \fatalpha} \colon \supp \, \fatalpha \rightarrow \Z_2$.
  The section $\chi^\fatalpha \sigma$ is supported within $U_{r}(\hat{\fatalpha})$.
  If $\fatalpha \neq \mathbf{0}$, then $\hat{\fatalpha}$ is not the empty map so that $P(\chi^\fatalpha \sigma)$ is also supported within $U_r(\hat{\fatalpha})$ by the proof of condition (3) above.
  
  Since $\supp \sigma \subseteq M \times (R'I^n)^\circ$ we can find a smooth and even function $\psi \colon \R \rightarrow [0, 1]$ that vanishes near $\R \setminus (-R',R')$ and such that $\psi^\mathbf{0}(x) \deff \prod_{j=1}^n \psi(x_j)$ is identically $1$ on $\supp \, \sigma$.
  The block form of $P$ implies
  \begin{align*}
      P(\chi^\mathbf{\fatalpha}\sigma) &= P(\chi^{\mathbf{\fatalpha}}\psi^\mathbf{0}\sigma) \\
      \begin{split}&= \prod_{j \in \supp \fatalpha} \psi \circ \pr_j \cdot P\left(\Bigl(\prod_{j \in \mathrm{Null}(\fatalpha)} \psi\circ \pr_j\Bigr)\cdot \chi^\fatalpha \sigma\right)\\
      & \qquad \quad + D(\hat{\fatalpha})\left(\prod_{j\in \supp \fatalpha} \psi \circ \pr_j\right) \cdot P\left(\Bigl(\prod_{j \in \mathrm{Null}\, \fatalpha} \psi\circ \pr_j\Bigr) \cdot \chi^\fatalpha \sigma\right).
      \end{split} 
  \end{align*}
  By Corollary \ref{BlockDescend - Cor}, $P(\fatalpha|_{\supp \fatalpha})$ is a block operator on $\Spinor_{g(\fatalpha|_{\supp \fatalpha})}$ that satisfy conditions (1) - (3) for the same constants as $P$ and $A \boxtimes D(\fatalpha|_{\supp \fatalpha})$ is a differential operator and hence only support-decreasing.
  This and the induction assumption, implies that the section 
  \begin{equation*}
      P\left(\Bigl(\prod_{j \in \mathrm{Null}(\fatalpha)} \psi\circ \pr_j\Bigr) \cdot \chi^\fatalpha \sigma\right)
  \end{equation*}
  must be supported within $\{|x_j| < R' \, : \, j \in \mathrm{Null}\, \fatalpha\} \cap U_r(\fatalpha|_{\supp \fatalpha})$.
  It follows that 
  \begin{equation*}
     \prod_{j \in \supp \fatalpha} \psi \circ \pr_j \cdot P\left(\prod_{j \in \mathrm{Null}(\fatalpha)} \psi\circ \pr_j \chi^\fatalpha \sigma\right) 
  \end{equation*}
  and 
  \begin{equation*}
      D(\hat{\fatalpha})\left(\prod_{j\in \supp \fatalpha} \psi \circ \pr_j\right) \cdot P\left(\prod_{j \in \mathrm{Null}\, \fatalpha} \psi\circ \pr_j \chi^\fatalpha \sigma\right)
  \end{equation*}
  are supported within $M \times (R'I^n)^\circ \cap U_r(\fatalpha|_{\supp \, \fatalpha})$, which finishes the proof of (2).
 \end{proof}
 
 \begin{cor}\label{InvariantDecomposition - Corollary}
   Let $P\in \PseudOpCl^1(\Spinor_g)$ be a block operator of type $A$ whose core is contained in $M \times RI^n$.
   For $\fatalpha \colon \mathbf{n}\rightarrow \Z_3$ and $\rho'>\rho>R$ define
   \begin{equation*}
     V_{\rho<\rho'}(\fatalpha) \deff M \times \{x \in \R^n \, : \, x_j \in (-\rho',\rho') \text{ if } j \in \mathrm{Null}(\fatalpha), \alpha_jx_j>\rho \text{ otherwise} \}.  
   \end{equation*}
   If $\sigma \in \Gamma_c(\Spinor_g)$ is supported within $V_{\rho<\rho'}(\fatalpha)$, then $P\sigma$ is also supported there.
 \end{cor}
 \begin{proof}
  If $\fatalpha = \mathbf{0}$, then $V_{\rho<\rho'}(0) = M \times \rho'I^n$ and the claim follows from Lemma \ref{CoreLowerBound - Lemma}.
  
  Assume now that $\fatalpha \neq \mathbf{0}$, so that $V_{\rho<\rho'}(\fatalpha) \subseteq U_\rho(\hat{\fatalpha})$.
  For each $u \in \Gamma_c(\Spinor_g)$ supported with in $V_{\rho<\rho'}(\fatalpha)$ we deduce with abuse of notation 
  \begin{align*}
     Pu = (P(\hat{\fatalpha})\boxtimes \id)(\chi^\fatalpha u) + \evenodd \boxtimes D(\hat{\fatalpha})(\chi^\fatalpha)(\chi^\fatalpha u)
  \end{align*}
  The second summand is a differential operator and hence support-decreasing.
  
  We know from Corollary \ref{BlockDescend - Cor} that $P(\hat{\fatalpha})$ is again a block operator whose core is contained in $M \times \rho' I^{|\mathrm{Null}(\fatalpha)|}$.
  Since $V_{\rho<\rho'}(\fatalpha) \cong \rho'I^{|\mathrm{Null}(\fatalpha)|} \times \R_\rho^{\hat{\fatalpha}}$, it follows from the previous Lemma that 
  \begin{equation*}
      \supp\, P(\fatalpha)(u(\placeholder,y)) \subseteq \rho'I^{|\mathrm{Null}(\fatalpha)|}
  \end{equation*}
  and hence
  \begin{equation*}
      \supp\, P(\fatalpha)( u) \subseteq \rho'I^{|\mathrm{Null}(\fatalpha)|} \times \R_\rho^{\hat{\fatalpha}} = V_{\rho<\rho'}(\fatalpha).
  \end{equation*}
 \end{proof}
 
 The technique used in the previous proofs can be used to prove that block operators are local in the following limited sense.
 
 \begin{cor}\label{BlockOpPartLocal - Cor}
   Let $P$ be a block operator whose core is contained in $M \times \rho I^n$.
   If $\rho' > \rho$ and $V \subseteq U_{\rho'}(\fateps)$ is an open set that corresponds under $U_{\rho'}(\fateps) \cong M \times \R^n(\fateps) \times \R^{\fateps}_{\rho'}$ to $M \times \R^n(\fateps) \times V^\fateps$, then $(Pu)|_V$ only depends on (the germ of) $u|_{\Bar{V}}$.
 \end{cor}
 \begin{proof}
  It suffices to show that if $u$ vanishes identically near $\Bar{V}$, then $Pu$ vanishes identically on $V$.
  Let $\phi \colon M \times \R^n \rightarrow \R_{\geq 0}$ be a smooth function such that $\phi^{-1}(\R_{>0}) = V$ and such that $\phi$ only depends on variables indexed by $\dom(\fateps)$.
  Let $\chi \colon \R \rightarrow [0,1]$ be a smooth even function such that $\chi$ is identically $1$ on $[-\rho,\rho]$ and whose support lies in $(-\rho',\rho')$. 
  Denote by $\chi^\fatalpha$ the partition of unity from Definition \ref{BlockPartOfUnity - Def}.
  
  Since $\supp \phi \subseteq U_{\rho'}(\fateps)$, we conclude that $\supp \, \phi$ and $V_{\rho<\rho'}(\hat{\fatalpha})$ have non-empty intersection only if $\hat{\fatalpha} \geq \fateps$.
  We know from Corollary \ref{InvariantDecomposition - Corollary} that $\supp(P\chi^\fatalpha u) \subseteq V_{\rho<\rho'}(\fatalpha)$, so we conclude
  $\phi \cdot P \chi^\fatalpha u = 0$ if $\hat{\fatalpha} \geq \fateps$ is false.
  On the other hand, if $\hat{\fatalpha} \geq \fateps$, then $V_{\rho<\rho'}(\fatalpha) \subseteq U_{\rho}(\fateps)$ and the block decomposition of $P$ implies
 \begin{align*}
     \phi \cdot P \chi^\fatalpha u &= P \phi \chi^\fatalpha u + [\phi,P]\chi^\fatalpha u = P\chi^\fatalpha\underbrace{\phi u}_{=0} - i^{-1}\underbrace{\symb_1(P)(\placeholder, \diff \phi)\chi^\fatalpha u}_{=0 \text{ on }V} \\ 
     &= 0 \qquad \text{ on } V
 \end{align*}
 because $u|_V = 0$.
 Since $\{\chi^\fatalpha\}_{\fatalpha \in \Z_3^n}$ is a partition of unity, the previous calculations imply $\phi P u|_V = 0$ from which we deduce $Pu |_V =0$.
 \end{proof}
 
 \begin{definition}
  A block operator $P \in \PseudOpCl^1(\Spinor_g)$ is called \emph{block operator of Dirac-type}, if $A = \evenodd$ and $D(\fateps) = \DiracSt_{\R^\fateps}$ 
  \nomenclature{$\DiracSt$}{Dirac operator of the euclidean metric}
  is the Dirac operator of the Euclidean metric on $\R^\fateps$ for all $\fateps \colon \dom \fateps \rightarrow \Z_2$, and if the principal symbol satisfies $\symb_1(P) = \symb_1(\Dirac_g)$. 
 \end{definition}
 
 The motivating examples for this definition are the Dirac operator of a block metric and the adjoint operator of a block map of pseudo Dirac operators.
 \begin{example}
  We apply the discussion of Example \ref{Spinor bundle product manifold - Example} to show that $\Dirac_g$ is a block operator (of Dirac type), if $g$ is a block metric on $M \times \R^n$.
  Assume that the block metric decomposes outside of $M \times rI^n$. 
  Then, since differential operators are only support decreasing, we can choose any $R>r$ so that condition (1) and (2) are trivially satisfied.
  It remains to show that $\Dirac_g$ decomposes appropriately on $U_r(\fateps)$.
  Since the defining sequence of the Dirac operator (\ref{Eq: Def Dirac Operator}) is natural with respect to Pin-structure preserving, isometric diffeomorphisms, we have
  \begin{align*}
      (\Phi_{\fateps})_\ast \Dirac_g|_{M \times U_r(\fateps)} &= \Dirac_{c(\fateps)_\ast g} = \Dirac_{g(\fateps) + \euclmetric_{\R^{|\dom \fateps|}}} = \Dirac_{g(\fateps)} \boxtimes \id + \evenodd \boxtimes \DiracSt_{\R^{|\dom \fateps|}}.
  \end{align*}
 \end{example}
 
 The decomposition of the spinor bundle in Example \ref{Spinor bundle product manifold - Example} has a weaker version.
 If $U\subseteq \R^n$ is open and $g \colon U \rightarrow \Riem(M)$ is a smooth map, then $TU$ is still an orthogonal complement of $TM$ in $T(M \times U)$ with respect to $\susp(g)$.
 This induces a splitting of the spinor bundle $\Spinor_{\susp(g)} \cong \Spinor_g \otimes \Cliff[n]$.
 
 \begin{definition}
  Let $P \colon U \rightarrow \PseudOp^m(\Spinor)$ be a smooth map, where the target carries the (subspace topology of the) amplitude topology. 
  It comes with a smooth map $g \colon U \rightarrow \Riem(M)$.
  Define on $\Spinor_{\susp(g)}$ the operators that are given under the isomorphism $\Spinor_{\susp(g)} \cong \Spinor_g \otimes \Cliff[n]$ by
 \begin{align*}
    P^\extension(\sigma\otimes e_I)(m,t) \deff \left(P_t\sigma(\cdot,t)\right)   \otimes e_I, \nomenclature{$P^\extension$}{tensorial extension of $P$}\\
    \susp(P) \deff P^\extension  + \sum_{k=1}^n \partial_{t_k} \cdot  \nabla_{\partial_{t_k}}^{\Spinor_{\susp \, g}}.
  \end{align*} 
  We call $\susp(P)$ 
  \nomenclature{$\susp(P)$}{(operator) suspension of the block map $P \colon \R^n \rightarrow \PseudDir$}
  the \emph{operator suspension of} $P$. 
 \end{definition}
 
 By the parametrised version of the Atiyah-Singer exterior tensor theorem \ref{FamilyASTensor - Theorem} we have $\susp(P) \in \PseudOpCl^m(\Spinor_{\susp(g)})$ if $m>0$.
 
 We are mostly interested in the case where $U = \R^n$ and $P$ is a smooth block map. 
 Note that this implies that the underlying map of Riemannian metrics is also a smooth block map.
 However, in some proofs, we need the operator suspension on more general open subsets.
 \begin{lemma}\label{SuspensionBlockOperator  - Lemma}
  Let $m>0$ and let $P \colon \R^n \rightarrow \PseudOp^m(\Spinor)$ be a smooth block map, where the target carries the (subspace topology of the) amplitude topology and underlying block map $g$.
  
  Then $P^\extension$ and  $\susp(P)$ are a block operators of type $\id$ and $\evenodd$, respectively.
  If $P$ takes values in $\PseudDir$, then $\susp(P)$ is a block operator of Dirac type.
 \end{lemma}
 \begin{proof}
  We will only prove the statement for $\susp(P)$ in the case that $P$ takes values in $\PseudDir$, for the other cases are similar.
  Since $P$ and $g$ are block maps, there is an $r>0$ such that $P$ and $g$ are independent of $x_i$ if $|x_i|>r$.
  Then $\susp(g)$ decomposes outside of $M \times rI^n$.
  Since $\susp(P) - \evenodd \otimes \DiracSt_{\R^n}$ is $\mathcal{C}^\infty(\R^n)$-linear, the operator restricts to open subsets of the form $M \times U$ with $U \subseteq \R^n$ open, and thus $P$ does it, too.
  
  It remains to show that $P$ decomposes appropriately on $U_r(\fateps)$.
  On this subset, a permutation of coordinates $c(\fateps) \colon M \times U_r(\fateps) \rightarrow M \times \R^n(\fateps) \times \R^\fateps_r$ is an isometry if the domain carries the metric $g|_{U_r(\fateps)}$ and the right hand side has $\susp(g)(\fateps) \oplus \euclmetric$.
  Thus the spinor bundle decomposes into 
  \begin{equation*}
     \Spinor_{\susp(g)(\fateps)} \boxtimes \Cliff[\fateps] \cong \Spinor_g \otimes \Cliff[n - |\fateps|] \boxtimes \Cliff[\fateps] 
  \end{equation*}
  over $(M\times \R^n(\fateps)) \times \R^\fateps_r$ via $\Phi_\fateps = \Spinor(c(\fateps))$.
  Since $P|_{U_r(\fateps)}$ is independent of the coordinates indexed by $\dom \fateps$, we derive the equation
  \begin{align*}
      (\Phi_\fateps)_\ast (\susp(P)|_{M \times U_r(\fateps)}) &= (\Phi_\fateps)_\ast \bigl(P \otimes \id\bigr) + (\Phi_\fateps)_\ast \bigl(\evenodd_{\Spinor_g} \otimes \DiracSt_{\R^n}\bigr) \\ 
      \begin{split}&= P_{c(\fateps)^{-1}(\placeholder)} \otimes \id + (\evenodd_{\Spinor_g} \otimes \DiracSt_{\R^n(\fateps)})\boxtimes \id_{\Cliff[\fateps]} \\
      & \qquad + \evenodd_{\Spinor_g} \otimes \evenodd_{\Cliff[n-|\fateps|][]} \boxtimes \DiracSt_{\R^\fateps_r}
      \end{split}\\
      \begin{split}&= \Bigl((P|_{M \times \R^n(\fateps)})_{c(\fateps)^{-1}(\placeholder)} \otimes \id_{\Cliff[n-|\fateps|]} + \evenodd_{\Spinor_g} \otimes \DiracSt_{\R^n(\fateps)} \Bigr) \boxtimes \id \\
      & \qquad+ \evenodd_{\Spinor_{\susp(g)(\fateps)}} \boxtimes \DiracSt_{\R^\fateps},
      \end{split}
  \end{align*}
  which shows that $\susp(P)$ decomposes as required.
 \end{proof}
 
 So far, we have only studied block operators on the level of vector bundles over manifolds.
 Before we are going to construct the block version of $\InvPseudDir$, we will  study its (functional) analytic properties in the next section.

 \section{Analytical Properties of Block Operators}\label{Analytical Properties of Block Operators - Section}

The purpose of this section is to derive the needed analytical properties of block operators.
We will show that each block operator realises to a bounded operator on positive Sobolev spaces and that the realisations of block operators of Dirac type satisfy analogous versions of the fundamental elliptic estimate.
Moreover, we will show that if they are self-adjoint and invertible near infinity, then they are also Fredholm operators.
Many of these properties rely on the structure of a block metric near infinity because we are working on non-compact manifolds.

In the following, we only consider pseudo differential operators of order $m \in \{0,1\}$.
\begin{lemma}\label{BlcokOperatorsInduceBound - Lemma}
 Let $g$ be a block metric on $M \times \R^n$ and let $P \in \PseudOpCl^m(\Spinor_g)$ be a block operator.
 Then $P$ extends to a bounded linear operator
 \begin{equation*}
     P \colon H^s(\Spinor_g) \rightarrow H^{s-m}(\Spinor_g)
 \end{equation*}
 for all $s \in \Z$.
\end{lemma}
\begin{proof}
 We prove this lemma by induction over $n$.
 For $n=0$, this is classical because $M$ is closed, see Theorem \ref{PseudosInduceBndMapsSobolev - Theorem}.
 For $n>0$, we pick appropriate $R > r >0$ such that $P$ restricts to $M \times (RI^n)^\circ$ and decomposes on $\{\eps x_i > r\}$ for all $i \in \{1,\dots,n\}$ and all $\eps \in \Z_2$.
 Let $\chi \colon \R \rightarrow [0,1]$ be a smooth, symmetric function that is identically $1$ on $(-r,r)$ and vanishes near $\R \setminus (-R,R)$.
 Let $\{\chi^\fatalpha\}_{\fatalpha \in \Z_3^n}$ be the partition of unity obtained from Definition \ref{BlockPartOfUnity - Def}.
 By the triangle inequality, it suffices to show that $P\cdot \chi^\fatalpha$ extends to a bounded operator.
 
 If $\fatalpha = \mathbf{0}$, then $P\chi^\mathbf{0}$ restricts to an operator on the relatively compact set $M \times (RI^n)^\circ$.
 Thus, by Theorem \ref{PseudosInduceBndMapsSobolev - Theorem}, $P \chi^\mathbf{0}$ extends to a bounded operator $H^s(\Spinor_g) \rightarrow H^{s-1}(\Spinor_g)$.
 
 If $\fatalpha \neq \mathbf{0}$, then there is a tupel $(i,\eps) \in \{1,\dots,n\} \times \Z_2$ such that $\supp \chi^\fatalpha$ is contained in $U_r((i,\eps))$.
 The induction hypothesis and \cite{palais1965seminar}*{Corollary XIV.1} applied to the discrete Hilbert chains $H^s(M\times \R^{n-1};\Spinor_{g(i,\eps)})$ and $H^s(\R;\Cliff[1][0])$ imply 
 \begin{equation*}
    P(i,\eps) \boxtimes \id + A \boxtimes D(i,\eps) \in Op^1(H^s(\Spinor_{g(i,\eps)} \boxtimes \Cliff[1][0])), 
 \end{equation*}
 which means that this operator induces a bounded linear map
 \begin{equation*}
     H^s(\Spinor_{g(i,\eps)} \boxtimes \Cliff[1][0]) \rightarrow H^{s-1}(\Spinor_{g(i,\eps)} \boxtimes \Cliff[1][0])
 \end{equation*}
 for all $s \in \Z$.
 Since $P$ decomposes into $P(i,\eps) \boxtimes \id + A \boxtimes D(i,\eps)$ on $\Spinor_g|_{U_r((i,\eps))} \cong \Spinor_{g(i,\eps)} \boxtimes \Cliff[1][0]|_{\eps \R_{>r}}$, it follows $P\chi^\fatalpha \in Op^1(\Spinor_{g})$. 
\end{proof}

Recall that a sequence of operators $P_n$ converges to $P$ in the Atiyah-Singer topology if the $P_i$ converge uniformly on compact subsets to $P$.
Thus, it is not clear that the fundamental ellitpic estimate carries over to elements in the Atiyah-Singer closure.
However, this is the case, if the operator $P$ increases the support in a controlled manner.

\begin{lemma}\label{EllipticEstimateCpt - lemma} 
   Let $M$ be a not necessarily closed manifold and $E,F \rightarrow M$ be vector bundles over $M$.
   Let $P \in \PseudOpCl^1(E,F)$ be an elliptic pseudo differential operator whose principal symbol is the principal symbol of an actual pseudo differential operator (i.e. smooth) and that has the following property: For all compact neighbourhoods $K \subseteq M$ there is a relatively compact open subset $W$ such that if $\supp u \subseteq K$ then $\supp Pu \subseteq W$.
   
   Then, for all compact $K \subseteq M$, there is a compact $L \subseteq M$ such that $P$ induces a bounded map $\Sobolev[K]{s}{E} \rightarrow \Sobolev[L]{s-1}{F}$, where $H_K^s(E)$ denotes the subspace of all section with support in $K$. 
   The operator $P$ furthermore satisfies the elliptic estimate:
   \begin{equation*}
    ||u||_s \leq C\left( ||u||_{s-1} + ||Pu||_{s-1}\right).
   \end{equation*}
 \end{lemma}
\begin{proof}
 The subspace $\Gamma_c(E)$ is dense in $H^s(E)$ for all $s \in \Z$ and the multiplication with a smooth compactly supported function induces a continuous endomorphism on all Sobolev spaces.
 Thus, if $\chi$ is a smooth compactly supported function that is identically $1$ on $K$, then we can approximate each element of $H_K^s(E)$ by an element of $\chi \cdot \Gamma_c(E)$.
 If $L$ is the compact subset from the assumption assigned to $\supp \, \chi$, then a simple approximation argument shows that $P \colon H_K^s(E) \rightarrow H^{s-1}_L(F)$ is well defined and continuous, because $H^{s-1}_L(F)$ is a closed subspace. 
 
 To prove the elliptic estimate we argue as follows.
 Pick smooth functions $\rho_0$ and $\rho_1$ such that $\rho_0$ is identically $1$ on $K$ and supported within $L^\circ$ and $\rho_1$ is identically $1$ on $\supp \rho_0$ and supported within $L^\circ$. 
 Pick further a pseudo-inverse $R$ for $P$, i.e., an actual pseudo differential operator of order $-1$ whose principal symbol is the inverse of $\symb_1(P)$. 
 Then $\rho_0 - \rho_0RP \rho_1 \in \PseudOpCl^0_L(E)$ and its principal symbol vanishes. 
 Since $\PseudOp^{-1}(E)$ lies dense in $\ker \symb_0$, we can approximate it by actual pseudo differential operators of order $-1$, see Lemma \ref{ASExactSequence - Lemma}, so there is a $Q \in \PseudOp^{-1}_L(E)$  such that $||\rho_0 - \rho_0RP\rho_1 - Q||_{s,s} < 1/2$.

 With that knowledge we can derive the required elliptic estimate from the elliptic estimate of actual pseudo differential operators. 
 For any $u \in \Sobolev[K]{s}{E}$ we have $\rho_j u = u$ and thus
 \begin{align*}
  ||u||_s &=|| (\rho_0 - \rho_0RP\rho_1)u + \rho_0RP \rho_1u||_s\\
  &\leq ||(\rho_0 - \rho_0RP\rho_1) u||_s + ||\rho_0RPu||_s \\
  &\leq ||(\rho - \rho_0RP\rho_1)u - Qu||_s + ||Qu||_s + ||\rho_0R||_{s-1,s}||Pu||_{s-1} \\
  &\leq 1/2 ||u||_s + ||Q||_{s-1,s}||u||_{s-1} + ||\rho_0R||_{s-1,s}||Pu||_{s-1},
 \end{align*}
 where we used the fact that $\rho_0 R$ is an actual pseudo differential operators of order $-1$ that maps sections to sections that are supported within $\supp \rho_0$, to deduce that $\rho_0R$ induces bounded maps between Sobolev spaces.
 The chain of inequalities can be rewritten as
 \begin{align*}
   ||u||_s &\leq 2\left( ||Q||_{s-1,s}||u||_{s-1} + ||\rho_0R||_{s-1,s}||Pu||_{s-1}\right) \\
           &\leq 2 \max\{ ||Q||_{s-1,s}, ||\rho_0R||_{s-1,s}\}\left(||u||_{s-1} + ||Pu||_{s-1}\right),
 \end{align*}
 which gives the second claim.
\end{proof} 

Since the constants in the fundamental elliptic estimate may depend on the compact subsets, the result does not immediately carry over to block operators.
To carry the result over, we use the following partition of unity which is due to Bunke \cite{bunke2009index}, who introduced it in special local coordinates for manifolds with corners.
 \begin{lemma}\label{BunkesPartOfOne - Lemma}
  Let $P \in \PseudOpCl^1(\Spinor_g)$ be a block operator of Dirac type whose core is contained in $M \times (-\rho,\rho)^n$ and that decomposes on all $U_{\rho - 1/2}(\fateps)$. 
  For $\fatalpha \in \Z_3^n$, set 
  \begin{equation*}
   V_\rho(\fatalpha) = \bigl\{(m,x) \in M \times \R^n \, : \, \alpha_j x_j > \rho - 1/2 \text{ if } \alpha_j \neq 0, \, x_j \in (-\rho,\rho) \text{ otherwise}\bigr\}.
  \end{equation*}
  There is a continuous partition of unity 
  $(\phi_0, \dots,\phi_{2n})$ that satisfies the following three properties:
  \begin{enumerate}
   \item[(i)] Each $\phi_l$ is almost everywhere differentiable, has bounded derivative that is zero at 
        infinity, which means $ \lim_{s \to \infty} \mathrm{sup}\{ |\mathrm{d}\phi_l|(x) \,|\, x \in M\times \R^n \setminus (-s,s)^n\} = 0$. 
   \item[(ii)] If $u$ is supported within $V_\rho(\fatalpha)$ with $\fatalpha \neq 0$, then $\lbrack P, \phi_l\rbrack u = \Cliffmult(\mathrm{grad}(\phi_l))(u)$.
   \item[(iii)] $\phi_0$ is compactly supported and for each $l>0$ there is a  unique tupel $(i,\eps) \in \{1,\dots,n\} \times \Z_2$ such that $\phi_l$ is supported within $\{\eps x_i > \rho\}$.   
  \end{enumerate}
 \end{lemma} 
 \begin{proof}
 Fix  $0 < a < (2\sqrt{n})^{-1}$. 
 Let $\chi \colon \R_{\geq 0} \rightarrow \lbrack 0, 1\rbrack$ be a smooth function whose zero set is $\R_{\leq a}$ and that is identically $1$ near $\R_{\geq 2a}$.
 We denote its mirror symmetric extension to $\R$ again with $\chi$. 
 For $i \in \{1,\dots, n\}$, we define 
 $\widetilde{\chi}_i \colon S^{n-1} \rightarrow \lbrack 0, 1\rbrack$ via $\widetilde{\chi}_i(x) = \chi(\pr_i(x))$ and set 
 \begin{equation*}
  \chi_i \deff \frac{\widetilde{\chi_i}}{\sum_{j=1}^{n} \widetilde{\chi}_j }.
 \end{equation*}
 The denominator is always non-zero because it vanishes at $x$ if and only if all summands vanish at $x$, which implies $||x||^2 < n\cdot a < 1$. 
 The $n$-tuple of functions $(\chi_1,\hdots,\chi_{n})$ has the following three evident properties:
 \begin{enumerate}
  \item $\chi_i(x) = 0$ if $|x_i| < a$,
  \item $\sigma^\ast \chi_i = \chi_{\sigma^{-1}(i)}$ for each permutation $\sigma \in \Sigma_n$,
  \item $\sum_j \chi_j = 1$. 
 \end{enumerate}
 For each $\fatalpha = (\alpha_1,\hdots,\alpha_n) \in \Z_3^n$, we define the set
 \begin{align*}
  A_\fatalpha &\deff \{(m,x) \in M \times \R^n \, : \, x_j \in \lbrack -\rho,\rho \rbrack \text{ if } \alpha_j = 0, \ \alpha_jx_j \geq \rho \} \\
  &\cong M \times \lbrack -\rho,\rho \rbrack^{|\mathrm{Null}(\fatalpha)|} \times \prod_{j \in \supp \fatalpha} \alpha_j\R_{\geq \rho}.
 \end{align*}
 The complement of $M \times (-\rho,\rho)^n$ is the union of all $A_\fatalpha$ with $\fatalpha \neq 0$.
 For $(i,\eps) \in \{1,\hdots , n\} \times \Z_2$, we define 
 \begin{equation*}
  \hat{\phi}^\eps_i |_{A_\fatalpha}(m,x) = \begin{cases} 
                                       \chi_i\left( \frac{\hat{x}}{||\hat{x}||}\right), & \text{if } \alpha_i \neq 0, \\
                                       0,    & \text{if } \alpha_i = 0,
                                     \end{cases} 
 \end{equation*}  
 where $\hat{x} = (|\alpha_1|x_1 - \alpha_1\rho, \dots, |\alpha_n|x_n - \alpha_n\rho)$. 
 The assignments agree on the overlaps $A_\fatalpha \cap A_{\bm{\beta}}$ yielding well defined functions on the complement of $M \times (-\rho,\rho)^n$. 
 Each $\hat{\phi}^\eps_i$ is continuous because it is continuous on all $A_\fatalpha$.
 For $0 < \delta < 1$, we pick a smooth function $\xi \colon \R_{\geq 0} \rightarrow \lbrack 0, 1\rbrack$ whose zero set is $\lbrack 0, \rho+\delta\rbrack$ and which is identically $1$ near $\R_{\geq \rho+1}$.
 Its mirror-symmetric extension to $\R$ is again denoted by $\xi$. 
 Set $\xi_i(x) := \xi(\pr_i(x))$ and extend $\hat{\phi_i^\eps}$ to $M \times \R^n$ via  $\phi_i^\eps \deff \xi_i \hat{\phi}_i^\eps$. 
  
 By construction, all $\phi_i^\eps$ sum up to $1$ on the complement of $M \times (-\rho+1,\rho+1)^n$ and the sum is bounded from above by 1 everywhere. 
 Each $\phi^\eps_i$ is supported within $\{\eps x_i \geq \rho + \delta \}$.
 Thus, by setting  $\phi_0 \deff 1 - \sum_{(i,\eps)}\phi_i^\eps$, we have constructed a continuous partition
 of unity. 
 
 The function $\phi_0$ is supported within $M \times (-\rho-1, \rho+1)^n$, which proves property (iii).
 
 The points at which the partition of unity might not be differentiable lie in $\bigcup_{\fatalpha \neq {\bm{\beta}}} A_\fatalpha \cap A_{{\bm{\beta}}}$, which is a zero-set.
 Outside this set, the differential is given by
 \begin{equation*}
  \mathrm{d}\phi_i^\eps = \hat{\phi}_i^\eps(x) \xi'(x_i) \diff x_i  + \xi(x_i)\mathrm{d}\chi_i \cdot \left(\frac{\hat{x} \cdot \hat{x}^T - ||\hat{x}||^2\id}{||\hat{x}||^3}\right) \mathrm{diag}(|\alpha_1|,\dots ,|\alpha_n|).
 \end{equation*}  
 Since $\mathrm{d}\chi_i$ is bounded, we conclude 
 \begin{equation*}
  |\mathrm{d}\phi_i^\eps| \in L^\infty(M \times \R^n) \text{ and } \lim_{s \to \infty} \mathrm{sup}\{ |\mathrm{d}\phi_i^\eps|(x) \,|\, x \in M\times \R^n \setminus (-s,s)^n\} = 0,
 \end{equation*}
 showing property (i). 
 
 By construction, the restriction of each $\phi^\eps_i$ to $V_\rho(\fatalpha)$ with $\fatalpha \neq 0$ depends only on the coordinates $x_j$ indexed by $j \in \supp \fatalpha$. 
 If $u \in \Gamma_c(\Spinor_g)$ is supported within $V_\rho(\fatalpha)$ with $\fatalpha \neq \mathbf{0}$, then the block form of $P$ implies
 \begin{align*}
  \lbrack P, \phi_i^\eps\rbrack u &= \left(P(\fatalpha) + \sum_{j \in \supp \fatalpha} e_j \cdot \partial_j\right)\phi_i^\eps u - \phi_i^\eps \cdot \left(P(\fatalpha) + \sum_{j \in \supp \fatalpha} e_j \cdot \partial_j\right)  u \\
  &= \phi^\eps_i(P(\fatalpha)-P(\fatalpha))u + \sum_{j \in \supp \fatalpha }\lbrack  e_j\cdot \partial_{x_j},\phi_i^\eps\rbrack u \\
  &= \grad{\phi_i^\eps} \cdot u = \Cliffmult(\grad{\phi_i^\eps})(u).
 \end{align*}
 This shows property (ii).
 \end{proof}
 
 We use this partition of unity to derive the fundamental elliptic estimates for block operators of Dirac type.
 The idea of the following two propositions is borrowed from \cite{ebert2017indexdiff}*{Prop 3.3}. 
 
 \begin{proposition}\label{EllipticEstBlock - Prop} 
  Let $P \in \PseudOpCl^1(\Spinor_g)$ be an elliptic block operator of Dirac type. 
  Then $P$ satisfies the fundamental elliptic inequality, which means that there is a $C>0$ such that 
  \begin{equation*}
    ||u||_1 \leq C\left( ||u||_0 + ||Pu||_0\right)
  \end{equation*}
  for all $u \in \Sobolev{1}{\Spinor_g}$ .
 \end{proposition}
 \begin{proof}
  The proof will be carried out by induction over $n$, the dimension of the Euclidean part of our block manifolds $M \times \R^n$.
  For $n = 0$ our base manifold is compact and the theorem follows from Lemma \ref{EllipticEstimateCpt - lemma}.
  
  Assume the statement for $n-1$ and pick an elliptic block operator $P \in \PseudOpCl^{1}(\Spinor_g)$ over a block metric $g \in \BlockMetrics[n]$.
  Pick $R>r>0$ such that $g$ decomposes outside $M \times (-r,r)^n$ and $M \times (-R,R)^n$ is a core for $P$.
  Since $P$ is of block form the classical elliptic estimate gives us 
  \begin{equation*}
    ||u||_1 \leq C_{0,R}\left( ||u||_0 + ||Pu||_0 \right)
  \end{equation*}
  for some $C_{0,R}$ and all sections $u$ with support in $M\times (-R,R)^n$.
  If $u$ is supported within $U_r(i,\eps) = \{\eps x_i > r\}$, then we can identify $\Spinor_g|_{U_r(i,\eps)}$ with $\Spinor_{g(i,\eps)}$ and $P$ with $P(i,\eps) \boxtimes \id + \evenodd \boxtimes \DiracSt$.
  Since $g$ is a product metric here, we can apply Fubini's theorem to get:
  \begin{align*}
   ||u||_1^2 &= \int_{\eps\R_{>r}} ||u(t)||_{1,g(i,\eps)}^2\diff t + \int_{\eps\R_{>r}} ||\partial_t u(t)||_{0,g(i,\eps)}^2 \diff t \\
   &\leq \int C^2\left( ||u(t)||^2_{0,g(i,\eps)} + ||P(i,\eps)u||^2_{0,g(i,\eps)} \right) + ||\partial_t u(t)||_{0,g(i,\eps)}^2 \diff t \\
   &\leq C^2 \int ||u(t)||^2_{0,g(i,\eps)} + ||P(i,\eps)u||^2_{0,g(i,\eps)}  + ||\partial_t u(t)||_{0,g(i,\eps)}^2 \diff t \\
   &= C^2\int ||u(t)||^2_{0,g(i,\eps)} + \langle P(i,\eps)^2u,u\rangle_{0,g(i,\eps)}  - \langle\partial^2_t u(t),u(t)\rangle_{0,g(i,\eps)} \diff t \\
   &= C^2 \int ||u(t)||^2_{0,g(i,\eps)} + \langle P(i,\eps)^2u - \partial^2_t u(t),u\rangle_{0,g(i,\eps)} \diff t \\
   &= C^2 \int ||u(t)||_{0,g(i,\eps)}^2 + \langle(P(i,\eps) \boxtimes \id + \evenodd \boxtimes \DiracSt)^2u(t),u(t) \rangle \diff t \\
   &= C^2 \int ||u(t)||^2_{0, g(i,\eps)} + ||(Pu)(t)||_{0,g(i,\eps)}^2 \diff t \\
   &= C^2(||u||_0^2 + ||Pu||_0^2),
  \end{align*}
  from which we immediately get the fundamental elliptic inequality for this case.
  
  Finally, we need to patch these inequalities together.
  We define the set $L \deff \{0\} \cup \{1,\dots,n\}\times \Z_2$ and let $(\phi_l)_{l \in L}$ be the partition of unity from Lemma \ref{BunkesPartOfOne - Lemma} adapted to $P$.
  We use this partition of unity to patch the inequalities together:
  \begin{align*}
   ||u||_1 &= ||\sum_L \phi_lu||_1 \leq \sum_L ||\phi_lu||_1 \\
   &\leq C_0(||\phi_0u||_0 + ||P\phi_0u||_0) + \sum_{L \setminus \{0\}} C_l(||\phi_lu||_0 + ||P\phi_l u||_0) \\
   &\leq (2n+1)\max\{C_l\}||u||_0 + \sum_L C_l(||\phi_lPu||_0 + ||\lbrack P, \phi_l \cdot\rbrack u||_0 ) \\
   &\leq (2n+1)\max\{C_l,||\lbrack P, \phi_l \cdot\rbrack||_{0,0}\}\left(||u||_0 + ||Pu||_0\right). 
  \end{align*}
  This concludes the induction step. 
 \end{proof}
 
 \begin{definition}\label{InvBlockOp - Def}
   A block operator $P \in \PseudOpCl^1(\Spinor_g)$ of Dirac type is \emph{invertible} if its induced operator $H^1(\Spinor_g) \rightarrow H^0(\Spinor_g)$ is an isomorphism.
\end{definition}
\begin{rem}
 By Proposition \ref{EllipticEstBlock - Prop} the defining condition is equivalent to require that $P$ is an invertible unbounded operator $H^0(\Spinor_{g_P})$.
 Recall that this means that the extension $P\colon \overline{\Gamma_c(\Spinor_{g_P})}^{||\placeholder||_P} \rightarrow H^0(\Spinor_g)$ is invertible, where the domain is the completion of $\Gamma_c(\Spinor_{g_P})$ with respect to the graph norm $||\placeholder||_P = ||\placeholder|| + ||P(\placeholder)||$.
\end{rem}
 
 If the block operator $P$ has invertible faces,  then the previous proposition can be improved.
 
 \begin{proposition}\label{EllipticEstimateInvInf - Prop}
   Let $P$ be a block operator of Dirac type whose faces $P(i,\eps)$ are invertible. 
  Then $P$ satisfies the elliptic estimate 
  \begin{equation*}
    ||u||_1 \leq C\left( ||Pu||_0 + ||\phi_0u||_0 + \sum_{l=1}^{2n} ||\lbrack P, \phi_l\cdot \rbrack u||_0 \right),
  \end{equation*}
  where $\{\phi_l\}$ is the partition of unity from Lemma \ref{BunkesPartOfOne - Lemma}. 
  Furthermore, $P$ has closed image and finite dimensional kernel.
  
  If $P$ is self-adjoint, then $P$ is an (unbounded) Fredholm operator with domain $H^1(M \times \R^n;\Spinor_g)$.
 \end{proposition}
 The proof of this proposition requires the main Fredholm Lemma.  
 \begin{lemma}[\cite{salamon2018functional}*{Lemma 4.3.9}]\label{FAAuxilary - Lemma}
  Let $X, Y$ be Banach spaces and $Z$ be a normed vector space.
  Let $A \colon X \rightarrow Y$ be a bounded linear map and let $K \colon X \rightarrow Z$ be a compact operator.
  If there is a constant $C$ such that for all $x \in X$ the inequality  
  \begin{equation*}
      ||x||_X \leq C(||Ax||_Y + ||Kx||_Z)
  \end{equation*}
  holds, then $A$ has closed image and finite dimensional kernel.
 \end{lemma}
 \begin{proof}[Proof of Proposition \ref{EllipticEstimateInvInf - Prop}] 
  Pick $R>r$ such that $M \times (-R,R)^n$ is a core of $P$ and that the block operator $P$ decomposes on $U_r(i,\eps)$ under the identification
  \begin{equation*}
      U_r(i,\eps) \deff \{(m,x) \in M \times \R^n \, : \, \eps x_i > r\} \cong M \times \R^n(i,\eps) \times \{\eps x_i > r\}
  \end{equation*}
  into the operator $P(i,\eps) \boxtimes \id + \evenodd \boxtimes \DiracSt$.
 
 For each $u \in \Gamma_c(M \times \R^n, \Spinor_g)$ with support in $U_r(i,\eps)$, we have 
 \begin{align*}
  ||Pu||_0^2 &= \int_{M\times \R^n} |Pu|^2 = \int_{M \times \R^n} | P(i,\eps) u + \partial_i\cdot \partial_{i} u|^2 \\
  &= \int_{M \times \R^n}\langle u, (P(i,\eps) + \partial_i\cdot \partial_i)^2u\rangle 
   = \int_{M \times \R^n}\langle u, P(i,\eps)^2 u - \partial_i^2u \rangle \\
  &= \int_{M \times \R^n}|P(i,\eps) u|^2 +  |\partial_i u|^2 \\
  &= \int_{\eps\cdot \R_{\geq r}} ||(P(i,\eps) u)(x)||^2_{L^2(M\times\R^n(i,\eps),\Spinor_g)} + \int_{M \times \R^n} |\partial_iu|^2 \\
  &\geq \int_{\eps\cdot \R_{\geq r}} c^{-2}_{(i,\eps)} ||u(x)||^2_{H^1(M\times\R^n(i,\eps),\Spinor_{g})}  + \int_{M \times \R^n} |\nabla_{\partial_i}u|^2 \\
  &\geq \mathrm{min}\{c^{-2}_{(i,\eps)},1\} ||u||^2_{H^1(M \times \R^n,\Spinor_g)}. 
 \end{align*}
 In this chain of equations, we used in the second line that $P$ is symmetric and that $P(i,\eps)$ anti-commutes with $e_i$. 
 In the fourth line, we applied Fubini's theorem,
 while in the fifth line, we used the invertibility of $P(i,\eps)$ and that $\nabla_{\partial_i} = \partial_i$. 
 In the last line we estimated against the minimum and identified the expression as the Sobolev-1-norm.
 
 On the other hand, if $u$ is supported within $M \times (-R,R)^n$, then $Pu$ is again supported within this relatively compact subset and the classical elliptic estimate implies 
 \begin{equation*}
  ||u||_1 \leq C(||u||_0 + ||Pu||_0).
 \end{equation*}

 We now use the special adapted partition of unity from Lemma \ref{BunkesPartOfOne - Lemma} with $\rho = R$ and the notation from there to patch the inequalities together.
 Using that $\phi_l$ is contained in some $U_r(i,\eps)$ for $l \neq 0$, the previous estimates yield 
 \begin{align*}
  ||u||_1 &= || \sum_L \phi_l u||_1 \leq \sum_L ||\phi_l u||_1 \\
  &\leq C_0(||\phi_0 u||_0 + ||P\phi_0 u||_0) + \sum_{L \setminus \{0\}} C_l||P\phi_l u||_0 \\
  &\leq C_0||\phi_0u||_0 + \sum_L C_l (||\phi_lPu||_0 + ||\lbrack P, \phi_l \cdot\rbrack u||_0) \\
  &\leq C_0||\phi_0u||_0 + \max\{C_l\cdot||\phi_l||_{0,0},1\} \cdot \left( \sum_L ||Pu||_0 + ||[P,\phi_l \cdot ]u||_0  \right)\\
  &\leq C\left( ||Pu||_0 + ||\phi_0u||_0 + \sum_{L} ||\lbrack P, \phi_l \cdot \rbrack u||_0 \right),
\end{align*}   
 where $C_l$ is the inverse of $\min\{c_{(i,\eps)^{-1}}\}$ for the pair $(i,\eps)$ corresponding to $l$.

 It turns out that $\phi_0$ and all $[P,\phi_l \cdot]$ are compact operators.
 Indeed, $\phi_0$ is compactly supported, so $\phi_0 \cdot $ is a compact operator by the classical Rellich-Lemma. 
 
 To see that $T_l \deff \lbrack P, \phi_l \cdot \rbrack$ is compact, note first that it is an element of $\PseudOpCl^1(\Spinor_g)$.
 Symbol calculus shows that its principal symbol vanishes.
 Pick a partition of unity as in Definition \ref{BlockPartOfUnity - Def} such that $\chi^\mathbf{0}$ is identically $1$ on $M \times [-R,R]^n$.
 We decompose it into $T_l = \sum_{\fatalpha \in \Z_3^n} T_l \chi^\fatalpha$. 
 
 By Lemma \ref{CoreLowerBound - Lemma} the summand $T_l\chi^\mathbf{0}$ is supported within $\supp \chi^\mathbf{0}$ in the sense that $\supp T_l \chi^\mathbf{0} u \subseteq \supp \chi^\mathbf{0}$ and $T_l\chi^\fatalpha u = 0$ if $\supp \chi^\fatalpha \cap \supp u = \emptyset$.
 Since the principal symbols of $T_l\chi^\mathbf{0}$ is zero, the Atiyah-Singer exact sequence, see Lemma \ref{ASExactSequence - Lemma}, implies that $T_l\chi^\mathbf{0}$ is compact.
 
 If $\fatalpha \neq 0$, then $\hat{\fatalpha} \deff \fatalpha|_{\supp \fatalpha}$ is not the empty map and  $\chi^\fatalpha u$ is supported within $U_R(\hat{\fatalpha})$. 
 Thus, property (ii) of the partition of unity $\{\phi_l\}$ yields
 \begin{equation*}
  T_l \chi^\fatalpha u = \mathrm{grad} \phi_l \cdot \chi^\fatalpha u. 
 \end{equation*}  
 By Lemma \ref{BunkesPartOfOne - Lemma} (i) and the generalised Rellich lemma for complete manifolds, see \cite{bunke2009index}*{p.50}, we conclude that $T_l \chi_\fatalpha$ is also compact.
 
 A finite sum of compact operators is still compact, so $T_l$ is a compact operator.
 Lemma \ref{FAAuxilary - Lemma} now yields that $P$ has finite kernel and closed image. 
 
 If $P$ is a self-adjoint (unbounded) operator on $L^2(M\times \R^n,\Spinor_g)$, then the maximal domain agrees with the minimal domain.
 The elliptic estimate implies that the graph norm of $P$ is equivalent to $||\placeholder||_1$, so the (minimal) domain of $P$ agrees with $\Sobolev{1}{M \times \R^n, \Spinor_g}$.  
 The kernel and the cokernel of a self-adjoint operator are isomorphic, so the cokernel is also finite dimensional and $P$ is a Fredholm operator.
 \end{proof}
 
 Usually, it is not easy to show that a (pseudo) differential operator is self-adjoint.
 On complete manifolds, there is, however, a general criterion, the boundness of the propagation speed.
 \begin{proposition}[\cite{higson2000analytic}*{Prop. 10.2.11}]
    Let $N$ be a complete Riemannian manifold and let $D$ be a symmetric differential operator on $M$.
    If the \emph{propagation speed}
    \begin{equation*}
        c_D \deff \sup \{||\symb_1(D)(x,\xi)|| \, : \, x \in M, \xi \in T^\vee_xM, ||\xi|| = 1\}
    \end{equation*}
    is finite, then $D$ is (essentially) self-adjoint, i.e. its closure is a self-adjoint operator. 
 \end{proposition}
 \begin{cor}\label{BlockDiracSA - Cor}
  Let $g \in \BlockMetrics[n]$ be a block metric on $M \times \R^n$. 
  Then $\Dirac_g$ is a self-adjoint operator.
  
  More generally, if $P \in \PseudOpCl^1(\Spinor_g)$ is a symmetric, block operator of Dirac type, such that $P - \Dirac_g$ extends to a bounded operator on $H^0(\Spinor_g) = L^2(\Spinor_g)$, then $P$ is a self-adjoint operator.
 \end{cor}
 \begin{proof}
  A block metric is complete by Proposition \ref{PropertiesBlockMetric - Prop} and the propagation speed obviously satisfies $c_{\Dirac_g}=1$.
  By the preceding proposition, $\Dirac_g$ is therefore self-adjoint.
  
  The second statement follows from Exercise 1.9.21 of \cite{higson2000analytic}.
 \end{proof}
 
 A common strategy to verify that a self-adjoint operator is invertible is to verify that it squares to a uniformly positive operator.
 We recall this concept. 
 \begin{definition}\label{boundedbelow - Definition}
  An unbounded operator $T$ on a Hilbert space $H$ is called \emph{non-negative}, if $\langle Tx,x\rangle \geq 0$ for all $x \in \dom T$.
  It is called \emph{uniformly positive}, if there is a constant $c>0$ such that $\langle Tx,x\rangle \geq c ||x||^2$.
  
  An unbounded operator is \emph{bounded from below} if there is a positive constant $c>0$ such that $||Tx|| \geq c ||x||$ for all $x \in \dom T$.
 \end{definition}
 Clearly, a closed unbounded operator $T$ is bounded from below if and only if $T^\ast T$ is uniformly positive.
 The next theorem is an immediate consequence from unbounded functional calculus, see \cite{werner2006funktionalanalysis}.
 \begin{lemma}
 Let $T$ be a self-adjoint unbounded operator on a separable Hilbert space.  
 \begin{enumerate}
  \item[(i)] If $T$ is uniformly positive, then $T$ is invertible.
  \item[(ii)] If $T$ is closed and invertible, then $T^\ast T$ is uniformly positive.  
 \end{enumerate}
\end{lemma}

Being uniformly positive is an open condition for block operators of Dirac type in the space $\mathrm{Hom}(H^1,H^0)$, where $H^j = H^j(\Spinor_g)$.
\begin{lemma}\label{TopInjOpen - Lemma}
 Let $P, Q \in \mathrm{Hom}(H^1,H^0)$ be bounded operators.
 Let $P$ be bounded from below and satisfy a fundamental elliptic estimate.
 If $Q$ is sufficiently close to $P$ in the operator norm of $\mathrm{Hom}(H^1,H^0)$, then $Q$ is also bounded from below.
\end{lemma}
\begin{proof}
 Given $c \in (0,1)$ such that $||Px||_0 > c||x||_0$ and a constant $C=C(P)$ making the fundamental elliptic estimate valid for $P$.
 We claim that all operators $Q$ satisfying the inequality
 \begin{equation*}
     ||P-Q||_{1,0} < \frac{c}{(c+1)C}
 \end{equation*}
 are bounded from below.
 
 Indeed, the fundamental elliptic estimate applied to $P$ and the assumed inequality yield
 \begin{align*}
     ||Qu||_0 &\geq \bigl| ||Pu||_0 - ||(Q-P)u||_0\bigr| \\
     &\geq ||Pu||_0 - ||P-Q||_{1,0}C(||u||_0 + ||Pu||_0) \\
     &= (1-||P-Q||_{1,0}\cdot C)||Pu||_0 - ||P-Q||_{1,0}C||u||_0. 
 \end{align*}
 Using uniform positivity and the assumed inequality a second time, we continue the estimation with
 \begin{align*}
    ||Qu||_0 &\geq c(1-C)||P-Q||_{1,0}||u||_0 - C||P-Q||_{1,0}||u||_0 \\
    &= \underbrace{\bigl(c - (c+1)C||P-Q||_{1,0}\bigr)}_{>0}||u||_0
 \end{align*}
 This proves the claim.
\end{proof}

 \subsubsection*{Gluing Techniques}
 
 In the forthcoming sections we will carry out several cut and paste arguments for block operators.
 Although we will start and end with a block operator on a complete manifold, during the procedures we will consider operators on non-complete manifolds. 
 Being self-adjoint or invertible are global properties of an operator that are usually lost when we restrict the operator to an open subset, so we cannot expect them to hold through the entire cut and paste procedure.
 On the other hand, being symmetric and bounded from below are properties that are preserved under such restrictions and the purpose of this section is to show that they are also preserved under certain gluing constructions.
 
 The first gluing technique for block operator roughly says that being symmetric is a sheaf-like property, provided the operators are suitably local.
 
 \begin{lemma}\label{SymmetrySheafLike - Lemma}
  Let $E \rightarrow M$ be a vector bundle over a not necessarily compact manifold $M$.
  Let $\{U_\alpha\}_{\alpha \in A}$ be a locally finite open cover and let $\{P_\alpha\}_{\alpha \in A}$ be a family of operators
  \begin{equation*}
      P_\alpha \colon \Gamma_c(U_\alpha, E|_{U_\alpha}) \rightarrow \Gamma_c(U_\alpha, E|_{U_\alpha}) \quad \text{ and }\quad  P_\alpha \in \PseudOpCl^1(E|_{U_\alpha})
  \end{equation*}
  such that $P_\alpha|_{U_\alpha \cap U_\beta} = P_\beta|_{U_\alpha \cap U_\beta} \colon \Gamma_c(U_{\alpha\beta},E_{U_{\alpha\beta}}) \rightarrow \Gamma_c(U_{\alpha\beta},E_{U_{\alpha\beta}})$ for all $\alpha,\beta \in A$.
  Assume further that there is a partition of unity $\{\chi_\alpha\}_{\alpha \in A}$ that is subordinated to the cover $\{U_\alpha\}_{\alpha \in A}$ such that the operators acts as derivation on the partition of unity, that is
  \begin{equation*}
      [P_\alpha,\chi_\beta \cdot ]u = i^{-1}\symb_1(P)(\placeholder,\diff \chi_\beta)\cdot u. 
  \end{equation*}
  Then each $P_\alpha$ is symmetric if and only if the operator
  \begin{equation*}
      P \deff \sum_{\alpha \in A} P_\alpha \circ \chi_\alpha
  \end{equation*}
  is symmetric.
 \end{lemma}
 \begin{proof}
  Assume that $P$ is symmetric, then also $P_a$ is symmetric for all $a \in A$ because for all $u \in \Gamma_c(U_a, E|_{U_a})$ we have
  \begin{equation*}
      Pu = \sum_{\alpha \in A} P_\alpha \chi_\alpha u = \sum_{\{\alpha : U_\alpha \cap U_a \neq \emptyset\}} P_a \chi_\alpha u = P_a u.
  \end{equation*}
  
  For the converse, let $u,v \in \Gamma_c(M,E)$ be two sections.
  Set $u_\alpha \deff \chi_\alpha u$ and $v_\beta \deff \chi_\beta v$.
  By bilinearity, it suffices to show that $(P_\alpha u_\alpha,  v_\beta) = (u_\alpha, P_\beta v_\beta)$ for all $\alpha,\beta \in A$.
  
  To verify this, we rely on the following two facts that are immediate consequences from the derivation property:
  \begin{itemize}
      \item[(I)] $\supp \, P_\beta u_\alpha \subseteq \supp \, \chi_\alpha$ for all $\alpha, \beta \in A$.
      \item[(II)] If $u|_{\{\chi_\alpha \neq 0\}} = 0$, then $ (P_\beta u)|_{\{\chi_\alpha \neq 0\}}=0$.
  \end{itemize}
  Indeed, property (I) follows from
  \begin{align*}
      P_\beta u_\alpha &= P_\beta \chi_\alpha u = \chi_\alpha P_\beta u + [P_\beta,\chi_\alpha]u = \chi_\alpha u + i^{-1}\symb_1(P_\beta)(\placeholder, \diff \chi_\alpha)u
  \end{align*}
  and that $\symb_1(P_\beta)(\placeholder, \diff \chi_\alpha) = \symb_1(P_\beta)(\placeholder, 0) = 0$ on $\{\chi_\alpha = 0\}$.
  Property (II) follows from
  \begin{align*}
      (\chi_\alpha P_\beta u)|_{\chi_\alpha \neq 0} = P_\beta\underbrace{\chi_\alpha u}_{=0} - i^{-1}\symb_1(P_\beta)(\placeholder,\diff \chi_\alpha) \underbrace{u|_{\chi_\alpha \neq 0}}_{=0}.
  \end{align*}
  
  Now, for each $\alpha \in A$, we pick a smooth function $\psi_\alpha$ that is identically $1$ near $\supp \, \chi_\alpha$ and supported within $U_\alpha$.
  Symmetry of $P$ now follows from the calculation 
  \begin{align*}
      (P_\alpha u_\alpha, v_\beta) &\overset{(\mathrm{I})}{=} (\psi_\alpha P_\alpha u_\alpha, v_\beta) = (P_\alpha u_\alpha, \psi_\alpha v_\beta) \\
      &= (u_\alpha, P_\alpha \psi_\alpha v_\beta) \\
      &= (u_\alpha, P_\beta \psi_\alpha v_\beta) \qquad \text{because } P_\alpha|_{U_{\alpha\beta}} = P_\beta|_{U_{\alpha\beta}} \\
      &= (u_\alpha, (P_\beta\psi_\alpha v_\beta)|_{\supp\, \chi_\alpha}) \\
      &= (u_\alpha,P_\beta v_\beta)
  \end{align*}
  in which we use (II) and linearity of $P_\beta$ to deduce that $(P_\beta\psi_\alpha v_\beta)|_{\supp\, \chi_\alpha}$ only depends on (the germ of) $\psi_\alpha v_\beta |_{\supp \, \chi_\alpha} = v_\beta|_{\supp \, \chi_\alpha}$.
 \end{proof}

The second gluing technique yields that two invertible block operators with matching faces can be glued together to an invertible operator, provided the ``gluing strip'' is sufficiently large. 
We carry out the gluing construction in a more general setup as we allow our underlying Riemannian manifolds to be non-complete.
This forces us, to work with the following weaker notion that is equivalent to invertibility if the operator in question is self-adjoint.
\begin{definition}\label{Uniform Bounded Below infty - Def}
 Let $P$ be a densely defined operator on a Hilbert space $H$. 
 We call $P$ \emph{(uniformly) bounded from below}, if there is a $c>0$ such that $||Pv|| \geq c ||v||$ for all $v \in \mathrm{dom}(P)$.
 More generally, an operator $P \in\PseudOpCl^1(E)$ is bounded from below, if $P \colon \Gamma_c(M,E) \rightarrow \Gamma(M,E) \subset L^2_{loc}(M,E)$ is bounded from below (meaning we also allow $||P\sigma||_0 = \infty$).
\end{definition}

We work in the following setup.
\begin{setup}\label{SetUp Manifold Parition}
 Let $H \subseteq M$ be a separating hypersurface that is closed as a subset but not necessarily compact.
 Set $M \setminus H = M^- \sqcup M^+$ and let $\Phi \colon H \times (-a,a) \hookrightarrow M$ be an open embedding that restricts to the inclusion on $H\times\{0\} = H$ and satisfies $\Phi(H \times (-a,0)) \subseteq M^-$ and $\Phi(H \times (0,a)) \subseteq M^+$.
 We further denote $M_{\pm a} \deff M^\pm \setminus \Phi(H\times (-a,a))$ and, for all $t,t_1,t_2 \in (-a,a)$, we set
 \begin{align*}
     &M_{<t} \deff M_{-a} \, \cup \, \Phi(H\times(-a,t)), \qquad  M_{>t} \deff \Phi(H \times (t,a)) \, \cup \, M_a, \quad \text{and} \\
     &M_{(t_1,t_2)} \deff \Phi(H \times (t_1,t_2)).
 \end{align*} 
\end{setup}
 
\begin{definition}
 Let $\mathcal{A} \subseteq \mathcal{C}^\infty(M,\R)$ be a sub-algebra of smooth functions on $M$. 
 An operator $P \in \PseudOpCl^1(E)$ is an $\mathcal{A}$\emph{-derivation}, if 
 \begin{equation*}
  [P, f\cdot ] = i^{-1}\symb_1(P)(\placeholder,\diff f)
 \end{equation*}
 for all $f \in \mathcal{A}$.
\end{definition}
 
 We are interested in the algebra $\mathcal{A} \deff (\pr_2 \circ \Phi^{-1})^\ast \mathcal{C}_c^\infty((-a,a),\R))$ considered as a sub-algebra of $\mathcal{C}^\infty(M,\R)$ (via extension by zero).
 
\begin{lemma}\label{LocalityLemma - Lemma}
 Let $E \rightarrow M$ be a vector bundle over a manifold as in setup \ref{SetUp Manifold Parition}.
 Assume that  $P \in \PseudOpCl^1(E)$ restricts to $M_{>b}$ and $M_{<-b}$ for some $0<b<a$ and that $P$ is an $\mathcal{A}$-derivation.
 If $\sigma \in \Gamma_c(M,E)$ is supported within $M_{(t_1,t_2)}$, then $P\sigma$ is also supported within $M_{(t_1,t_2)}$ for all $-a<t_1<t_2<a$.
 In particular, $(P\sigma)|_{M_{(t_1,t_2)}}$ only depends on the germ of $\sigma|_{\overline{M}_{(t_1,t_2)}}$.
 
 Analogous statements hold for $M_{>t}$ and $M_{<t}$ for all $t \in (-b,b)$.
\end{lemma}
\begin{proof}
 We only prove the statement for $M_{(t_1,t_2)}$. 
 The other proofs are quite similar.
 Assume $\sigma$ is supported within this set. 
 Since the support is compact, there are $t_1 < t_1' < t_2' < t_2$ such that $\sigma$ is also supported within $M_{(t_1't_2')}$.
 Let $\chi \colon (t_1,t_2) \rightarrow [0,1]$ be identically $1$ on $(t_1',t_2')$ and compactly supported. 
 To ease the notation, we denote $\chi \circ \pr_2 \circ \Phi^{-1}$ (and its extension by zero) again with $\chi$.
 We have
 \begin{align*}
     P\sigma &= P \chi \sigma = \chi P\sigma + [P,\chi \cdot]\sigma \\
     &= \chi P\sigma + i^{-1}\symb_1(P)(\placeholder,\chi') \sigma \\
     &= \chi P \sigma,
 \end{align*}
 because $\chi \equiv 1$ on $\supp \sigma$. 
 Thus, $P\sigma$ is contained in $\supp \chi \subseteq M_{(t_1,t_2)}$.
 
 To prove that $P\sigma|_{M_{(t_1,t_2)}}$ only depends on the germ of $\sigma|_{\overline{M_{(t_1,t_2)}}}$, it suffices to show that if $\sigma$ is identically zero near $\overline{M_{(t_1,t_2)}}$, then $P\sigma$ vanishes identically on $M_{(t_1,t_2)}$.
 But if $\sigma$ vanishes near $\overline{M_{(t_1,t_2)}}$, then, by compactness of $\supp \sigma$, there are $t_1'<t_1$ and $t_2'>t_2$ such that $\supp(\sigma) \cap M_{(t_1',t_2')} = \emptyset$.
 Thus $\supp(\sigma)$ is supported within $M_{<t_1'}$ and $M_{>t_2'}$, so $P\sigma$ must be supported there, too.
\end{proof}

\begin{definition}\label{SetUp2}
 Let $R<b<a$ and let $P_1$, $P_2$ be operators on $M_{<a}$ and $M_{>-a}$, respectively.
 Assume that $P_1$ and $P_2$ restrict to $M_{<-b}$ and $M_{>b}$, respectively, and that they are $\mathcal{A}$-derivations. 
 By the previous Lemma, the two operators restrict to $M_{(-R,R)}$ and we assume that the restrictions agree.
 This allows us to patch the two operators together: Let $\chi \colon \R \rightarrow [0,1]$ is a smooth function that is identically $1$ near $\R_{\leq -R}$ and identically $0$ near $\R_{\geq R}$. 
 We extend he function via $0$ and $1$ to $M$ and denote the result again with $\chi$.
 We define
 \begin{equation*}
    P_1 \cup P_2 \deff P_1 (1-\chi) + P_2 \chi. \nomenclature{$P_1 \cup P_2$}{Convex-combination of two sufficiently local operators}
 \end{equation*}
\end{definition}

Since $P_1$ and $P_2$ agree on $M_{(-R,R)}$ and act as $\mathcal{A}$-derivations, the result does not depend on the choice of $\chi$.
An application of Lemma \ref{SymmetrySheafLike - Lemma} shows that $\cup$ preserves symmetry.
\begin{cor}\label{SymmetryUnion - Cor}
 The operator $P_1 \cup P_2$ is symmetric if and only if the operators $P_1$ and $P_2$ are symmetric.
\end{cor}
\begin{proof}
 The only if direction is easy because restrictions of symmetric operators are symmetric and $P_1$ agrees with $P_2$ on $M_{(-R,R)}$.
 
 For the other direction, we will prove that $M_{<R}$ and $M_{> -R}$ together with $\chi$ and $1-\chi$ from Definition \ref{SetUp2} satisfy the condition of Lemma \ref{SymmetrySheafLike - Lemma} for $P_1$ and $P_2$.
 To this end, fix a function $\eta \colon (-a,a) \rightarrow [0,1]$ with compact support and $\eta|_{(-b,b)} \equiv 1$.
 Let $\eta^{-1}$, $\eta^0$, and $\eta^1$ be the functions obtained from the construction in Definition \ref{BlockPartOfUnity - Def} but considered as function on $M$.
 They form a partition of unity.
 
 For all compactly supported sections $\sigma \in \Gamma_c(M_{<R},E)$, the section $\eta^{-1}\sigma$ is supported within $M_{<-b}$, so that $P_1\eta^{-1}\sigma$ is also supported there by assumption.
 This implies
 \begin{align*}
     P_1\chi \sigma &= P_1\chi \eta^{-1}\sigma + P_1 \chi \eta^{0}\sigma = P_1\eta^{-1} + P_1\chi\eta^0 \sigma\\
                    &= \chi P_1\eta^{-1} \sigma + \chi P_1\eta^0 \sigma + [P_1,\chi]\eta^0 \sigma \\
                    &= \chi P_1(\eta^{-1} + \eta^0)\sigma + i^{-1}\symb_1(P_1)(\placeholder,\chi')\cdot \eta^0 \sigma \\
                    &= \chi P_1 \sigma + i^{-1}\symb_1(P_1)(\placeholder,\chi') \sigma.
 \end{align*}
 Note that the last equation follows from $\chi' \equiv 0$ on $\supp \eta^{-1}$.
 Of course, this implies that $P_1$ acts on $1-\chi$ as a derivation.
 
 The corresponding statement for $P_2$ and $(1-\chi)$ can be derived in an analogous manner. 
\end{proof}

\begin{lemma}\label{LowerBndGluing - Lemma}
 Let $P_1$ and $P_2$ be operators on $M_{<a}$ and $M_{>-a}$ as in setup \ref{SetUp2}.
 Assume that 
 \begin{equation*}
     ||\symb_1(P_1)(\placeholder,\diff t)||_{\infty,M_{(-R,R)}} \deff \{||\symb_1(P_1)(x,\diff t)||_{\op} \, : \, x \in M_{(-R,R)}\} \leq C.
 \end{equation*}
 If there is a section $u \in \Gamma_c(M,E)$ such that $||(P_1\cup P_2) u||_0 < \eps ||u||_0$, then there is a $v \in \Gamma_c(M,E)$ that is supported within $M_{<R}$ that satisfies 
 \begin{equation*}
     ||P_1v||_0^2  < 4\left(\eps^2 + \frac{C^2}{R^2}\right)||v||_{0}^2   
 \end{equation*}
 or is supported within $M_{>-R}$ and satisfies an analogous inequality with $P_1$ replaced by $P_2$.
\end{lemma}
\begin{proof}
 Abbreviate $P_1 \cup P_2$ to $Q$.
 Let $u \in \Gamma_c(M,E)$ with $||Qu||_0^2 < \eps^2 ||u||_0^2$.
 Assume without loss of generality, that $||u||_{0,M_{\leq 0}}^2 \geq 1/2||u||_0^2$.
 Let $\eta \colon \R \rightarrow [0,1]$ be a smooth function that satisfies $\eta \equiv 1$ on $\R_{\leq 0}$, $\eta \equiv 0$ on $\R_{\geq 1}$, and $|\eta'|\leq 1.2$. 
 Set $\eta_R(t) = \eta(R^{-1}t)$ and extend this function as described above to $M$. 
 We denote the result again with $\eta_R$.
 
 The norm splits into
 \begin{align*}
     ||P_1\eta_Ru||_0^2 = ||P_1\eta_Ru||_{0,M_{\leq -R}}^2 + ||P_1u||_{0,M_{(-R,R)}}^2  
 \end{align*}
 and, by Lemma \ref{LocalityLemma - Lemma}, the summands only depend on the germs of the restrictions to these domains. 
 For example,
 \begin{equation*}
    ||P_1\eta_Ru||_{0,M_{\leq -R}}^2 = ||P_1\eta_Ru|_{M_{\leq -R}}||_{0,M_{\leq -R}}^2 = ||P_1u||_{0,M_{\leq -R}}^2. 
 \end{equation*}
 Since also $Q$ satisfies the assumptions of Lemma \ref{LocalityLemma - Lemma}, we have
 \begin{equation*}
     (P_1\eta_Ru)|_{M\leq-R} = P_1(\eta_Ru|_{M\leq-R}) = (P_1u)|_{M\leq-R} = (Qu)|_{M\leq-R}
 \end{equation*}
 and because $P_1$ and $P_2$ agree on $M_{(-R,R)}$ we have 
 \begin{align*}
     (Q\eta_Ru)|_{M_{(-R,R)}} &= Q\eta_R(u|_{M_{(-R,R)}}) = P_1\eta_R(u|_{M_{(-R,R)}}) \\
     &= \eta_RP_1(u|_{M_{(-R,R)}}) + [P,\eta_R \cdot](u|_{M_{(-R,R)}}) \\
     &= (\eta_RQu)|_{M_{(-R,R)}} + i^{-1}\symb_1(P_1)(\placeholder,\eta_R'\diff t)(u)|_{M_{(-R,R)}}.
 \end{align*}
 This implies
 \begin{align*}
     ||P_1 \eta_R u||_0^2 &= ||Qu||^2_{0,M_{\leq -R}} + ||\eta_RQu + \symb_1(P_1)(\placeholder,\eta_R'\diff t)(u)||_{0,M_{(-R,R)}}^2 \\
     &\leq ||Qu||^2_{0,M_{\leq -R}} + 2||\eta_RQu||^2_{0,M_{(-R,R)}} + 2||\symb_1(P_1)(\placeholder,\eta_R'\diff t)(u)||_{0,M_{(-R,R)}}^2 \\
     &\leq 2||Qu||^2_{0,M_{<R}} + 2||\symb_1(P_1)(\placeholder,\eta_R'\diff t)(u)||_{0,M_{(-R,R)}}^2 \\
     &< 2 \eps^2 ||u||^2_0 + 2\frac{C^2}{R^2} ||u||_0^2, 
 \end{align*}
 where we used in the last inequality that $||Qu||_0^2 < \eps^2||u||^2_0$ and that
 \begin{align*}
     &\quad \  ||\symb_1(P_1)(\placeholder,\eta_R'\diff t)(u)||_{0;M_{(-R,R)}} \\
     &\leq ||\symb_1(P_1)(\placeholder,\eta_R'\diff t)||_{0,0;M_{(-R,R)}} ||u||_{0,M_{(-R,R)}} \\
     &\leq ||\symb_1(P_1)(\placeholder,\eta_R'\diff t)||_{\infty;M_{(-R,R)}} ||u||_0\\
     &\leq \frac{|\eta'|_\infty}{R}||\symb_1(P_1)(\placeholder,\diff t)||_{\infty;M_{(-R,R)}} ||u||_0\\
     &\leq \frac{C}{R}||u||_0
 \end{align*}
 From $1/2||u||_0^2 \leq ||u||_{0,M_{\leq 0}}^2$ it follows that $||u||_0^2 \leq 2||\eta_Ru||_0^2$.
 Putting all inequalities together, we get
 \begin{equation*}
     ||P_1\eta_R u||_0^2 < 4\left(\eps^2 + C^2/R^2 \right)||\eta_Ru||_0^2,
 \end{equation*}
 so the claim follows for $v \deff \eta_R u$.
\end{proof}
The following implication is an indispensable tool in the following sections because it allows us to glue invertible block operators of Dirac type together.
Note that, in this set up, we are allowed to choose the length of the ``gluing-strip'' to be arbitrarily large without increasing the propagation speed of the operator.
\begin{cor}\label{Gluing Lower Bounds - Cor}
 Let $P_1$ and $P_2$ be as in the previous lemma.
 If they are bounded from below with lower bounds $c_1$ and $c_2$ respectively and if $R > \max\{C^2/c_1^2,C^2/c_2^2\}$, then $P_1 \cup P_2$ is bounded from below with lower bound 
 \begin{equation*}
     c_{P_1\cup P_2}^2 = \frac{\min\left\{{c_1^2 -C^2/R^2}, c_2^2 - C^2/R^2\right\}}{4}.
 \end{equation*}
\end{cor}
\begin{proof}
 By assumption on $R$, the constant $c_{P_1 \cup P_2}$ is positive.
 Assume that there is a section $u \in \Gamma_c(M,E)$ such that $||(P_1 \cup P_2)(u)||_0 <  c_{P_1\cup P_2}||u||_0$, then the previous Lemma implies the existence of $v \in \Gamma_c(M,E)$ that is, without loss of generality, supported within $M_{<R}$ and satisfies 
 \begin{align*}
     ||P_1v||_0^2 &< 4\left(c_{P_1\cup P_2}^2 + C^2/R^2\right)||v||_0^2 < c_1^2||v||_0^2,
 \end{align*} 
 which contradicts that $P_1$ is bounded from below by $c_1$.
\end{proof}

 \section{Foundations of the Operator Concordance Set}\label{Section - Foundations of the Operator Concordance Set}

With the theory developed in the previous section we are finally in the position to define the operator concordance set $\InvBlockDirac$ as the cubical subset of all block Dirac operators.

But, first we give a combinatorial model for the classifying space for real $K$-theory, by considering the cubical set of block maps into $\PseudDir[ ]$ and $\InvPseudDir[ ]$ instead of the spaces themselves.

\begin{definition}
 Let $\PseudDir[\bullet]$ be the cubical set whose $n$-cubes consists of smooth block maps $P \colon \R^n \rightarrow \PseudDir[]$. 
 The connecting maps are given by
 \begin{equation*}
     \face{i}P = \lim_{R \to \infty} P \circ \CubeIncl[R\eps]{i} \qquad \text{ and } \qquad \degen{i} P = P \circ \CubeProj{i}.
 \end{equation*}
 Let $\InvPseudDir[\bullet]$ be the cubical subset consisting of all smooth block maps with values in $\InvPseudDir$.
\end{definition}

Recall that $\PseudDir$ is an affine bundle over $\Riem(M)$ with fibre $\PseudDir_g$, the set of all self-adjoint, odd, Clifford-linear pseudo differential operators of order $1$ whose principal symbol agree with the principal symbol of the Dirac operator $\Dirac_g$. 
The map $\PseudDir \rightarrow \Riem(M)$ assigns to an operator $P$ its underlying metric $g$.
This map induces a map of cubical sets $\PseudDir[\bullet] \rightarrow \SingMet$.

Now, we are going to construct the combinatorial reference space.

\begin{definition}\label{BlockDiracOperator - Def}
  Let $g$ be a block metric on $M \times \R^n$.
  A block operator $P \in \PseudOpCl^1(\Spinor_g)$ of Dirac type is called a \emph{block Dirac operator} if it is symmetric and $P - \Dirac_g$ extends to a bounded endomorphism on $H^0(\Spinor_g) = L^2(\Spinor_g)$.
\end{definition}

An immediate consequence of Corollary \ref{BlockDiracSA - Cor} is that a block Dirac operator is automatically essentially self-adjoint.

For the definition of the operator concordance set it will be important to trace permutation of coordinates.
A permutation $\sigma \in \mathrm{Aut}(\{1,\dots,n\}) = S_n$ induces a linear map
\begin{equation*}
     \sigma \colon (x_1,\dots, x_n) \mapsto (x_{\sigma^{-1}(1)},\dots, x_{\sigma^{-1}(n)}). 
\end{equation*}
If $g$ is a block metric on $M \times \R^n$, then $\sigma^\ast g$ is also a block metric on $M \times \R^n$.
Since $\R^n$ carries up to equivalence only one $\mathrm{Pin}^-$-structure, permutations are Pin structure preserving.
Functoriality of the spinor bundle construction yields a map, which we denote with
\begin{equation*}
      \Spinor(\sigma) \colon \Spinor_{\sigma^\ast g} \rightarrow \Spinor_{g} \quad \text{ or equivalently } \quad \Spinor(\sigma) \colon \Spinor_g \rightarrow \Spinor_{\sigma_\ast g}.
\end{equation*}
We are particularly interested in cyclic permutations.
\begin{definition}
    For positive integers $a < b$ let $\cycl(a,b)$ 
    \nomenclature{$\cycl(a,b)$}{Cyclic permuations of elements between $a$ and $b$}
    be the cyclic permutation given by $a\mapsto a+1$ and $b\mapsto a$ and let $\cycl(b,a) = \cycl(a,b)^{-1}$.
   We denote the induced maps on spinor bundles with $\Cycl(a,b)$ 
   \nomenclature{$\Cycl(a,b)$}{Induced map of $\cycl(a,b)$ on spinor bundles}
   and $\Cycl(b,a)$, respectively.
\end{definition}

\begin{definition}
  Let $\BlockDirac$ be the cubical set whose set of $n$-cubes is given by
  \begin{equation*}
      \BlockDirac[n] \deff \{P \in \PseudOpCl^1(\Spinor_g) \, : \, g \in \BlockMetrics[n], \, P\text{ block Dirac operator}\}
      \nomenclature{$\BlockDirac[n]$}{Cubical Set of block Dirac operators}
  \end{equation*}
  and whose connecting maps are given by
  \begin{align*}
      \face{i} P &\deff \lim_{R \to \infty}{\CubeIncl[R\eps]{i}}^\ast P(i,\eps) \in \PseudOpCl^1(\Spinor_{\face{i}g}), \\
      \degen{i}P &\deff \Cycl(i,n+1)_\ast(P \boxtimes \id + \evenodd \boxtimes \DiracSt) \in \PseudOpCl^1(\Spinor_{\degen{i}g})
  \end{align*}
\end{definition}
\begin{rem}
 A comment to the (subtle) notation.
 Recall that if a block operator $P$ with underlying block metric $g$ decomposes on $M \times U_r(i,\eps) = \{(m,x) \, : \, \eps x_i > r\}$ into $P(i,\eps)\boxtimes \id + \evenodd \boxtimes \DiracSt$ under the isometry $U_r(i,\eps) \cong M \times \{\eps x_i = r\} \times \eps\R_{>r}$ given by permutation, then the operator $P(i,\eps)$ restricts to a block operator on $\Spinor_{g(i,\eps)} \rightarrow M \times \{\eps x_i = R\}$ for all sufficiently large $R$ and is independent of $R$.
 The map 
 \begin{equation*}
     \CubeIncl[R\eps]{i} \colon (M \times \R^{n-1},\face{i}g) \rightarrow \bigl(M \times \{\eps x_i = R\},g(i,\eps)\bigr)
 \end{equation*}
 is a Pin-structure preserving isometry and thus provides a bundle isometry $\Spinor(\CubeIncl[R\eps]{i})$ between the corresponding spinor bundles that induces an isomorphism between the corresponding sections with compact support.
 We use this isomorphism to pull $P(i,\eps)$, an operator on $\Gamma_c(\Spinor_{g(i,\eps)})$, back to ${\CubeIncl[R\eps]{i}}^\ast P(i,\eps)$, an operator on $\Gamma_c(\Spinor_{\face[\eps]{i}g})$.
 
 For the degeneracy map, we use that the cyclic permutation 
 \begin{equation*}
    \cycl(i,n+1) \colon \left(M \times \{\eps x_i > r\},\degen{i}g\right) \rightarrow \left(M \times \R^n \times \eps \R_{>r}, g \oplus \diff x_{n+1}^2\right) 
 \end{equation*}
 is an isometry and the natural isomorphism $\Spinor_g \boxtimes \Cliff[1][0] \cong \Spinor_{g \oplus \diff x_{n+1}^2}$. 
\end{rem}

\begin{lemma}\label{ConnectingMapsWellDef - Lemma}
 The connecting maps of $\BlockDirac$ are well-defined and satisfy the cubical identities. 
\end{lemma}

The proof requires the following simple functional analytic result.
\begin{lemma}\label{OperatornormTensoProduct - Lemma}
 Let $H_1$, $H_2$ be Hilbert spaces and $A \colon H_1 \rightarrow H_1$ be an unbounded operator.
 Then $A \otimes \id$ is a bounded operator on $H_1 \otimes H_2$ if and only if $A$ is bounded on $H_1$.
 In this case, the operator norms satisfy $||A \otimes \id||_{\op} = ||A||_{\op}$.
\end{lemma}
\begin{proof}
 Since the algebraic tensor product $H_1 \odot H_2$ lies dense in $H_1 \otimes H_2$, so does $\dom(A \otimes \id) \deff \dom A \odot H_2$.
 If $A$ is bounded, then we can continuously extend it to an operator on $H_1 \otimes H_2$ and the operator norm of this extension necessarily agrees with $||A||_{\op}$.
 
 Conversely, assume that $A \otimes \id$ is bounded.
 Let $v\in H_2$ be a unit-length vector.
 The map $\iota_v \colon H_1 \rightarrow H_1 \otimes H_2$ given by $x \mapsto x \otimes v$ is linear and bounded.
 The map $x \otimes y \mapsto \langle y,v\rangle \cdot x$ extends from $H_1 \odot H_2$ to a bounded linear map $p_v \colon H_1 \otimes H_2 \rightarrow H_1$.
 From $A = p_v \circ A \otimes \id \circ \iota_v$ we deduce that $A$ is a bounded operator.
\end{proof}

\begin{proof}[Proof of Lemma \ref{ConnectingMapsWellDef - Lemma}]
 By Corollary \ref{BlockDescend - Cor}, $\face{i}P$ is a again a block operator.
 Symbol calculus shows that the principal symbol of $\face[\eps]{i} P$ agrees with the principal symbol of $\Dirac_{\face{i}g}$, so $\face{i}P$ is a block operator of Dirac type.
 In a similar fashion, but using Lemma \ref{BlockSuspends - Lemma} instead, we see that $\degen{i}P$ is a block operator of Dirac type.
 
 We need to show that $\face{i}P$ and $\degen{i}P$ are symmetric operators if $P$ is symmetric.
 We start with $\degen{i}P = \Cycl(i,n+1)_\ast(P \boxtimes \id + \evenodd \boxtimes \DiracSt)$.
 Since $\Cycl(i,n+1)$ induces an isometry between the corresponding Hilbert spaces of square integrable sections, the push-forward $\Cycl(i,n+1)_\ast$ maps symmetric operators to symmetric operators.
 Thus, it suffices to prove that $P \boxtimes \id + \evenodd \boxtimes \DiracSt$ is a symmetric operator on 
 \begin{equation*}
     L^2(M\times \R^{n+1};\Spinor_{g \oplus \diff x_{n+1}^2}) \cong L^2(M \times \R^n;\Spinor_g)\otimes L^2(\R,\Cliff[1][0]),
 \end{equation*}
 which is true because the tensor product of two symmetric operators is again symmetric and all four operators in question $P$, $\id$, $\evenodd$ and $\DiracSt$ are symmetric.
 
 To argue that $\face{i}P$ is symmetric, we proceed as follows:
 On $U(i,\eps) \deff U_R(i,\eps)$, for some sufficiently large $R>0$, we can identify the symmetric operator $P|_{U(i,\eps)}$ with $P(i,\eps) \boxtimes \id + \evenodd \boxtimes \DiracSt_{\R}$ via an isometric bundle morphism.
 The second summand is symmetric, so $P(i,\eps) \boxtimes \id$ must be a symmetric operator on $\Spinor_{g(i,\eps)\oplus \diff x_n^2} \rightarrow (M \times \R^n)(i,\eps) \times \eps\R_{>R}$.
 This implies that $P(i,\eps)$ is symmetric on $\Spinor_{g(i,\eps)}$ and consequently that $\face{i}P = \CubeIncl[R\eps]{i}^\ast P(i,\eps)$ is symmetric because $\CubeIncl[R\eps]{i} \colon (M \times \R^{n-1},\face{i}g) \rightarrow (M\times \R^n(i,\eps),g(i,\eps))$ is an isometry. 
 
 From 
 \begin{align*}
     &\quad \ \degen{i}P - \Dirac_{\degen{i}g} \\
     &= \Cycl(i,n+1)_\ast\left(P \boxtimes \id + \evenodd \boxtimes \DiracSt\right) - \Cycl(i,n+1)_\ast\left(\Dirac_g \boxtimes \id + \evenodd \boxtimes \DiracSt\right) \\
     &= \Cycl(i,n+1)_\ast((P - \Dirac_g)\boxtimes \id)
 \end{align*}
 and Lemma \ref{OperatornormTensoProduct - Lemma} it follows that the difference $\degen{i}P - \Dirac_{\degen{i}g}$ extends to a bounded operator on $H^0(\Spinor_{\degen{i}g})$ with operator norm $||P - \Dirac_g||_{0,0}$.
 
 To see that $\face{i}P - \Dirac_{\face[\eps]{i}g}$ extends to a bounded operator on $H^0(\Spinor_{\face{i}g})$ we argue as follows.
 For $R>0$ sufficiently large the operator $P - \Dirac_g$ restricts to $\Gamma_c(U_R(i,\eps),\Spinor_g)$. 
 Since $P -\Dirac_g$ extends to a bounded operator on $H^0(\Spinor_g)$, the restriction extends to a bounded operator on $H^0(U_R(i,\eps),\Spinor_g)$, which is isometric equivalent to $H^0(M\times \{\eps x_i = R\},\Spinor_{g(i,\eps)}) \otimes L^2(\eps\R_{>R},\Cliff[1][0])$.
 Under this identification, the restriction of $P - \Dirac_g$ agrees with $\left(P(i,\eps) - \Dirac_{g(i,\eps)}\right) \boxtimes \id$.
 By Lemma \ref{OperatornormTensoProduct - Lemma}, this implies that $P(i,\eps) - \Dirac_{g(i,\eps)}$ is bounded.
 Thus, $\face{i}P - \Dirac_{\face{i}g} = \CubeIncl[R\eps\,\ast]{i}\left(P(i,\eps) - \Dirac_{g(i,\eps)}\right)$ is a bounded operator on $H^0(\Spinor_{\face{i}g})$.
 
 Lastly, we need to show, that the connecting maps satisfy the cubical identities.
 For $i < j$ and $R > 0$ sufficiently large, we have
 \begin{align*}
  \face{i}\face[\omega]{j}P &=  \face{i}\left(\CubeIncl[\omega R]{j}^\ast P(j,\omega)\right) \\
  &= \CubeIncl[\eps R]{i}^\ast\left(\CubeIncl[\omega R]{j}^\ast \left(P(j,\omega)\right)(i,\eps)\right) \\
  &= \left(\CubeIncl[\omega R]{j}\CubeIncl[\eps R]{i}\right)^\ast P(i,j,\eps,\omega) \\
  &= \left(\CubeIncl[\eps R]{i}\CubeIncl[\omega R]{j-1}\right)^\ast P(i,j,\eps,\omega) \\
  &= \CubeIncl[\omega R]{j-1}^\ast\left(\CubeIncl[\eps R]{i}^\ast \left(P(i,\eps)\right)(j-1,\omega)\right) \\
  &= \face[\omega]{j-1}\face{i}P.
 \end{align*}
 For $i < j$ and a sufficiently large $R>0$, we calculate 
 \begin{align*}
  \face{i}\degen{j}(P) &= \CubeIncl[R\eps]{i}^\ast(\degen{j}(P)(i,\eps)) \\
  &= \CubeIncl[R\eps]{i}^\ast\Cycl(j,n+1)_\ast\left(P \boxtimes \id + \evenodd \boxtimes \DiracSt\right)(\cycl(j,n+1)^{-1}(i),\eps) \\
  &= \CubeIncl[R\eps]{i}^\ast\Cycl(j,n+1)_\ast(P(i,\eps) \boxtimes \id + \evenodd \boxtimes \DiracSt) \quad \quad \quad \text{ for } i < j \\
  &= \left(\CubeIncl[R\eps]{i}^{-1}\right)_\ast \Cycl(j,n+1)_\ast(P(i,\eps) \boxtimes \id + \evenodd \boxtimes \DiracSt) \\
  &= \Cycl(j-1,n)_\ast \left(\CubeIncl[R\eps]{i}^{-1}\right)_\ast \left(P(i,\eps) \boxtimes \id + \evenodd_{\Spinor_{g(i,\eps)}} \boxtimes \DiracSt\right) \\
  &= \Cycl(j-1,n)_\ast(\face{i} P \boxtimes \id + \evenodd_{\Spinor_{\face{i}g}} \boxtimes \DiracSt) \\
  &= \degen{n-1} \face{i}P,
 \end{align*}
 and analogously $\face{i}\degen{j}P = \degen{j}\face{i-1}P$ for $i>j$. 
 Finally, with a slight abuse of notation in the last lines, we calculate
 \begin{align*}
  \face{i}\degen{i}P &= \CubeIncl[R\eps]{i}^\ast\left((\degen{i}P)(i,\eps)\right)\\
  &= \CubeIncl[R\eps]{i}^\ast\left(\Cycl(i,n+1)_\ast(P \boxtimes \id + \evenodd \boxtimes \DiracSt)(i,\eps)\right) \\
  &= \CubeIncl[R\eps]{i}^{-1}_\ast \Cycl(i,n+1)_\ast\left((P \boxtimes \id + \evenodd \boxtimes \DiracSt)(\cycl(i,n+1)^{-1}(i),\eps)\right) \\
  &= \CubeIncl[R\eps]{i}^{-1}_\ast \Cycl(i,n+1)_\ast\left((P \boxtimes \id + \evenodd \boxtimes \DiracSt)(n+1,\eps)\right) \\
  &= \CubeIncl[R\eps]{n+1}^{-1}_\ast\left((P \boxtimes \id + \evenodd \boxtimes \DiracSt)(n+1,\eps)\right) = P. 
 \end{align*}
 For $i \leq j$, we use the identity
 \begin{equation*}
   \cycl(i,n+2) \circ \cycl(j,n+1) = \cycl(j+1,n+2) \circ \cycl(i,n+1) \circ \tau_{n+1,n+2}
 \end{equation*}
 in which $\tau_{n+1,n+2}$ denotes the transposition that interchanges $n+1$ and $n+2$, to conclude 
 \begin{align*}
  \degen{i}\degen{j}(P) &= \degen{i}(\Cycl(j,n+1)_\ast(P \boxtimes \id + \evenodd \boxtimes \DiracSt)) \\
  \begin{split}
   &= \Cycl(i,n+2)_\ast\left(\Cycl(j,n+1)_\ast(P \boxtimes \id + \evenodd_{\Spinor_g}\boxtimes \DiracSt) \boxtimes \id \right. \\
   &\quad \left. + \evenodd_{\Spinor_{\degen{j}(g)}} \boxtimes \DiracSt \right) 
  \end{split}\\
  \begin{split}
   &= \Cycl(i,n+2)_\ast \circ \Cycl(j,n+1)_\ast \biggl( P \boxtimes \id \boxtimes \id  \\
   &\quad + \left. \evenodd_{\Spinor_g} \boxtimes \DiracSt \boxtimes \id + \evenodd_{\Spinor_{\degen{n+1}(g)}} \boxtimes \DiracSt \right)
  \end{split} \\
  \begin{split}
   &= \Cycl(i,n+2)_\ast \circ\Cycl(j,n+1)_\ast\biggl( P \boxtimes \id \boxtimes \id  \\
   & \quad \left.+ \evenodd_{\Spinor_g} \boxtimes \DiracSt \boxtimes \id + \evenodd_{\Spinor_{g}} \boxtimes \evenodd \boxtimes \DiracSt \right)
  \end{split}\\
  \begin{split}
   &= \Cycl(j+1,n+2)_\ast \circ \Cycl(i,n+1)_\ast \circ {\tau_{n+1,n+2}}_\ast \biggl(  P \boxtimes \id \boxtimes \id  \\
   & \quad \left. + \evenodd_{\Spinor_g} \boxtimes \DiracSt \boxtimes \id + \evenodd_{\Spinor_{g}} \boxtimes \evenodd \boxtimes \DiracSt \right)
  \end{split}\\
  \begin{split}
   &= \Cycl(j+1,n+2)_\ast \circ \Cycl(i,n+1)_\ast \biggl(  P \boxtimes \id \boxtimes \id  \\
   & \quad \left. + \evenodd_{\Spinor_g} \boxtimes \DiracSt \boxtimes \id + \evenodd_{\Spinor_{g}} \boxtimes \evenodd \boxtimes \DiracSt \right) 
  \end{split}\\
  &= ... = \degen{j+1} \degen{i}(P).
 \end{align*}
 In the firth line, we can drop ${\tau_{n+1,n+2}}_\ast$ in the second factor to go to the sixth line because 
 \begin{equation*}
     \DiracSt \boxtimes \id + \evenodd \boxtimes \DiracSt = \partial_{x_{n+1}}\cdot \frac{\partial}{\partial{x_{n+1}}} +  \partial_{x_{n+2}} \cdot \frac{\partial}{\partial{x_{n+2}}} 
 \end{equation*}
 is symmetric in the $(n+1)$-th and $(n+2)$-th coordinate.
\end{proof}

We come now to one of the most essential definitions of this thesis.
\begin{definition}
 Let $\InvBlockDirac$ 
 \nomenclature{$\InvBlockDirac$}{Cubical Set of invertible block Dirac operators}
 be the sequence of subsets of $\BlockDirac$, where $\InvBlockDirac[n]$ consists of all block Dirac operator $P \in \BlockDirac[n]$ that extend to an invertible operator $P \colon H^1(\Spinor_{g_P}) \rightarrow H^0(\Spinor_{g_P})$.
\end{definition}

The proof that the face maps of $\BlockDirac$ restricts to $\InvBlockDirac$ requires the following lemma that roughly says that positivity descends to hypersurfaces.

\begin{lemma}\label{PositivityDescendHyperSurf - Lemma}
 Let $\Spinor_g \rightarrow M$ be a spinor bundle over a  Riemannian manifold $(M,g)$.
 For $Q \in \PseudOpCl^1(E)$ symmetric and odd, define $P \deff Q \boxtimes \id + \evenodd \boxtimes \DiracSt$ on $\Spinor \boxtimes \Cliff[1][0] \rightarrow M \times (-a,a)$ for all $0<a\leq \infty$.
 
 If $||Pu||_0 \geq c||u||_0$ for all $u \in \Gamma_c(M \times (-a,a),\Spinor_g \boxtimes \Cliff[1][0])$, then 
 \begin{equation*}
     ||Qv||_0^2 \geq \max\{0, c-4/a^2\} ||v||_0^2
 \end{equation*}
 for all $v \in \Gamma_c(M,\Spinor_g)$. 
\end{lemma}
\begin{proof}
 Let $v \in \Gamma_c(M,\Spinor_g)$ and $f \in \Gamma_c((-a,a),\Cliff[1][0])$.
 The operator $P$ is again symmetric by the proof of Lemma \ref{ConnectingMapsWellDef - Lemma} and odd.
 Using Fubini's theorem and that $P$ is bounded from below, we get
 \begin{align*}
     ||Pv \otimes f||_0^2 &= \langle P^2(v \otimes f), v\otimes f\rangle \\
     &= \langle Q^2(v) \otimes f - v \otimes \Delta_{\R}(f), v \otimes f \rangle \\
     &= ||Q(v) \otimes f||_0^2 + ||v \otimes f'||_0^2 \\
     &= ||Q(v)||_0^2 ||f||_0^2 + ||v||_0^2 ||f'||_0^2 \\
     &\geq c ||v||_0^2||f||_0^2 = c||v \otimes f||_0^2.
 \end{align*}
 The function 
 \begin{equation*}
     \varphi \colon (-1,1) \rightarrow \R \quad \text{ given by } \quad \varphi(t) = \begin{cases}
       1 + t, & \text{if } t \leq 0, \\
       1 - t, & \text{if } t \geq 0,
     \end{cases}
 \end{equation*}
 satisfies $||\varphi||_0^2 = 2/3$ while $||\varphi'||_0^2 = 2$.
 Pick a compactly supported function $f$ that is sufficiently $||\placeholder||_1$-close to $\varphi$ such that
 $||f'||_0^2/||f||_0^2 \leq 4$.
 
 By the transformation formula, the assignment $f \mapsto a^{-1/2}f(a^{-1}\cdot)$ is an isometry $L^2(-1,1) \rightarrow L^2(-a,a)$.
 Thus, $||f_a'||_0^2/||f_a||_0^2 = a^{-2}||f'||_0^2/||f||_0^2\leq 4/a^2$.
 
 In conclusion, if we plug into $f_a$ in the previous inequality, we obtain
 \begin{align*}
     ||Qv||_0^2 &\geq \left(c - ||f_a'||_0^2/||f_a||_0^2\right) ||v||_0^2 \\
     &\geq \left(c - 4/a^2\right) ||v||_0^2,
 \end{align*}
 which is positive if $c > 4/a^2$.
\end{proof}

\begin{lemma}
 The connecting maps of $\BlockDirac$ send invertible elements to invertible elements.
 Thus, $\InvBlockDirac$ is a cubical subset.
\end{lemma}
 \begin{proof}
  Fix a sufficiently large $R > 0$ such that $g = g_P$ decomposes into $g(i,\eps) \oplus \diff x_i^2$ on $U_R{(i,\eps)}$ and, by abusing notation, $P$ decomposes into $P(i,\eps)\boxtimes \id + \evenodd \boxtimes \DiracSt_{\R}$.
  Since $P$ is bounded from below, its restriction $P|_{U_R(i,\eps)}$ is also bounded from below.
  Lemma \ref{PositivityDescendHyperSurf - Lemma} implies that $P(i,\eps)$, and hence $\face{i}P$ too, is bounded from below.
  Since $\face{i}P$ is self-adjoint, it must be invertible.
 
To see that $\degen{i}(P)$ is invertible, it suffices to check that
$P \boxtimes \id + \evenodd \boxtimes \DiracSt$ is invertible.
Since $P \boxtimes \id$ and $\evenodd \boxtimes \DiracSt$ anti-commute we have 
\begin{equation*}
 \left(P \boxtimes \id + \evenodd \boxtimes \DiracSt\right)^2 = P^2 \boxtimes \id - \id \boxtimes \Delta.
\end{equation*}
Both summands are non-negative operators. Since $P$ is invertible, $P^2$ is positive and therefore the sum is positive and therefore invertible. 
\end{proof}

We would like to know whether $\BlockDirac$ and $\InvBlockDirac$ are Kan sets.
For $\BlockDirac$ even more is true.

\begin{proposition}\label{BlockDiracConbContr - Prop}
  The cubical set $\BlockDirac$ is combinatorially contractible.
\end{proposition}

The proposition follows immediately from the following elementary lemma.

\begin{lemma}\label{Operator Suspension Union - Lemma}
 Let $P_{(j,\omega)}$ and $P_{(k,\eta)}$ be two block Dirac operators on $M \times \R^{n-1}$. 
 Assume that $j \leq k$, that $\face[\omega]{j}P_{(k,\eta)} = \face[\eta]{k-1}P_{(j,\omega)}$ and that the two operators decompose outside of $M \times \rho I^n$. 
 Then the operators $\degen{j}P_{(j,\omega)}$ and $\degen{k}P_{(k,\eta)}$ agree on the set $\{\omega x_j > \rho, \eta x_k > \rho\}$. 
\end{lemma}
\begin{proof}
 Abbreviate $\{\omega x_j > \rho, \eta x_k > \rho\}$ to $U$.
 If $j = k$ there is nothing to prove for $U$ is either empty or $(j,\omega) = (k,\eta)$. 
 In the other cases, it follows from Lemma \ref{Metric Suspension Union - Lemma} that the underlying metrics agree on $U$.
 The proof of Lemma \ref{BlockSuspends - Lemma} shows that $\sigma_jP_{(j,\omega)}$ restricts to $U$.
 Furthermore, for $j<k$, the calculation
 \begin{align*}
  \degen{j}P_{(j,\omega)}|_U &= \degen{k} \face[\eta]{k} \degen{j} P_{(j,\omega)} |_U \\
  &= \degen{j} \degen{k-1} \face[\eta]{k-1} P_{(j,\omega)} |_U\\
  &= \degen{j} \degen{k-1} \face[\omega]{j}P_{(k,\eta)} |_U \\
  &= \degen{j} \face[\omega]{j} \degen{k} P_{(k,\eta)} |_U = \degen{k} P_{(k,\eta)}|_U
 \end{align*}
 shows the claim. 
 Here, we used the block form of $P_{(j,\omega)}$ in the first line, cubical identities in the second line, and the assumed compatibility in the third line.
\end{proof}

\begin{proof}[Proof of Proposition \ref{BlockDiracConbContr - Prop}]
 Let $\partial \StandCube{n} \rightarrow \BlockDirac$ be given by the family of block Dirac operators 
 \begin{equation*}
  \{P_{(j,\omega)} \in \BlockDirac[n-1] \, : \, \face[\omega]{j} P_{(k,\eta)} = \face[k-1]{\eta} P_{(j,\omega)} \text{ for } j<k\},
 \end{equation*}
 where $(j,\omega) \in \{1,\dots, n\} \times \Z_2$.
 Since we are given only finitely many block Dirac operators, we can pick a sufficiently large $\rho>0$ such that the core of each $P_{(j,\omega)}$ is a subset of $M \times \rho I^{n-1}$.
 All $\degen{j}P_{(j,\omega)}$ are again block Dirac operators whose cores are contained in $M \times \rho I^n$, so they restrict to $\{\omega x_j > \rho\}$ and all of their intersections.
 The compatibility assumption yields that the operators agree on the intersections by Lemma \ref{Operator Suspension Union - Lemma}.
 
 Thus, these block Dirac operators can be patched together with a suitable choice of the partition of unity from Definition \ref{BlockPartOfUnity - Def} to an operator $P$ on $\Spinor_g \rightarrow M \times \R^n\setminus \rho I^n$, where $g$ is the union of all (compatible) underlying Riemannian metrics $g_{P_{(j,\omega)}}$.
 Extend the metric $g$ to a block metric on $M \times \R^n$ as in the proof of Proposition \ref{BlockMetric Combi Contrac - Prop}
 and denote the extension again with $g$.
 
 Let $\chi \colon M \times \R^n \rightarrow [0,1]$ be a smooth function that is identically $1$ on an open neighbourhood of $M \times \rho I^n$ and identically zero on $M \times \R^n \setminus (\rho +2)I^n$. 
 Let $\psi \colon M \times \R^n \rightarrow [0,1]$ be a function that is identically $1$ on $\supp (1-\chi)$ and supported within the complement $M \times \rho I^n$. 
 The operator $\mathtt{P} \deff \psi P(1-\chi) +  \Dirac_g \chi$ is a block operator on $M \times R^n$ that satisfies $\face[\omega]{j} \mathtt{P} = P_{(j,\omega)}$ because it agrees with $P$ away from $M \times (\rho + 2)I^n$.
 Symbol calculus and the choice of $\chi$ and $\psi$ imply that the principal symbol of $\mathtt{P}$ agrees with the principal symbol of $\Dirac_g$.
 From $\mathtt{P} - \Dirac_g = \psi(P - \Dirac_g)(1-\chi)$ follows that $P$ differs from $\Dirac_g$ by a bounded operator $H^0$. 
 If $\mathtt{P}$ were symmetric, then it would be a block Dirac operator and a filler.
 Although $\mathtt{P}$ might not be symmetric, its formal adjoint exists and is given by
 \begin{equation*}
     \mathtt{P}^\ast = \chi\Dirac_g   + (1 - \chi)P\psi. 
 \end{equation*}
 It is again a block operator of Dirac type.
 Thus a filler is given by $P_{\mathrm{fill}} \deff (\mathtt{P}^\ast + \mathtt{P})/2$, so $\BlockDirac$ is combinatorially contractible. 
\end{proof}

We now describe the comparison map $\PseudDir[\bullet] \rightarrow \BlockDirac$.
To this end, recall that for a block map $P \colon \R^n \rightarrow \PseudDir$ with underlying block map of Riemannian metrics $h \colon \R^n \rightarrow \Riem(M)$, the spinor bundle decomposes $\Spinor_{\susp\, h} \cong \Spinor_h \otimes \Cliff[n][0]$ and that we have defined the \emph{operator suspension} $\susp(P)$ as the operator that corresponds under the isomorphism $\Gamma_c(\Spinor_{\susp\, h}) \cong \Gamma_c(\Spinor_h \otimes \Cliff[n][0])$ to 
\begin{equation*}
    \susp(P) = P \otimes \id + \sum_{k=1}^n \partial_{t_k} \cdot \nabla_{\partial_{t_k}}^{\Spinor_{\susp \, h}}.
\end{equation*}

\begin{lemma}
 The suspensions are well defined, affine maps and assemble to a cubical map $\PseudDir[\bullet] \rightarrow \BlockDirac$
\end{lemma}
\begin{proof}
 
 We already know from Lemma \ref{SuspensionBlockOperator  - Lemma} that $\susp(P)$ is a block operator of Dirac type.
 In remains to show that $\susp(P)$ is symmetric and differs from $\Dirac_{\susp\, h}$ by a bounded operator.
 In the following, we will abbreviate $\susp \, h$ with $g$.
 
 To show that $\susp(P)$ is symmetric, we will show that the summands $P^\extension = P \otimes \id$ and $D \deff \sum \partial_{k} \cdot \nabla_{k}^{\Spinor_{g}}$ are symmetric operators.
 It is enough to consider compactly supported sections of the form $\sigma(m,t) \otimes e_I$. 
 Using Fubinis theorem we derive symmetry of $P^\extension$ by the following lines:
 \begin{align*}
  \left( P^\extension(\sigma_1 \otimes e_I), \sigma_2\otimes e_J\right) &= \int_{\R^n} \left( P_t\sigma_1(\cdot,t)\otimes e_I, \sigma_2(\cdot,t)\otimes e_J\right)_{\Sobolev{0}{M,\Spinor_{h(t)}} \otimes \Cliff} \diff t \\
  &= \int_{\R^n}  \left( P_t\sigma_1(\cdot,t) e_I, \sigma_2(\cdot,t)\right)_{\Sobolev{0}{M,\Spinor_{h(t)}}} \langle e_I, e_J\rangle_{\Cliff} \diff t \\
  &= \int_{\R^n} \left( \sigma_1(\cdot,t), P_t\sigma_2(\cdot,t)\right)_{\Sobolev{0}{M,\Spinor_{h(t)}}} \langle e_I, e_J\rangle_{\Cliff} \diff t \\
  &= \dots = \left( \sigma_1 \otimes e_I, P^\extension(\sigma_2 \otimes e_J)\right). 
 \end{align*}
 The operator $D$ is symmetric as it is a $\Spinor_g$-twisted Dirac operator on $\Cliff[n][0] \rightarrow M \times \R^n$ in the sense of \cite{LawsonMichelsonSpin}.
 
 Symbol calculus and Theorem \ref{App: ExternalTensor is Continuous - Theorem} show that
 \begin{equation*}
     \susp(P) - \Dirac_g = (P - \Dirac_h) \otimes \id \in Op^0(H^0(\Spinor_g)))
 \end{equation*}
 because $P - \Dirac_h$ is a block map with values in pseudo differential operators of order zero.
 The operator $\susp(P) - \Dirac_g$ thus decomposes in the same manner as a block operator\footnote{The only thing that prevents the operator to be a block operator is that it only lives in $Op^0$ instead of $\PseudOpCl^0$.}. 
 The proof of Lemma \ref{BlcokOperatorsInduceBound - Lemma} does not require that the operator can be approximated by actual pseudo differential operators of order zero, so $\susp(P) - \Dirac_g$ induces a bounded endomorphism on $H^0(\Spinor_g)$.
 
 We have seen in the proof of Lemma \ref{SuspensionBlockOperator  - Lemma} that $P$ decomposes on $U_r(\fateps)$, provided $r$ is sufficiently large, as follows
 \begin{align*}
  (\Phi_\fateps)_\ast \susp(P)|_{U(\fateps)} &= \susp(P(\fateps)) \boxtimes \id + \evenodd \boxtimes \sum_{k=1}^{|\dom \fateps|} \left((\Phi_\fateps)_\ast \nabla^{\Spinor_g}\right)_{\partial_{t_k}} \\
  &= \susp(P(\fateps)) \boxtimes \id + \evenodd \boxtimes \sum_{k=1}^{|\dom \fateps|} \partial_{t_k} \cdot \frac{\partial}{\partial t_k},
 \end{align*}
 which implies $\face{i}\susp(P) =\susp(\face{i}P)$.
 
 Using the isomorphism $\Spinor_{\degen{n+1}g} \cong \Spinor_g \boxtimes \Cliff[1]$ we can identify $\susp(\degen{n+1}P)$ with $\susp(P) \boxtimes \id + \evenodd \boxtimes \DiracSt = \degen{n+1} \susp(P)$, which shows that $\susp$ commutes with $\degen{n+1}$. It is a bit cumbersome but not hard to prove that 
 \begin{equation*}
  \Cycl(i,n)_\ast\susp(P) = \susp(\cycl(i,n)_\ast P) = \susp(P \circ \cycl(i,n)^{-1}) 
 \end{equation*}
 for every $n\geq 1$, which yields $\degen{i}\susp(P) = \susp(\degen{i}P)$.
\end{proof}

As in the case of psc metrics, it is not true that the operator suspension sends block maps of invertible operators to an invertible block operator. 
In order to correct this defect, we will construct a smaller model for $\InvPseudDir[\bullet]$ in the next section.

 \section{The Comparison Map}\label{Section - The Comparison Map}

The operator suspension of a smooth, invertible pseudo Dirac operator valued block map does not need to be invertible.
To repair this defect, we will construct a weakly equivalent cubical subset $B_\bullet$ of $\InvPseudDir[\bullet]$, which is Kan and such that the image of the operator suspension map lies in $\ConcSet$ and $\InvBlockDirac$.
Essentially, $\AuxInvDir$ is the subset of smooth block maps that are sufficiently ``slowly'' parameterised.

Before we start, recall that an inner product $h$ on a finite dimensional vector space $V$ induces an isomorphism $h^\flat \colon V \rightarrow V^\vee$ via $v \mapsto h(v,\placeholder)$ whose inverse is denoted by $h_\sharp$. 
If $b$ is another bilinear form on $V$, we denote its induced endomorphism $h_\sharp \circ b$ with $b^\op$.

For each block map of Riemannian metrics $h \colon \R^n \rightarrow \Riem(M)$, the spinor bundle decomposes into $\Spinor_{\susp \, h} \cong \Spinor_h \otimes \Cliff[n][0]$, where $\Spinor_h \rightarrow M \times \R^n$.
We need to properly investigate the spinor connection of a suspension metric $\susp \, h$.
We denote the standard basis on $\R^n$ with $\partial_j = \partial_{t_j}$.
\todo[inline]{Rewrite the beginning. See Thomas remarks!}

\begin{lemma}\label{CalcChristoffelSymb - Lemma}
 Let $g_{0} \in \Riem(M)$ be a Riemannian metric. 
 Denote the constant block map with values $g_0$ and its suspension with the same letter.
 For each block map of Riemannian metrics $h \colon \R^n \rightarrow \Riem(M)$, the spinor connection on $\Spinor_{\susp \, h}$ satisfies
 \begin{equation*}
     \preGauge(\susp\, h,g_0)_\ast(\nabla^{\Spinor_{\susp\, h}})_{\partial_i} = \partial_i
 \end{equation*}
 for all $1 \leq i \leq n$.
 
 Informally speaking, the Christoffel endomorphism $\Gamma_i = \partial_i - \nabla^{\Spinor_{\susp \, h}}_{\partial_i}$ vanishes for all $1 \leq i \leq n$.
\end{lemma}
\begin{proof}
 Write $g$ for $\susp \, h$, abuse notation with $g_0 \deff \susp \, g_0$ and let $\preGauge(g,g_0) \colon \Spinor_g \rightarrow \Spinor_{g_0}$ be the pre-gauge map.
 The push-forward $\preGauge(g,g_0)_\ast(\nabla^{\Spinor_g})$ is a connection $\Spinor_{g_0}$.
 Since $g_0 = g_{0} \oplus \euclmetric_{\R^n}$ is a product metric, the space of section $\Gamma(M\times \R^n,\Spinor_{g_0})$ can be identified with $\mathcal{C}^\infty\bigl(\R^n,\Gamma(M,\Spinor_{g_{0}}\otimes \Cliff[n][0])\bigr)$, so that
 all $\partial_{i}$ are well defined differential operators.
 Thus, the difference
 \begin{equation*}
     \Gamma_i \deff \partial_i - \preGauge(\susp\, h,g_0)_\ast(\nabla^{\Spinor_\susp{\, h}})_{\partial_i}
 \end{equation*}
 is a bundle endomorphism of $\Spinor_{g_0}$, the Christoffel endomorphism for $\preGauge(g,g_0)_\ast (\nabla^{\Spinor_{g}})_{\partial_i}$. 
 This means, that if we pick an orthogonal frame of $TM$ complementing the standard basis of $\R^n$, then the matrix representation of that endomorphism with respect to the chosen frame yield the Christoffel symbols of $\preGauge(g,g_0)_\ast(\nabla^{\Spinor_g})$.
 (Of course, the identity $\Gamma_i=0$ is independent of the chosen frame on $TM$.)
 
 The Christoffel symbols for $\preGauge(g,g_0)_\ast(\nabla^{\Spinor_g})_{\partial_i}$ with respect to a frame induced by $(e_{1-d},\dots,e_0$, $\partial_1 \dots, \partial_n)$ agree with the Christoffel symbols for $\nabla^{\Spinor_g}_{\partial_i}$ with respect to the frame induced by $\preGauge(g_0,g)(e_{(1-d)}), \dots,\preGauge(g_0,g)(e_0)$, $\partial_1,\dots,\partial_n$.
 
 The Christoffel symbols of the spinor connection with respect to an arbitrary orthonormal frame $(e_i)$ on the manifold and the induced orthonormal frame $(\sigma_\alpha)$ on the spinor bundle is given by the formula \cite{LawsonMichelsonSpin}*{Theorem II.4.14}
 \begin{equation*}
     \nabla_{e_i}^{\Spinor_g} \sigma_\alpha = \frac{1}{2}\left(\sum_{j<k}\Gamma_{ij}^k e_je_k\right) \cdot \sigma_\alpha,
 \end{equation*}
 where $\Gamma_{ij}^k$ are the Christoffel symbols of the Levi-Cevita connection with respect to the orthonormal frame $(e_i)$.
 
 To determine $\Gamma_i = \Gamma^{L.C.,g}_i$ at $(m_0,t_0)$, the Christoffel endomorphism of the Levi-Cevita connection on $M_{t_0} = M \times \{t_0\}$, we will choose a local orthonormal frame that is particularly adapted to this point.
 Let $(y_{1-d},\dots,y_0) \colon U \subset M \rightarrow V \subseteq \R^d$ be geodesic coordinates for $h(t_0)$ centered at $m_0$.
 These coordinates map $m_0$ to $0$.
 Let $\partial_{y_{(1-d)}}, \dots, \partial_{y_0}$ be the standard basis. 
 The push-forward of $g = \susp(h)$ with respect to this coordinate system is uniquely determined by the function $h \colon V \rightarrow \R^{d\times d}$, which is given by $h_{\alpha \beta}(v) = g(v)(\partial_{y_\alpha}, \partial_{y_\beta})$.
 Observe that $h(0) = \id$ and $dh = 0$. 
 Applying the Gram-Schmidt algorithm to $(\partial_{y_{-(d-1)}}, \dots , \partial_{y_0})$ with respect to $h(t_0)$, we get an orthonormal frame $f_{-(d-1)}, \dots f_0$ on $TM \times \{t_0\}$.
 In geodesic coordinates, the metric satisfies $h(v,t_0) = \id + \mathcal{O}(||v||^2)$ around the origin.
 Induction over the Gram-Schmidt algorithm shows the orthonormal basis $(f_j)$ is tangent to $\partial_{y_j}$ at the origin, i.e., $f_j(v) - \partial_{y_j} = 0 + \mathcal{O}(||v||^2)$.

  There are unique smooth functions $a_{\alpha,j} \colon V  \rightarrow \R$ such that
 \begin{equation*}
     f_j(v) = \sum_{\alpha=1-d}^0 a_{\alpha,j}(v)\partial_{\alpha}.
 \end{equation*}
 Since the $(f_j)$ and $(\partial_{y_j})$ are tangent to each other, we have $a_{\alpha,j} = \delta_{\alpha,j}$ and $\diff a_{j,\alpha}=0$.
 
 In these local coordinates, we can extend the frame $f_{1-d},\dots,f_{0}$ from $V = V\times{t_0}$ to a frame on $V \times \R^n$ via 
 \begin{equation*}
     f_j(v,t) \deff h(v,t)^{-1/2} \cdot f_j(v,t_0).
 \end{equation*}
 (To be precise, $h(v,t)$ denotes the push-forward of the metric $h(t)|_U = g|_{TU \times \{t\}}$ with respect to the chosen geodesic coordinates $(y_{(1-d)},\dots, y_0$)).
 
 To calculate the Christoffel symbols of the Levi-Cevita connection of $g = \susp(h)$ with respect to the frame $(e_{-(d-1)}, \dots e_n)=$ $(f_{-(d-1)}, \dots, f_0 , \partial_{t_1}, \dots, \partial_{t_n})$ we use the Koszul formula for the Levi-Cevita connection
\begin{equation*}
\begin{split}
 2g(\nabla_XY,Z) = \  &Xg(Y,Z) + Yg(Z,X) - Zg(X,Y) \\
                  &-g(Y,[X,Z]) - g(Z,[Y,X]) + g(X,[Z,Y]).
\end{split}
\end{equation*}
Let $\Gamma_{ij}^k = g(\nabla_{e_i}e_j,e_k)$. 
We have chosen an orthonormal frame, so $g(e_i,e_j) \equiv \delta_{ij}$ and the upper line in the Koszul formula vanishes identically.
Thus, only the lower line contributes.
 
 The commutators at $(0,t_0)$ are given by:
\begin{align*}
 [f_i,f_j] &= [h(t)^{-1/2}(\sum a_{\alpha,i}\partial_{y_\alpha}), h(t)^{-1/2}(\sum a_{\beta,j}\partial_{y_\beta}) ] \\
 &= \sum_{\alpha, \beta} [ h(t)^{-1/2}( a_{\alpha,i}\partial_{y_\alpha}), h(t)^{-1/2}( a_{\beta,j}\partial_{y_\beta})] \\
 &= \sum_{\alpha, \beta} a_{\beta,j}\partial_{y_j}a_{\alpha,i} - a_{\alpha,i}\partial_{y_i} a_{\beta,j} + a_{\alpha,i}a_{\beta,j} [ h(t)^{-1/2} \partial_{y_\alpha}, h(t)^{-1/2}\partial_{y_\beta}] \\
 &= \sum_{\alpha, \beta}-1/2 \cdot a_{\alpha,i}a_{\beta,j}(0,t_0)\bigl( \partial_{y_\alpha} h(0,t_0) \cdot \partial_{y_\beta} - \partial_{y_\beta} h(0,t_0) \cdot \partial_{y_\alpha}\bigr)=0, 
\end{align*}
\begin{align*}
 [\partial_{t_i}, f_j] &= \sum_\alpha [\partial_{t_i}, h^{-1/2}(t)_v a_{\alpha, j}(v)\partial_{y_\alpha}] \\
 &= \sum_\alpha \partial_{t_i}a_{\alpha, j}(v) + a_{\alpha, j}[\partial_{t_i}, h^{-1/2}(t)_v \partial_{y_\alpha}] \\
 &= 0 + \sum_{\alpha} a_{\alpha, j}(0,t_0) \cdot (-1/2 \partial_{t_i}h(t)|_{(0,t_0)} \cdot \partial_{y_\alpha}) \\
 &= -1/2 (\partial_{t_i}h)(0,t_0) \partial_{y_j},
\end{align*}
and, of course,
\begin{equation*}
 [\partial_{t_i}, \partial_{t_j}] = 0.
\end{equation*}

We can now calculate the Christoffel symbols $\Gamma_{ij}^k$ at $(0,t_0)$.
If two of the three indices $i,j,k$ are bigger than zero, then $\Gamma_{ij}^k =0$ because either the commutator vanishes, or the commutator is orthogonal to all $\partial_{t_m}$.
We deal with the remaining four cases individually:
For $i,j,k \leq 0$ we have $2\Gamma_{ij}^k = 0$, because all commutators vanish.
For $i,j \leq0$ and $k>0$, we have
\begin{align*}
  2\Gamma_{ij}^k &= g(e_i, [e_k,e_j]) - g(e_j,[e_i,e_k]) = g(f_i,[\partial_{t_k},f_j]) - g(f_j,[f_i,\partial_{t_k}]) \\
  &= g(\partial_{y_i}, -1/2\partial_{t_k}h \cdot \partial_{y_j}) + g(\partial_{y_j}, -1/2\partial_{t_k}h \cdot \partial_{y_i}) \\
  &= -1/2 \partial_{t_k} h_{ij} -1/2\partial_{t_k}h_{ji} = - \partial_{t_k}h_{ij}, 
\end{align*}  
where we used that $f_\alpha(0,t_0) = \partial_{y_\alpha}$ in the second line.
For $i,k \leq 0$ and $j > 0$ we derive analogously $2\Gamma_{ij}^k = \partial_{t_j}h_{ik}$.
For $j,k \leq 0$ and $i>0$ we have: $2\Gamma_{ij}^k = 1/2\partial_{t_i}h_{jk} - (1/2)\partial_{t_i}h_{kj} = 0$.

This implies for the Christoffel symbol of the spinor connection with respect to the frame $(e_{-(d-1)}, \dots e_n)=$ $(f_{-(d-1)}, \dots, f_0 , \partial_{t_1}, \dots, \partial_{t_n})$ in the case $i>0$ the identity
\begin{align*}
    \Gamma_{\partial_i}^{\Spinor_g} &= \Gamma_i^{\Spinor_g} = \frac{1}{2}\sum_{1-d \leq j < k \leq n} \Gamma_{ij}^ke_je_k = 0,
\end{align*}
and the lemma is proven.
\end{proof}

Let $P \colon \R^n \rightarrow \PseudOp^1(\Spinor)$ be a block map with underlying block map $h$.
Recall that $P^\extension$ is the operator which corresponds under the decomposition $\Spinor_{\susp \, h} \cong \Spinor_h \otimes \Cliff[n][0]$ to $P \otimes \id$.
The commutators $\lbrack P^\extension, \nabla_{\partial_{t_k}}^{\Spinor_{\susp h}}\rbrack$ is $\mathcal{C}^\infty_c(\R^n,\R)$-linear, so it restricts to an operator on $\Gamma(M\times t, \Spinor_{\susp h}) = \Gamma(M\times t,\Spinor_{h_t} \otimes \Cliff)$ that we denote with  $\lbrack P^\extension, \nabla_{\partial_{t_k}}^{\Spinor_{\susp h}}\rbrack\restrict_{M \times t}$. 

\begin{lemma}\label{CommutatorWithExtension - Lemma}
 The assignment 
 \begin{equation*}
     t \mapsto [P^\extension,\nabla^{\Spinor_{\susp \, h}}_{\partial_j}] \restrict_{M\times t}
 \end{equation*}
 is the extension of a smooth block map $\R^n \rightarrow \PseudOp^1(\Spinor)$ for all $1 \leq j \leq n$.
\end{lemma}
\begin{proof}
 Write $g$ for $\susp \, h$. 
 Pick a Riemannian metric $g_{00} \in \Riem(M)$ and denote the constant map with value $g_{00}$ by the same letter. 
 Set $g_0 = \susp(g_{00}) = g_{00} \oplus \euclmetric_{\R^n}$.
 
 We prove the equivalent statement that
 \begin{equation*}
     t \mapsto \preGauge(g,g_0)_\ast\left([P,\nabla^{\Spinor_{\susp \, h}}_{\partial_j}]\right) \restrict_{M\times t}
 \end{equation*}
 is a smooth block map $\R^n \rightarrow \PseudOp^1(\Spinor_{g_0})$.
 
 Note that
 \begin{equation*}
     \preGauge(g,g_0)_\ast P^{\extension} = \left(\preGauge(h,g_{00})_\ast P\right)^\extension
 \end{equation*}
 because both operators agree on the dense subspace spanned by sections of the form $u \otimes f$, with $u \in \Gamma(\Spinor_{g_{00}})$ and $f \in \Gamma_c(\Cliff[n][0])$.
 Denote $\preGauge(h,g_{00})_\ast P$ with $Q$.
 Lemma \ref{CalcChristoffelSymb - Lemma} implies
 \begin{equation*}
     [Q^\extension, \preGauge(g,g_0)_\ast\left(\nabla^{\Spinor_g}\right)_{\partial_j}]\restrict_{M \times t} = (\partial_j Q)^\extension\restrict_{M \times t},
 \end{equation*}
 where $\partial_j Q$ denotes the partial differential of the smooth block map $Q\colon \R^n \rightarrow \PseudOp^1(\Spinor_{g_{00}})$.
 Since the partial differential of a smooth block map is again a smooth block map, the claim follows.
\end{proof}

From now on, we will not notationally distinguish $\lbrack P^\extension, \nabla_{\partial_{t_k}}^{\Spinor_{\susp h}}\rbrack$ from its restriction to $M \times \{t\}$ if it does not lead to confusion.
Also, we will write $M_t$ for $M \times \{t\}$.

\begin{definition}\label{AuxInvDir - Definition}
 Let $\AuxInvDir$ be the following cubical subset of $\InvPseudDir[\bullet]$:
 A block map $P \in \InvPseudDir[n]$ with underlying block map of metrics $h$ belongs to $\AuxInvDir[n]$ if the following conditions are satisfied
 \begin{align}
   &\sum_{j=1}^n ||(\partial_k h)^\op||_{\mathrm{op},h_t} < \min \{d^{-1}, 1/8 \cdot 2^{-(d+1)} \lowBnd(t)\}, \label{eq: AuxDir1} \nomenclature{$(\partial_k h)^\op$}{Operator version of $\partial_k h$ that is ${h_t}_\sharp \circ \partial_k h$}\\
   &\sum_{j=1}^n || \lbrack P^\extension, \nabla^{\Spinor_{\susp h}}_{\partial_j}\rbrack||_{1,0;M_t} < 1/32 \label{eq: AuxDir2}\, \lowBnd(t)^2.
 \end{align}  
 Here, $d = \dim M$ and $\lowBnd(t) = \lowBnd_P(t) \deff \inf \{ ||P_tv||_{0,h_t} \, : \, ||v||_{1,h_t} = 1\}$.
 The function $\lowBnd_P$ is uniformly positive because $t \mapsto P_t$ is a block map with values in invertible operators. 
\end{definition}

\begin{lemma}
 $\AuxInvDir$ is a cubical subset of $\InvPseudDir[\bullet]$.
\end{lemma}
\begin{proof}
 Given $U = \R_{<0}$ and $X = \InvPseudDir[ ]$, which is an open subset of the affine Fréchet space $\PseudDir[ ]$. 
 Then 
 \begin{equation*}
  f_n(P) \deff \sum_{k=1}^n ||(\partial_{t_k}h)^\op||_{\op,h_t} - \min \{d^{-1}, 1/8 \cdot 2^{-(d+1)} \lowBnd(t)\} 
 \end{equation*}
 and 
 \begin{equation*}
  g_n(P) \deff \sum_{k=1}^n ||\lbrack P^\extension, \nabla_{\partial_k}^{\Spinor_{\susp h}}\rbrack||_{1,0;M_t} - 1/32 \, \lowBnd(t)^2
 \end{equation*}
 are local and stable sequences of continuous maps $\mathcal{C}^\infty(\R^n,X) \rightarrow \mathcal{C}^0(\R^n,\R)$. 
 By definition, $\AuxInvDir[n] = f_n^{-1}(U) \cap g_n^{-1}(U)$. 
 Since intersections of two cubical subsets are again cubical subsets, the claim follows from Lemma \ref{local stable criterium - Lemma}.
\end{proof}

Next we will show that $\susp|_{\AuxInvDir}$ takes values in $\InvBlockDirac$ by showing that $\susp(P)^2$ is invertible for all $P \in \AuxInvDir[n]$. 

\begin{lemma}\label{susp quad - Lemma}
 Let $P \colon \R^n \rightarrow \PseudDir^\times$ be a smooth block map and set $g \deff \susp(h)$. 
 For each $u \in \Gamma_c(M \times \R^n,\Spinor_h)$, we have 
 \begin{align*}
   \susp(P)^2(u) = (P^2)^\extension(u) &- \sum_{j=1}^n \partial_j \cdot \lbrack P^{\extension} , \nabla^{\Spinor_g}_{\partial_{j}} \rbrack (u) \\ 
   &- \sum_{j=1}^n \left(\nabla_{\partial_j}^{\Spinor_g}\right)^2(u) + \sum_{1 \leq j < k\leq n}\mathfrak{R}^{\Spinor_g}(\partial_j,\partial_k)(u). 
 \end{align*}
\end{lemma}
\begin{proof}
 We abbreviate $\nabla_{\partial_j}^{\Spinor_g}$ to $\nabla_j$ and $\sum_{j=1}^n \partial_j \cdot \nabla_j$ to $\Dirac$. 
 Clearly,
 \begin{equation*}
  \susp(P)^2 = (P^\extension)^2 + P^\extension \Dirac + \Dirac P^\extension + \Dirac^2
 \end{equation*}
 and $(P^2)^\extension = (P^\extension)^2$. 
 Recall that $P^\extension$ is an odd operator. 
 This yields
 \begin{align*}
  P^\extension \Dirac u &= \sum_{j=1}^n P^\extension(\partial_j \cdot \nabla_j u) \\
                        &= \sum_{j=1}^n (-1) \partial_j \cdot P^\extension(\nabla_j u)
 \end{align*}
 implying
 \begin{align*}
  P^\extension \Dirac u + \Dirac P^\extension u &= \sum_{j=1}^n \partial_j \cdot \left( - P^\extension \nabla_j u + \nabla_j P^\extension u\right) \\
  &= -\sum_{j=1}^n \partial_j \cdot \lbrack P^\extension, \nabla_j \rbrack u.
 \end{align*}
 The last summand is given by
 \begin{align*}
   \Dirac^2 &= \sum_{j=1}^n\sum_{k=1}^n \partial_j\cdot \nabla_j(\partial_k \cdot \nabla_k u) \\
            &= \sum_{j,k} \partial_j\partial_k \cdot \nabla_j\nabla_k u + \partial_j \cdot \underbrace{\nabla_j\partial_k}_{=0} \cdot \nabla u \\
            &= \sum_{j=1}^n \partial_j^2 \cdot \nabla_j \nabla_j u + \sum_{j<k} \partial_j \partial_k \cdot \left( \nabla_j \nabla_k - \nabla_k \nabla_j\right)u \\
            &= \sum_{j=1}^n -\nabla_j^2 u + \sum_{j < k} \partial_j \partial_k \cdot R^{\nabla^{\Spinor_g}}(\partial_j,\partial_k)u \\
            &= \sum_{j=1}^n -\nabla_j^2 u + \sum_{j < k}\mathfrak{R}(\partial_j, \partial_k)u.
 \end{align*} 
 Putting everything together gives the claim. 
\end{proof}

\begin{theorem}\label{OperatorSuspensionFaktorisation - Theorem}
 For each $P \in \AuxInvDir[n]$, the operator $\susp(P)$ is invertible. Consequently, the diagram 
 \begin{equation*}
  \xymatrix@C+1em{\AuxInvDir \ar@{^{(}->}[r] \ar[d]_\susp & \InvPseudDir[\bullet] \ar[d]^{\susp} \\ \InvBlockDirac \ar@{^{(}->}[r] & \BlockDirac}
 \end{equation*}
 commutes.
\end{theorem}

The strategy of the proof is actually quite easy. By Lemma \ref{susp quad - Lemma}, $\susp(P)^2$ has a strictly positive summand, a non-negative one, and two summands we need to control. 
We will show that the conditions we have imposed on $\AuxInvDir$ guarantee that the two remaining terms have small operator norm so that $\susp(P)^2$ remains strictly positive.

In order to prove this we need to derive several analytical results that allow us to control the not necessarily non-negative summands of $\susp(P)^2$.

\begin{lemma}
 Let $U \subseteq \R^n$ be open. For a smooth map $g \colon U \rightarrow \Riem(M)$ of Riemannian metrics, let $R$ be the Riemannian curvature tensor of $\susp (g)$. If $\partial_i = \partial_{t_i}$ is the standard basis of $\R^n$, then
 \begin{equation*}
  R(\partial_{i},\partial_j) = \frac{1}{4}\lbrack (\partial_jg)^\op, (\partial_i g)^\op\rbrack.
 \end{equation*}
\end{lemma}
\begin{proof}
 We prove this formula in local coordinates. 
 Pick a local holonomic frame of $M\times \R^n$ that extends $\partial_1, \dots, \partial_n$. This ``spacial'' frame is indexed by $1-d \leq j \leq 0$. Using Einstein's sum convention, the Christoffel symbols of the Levi-Cevita connection of $\susp(g)$ in these coordinates are given by
 \begin{align*}
  2\Gamma_{t_i j}^k &= \susp(g)^{k\alpha}(\susp(g)_{t_i\alpha,j} + \susp(g)_{\alpha j,t_i} - \susp(g)_{t_ij,\alpha}) \\
  &= \susp(g)^{k\alpha}\susp(g)_{\alpha j,t_i} \\
  &= \begin{cases}
     0 & \text{if } \mathrm{max}\{j,k\} > 0, \\
     \sum_{\alpha<0} g^{k\alpha}g_{j\alpha,t_i} & \text{if } \max\{j,k\}\leq 0.
  \end{cases} 
 \end{align*}
 Thus, we can identify the endomorphism $2\Gamma_{t_i}$ with $g^{-1}\partial_{t_i}g = (\partial_{t_i}g)^\op$. 
 Abbreviating $t_i$ and $t_j$ with $i$ and $j$, respectively, we get for the Riemannian curvature of $\susp(g)$ the following expression:
 \begin{align*}
  R(\partial_i,\partial_j) &= (\partial_i + \Gamma_i)(\partial_j + \Gamma_j) - (\partial_j + \Gamma_j)(\partial_i + \Gamma_i) \\
  &= \partial_i\Gamma_j - \partial_j\Gamma_i + \Gamma_i\Gamma_j - \Gamma_j\Gamma_i \\
  \begin{split}
   &= -1/2g^{-1}(\partial_ig)g^{-1}(\partial_jg) + 1/2g^{-1}(\partial_jg)g^{-1}\partial_i g \\
   &\qquad +1/4 g^{-1}(\partial_ig) g^{-1} (\partial_jg) - 1/4 g^{-1}(\partial_j g)g^{-1}(\partial_i g)
  \end{split}\\
  &= 1/4\lbrack g^{-1}\partial_jg,g^{-1}\partial_i g\rbrack = 1/4 \lbrack (\partial_jg)^\op, (\partial_ig)^\op\rbrack.
 \end{align*} 
 and the lemma is proven. 
\end{proof}

\begin{cor}\label{Spinor Curvature Norm Estimation - Cor}
 The operator norm of the bundle endomorphism $\mathfrak{R}_t(\partial_j,\partial_k)\in \mathrm{End}(\Spinor|_{M_t})$ is bounded by
 \begin{equation*}
  2||\mathfrak{R}_t(\partial_j,\partial_k)||_{\op,t} \leq d^2||(\partial_jg)^\op||_{\op,t}||(\partial_kg)^\op||_{\op,t}
\end{equation*}
 for every $t \in U$, where $d = \dim M$.
 This implies
 \begin{equation*}
  \sum_{1 \leq j < k \leq n}\left|\left|\mathfrak{R}_t(\partial_j,\partial_k)\right|\right|_{\op,t} \leq \frac{d^2}{2}\left(\sum_{j=1}^n ||(\partial_jg)^\op||_{\op,t} \right)^2.
 \end{equation*}
 Here, $||\placeholder||_{\op,t}$ denotes the operator norm on $\Spinor_{g(t)}$ on the left hand side and the operator norm on $TM_t$ with respect to $g(t)$ on the right hand side.
\end{cor}
\begin{proof}
 Pick a local orthonormal frame $e_{1-d},\dots,e_0$ on $TM$ that completes $\partial_1,\dots,\partial_n$ to a local orthonormal frame on $TM \oplus \R^n$
 Then, by \cite{LawsonMichelsonSpin}*{Theorem II.4.15}, the curvature operator of the spinor connection is given by
 \begin{equation*}
  \mathfrak{R}(\partial_i,\partial_j) = \frac{1}{2}\sum_{k<l\leq 0} \langle R_{\partial_i,\partial_j}(e_k),e_l\rangle e_ke_l\cdot.
\end{equation*}  
 This yields
 \begin{align*}
  2||\mathfrak{R}(\partial_i,\partial_j)||_{\op,t} &= \left|\left| \sum_{k<l\leq 0} \langle R_{\partial_i\partial_j}(e_k),e_l\rangle e_ke_l\right|\right|_{\op, t} \\
  &\leq \sum_{k<l\leq 0} |\langle R_{\partial_i \partial_j}(e_k),e_l\rangle|\underbrace{||e_ke_l \cdot||_{\op,t}}_{=1} \\
  &\leq \sum_{k<l\leq 0}||R_{\partial_i\partial_j}||_{\op, TM_t} \\
  &= \sum_{k<l\leq 0} 1/4\cdot || \lbrack (\partial_jg)^\op, (\partial_i g)^\op\rbrack||_{\op, t} \\
  & \leq 2d^2||(\partial_ig)^\op||_{\op,t}||(\partial_j g)^\op||_{\op,t}.
 \end{align*}
\end{proof}

Since $\Spinor_{\susp \, h}|_{M_t} \cong \Spinor_{h(t)} \otimes \Cliff$ we have two canonical connections on this restriction, namely $\nabla^{\Spinor_{\susp h}|_{M_t}}$ and $\nabla^{\Spinor_{h(t)} \otimes \id}$. 
We denote the first connection with $\nabla^\susp$ and the second with $\nabla \otimes \id$. 
Of course, since $M$ is compact, these two connections yield equivalent Sobolev norms $||\placeholder||_{1,\nabla \otimes \id}$ and $||\placeholder||_{1,\nabla^\susp}$ on $\Spinor_{\susp h}|_{M_t}$. The next lemma helps to find explicit constants.
\begin{lemma}\label{TensorSuspNormEqui - Lemma}
 The bundle endomorphism valued 1-form 
 \begin{equation*}
  \mathcal{E} \deff \nabla \otimes \id - \nabla^\susp \colon TM \rightarrow \mathrm{End}\left(\Spinor_{\susp\,h}|_{M_t}\right)
\end{equation*}  
induces a bounded operator $L^2(M_t;TM_t) \rightarrow L^2(M_t;\mathrm{End} (\Spinor_{\susp \, h}|_{M_t})$ that satisfies the following operator norm estimation:
 \begin{equation*}
  ||\nabla \otimes \id - \nabla^\susp||_{0,0} \leq d/4 \sum_{k=1}^n ||\partial_{t_k}h^\op||_{\op,h(t)}. 
 \end{equation*}  
\end{lemma}
\begin{proof}
 Let $e_{-(d-1)}, \dots, e_0$ be a local orthonormal frame of $M$ with respect to $h(t)$. 
 With respect to this frame, $\mathcal{E}$ can be expressed in terms of the Christoffel symbols of $\nabla \otimes \id$ and $\nabla^\susp$, which we denote as $\Gamma^h$ and $\Gamma^\susp$, respectively. 
 Let $\Gamma_{ij}^k$ denote the Christoffel symbols of $\susp(h)$ on $M \times \R^n$, and denote $\partial_{t_k}$ by $e_k$ for $1\leq k \leq n$.
 For all $1-d \leq i \leq 0$ we have by \cite{LawsonMichelsonSpin}*{Theorem II.4.14} and the proof of Lemma \ref{CalcChristoffelSymb - Lemma} the following expression
 \begin{align*}
  \mathcal{E}(e_i) &= \Gamma_i^{h} - \Gamma^\susp_i = \frac{1}{2}\left(\sum_{j<k\leq 0} \Gamma_{ij}^k e_je_k - \sum_{-(d-1)\leq j<k\leq n} \Gamma_{ij}^k e_j e_k\right) \\
  &= \frac{1}{2} \sum_{j < k, \, k>0} - \Gamma^k_{ij}e_je_k \\
  &= \frac{1}{4}\sum_{j\leq 0 < k}  \partial_{t_k}h_{ij}e_j \partial_{t_k}  = \frac{1}{4}\sum_{k>0} (\partial_{t_k}h)^\op(e_i)  \partial_{t_k}  .
 \end{align*}
 Using this identity, we can estimate the operator norm of $\mathcal{E}$ as follows:
 \begin{align*}
  4||\mathcal{E}||_{0,0} &\leq \sum_{i=1-d}^0 4||\mathcal{E}(e_i)||_{\mathrm{End}(\Spinor_{\susp(h)})}\\
   &\leq \sum_{i=1-d}^0\sum_{k=1}^n ||\partial_{t_k}h^\op(e_i)\partial_{t_k} \cdot ||_{\mathrm{End}(\Spinor_{\susp(h)})} \\
  &= \sum_{i=1-d}^0\sum_{k=1}^n ||\partial_{t_k}h^\op(e_i)||_{TM}\\
  &\leq d \sum_{k=1}^n ||\partial_{t_k}h^\op||_{\op,h(t)}.
 \end{align*}
 Here, $||\placeholder||_{\mathrm{End}(\Spinor_{\susp(h)})}$ denotes the usual $C^\ast$-algebra operator norm coming from the inner product on the fibres of $\Spinor_{\susp(h)}$ and $||\placeholder||_{TM_t}$ denotes the norm on $TM$ induced by $h(t)$.
 Note that we used in the third line that the endomorphim is given by the Clifford multiplication, so we could use $C^\ast$-identity and the Clifford identity.
\end{proof}
\begin{cor}\label{EquiSobolevNorms - Cor}
 If $||\placeholder||_{1,\nabla^\susp}$ and $||\placeholder||_{1,\nabla \otimes \id}$ denote the Sobolev norms on $\Spinor_{\susp h}|_{M_t}$ with respect to $\nabla^\susp$ and $\nabla \otimes \id$, respectively, then they are equivalent by the following constants:
 \begin{align*}
   ||\placeholder||_{1,\nabla^\susp} &\leq \left( 1 + d \sum_{k=1}^n ||\partial_{t_k} h||_{\op}\right) ||\placeholder||_{1,\nabla \otimes \id}\\
   & \leq \left( 1 + d \sum_{k=1}^n ||\partial_{t_k} h||_{\op}\right)^2 ||\placeholder||_{1,\nabla^\susp}.
 \end{align*}
\end{cor}
\begin{proof}
 Using the generalised triangle inequality and Lemma \ref{TensorSuspNormEqui - Lemma}, we estimate
 \begin{align*}
     ||\placeholder||_{1,\nabla^{\susp}} &\deff \left( ||\placeholder||_0^2 + ||\nabla^\susp(\placeholder)^2||_0 \right)^{1/2} \\
      &= \left( ||\placeholder||_0^2 + ||\nabla\otimes \id(\placeholder) + (\nabla^\susp - \nabla \otimes \id)(\placeholder)^2||_0 \right)^{1/2} \\
      &\leq \left( ||\placeholder||_0^2 + ||\nabla\otimes \id(\placeholder)^2||_0 \right)^{1/2} + ||(\nabla^\susp -\nabla \otimes \id)(\placeholder)||_0 \\
      &\leq ||\placeholder||_{1,\nabla\otimes \id} + ||\nabla^\susp - \nabla \otimes \id||_{0,0} \cdot ||\placeholder||_0\\
      &\leq \left( 1 + d \sum_{k=1}^n ||\partial_{t_k} h||_{\op}\right) ||\placeholder||_{1,\nabla \otimes \id}.
 \end{align*}
 The proof of the second estimate is similar.
\end{proof}

\begin{proof}[Proof of Theorem \ref{OperatorSuspensionFaktorisation - Theorem}]
 Let $P \in \AuxInvDir[n]$ with underlying block map $h \in \SingMet[M][n]$. 
 Set $g = \susp(h)$. 
 We will prove the theorem by showing that $\susp(P)^2$ is invertible.
 Recall from Lemma \ref{susp quad - Lemma} that
 \begin{align*}
  \susp(P)^2 &= (P^\extension)^2 + \{P^{\extension}, \Dirac\} + \Dirac^2 \\
             &= (P^2)^\extension - \sum_{j=1}^n \partial_j \cdot \lbrack P^\extension, \nabla_j^{\Spinor_g} \rbrack - \Delta + \sum_{1\leq j < k} \mathfrak{R}(\partial_j,\partial_k).
 \end{align*}
 Note that 
 \begin{equation*}
   - \Delta = \sum_{j=1}^n (\nabla^{\Spinor_g}_{\partial_j})^2
 \end{equation*}
 is a non-negative operator. The proof is analogous to the one for the Dirac operator presented in \cite{LawsonMichelsonSpin}*{Proposition II.5.3}.
 Since $P$ is invertible, $P^2$ is strictly positive.
 We have, for each compactly supported smooth section, the following estimate
 \begin{align*}
  \langle (P^\extension)^2 u,u\rangle &= ||P^\extension(u)||_0^2 = \int_{\R^n} ||P_tu||^2_{0,\Spinor_g|_{M_t}} \diff t \\
  &\geq \int_{\R^n} \lowBnd(t)^2 ||u||^2_{1,\Spinor_g|_{M_t},\nabla\otimes \id} \diff t \\
  &> \frac{1}{4}\int \lowBnd(t)^2||u||^2_{1,\Spinor_g|_{M_t},\nabla^\susp} \diff t. 
 \end{align*}
 Here, we used in the first estimate that $\Sobolev{1}{\Spinor_g|_{M_t},\nabla\otimes \id}$ is isometric isomorphic to $\Sobolev{1}{\Spinor_{h(t)},\nabla} \otimes \Cliff$, so that we can apply $||P_tu||_0 \geq \lowBnd(t)^2 ||u||_1$ component-wise. 
 The second estimate follows from the equivalence of norms (Corollary \ref{EquiSobolevNorms - Cor}) and condition (\ref{eq: AuxDir1}) of $\AuxInvDir$.
 
 In the following estimate we denote by $||\placeholder||_{1;t}$ the Sobolev norm on $\Spinor_g|_{M_t}$ with respect to $\nabla^{\susp}$: 
 \begin{align*}
  \left|\langle \sum_{j=1}^n\lbrack P^{\extension}, \Dirac\rbrack u,u\rangle\right| &\leq \int_{\R^n} \left|\sum_{j=1}^n \langle  \partial_j \cdot \lbrack P^\extension, \nabla^{\Spinor_g}_{\partial_j} \rbrack u_t, u_t\rangle\right|_{0;t} \diff t \\
  &\leq \int_{\R^n} ||\sum_j \partial_j \cdot \lbrack P^\extension, \nabla^{\Spinor_g}_{\partial_j} \rbrack_t u_t||_{0;t} ||u_t||_{0;t} \diff t \\
  &\leq \int_{\R^n} ||\sum_j \partial_j \cdot \lbrack P^\extension, \nabla^{\Spinor_g}_{\partial_j} \rbrack_t||_{1,0;t} ||u_t||_{1;t}^2 \diff t \\
  &\leq \int_{\R^n} \sum_{j=1}^n \underbrace{||\partial_j \cdot||_{0,0;t}}_{=1} ||\lbrack P^\extension, \nabla^{\Spinor_g}_{\partial_j}\rbrack||_{1,0;t} ||u_t||^2_{1;t} \diff t \\
  &< \int_{\R^n} 1/32 \cdot \lowBnd(t)^2 ||u_t||^2_{1;t} \diff t.
 \end{align*}
  
 We estimate the fourth summand with the help of Corollary \ref{Spinor Curvature Norm Estimation - Cor} in a similar manner and notation as before
 \begin{align*}
  |\langle \sum_{1\leq j<k\leq n}\mathfrak{R}(\partial_j,\partial_k)(u),u\rangle| &\leq \int_{\R^n} |\langle \sum_{j<k}\mathfrak{R}_t(\partial_j,\partial_k)(u_t),u_t\rangle| \diff t \\
  &\leq \int_{\R^n}||\sum_{j<k}\mathfrak{R}_t(\partial_j,\partial_k)||_{0,0;t} ||u_t||_{0;t}^2 \diff t \\
  &\leq \int_{\R^n}||\sum_{j<k}\mathfrak{R}_t(\partial_j,\partial_k)||_{\op,t} ||u_t||_{1;t}^2 \diff t \\
  &\leq \int_{\R^n} d^2/2 \left( \sum_{j=1}^n ||(\partial_j h)^\op||_{\op,t}\right)^2 ||u||_{1;t}^2\diff t \\
  &< \int_{\R^n} 1/64\cdot 4^{-(d+1)} \cdot d^{2}/2 \cdot \lowBnd(t)^2 ||u_t||^2_{1;t}  \diff t\\
  &< \int_{\R^n} 1/32 \lowBnd(t)^2 ||u_t||^2_{1;t} \diff t. 
 \end{align*}
 Plugging everything together, we get:
 \begin{align*}
  &\quad \  |\langle \susp(P)^2u,u\rangle| \\
  &\geq \langle (P^\extension)^2 u,u\rangle - \langle\Delta u,u \rangle - |\langle\{P^\extension,\Dirac\}u,u\rangle| - |\langle \sum_{1\leq j < k \leq n} \mathfrak{R}(\partial_j,\partial_k)u,u\rangle| \\
  &\geq \int_{\R^n} (1/4 - 2/32) \lowBnd(t)^2||u_t||_{1;t}^2 \diff t \\
  &\geq 3/16  \min\{\lowBnd(t)^2 \, : \, t \in \R^n\} ||u||_0^2,
 \end{align*}
 which is equivalent to
 \begin{align*}
  ||\susp(P)u||_0 \geq \sqrt{3}/4\min\{\lowBnd(t) \,:\, t \in \R^n\}||u||_0.
\end{align*}  
 Since $\susp(P)$ is a self-adjoint operator, this inequality implies that $\susp(P)$ is invertible.
\end{proof}
Finally, we wish to show that $\AuxInvDir$ is a Kan set and that the inclusion is a weak homotopy equivalence. 
To this end, we need to know how the anti-commutator behaves under rescaling. 

\begin{lemma}\label{Rescaling Behaviour - Lemma}
 Let $D \subseteq \R^n$ be an open subset of $\R^n$ and let $\varphi \colon D \rightarrow D$ and $P \colon D \rightarrow \PseudDir[ ]$ be a smooth maps. 
 The map $P$ comes with an underlying smooth map $h \colon D \rightarrow \Riem(M)$.
 Then we have
 \begin{gather*}
  \lbrack (P \circ \varphi)^\extension, \nabla_{\partial_{t_i}}^{\Spinor_{\susp(h\circ \varphi)}}\rbrack = \sum_{k=1}^n \lbrack P^{\extension}, \nabla_{\partial_{t_k}}^{\Spinor_{\susp(h)}}\rbrack \circ \varphi \cdot \partial_{t_i} \varphi_k, \\
  \partial_{t_i}(h \circ \varphi) = \sum_{k=1}^n (\partial_{t_k}h) \circ \varphi \cdot \partial_{t_i} \varphi_k.
 \end{gather*}
\end{lemma}
\begin{proof}
 The second equation is simply the chain rule for differentials, so it remains to prove the first equation.
 Let $g_{00} \colon D \rightarrow \Riem(M)$ be the constant map with value $g_{00}$. 
 Set further $g_0 = \susp(g_{00})$. 
 Recall from the proof of Lemma \ref{CommutatorWithExtension - Lemma} the identity
 \begin{equation*}
  \preGauge(\susp (h),g_0)_\ast(P^\extension) = (\preGauge(h,g_{00})_\ast P)^\extension.
 \end{equation*}
 In the following, we will abbreviate $\preGauge(\susp(h),g_0)$ to $\preGauge(\susp(h))$ and $\preGauge(h,g_{00})$ to $\preGauge(h)$.
 Using Lemma \ref{CalcChristoffelSymb - Lemma}, we calculate 
 \begin{align*}
  &\qquad \preGauge(\susp(h \circ \varphi))_\ast\left(\lbrack (P \circ \varphi)^\extension, \nabla_{\partial_{t_i}}^{\Spinor_{\susp(h \circ \varphi)}} \rbrack\right) \\
  &= \lbrack (\preGauge(\susp (h\circ \varphi))_\ast(P \circ\varphi))^\extension, (\preGauge(h \circ \varphi)_\ast\nabla^{\Spinor_{\susp(h\circ \varphi)}})_{\partial_{t_i}}\rbrack \\
  &= \lbrack (\preGauge(h \circ \varphi)_\ast (P \circ\varphi))^\extension, \partial_{t_i}\rbrack = \lbrack \bigl((\preGauge(h)_\ast(P))\circ \varphi\bigr)^\extension, \partial_{t_i} \rbrack\\
  &\underset{\mathrm{rule}}{\overset{\mathrm{chain}}{=}} \Bigl(\partial_{t_i}\bigl((\preGauge(h)_\ast P) \circ \varphi\bigr)\Bigr)^\extension = \sum_{j=1}^n ((\partial_{t_k}\preGauge(h)_\ast P)\circ \varphi)^\extension \cdot \partial_{t_i} \varphi_k.
 \end{align*}
 A similar calculation shows 
 \begin{align*}
   \preGauge(\susp(h \circ \varphi))_\ast &\left( \sum_{j=1}^n (\lbrack P^\extension, \nabla_{\partial_{t_k}}^{\Spinor_{\susp h}}\rbrack)\circ \varphi \cdot \partial_{t_i} \varphi_k\right) \\
   &= \sum_{j=1}^n ((\partial_{t_k}\preGauge(h)_\ast P)\circ \varphi)^\extension \cdot \partial_{t_i} \varphi_k,
 \end{align*}
 so the lemma follows.
\end{proof}
\begin{proposition}\label{AuxPseduDir Is Kan - Prop}
 The cubical set $\AuxInvDir$ is a Kan set.
\end{proposition}
\todo[inline]{Draw a picture and give a rough proof idea.}
 The rough idea of the proof is the following one:
 A given cubical $n$-horn $\CubeBox$ gives rise to a well defined map $P_0$ on the complement of an infinte half-cylinder $Q_K$ with base $K$, where $K$ is a compact, convex set with smooth boundary.
 By the block form of the operator-valued maps in the given horn, the operator-valued map $P_0$ is the elongation of $P_0|_{\{\eps x_i < S\}}$ for some sufficiently large $S$.
 We need to "extend" $P_0|_{\{\eps x_i < S\}}$ to $Q_K$ and we do this by modifying $P_0|_{\{\eps x_i < S\}}$ on a compact set of $\{\eps x_i \leq S\}$.
 More precisely, we do this by pushing down the argument of $P_0$ in $\eps x_i$-direction as closer we get to $Q_K$ so that we can extend the modification to $Q_K$ by $P_{(i,-\eps)} \circ \CubeProj{i}$.
\begin{proof}
 Given an arbitrary cubical $n$-horn
 \begin{equation*}
  \CubeBox = \{P_{(j,\omega)} \in \AuxInvDir[n-1] \, : \, \face[\omega]{j}P_{(k,\eta)} = \face[\eta]{k-1}P_{(j,\omega)}; \, j<k; (j,\omega), (k,\eta) \neq (i,\eps)\}.
 \end{equation*}
 Since that horn contains only finitely many block maps, there is a sufficiently large $\rho>0$ such that all block maps $P_{(j,\omega)}$ decompose outside of $\rho I^{n-1}$.
 We restrict $\degen{j}P_{(j,\omega)}$ to a map $\{\omega x_j > \rho \} \rightarrow \InvPseudDir[ ]$.
 By the compatibility requirement, these restrictions agree on the intersection of their domains and we can therefore glue them together to a map
 \begin{equation*}
   P_0 \colon \left(\quader[n][i,\eps][\rho]\right)^c\rightarrow \InvPseudDir[ ],
 \end{equation*}
 where $\quader[n][i,\eps][\rho] \deff \{x \in \R^n \, : \, x_j \in [-\rho,\rho] \text{ if } j \neq i,\,  \eps x_i \geq -\rho\}$ is the cuboid with diameter $\rho$. 
 Let $K \subseteq \R^{n-1}$ be a compact, convex subset that contains $\rho I^{n-1}$, that is point symmetric at the origin, and that has a smooth boundary. 
 Let $Q_K \deff \{x \in \R^n \, : \, \CubeProj{i}(x) \in K, \eps x_i \geq -\rho \}$ be the cuboid with base $K$. 
 Restricting $P_0$ yields 
 \begin{equation*}
  P_0 \colon Q_K^c \cap \{\eps x_i < \rho + R_1\} \rightarrow \InvPseudDir[ ],
 \end{equation*}
 for a fixed choice of $R_1>1$. 
 
 For the chosen $K$, there is a unique norm $||\placeholder ||_K$ whose unit ball is $K$. 
 Since $\partial K$ is smooth, the norm is smooth on $\R^{n-1}\setminus \{0\}$.
 For a sufficiently large $R_2 = R_2(\rho ,R_1)$ to be determined later, pick $\chi \colon \R \rightarrow \lbrack 0,1\rbrack$ that is identically $1$ on $\R_{\leq 1.2}$, vanishes on $\R_{\geq R_2-1}$, and whose first and second derivative satisfy $|\chi'| \leq 1.3/R_2$ and $|\chi''| \leq 10/R_2^2$.
 Furthermore, pick a monotonically increasing function $q \colon \R \rightarrow \R$ that is the identity on $\R_{\leq \rho}$, the constant map with value $\rho + R_1 - 1$ on $\R_{\rho + R_1}$, and whose derivatives satisfy $q' \leq 1$ and $|q''| \leq 1.2/R_1$.
 The map $q$ defines a map 
 \begin{align*}
  q &\colon \R^n \rightarrow \R^n \\
  (x_1, \dots, x_n) &\mapsto (x_1, \dots, x_{i-1}, \eps q(\eps x_i), x_{i+1}, \dots, x_n).
 \end{align*}
 Finally, we define
 \begin{align*}
  H \colon \left(Q_K^c  \cap \{\eps x_i < \rho + R_1\}\right) \times [0,1] &\rightarrow  Q_K^c \cap \{\eps x_i < \rho + R_1\}, \\
  (x,t) &\mapsto x - 2\eps(\rho +R_1)t e_i 
 \end{align*}
 and set 
 \begin{align*}
  P \colon \R^n &\rightarrow \InvPseudDir[ ], \\
  x &\mapsto \begin{cases}
   P_0 \circ H\bigl(q(x), \chi(||\CubeProj{i}(x)||_K)\bigr), & \text{if } x \in Q_K^c, \\
   P_{(i,-\eps)}(\CubeProj{i}(x)), & \text{if } x \in Q_K.
  \end{cases}
 \end{align*}
 Note that $P$ is well defined and smooth because $P_0$ agrees with $P_{(i,-\eps)}\circ \CubeProj{i}$ on the set $\{-\eps x_i > \rho\}$ and $\partial Q_K$ lies in the preimage of $\{-\eps x_i > \rho\}$ under the map $H \circ (q(\placeholder ), \chi(||\CubeProj{i}(\placeholder )||_K))$. 
 
 The map $P$ is a block map.
 Indeed, $\chi(||\CubeProj{i}(x)||_K) = 0$ on $\{\omega x_j > R_2\}$ for all $j \neq i$, which yields
 \begin{equation*}
  P(x) = P_0(H(q(x),0)) = P_0(q(x)) = P_0(x) = P_{(j,\omega)}(p_j(x))
 \end{equation*}
 because $P_{(j,\omega)(\CubeProj{j}(\placeholder))}$ is independent of $x_i$ if $|x_i|> \rho$.
 For all $x$ with $-\eps x_i > R_2$, we have
 \begin{align*}
  P(x) &= P_0(q(x) - 2 \eps(\rho+R_1)\chi(||\CubeProj{i}(x)||_K) e_i) \\
       &= P_0(x - 2 \eps(\rho+R_1)\chi(||\CubeProj{i}(x)||_K) e_i) \\
       &= P_{(i,-\eps)}(\CubeProj{i}(x))
 \end{align*}
 because, in this case, $-\eps x_i + 2\eps(\rho + R_1)t > R_2$ for every $t \geq 0$.
 Lastly, for all $x$ with $\eps x_i > R_2 > \rho + R_1$ we have 
 \begin{equation*}
   P(x) = P(x_1,\dots, x_{i-1},\eps(\rho + R_1 - 1), x_{i+1}, \dots, x_n)
 \end{equation*}
 by the definition of $q$.
 
 The proof also shows that $\face[\omega]{j}P = P_{(j,\omega)}$ for all $(j,\omega) \neq (i,\eps)$.
 This means that $P$ is a filler for $\CubeBox$ in $\InvPseudDir[n]$.
 
 It remains to show that $P$ lies in $B_n$.
 It is clear that $P$ satisfies the defining conditions of $B_n$ in an open neighbourhood of $Q_K$, for example on $\chi^{-1}(\R_{<1.2})$. 
 We will show that $P$ satisfies the defining conditions of $B_n$ on $\chi^{-1}(\R_{>1.1})$, too.
 Clearly, we have
 \begin{equation*}
  \lowBnd_P(t) = \lowBnd_{P_0}(H(q(t), \chi(||\CubeProj{i}(t)||_K)).
 \end{equation*}
 
 Actually, this equation holds for every smooth self map $\varphi \colon Q_K^c \rightarrow Q_K^c$.
 For every such self map $\varphi$ and each smooth map $Q \colon Q_K^c \rightarrow \InvPseudDir[ ]$ the underlying map of Riemannian metrics $h_Q$ satisfies $h_{Q \circ \varphi} = h_Q \circ \varphi$.
 By Lemma \ref{Rescaling Behaviour - Lemma}, we find the following upper bound for its derivatives
 \begin{align*}
  ||(\partial_jh_{Q \circ \varphi})^\op||_{\op, h_{Q \circ \varphi}} &= ||\sum_{k=1}^n (\partial_{k}h_Q)^\op\circ \varphi \cdot \partial_{j} \varphi_k||_{\op, h_{Q\circ \varphi}} \\
  &\leq \sum |\partial_{j} \varphi_k| \cdot ||(\partial_{k} h_Q)^\op \circ \varphi||_{\op, h \circ \varphi} \\
  &= \sum_{k=1}^n |\partial_{j} \varphi_k| \cdot ||(\partial_{k} h_Q)^\op||_{\op, h} \circ \varphi. 
\end{align*}  

We would like to apply this estimate to $Q = P_0$ and $\varphi = H \circ (q, \chi(||\CubeProj{i}(\placeholder )||_K)$.
For $j \neq i$, the derivatives of $\varphi_k$ are given by
\begin{align*}
 \partial_j\varphi_k = \begin{cases}
  0, & \text{if } k \notin\{i,j\}, \\
  1, & \text{if } k = j, \\
  -2\eps(\rho + R_1)\chi' \cdot (\partial_j|| \placeholder ||_K)\circ \CubeProj{i}, & \text{if } j<i, k=i, \\
  -2\eps(\rho + R_1)\chi' \cdot (\partial_{j-1} || \placeholder ||_K) \circ \CubeProj{i}, & \text{if } j>i, k=i.
 \end{cases}
\end{align*}
For $j=i$, they are given by
\begin{align*}
 \partial_i \varphi_k = \begin{cases}
  q'(\eps \cdot \placeholder), & \text{if } k=i,\\
  0, & \text{if } k \neq i.
 \end{cases}
\end{align*}
Plugging these results into the previous estimate yields 
\begin{align*}
 ||(\partial_j h_P)^\op||_{\op, h_P} \leq \begin{cases}
  ||(\partial_j h_{P_0})^\op||_{\op, P_0}\circ \varphi + |\partial_j \varphi_i| \cdot ||(\partial_{t_i}h)^\op||_{\op, h_{P_0}} \circ \varphi, & \text{if } j \neq i, \\
  q'(\eps \cdot \placeholder) \cdot || (\partial_{t_i}h_{P_0})^\op||_{\op, h_{P_0}}\circ \varphi, & \text{if } j= i,
 \end{cases}
\end{align*}
which gives
 \begin{align*}
  \begin{split}
    \sum_{j=1}^n ||(\partial_{j}h_P)^\op||_{\op, h_P} \leq &\sum_{j=1}^n ||\partial_{j}h_{P_0}^\op||_{\op, h_{P_0}} \circ \varphi \\
    &+ \left(\sum_{j=1}^{n-1} \frac{2(\rho + R_1)\cdot 1.3}{R_2} \cdot |\partial_{j}|| \placeholder ||_K| \circ \CubeProj{i} \cdot ||\partial_{i}h_{P_0}^\op||_{\op,h_{P_0}}\right) \circ \varphi. 
  \end{split}
\end{align*}

Applying Lemma \ref{Rescaling Behaviour - Lemma} to $[P^\extension, \nabla_{\partial_{j}}]$ we obtain, by a similar calculation, the estimate
\begin{align*} 
 \begin{split}
 &\quad \ \sum_{j=1}^n ||\lbrack P^\extension, \nabla_j^{\Spinor_{\susp(h_P)}}\rbrack||_{1,0,h_P} \\
 &\leq  \sum_{j=1}^n || \lbrack P_0^\extension, \nabla_{j}^{\Spinor_{\susp(h_{P_0})}}\rbrack||_{1,0,h_{P_0}} \circ \varphi \\
 & \quad \ + \left(\sum_{j=1}^{n-1} \frac{2.6(\rho + R_1)}{R_2}|\partial_j || \placeholder ||_K|\circ \CubeProj{i} \cdot  ||[P_0^\extension, \nabla_{i}^{\Spinor_{\susp(h_{P_0})}}]||_{1,0,h_{P_0}}\right) \circ \varphi. 
 \end{split}
\end{align*}

Since $\partial_i h^\op_P$ and $[P^\extension, \nabla_{\partial_{i}}]$ have smooth extensions from $\chi^{-1}(\R_{>1.1})$ to $\chi^{-1}(\R_{\geq 1.05})$ their norms are uniformly bounded by $B_1$, say.
Also all $|\partial_j || \placeholder ||_K|$ are uniformly bounded by $B_2$, say.
Thus, if we choose $R_2 > 2.6(\rho + R_1)B_1B_2\mathcal{I}^{-1}$, where $\mathcal{I} = \min\{\mathcal{I}_1,\mathcal{I}_2\}$ is the minimum of
\begin{equation*}
 \mathcal{I}_1 \deff \mathrm{inf}\left\{\frac{1}{32} \, \lowBnd_{P_0}(t)^2 - \frac{1}{2}\sum_{j=1}^n ||[P_0^\extension, \nabla_{j}^{\Spinor_{\susp h_{P_0}}}] ||_{1,0,h_{P_0}} \, : \, t \in Q_K^c\right\}
\end{equation*}
and
\begin{equation*}
 \mathcal{I}_2 \deff \left\{ \mathrm{min}\{d^{-1}, 2^{-(d+4)} \lowBnd(t)\} - \sum_{k=1}^n ||\partial_kh^\op||_{\op,h_{P_0}(t)} \, : \, t\in Q_K^c\right\}
\end{equation*}
then $P$ satisfies the defining equations of $\AuxInvDir$.
Note that $\mathcal{I} >0$ because $P_0|_{Q_K^c}$ is a union of restriction of degenerate elements of $B_n$ and the defining conditions of $B_n$ are local.

Thus, $P \in \AuxInvDir[n]$ is a filler for the given cubical $n$-horn and $\AuxInvDir$ is therefore a Kan set.
\end{proof}
\begin{cor}\label{AuxDir is weak equivalent -  Cor}
 The inclusion $\AuxInvDir \hookrightarrow \InvPseudDir[\bullet]$ is a weak homotopy equivalence.
\end{cor}
\begin{proof}
 We first prove that the inclusion induces a surjection on homotopy groups.
 Let $[P] \in \pi_n(\InvPseudDir[\bullet],P_0)$ be represented by a smooth block map $P \colon \R^n \rightarrow \InvPseudDir[ ]$ that is constant with value $P_0$ outside some cube $rI^n$.
 For each $\rho \geq 1$ the rescaled map ${}_\rho P(x) \deff P(\rho^{-1}x)$ is still constant near infinity with value $P_0$ and thus represents an element in $\pi_n(\InvPseudDir[\bullet])$, too.
 By Lemma \ref{Rescaling Behaviour - Lemma} there is a $\rho_0 \geq 1$ such that ${}_\rho P \in \pi_n(\AuxInvDir)$ for all $\rho \geq \rho_0$ because the left hand side of the defining inequalities for $\AuxInvDir$ rescale with $1/\rho$, while the right hand side is scaling invariant.
 
 We need to show ${}_\rho P$ and $P = {}_1P$ are homotopic relative boundary.
 A homotopy $H \in \InvPseudDir[n+1]$ relating ${}_\rho P$ and ${}_1 P$ is given by 
 \begin{equation*}
  H(t,x) \deff {}_{(1-\chi(t))\rho + \chi(t)}P(x),
\end{equation*}   
where $\chi \colon \R \rightarrow [0,1]$ is a smooth function that is identically zero on $\R_{\leq -1}$ and identically $1$ on $\R_{\geq 1}$. 
Indeed, $H$ satisfies
\begin{align*}
 \face[-1]{1}H(x) &= \lim_{R \to \infty} H(-R,x) = {}_\rho P(x), \\
 \face[1]{1}H(x) &= \lim_{R \to \infty}H(R,x) = {}_1P(x), \\
 \face{i}H(t,x) &= \lim_{R \to \infty}H(t,x_1,\dots,x_{i-1},\eps R,x_i, \dots, x_{n-1}) \\
                &= \lim_{R \to \infty} P\left(((1-\chi(t))\rho + \chi(t))^{-1}\cdot (x_1, \dots,\eps R,\dots, x_{n-1})\right) \\
                &= P_0 = \degen{1}\face{i-1}P(t,x). 
\end{align*}

To prove injectivity, we follow the common strategy that uniqueness is just a relative form of existence.
For each $[P] \in \pi_n(\AuxInvDir[\bullet],P_0)$ with $0=[P] \in \pi_n(\InvPseudDir[\bullet],P_0)$ there exists an element $H \in \InvPseudDir[n+1]$ that satisfies
\begin{align*}
 \face[1]{1}H &= P_0, \\
 \face[-1]{1}H &= P,\\
 \face{j}H &= \degen{1}\face{j-1}P = P_0.
\end{align*}

As before, there is a sufficiently large $\rho \geq 1$ such that ${}_\rho H \in \AuxInvDir[n+1]$.
Note that ${}_\rho H$ is a homotopy between ${}_\rho P$ and $P_0$ in $\AuxInvDir$, so it remains to find a homotopy between $P$ and ${}_\rho P$ in $\AuxInvDir$. 
An example for such a homotopy is $\mathcal{H} \deff \degen{1}P \circ {}_R\varphi$, where
\begin{align*}
 {}_R\varphi \colon \R^{n+1} &\rightarrow \R^{n+1} \\
 (t,x) &\mapsto \left(t, ((1-\chi(R^{-1}t)\rho + \chi(R^{-1}t))^{-1}x\right),
\end{align*} 
provided $R$ is sufficiently large.
Abbreviate $(1-\chi(t))\rho + \chi(t)$ to $\mu(t)$ and set ${}_R\mu(t) = \mu(R^{-1}t)$.
Furthermore, abbreviate $(t,x)$ to $y$ so that $t=y_1$ and $x_i = y_{i+1}$.
Lemma \ref{Rescaling Behaviour - Lemma}, applied to $y_1=t$ yields 
\begin{align*}
 \left\lbrack \mathcal{H}^\extension, \nabla_{\partial_{y_1}}^{\Spinor_{\susp (g_{\mathcal{H}}) }}\right\rbrack (y) &= \left\lbrack \degen{1} P^\extension \circ {}_R\varphi, \nabla_{\partial_{y_1}}^{\Spinor_{\susp (\degen{1}g_{P} \circ {}_R\varphi) }} \right\rbrack (y) \\
 &= \sum_{k=1}^{n+1} [\degen{1}P^\extension, \nabla_{\partial_{y_k}}^{\Spinor_{\susp (\degen{1} g_P)}}] \circ {}_R\varphi \cdot \partial_{y_1} {}_R\varphi_k(y) \\
 &= \frac{1}{R} \sum_{k=2}^{n+1} [\degen{1}P^\extension, \nabla_{\partial_{ y_k}}^{\Spinor_{\susp (\degen{1} g_P)}}] \circ {}_R \varphi \cdot {}_R\mu^{-2} \cdot ((\rho -1)\cdot {}_R(\chi'))y_k
\end{align*}
where we used $[\degen{1}P^\extension,\nabla^{\Spinor_{\susp(\degen{1}g_P)}}_{\partial_{y_1}}] = [\degen{1}P^\extension,\partial_{y_1}] = (\partial_{y_1}\degen{1}P)^\extension = 0$.
Applied to $y_i = x_{i-1}$ for $i>1$, Lemma \ref{Rescaling Behaviour - Lemma} yields
\begin{align*}
 \left\lbrack \mathcal{H}^\extension, \nabla_{\partial_{y_i}}^{\Spinor_{\susp (g_{\mathcal{H}})}}\right\rbrack (y) = \dots &= \sum_{k=2}^{n+1} [\degen{1}P^\extension, \nabla_{\partial_{y_k}}^{\Spinor_{\susp (\degen{1} g_P)}}] \circ {}_R \varphi \cdot {}_R\mu^{-1} \delta_{ik} \\
 &= {}_R\mu^{-1} \cdot [\degen{1}P^\extension, \nabla_{\partial_{y_i}}^{\Spinor_{\susp (\degen{1} g_P)}}] \circ {}_R \varphi\\
 &= {}_R\mu^{-1} \cdot \degen{1}\left([P^\extension, \nabla_{\partial_{x_{i-1}}}^{\Spinor_{\susp (g_{P})}} ]\right) \circ {}_R\varphi. 
\end{align*}
We assumed that $P = P_0$ outside of $rI^n$, so $[\degen{1}P^\extension, \nabla_{\partial_{y_k}}] \circ {}_R \varphi$ is supported within $\R \times \rho r I^n$ for $k>1$ and therefore 
it is uniformly bounded.
Thus, we can choose $R$ sufficiently large such that $[\mathcal{H}^\extension, \nabla_{\partial_{t}}]$ is arbitrarily small.
Since $\lowBnd_\mathcal{H} = \lowBnd_{\degen{1}P} \circ {}_R \varphi = \lowBnd_P \circ \CubeProj{1} \circ {}_R\varphi = \degen{1}(\lowBnd_P) \circ {}_R\varphi$ and since ${}_R\mu^{-1} \leq 1$, we have
\begin{equation*}
    \sum_{k=2}^{n+1} \left|\left|\left\lbrack \mathcal{H}^\extension, \nabla_{\partial_{y_k}}^{\Spinor_{\susp( g_{\mathcal{H}})}}\right\rbrack (y)\right|\right|_{1,0,M_{y}} < 1/32 \lowBnd_{\mathcal{H}}(y)^2, 
\end{equation*}
so we can choose $R>1$ so large that $\mathcal{H}$ satisfies the second defining inequality (\ref{eq: AuxDir2}) of $\AuxInvDir$.

The same strategy shows that the underlying map of metrics $g_{\mathcal{H}}$ satisfies the first defining inequality (\ref{eq: AuxDir1}) after enlarging $R$, if necessary.
\end{proof}

\section{Kan Property of the Operator Concordance Set}\label{Section - Kan Property of the Operator Concordance Set}

The main theorem of this section is the analogous result of Theorem \ref{Concordance Set is Kan - Theorem} for invertible block Dirac operators.

\begin{theorem}\label{operatorfundamentallemma}
 The operator concordance set $\InvBlockDirac$ is a Kan set for all closed, spin manifolds $M$ that admit at least one invertible pseudo Dirac operator.
\end{theorem}
The proof of this theorem relies on an operator theoretic analog of Lemma \ref{Metric Modification - Lemma}.

\begin{lemma}[Operator Modification Lemma]\label{Operator Modification - Lemma}
 Let $P$ be a block Dirac operator on $M \times \R^n$ with underlying block metric $g$. 
 Assume that the core of $P$ is contained in $M \times \rho/2 \cdot I^n$ for $\rho > 20$ (in particular, $g$ decomposes on the complement).
 Assume further that $\face[\omega]{j} P$ is invertible if $(j,\omega) \neq (1,1)$.
 Let $H$ and $\Phi$ as in Theorem \ref{Embedding Theorem} and denote ${}_R\Phi \deff R \Phi \circ (\id \times R^{-1})$. 
 Then
 \begin{itemize}
     \item[(o)] the operator $P$ restricts to $M \times R\Phi(H \times (a,b))$ for all $1 \leq a < b \leq 2$ and all $R \geq 1$. 
 \end{itemize}
 
 For all $R \geq 1$, there is a Riemannian metric $G_R$ on $M \times R\Phi(H \times [1,2])$, a compact subset $K = K(R)$ of $M \times \R^n$, a relative compact, open subset $L=L(R) \supseteq K$ of $M \times \R^n$, and a pseudo Dirac operator $Q\deff Q(R)$ on $\Gamma_c(R\Phi(H \times (1,2)),\Spinor_{G_R})$ that satisfy:
 \begin{enumerate}
  \item[(i)] The operator ${}_R\Phi^\ast Q$ restricts to $M \times H \times (a,b)$ for all $R \leq a <b \leq 2R$ such that ${}_R\Phi^\ast Q$ restricts to ${}_R\Phi^{-1}(L)$ and to the complement of ${}_R\Phi^{-1}(K)$.
  \item[(ii)] The restriction of ${}_R\Phi^\ast Q$ and ${}_R\Phi^\ast P$ to $M \times H \times (1.8R,2R)$ and to the complement of $K$ agree.
  \item[(iii)] The metric ${}_R\Phi^\ast G_R$ is a product metric on $M \times H \times [R,1.2R]$. 
  Under the decomposition $\Spinor_{{}_R\Phi^\ast G_R}|_{M \times H \times [R,1.2R]} \cong \Spinor_{{}_R\Phi^\ast G_R}\restrict_{M \times H \times {1.1R}} \boxtimes Cl_{1,0}$, the operator ${}_R\Phi^\ast Q$ decomposes into
  \begin{equation*}
   {}_R\Phi^\ast Q = {}_R\Phi^\ast Q\restrict_{{M \times H \times {1.1R}}} \boxtimes \id + \evenodd \boxtimes \DiracSt_{\R}.
  \end{equation*}   
  \item[(iv)] There is an $R_0 \geq 1$ and a constant $c > 0$, such that $Q$ is bounded from below with lower bound $\geq c$ for all $R \geq R_0$. 
  \item[(v)] If $P - \Dirac_{g_P}$ is bounded, then the operator $Q - \Dirac_{G_R}$ is bounded. 
\end{enumerate}  
\end{lemma}
\begin{proof}
 We start with the proof of $(o)$.
 The block form of $P$ implies that it restricts to $V_{R\rho}(\fatalpha)$, see Lemma \ref{CoreLowerBound - Lemma}, and that $P|_{V_{R\rho(\fatalpha)}}$ can be identified with $P(\hat{\fatalpha}) \boxtimes \id + \evenodd \boxtimes \DiracSt_{\R^{\hat{\fatalpha}}}$ via a permutation of coordinates.
 The same permutation of coordinates gives a diffeomorphism between $R\Phi(H \times (1,2)) \cap V_{R\rho}(\fatalpha)$ and $(-R\rho,R\rho)^{|\mathrm{Null \, \fatalpha}|} \times R\Phi^{\hat{\fatalpha}}(H^{\hat{\fatalpha}} \times (1,2))$.
 The first factor contains the core of $P(\hat{\fatalpha})$.
 On the second factor, $P$ acts as the differential operator $\DiracSt_{\R^{\hat{\fatalpha}}}$.
 Thus, $P$ restricts to $M \times V_{R\rho}(\fatalpha) \cap R\Phi(H \times (a,b))$ for all $1\leq a < b \leq 2$.
 
 As $\{V_{R\rho}(\fatalpha) \, : \, \fatalpha \neq \mathbf{0}, (1,0\dots,0)\}$ covers $R\Phi(H \times (1,2))$, the operator $P$ restricts to a symmetric operator on $M \times R\Phi(H \times (a,b))$.
 
 The same strategy shows that $P$ also restricts to a symmetric operator on the complement of $V_{R\rho}(\mathbf{0})$. 
 This proves $(o)$. 
 
 For all $R \geq 1$, let $G_R$ be the metric from the Metric Modification Lemma \ref{Metric Modification - Lemma}.
 Define the operator $Q$ on $\Gamma_c(M \times R\Phi(H \times (1,2)) , \Spinor_{G_R})$ via
 \begin{equation*}
     Q \deff Q(R) \deff \fullGauge(g_P,G_R)_\ast(P - \Dirac_{g_P}) + \Dirac_{G_R}. 
 \end{equation*}
 
 Since $Q$ is obtained from $P$ by pushing it forward with a vector bundle map and adding a differential operator to $P$, we conclude $Q \in \PseudOpCl^1(\Spinor_{G_R})$.
 The operator $Q$ also restricts to $M \times R\Phi(H \times (a,b))$ for all $1\leq a < b\leq 2$, to the intersection of $L(R) \deff V_{R(\rho+5)}(\mathbf{0})$ and $M \times R\Phi(H \times (a,b))$, and to the complement of $K(R) \deff \overline{V_{R(\rho+4)}(\mathbf{0})}$ in $M \times R\Phi(H \times (a,b))$, which proves $(i)$.
 
 The gauge map $\fullGauge(g_P,G_R)$ is an isometry between spaces of square integrable spinors, so $Q$ is symmetric.
 We will calculate below that the principal symbol of $Q$ and $\Dirac_{G_R}$ agree, thus $Q$ is a symmetric pseudo Dirac operator.
 
 To prove $(ii)$, recall that $g_P$ and $G_R$ agree on $M \times R\Phi(H \times (1.8,2))$, so that $\fullGauge(g_P,G_R) = \id$ and $\Dirac_{g_P} = \Dirac_{G_R}$ on this subset. 
 Thus, $P = Q$ on this subset.
 The same is true on the complement of $K(R) = \overline{V_{R(\rho+4)}(\mathbf{0})}$.
 
 Next we prove $(iii)$.
 Consider the restriction of ${}_R\Phi^\ast Q$ to ${}_R\Phi^{-1}(V_{R\rho}(\fatalpha))$.
 Recall from the Metric Modification Lemma \ref{Metric Modification - Lemma} that, on these subsets, the metrics decompose into
 \begin{align*}
     {}_R\Phi^\ast g_P = R^\ast\face[\hat{\fatalpha}]{ }g \oplus ({}_R\Phi^{\hat{\fatalpha}})^\ast \euclmetric \quad \text{and} \quad {}_R\Phi^\ast G_R = R^\ast\face[\hat{\fatalpha}]{ }g \oplus k_R^{\hat{\fatalpha}}. 
 \end{align*}
 Thus, the gauge map $\fullGauge({}_R\Phi^\ast g_P,{}_R\Phi^\ast G_R)$ decomposes into $\id \boxtimes \fullGauge(({}_R\Phi^{\hat{\fatalpha}})^\ast\euclmetric, k_R^{\hat{\fatalpha}})$.
 Since the gauge map is natural with respect to linear maps, we get
 \begin{align*}
     {}_R\Phi^\ast Q &= \fullGauge({}_R\Phi^\ast g_P, {}_R\Phi^\ast \euclmetric)_\ast ({}_R\Phi^\ast P - \Dirac_{{}_R\Phi^\ast g_P}) + \Dirac_{{}_R\Phi^\ast G_R} \\
     \begin{split}
         \begin{split}
         &= \id \boxtimes \fullGauge(({}_R\Phi^{\hat{\fatalpha}})^\ast \euclmetric, k_R^{\hat{\fatalpha}})_\ast\bigl(P(\hat{\fatalpha})\boxtimes \id + \evenodd \boxtimes \Dirac_{({}_R\Phi^{\hat{\fatalpha}})^\ast \euclmetric} \\
         & \qquad \qquad \qquad \qquad \qquad \qquad -\Dirac_{R^\ast \face[\hat{\fatalpha}]{ }g } \boxtimes \id - \evenodd \boxtimes \Dirac_{({}_R\Phi^{\hat{\fatalpha}})^\ast \euclmetric}\bigr)
         \end{split}\\
         &\quad  + \Dirac_{R^\ast \face[\hat{\fatalpha}]{ }g} \boxtimes \id + \evenodd \boxtimes \Dirac_{k_R^{\hat{\fatalpha}}}
     \end{split}\\
     &= P(\hat{\fatalpha}) \boxtimes \id + \evenodd \boxtimes \Dirac_{k_R^{\hat{\fatalpha}}}.
 \end{align*}
 Since $k_R^{\hat{\fatalpha}} = k_R^{\hat{\fatalpha}}\restrict_{H \times 1.1R} \oplus \diff t^2$ on $H^{\fatalpha} \times (R,1.2R)$, the Dirac operator $\Dirac_{k_R^{\hat{\fatalpha}}}$ decomposes accordingly on this set, hence $Q$ does, too.
 These decompositions are compatible on the intersections of different $V_{R\rho}(\fatalpha)$, so we get a global decomposition
 \begin{equation*}
     {}_R\Phi^\ast Q = {}_R\Phi^{\ast} Q \restrict_{M \times H \times 1.1R} \boxtimes \id + \evenodd \boxtimes \DiracSt_{\R},
 \end{equation*}
 which proves $(iii)$.
 
 We can read off the principal symbol of ${}_R\Phi^\ast Q$ from the (local) decompositions.
 On ${}_R\Phi^{-1}(V_{R\rho}(\fatalpha))$, we have
 \begin{align*}
     \symb_1({}_R\Phi^\ast Q) &= i\Cliffmult_{\face[\hat{\fatalpha}]{ }g}(\face[\hat{\fatalpha}]{ }g_\sharp(\placeholder)) \boxtimes \id + \evenodd \boxtimes i\Cliffmult_{k_R^{\hat{\fatalpha}}}({k_R^{\hat{\fatalpha}}}_\sharp(\placeholder)) \\
     &= \symb_1(\Dirac_{{}_R\Phi^\ast G_R}),
 \end{align*}
 so $Q$ also has the same principal symbol as $\Dirac_{G_R}$.
 
 We continue with the proof of $(iv)$.
 Since the gauge map $\fullGauge(g_P,G_R)$ induces an isometry between $H^0(\Spinor_{g_P})$ and $H^0(\Spinor_{G_R})$, the operator $\fullGauge(g_P,G_R)_\ast(P)$ is bounded from below with the same lower bound as $P$.
 Being bounded from below is an open condition for elliptic operators in $\mathrm{Hom}(H^1(\Spinor_{G_R}),H^0(\Spinor_{G_R}))$, 
 see Lemma \ref{TopInjOpen - Lemma}, so it suffices to show that $Q$ and $\fullGauge(g_P,G_R)_\ast P$ are sufficient close.
 
 The difference of these two operators is a differential operator that is supported within $M \times R\Phi(H \times [1,1.8])$.
 
 The previous calculation shows that 
 \begin{align*}
     (R\Phi)^\ast(Q - \fullGauge(g_P,G_R)_\ast(P)) = \left( \evenodd \boxtimes \Dirac_{R^2k^{\hat{\fatalpha}}} - \fullGauge((R\Phi^{\hat{\fatalpha}})^\ast \euclmetric, R^2k^{\hat{\fatalpha}})_\ast \Dirac_{(R\Phi^{\hat{\fatalpha}})^\ast\euclmetric}\right)
 \end{align*}
 on $R\Phi^{-1}(V_{R\rho}(\fatalpha))$, which we abbreviate to $\evenodd \boxtimes T$.
 Recall that $R\Phi^{-1}(V_{R\rho}(\fatalpha)) = M \times (-\rho,\rho)^{\mathrm{Null}(\fatalpha)} \times H^{\hat{\fatalpha}} \times (1,2)$ for $\fatalpha \neq \mathbf{0},(1,0\dots,0)$.
 Hence, for each $\sigma \in \Gamma_c(R\Phi^{-1}(V_{R\rho}(\fatalpha)),\Spinor_{\face[\hat{\fatalpha}]{ }g} \boxtimes \Spinor_{R^2k^{\fatalpha}})$, using Fubini's theorem we have
 \begin{align*}
     &\qquad ||\evenodd \boxtimes T(\sigma)||_0^2 = \langle \id \boxtimes T^2(\sigma),\sigma\rangle \\
     &= \int_{M \times (-R\rho,R\rho)^{|\mathrm{Null}(\fatalpha)|}} \int_{H^{\hat{\fatalpha}} \times (1,2)} ||\id \boxtimes T^2 \sigma(t,x)||^2_{\Spinor_{\face[\hat{\fatalpha}]{ }g} \boxtimes \Spinor_{R^2k^{\fatalpha}}} \diff \vol_{R^2k^{\hat{\fatalpha}}}(x) \diff \vol_{\face[\hat{\fatalpha}]{ }g}(t) \\
     &\leq \int \int ||T||_{1,0}^2\left(||\sigma(t,x)||^2 + ||\nabla^{\Spinor_{R^2k^{\hat{\fatalpha}}}} \sigma(t,x)||^2\right) \diff \vol(x) \diff \vol(t) \\
     &\leq ||T||_{1,0}^2 ||\sigma||_{1}^2,
 \end{align*}
 where the first inequality is the operator norm inequality applied to the section $\sigma(t,\placeholder)$.
 This gives $||\evenodd \boxtimes T||_{1,0} \leq ||T||_{1,0}$.
 
 It follows from \cite{hijazi1986conformal}*{Prop 4.3.1 and 4.2.1} that
 \begin{equation*}
     \fullGauge(R^2g,g)_\ast \Dirac_{R^2g} = R^{-1} \Dirac_g \quad \text{and} \quad \fullGauge(R^2g,g)_\ast \nabla^{\Spinor_{R^2g}} = \nabla^{\Spinor_{g}}
 \end{equation*}
 for all Riemannian metrics $g$ and all values $R > 0$.
 The latter equation implies that the Gauge map induces an isometry between all Sobolev spaces.
 Furthermore, we have $\fullGauge\bigl(R^2k^{\hat{\fatalpha}},k^{\hat{\fatalpha}}\bigr) \circ \fullGauge\bigl((R\Phi^{\hat{\fatalpha}})^\ast\euclmetric, R^2k^{\hat{\fatalpha}}\bigr) = \fullGauge\bigl((R\Phi^{\hat{\fatalpha}})^\ast\euclmetric, k^{\hat{\fatalpha}}\bigr)$.
 Putting these facts together yields
 \begin{align*}
     ||T||_{1,0} &= || \Dirac_{R^2k^{\hat{\fatalpha}}} - \fullGauge((R\Phi^{\hat{\fatalpha}})^\ast\euclmetric, R^2k^{\hat{\fatalpha}})_\ast \Dirac_{(R\Phi^{\hat{\fatalpha}})^\ast\euclmetric} ||_{1,0} \\
     &= ||\fullGauge(R^2k^{\hat{\fatalpha}}, k^{\hat{\fatalpha}} )_\ast\left( \Dirac_{R^2k^{\hat{\fatalpha}}} - \fullGauge((R\Phi^{\hat{\fatalpha}})^\ast\euclmetric, R^2k^{\hat{\fatalpha}})_\ast \Dirac_{(R\Phi^{\hat{\fatalpha}})^\ast\euclmetric}\right) ||_{1,0} \\
     &=|| R^{-1}\Dirac_{k^{\hat{\fatalpha}}} - \fullGauge((R\Phi^{\hat{\fatalpha}})^\ast\euclmetric, k^{\hat{\fatalpha}})_\ast(\Dirac_{(R\Phi^{\hat{\fatalpha}})^\ast\euclmetric}) ||_{1,0} \\
     &=|| \Dirac_{k^{\hat{\fatalpha}}} - \fullGauge({\Phi^{\hat{\fatalpha}}}^\ast\euclmetric, k^{\hat{\fatalpha}})_\ast(\Dirac_{{\Phi^{\hat{\fatalpha}}}^\ast\euclmetric}) ||_{1,0} /R \xrightarrow{R \to \infty} 0.
 \end{align*}
 It follows that $Q$ and $\fullGauge(g_P,G_R)_\ast(P)$ are arbitrarily close if $R$ is sufficient large.
 Together with Lemma \ref{TopInjOpen - Lemma}, we find an $R_0 \geq 1$ such that for all $R \geq R_0$, the operator $Q$ is bounded from below with lower bound $c>0$, that only depends on the lower bound of $P$ and some chosen upper bound of the distance of $Q$ to $\fullGauge(g_P,G_R)_\ast(P)$.
 In fact, the proof of Lemma \ref{TopInjOpen - Lemma} shows, that the constant $c$ can be made arbitrarily close to the lower bound of $P$ at the cost of increasing $R_0$.
 
 It remains to prove $(v)$.
 On $V_{R\rho}(\fatalpha) \cap M \times R\Phi(H \times (1,2))$, we have
 \begin{align*}
     {}_R\Phi^\ast (Q - \Dirac_{G_R}) = (P(\hat{\fatalpha})  - \Dirac_{\face[\hat{\fatalpha}]{ }g}) \boxtimes \id, 
 \end{align*}
 which is bounded by assumption of $P$.
 Hence, ${}_R\Phi^\ast(Q - \Dirac_{G_R})$ is bounded by a partition of unity argument.
\end{proof}

\begin{proof}[Proof of Theorem \ref{operatorfundamentallemma}]
 By assumption, $M$ has at least one invertible pseudo Dirac operator so that $\InvBlockDirac[0]$ (and hence $\InvBlockDirac$) is non-empty.
 A given cubical $n$-horn $\CubeBox \rightarrow \InvBlockDirac$ is represented by the following set of invertible block Dirac operators
 \begin{equation*}
     \{P_{(j,\omega)} \in \InvBlockDirac[n-1] \, : \, \face[\omega]{j}P_{(k,\eta)} = \face[\eta]{k-1}P_{(j,\omega)} \text{ for } j<k; (j,\omega),(k,\eta) \neq (i,\eps) \}.
 \end{equation*}
 The problem is symmetric in $(i,\eps)$, so we assume $(i,\eps) = (1,1)$.
 
 We pick a sufficient large $\rho>20$ such that the cores of all $P_{(j,\omega)}$ are contained in $M \times \rho/2I^{n-1}$.
 The operators $\degen{j}P_{(j,\omega)}|_{\{\omega x_j > \rho\}}$ restrict to each $\{\eta x_k > \rho\}$ if the intersection is not empty and the restrictions agree there, see Lemma \ref{Operator Suspension Union - Lemma}.
 Thus, the union of these operators form an operator $P$ on the complement of the half-infinite cuboid $\R^n \setminus \quader[n][1,1][\rho] = \bigcup_{((j,\omega)\neq (i,\eps))} \{\omega x_j > \rho\}$ that restricts to $\R^n \setminus \quader[n][1,1][\rho']$ for all $\rho' > \rho$.
 
 Applying Corollary \ref{Gluing Lower Bounds - Cor} inductively, we conclude that the restriction of $P$ to $\R^n \setminus \quader[n][1,1][\rho']$ is bounded from below for all sufficient large $\rho'$.
 
 To meet the assumptions of the Operator Modification Lemma \ref{Operator Modification - Lemma}, we extend $P$ to a block Dirac operator $\mathtt{P}$ on $M \times \R^n$.
 By increasing $\rho$ if necessary, we may assume that the core of $\mathtt{P}$ is contained in $M \times \rho/2I^n$.
 
 Let $H$, $\Phi$, and $U_{left} \subseteq \R^n \setminus \quader[n][1,1][\rho]$ be as in the Embedding Theorem \ref{Embedding Theorem}.
 By the Operator Modification Lemma \ref{Operator Modification - Lemma}, $\mathtt{P}$ restricts to $P$ on $M \times R\Phi(H \times (a,b))$ for all $1 \leq a < b \leq 2$ and acts like a differential operator on the cylinder variable; more precisely, it is a ${(R\Phi)^{-1}}^\ast \pr_2^\ast \mathcal{C}_c^\infty(a,b)$-derivation.
 Thus, $\mathtt{P}$ also restricts to $P$ on $M \times \bigl(U_{left} \cup R\Phi(H\times(1.8,2))\bigr)$.
 
 By the Operator Modification Lemma \ref{Operator Modification - Lemma}, there is a $Q=Q(R)$ on $\Gamma_c(M \times R\Phi(H \times (1,2);\Spinor_{G_R})$ for all $R\geq1$ that agrees with $P$ on $M \times R\Phi(H\times(1.8,2))$ and decomposes on $M \times R\Phi(H\times(1,1.2))$ as follows
 \begin{equation*}
     \Phi_{R}^\ast Q \deff \Phi_R^\ast Q\restrict_{M \times H \times 1.1R} \boxtimes \id + \evenodd \boxtimes \DiracSt_{\R}.
 \end{equation*}
 Let ${}_R\Theta$ and ${}_R\Psi_n$ be as in the proof of Theorem \ref{Concordance Set is Kan - Theorem}.
 To ease the notation, we abbreviate $\id \times {}_R\Theta$ and $\id \times {}_R\Psi_n$ to ${}_R\Theta$ and ${}_R\Psi_n$, respectively.
 
 The operator $({}_R\Theta \circ{}_R\Psi_n)_\ast(Q)$ is again a symmetric pseudo Dirac operator.
 We chose ${}_R\Theta$ in the proof of Theorem \ref{Concordance Set is Kan - Theorem} such that ${}_R\Theta \circ {}_R\Psi_n \circ {}_R\Phi(m,h,s) = (m,R\Phi(h,1)-(s-R)e_1)$ on $M \times H \times R[1,1+\epsilon))$ for all $1/10 > \epsilon>0$.
 Since ${}_R\Theta \circ {}_R\Psi_n \circ {}_R\Phi(H \times R[1,1+\epsilon)) = \{R(\rho+1-\epsilon)< x_1 \leq R(\rho+1)\}$, we obtain the decomposition
 \begin{equation*}
     ({}_R\Theta {}_R\Psi_n)_\ast(Q) = Q\restrict_{\{x_1 = R(\rho+1)\}} \boxtimes \id + \evenodd \boxtimes \DiracSt_{\R} 
 \end{equation*}
 on $\{R(\rho+1-\epsilon)< x_1 \leq R(\rho+1)\}$.
 
 By $(iv)$ of the Operator Modification Lemma \ref{Operator Modification - Lemma}, the operator $Q$ and hence $({}_R\Theta \circ {}_R\Psi_n)_\ast(P)$ is bounded from below with lower bound $c>0$ that is independent of $R$, provided $R$ is sufficiently large.
 Lemma \ref{PositivityDescendHyperSurf - Lemma} implies that the restriction $Q\restrict_{\{x_1 = R(\rho+1)\}}$ is bounded from below with lower bound bigger than $c/2$, for all sufficient large $R$.
 Thus, the constant extension of $Q\restrict_{\{x_1 = R(\rho+1)\}} \boxtimes \id + \evenodd \boxtimes \DiracSt_{\R}$ to $\{R(\rho+1-\epsilon) < x_1\}$ is bounded from below with lower bound bigger than $c/2$.
 
 It follow from Corollary \ref{Gluing Lower Bounds - Cor} and Lemma \ref{SymmetryUnion - Cor}, that
 \begin{equation*}
    \overline{Q} \deff ({}_R\Theta \circ {}_R\Psi_n)_\ast(Q) \cup Q\restrict_{\{x_1 = R(\rho+1)\}} \boxtimes \id + \evenodd \boxtimes \DiracSt_{\R} 
 \end{equation*}
 is a pseudo Dirac operator on $M \times \bigl({}_R\Theta \circ {}_R\Psi_n \circ {}_R\Phi(H \times R(1,2)) \cup \{R(\rho+1-\epsilon)<x_1\}\bigr)$ that is symmetric and bounded from below, provided $R$ is sufficiently large.
 
 Of course, $\overline{Q}$ restricts to $M \times R\Phi(H \times (1.8,2))$ and agrees there with $P$.
 Applying Corollary \ref{Gluing Lower Bounds - Cor} and \ref{SymmetryUnion - Cor} a second time, the operator $P_{\mathrm{fill}} \deff P \cup \overline{Q}$ is a symmetric pseudo Dirac operator on $M \times \R^n$ that is bounded from below, provided $R$ is sufficiently large.
 
 It remains to show that $P_{\mathrm{fill}}$ is a block Dirac operator. 
 Let $\mathtt{R}_1 = \mathtt{R}_1(R) > R(\rho+1)$ be sufficiently large such that $M \times \mathtt{R}_1I^n$ contains $L(R)$ and ${}_R\Theta \circ {}_R\Psi_n(L(R))$ and such that ${}_R\Theta$ and ${}_R\Psi_n$ restrict to the identity on the complement of $\quader[n][1,1][\mathtt{R}_1]$.
 Then $M \times \mathtt{R}_1I^n$ serves as a core for $P_{\mathrm{fill}}$.
 Since $P_{\mathrm{fill}}$ restricts to $P$ on $M \times \R^n\setminus \quader[n][1,1][\mathtt{R}_1]$ and since, by construction, 
 \begin{equation*}
   P_{\{x_1 > R(\rho+1)\}} = ({}_R\Theta \circ {}_R\Psi_n)_\ast(Q)\restrict_{\{x_1 = R(\rho+1)\}} \boxtimes \id + \evenodd \boxtimes \DiracSt_{\R},  
 \end{equation*}
 we deduce that $P_{\mathrm{fill}}$ has block form.
 Finally, $P_{\mathrm{fill}} - \Dirac_{g_{P_{\mathrm{fill}}}}$ is bounded, because it is true for all the three summands whose union forms $P_{\mathrm{fill}}$.

 Thus, $P_{\mathrm{fill}}$ a filler for the given $n$-horn, so $\InvBlockDirac$ is a Kan set. 
\end{proof}

 \section{The Symbol Map is a Kan Fibration}\label{SymbolMapKanFib - Section}

Recall that every pseudo Dirac operator $P$ comes with an underlying Riemannian metric $g_P$ that determines the principal symbol of $P$ (and also determines the spinor bundle on which $P$ acts). 
We demand that $g_P$ is a block metric and that the decompositions of $P$ and $g_P$ are compatible.
In Section \ref{Section - Foundations of the Operator Concordance Set} we have established that the map that assigns to a pseudo Dirac operator its underlying metric assembles to a cubical map $\BlockDirac \rightarrow \BlockMetrics$ which we refer to as the symbol map.

We use the Kan property of $\InvBlockDirac$ from the previous subsection to establish the main result of this section.
\begin{theorem}\label{SymbolMapKanFib - Theorem}
  The symbol map $\InvBlockDirac \rightarrow \BlockMetrics$ is a Kan fibration.
\end{theorem}
This result is remarkable because it implies that all index theoretic properties are stored in a single fibre $\InvBlockDirac_{g_0}$ because $\BlockMetrics$ is combinatorially contractible, see Proposition \ref{CombContract - Prop}.
The fibre $\InvBlockDirac_{g_0}$ is a Kan set by formal reasons, see Lemma \ref{KanFormalImplication - Lemma}.
When we are determining the homotopy groups of $\InvBlockDirac$, we will make explicit use of the Kan property of the fibre $\InvBlockDirac_{g_0}$, because we will compare the fibre to the space of Fredholm operators on a \emph{single} Hilbert space determined by $g_0$.

Of course, Theorem \ref{SymbolMapKanFib - Theorem} implies Theorem \ref{operatorfundamentallemma} by formal reasons.
However, the proof of Theorem \ref{SymbolMapKanFib - Theorem} relies on $\InvBlockDirac$ being a Kan set.
It goes roughly as follows:
For a given $i$-horn of invertible block Dirac operator and a given block metric, we find  a filler by the Kan property of $\InvBlockDirac$. 
So far, the underlying metric may disagree with the given one.
To correct this mismatch, we will modify the filler so that they do agree.

The modification procedure will rely on two tools.
The first tool is a fibration-kind result inspired from the result of \cite{ebert2017indexdiff}*{Section 4.2} although the methods are quite different. 
It allows us to ``rotate'' the principal symbol of the filler to the required one.
The second tool is based on the gluing constructions from Section \ref{Analytical Properties of Block Operators - Section}.

The first tool is based on the following analytical result and its consequences.
\begin{lemma}\label{GromovLawsonEstimate - Lemma}
 Let $P \in \BlockDirac[n]$ such that all faces $\face{i} P$ are invertible. 
 Then there are constants $c, \texttt{C}>0$, each depending only on the faces $\face{i} P$, and a compact subset $K \subseteq M \times \R^n$, such that all  sections $0 \neq u \in L^2(M \times \R^n; \Spinor_{g_P})$ with $||Pu||_0 \leq \gamma ||u||_0$ for $\gamma < c$ satisfy
 \begin{equation}\label{eq: GromovLawsonEstimate}
  \frac{||u||_{0,K}}{||u||_0} \geq \frac{10}{11}\left(\frac{9}{10} - \mathtt{C}\gamma \right)>0,
 \end{equation}
 where $||u||_{0,K} \deff ||1_K \cdot u||_{0}$.
\end{lemma}
\begin{proof}
 For this proof, it will be more convenient to use the following variant of the Sobolev $1$-norm
 \begin{equation*}
  ||u||_1 \deff ||u||_0 + ||\nabla u||_0,
 \end{equation*}
 which is equivalent to the (usual) Sobolev 1-norm.
 More generally, for each measurable $K \subset M \times \R^n$ we will denote, just for this proof,
 $||u||_{1,K} \deff ||u||_{0,K} + ||\nabla u||_{0,K}$.
 For each compactly supported smooth function $\varphi$, we can estimate $||\varphi u||_1$ against $||u||_{1,\supp \varphi}$ as follows:
 \begin{align*}
  ||\varphi u||_1 &= ||\varphi u||_0 + ||\nabla \varphi u||_0 \\
  &= ||\varphi u||_0 + ||\diff \varphi \otimes u + \varphi \nabla u||_0 \\
  &\leq ||\varphi u||_0 + ||\grad \varphi||_\infty \cdot ||u||_{0,\supp \varphi} + ||\varphi \nabla u||_{0} \\
  &\leq ||\varphi||_\infty ||u||_{1,\supp \varphi} + ||\grad \varphi||_\infty \cdot ||u||_{0,\supp \varphi}.
\end{align*}

Let $\chi\colon \R \rightarrow [0,1]$ be a smooth, even function that is identically $1$ on $(-r,r)$ and such that $\chi^{-1}(\R_{>0}) = (-r',r')$ for some $r'>r>0$.
We assume that $r$ and $r'$ are sufficiently large such that the interior of $M\times rI^n$ contains the core of $P$.
Let $\{\chi^\fatalpha\}$ be the associated partition of unity from Definition \ref{BlockPartOfUnity - Def} and set $K \deff \supp \chi^\mathbf{0}$.
Recall from the proof of Proposition \ref{EllipticEstimateInvInf - Prop} that invertibility of $\face{i}P$ implies, for all sections $u$ supported within $\{\eps x_i > R\}$, the inequality
\begin{equation*}
 ||Pu||_0 \geq c(i,\eps)||u||_1, 
\end{equation*}
 where $c(i,\eps)$ is the product of the operator norm of $(\face{i} P)^{-1} \colon  L^2(M \times \R^{n-1}; \Spinor_{g_{\face{i}P}}) \rightarrow H^1(M \times \R^{n-1}; \Spinor_{g_{\face{i}P}})$ and a constant that comes from the equivalence of the two variants of the Sobolev 1-norms on $H^1(M \times \R^{n-1}; \Spinor_{g_{\face{i}P}})$. 
The second constant only depends on $g_{\face{i}P}$, too.

Recall further that the block form implies $\supp P\chi^\fatalpha u \subseteq \supp \chi^\fatalpha$ by Corollary \ref{InvariantDecomposition - Corollary}.
Indeed, if we use the notation from Corollary \ref{InvariantDecomposition - Corollary}, then $\supp \chi^{\fatalpha}$ is the closure of $V_{r < r'}(\fatalpha)$ and the claim follows from applying Corollary \ref{InvariantDecomposition - Corollary} to $V_{r-\epsilon < r'+\epsilon}$ for all sufficiently small $\epsilon>0$.

 For $\fatalpha \neq \mathbf{0}$, this implies the inequality
\begin{align*}
 ||\chi^\fatalpha u||_1 &\leq \max\{c(i,\eps)^{-1}\} ||P \chi^\fatalpha u||_{0,\supp \chi^\fatalpha} \\
 &\leq \max\{c(i,\eps)^{-1}\} \left(||\chi^\fatalpha Pu||_0 + ||\grad \chi^\fatalpha||_\infty ||u||_{0,\supp \chi^\fatalpha}\right) \\
 &\leq \max\{c(i,\eps)^{-1}\} \left(||Pu||_0 + ||\grad \chi^\fatalpha||_\infty ||u||_0\right).
\end{align*}
Putting these inequalities together yields
 \begin{align*}
 ||u||_1 &\leq ||\chi^\mathbf{0} u||_1 + \sum_{\fatalpha \neq \mathbf{0}} ||\chi^\fatalpha u||_1 \\
 \begin{split}&\leq ||u||_{1,K} + ||\grad{ \chi^\mathbf{0}}||_\infty ||u||_{0,K} \\
 &\qquad \qquad \ + \underbrace{3^n\max\{c(i,\eps)^{-1}\}}_{=: \mathtt{C}}\left( ||Pu||_0 + \max\{||\grad{ \chi^\fatalpha}||_\infty\} ||u||_0\right).
 \end{split}
 \end{align*}
 
 We now assume that $||Pu||_0 \leq \gamma ||u||_0$.
 The previous estimate then yields 
 \begin{align*}
  ||u||_0 &\leq ||u||_0 + ||\nabla u||_0 - ||\nabla u||_{0,K} = ||u||_1 - ||\nabla u||_{0,K} \\
  \begin{split}
     &\leq ||u||_{1,K} + ||\grad{\chi^\mathbf{0}}||_\infty ||u||_{0,K} - ||\nabla u||_{0,K}\\
     &\qquad \qquad \  +\mathtt{C}\left(\gamma ||u||_0 + \max\{||\grad{ \chi^\fatalpha}||_\infty\} ||u||_0\right)  
  \end{split}\\
  &\leq \left(1 + ||\grad{ \chi^\mathbf{0}}||_\infty\right)||u||_{0,K} + \mathtt{C}\left(\gamma ||u||_0 + \max\{||\grad{ \chi^\fatalpha}||_\infty\} ||u||_0\right), 
 \end{align*}
 or equivalently,
 \begin{equation*}
  \bigl(1 - \mathtt{C}(\gamma + \max\{|| \grad{ \chi^\fatalpha} ||_\infty\})\bigr)||u||_0 \leq (1+ ||\grad{ \chi^\mathbf{0}}||_\infty)||u||_{0,K}.
 \end{equation*}
 
 By replacing $\chi^\fatalpha$ with $\chi^\fatalpha(\rho^{-1} \,\cdot\,)$, where $\rho$ is sufficiently large, we may assume\footnote{Keep in mind that rescaling the partition of unity replaces $K$ by $\rho K$.}   
 $|| \grad{ \chi^{\fatalpha}}||_\infty < \min\{1/10, 1/10 \mathtt{C}^{-1}\}$.
 For $u \neq 0$, we end up with the desired inequality:
 \begin{equation*}
  \frac{||u||_{0,K}}{||u||_0} \geq \frac{10}{11}\left(\frac{9}{10} - \mathtt{C}\gamma\right).
 \end{equation*}
 This lower bound is positive if and only if 
 \begin{equation*}
  \gamma < 3^{2-n} \frac{\min\{c(i,\eps)\}}{10} =: c.   \qedhere
 \end{equation*}
\end{proof}

This inequality has important implications concerning the spectrum of $P$.
Recall that a sequence of elements $(v_n)_{n\in \N}$ in a Hilbert space $\mathcal{H}$ is a sequence of \emph{approximate eigenvectors for} $P$ if $||v_n||=1$ and there is a (necessarily unique) $\lambda \in \C$ such that the sequence $(P - \lambda \id)v_n$ converges to zero. 
In that case, $\lambda$ is called an \emph{approximate eigenvalue} (associated to $(v_n)_{n \in \N}$).
Note that the spectrum of a self-adjoint operator consists of approximate eigenvalues only.
\begin{cor}\label{ApproxEigenvalues - Cor}
 Let $P$ be as in the previous lemma. 
 Then each sequence of approximate eigenvectors of $P$ with approximate eigenvalue $\lambda \in (-c,c)$ has a convergent subsequence. In particular, $\lambda$ is an eigenvalue.
\end{cor} 
\begin{proof}
 For a given spectral value $\lambda \in (-c,c)$ of $P$ take a sequence of approximate eigenvalues $(v_n)_{n\in \N}$ associated to it. Since $||Pv_n||_0 < c||v_n||_0 = c$ for almost all $n\in \N$, the fundamental elliptic estimate implies $v_n \in H^1(\Spinor_{g_P})$ and that $(||v_n||_1)_{n \in \N}$ is bounded.
 
 Pick a smooth compactly supported function $\varphi$ that is identically $1$ on $K$, the compact set from Lemma \ref{GromovLawsonEstimate - Lemma}.
 Then $||v_n||_{1,K} \leq ||\varphi v_n||_1 \leq ||\varphi||_{1,1}||v_n||_1$ and the classical Rellich lemma implies that $(\varphi v_n)_{n \in \N}$ has a convergent subsequence with respect to $||\placeholder||_{0}$. 
 Thus $(v_n)_{n\in \N}$ has a convergent subsequence with respect to $||\placeholder||_{0,K}$.
 But inequality \ref{eq: GromovLawsonEstimate} implies that the convergent subsequence also converges with respect to $||\placeholder ||_0$.
 This limit is an eigenvector and $\lambda$ its associated eigenvalue.   
\end{proof}
As an application, we prove the for us very important
\begin{proposition}[Sausage Lemma]
 Let $P$ and $c$ as in Lemma \ref{GromovLawsonEstimate - Lemma}. 
 Then $\mathrm{spec}(P) \cap (-c,c)$ consists of eigenvalues only and the linear hull $\bigoplus_{|\lambda| < c_0} \mathrm{Eig}(P,\lambda)$ of eigenspaces forms a finite dimensional subvector space for all $c_0 < c$. 
\end{proposition}
\begin{proof}
 The first claim was part of the previous corollary.
 
 To prove the second claim, we first show that every eigenspace is finite dimensional.
 Obviously, $\mathrm{Eig}(P,\lambda)$ is a closed subspace of $L^2(\Spinor_{g_P})$, hence a Hilbert space.
 For a given sequence of normalised eigenvectors $(v_n)_{n\in \N}$, the previous corollary implies the existence of a convergent subsequence.
 Thus, the unit ball of $\mathrm{Eig}(P,\lambda)$ is compact, so $\mathrm{Eig}(P,\lambda)$ is finite dimensional.
 
 We now prove that $P$ has only finitely many eigenvalues in $[-c_0,c_0]$ for all $c_0 < c$.
 Assume the contrary. 
 Then there is a sequence $(\lambda_n)_{n \in \N}$ of pairwise different eigenvalues lying in $[-c_0,c_0]$.
 For each $\lambda_n$, pick a normalised eigenvector $v_n$. 
 For $\lambda \deff \limsup_n \lambda_n \in [-c_0,c_0]$, we can pick a subsequence of $(v_n)_{n\in \N}$ that serves as a sequence of approximate eigenvectors for $\lambda$.
 But Corollary \ref{ApproxEigenvalues - Cor} implies that this subsequence has a convergent subsequence $(v_{n_k})_{k \in \N}$. 
 The limit $w$ is an eigenvector corresponding to $\lambda$ and eigenvectors corresponding to different eigenvalues are perpendicular, so we end up with 
 \begin{equation*}
  0 = \langle v_{n_k}, w\rangle \xrightarrow{k \to \infty} \langle w, w \rangle = 1,
 \end{equation*} 
 a contradiction.
\end{proof}
The Sausage Lemma offers the possibility to use spectral-flow-style arguments.
\begin{lemma}
 Let $P \colon X \rightarrow \BlockDirac[n]$ be a map, such that ${\fullGauge_{g_{P},g_0}}_\ast(P) \colon X \rightarrow \mathrm{Hom}(H^1(\Spinor_{g_0}),H^0(\Spinor_{g_0}))$ is continuous and such that $\face[\eps]{i}P$ is a constant map of invertible operators.  
 Assume that there is a constant $c>0$ such that 
 \begin{equation*}
  \mathrm{spec}(P_y) \cap (-c,c) =: \{\lambda_0(y) \leq \dots \leq \lambda_{N_y}(y)\}
 \end{equation*}
 consists of finitely many eigenvalues for all $y \in X$. 
 Then, for each $x \in X$, there is an open neighbourhood $U$ of $x$ such that $N_{(\placeholder)}$ is constant on $U$ and all $\lambda_j \colon U \rightarrow \R$ are continuous functions.
\end{lemma}
\begin{proof}
 For every $c_0$ that lies strictly between $\max\{|\lambda_0(x)|,|\lambda_N(x)|\}$ and $c$, choose a smooth, monotonically increasing function $b \colon \R \rightarrow [-2c,2c]$ that is the identity on $[-c_0,c_0]$.
 The bounded transform $\hat{P} \deff b(\fullGauge_{g_{P},g} \circ P \circ \fullGauge_{g_0,g_{P}})$ is a continuous map of bounded, self-adjoint operators $X \rightarrow B(L^2(\Spinor_{g_0}))$ \cite{ebert2017indexdiff}*{Proposition 3.7} whose values have the same eigenvalues between $(-c_0,c_0)$ as the values of $P$.
 The statement now follows from the fact that eigenvalues of bounded operators depend continuously on the operator with respect to the operator norm \cite{booss1993elliptic}*{Lemma 17.1}.
\end{proof}
\begin{rem}
 The functions of eigenvalues $\lambda_j$ are only defined locally near $x$ because they may leave the interval $(-c,c)$. 
 Also, the number of eigenvalues may not be locally constant.
\end{rem}

 Recall from Theorem \ref{App: ExternalTensor is Continuous - Theorem} that a smooth path $P \colon (0,1) \rightarrow \PseudOpCl^1(\Spinor)$ with underlying path $g_P$ of Riemannian metrics has an extension $P^\extension \in Op^1(L^2(\Spinor_{\susp(g_P)}))$ whose values are not necessarily in the subspace $\PseudOpCl^1(\Spinor_{\susp(g_P)})$.
 However, this is a nice property to have, so we will give it a name.
 \begin{definition}
  Let $X\subseteq \R^n$ be open and let $P \colon X \rightarrow \PseudOpCl^1(\Spinor)$ be a smooth map with underlying map $g_P$ of Riemannian metrics. Then $P$ called \emph{suspensionable} if the extension $P^\extension$ given by
  \begin{equation*}
      (P^\extension\sigma)(n,x) = P_x(\sigma(\placeholder,x) )
  \end{equation*}
  is an element of $\PseudOpCl^1(\Spinor_{\susp(g_P)})$.
 \end{definition}
 The first tool is the following fibration-like result.
\begin{proposition}[Path Lifting Property]\label{PathLiftingProp - Prop}
 Consider $\BlockMetrics[n]$ as a topological space equipped with the smooth Fréchet topology and $\InvBlockDirac[n]$ as a topological space equipped with the subspace topology of $\PseudOpCl^1(\Spinor)$. 
 Let $\gamma \colon [0,1] \rightarrow \BlockMetrics[n]$ be a smooth curve that is constant near the end points and satisfies 
 \begin{equation*}
     \gamma_t|_{\{\eps x_i > R\}} = \degen{i}\face[\eps]{i}\gamma_0|_{\{\eps x_i > R\}}
 \end{equation*}
 for some $R>0$, all $(i,\eps) \in \{1,\dots,n\}\times \Z_2$, and all $t \in \lbrack 0,1\rbrack$.
 
 Then each commutative outer square 
 \begin{align*}
  \xymatrix{ \{0\} \ar[d] \ar[r]^-{Q_{0}} & \InvBlockDirac[n] \ar[d] \\ 
  \lbrack 0,1\rbrack \ar[r]^\gamma \ar@{-->}[ru]^{Q} & \BlockMetrics[n] }
 \end{align*}
 can be filled with a smooth map that is constant near the boundary and such that the upper triangle commutes and the lower right triangle commutes up to homotopy relative to the boundary. 
 In fact $g_{Q_{(\, \cdot \, )}}$ and $\gamma_{(\, \cdot \,)}$ differs by a (non-diffeomorphical) reparametrisation.
 Furthermore, we may assume that the lift is suspensionable and satisfies $\face[\omega]{j}Q_t = \face[\omega]{j}Q_{0}$ for all $(j,\omega) \in \{1,\dots,n\} \times \Z_2$.
\end{proposition}
In contrast to the proof in \cite{ebert2017indexdiff}*{Section 4.2} our argument is quite geometrical. 
The idea is the following one: For an ordered pair of two Riemannian metric $(g_1,g_2)$ and a non-zero covector $\xi \in T^\vee M$, there is a unique $A \in \GL^+(d)$ that stretches and rotates ${g_1}_\sharp(\xi)$ to ${g_2}_\sharp(\xi)$ in the plane spanned by the two vectors and that acts identically on the orthogonal complement.
We find a lift $b_{g_1,g_2}(\xi)$ in $\Cliff[d][0]$ such that its induced linear map $b_{g_1,g_2}(\xi)(\placeholder)b_{g_1,g_2}(\xi)^\ast$ on $\R^d$ agrees with $A$.
We quantise the symbol to get a pseudo differential operator $B_{g_1,g_2}$.
If $B_{\gamma_0,\gamma_t}PB_{\gamma_0,\gamma_t}^\ast$ were a curve of invertible operators, then it would be the desired lift.
Unfortunately, it is not clear that this curve is invertible, so we need to stop the curve from time to time and ``push the spectral values away from zero''.

Let us now introduce the required concepts.
\begin{definition}
 For $h \in \BlockMetrics[n]$ and $R>0$, denote by $\BlockMetrics[n]_{h,R}$ the space of all block metrics that agree with $h$ outside of $M \times RI^n$ equipped with the smooth compact open topology.
\end{definition}

\begin{lemma}\label{RotationQuantification - Lemma}
 For all $h \in \BlockMetrics[n]$ and all $R >0$, there is a smooth map $B \colon \BlockMetrics[n]_{h,R} \times \BlockMetrics[n]_{h,R} \rightarrow \PseudOp^0(\Spinor,\Spinor)$ that satisfies the following three conditions:
 \begin{enumerate}
  \item[(i)] The following diagram commutes
  \begin{equation*}
   \xymatrix{ \BlockMetrics[n]_{h,R} \times \BlockMetrics[n]_{h,R} \ar[r]^-B \ar[dr]_{\id} & \PseudOp^0(\Spinor, \Spinor) \ar@<1ex>[d]^{\mathrm{target}} \ar@<-1ex>[d]_{\mathrm{source}}  \\ & \BlockMetrics[n]_{h,R} \times \BlockMetrics[n]_{h,R}. }
  \end{equation*}
  \item[(ii)] For all $g_1, g_2 \in \BlockMetrics[n]_{h,R}$, the operator $B_{g_1,g_2} - \preGauge_{g_1,g_2}$ is compactly supported.
  \item[(iii)] For all $g_1,g_2 \in \BlockMetrics[n]_{h,R}$, we have
  \begin{equation*}
   \symb_0(B_{g_1,g_2}) \circ g_{1\, \sharp}(\, \cdot \,) \cdot_{g_1} \circ \, \symb_0(B_{g_1,g_2}^\ast) =  g_{2 \, \sharp}(\, \cdot \,) \cdot_{g_2}, 
  \end{equation*}   
              where $\cdot_{g_j}$ is the Clifford left action on $\Spinor_{g_j}$.
 \end{enumerate}
\end{lemma}
\begin{proof}
 Let $\pi \colon \mathrm{Pos} \rightarrow M\times \R^n$ be the bundle whose fibres consist of positive definite, symmetric bilinear forms and $\Spinor \rightarrow \mathrm{Pos}$ be the fibre bundle whose fibre of $g_x$ is $\Spinor_{g_x}$.
 Similarly, we define $Cl(T(M\times \R^n), \placeholder ) \rightarrow \mathrm{Pos}$ to be the bundle whose fibre of $g_x$ is the Clifford algebra $Cl(T_x(M \times \R^n),g_x)$.
 We denote with $||\placeholder|| \deff ||\placeholder||_{g_x}$ the norm on $Cl(T_x(M \times \R^n),g_x)$ induced by $g_x$ and we denote with $\cdot_{g_x}$ the Clifford multiplication and left action of $Cl(M,g_x)$ on $\Spinor_{g_x}$.
 Recall furthermore that 
 \begin{equation*}
    \Phi_{g_{1,x},g_{2,x}} \deff Cl(\tau_{g_{1,x},g_{2,x}}) \colon Cl(T_x(M \times \R^n),g_{1,x}) \rightarrow Cl(T_x(M \times \R^n),g_{2,x}) 
 \end{equation*}
 denotes the algebra isomorphism induced by the pre-gauge map on the tangent spaces.
 Using this notation, we set
 \begin{align*}
    \begin{split}
      b_{g_1,g_2}(\xi) := &  
  \frac{||\Phi_{g_2,g_1}g_{2\, \sharp}(\xi)||^{1/2}}{||g_{1\, \sharp}(\xi)||^{1/2}} \cdot \\
  &\quad \frac{\Phi_{g_2,g_1}g_{2\, \sharp}(\xi)/||\Phi_{g_2,g_1}g_{2\, \sharp}(\xi)|| + g_{1\, \sharp}(\xi)/||g_{1\, \sharp}(\xi)||}{||\Phi_{g_2,g_1} g_{2\, \sharp}(\xi)/||\Phi_{g_2,g_1} g_{2\, \sharp}(\xi)|| + g_{1\, \sharp}(\xi)/||g_{1\, \sharp}(\xi)|| \, ||} \cdot_{g_1} 
  \frac{g_{1 \, \sharp}(\xi)}{||g_{1\, \sharp}(\xi)||} \cdot_{g_1} .
  \end{split} 
 \end{align*}
 The map $b$ takes values in even and $\Cliff$-right linear symbols, more precisely,
 \begin{equation*}
     b_{g_1,g_2} \in \{\sigma \in \Symb_0(\Spinor_{g_1},\Spinor_{g_2}) \, : \, \sigma \text{ is even and } \Cliff[d][0]\text{-right linear} \},
 \end{equation*}
 see Definition \ref{SetOfSymbols - Def} for a reminder of $\Symb_0(\Spinor_{g_1},\Spinor_{g_2})$.
 
 Denote by $\mathrm{Pos}^2 = \mathrm{Pos} \times_M \mathrm{Pos}$ the fibre bundle of pairs of positive definite, symmetric, bilinear forms on $T(M \times \R^n)$.
 We consider it as fibre bundle over $\mathrm{Pos}$, whose bundle projection is the projection to the first component. 
 Furthermore, let
 \begin{equation*}
  F \deff \{(v,w) \in T(M\times \R^n) \setminus \{0\} \oplus T(M \times \R^n) \setminus \{0\} \, : \, v \neq - \lambda w \text{ for all  } \lambda > 0\} 
 \end{equation*}
 be the bundle of pairs of tangent vectors whose cones are not antipodal half lines.
 
 The map $b \colon \BlockMetrics[n]_{h,R} \times \BlockMetrics[n]_{h,R} \rightarrow \Symb_0^\times(\Spinor,\Spinor)$ is continuous because it is induced by the following sequence of bundle maps
 \begin{equation*}
     \xymatrix{\pi^\ast T^\vee(M \times \R^n)\setminus \{0\} \times_M \mathrm{Pos}^2 \ar[r] & \pi^\ast F \ar[r] & \Cliff[ ][ ](M \times \R^n,\placeholder) \ar[r] & \mathrm{Hom}(\Spinor,\Spinor). } 
 \end{equation*}
 
 The first map is given by $(g_1,g_2) \mapsto (g_{1\, \sharp}(\xi), \Phi_{g_2,g_1}(g_{2\, \sharp}(\xi)))$. 
 Notice that this map indeed takes values in $\pi^\ast F$ because $g_{1\, \sharp}(\xi) = - \Phi_{g_2,g_1}(\lambda g_{2 \, \sharp}(\xi))$ for some $\lambda > 0$ implies that $g_1$ is not positive definite:
 \begin{align*}
  g_1\bigl(g_{1 \, \sharp}(\xi), g_{1 \, \sharp}(\xi)\bigr) &= - \lambda g_1\bigl(g_{1\, \sharp}(\xi),\Phi_{g_2,g_1}(g_{2\, \sharp}(\xi))\bigr) \\
  &= - \lambda g_1\bigl({g_1}_\sharp(\xi), \Phi_{g_2,g_1}({g_{2}}_\sharp \circ g_1^\flat \, {g_1}_\sharp(\xi)\bigr) \\
  &= - \lambda g_1\bigl({g_1}_\sharp(\xi), \Phi_{g_2,g_1}^3 ({g_1}_\sharp(\xi))\bigr) < 0 
 \end{align*}
 because $\Phi_{g_2,g_1}|_{T(M\times \R^n)} = \tau_{g_1,g_2}$ is the pre-gauge map, which is positive definite and self-adjoint with respect to $g_1$.
  The second map is given by 
  \begin{equation*}
   (v,w) \mapsto \frac{||w||^{1/2}_{g_1}}{||v||^{1/2}_{g_1}} \frac{w/||w||_{g_1} + v/||v||_{g_1}}{|| w/||w||_{g_1} + v/||v||_{g_1} ||_{g_1}} \cdot_{g_1} \frac{v}{||v||_{g_1}} \cdot_{g_1}, 
  \end{equation*}
  which is clearly continuous.
 The third map is just the action of the Clifford algebra bundle on the spinor bundle, which is even isometric. 
 Note, in particular, that $b$ takes values in invertible symbols.
   
 Away from $M \times RI^n$, the map $b$ is the constant map with value $\id_{\Spinor_{h_x}}$.
 Recall from Lemma \ref{SymbolExactSequence - Lemma} the continuous (with respect to the amplitude topology) section $\rho \colon \Symb_0(\Spinor_{g_1}) \rightarrow \PseudOp^0(\Spinor_{g_1})$.
 We set $\tilde{B}_{g_1,g_2} \deff \rho(b_{g_1,g_2}) - \rho(b_{g_1,g_1}) + \id$.
 So far, $\tilde{B}_{g_1,g_2}-\id$ might not have compact support. 
 To achieve this, we pick smooth compactly supported functions $\varphi, \psi \colon M \times \R^n \rightarrow [0,1]$ such that $\varphi \equiv 1$ on $M \times RI^n$ and $\psi \equiv 1$ on $\supp \varphi$. 
 The operator $\hat{B}_{g_1,g_2} \deff \psi \tilde{B}_{g_1,g_2} \varphi + \id \cdot (1-\varphi)$ has the same principal symbol as $\tilde{B}_{g_1,g_2}$ and differs from the identity by a compactly supported operator.
 Recall the function $\alpha_{g_2,g_1}^2 = \diff \mathrm{vol}(g_2)/\diff \mathrm{vol}(g_1)$. 
 Finally, we define the desired map
 \begin{align*}
  B \colon \BlockMetrics[n]_{h,R} \times \BlockMetrics[n]_{h,R} &\rightarrow \PseudOp^0(\Spinor,\Spinor), \\
 (g_1,g_2) &\mapsto \preGauge_{g_1,g_2} \circ \alpha_{g_1,g_2} \hat{B}_{g_1,g_2} \in \PseudOp^0(\Spinor_{g_1},\Spinor_{g_2}).
 \end{align*}
 
 By the definition of the topology of $\PseudOp^0(\Spinor,\Spinor)$, the map $B$ is smooth if and only if $\hat{B}$ is smooth, which follows from the construction.
 It is clear, that $B$ satisfies $(i). $
 
 The composition of compactly supported pseudo differential operators with vector bundle morphisms are still compactly supported. 
 Furthermore, $g_1 = g_2$ on the complement of $M \times RI^n$ so that $\alpha_{g_1,g_2}=1$ on this set.
 Hence, $B_{g_1,g_2} - \preGauge_{g_1,g_2}$ is compactly supported, which is $(ii)$.
 
 We now prove $(iii)$. 
 Using the relation $(v+w)\cdot w = v \cdot (v+w)$ for unit-length vectors $v,w$ in the Clifford algebra, we derive
 \begin{equation*}
  b_{g_1,g_2}(\xi) \cdot_{g_1} g_{1\, \sharp}(\xi) \cdot_{g_1} b_{g_1,g_2}(\xi)^\ast \cdot_{g_1} = g_{2 \, \sharp}(\xi) \cdot_{g_1}.
 \end{equation*}  
 Recall that the pre-Gauge map is an isometric vector bundle map and $Cl(M,g_1)$-$Cl(M,g_2)$-equivariant. 
 It induces a linear map $Op(\preGauge_{g_1,g_2}) \colon L^2(\Spinor_{g_1}) \rightarrow L^2(\Spinor_{g_2})$.
 The adjoints of the two maps are related by the formula
 \begin{equation*}
  Op(\preGauge_{g_1,g_2})^\ast = Op(\alpha_{g_2,g_1}^2 \preGauge_{g_1,g_2}^\ast) = Op(\alpha_{g_2,g_1}^2 \preGauge_{g_2,g_1})
 \end{equation*}
 because the underlying volume induced by the metrics is different.
 Symbol calculus now implies the third claim
 \begin{equation*}
  \symb_0(B_{g_1,g_2}) \circ g_{1\, \sharp}(\, \cdot \,) \cdot_{g_1} \circ \, \symb_0(B_{g_1,g_2}^\ast) = g_{2 \, \sharp}(\, \cdot \,) \cdot_{g_2}. \qedhere
 \end{equation*} 
\end{proof}
 
 Since compactly supported pseudo differential operators have compactly supported adjoint, the previous proposition has the following important consequence.

\begin{cor}\label{BPB is block - Cor}
 For every $g \in \BlockMetrics[n]_{h,R}$ and every block Dirac operator $P$ with underlying metric $g_P \in \BlockMetrics[n]_{h,R}$ the operator $B_{g_P,g} \circ P \circ B_{g_P,g}^\ast$ is a block Dirac operator with underlying metric $g$. 
\end{cor}
\begin{proof}
 The composition $B_{g_P,g}PB_{g_P,g}^\ast$ lies in $\PseudOpCl^1(\Spinor_g)$ because $B_{g_P,g}$ is an actual pseudo differential operator of order zero, $\PseudOpCl^1(\Spinor_g)$ is local, and compositions on actual pseudo differential operators are actual pseudo differential operators.
 The composition is clearly a symmetric operator.
 
 Next we show that the composition is of block form.
 Let $r_2>r_1>R$ be a real numbers such that $M \times r_1I^n$ contains the core of $P$, the support of $B_{g_P,g} - \preGauge_{g_P,g}$ and $B_{g_P,g}^\ast - \preGauge_{g_P,g}^\ast$ and such that $M \times rI^n$ and such that $P$ decomposes away from $M \times r_1I^n$.
 Lemma \ref{CoreLowerBound - Lemma} applied to $P$ now implies that $M \times r_2I^n$ serves as a core for $B_{g_P,g}PB_{g_P,g}^\ast$.
 
 Since $g_P=g=h$ away from $M \times r_1I^n$, we deduce that $B_{g_P,g}PB_{g_P,g}^\ast$ decomposes there as $P$.
 
 The operator $B_{g_P,g}\Dirac_{g_P}B_{g_P,g}^\ast - \Dirac_g$ is a pseudo differential operator of order zero with compact support, hence a bounded operator.
 The calculation 
 \begin{align*}
  B_{g_P,g}PB_{g_P,g}^\ast-\Dirac_g &= B_{g_P,g}PB_{g_P,g}^\ast- B_{g_P,g}\Dirac_{g_P}B_{g_P,g}^\ast + B_{g_P,g}\Dirac_{g_P}B_{g_P,g}^\ast -\Dirac_g\\
  &= B_{g_P,g}\left(P - \Dirac_P\right)B_{g_P,g}^\ast + \left(B_{g_P,g}\Dirac_{g_P}B_{g_P,g}^\ast -\Dirac_g\right)
 \end{align*}
 implies that $B_{g_P,g}PB_{g_P,g}^\ast - \Dirac_g$ is the sum of two bounded operators, hence bounded.
 
 In conclusion, $B_{g_P,g}PB_{g_P,g}^\ast$ is a block Dirac operator with underlying metric $g$.
\end{proof}

Unfortunately, the information that the eigenvalues of a continuous family of block Dirac operators are continuous functions is not good enough for the arguments to come. 
More refined information for the case $g \mapsto B_{g_P,g}PB_{g_P,g}^\ast$ is provided by the next lemma. 
To formulate it, we recall some notations taken from \cite{kato1966perturbation}*{Section IV.2}.
\begin{definition}
 Let $S, T \colon X \rightarrow Y$ be two closed, unbounded operators between the Banach spaces $X$ and $Y$.
 We define $\delta(S,T)$ to be the smallest number $\delta$ that satisfies for all $u \in \mathrm{dom}(S)$ the inequality
 \begin{equation*}
  \mathrm{dist}\bigl((u,Su);\mathrm{Graph}(T)\bigr) \leq \delta (||u||_X + ||Su||_Y) \overset{\text{def}}{=} \delta \cdot ||(u,Su)||_{X \oplus Y}  
 \end{equation*}  
 and $\hat{\delta}(S,T) \deff \max\{\delta(S,T); \delta(T,S)\}$.
\end{definition}
It is shown in $\mathrm{loc}.$ $\mathrm{cit}.$ that $\hat{\delta}$ is a (incomplete) metric on the set $C(X,Y)$ of all closed, unbounded operators from $X$ to $Y$.
The distance $\hat{\delta}$ behaves nicely with compositions of invertible bounded operators.
\begin{lemma}\label{Estimation GrassmannDist - Lemma}
 Let $T \in C(X,Y)$ and let $A \colon X \rightarrow X$, $B\colon Y \rightarrow Y$ be bounded invertible operators such that $A$ and $A^{-1}$ preserve the domain of $T$.
 
 Then 
 \begin{align*}
  &\delta(T,BT) \leq ||\id - B||, &\delta(BT,T) \leq ||\id - B||,  \\
  &\delta(T,TA) \leq ||\id - A^{-1}||, &\delta(TA,T) \leq ||\id - A||.
 \end{align*}
\end{lemma}

\begin{proof}
 Let $u \in \dom T$ be non-zero. For $v = u \in \dom BT$ we have 
 \begin{align*}
  \mathrm{dist}\bigl((u,Tu);\mathrm{Graph}(BT)\bigr) \leq ||(u,Tu)-(v,BTv)|| &= ||(0,(\id -B)Tu)|| \\
  &\leq ||\id - B||\cdot||u||.
 \end{align*}
 If we choose $v = A^{-1}u \in \dom T$, then
 \begin{align*}
  \mathrm{dist}\bigl((u,Tu);\mathrm{Graph}(TA)\bigr) \leq ||(u,Tu) - (v,TAv)|| &= ||((\id - A^{-1})u,0)|| \\
  &\leq ||\id - A^{-1}||\cdot ||u||.
 \end{align*}
 The other two inequalities are derived in a similar manner.
\end{proof}

Important for us is the fact that spectral gaps of operators are, in a more general sense, lower semi-continuous (they cannot suddenly decrease but they may suddenly increase) with respect to the topology generated by $\hat{\delta}$.
\begin{theorem}\cite{kato1966perturbation}*{Theorem IV.3.1}\label{Kato SpectralDistance - Theorem}
 Let $T \in C(X,Y)$ be a closed operator and let $R \colon \zeta \mapsto (T-\zeta)^{-1} \in B(Y,X)$ be its resolvent function. 
 For every compact subset $\Gamma$ of the resolvent set $\rho(T)$ of $T$, there is a $\delta > 0$ such that if $\hat{\delta}(S,T) < \delta$, then $\Gamma$ is a subset of $\rho(S)$.
 An explicit choice for $\delta$ is 
 \begin{equation*}
  \delta = \min_{\zeta \in \Gamma} \frac{1}{2} (1+|\zeta|^2)^{-1}(1+||R(\zeta)||^2)^{-1/2}. 
 \end{equation*}  
\end{theorem}

\begin{proof}[Proof of the Path Lifting Property]
 Let the following commutative square be given
 \begin{align*}
  \xymatrix{ \{0\} \ar[d] \ar[r]^-{Q_0} & \InvBlockDirac[n] \ar[d] \\ 
  \lbrack 0,1\rbrack \ar[r]^\gamma  & \BlockMetrics[n]. }
 \end{align*}
 By assumption, $\gamma$ factors through $\BlockMetrics_{\gamma_0,R}$ for some $R>0$.
 Let $0<c<1$ be a constant smaller than the corresponding constant from the Sausage Lemma for the faces $\{\face{i}Q_0\}$.  
 Choose a positive number $\mu_0 < c/8$. 
 Since $\gamma$ is uniformly continuous there is an $N \in \N$ such that all $s,t \in [0,1]$ with $|s-t| < 1/N$ satisfy
 \begin{equation}\label{eq: B-preGauge}
  ||B(\gamma_s,\gamma_t) - \fullGauge_{\gamma_s,\gamma_t}||_{0,0} = || B(\gamma_s,\gamma_t)\fullGauge_{\gamma_s,\gamma_t}^\ast - \id ||_{0,0} < \mu_0/8. 
 \end{equation}   
 
 We will prove the following slightly stronger statement.
 For every $P \in \InvBlockDirac[n]$ that contains $(-2\mu_0,2\mu_0)$ within its resolvent set, we construct a smooth path $Q \colon [0,1/N] \rightarrow \InvBlockDirac[n]$ that is constant near the boundary, lifts the restriction of $\gamma|_{[0,1/N]}$ with initial value $Q_0 = P$ up to reparametrisation, satisfies $\face[\omega]{j}Q_t = \face[\omega]{j}Q_0$, and whose endpoint contains $(-2\mu_0,2\mu_0)$ in its resolvent set.
 
 We will apply this construction iteratively over the partition $[k/N,(k+1)/N]$ starting with the given operator $Q_0$, which satisfies our construction assumptions.
 The union of all these curves will then be the desired homotopy lift.

 Let us now describe the construction in more detail.
 Fix, once and for all, a smooth, odd, injective function $\varphi \colon \R\setminus \{0\} \rightarrow \R \setminus [-3\mu_0,3\mu_0]$  that is the identity on the complement of $(-c/2,c/2)$.
 The construction is divided into two steps.
 The first step is to extend the inital value $Q_0$ to a path, the second step is adjusting the eigenvalues of the endpoint $Q_{1/N}$.
 
 Let us start with the first step:
 Abbreviate $B_{\gamma_0,\gamma_s}$ to $B_s$ and consider on $[0,1/N]$ the map $t \mapsto B_tPB_t^\ast$.
 By Lemma \ref{RotationQuantification - Lemma} and Corollary \ref{BPB is block - Cor} this map is a smooth map with values in $\BlockDirac[n]$.
 
 To see that $B_tPB_t^\ast$ is invertible, we argue as follows. 
 The distance between the spectrum of $Q_0$ and $\Gamma \deff [-\mu_0,\mu_0]$ is bigger than $\mu_0$.
 Thus we can estimate the constant $\delta$ in Theorem \ref{Kato SpectralDistance - Theorem} by
 \begin{align*}
  \delta &= \min_{\zeta \in \Gamma} \frac{1}{2} (1 + |\zeta|^2)^{-1}(1 + ||R(\zeta)||^2)^{-1/2} \\
  &\geq \frac{1}{2} (1 + \mu_0^2)^{-1}(1 + \mu_0^{-2})^{-1/2} \\
  &\geq \frac{1}{2}\mu_0(1+\mu_0^2)^{-3/2} > \frac{1}{3} \mu_0.
 \end{align*}
 The last inequality holds because $\mu_0 < 1/8$.
 
 On the other hand Lemma \ref{Estimation GrassmannDist - Lemma} and inequality (\ref{eq: B-preGauge}) together with the estimate from the geometric series yield for all $t \in [0,1/N]$ the estimate
 \begin{align*}
  \hat{\delta}(\fullGauge_t^\ast B_tPB_t^\ast \fullGauge_t,P) &\leq \hat{\delta}(\fullGauge_t^\ast B_tPB_t^\ast \fullGauge_t, \fullGauge_t^\ast B_tP) + \hat{\delta}(\fullGauge_t^\ast B_tP,P) \\
  &\leq \max\{||B^\ast_t - \fullGauge_{\gamma_t,\gamma_0}||_{0,0},||(B^\ast_t)^{-1} - \fullGauge_{\gamma_0,\gamma_t}||_{0,0}\} + ||B_t - \fullGauge_{\gamma_0,\gamma_t}||_{0,0}\\
  &\leq \mu_0/8 \cdot (1 - \mu_0/8)^{-1} + \mu_0/8 < \mu_0/3 < \delta,
 \end{align*}
 so all $B_t P B_t^\ast$ have $[-\mu_0,\mu_0]$ in their resolvent set.
 
 Compose $t \mapsto B_t P B_t^\ast$ with a smooth, monotone increasing, surjective self map on $[0,1/N]$ such that the composition is a map $[0,1/N] \rightarrow \InvBlockDirac[n]$ that is constant near the boundary.
 By abuse of notation, we will denote the composition again with $BPB^\ast$.
 
 Now comes the second step.
 Let $\lambda_0^+ \leq \dots \leq \lambda_{m^+}^+$ be all positive eigenvalues and let $\lambda_0^- \geq \dots \geq \lambda^-_{m^-}$ be all negative eigenvalues of $BQ_0B^\ast$ inside $(-c,c)$.
 Pick eigenvectors $v_j^\pm$ corresponding to $\lambda_j^\pm$ such that they form an orthonormal system and pick smooth compactly supported sections $w_j^\pm$ that satisfy $||v_j^\pm - w_j^\pm||_0 < \mu_0/C \cdot 10^{-j}$, where $C \geq 1$ is a constant to be determined later.
 Let $\Phi \colon [1/N, 1 + 1/N] \times \R\setminus \{0\} \rightarrow \R\setminus \{0\}$ be the homotopy between $\id$ and $\varphi$ given by convex combination and compose the first entry with a smooth, monotonically increasing, surjective self-map on $[1/N, 1+1/N]$ such that $\Phi(\,\cdot \, ,x) = x$ near $1/N$ and $\Phi(\,\cdot\, ,x) = \varphi(x)$ near $1+1/N$.
 Denote the result again with $\Phi$.
 Finally, define 
 \begin{gather*}
   \Theta \colon [1/N, 1+ 1/N] \rightarrow \InvBlockDirac[n][\Spinor_{g_{B_{1/N}Q_{0}B_{1/N}^\ast}}], \\
   t \mapsto B_{1/N}Q_{0}B_{1/N}^\ast + \sum_{j=0}^{m^\pm} (\Phi(t,\lambda_j^\pm)-\lambda_j^\pm)\langle \cdot, w_j^\pm\rangle w_j^\pm 
 \end{gather*}
 and the auxiliary map
 \begin{gather*}
   \widetilde{\Theta} \colon [1/N, 1+1/N] \rightarrow \PseudOpCl^1(\Spinor_{g_{B_{1/N}Q_{0}B_{1/N}^\ast}}), \\
   t \mapsto B_{1/N}Q_{0}B_{1/N}^\ast + \sum_{j=0}^{m^\pm} (\Phi(t,\lambda_j^\pm)-\lambda_j^\pm)\langle \cdot, v_j^\pm\rangle v_j^\pm. 
 \end{gather*}
 
 For all $t \in [1/N,1+1/N]$, the two operators $\Theta_t$ and $\widetilde{\Theta}_t$ differ from $B_{1/N}Q_{0}B_{1/N}^\ast$ by a self-adjoint infinitely smoothing operator of finite rank. 
 They are therefore self-adjoint pseudo differential operators of order $1$ with the same principal symbol as the Dirac operator of the underlying metric.
 They also differ from this Dirac operator by a bounded operator, because $B_{1/N}Q_{0}B_{1/N}^\ast$ does it.
 Since all $w_j^\pm$ are compactly supported, $\Theta_t$ is of block form and has the same faces as $B_{1/N} Q_{0} B_{1/N}^\ast$; hence, $\Theta_t$ is a block Dirac operator.
 
 We prove invertibility of $\Theta$ by showing that $\widetilde{\Theta}$ is invertible and that $\Theta$ is a sufficiently small perturbation of $\widetilde{\Theta}$.
 On the positive subspace of $B_{1/N}Q_{0}B_{1/N}$, i.e., the one on which $B_{1/N}Q_{0}B_{1/N}^\ast$ is positive definite,
 we have, by construction,
 \begin{equation*}
  \widetilde{\Theta_t} \geq \Phi(t,\lambda_0^+) \id.
 \end{equation*}   
 A similar result holds for the negative subspace, so the interval $(\Phi(t,\lambda_0^-),\Phi(t,\lambda_0^+))$ does not intersect the spectrum of $\widetilde{\Theta_t}$.
 In particular, each $\widetilde{\Theta_t}$ is invertible.
 
 It remains to show that $\Theta$ is a sufficiently small perturbation of $\widetilde{\Theta}$.
 First observe
 \begin{align*}
  &\quad \ ||\langle \, \cdot \, , w_j^\pm\rangle w_j^\pm - \langle \, \cdot \, , v_j^\pm\rangle v_j^\pm||_{0,0}\\
   &\leq ||\langle \, \cdot \, , w_j^\pm - v_j^\pm \rangle||_{\mathrm{dual}} ||w_j^\pm||_0 + ||\langle \, \cdot \, , v_j^\pm\rangle ||_{\mathrm{dual}} ||v_j^\pm-w_j^\pm||_0 \\
  &\leq \left( 2 + ||v_j - w_j||_0\right)||v_j + w_j||_0 < 2.5 \mu_0/C 10^{-j}. 
 \end{align*}
 This implies
 \begin{align*}
  ||\Theta_t - \widetilde{\Theta_t}||_{0,0} &\leq \sum_{j=0}^{M,N} \underbrace{|\Phi(t,\lambda_j^\pm)-\lambda_j^\pm|}_{\leq c} 2.5\mu_0/C 10^{-j}\\
  &\leq c \cdot 2.5 \cdot 2 \cdot 10/9 \cdot \mu_0/C\\
  &\leq 6  \mu_0/C.
 \end{align*}
 Thus, the spectrum of $\Theta_t$ does not intersect $(\Phi(t,\lambda_0^-) + 6\mu_0/C,\Phi(t,\lambda_0^+) - 6\mu_0/C)$. 
 If we choose $C$ sufficiently large, for example $C = 12$ and estimate $\tilde{\Theta}_t$ against its lowest positive eigenvalue $\Phi(t,\lambda_0^+)$ on its positive subspace, we derive the following chain of inequalities
 \begin{align*}
  \Theta_t &= \Theta_t + \tilde{\Theta}_t - \tilde{\Theta}_t \geq \left(\Phi(t,\lambda_0^+) - ||\tilde{\Theta}_t - \Theta_t||_{0,0} \right)\id \\
  &\geq \left(\lambda_0^+ - ||\tilde{\Theta}_t - \Theta_t||_{0,0}\right)\id \geq \left(2\mu_0 - ||\tilde{\Theta}_t - \Theta_t||_{0,0}\right)\id \\
  &\geq (2\mu_0 - 6\mu_0/C) \id \geq \mu_0 \id.
 \end{align*}  
 A similar result holds on the negative subspace, so $\Theta_t$ is invertible for all $t \in [1/N,1+1/N]$.
 Using $\Phi(1+1/N,\lambda_0^+) = \varphi(\lambda_0^+) \geq 3\mu_0$ and
 \begin{align*}
  \Theta_{1+1/N} &= \Theta_{1+1/N} - \tilde{\Theta}_{1+1/N} + \tilde{\Theta}_{1+1/N} \\
  &\geq \left(3\mu_0 - ||\Theta_{1+1/N} - \tilde{\Theta}_{1+1/N}||_{0,0}\right)\id \\
  &\geq (3\mu_0 - \mu_0/2) \id > 2\mu_0 \id,  
 \end{align*}
  we see that the spectrum of $\Theta_{1/N}$ does not intersect $(-2\mu_0,2\mu_0)$.
 
 The desired homotopy lift is now given by
 \begin{equation*}
  \xymatrix@+3.3em{Q \colon [0,1/N] \ar[r]^{t \mapsto (N+1)t} & [0,1+1/N] \ar[r]^{B_{(\cdot )}Q_{0}B_{(\cdot)} \cup \Theta} & \InvBlockDirac[n]. }
 \end{equation*}
 
 Since $Q_{1/N}$ satisfies the assumption of the described construction, we can construct a lift on $[1/N,2/N]$ starting with $Q_{1/N}$ instead of $Q_{0}$ (but with the same $\mu_0$).
 If we repeat this construction $(N-2)$-times further and glue the results together, we end up with the desired (homotopy) lift $Q \colon [0,1] \rightarrow \InvBlockDirac[n]$. 
 Since the homotopy between $g_Q$ and $\gamma$ is given by reparametrisation relative endpoints, the map $Q$ is indeed a homotopy lift.
 
 It remains to prove that our construction is a suspensionable path.
 More precisely, if $Q$ also denotes the constant extension of the constructed lift, we need to show that $Q^\extension \in \PseudOpCl^1(M\times \R^n \times \R;\Spinor_{\susp(\gamma)})$.
 Since $\PseudOpCl^1(\Spinor)$ is local, see Proposition \ref{PropertiesASPseudos - Proposition}, it suffices to prove the statement only for each piece of the construction.
 
 To verify $\bigl(B_{(\cdot)}Q_0B_{(\cdot)}^\ast\bigr)^\extension \in \PseudOpCl^1(\Spinor)$, take a sequence of actual pseudo differential operators $P_n \in \PseudOp^1(\Spinor_{\gamma_0})$ that approximates $Q_0$ with respect to the Fréchet structure on $Op^1$ (the one that induces the Atiyah-Singer topology on $\PseudOpCl^1(\Spinor_{\gamma_0})$).
 Since $B - \preGauge$ is compactly supported, Proposition \ref{App: Composition PSiDOs - Prop} and a partition of unity argument shows that $P \mapsto BPB^\ast$ is continuous with respect to the amplitude topology.
 But $(\placeholder)^\extension$, which is denoted with $\placeholder \tilde{\boxtimes} \id$ in the appendix, is a continuous linear map, see Theorem \ref{App: ExternalTensor is Continuous - Theorem}, so we conclude
 \begin{align*}
  \lim_{n \to \infty}\bigl(BP_nB^\ast\bigr)^\extension &= B^\extension \lim_{n\to\infty} P_n^\extension (B^\ast)^\extension\\
  &= \bigl(BQ_0B^\ast\bigr)^\extension \in \PseudOpCl^1(M \times \R^{n}\times \R,\Spinor_{g_{\susp(BQ_0B^\ast)}}).
 \end{align*}
 
 For $\tilde{\Theta}$ it is even easier to verify because it is the sum of a block Dirac operator (interpreted as constant curve) and a smooth curve $W$ of infinitely smoothing operators.
 We conclude
 \begin{equation*}
  \tilde{\Theta}^\extension = \bigl(B_{1/N}Q_0B_{1/N}^\ast\bigr)^\extension  + W^\extension \in \PseudOpCl^1(M\times \R^n \times \R;\Spinor_{\susp(\gamma_{1/N})}) 
 \end{equation*}  
 from the Theorem \ref{ASTensor - Theorem} and Theorem \ref{App: ExternalTensor is Continuous - Theorem}. 
\end{proof} 

The next tool allows to glue two matching invertible pseudo Dirac operators together.
To introduce this tool, we need to discuss shift operators on pseudo differential operators.
\begin{definition}
 Let $E \rightarrow M \times \R^n$ be a vector bundle.
 For $v \in \R^n$, define the \emph{shift operator}
 \nomenclature{$\Shift_v$}{Shift operator on sections $v$-direction}
 \begin{align*}
  \Shift_v \colon \Gamma_c(M\times \R^n;E) &\rightarrow \Gamma_c(M \times \R^n; E), \\
  \Shift_v(\sigma)(m,x) &= \sigma(m,x-v).
 \end{align*}
\end{definition}
The shift operators induce linear maps on the spaces of pseudo differential operators via push-forward
\begin{align*}
 \Shift_{v \, \ast} \colon \PseudOp^k(M \times \R^n; E) &\rightarrow \PseudOp^k(M \times \R^n; E), \\
 P &\mapsto \Shift_v \circ P \circ \Shift_{-v}.
\end{align*}
In case that $E$ is the tangent or the spinor bundle, there is a less ad-hoc way to define the shift operator.
The translation $\mathrm{trans}_v \colon (m,x) \mapsto (m,x+v)$ is a diffeomorphism on $M \times \R^n$.
Its differential gives an isometry 
\begin{equation*}
 T \mathrm{trans}_v \colon (T(M \times \R^n),g) \rightarrow (T(M \times \R^n),\mathrm{trans}_{-v}^\ast g).
\end{equation*}
Clearly, $\mathrm{trans}_{-v}^\ast g$ and (the ad-hoc version of) $\Shift_{v}g$ agree.
Since the spinor bundle construction is functorial with respect to Pin structure preserving isometries, we obtain an  isomorphism of Clifford-module bundles $\Spinor_g \rightarrow \Spinor_{\Shift_{v}g}$ over $\mathrm{trans}_v$.
The fibres $(\Spinor_g)_{(m,x-v)}$ and $(\Spinor_{\Shift_v(g)})_{(m,x)}$ are the same, so the push forward with this Clifford bundle isomorphism on the level of sections is again given by
\begin{align*}
  \Shift_v \colon \Gamma_c(M\times \R^n;\Spinor_g) &\rightarrow \Gamma_c(M \times \R^n; \Spinor_{\Shift_v g}), \\
  \Shift_v(\sigma)(m,x) &= \sigma(m,x-v).
 \end{align*}
Since the Sobolev norms, as defined in Appendix \ref{Sobolev Spaces - Chapter}, depend only on the local geometry of the bundles in considerations, all shift operators induce isometries between the corresponding Sobolev spaces
\begin{equation*}
 \Shift_v \colon H^k(\Spinor_g) \rightarrow H^k(\Spinor_{\Shift_{v \, \ast} g}).
\end{equation*}
Together with the observation made above that conjugations with $\Shift_v$ preserves the set of (actual) pseudo differential operators, we conclude that these conjugations induce linear maps on the set of pseudo differential operators that are continuous with respect to the Atiyah-Singer topology. 
Continuity with respect to the Atiyah-Singer topology allows us to extend $\Shift_{v\, \ast}$ continuously to the Atiyah-Singer closure.

We will write $\Shift_{j,R}$ instead of $\Shift_{Re_j}$.

\begin{definition}\label{Operator Addition - Definition}
 Let $P,Q \in \BlockDirac[n]$ be two block Dirac operators whose cores are contained in $M \times \rho I^n$ for some $\rho>0$ and that satisfy $\face[1]{i}P= \face[-1]{i} Q$ for some $1\leq i \leq n$.
 Then we define, for all $R>\rho$, the block metric
 \begin{equation*}
     g_P +_{i,R} g_Q \deff {\Shift_{i,-2R}}_\ast(g_P) \cup {\Shift_{i,2R}}_\ast(g_Q) 
 \end{equation*}
 and the block operator
 \begin{equation*}
     P +_{i,R} Q \deff {\Shift_{i,-2R}}_\ast(P) \cup {\Shift_{i,2R}}_\ast(Q) 
 \end{equation*}
 on $\Spinor_{g_p +_{i,R} g_Q}$, where $\cup$ was defined in Definition \ref{SetUp2}.
\end{definition}
\begin{proposition}\label{GluingBlockOperators - Prop}
 With the notation from the previous definition, the operator $P +_{i,R} Q$ is a block Dirac operator, whose faces are given by
 \begin{equation*}
     \face[\omega]{j}(P +_{i,R} Q) \deff \begin{cases}
      \face[-1]{i} P, & \text{if } j=i, \omega = -1, \\
      \face[1]{i} Q, & \text{if } j=i, \omega = 1, \\
      \face[\omega]{j}P +_{i-1,R} \face[\omega]{j}Q, & \text{if } j < i, \\
      \face[\omega]{j}P +_{i,R} \face[\omega]{j}Q, & \text{if } j > i.
    \end{cases}
 \end{equation*}
 If $P$ and  $Q$ are invertible, then there is an $R_0$, depending on $P$ and $Q$, such that $P +_{i,R} Q$ is also invertible for all $R > R_0$.
\end{proposition}
\begin{proof}
 It is easy to see that $P +_{i,R} Q \in \PseudOpCl^{1}(\Spinor_{g_{P +_{i,R} Q}})$ because the space is local, see \cite{atiyah1968indexI}.
 Furthermore, it is easy to see that it is a block operator with the claimed faces.
 Symbol calculus implies $\symb_1(P +_{i,R} Q) = \Cliffmult(g_{P +_{i,R} Q \, \sharp}(\placeholder))$. 
 Corollary \ref{SymmetryUnion - Cor} implies that $P +_{i,R} Q$ is symmetric.
 
 Since the shift operator is induced by the pushforward of a diffeomorphism, we have ${\Shift_{i,R}}_\ast\Dirac_{g} = \Dirac_{\Shift_{i,R}(g)}$ so that $\Dirac_{g_P +_{i,R}g_Q} = \Dirac_{g_P} +_{i,R} \Dirac_{g_Q}$.
 It follows that
 \begin{equation*}
     P +_{i,R} Q - \Dirac_{g_{ P +_{i,R} Q}} =  P +_{i,R} Q - \Dirac_{g_P} +_{i,R} \Dirac_{g_Q} = (P - \Dirac_{g_P}) +_{i,R} (Q - \Dirac_{g_Q}) 
 \end{equation*}
 is a convex combination of two bounded operators and hence itself bounded.
 
 Assume now that $P$ and $Q$ are invertible, then there are positive constants $c_P$ and $c_Q$ respectively, such that $||P\sigma||_0 \geq c_P ||\sigma||_0$ and $||Q\sigma||_0 \geq c_Q||\sigma||_0$.
 By Corollary \ref{Gluing Lower Bounds - Cor} there is an $R_0$, such that $P +_{i,R} Q$ is bounded from below for all $R \geq R_0$. 
 Since $P+_{i,R}Q$ is symmetric and differs from $\Dirac_{g_{P+_{i,R}Q}}$ by a bounded operator, $P+_{i,R}Q$ is also self-adjoint and therefore invertible.
\end{proof}

Finally, we are able to prove the main theorem of this section.

\begin{proof}[Proof of the Fibration Theorem \ref{SymbolMapKanFib - Theorem}]
 We need to show that every commutative outer square 
 \begin{equation*}
  \xymatrix{\CubeBox \ar[r] \ar[d] & \InvBlockDirac \ar[d] \\ \StandCube{n} \ar@{-->}[ru] \ar[r] & \BlockMetrics.}
 \end{equation*}
 has a dotted lift. In other words, for a given set of elements
 \begin{equation*}
  \{P_{(j,\omega)} \in \InvBlockDirac[n-1] \, : \, \face[\omega]{j}P_{(k,\eta)} = \face[\eta]{k}P_{(j,\omega)} \text{ for } j<k, \, (j,\omega),(k,\eta) \neq (i,\eps)   \}
 \end{equation*}
 and a block metric $g \in \BlockMetrics[n]$ such that $g_{P_{(j,\omega)}} = \face[\omega]{j}g$, we need to find a $P \in \InvBlockDirac[n]$ such that $\face[\omega]{j} P = P_{(j,\omega)}$ and $g_P = g$.
 The problem is completely symmetric in $(i,\eps)$, so we will assume that $(i,\eps) = (n,1)$.
 Since $\InvBlockDirac$ is a Kan set, there is a $\hat{P} \in \InvBlockDirac[n]$ such that $\face[\omega]{j}\hat{P} = P_{(j,\omega)}$.
 However, in general $g_{\hat{P}} \neq g$ and $g_{\face[1]{n}\hat{P}} = \face[1]{n}g_{\hat{P}} \neq \face[1]{n}g$.
 
 We will first modify $\hat{P}$ such that the resulting invertible block Dirac operator $P^{\mathrm{aux}}$ additionally satisfies  $g_{\face[1]{n}P^{\mathrm{aux}}} = \face[1]{n}g$.
 To this end, pick a smooth function $\chi \colon [-1,1] \rightarrow [0,1]$ that is identically $0$ near $-1$ and identically $1$ near $1$.
 The map 
 \begin{equation*}
   \gamma \colon t \mapsto  (1- \chi(t))g_{\face[1]{n} \hat{P}} + \chi(t) \face[1]{n}g 
 \end{equation*} 
 is a smooth curve $\gamma \colon [-1,1] \rightarrow \BlockMetrics[n-1]_{\face[1]{n}\hat{P},R}$ for some sufficiently large $R>0$ because of the assumed compatibility conditions $\face[\omega]{j} g_{\face{i} \hat{P}} = \face[\omega]{j} \face{i}g$ for all $(j,\omega) \in \{1,\dots,n-1\} \times \Z_2$.
 
 The Path Lifting Proposition \ref{PathLiftingProp - Prop} applied to $Q_{-1} := \face[1]{n}\hat{P}$ and $\gamma$ yield a smooth, suspensionable curve of invertible block Dirac operators $Q \colon [-1,1] \rightarrow \InvBlockDirac[n-1]$, constant near the boundary, such that $g_{Q_1} = \face[1]{n}g$.
 Extend it constantly to a smooth map $Q \colon \R \rightarrow \InvBlockDirac[n-1]$.
 By re-parametrising $Q$, if necessary, we may assume that
 \begin{equation*}
   \tilde{Q} \deff \susp(Q) = Q^\extension  + \evenodd\boxtimes \partial_n \cdot \frac{\partial}{\partial x_n}
 \end{equation*}
 is an invertible operator on $L^2(M \times \R^n; \Spinor_{\susp(\gamma)})$.
 By Proposition \ref{GluingBlockOperators - Prop}, there is a sufficiently large $R > 0$ such that the block Dirac operator 
 \begin{equation*}
  P^{\mathrm{aux}} \deff {\Shift_{n,2R}}_\ast(\hat{P} +_{n,R} \tilde{Q})
 \end{equation*}  
 is invertible.
 Since the faces of $Q$ are independent of the curve parameter, we have $\face[\omega]{j} \tilde{Q} = \face[\omega]{j} \degen{n} \face[1]{n} \hat{P}$ for all $(j,\omega) \neq (n,1)$. 
 We conclude that $\face[\omega]{j}\left(P^{\mathrm{aux}} \right) = \face[\omega]{j}P = P_{(j,\omega)}$.
 By construction, $\face[1]{n}\left(P^\mathrm{aux} \right) = \face[1]{n}\tilde{Q} = Q_1$.
 Thus, $P^{\mathrm{aux}}$ is an invertible block Dirac operator that satisfies $\face[\omega]{j}g_{P^{\mathrm{aux}}} = \face[\omega]{j}g$ for all $(j,\omega) \in \{1,\dots, n\} \times \Z_2$.
 
 Proceeding as before, the smooth curve
 \begin{equation*}
  \Gamma \colon [-1,1] \rightarrow \BlockMetrics[n]_{g_{P^{\mathrm{aux}}}}, \qquad t \mapsto  (1- \chi(t))g_{P^{\mathrm{aux}}} + \chi(t) g
 \end{equation*}
 satisfies $\face[\omega]{j} \Gamma_t = \face[\omega]{j} \Gamma_{-1}$ for all $t \in [-1,1]$ and all face maps.
 Apply the path-lifting lemma to $P^{\mathrm{aux}}$ and $\Gamma$ to find a homotopy lift $\mathcal{P}$ of $\Gamma$. 
 Since the homotopy between $\Gamma$ and $g_{\mathcal{P}}$ is relative to the boundary, the operator $P_{\mathrm{fill}} \deff \mathcal{P}_1$ is the desired filler. 
\end{proof}

\section{Operator Theoretic Addition}\label{Operator Theoretic Addition - Section}

This section is the operator-theoretic analog to Section \ref{Geometric Addition - Section}.
We will prove a version of the ``angle rotation'' for operators.
Once this result is established, we deduce the operator-theoretic analogs of the results of Section \ref{Geometric Addition - Section} from the pictures therein.
The construction also produces an underlying block metric $g^\angle_R$ that can be thought, up to a shift, as the elongation of the resulting metric of the construction described in Proposition \ref{Construction-AngleRotation - Prop}.
For the construction, recall the shift operator in the context of spinor bundles from Section \ref{SymbolMapKanFib - Section}.

\begin{proposition}
  Let $n \geq 1$ and let $P$ be a block Dirac operator with underlying block metric $g$ on $M \times \R^n$ whose core is contained in $M \times \rho I^n$. 
  Then, for all $R \geq 1$, there is a block Dirac operator $P_R^\angle$ 
  \nomenclature{$P_R^\angle$}{operator that is obtained from $P$ by ``angle-rotation'' with rotation speed $1/R$}
  on $M \times \R^{n+1}$ with underlying block metric $g^\angle_R$ satisfying the following properties:
  \begin{itemize}
      \item[(i)] The faces of $P^\angle_R$ are given by 
      \begin{equation*}
          \face{i}(P_R^\angle) = \begin{cases}
            \degen{n+1}\face{i}P, & \text{if } i<n, \\
            \degen{n+1}\face{i}P, & \text{if } i \geq n, \eps = 1, \\
            P, & \text{if } i \geq n, \eps = -1. 
          \end{cases}
      \end{equation*}
      \item[(ii)] If $P$ is invertible, so is $P^\angle_R$ for all sufficiently large $R$.
  \end{itemize}
\end{proposition}

\begin{proof}
 As in the proof of Proposition \ref{Construction-AngleRotation - Prop}, let $K$ be a compact, convex, point-symmetric body inside $I^2$ with smooth boundary $\partial K$ that further agrees with $I^2$ inside $\{|x_1| \leq 1/2\}$ and $\{|x_2| \leq 1/2\}$.
 Let $\gamma \colon [0,l] \rightarrow \partial K \cap \R_{\geq 0}^2$ be the smooth curve parametrised by arclength that satisfies $\gamma(0) = e_1$, let $v \colon [0,l] \rightarrow \R^2$ be the normalised vector field that is perpendicular to $\gamma'$ and satisfies $\det(v,\gamma') > 0$, and let $\kappa$ be the curvature of $\gamma$.
 Recall that this sign convention implies $v(0) = e_1$ and that $v' = \kappa \gamma'$.
 Denote by $\gamma_R$, $v_R$, and $\kappa_R$ the corresponding objects obtained by replacing $K$ by $R\cdot K$.
 
 Since $K \cap \{|x_j| \leq 1/2\}$ is a rectangle, we can form the elongation
 \begin{align*}
     \mathbf{K} \deff K \cap \R_{\geq 0}^2 \cup \{(x_1,x_2) \in \R^2 \setminus \R_{>0}^2 \, : x_1 \in K \text{ or } x_2 \in K\} \cup \R_{\leq 0}^2
 \end{align*}
 and extend $\gamma$, $v$, $\kappa$ in a constant manner to smooth maps $\R \rightarrow \R^2$. 
 The same holds true for $R\cdot K$ and the rescaled analogs $\gamma_R$, $v_R$, and $\kappa_R$.
 
 The extended objects introduce coordinates on $\R^2 \setminus R\mathbf{K}$ via
 \begin{align*}
     \Xi_R \colon (0,\infty) \times \R &\rightarrow \R^2 \setminus R\mathbf{K} \\
     (r,\varphi) &\mapsto \gamma_R(\varphi) + r\cdot v_R(\varphi),
 \end{align*}
such that
 \begin{equation*}
     \left(\Xi_R^\ast\euclmetric\right)_{(r,\varphi)} = \diff r^2 + (1 + r\kappa_R(\varphi))^2 \diff \varphi^2.
 \end{equation*}
 
 Now that we have recalled the required notation, we can start with the construction.
 In order to avoid case distinctions, we assume that the operator $P$ is invertible to begin with.
 This is no restriction because each step in the construction can be carried out without the invertibility assumption (of course, the resulting operator then has no chance to be invertible).
 
 We fix a lower bound $c_P$ for $P$ that we assume to be smaller than $1$.
 We abbreviate $M \times \R^{n-1}$ to $N$ so that $M \times \R^n = N \times \R$ (the last coordinate is the distinguished one).
 For a fixed $\rho' > \max\{32/c_P,\rho\}$, the operator $\Shift_{n,2\rho}(P)$ restricts to $N \times \R_{>0}$ and we denote the restriction with $Q$.
 
 Define on $N \times \R_{>0} \times \R$ (in $\Xi_R$-coordinates) the metric 
 \begin{equation*}
     (g_R)_{(n,r,\varphi)} \deff \left(\Shift_{r,2\rho'}(g)\right)_{(n,r)} + (1+ r\kappa_R(\varphi))^2\diff \varphi^2 
 \end{equation*}
 and on $\Spinor_{g_R} \rightarrow N\times \R_{>0} \times \R$ the operator
 \begin{equation*}
     \mathtt{P}_R \deff \fullGauge(\Shift_{r,2\rho'}(g)\oplus \diff \varphi^2,g_R)_\ast (Q\boxtimes \id - \Dirac_{\Shift_{r,2\rho'}g} \boxtimes \id) + \Dirac_{g_R}.
 \end{equation*}
 
 Note that $(\Xi_R)_\ast(g_R) = \left(\Shift_{n,2\rho'+R}\right)(g) \oplus \diff x_{n+1}^2$ on $N \times \Xi_R((0,\rho')\times \R)$, so that the metrics $(\Xi_R)_\ast(g_R)$ and $\face[-1]{n}g \oplus \euclmetric_{\R^2}$ agree there.
 Since $\rho' > \rho$, the operators $Q$ and $\Dirac_{\Shift_{n,2\rho' \ast}g}$ decompose on $N \times (0,\rho')$ into
 \begin{align*}
     Q = \face[-1]{n}P \boxtimes \id + \evenodd \boxtimes \DiracSt_{(0,\rho')} \quad \text{ and } \quad
     \Dirac_{\Shift_{n,2\rho' \ast}g} = \Dirac_{\face[-1]{n}g} \boxtimes \id + \evenodd \boxtimes \DiracSt_{(0,\rho')},
 \end{align*}
 where $\DiracSt_{(0,\rho')}$ denotes the Dirac operator of the Euclidean metric on $(0,\rho')$.
 Thus, the operators $(\Xi_R)_\ast(\mathtt{P}_R)$ and $\face[-1]{n}P \boxtimes \id + \evenodd \boxtimes \DiracSt_{\R^2}$ restrict to $N \times \Xi_R((0,\rho') \times \R)$ and agree there.
 They also act as derivations in $r$-direction, so we can glue these operators together\footnote{That means that we are in setup of Definition \ref{SetUp2} although we refrain from spelling out the constants $a$ and $b$ explicitly. For the $R$ in Definition \ref{SetUp2}, we choose $R=4/c_P$ so that $2R < \rho'/2$. For completeness, the remaining data is the following: The ambient manifold $M$ of Definition \ref{SetUp2} is here $M \times \R^{n+1}$, the separating hypersurface is $N \times \Xi_R({\rho'/2}\times \R)$, the operator $P_1$ is $\face[-1]{n}P \boxtimes \id + \evenodd \boxtimes \DiracSt_{\R^2}$ on $N \times \bigl(RK \cup \Xi_R((0,\rho') \times \R)\bigr)$, and $P_2$ is $(\Xi_R)_\ast(\mathtt{P}_R)$}.
 We define 
 \begin{equation*}
     P_R^\angle \deff \Shift_{-(R+2\rho')(e_n + e_{n+1})\, \ast}\bigl( \face[-1]{n}P\boxtimes \id + \evenodd \boxtimes \DiracSt_{\R^2} \cup \Xi_{R,\ast}(\mathtt{P}_R)\bigr)
 \end{equation*}
 with underlying metric $g_R^\angle \deff \Shift_{-(R+2\rho')(e_n + e_{n+1})\, \ast}(g_R \cup \face[-1]{n}g \oplus \euclmetric_{\R^2})$ on $M \times \R^{n+1}$.
 
 Note that $P_R^\angle \in \PseudOpCl^1(\Spinor_{g_R^\angle})$ because the space $\PseudOpCl^1(\Spinor_{g_R^\angle})$ is local and $P_R^\angle$ is obtained from $P$ by adding other pseudo differential operators (in the Atiyah-Singer sense), tensoring it with the identity, and pushing operators forward with vector bundle maps and diffeomorphisms.
 
 The operator $P_R^\angle$ is also a block operator of Dirac type.
 Indeed, its core is contained inside $M \times (R + 5\rho')I^n$ as one can check for sections that are either supported within the union of $N \times R\cdot \mathbf{K}$ and the gluing-strip $N \times \Xi_R((0,\rho') \times \R)$, or supported outside of $N \times R \cdot \mathbf{K}$ separately.
 
 The operator $P_R^\angle$ is symmetric because $\mathtt{P}_R$ acts on $N \times (0,\rho') \times \R$ as a differential operator in $r$-direction, so we can apply Corollary \ref{SymmetryUnion - Cor}.
 
 The operator $(P - \Dirac_g) \boxtimes \id$ induces a bounded operator on $L^2(\Spinor_{g \oplus \diff \varphi^2})$, which also implies that its principal symbol vanishes by Lemma \ref{ASExactSequence - Lemma}.
 Hence, $\mathtt{P}_R$ differs from the Dirac operator of its underlying metric $g_R$ by a bounded operator and has the same  principal symbol as the Dirac operator of its underlying metric.
 The same is true for $\face[-1]{n}P \boxtimes \id + \evenodd \boxtimes \DiracSt_{\R^2}$, so its also true for the union $(\Xi_R)_\ast\mathtt{P}_R \cup (\face[-1]{n}P \boxtimes \id + \evenodd \boxtimes \DiracSt_{\R^2})$.
 
 Visual reasons, see figure \ref{fig:OperatorConstruction} 
 and the fact that we modify the operator $P\boxtimes \id$ only in the last two coordinates imply that $P^\angle_R$ has the faces as claimed. 
 This proves (i).
 
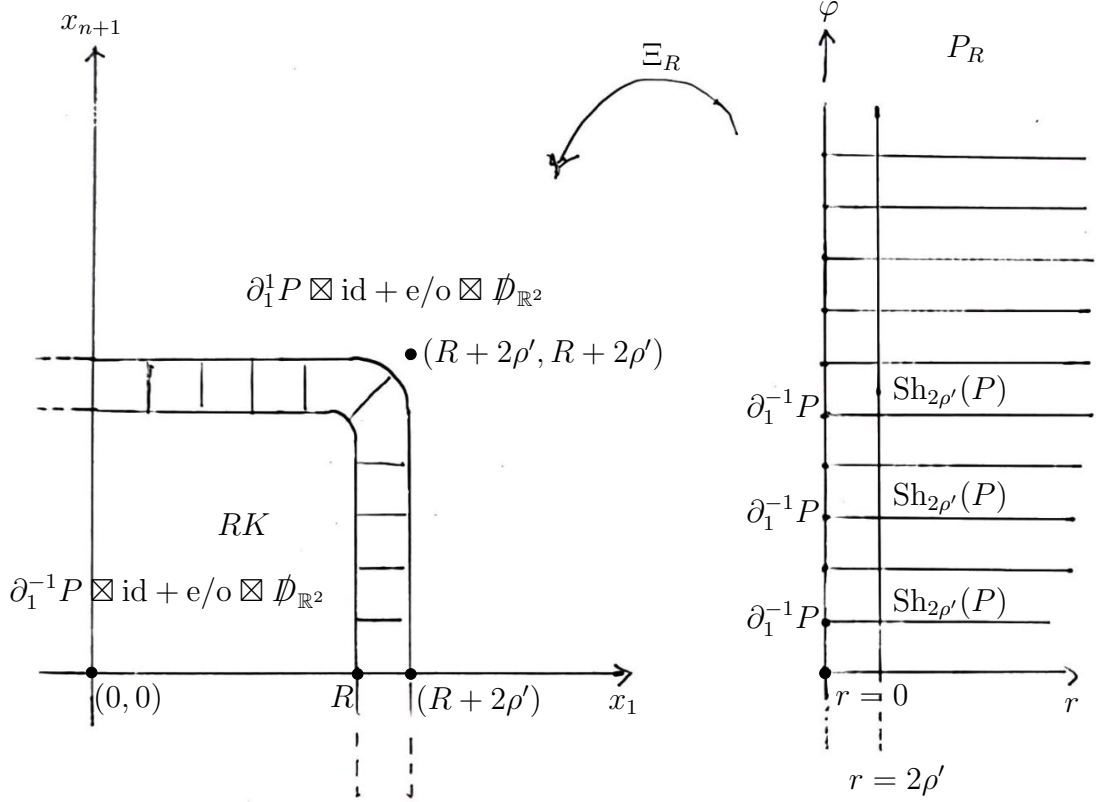
\begin{figure}
    \centering
    \begin{tikzpicture}
		\node at (0,0) {\includegraphics[width=\textwidth]{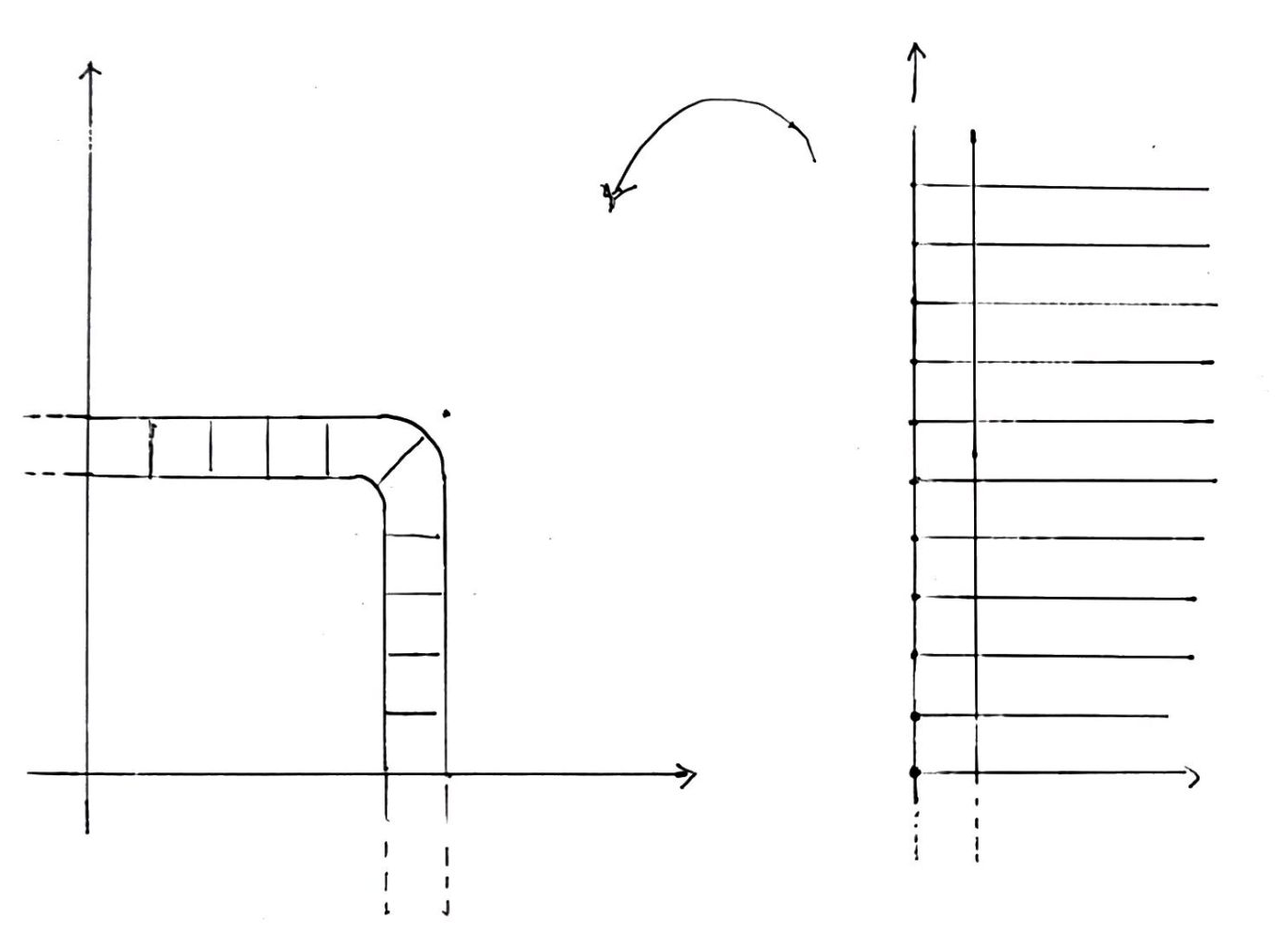}};
		\node at (-6.5,5) {$x_{n+1}$};
		\node at (.5,-4) {$x_1$};
		\node at (6.4,-4) {$r$};
		\node at (3.2,5.2) {$\varphi$};
		\node at (1,4.6) {$\Xi_R$};
		\node at (-4.5,-1.6) {$RK$};
		\node at (-5.5,-2.5) {$\face[-1]{1}P \boxtimes \id + \evenodd \boxtimes \DiracSt_{\R^2}$};
		\node at (-2.5,1.5) {$\face[1]{1}P \boxtimes \id + \evenodd \boxtimes \DiracSt_{\R^2}$};
		\node at (5,4.7) {$P_R$};
		\node at (4.1, -5) {$r=2\rho'$};
		\node at (2.6,-2.85) {$\face[-1]{1}P$}; \node at (4.8,-2.65) {$\Shift_{2\rho'}(P)$};
		\node at (2.6,-1.45) {$\face[-1]{1}P$}; \node at (4.8,-1.25) {$\Shift_{2\rho'}(P)$};
		\node at (2.6,-0.05) {$\face[-1]{1}P$}; \node at (4.8,0.15) {$\Shift_{2\rho'}(P)$};
		\filldraw[black] (-6.5,-3.57) circle (2pt) node[anchor=north]{$\ \qquad (0,0)$};
		\filldraw[black] (-3,-3.59) circle (2pt) node[anchor=north]{$R \quad$};
		\filldraw[black] (-2.3,-3.59) circle (2pt) node[anchor=north]{$\qquad \qquad \ (R+2\rho')$};
		\filldraw[black] (3.15,-3.58) circle (2pt) node[anchor=north]{$\qquad \quad r=0$};
		\filldraw[black] (-2.3,0.65) circle (2pt) node[anchor=west]{$(R+2\rho',R+2\rho')$};
	\end{tikzpicture}
    \caption{The construction of $\bigl( \face[-1]{n}P\boxtimes \id + \evenodd \boxtimes \DiracSt_{\R^2} \cup \Xi_{R,\ast}(\mathtt{P}_R)\bigr)$ and $\mathtt{P}_R$. We need to shift the operator towards $-2\rho(e_n + e_{n+1})$ so that the faces agree with the given operator.
    The right strip shows $\mathtt{P}_R$ in $\Xi_R$-coordinates although we did not draw the Dirac operator summand.}
    \label{fig:OperatorConstruction}
\end{figure}
 
 We will prove in Corollary \ref{AuxilaryOperatorHasLowerBound - Cor} below that $||\mathtt{P}_R\sigma||_0 \geq c_P/3||\sigma||_0$ if $R$ is sufficiently large.
 It follows that $(\Xi_R)_\ast(\mathtt{P}_R)$ is bounded from below on $\Spinor_{(\Xi_R)_\ast(g_R)} \rightarrow N \times \R^2 \setminus R\cdot \mathbf{K}$ with lower bound $c_P/3$.
 Since $\rho'> 32/c_P$, Corollary \ref{Gluing Lower Bounds - Cor} (with gluing strip length $4/c_{P}$) applied to $(\Xi_R)_\ast(\mathtt{P}_R)$ and $\face[-1]{n}P \boxtimes \id + \evenodd \boxtimes \DiracSt_{\R^2}$ implies that $P^\angle_R$ is bounded from below by
 \begin{align*}
     c_{P_R^\angle}^2 &= 1/4\cdot \min\{c^2_{(\Xi_R)_\ast(\mathtt{P}_R)} - c^2_P/16, c_{\face[-1]{1}P} - c^2_P/16\} \\
                      &=1/4 \cdot(c^2_{(\Xi_R)_\ast(\mathtt{P}_R)} - c_P^2/16) \\
                      &\geq 1/4(1/9 - 1/16) c_P^2 > c_P^2/100.
 \end{align*}
 As $P_R^\angle$ is a block Dirac operator, in particular self-adjoint, it is invertible, which proves (ii).
\end{proof}

\begin{lemma}
 Using the notation from the proof of the previous proposition, the operator $D \deff \fullGauge(g_R,\Shift_{r,2\rho'}(g) \oplus \diff  \varphi^2)_\ast(\Dirac_{g_R}) - \Dirac_{\Shift_{r,2\rho'}(g)} \boxtimes \id$ is given by
 \begin{equation*}
     D = (1+r\kappa_R)^{-1}\partial_\varphi \cdot \frac{\partial}{\partial \varphi} -\frac{1}{2}(1+r\kappa_R)^{-2}r(\partial_\varphi \kappa_R) \partial_\varphi \cdot  .
 \end{equation*}
\end{lemma}
\begin{proof}
 Abbreviate $\fullGauge(g_R,\Shift_{r,2\rho'}(g) \oplus \diff\varphi^2)$ to $\fullGauge$.
 We need to express $D$ in the coordinates induced by $\Xi_R$.
 To this end, we fix a local $(\Shift_{r,2\rho'}(g)\oplus \diff \varphi^2)$-orthonormal frame $e_1,\dots,e_n,e_{n+1}$ such that $e_1,\dots, e_n$ are independent of $\varphi$ and $e_{n+1} = \partial_\varphi$.
 
 Recall the pre-gauge map $\tau_{g_R, \Shift_{r,2\rho'}(g)\oplus \diff \varphi^2}$ from Definition \ref{TangentPreGauge - Definition} and its induced map $\Phi = \Phi_{g_R,\Shift_{r,2\rho'}(g) \oplus \diff \varphi^2}$ between the corresponding Clifford algebras. 
 Define the $g_R$-orthonormal frame 
 \begin{equation*}
     f_j \deff \tau_{\Shift_{r,2\rho'}(g) \oplus \diff \varphi^2,g_R}(e_j) = \begin{cases}
     e_j, & \text{if } j < n+1, \\
     (1+r\kappa_R)^{-1}\partial_\varphi, & \text{if } j = n+1.
    \end{cases}
 \end{equation*}
 Recall from the proof of Lemma \ref{DiracOperatorSmooth - Theorem} the formula
 \begin{equation*}
     \fullGauge_\ast(\Dirac_{g_R}) = \preGauge_\ast(\Dirac_{g_R}) -\alpha^{-1} \Phi(\mathrm{grad}_{g_R}\alpha) \cdot\ ,
 \end{equation*}
 where 
 \begin{equation*}
    \alpha^2 = \alpha^2_{g_R,\Shift_{r,2\rho'}(g)\oplus \diff \varphi^2} = \frac{\diff \vol_{g_R}}{\diff \vol_{\Shift_{r,2\rho'}(g)\oplus \diff \varphi^2}} = (1+r\kappa_R).
 \end{equation*}
 Differentiating $\alpha$ yields
 \begin{equation*}
     \diff \alpha = \frac{1}{2}(1+r\kappa_R)^{-1/2}\left(\kappa_R \diff r + r \partial_\varphi \kappa_R \diff \varphi\right).
 \end{equation*}
 hence, the gradient is given by
 \begin{align*}
     \mathrm{grad}_{g_R}(\alpha) &= \mathrm{grad}_{\Shift_{r,2\rho'}(g) + (1+r\kappa_R)^2 \diff \varphi^2}(\alpha) \\
     &= \frac{1}{2}(1+r\kappa_R)^{-1/2}\left(\kappa_R\mathrm{grad}_{\Shift_{r,2\rho'}(g)}(r) + r (\partial_\varphi \kappa_R) \cdot (1+r\kappa_R)^{-2} \partial_\varphi \right)
 \end{align*}
 so that the first summand is given by
 \begin{align*}
     \alpha^{-1}\Phi(\mathrm{grad}_{g_R}\alpha) = \frac{1}{2}(1+r\kappa_R)^{-1}\left(\kappa_R \mathrm{grad}_{\Shift_{r,2\rho'}(g)}(r) + r (\partial_\varphi \kappa_R) \cdot (1+r\kappa_R)^{-1} \partial_\varphi\right).
 \end{align*}
 
 The determination of the second summand requires calculations of Christoffel symbols.
 From \cite{LawsonMichelsonSpin}*{Theorem II.4.14} follows
 \begin{align*}
     \preGauge_\ast(\Dirac_{g_R}) &= \preGauge_\ast\left(\sum_{i=1}^{n+1} f_i \cdot \nabla_{f_i}^{\Spinor_{g_R}}\right) \\
     &= \sum_{i=1}^{n+1} e_i \cdot \preGauge_\ast\left(\nabla^{\Spinor_{g_R}}\right)_{f_i} \\
     &= \sum_{i=1}^{n+1} e_i \cdot \left(\partial_{f_i} + \frac{1}{2}\sum_{j<k} {}_f^{g_R}\Gamma_{ij}^k e_je_k\right),
 \end{align*}
 where ${}_f^{g_R}\Gamma_{ij}^k$ are the Christoffel symbols for the Levi-Cevita connection of $g_R$ expressed in the local orthonormal frame $f_1,\dots,f_{n+1}$, and $\partial_{f_j}$ is the usual partial differential in $f_j$-direction under the identification 
 \begin{align*}
     \Gamma_c(M\times \R^n,\Spinor_{\Shift_{r,2\rho'}(g) \oplus \euclmetric_{\R^n}}) &= \Gamma_c(M\times \R^n,\Spinor_{\Shift_{r,2\rho'}(g)} \boxtimes \Cliff[n][0]) \\
     &\cong \mathcal{C}_c^\infty(\R^n,\Gamma(M,\Spinor_{\Shift_{r,2\rho'}(g)})\otimes \Cliff[n][0]).
 \end{align*}
 
 We also denote $f_{n+1}$ by $f_\varphi$ to emphasise the distinguished coordinate.
 The vector fields $f_1,\dots,f_n$ are independent of $\varphi$, so the commutator $[f_j,f_\varphi]$ is given by
 \begin{align*}
     [f_j,f_\varphi] = [f_j, (1+r\kappa_R)^{-1}\partial_\varphi] &= - (1+r\kappa_R)^{-1}\partial_{f_j}(r) \cdot \kappa_R \partial_\varphi \\
     &= -(1+r\kappa_R)^{-1}\partial_{f_j}(r) \kappa_R f_\varphi
 \end{align*}
 The Koszul formula for the Christoffel symbols implies
 \begin{equation*}
     {}^{g_R}{{}_f\Gamma_{ij}^k} = {}^{\Shift_{r,2\rho'}(g)}{{}_e\Gamma_{ij}^k}, \quad {}^{g_R}{{}_f\Gamma_{ij}^\varphi} = {}^{g_R}{{}_f\Gamma_{\varphi j}^k} = 0, \quad \text{and} \quad {}^{g_R}{{}_f\Gamma_{\varphi j}^\varphi} = (1+r\kappa_R)^{-1}\kappa_R \partial_{f_j}(r).
 \end{equation*}
 Thus, the second summand is given by
 \begin{align*}
     \preGauge_\ast(\Dirac_{g_R}) - \Dirac_{\Shift_{r,2\rho'}(g)} \boxtimes \id &= \preGauge_\ast\left(\sum_{j=1}^{n+1}f_j \cdot \nabla_{f_j}^{\Spinor_{g_R}}\right) - \sum_{j=1}^n e_j \cdot \nabla_{e_j}^{\Spinor_{\Shift_{r,2\rho'}(g)}} \\
     &= \partial_\varphi \cdot \left( \partial_{f_\varphi} + \frac{1}{2}\sum_{j=1}^n (1+r\kappa_R)^{-1}\kappa_R \partial_{f_j}(r) e_j\partial_\varphi \cdot \right) \\
     &= (1+r\kappa_R)^{-1}\left(\partial_\varphi \cdot \frac{\partial}{\partial \varphi} + \frac{1}{2}\kappa_R\mathrm{grad}_{\Shift_{r,2\rho'}(g)}(r) \cdot \right).
 \end{align*}
 If we put everything together, we end up with
 \begin{align*}
     D = \fullGauge_\ast(\Dirac_{g_R}) - \Dirac_{\Shift_{r,2\rho'}(g)}\boxtimes \id = (1+r\kappa_R)^{-1}\partial_\varphi \cdot \frac{\partial}{\partial \varphi} -\frac{1}{2}(1+r\kappa_R)^{-2}r(\partial_\varphi \kappa_R) \partial_\varphi \cdot 
 \end{align*}
 as claimed.
\end{proof}

\begin{cor}\label{AuxilaryOperatorHasLowerBound - Cor}
 Abbreviate $\fullGauge(g_R,\Shift_{r,2\rho'}(g) \oplus \diff \varphi)$ to $\fullGauge$.
 For $\Bar{\rho}>5\rho'$, define $N_{\Bar{\rho}} \deff N \times (0,\Bar{\rho}) \times \R$.
 Then there is a constant $C$ only depending on $\Bar{\rho}$ such that, on $N_{\Bar{\rho}}$, we have
 \begin{align*}
     \left|\left| \fullGauge(g_R,\Shift_{r,2\rho'}(g)\oplus \diff \varphi^2)_\ast(\mathtt{P}_R) - Q \boxtimes \id - \evenodd \boxtimes \DiracSt_{\R}\right|\right|_{1,0} \leq C/R.
 \end{align*}
 In particular, $||\mathtt{P}_R \sigma||_0 \geq c_P/12||\sigma||_0$ for all sufficiently large $R$ and all sections supported within $N \times \R_{>0} \times \R$.
\end{cor}
\begin{proof}
 The previous lemma implies that
 \begin{align*}
     &\quad \ \fullGauge_\ast(\mathtt{P}_R) - Q \boxtimes \id - \evenodd \boxtimes \DiracSt_{\R} = D - \partial_\varphi \cdot \frac{\partial}{\partial \varphi} \\
     &= -(1+r\kappa_R)^{-1}r\kappa_R \partial_\varphi \cdot \frac{\partial}{\partial \varphi} -\frac{1}{2}(1+r\kappa_R)^{-2}r(\partial_\varphi \kappa_R) \partial_\varphi \cdot  .
 \end{align*}
 Since $\kappa_R = 1/R\cdot \kappa(\placeholder/R)$, we estimate the operator norm of $D - \partial_\varphi \cdot \partial/\partial \varphi$ on $N_{\Bar{\rho}}$  as follows
 \begin{align*}
     \left|\left|D - \partial_\varphi\cdot \frac{\partial}{\partial\varphi}\right|\right|_{1,0} &\leq  \left|\left|\frac{r\kappa_R}{1+r\kappa_R}\right|\right|_{0,0}\cdot \underbrace{\left|\left|\partial_\varphi \cdot \frac{\partial}{\partial \varphi}\right|\right|_{1,0}}_{=1} 
     + \frac{1}{2}\left|\left|\frac{r(\partial_\varphi \kappa_R)}{(1+r\kappa_R)^{2}} \right|\right|_{0,0}\\
     &\leq \frac{1}{R}\left|\left|\frac{r(\kappa\circ(1/R))}{(1+r\kappa_R)}\right|\right|_{\infty,N_{\Bar{\rho}}} 
     + \frac{1}{R^2} \left|\left|\frac{r(\partial_\varphi\kappa)\circ(1/R)}{(1+r\kappa_R)^{2}}\right|\right|_{\infty,N_{\Bar{\rho}}}\\
     &\leq C(\Bar{\rho})/R.
 \end{align*}
 
 Since $Q \boxtimes \id + \evenodd \boxtimes \DiracSt_{\R}$ is a block Dirac operator, it satisfies the fundamental elliptic estimate.
 This operator is further bounded from below by $c_P$.
 By Lemma \ref{TopInjOpen - Lemma}, we find a constant $R_0$, only depending on $Q$, such that for all $R > R_0$ the operators $\fullGauge_\ast(\mathtt{P}_R)$ and $Q \boxtimes \id + \evenodd \boxtimes \DiracSt_{\R}$ are sufficiently close on $N_{\Bar{\rho}}$ such that $||\fullGauge_\ast(\mathtt{P_R})\sigma||_0 > c_P/3||\sigma||_0$ for all compactly supported sections inside of $N_{\Bar{\rho}}$.
 The gauge map $\fullGauge(g_R,\Shift_{r,2\rho'}(g)\oplus \diff \varphi^2)$ is a bundle map over the identity that induces an isometry between the corresponding Hilbert spaces of square integrable spinors, so the same holds for $\mathtt{P}_R$.
 
 On $N_{>5\rho'} \deff N \times \R_{>5\rho'} \times \R$, the metrics $g_R$ and $\Shift_{r,2\rho'}(g) \oplus \diff \varphi^2$ decompose into
 \begin{equation*}
     g_R |_{N_{>5\rho'}} = \face[1]{n} g \oplus \Xi_R^\ast\euclmetric \quad \text{ and } \quad  \Shift_{r,2\rho'}(g)|_{N_{>5\rho'}} = \face[1]{n} g \oplus \diff r^2 \oplus \diff \varphi^2
 \end{equation*}
 so that $\fullGauge(g_R,\Shift_{r,2\rho'}(g) \oplus \diff \varphi^2) = \id_N \boxtimes \fullGauge(\Xi_R^\ast \euclmetric, \diff r^2 \oplus \diff \varphi^2)$.
 The operator $Q$ decomposes there into 
 \begin{equation*}
     Q|_{N_{>5\rho'}} = \face[1]{n}Q \boxtimes \id + \evenodd \boxtimes \DiracSt_{\R_{>5\rho'}}.
 \end{equation*}
 Tracing through the definitions now yields
 \begin{align*}
     \fullGauge_\ast(\mathtt{P}_R) &= (Q - \Dirac_{\Shift_{r,2\rho'}(g)}) \boxtimes \id_{\R} + \fullGauge_\ast(\Dirac_{g_R}) \\
     &= \face[1]{n}Q \boxtimes \id_{\R_{>5\rho' \times \R}} + \evenodd \boxtimes \fullGauge(\Xi_R^\ast\euclmetric,\diff r^2 \oplus \diff \varphi^2)_\ast\bigl(\Xi_R^\ast(\DiracSt_{\R^2})\bigr),
 \end{align*}
 from which we deduce that $\fullGauge_\ast(\mathtt{P_R})^2 \geq c_{\face[n]{1}P} \id$ is uniformly positive.
 Hence $\fullGauge_\ast(\mathtt{P}_R)$ is bounded from below with lower bound $c_{\face[n]{1}P}$.
 By choosing the difference $\Bar{\rho} - 5\rho'$ sufficiently large, we can apply Corollary \ref{Gluing Lower Bounds - Cor} to deduce that $||\fullGauge_\ast(\mathtt{P}_R)\sigma||_0 \geq c_P/12||\sigma||_0$ for all compactly supported sections $\sigma$ inside $N \times \R_{>0} \times \R$, which proves the addendum.
\end{proof}

\subsubsection*{Operator Theoretic Addition on the Operator Concordance Set}

Since $\InvBlockDirac$ is a Kan set, the (combinatorial) homotopy groups are well-defined and carry a combinatorical group structure.
This group structure is formally defined through the Kan condition and, as we do not keep track of the fillers, the group structure is only defined on homotopy classes.
We have constructed several geometric group structures on $\pi_n(\ConcSet)$ in subsection \ref{Geometric Addition - Section} that almost work on the level of representatives.
The analogous operations for invertible block Dirac operators have already been introduced in Definition \ref{Operator Addition - Definition}.  
Recall that, for two invertible block operators $P$, $Q \in \InvBlockDirac$ with $\face[1]{j}P = \face[-1]{j}Q$ whose cores are contained in $M \times \rho I^n$, the operator 
\begin{equation*}
     P +_{j,R} Q \deff \Shift_{j,R}(P)|_{\{x_j \leq 0\}} \cup \Shift_{j,R}(Q)|_{\{x_j \geq 0\}} 
     \nomenclature{$ P +_{j,R} Q$}{Operator Addition in $j$-th coordinate direction}
 \end{equation*}
 is a block Dirac operator over the block metric
 \begin{align*}
     g_{P+_{j,R}Q} &\deff \Shift_{j,R}(g_P)|_{\{x_j \leq 0\}} \cup \Shift_{j,R}(g_Q)|_{\{x_j \geq 0\}}  \\
     &= g_P +_{j,R} g_Q
 \end{align*}
 for all $R > \rho$.
 If $R$ is sufficiently large, then the resulting block Dirac operator is also invertible.

The next theorem is the operator-theoretic analogue of Theorem \ref{Geometric Addition - Thm}.
\begin{theorem}\label{Operatortheoretic Addition - Thm}
  There are $n$-many group structures on $\pi_n(\InvBlockDirac)$ that are defined as follows:
  If $P, Q$ are representatives of $[P], [Q] \in \pi_n(\InvBlockDirac)$ and if $R$ is sufficiently large, then we set 
  \begin{equation*}
      [P] +_j [Q] \deff [P +_{j,R} Q].
  \end{equation*}
  If $r_j$ denotes the reflection at the hyperplane $\{x_j = 0\}$, then the inverse element of $[P]$ for $+_j$ is given by $[r_j^\ast P]$.
  The structures are Eckmann-Hilton related to each other.
  Furthermore, $+_1$ agrees with the group structure from cubical set theory.
\end{theorem}
\begin{proof}
 The entire proof is completely analogous to the proof of Theorem \ref{Geometric Addition - Thm} except for a small technical nuisance in the proof that $+_{j,R}$ is independent of $R$.
 
 For two operators $P,Q$, let $R_0 < R_1$ be sufficiently large such that $P +_{j,R_k}Q$ is defined.
 We would like to construct a concordance between them.
 
 If $\chi \colon \R \rightarrow [0,1]$ is a smooth function with $\chi|_{\R\leq -1} = 0$ and $\chi|_{\R_{\geq S}} = 1$ with $|\chi'|\leq /1/S$ and $|\chi''|\leq 4/S^2$, then $P +_{j,R_\chi}Q$ is a path over $g_P +_{j,R_\chi} g_Q$ joining $P+_{j,R_0}Q$ and $P +_{j,R_1}Q$.
 
 If $\susp(P +_{j,R_\chi}Q)$ were an element in $\PseudOpCl^1(\Spinor_{\susp(g_P +_{j,R}g_Q)})$, then, provided $S$ is sufficiently large, the operator $\susp(P +_{j,R_\chi}Q)$ would be an invertible block Dirac operator and, hence, serve as a concordance.
 We cannot apply Theorem \ref{FamilyASTensor - Theorem} directly to the curve $P +_{j,R_\chi} Q$, so we have to employ the following ``Swiss army knife trick'':
 Consider the diffeomorphisms $\Shift_{j,\pm2R_\chi}$ on $\R^{n+1}$ given by
 \begin{equation*}
     \Shift_{j,R_\chi} \colon x \mapsto x + \Bigl( (1-\chi(x_{n+1}))R_0 + \chi(x_{n+1})R_1 \Bigr)\cdot e_j.
 \end{equation*}
 Theorem \ref{ASTensor - Theorem} implies that $P \boxtimes \id \in \PseudOpCl^1(\Spinor_{g_P \oplus \diff x_{n+1}^2})$ and similarly $Q \boxtimes \id \in \PseudOpCl^1(\Spinor_{g_Q \oplus \diff x_{n+1}^2})$.
 It follows that $\Shift_{j,-R_\chi}(P \boxtimes \id) \in \PseudOpCl^1(\Spinor_{\Shift_{j,R_\chi}(g) + \diff x_{n+1}^2})$.
 A similar result holds for $\Shift_{j,R_\chi}(Q \boxtimes \id)$.
 Since $\PseudOpCl^1$ is local, we conclude that
 \begin{equation*}
    \Shift_{j,-R_\chi}(P \boxtimes \id) \cup \Shift_{j,R_\chi}(Q \boxtimes \id) \in \PseudOpCl^1(\Spinor_{\susp(g_P +_{j,R_\chi} g_Q)})
 \end{equation*}
 and hence $\susp(P+_{j,R_\chi}Q) \in  \PseudOpCl^1(\Spinor_{\susp(g_P +_{j,R_\chi} g_Q)})$ because the difference between the two operators is a differential operator.
\end{proof}

\subsubsection*{Index Additivity for Block Dirac Operators}

As an application of Proposition \ref{GluingBlockOperators - Prop}, we prove an additivity formula for the Fredholm index of a block Dirac operator with invertible faces.
There are several theorems of this kind in the existing literature, see \cite{booss1993elliptic}, \cite{bunke1995Callias}, \cite{botvinnik2017infinite}.
The new features of our version is that it works for pseudo differential operators and, more importantly, for separating hypersurfaces that are non-compact.

Although Bunke \cite{bunke1995Callias} was the first who proved such an additivity theorem, we follow in notation and strategy \cite{botvinnik2017infinite}.

For $\alpha \in \{0,\dots, 3\}$, let $P_\alpha \in \BlockDirac[n]$ be block Dirac operators with underlying block metrics $g_\alpha$. 
For a fixed but arbitrary $i \in \{1,\dots,n\}$, assume that $\face[1]{i}P_\alpha = \face[-1]{i}P_\beta$ for $(\alpha,\beta) \in \{(0,1),(2,3),(0,3),(2,1)\}$ (this implies that the corresponding faces of the underlying metrics agree) and that all other faces are invertible.
Let $R_0$ be a positive number such that the cores of all $P_j$ lie in the interior of $M \times R_0I^n$.
For all $R>R_0$ define $P_{\alpha\beta,R} \deff P_\alpha +_{R,i} P_\beta$.
These are block Dirac operators over $\Spinor_{g_{\alpha \beta,R}}$, where $g_{\alpha \beta,R} = g_\alpha +_{R,i} g_\beta$.

\begin{theorem}[Index Additivity]
 There is an $R_1 > R_0$ such that, for all $R > R_1$, the block Dirac operators $P_{\alpha\beta,R}$ have invertible faces and hence are unbounded Fredholm operators.
 Their indices are independent of $R$, and, moreover, they satisfy the following additivity formula:
 \begin{equation}\label{eq: Index Additivity}
     \ind(P_{01}) + \ind(P_{23}) = \ind(P_{03}) + \ind(P_{21}) \in KO^{-(d+n)}(\mathrm{pt}).
 \end{equation}
\end{theorem}
\begin{proof}
 Since the face operators commute with the gluing construction, we find, by Corollary \ref{Gluing Lower Bounds - Cor}, a sufficiently large $R_1$ such that $P_{\alpha\beta,R}$ has invertible faces for all $R > R_1$.
 Proposition \ref{EllipticEstimateInvInf - Prop} implies that $P_{\alpha\beta,R}$ is a Fredholm operator.
 
 The assignment 
 \begin{equation*}
     s \mapsto \fullGauge(g_{\alpha\beta,(1-s)R_1+s R_2}, g_{\alpha\beta,R_1})_\ast \circ  P_{\alpha\beta,(1-s)R_1+s R_2} \circ \fullGauge(g_{\alpha\beta,R_1},g_{\alpha\beta,(1-s)R_1+s R_2})_\ast
 \end{equation*}
 is a path of unbounded Fredholm operators between $P_{\alpha\beta,R_1}$ and an operator that is unitarily equivalent to $P_{\alpha\beta,R_2}$.
 Since the index is invariant under homotopy and unitary equivalence, we have $\ind(P_{\alpha\beta,R_1}) = \ind(P_{\alpha\beta,R_2})$.
 
 In the following, we will drop the suffix $R$ from the Hilbert spaces and operators if it does not lead to confusion. 
 Set $H_{\alpha\beta} \deff L^2(\Spinor_{g_{\alpha\beta}})$.
 Let $H_{\alpha \beta}^\op$ be the \emph{opposite} Clifford module with Clifford action $-\CliffRight$, the negative of the canonical right action, and equipped with the opposite grading $-\evenodd$.
 The operator $P_{\alpha\beta}^\op$ is the ``same'' as $P_{\alpha\beta}$, but acting on $H_{\alpha\beta}^\op$.
 
 By \cite{ebert2017indexdiff}*{Lemma 2.20} the claimed formula is equivalent to
 \begin{equation}\label{eq: AdditvityVanishing}
     \ind(P_{01}) + \ind(P_{23}) + \ind(P_{03}^\op) + \ind(P_{21}^\op) = 0.
 \end{equation}
 The left hand side can be written as the index of a single operator $P$.
 Indeed, let $H = H_{01} \oplus H_{23} \oplus H_{03} \oplus H_{21}$ with Clifford action $r$ and grading $\iota$ be given by
 \begin{equation*}
     r = \begin{bmatrix}
      \CliffRight_{01} & & & \\
      & \CliffRight_{23} & & \\
      & & -\CliffRight_{03} & \\
      & & & -\CliffRight_{21}
     \end{bmatrix}_, \qquad
     \iota = \begin{bmatrix}
      \evenodd_{01} & & & \\
      & \evenodd_{23} & & \\
      & & -\evenodd_{03} & \\
      & & & -\evenodd_{21}
     \end{bmatrix}_.
 \end{equation*}
 The operator $P$ is then given by
 \begin{equation*}
     P = \begin{bmatrix}
      P_{01} & & & \\
      & P_{23} & & \\
      & & P_{03} & \\
      & & & P_{21}
     \end{bmatrix}_.
 \end{equation*}
 
 For later purpose, we merge the Sobolev spaces $H^1 \deff H^1_{01} \oplus H_{23}^1 \oplus H^1_{03} \oplus H^1_{21} \subseteq H$.
 
 Pick $\lambda_0,\mu_0 \colon \R \rightarrow \R$ with $\supp \lambda_0 \subseteq \R_{\geq -1}$, $\supp \mu_0 \subseteq \R_{\leq 1}$ and $\lambda_0^2 + \mu_0^2 = 1$.
 We may assume that $|\lambda_0'|,|\mu_0'| \leq 1$.
 Define $\lambda_R,\mu_R \colon \R^n \rightarrow \R$ by $\lambda_R(x) = \lambda_0(R^{-1}x_i)$ and $\mu_R(x) = \mu_0(R^{-1}x_i)$.
 The operator 
 \begin{equation*}
     J_0 \deff \begin{bmatrix}
      &&-\mu_R & -\lambda_R \\
      &&-\lambda_R &\mu_R \\
      \mu_R & \lambda_R && \\
      \lambda_R & -\mu_R &&
     \end{bmatrix}
 \end{equation*}
 is a well-defined operator on $H = H_R$ because, for example, for the upper right corner, the restrictions agree
 \begin{equation*}
     \Spinor_{g_{21}}|_{\supp \, \lambda_R} =  \Spinor_{g_1}|_{\supp \lambda_R} = \Spinor_{g_{01}}|_{\supp \lambda_R}
 \end{equation*}
 and we can use these identities to transplant section from  $H_{21}$ to $H_{01}$.
 Similar results hold for the other entries.
 
 The operator $J \deff J_0 \iota$ is then a self-adjoint, Clifford-linear, odd (with respect to $\iota$), and bounded involution in $\mathrm{Fred}^{d+n,0}(H)$. 
 Note that, for algebraic reasons, $J$ cannot be essentially positive or essentially negative.
 To verify equation (\ref{eq: AdditvityVanishing}), it is enough to find a path in $\mathrm{Hom}(H^1,H)$ between $P$ and $J$ that take values in unbounded, self-adjoint, odd, Clifford-linear Fredholm operators.
 We claim that such a path is given by 
 \begin{equation*}
     P_s \colon [0,\pi/2] \rightarrow \mathrm{Hom}(H^1,H),\quad s\mapsto \cos(s)P + \sin(s)J.
 \end{equation*}
 Clearly, each $P_s$ is self-adjoint, odd, Clifford-linear.
 It remains to show that $P_s$ is a Fredholm operator, provided $R$ is sufficiently large.
 
 The anti-commutator of $P$ and $J$ is the operator
 \begin{equation*}
     \{P,J\} \deff PJ + JP = \begin{bmatrix}
      &&-\mu_R' & -\lambda_R' \\
      &&-\lambda_R' &\mu_R' \\
      \mu_R' & \lambda_R' && \\
      \lambda_R' & -\mu_R' &&
     \end{bmatrix}\iota \partial_{x_i}
 \end{equation*}
 because $\lambda_R$ and $\mu_R$ are locally constant away from $\{-R < x_i <R\}$, and $P$ acts as a differential operator in $x_i$-direction on the strip $\{-R < x_i <R\}$.
 It is a bounded operator and, by increasing $R$, we may choose its operator norm to be arbitrarily small (but positive).
 However, if $n>1$, the anti-commutator is \emph{not} compact anymore, and this is the very reason why we now need to deviate from the proof presented in \cite{botvinnik2017infinite}.
 
 The proof that $P_s$ is Fredholm is similar to the proof of Proposition \ref{EllipticEstimateInvInf - Prop}. 
 However, $P_s$ is not a block Dirac operator and mixes sections of different bundles, so we cannot formally deduce the Fredholm property from Proposition \ref{EllipticEstimateInvInf - Prop}.
 However, due to the similarities, we will not repeat every detail.
 
 Let us fix some notation: An element $u \in H$ can be thought as a tuple with entries $u_{\alpha\beta} \in H_{\alpha\beta}$ and we will not notionally distinguish a vector with a single non-zero entry from its entry.
 
 We have, for all $u \in H^1$, the equation
 \begin{equation}\label{eq: NormRotatedOperator}
     ||P_su||^2_0 = \cos(s)^2 ||Pu||_0^2 + \sin(s)^2 ||u||_0^2 + \cos(s)\sin(s)\langle\{P,J\}u,u\rangle
 \end{equation}
 and analogously
 \begin{equation*}
    ||P_su_{\alpha\beta}||^2_0 = \cos(s)^2 ||P_{\alpha\beta}u_{\alpha\beta}||_0^2 + \sin(s)^2 ||u_{\alpha\beta}||_0^2 + \cos(s)\sin(s)\langle\{P,J\}u_{\alpha\beta},u_{\alpha\beta}\rangle. 
 \end{equation*}
 If $u_{\alpha\beta}$ is supported within $\{\omega x_j > R\}$ such that $P_{\alpha\beta}$ decomposes there, then we have, by the invertibility of $\face[\omega]{j}P_{\alpha\beta}$, the estimate
 \begin{align*}
     ||P_su_{\alpha\beta}||_0^2 &\geq c(j,\omega)^2 \cos(s)^2||u_{\alpha\beta}||^2_1 + \sin(s)^2 ||u_{\alpha\beta}||_0^2 + \cos(s)\sin(s) \langle\{P,J\}u_{\alpha\beta},u_{\alpha\beta}\rangle \\
     &\geq \underbrace{\left(c(j,\omega)^2 - ||\{P,J\}||_{0,0}\right)}_{>0}||u_{\alpha\beta}||_1^2
 \end{align*}
 because the constant $1>c(j,\omega)$ can be chosen independently of $R$, see Corollary \ref{Gluing Lower Bounds - Cor}.
 
 Assume now that $u_{\alpha\beta}$ is contained in a fixed compact neighbourhood $K$ whose interior containes the core of $P_{\alpha\beta}$.
 For $s \in [0,\pi/2)$, the classical elliptic estimate Lemma \ref{EllipticEstimateCpt - lemma} and the relation between $||P_su_{\alpha\beta}||_0^2$ and $||P_{\alpha\beta}u_{\alpha\beta}||_0^2$ described by equation (\ref{eq: NormRotatedOperator}) yield
 \begin{align*}
     ||u_{\alpha\beta}||_1^2 &\leq C\left( ||u_{\alpha\beta}||_0^2 + ||P_{\alpha\beta}u_{\alpha\beta}||_0^2\right) \\
     \begin{split}
          &= C\bigl(||u_{\alpha\beta}||_0^2 + \cos(s)^{-2}||P_su_{\alpha\beta}||_0^2 - \tan(s)^{-2}||u_{\alpha\beta}||_0^2 \\ 
          & \qquad \ - \tan(s)\langle\{P,J\}u_{\alpha\beta},u_{\alpha\beta}\rangle\bigr)
     \end{split} \\
     &\leq \cos(s)^{-2}(1 + ||\{P,J\}||_{0,0}) \cdot C(||u_{\alpha\beta}||_0^2 + ||Pu_{\alpha\beta}||^2).
 \end{align*}
 
 We now use the partition of unity $\{\phi_l\}_{l \in L}$ from Lemma \ref{BunkesPartOfOne - Lemma} to patch these inequalities together.
 As in the proof of Proposition \ref{EllipticEstimateInvInf - Prop}, we have
\begin{align*}
    ||u||_1^2 &= \sum_{\alpha\beta} ||u_{\alpha\beta}||_1^2 \leq (2n+1)\sum_{\alpha\beta}\sum_{l\in L} ||\phi_l u_{\alpha\beta}||_1^2 \\
    &\leq(2n+1) \sum_{\alpha\beta}\left( C_{\alpha\beta,0}(s) (||\phi_0u_{\alpha\beta}||_0^2 + ||P_s\phi_0u_{\alpha\beta}||_0^2) + \sum_{L \setminus 0} C_{\alpha\beta,l}(s) ||P_s\phi_lu_{\alpha\beta}||_0^2 \right) \\
    &\leq C(s) \sum_{\alpha\beta}\left( ||\phi_0u_{\alpha\beta}||_0^2 + ||P_s u_{\alpha\beta}||_0^2 + \sum_{L} ||[P_s,\phi_l \cdot]u_{\alpha\beta}||_0^2\right) \\
    &= C(s)\left( ||P_s u||_0^2 + ||\phi_0u||_0^2 + \sum_L ||[P_s,\phi_s \cdot]u||_0^2 \right).
\end{align*}
Since $[P_s,\phi_l \cdot] = \cos(s)[P,\phi_l \cdot] = \cos(s)\grad{ \phi_l} \cdot$ we end up with
\begin{equation*}
    ||u||_1^2 \leq C(s) \left( ||P_s u||_0^2 + ||\phi_0u||^2 + \sum_L ||\grad{\phi_l} \cdot u||_0^2\right).
\end{equation*}
As before, the operators $\phi_0 \cdot$ and $\grad{\phi_l}\cdot$ are compact operators because block matrices of compact operators are compact operators.
Lemma \ref{FAAuxilary - Lemma} now yields that $P_s$ has closed image and finite dimensional kernel.
Since $P_s$ is self-adjoint, it must be Fredholm.

In conclusion, $s \mapsto P_s$ is a continuous path between $P_0$ and the invertible operator $J = P_1$.
Homotopy invariance of the index implies $\ind(P) = \ind(J) = 0$ and the theorem is proven.
\end{proof}
 As in \cite{gromov1983positive}, we deduce the following consequence.
\begin{cor}\label{BunkesIndexThm - Cor}
 Let $(M^d,g_0)$ be a closed spin manifold such that $\Dirac_{g_0}$ is invertible and let $\chi \colon \R \rightarrow [0,1]$ with $\chi(t) = 0$ for $t \ll 0$ and $\chi(t) = 1$ for $t \gg 0$.
 The index difference map $\pi_n(\InvBlockDirac) \rightarrow KO^{-(d+n+1)}(\mathrm{pt})$ given by
 \begin{equation*}
     P \mapsto \mathrm{index}\bigl(\chi \cdot \susp(\fullGauge(g_P,g_\chi)_\ast P) + (1-\chi) \cdot \Dirac_{\susp(g_\chi)}\bigr), 
 \end{equation*}
 where $g_\chi = \chi \cdot g_P + (1-\chi)(g_0 \oplus \euclmetric_{\R^n})$ is the interpolating block map of Riemannian metrics,
 is a well-defined map, in particular, independent of the interpolation function $\chi$.
 Furthermore, the index map is a well-defined group homomorphism. 
\end{cor}

\section{Comparison to real \textit{K}-Theory}\label{Section - Comparison to K-theory}

After the effort we have made to construct $\InvBlockDirac$ and to prove several properties of it, we would like to know whether it is a classifying space for real $K$-theory.
Fortunately, this is the case, so we can give a precise formulation of Theorem \ref{BlockModelKTheory - MainThm}:
\begin{theorem}\label{OperatorSuspWeakEquiv - Theorem}
 Let $(M^d,g_0)$ be a spin manifold such that $\Dirac_{g_0}$ is an invertible operator. Then the suspension
 \begin{equation*}\label{SuspWhe - Theorem}
     \xymatrix{\InvPseudDir[\bullet] && B_\bullet \ar@{_{(}->}[ll]_-{\incl}^-{\simeq} \ar[rr]^-{\susp_\bullet} &&\InvBlockDirac }
 \end{equation*}
 is a weak homotopy equivalence.
\end{theorem}

Since the canonical map $\InvBlockDirac \rightarrow \BlockMetrics$ is a Kan fibration and the base space is combinatorially contractible, it suffices to show that the restriction to the fibre
\begin{equation*}
    \xymatrix{\InvPseudDir[\bullet]_{g_0} && B_{\bullet,g_0} \ar@{_{(}->}[ll]_-{\incl}^-{\simeq} \ar[rr]^-{\susp_\bullet} && \InvBlockDirac_{g_0} }
\end{equation*}
is a weak homotopy equivalence.

The proof will be carried out in two steps.
The first step is to compare the two different cubical sets of invertible pseudo Dirac operators to $KK$-theory.
More precisely, we will construct a diagram
\begin{equation}\label{eq: DiagramKKcomp}
    \xymatrix{ \pi_n(\InvPseudDir[\bullet]_{g_0}) \ar[r]^-\cong \ar[d]_{\pi_n(\susp_\bullet)} & KK(Cl_{0,d}, \mathcal{C}_0(\R^{n+1})) \ar[d]^{\tau_{Cl_{0,n+1}}(\placeholder)\sharp \alpha_{n+1}}_{\cong} \\ 
    \pi_n(\InvBlockDirac_{g_0}) \ar[r] & KK(Cl_{0,d+n+1},\R)}
\end{equation}
We will show in Theorem \ref{KKCommutative - Theorem} below that the diagram commutes, from which we deduce that $\pi_n(\susp_\bullet)$ is injective.

The second step is to show that $\pi_n(\susp_\bullet)$ factors as follows:
\begin{equation*}
    \xymatrix{\pi_n(\InvPseudDir[\bullet]_{g_0}) \ar@{^{(}->}[r] & KO^{-(d+n+1)}(\mathrm{pt}) \ar@{->>}[r] & \pi_n(\InvBlockDirac_{g_0}).}
\end{equation*}
If $n\geq 1$, then these maps are group homomorphisms, so, as all $KO$-theory groups are cyclic, $\pi_n(\susp_\bullet)$ must also be a surjective map.
If $n=0$, we show surjectivity of $\pi_0(\susp_\bullet)$ ``by hand''.

The idea to compare $\pi_n(\susp_\bullet)$ against $KK$-theory is inspired from the PhD thesis of Lukas Buggisch \cite{buggisch2019spectral}, where he shows that the index difference of Hitchin and Gromov-Lawson agree in the context of higher index theory.
More precisely, he establishes, for all smooth, closed manifolds $X$, a commutative diagram
\begin{equation*}
    \xymatrix{[X,\Riem^+(M)] \ar[rr]^-{\mathrm{inddif}_H} \ar[rrd]_-{\mathrm{inddif}_{GL}} && KK\bigl(\Cliff[0][d],\mathcal{C}_0(\R\times X,C^\ast \pi_1(M))\bigr) \ar[d]^{\tau_{\Cliff[0][1]}(\placeholder) \sharp \alpha_1} \\
    && KK\bigl(\Cliff[0][d+n+1],\mathcal{C}(X,C^\ast \pi_1(M))\bigr).}
\end{equation*}
Unfortunately, his results and proofs do not allow an inductive treatment. 
Furthermore, some of the proofs there rely on properties of the Dirac operator that do not carry over to general pseudo Dirac operators, like the Lichnerowicz formula.
This is the reason why we present a complete proof.

\subsubsection*{Kasparov KK-groups}

We recall the basic facts of $KK$-theory that are needed for our proof. 
We refrain from giving all details as they can be found, for example, in \cite{kasparov1980operator} and \cite{blackadar1998k}.

For separable, graded $C^\ast$-algebras $A, B$  (real, Real, or complex) the abelian group $KK(A,B)$ is defined as the set of homotopy classes of Kasparov ($A$-$B$)-modules $(E,\rho, F)$. 
Recall that a Kasparov module consists of a graded, separable Hilbert $B$-right-module $E$, a graded $\ast$-homomorphism $\rho \colon A \rightarrow \mathbf{Lin}_B(E)$ into the space of adjoinable, continuous $B$-right linear endomorphism of $E$, and an odd operator $F \in \mathbf{Lin}_B(E)$ that satisfies the following identities modulo compact operators:
\begin{equation*}
    \rho(a)(F-F^\ast) \equiv 0, \quad \rho(a)(F^2-\id) \equiv 0, \quad [\rho(a),F] \equiv 0 \qquad \mod \mathcal{K}_B.
\end{equation*}
Here, $\mathcal{K}_B$ is the set of $B$-compact operators, which is defined as the closed linear hull of all $\theta_{e_1,e_2} \deff e_1\cdot (\placeholder,e_2)$ for all $e_j \in E$.
Recall further that a homotopy between two Kasparov ($A$-$B$)-modules $(E_0,\rho_0,F_0)$ and $(E_1,\rho_1,F_1)$ is a Kasparov ($A$-$\mathcal{C}([0,1],B)$)-module that, up to unitary equivalence, restricts to the given ones.

These groups generalise real $K$-homology and real $K$-theory of $C^\ast$-algebras in the sense that
\begin{equation*}
    KO^d(A) \cong KK(A \otimes Cl_{0,d},\R) \quad \text{ and } \quad KO_d(A) \cong KK(\R, A \otimes Cl_{0,d}).
\end{equation*}
In fact, Kasparov proved that the groups $KK(A \otimes \Cliff[p_1][q_2], B \otimes \Cliff[p_2][q_2])$ only depend on $d = (p_1 - q_2) - (q_1 - p_2)$, see \cite{kasparov1980operator}*{Theorem 5.4}, and that the canonical homomorphism
\begin{align*}
    \tau_{\Cliff[p][q]} \colon KK(A,B) &\rightarrow KK(A \otimes \Cliff[p][q], B\otimes\Cliff[p][q]), \\
    [E,\rho,F] &\mapsto [E \otimes \Cliff[p][q], \rho\otimes \id, F \otimes \id]
\end{align*}
is an isomorphism \cite{kasparov1980operator}*{Proof of Theorem 5.2}.

The connection to our set up is provided by the following isomorphism  
\begin{align*}
    [\R^{n+1} \cup \{\infty\},\{\infty\}; \,  \mathrm{Fred}^{d,0}_0(H), GL_0(H)] &\rightarrow KK(Cl_{0,d},\mathcal{C}_0(\R^{n+1})) \\
    F &\mapsto [\mathcal{C}_0(\R^{n+1},H),\check{\rho}, F], 
\end{align*}
see \cite{ebert2017indexdiff}*{App. A}. 
Here, $H$ is a graded, ample Hilbert $Cl_{d,0}$-right-module with right action $\rho$ and grading $\iota$. 
The $Cl_{0,d}$ left-action $\check{\rho}$ is the unique one that is defined on the level of vectors by $\check{\rho}(v) = \rho(v) \circ \iota$.
The space $\mathrm{Fred}_0^{d,0}(H)$ is the weakly equivalent subspace consisting of all operators with $F^2 - \id$ compact and $GL_0(H)$ is the subspace of involutions \cite{ebert2017indexdiff}*{App. A}.

For later purposes, it will be useful to describe an explicit isomorphism between $KK(\Cliff[0][p],\R)$ and $KO^{-p}(\mathrm{pt})$.
 Let $\hat{\mathfrak{M}}_{r,s}$ be the abelian monoid of isomorphism classes of finitely generated, graded $\Cliff[r][s]$-modules.
 Recall that the monoids $\hat{\mathfrak{M}}_{r,s}$ and $\hat{\mathfrak{M}}_{s,r}$ are canonically isomorphic.
 Indeed, if $\rho \colon \Cliff[r][s] \rightarrow \mathrm{End}(W,\iota)$ is a graded algebra homomorphism, where $\iota$ is the grading, then the map $\R^{r,s} \ni v \mapsto \rho(v) \circ \iota$ extends to a graded algebra homomorphism $\check{\rho} \colon \Cliff[s][r] \rightarrow \mathrm{End}(W,\iota)$.
 The assignment $\rho \mapsto \check{\rho}$ provides an example of this isomorphism.
 
 Recall that $KO^{-p}(\mathrm{pt}) \cong \hat{\mathfrak{M}}_{p,0}/\incl^\ast \hat{\mathfrak{M}}_{p+1,0} \cong \hat{\mathfrak{M}}_{0,p}/\incl^\ast \hat{\mathfrak{M}}_{0,p+1}$.
 The first isomorphism can be found in \cite{LawsonMichelsonSpin}*{Section 1.9} or \cite{eschenburg2021bott}, the second one is induced by the isomorphism just constructed.
\begin{lemma}\label{PictKOHomology - Lemma}
 The map
 \begin{equation*}
     KK(\Cliff[0][p],\R) \rightarrow KO^{-p}(\mathrm{pt}) \quad \text{ given by } \quad  [E,\rho, F] \mapsto [\ker(F)]
 \end{equation*}
 is an isomorphism.
\end{lemma}
 \begin{proof}
  We may assume that $F$ is self-adjoint and graded-commutes with $\Cliff[0][p]$, see \cite{blackadar1998k}*{Section 17.4}. 
  Since $\Cliff[0][p]$ is a unital algebra, the identity 
  \begin{equation*}
      \rho(1)\cdot(F^2-\id) = (F^2-\id) \equiv 0 \quad \mod \mathcal{K}_{\Cliff[0][p]}
  \end{equation*}
  implies that $F$ is a Fredholm operator.
  The kernel is therefore a finite-dimensional $\Cliff[0][p]$-module.
  
  To see that the map is well-defined, we work with $KK_{oh}(\Cliff[0][p],\R)$ instead of $KK(\Cliff[0][p],\R)$.
  The former group is the group of equivalence classes of Kasparov modules of the equivalence relation generated by operator homotopy and by adding degenerate $KK$-modules.
  Recall that an operator homotopy is a $\Cliff[0][p]$-$\mathcal{C}([0,1],\R)$-module $(E,\rho,F_t)$, where the Hilbert space $E$ and the representation $\rho$ are independent of $t \in I$.
  The canonical map $KK_{oh}(\Cliff[0][p],\R) \rightarrow KK(\Cliff[0][p],\R)$ is an isomorphism, see \cite{blackadar1998k}*{p.148}.
  
  The operators of degenerate modules have trivial kernel, and two operator homotopic operators $F_0$, $F_1$ represent the same element in $\hat{\mathfrak{M}}_{0,p}/\incl^\ast \hat{\mathfrak{M}}_{0,p+1}$, see \cite{LawsonMichelsonSpin}*{Theorem III.10.8}.
  In conclusion, the map is well-defined.
  
  The map clearly surjective: Each element $[M,\rho]$ is the image of $[M,\rho,0] + x_{\mathrm{deg}}$ (the degenerate module $x_{\mathrm{deg}}$ is needed to assure that our Hilbert space is ample).
  
  To see that the map is also injective, assume that the $\Cliff[0][p]$-module structure on $\ker(F)$ extends to a $\Cliff[0][d+1]$-module structure.
  If $\pr_{\ker(F)} \colon E \rightarrow \ker(F)$ denotes the orthogonal projection onto the kernel, then it is easy to check that $[E,\rho, F + \rho(e_{n+1})\circ \pr_{\ker(F)}]$ is a degenerate Kasparov module.
  Furthermore, $[E,\rho, F_t]$ with $F_t \deff F + t\rho(e_{n+1})\circ \pr_{\ker(F)}$ is an operator homotopy between these two Kasparov modules.
  Thus, $[E,\rho,F]$ already represents the zero element in $KK(\Cliff[0][p],\R)$, so the map is injective.
 \end{proof}

The advantage of Kasparov's $KK$-groups is the existence of a Kasparov product, which is a bilinear pairing 
\begin{equation*}
    \sharp_D \colon KK(A_1,B_1 \otimes D) \otimes KK(A_2 \otimes D,B_2) \rightarrow KK(A_1 \otimes A_2, B_1 \otimes B_2).
\end{equation*}
There are elements 
\begin{equation*}
    \alpha_n = [L^2(\R^n, \Cliff[n][0]), \check{\CliffRight}, \DiracSt_{\R^n}(1 + \DiracSt_{\R^n}^2)^{-1/2} ] \in KK(\mathcal{C}_0(\R^n, \Cliff[0][n]), \R) = KO_n(\R^n)
\end{equation*}
and 
\begin{equation*}
    \beta_n = [\mathcal{C}_0(\R^n,\Cliff[0][n]), \lambda, \frac{\lambda(x)}{(1+||x||^2)^{1/2}}] \in KK(\R, \mathcal{C}_0(\R^n,\Cliff[0][n])) = KO^{-n}(\R^n),
\end{equation*}
where $\lambda$ denotes left multiplication, that satisfy
\begin{equation*}
    \alpha_n \sharp \beta_n = 1 \in KK\bigl(\mathcal{C}_0(\R^n,\Cliff[0][n]), \mathcal{C}_0(\R^n,\Cliff[0][n])\bigr) \text{ and } \beta_n \sharp \alpha_n = 1 \in KK(\R,\R).
\end{equation*}
A detailed proof for the complex case with $n=1$ can be found in \cite{cuntz2017ktheory}*{p.99-105}.
The real case for each $n$ can be reduced to the cited one because the Kasparov product commutes with complexification \cite{kasparov1980operator}*{Theorem 4.4.2} and the complexified elements $\alpha_n^{\C} = \alpha_1^{\C} \sharp_{\C} \dots \sharp_{\C}\alpha_1^{\C}$ as well as $\beta_n^{\C} = \beta_1^{\C} \sharp_{\C} \dots \sharp_{\C} \beta_1^{\C}$, see \cite{kasparov1980operator}*{p.546-547}.
Kasparov proved that 
\begin{equation*}
    (\placeholder) \sharp \alpha_n \colon  KK(A, B\otimes \mathcal{C}_0(\R^n;\Cliff[0][n])) \rightarrow KK(A,B)
\end{equation*}
is an isomorphism\footnote{Here, we applied the Kasparov product with $D = \mathcal{C}_0(\R^n;\Cliff[0][n])$, and $A_2 = B_2 = \R$.}, \cite{kasparov1980operator}*{Theorem 5.7}.

In general, it is hard to compute Kasparov products.
To reduce this difficulty, we will use unbounded Kasparov modules introduced by Baaj and Julg \cite{baaj1983theorie}, which we refer to as \emph{spectral triples} following the convention of \cite{cuntz2007topological}.

\begin{definition}
  Let $A$ and $B$ be separable, graded, real $C^\ast$-algebras.
  An $A$-$B$ \emph{spectral triple} $(E,\rho,D)$ consists of a graded Hilbert $B$-right module $E$, together with a graded $\ast$-homomorphism $\rho \colon A \rightarrow \mathbf{Lin}_B(E)$ into the set of continuous, adjoinable, $B$-right linear operators, and an odd, self-adjoint and regular $B$-right linear map $D$ that satisfies:
  \begin{itemize}
      \item[(1)] The graded commutator $[\rho(a),D]$ extends to a bounded operator for all elements $a$ in a dense subset of $A$.
      \item[(2)] The operator $\rho(a)(1+D^2)^{-1}$ is compact\footnote{Compact in the Hilbert $C^\ast$-module sense, $\mathrm{i.e.}$,  $\rho(a)(1+D^2)^{-1} \in \mathcal{K}_B \subset \mathbf{Lin}_B(E)$} for all $a \in A$. 
  \end{itemize}
  A \emph{homotopy} between two $A$-$B$ spectral triples $(E_j,\rho_j,D_j)$ with $j=0,1$ is an $A$-$\mathcal{C}([0,1];B)$ spectral triple $(E,\rho,D)$ that restricts to given ones.
  We will simply use the term \emph{spectral triple} if the $C^\ast$-algebras $A$ and $B$ are clear from the context and denote the set of homotopy classes of all spectral triples by $\Psi(A,B)$.
\end{definition}
Baaj and Julg proved that the bounded transform
\begin{equation*}
    b \colon \Psi(A,B) \rightarrow KK(A,B) \quad \text{given by} \quad (E,\rho,D) \mapsto \bigl(E,\rho,D/(1+D^2)^{1/2}\bigr)
\end{equation*}
is surjective \cite{baaj1983theorie}.
The following result of Kucerovski gives a sufficient criterium for a spectral triple to represent a Kasparov product.

\begin{theorem}[\cite{kucerovsky1997kk}*{Lemma 10, Theorem 13}]\label{Kucerovski - Theorem}
 Let $(E_1,\rho_1,D_1) \in \Psi(A,C)$ and $(E_2,\rho_2,D_2) \in \Psi(C,B)$ be two spectral triples. 
 Assume $E = (E_1 \otimes_{\rho_2} E_2, \rho_1 \otimes \id, D) \in \Psi(A,B)$ satisfies the following properties:
 \begin{itemize}
     \item[(i)] For all $x$ in a dense subset of $\rho_1(A)E_1$ of pure degree, the operators 
     \begin{equation*}
         DT_x - (-1)^{\mathrm{deg}(x)}T_xD_2 \quad \text{ and } \quad  D_2T_x^\ast - (-1)^{\mathrm{deg}(x)}T_x^\ast D
     \end{equation*}
     extend to bounded operators, where $T_x \deff x \otimes (\placeholder)$.
     \item[(ii)] $\mathrm{dom}\, D \subseteq \mathrm{dom}\, D_1$.
     \item[(iii)] There is a $c\in \R$ such that 
     \begin{equation*}
         \langle D x, (D_1\otimes_{\rho_2}\id) x \rangle_B + \langle (D_1\otimes_{\rho_2}\id) x, Dx \rangle_B \geq c \langle x, x \rangle_A
     \end{equation*}
     for all all $x \in \mathrm{dom} D$.
 \end{itemize}
 Then the Kasparov module $(E_1 \otimes_{\rho_2} E_2, \rho_1 \otimes \id, D/(1+D^2)^{1/2})$ represents the Kasparov product of $(E_1,\rho_1,D_1/(1+D_1^2)^{1/2})$ and $(E_2,\rho_2,D_2/(1+D_2^2)^{1/2})$. 
\end{theorem}

\subsubsection*{Construction of the comparison maps}

We will obtain the horizontal maps in diagram (\ref{eq: DiagramKKcomp}) by constructing spectral triples.
First, recall that a dice function adapted to $\rho > 0$ is a function 
\begin{equation*}
    \dice_\rho \colon \R^n \rightarrow \R_{\geq 0}
\end{equation*}
that $\partial_j \dice_{\rho}(x)\neq 0$ only if $|x_j| > \rho$.
We may choose (and do so) the function to be proper, to vanish near the origin, and to have a bounded total differential, see Section \ref{Section - Construction of Special Submanifolds} for a construction.

$ $

\noindent\textbf{Construction of} $\pi_n(\InvPseudDir[\bullet]_{g_0}) \rightarrow \Psi(\Cliff[0][d],\mathcal{C}_0(\R^{n+1}))$:

\noindent Pick a smooth  function $\chi \colon \R \rightarrow [0,1]$ with $\chi(t) = 1$ for $t \gg 0$ and $\chi(t) = 0$ for $t \ll 0$.
For a given element $P \in \Omega^n_{\Dirac_{g_0}}(\InvPseudDir[\bullet]_{g_0})$, which is a smooth block map $P \colon \R^n \rightarrow \InvPseudDir_{g_0}$ that restricts to $\Dirac_{g_0}$ outside of $\rho I^n$ for some $\rho > 0$, define the block map
\begin{equation*}
    Q \colon \R^{n+1} = \R^n \times \R \rightarrow \InvPseudDir_{g_0} \quad \text{ via } \quad (x,t) \mapsto \chi(t)P + (1-\chi(t))\Dirac_{g_0}.
\end{equation*}
For the chosen $\rho$, define 
\begin{equation*}
    h = h_\rho \colon \R^{n+1} \rightarrow \R_{\geq 1} \quad \text{ via } \quad x \mapsto (1 + \dice_\rho(x)^2)^{1/4}.
\end{equation*}

\begin{lemma}
  Let $\CliffRight$ be the canonical $\Cliff[d][0]$-right action on $\Spinor_{g_0}$ and let $\check{\CliffRight}$ be the induced $\Cliff[0][d]$-left action. 
  The triple $(\mathcal{C}_0(\R^{n+1},L^2(M,\Spinor_{g_0})), \check{\CliffRight}, hQh)$ is a $\Cliff[0][d]$-$\mathcal{C}_0(\R^{n+1})$ spectral triple.
\end{lemma}
\begin{proof}
 The canonical isomorphism 
 \begin{equation*}
     L^2(\Spinor_{g_0}) \otimes_{\R} \mathcal{C}_0(\R^{n+1}) \cong \mathcal{C}_0(\R^{n+1};L^2(\Spinor_{g_0})) \quad \text{ induced by } \quad f \otimes \sigma \mapsto f \cdot \sigma ,
 \end{equation*}
 where $\mathcal{C}_0(\R^{n+1})$ is trivially graded, identifies the right hand side as a  Hilbert $\Cliff[0][d]$-$\mathcal{C}_0(\R^{n+1})$-bimodule.
 
 Considered as an operator on $\mathcal{C}_0(\R^{n+1},L^2(\Spinor_{g_0}))$, the block map $hQh$ is self-adjoint as a $\mathcal{C}_0(\R^{n+1})$-right linear map because it is point-wise self-adjoint.
 
 The operator-valued map $Q$ is $\Cliff[d][0]$-right linear, so it graded commutes with the associated $\Cliff[0][d]$-left action.
 Indeed, it suffices to prove it for vectors $\R^d \subseteq \Cliff[0][d]$, which is easily calculated
 \begin{align*}
     Q \circ \check{\CliffRight}(v) &= Q \circ \CliffRight(v) \circ \evenodd = \CliffRight(v) \circ Q \circ \evenodd \\
     &= - \CliffRight(v) \circ \evenodd \circ Q = -\check{\CliffRight}(v) \circ Q.
 \end{align*}
 Thus, the first condition of a spectral triple is satisfied.
 
 Since $Q$ is a block map that takes values in invertible operators outside a compact subset of $\R^{n+1}$, there is a $c_Q > 0$ such that $||Q_x \sigma||_0^2 \geq c_Q^2||\sigma||_0^2$ for all $\sigma \in L^2(\Spinor_{g_0})$ and all $x$ in the complement of this compact subset.
 Together with $hQ = Qh$, this implies 
 \begin{align*}
     (hQh)^2 &= h^4 Q^2 = (1 + \dice_\rho(x)^2) Q^2 \geq c_Q^2(1 + \dice_\rho(x)^2)
 \end{align*}
 for all $x$ outside the given compact subset.
 Since the manifold $M$ is closed and since $hQh \colon \R^{n+1} \rightarrow \PseudOp^1(\Spinor_{g_0})$ takes values in elliptic operators, $(1 + hQh)^{-1}$ is a (norm-)continuous map of compact operators that vanishes near $\infty$ by the previous inequiality.
 Thus, $(1 + hQh)^{-1}$ is a $\mathcal{C}_0(\R^{n+1})$-compact operator, so the second condition of a spectral triple is also satisfied.
\end{proof}

 Different choices yield homotopic spectral triples\footnote{to spell it out: Different choices of $\dice_\rho$, the representative $P$ of $[P]$, and the interpolating path between $P$ and $\Dirac_{g_0}$ can be joined by paths and these paths produce operator homotopies of spectral triples.}, 
 so we get a well-defined map
\begin{align*}
    \pi_n(\InvPseudDir[\bullet]_{g_0}) \rightarrow \Psi(\Cliff[0][d], \mathcal{C}_0(\R^{n+1})) \quad \text{ given by } \quad [P] \mapsto [hQh].
\end{align*}
\begin{lemma}\label{UpperArrowIso - Lemma}
 The composition with the bounded transform 
 \begin{equation*}
     \xymatrix{\pi_n(\InvPseudDir[\bullet]_{g_0}) \ar[r] & \Psi(\Cliff[0][d], \mathcal{C}_0(\R^{n+1})) \ar[r]^b & KK(\Cliff[0][d], \mathcal{C}_0(\R^{n+1}))  }
 \end{equation*}
 is an isomorphism if $n \geq 1$ and a bijection if $n=0$.
\end{lemma}
\begin{proof}
 Abbreviate $\R^{n+1} \cup \{\infty\}$ to $S^{n+1}$ and $\mathrm{Fred}^{d,0}(L^2(\Spinor_{g_0}))$ to $\mathrm{Fred}^d$.
 Consider the following diagram
 \begin{equation*}
     \xymatrix{ &&  [S^{n+1}, \{\infty\}; \, \mathrm{Fred}_0^{d},GL_0] \ar[r]^-\cong \ar@{^{(}->}[d]^{\mathrm{incl}} &  KK(\Cliff[0][d], \mathcal{C}_0(\R^{n+1}))\\
     \pi_n(\InvPseudDir[\bullet]_{g_0}) \ar[rr]_-{P \mapsto b(Q)} \ar[rru]^-{P \mapsto b(hQh)} && [S^{n+1},\{\infty\}; \, \mathrm{Fred}^d, GL]. }
 \end{equation*}
 The inclusion induces a bijection on homotopy classes by a spectral deformation argument, see \cite{ebert2017indexdiff}*{Appendix A} or Atiyah-Singer \cite{atiyah1969indexskew}*{p.12}. 
 The lower horizontal map represents the index difference, see Theorem \ref{InvPseudDir Classifying Space - Thm}, 
 so it is also an isomorphism, respectively, a bijection.
 Lastly, the upper horizontal map is an isomorphism that is nicely described in \cite{ebert2017indexdiff}*{Appendix A}.
 If $H = (H,\rho,\iota)$ is a Hilbert space with an ample $\Cliff[d][0]$-right action $\rho$ and grading $\iota$, then the map is given by
 \begin{equation*}
     \bigl(F \colon (S^{n+1},\infty) \rightarrow (\mathrm{Fred}^{d,0},GL_0)\bigr) \quad \mapsto \quad \bigl(\mathcal{C}_0(\R^{n+1},H), \Check{\rho}, F \bigr) \in \Psi(\Cliff[0][d],\mathcal{C}_0(\R^{n+1})). 
 \end{equation*}
 From this description, it is easy to see that the composition agrees with the composition from the statement.
 
 It remains to prove that the diagram commutes.
 Since $h$ and $Q$ commute and $h$ is a positive function, the operators $Q$ and $hQh$ have the same kernel and the same positive and negative subspace.
 Functional calculus now implies that  
 \begin{equation*}
     s \mapsto s b(Q) + (1-s)b(hQh)
 \end{equation*}
 is a homotopy between $\incl \circ b(Q)$ and $b(hQh)$ with values in Fredholm operators because the kernel is independent of $s$ and finite dimensional.
\end{proof}

\noindent\textbf{Construction of }$\pi_n(\InvBlockDirac_{g_0}) \rightarrow \Psi(\Cliff[0][d+n+1],\R)$:

\noindent Pick a smooth function $\chi \colon \R \rightarrow [0,1]$ that satisfies $\chi(t) = 1$ for $t \gg 0$ and $\chi(t) = 0$ for $ t \ll 0$.
For each element $P \in \Omega_{g_0}^n\InvBlockDirac_{g_0}$, we define a block map of block Dirac operators
\begin{equation*}
    Q \colon \R \rightarrow \BlockDirac[n]_{g_0} \quad \text{ by } \quad t \mapsto \chi(t) P + (1 - \chi(t)) \Dirac_{g_0 \oplus \euclmetric_{\R^n}}.
\end{equation*}
Its suspension 
\begin{equation*}
    \susp(Q) = \chi \circ \mathrm{pr}_{n+1} \cdot P^\extension + (1 - \chi \circ \mathrm{pr}_{n+1}) \cdot \Dirac_{g_0 \oplus \euclmetric_{\R^n}} + \partial_{x_{n+1}} \cdot \frac{\partial}{\partial x_{n+1}}
\end{equation*}
is a block Dirac operator on $\Spinor_{g_0 \oplus \euclmetric_{\R^{n+1}}} \rightarrow M \times \R^{n+1}$.
For $\rho > 0$ such that the core of $\susp(Q)$ is contained in $\rho I^{n+1}$, define the function
\begin{equation*}
    h \deff h_\rho \deff (1 + \dice_\rho^2)^{1/4} \colon \R^{n+1} \rightarrow \R.
\end{equation*}

The action $\check{\CliffRight}$ turns $L^2(M \times \R^{n+1};\Spinor_{g_0 \oplus \euclmetric})$ into a Hilbert $\Cliff[0][d+n+1]$-$\R$-bimodule on which $h\susp(Q)h$ acts as an odd, symmetric, unbounded operator whose domain is the dense subspace of all compactly supported smooth sections.

\begin{lemma}
 The triple $(L^2(M\times \R^{n+1},\Spinor_{g_0 \oplus \euclmetric}), \check{\CliffRight}, h \susp(Q)h )$ is a $\Cliff[0][d+n+1]$-$\R$ spectral triple.
\end{lemma}
\begin{proof}
 We need to show that 
 \begin{itemize}
     \item[(o)] the closure of operator $h\susp(Q)h$ is self-adjoint, 
     \item[(i)] that the graded commutator $[\check{\CliffRight}(v),h\susp(Q)h]$ extends to a bounded operator for all $v\in \Cliff[0][n+d+1]$, and
     \item[(ii)] that $\check{\mathbf{r}}(v)(1+(hQh)^2)^{-1}$ is compact for all $v \in \Cliff[0][n+d+1]$.
 \end{itemize}  
 
 The prove of (o) is similar to the proof of \cite{higson2000analytic}*{Proposition 10.2.10}:
 Since $\susp(Q)$ decomposes outside of $\rho I^{n+1}$, it acts as differential operator on $\dice_\rho$ and $h = h_\rho$.
 Thus,
 \begin{align*}
     [h\susp(Q)h,\log(1 + \dice_\rho)] &= [h^2\susp(Q),\log(1+\dice_\rho)] + [\grad{h} \cdot, \log(1+\dice_\rho)] \\
     &= h^2[\susp(Q), \log(1+\dice_\rho)] \\
     &= h^2 \grad{\log(1+\dice_\rho)}\cdot \\
     &= (1+\dice_\rho^2)^{1/2}(1+\dice_\rho)^{-1} \grad{\dice_\rho} \cdot .
 \end{align*}
 We chose $\dice_\rho$ to be proper and to have bounded derivative, so the commutator $[h\susp(Q)h,\log(1+\dice_\rho)]$ extends to a bounded operator.
 With the help of $\log(1+\dice_\rho)$ we find smooth functions $\varphi_n \colon M\times \R^{n+1} \rightarrow [0,1]$ such that $\varphi_n \rightarrow 1$ uniformly on a compact set and such that $[h\susp(Q)h,\varphi_n] \rightarrow 0$ in norm.
 
 Assume now that $u \in L^2(\Spinor_{g_0 \oplus \euclmetric_{\R^{n+1}}})$ lies in the maximal domain.
 Since $\susp(Q)$ is a block operator, it maps sections supported in the interior $M \times RI^{n+1}$ to sections that are also supported there, hence
 \begin{align*}
     h\bigl(\susp(Q)\bigr)h fu = f\cdot h\bigl(\susp(Q)\bigr)h u + [h\bigl(\susp(Q)\bigr)h,f]u
 \end{align*}
 is also bounded for all compactly supported functions $f$.
 The operator $h(\susp(Q)-\Dirac_{g_0\oplus\euclmetric_{\R^{n+1}}})h$ is bounded on compact sets, so the proof of \cite{higson2000analytic}*{Lemma 10.2.5} implies that $fu$ is also in the minimal domain of $h\susp(Q)h$, because we can choose Friedrichs' mollifiers in a manner so that they map sections supported in $\supp(f)$ to sections supported in a fixed compact neighbourhood $L$ of $\supp(f)$.
 In particular, $\varphi_nu$ belongs to the minimal domain of $h\susp(Q)h$.
 Since $u$ is square integrable, we have $\varphi_nu\xrightarrow{n \to \infty} u$ in norm, while
 \begin{equation*}
     h\susp(Q)h\varphi_n u = \varphi_n h\susp(Q)h u + [h\susp(Q)h,\varphi_n]u \xrightarrow{n \to \infty} u,
 \end{equation*}
 so $u$ belongs to the minimal domain of $h\susp(Q)h$.
 
 Condition (i) is satisfied because $[\check{\CliffRight}(v),h\susp(Q)h] = 0$.
 
 The proof of condition (ii) is borrowed from \cite{ebert2016elliptic}*{Theorem 2.40}.
 The difference to the proof presented there is that we avoid to use coercive functions.
 
 Note that the restriction of $\susp(Q)$ to $\{\eps x_i > \rho\}$ is an operator that is bounded from below in the sense of Definition \ref{boundedbelow - Definition} for all $(i,\eps) \in \{1,\dots,n+1\}\times \Z_2$.
 Corollary \ref{Gluing Lower Bounds - Cor} applied inductively to all restrictions of $\susp(Q)$ implies $||\susp(Q)u||_0 \geq c ||u||_0$ for some $c > 0$ and all sections $u \in \Gamma_c(M \times \R^{n+1},\Spinor_{g_0 \oplus \euclmetric})$ that are supported in the complement of $\rho' I^{n+1}$ for all sufficiently large $\rho'\geq \rho$.
 To keep the notation readable, we assume without loss of generality $\rho' = \rho$.
 Thus, if $\chi$ is a smooth function with compact support in the complement of $\rho I^{n+1}$, then we have
 \begin{align*}
     ||h\susp(Q)h \chi u||_0 \cdot ||\chi u||_0 \geq c ||h\chi u||_0 \cdot ||\chi u||_0 &\geq c \inf\{h(t) \, : \, t \in \supp \chi\} ||\chi u||_0^2 \\
     &=: \mathrm{low} \cdot ||\chi u||_0^2.
 \end{align*}
 
 Since $\susp(Q)$ is a block Dirac operator over $g_0\oplus \euclmetric_{\R^{n+1}}$, the operator $\susp(Q) - \Dirac_{g_0 \oplus \euclmetric_{\R^{n+1}}}$ is bounded, which implies that the commutator $[\susp(Q),\chi]$ is a bounded operator. 
 By Lemma \ref{CoreLowerBound - Lemma}, $\susp(Q)u$ is supported in $M \times RI^{n+1}$ if $u$ is supported there for all sufficient large $R>0$.
 Thus, if $\supp \chi$ is contained in $RI^{n+1}$, then the operator $[\susp(Q),\chi]$ is supported in $M \times RI^{n+1}$, too.
 It follows that $[h\susp(Q)h,\chi] = h[\susp(Q),\chi]h$ also extends to a bounded operator whose operator norm is bounded by
 \begin{equation*}
     ||[h\susp(Q)h,\chi]||_{0,0} \leq \sup\{h(t)^2 \, : \, t \in \supp([\susp(Q),\chi])\} \cdot ||[\susp(Q) , \chi]||_{0,0}.
 \end{equation*}
 Combining the two inequalities yields
 \begin{align*}
     ||\chi u||_0 &\leq \mathrm{low}^{-1} ||h\susp(Q)h\chi u||_0 \\
     &\leq \mathrm{low}^{-1}(||[h\susp(Q)h,\chi]u||_0 + ||\chi h\susp(Q) h u||_0 ) \\
     &\leq \mathrm{low}^{-1}\left( ||[h\susp(Q)h,\chi]||_{0,0} \cdot ||u||_0 + ||\chi||_{\infty}||h\susp(Q)hu||_0\right) \\
     &\leq \mathrm{low}^{-1}\max\{||[h\susp(Q)h,\chi]||_{0,0}, ||\chi||_\infty\} \cdot\left(||u||_0 + ||h\susp(Q)h||_0\right) \\
     &\leq \mathrm{low}^{-1} \sqrt{2} \max\{||[h\susp(Q)h,\chi]||_{0,0}, ||\chi||_\infty\} \left(||u||_0^2 + ||h\susp(Q)h||_0^2\right)^{1/2}. 
 \end{align*}
 If we interpret $Q$ as a complex linear operator on the complexification $\Spinor_{g_0 \oplus \euclmetric}\otimes \C$, we can extend the previous chain of inequalities
 \begin{align*}
     ||\chi u||_0 &\leq\mathrm{low}^{-1} \sqrt{2} \max\{||[h\susp(Q)h,\chi]||_{0,0}, ||\chi||_\infty\} \left(||u||_0^2 + ||h\susp(Q)h||_0^2\right)^{1/2} \\
     &= \mathrm{low}^{-1} \sqrt{2} \max\{||[h\susp(Q)h,\chi]||_{0,0}, ||\chi||_\infty\} ||(h\susp(Q)h + i)u||_0
 \end{align*}
 for all real sections $u \in \Gamma_c(\Spinor_{g_0 \oplus \euclmetric}\otimes \C)$.
 Since $h\susp(Q)h$ is self-adjoint, its spectrum is real so that $(h\susp(Q)h\pm i)^{-1}$ exists and its operator norm is bounded by $1$.
 Thus, by substituting $u = (h\susp(Q)h+i)^{-1}v$, the previous inequality becomes
 \begin{equation*}
     || \chi(h\susp(Q)h+i)^{-1}v||_0 \leq \mathrm{low}^{-1}\sqrt{2}  \max\{||[h\susp(Q)h,\chi]||_{0,0}, ||\chi||_\infty\} ||v||_0,
 \end{equation*}
 which implies
 \begin{align*}
     ||\chi \cdot (1+(h\susp(Q)h)^{-2})||_{0,0} &\leq ||\chi (h\susp(Q)h + i)^{-1}||_{0,0} \cdot ||(h\susp(Q)h + i)^{-1}||_{0,0} \\
     &\leq \mathrm{low}^{-1} \sqrt{2} \max\{||[h\susp(Q)h,\chi]||_{0,0}, ||\chi||_\infty\}.
 \end{align*}
 
 Pick $\tilde{\eta}_0 \colon \R \rightarrow [0,1]$ that is compactly supported and identically $1$ on the interval $[0,a]$, where $a \deff \max\dice_\rho((\rho+1)I^{n+1})$.
 If $||\diff \, \dice_{\rho}||_\infty \leq B$, then the mean value theorem gives $|\dice_\rho(x)| \leq B||x||$ for all $x \in \R^{n+1}$.
 For all $m \in \N$, pick a smooth function $\tilde{\eta}_m \colon \R_{\geq 0} \rightarrow [0,1]$ that satisfies 
 \begin{align*}
     \tilde{\eta}_m|_{[0,a+(m-1)m/2]} = 1, \  \tilde{\eta}_m|_{[a+m(m+1)/2,\infty)} = 0, \  \text{and } \left|\tilde{\eta}_m'|_{a+[(m-1)m/2,m(m+1)/2]}\right| \leq 2/m,  
 \end{align*}
 where the interval $[(m-1)m/2,m(m+1)/2]$ has length $m$.
 Set $\eta_m \deff \tilde{\eta}_m \circ \dice_\rho$.
 Since $\eta_m$ is constant on $\rho I^{n+1}$, we have $[\susp(Q),h] = \grad{h}$.
 Recall that $h_\rho = (1 + \dice_\rho^2)^{1/4}$.
 We estimate
 \begin{align*}
     \bigl|\bigl|[h\susp(Q)h,\eta_m]\bigr|\bigr|_{0,0} &= \bigl|\bigl|[h^2\susp(Q),\eta_m] + [h\grad{h}\cdot,\eta_m]\bigr|\bigr|_{0,0} \\
     &=||h^2\grad{\eta_m}\cdot||_{0,0} = ||h^2 \tilde{\eta}_m' \circ \dice_\rho \cdot  \grad{\dice_\rho}||_{0,0} \\
     &\leq ||h^2||_{\infty,\supp \eta_m} \cdot ||\eta'_m \circ \dice_\rho \cdot \grad{\dice_\rho}||_{0,0} \\
     &\leq \left(1+\left(\frac{a+m(m+1)}{2}\right)^2\right)^{1/4} \frac{2B}{m}  \\
     &\leq \mathrm{C}_a \cdot B.
 \end{align*}
 
 The operator $\eta_0 (1+(h\susp(Q)h)^2)^{-1}$ is compact because it factors through $H^1_{\supp \eta_0}$, the space of elements in $H^1$ supported within $\supp \eta_0$,
 so it remains to check that $(1-\eta_0)(1+(h\susp(Q)h)^2)^{-1}$ is compact.
 Clearly, the sequence converges pointwise: 
 \begin{equation*}
     \eta_m(1-\eta_0)(1+(h\susp(Q)h)^2)^{-1} \xrightarrow{m \to \infty} (1-\eta_0)(1+(h\susp(Q)h)^2)^{-1}. 
 \end{equation*}
 The limit is also the norm limit of $\eta_m(1-\eta_0)(1+(h\susp(Q)h)^2)^{-1}$ by the following estimate:
 \begin{align*}
     & \quad \  ||(\eta_l - \eta_m)(1-\eta_0)(1+(h\susp(Q)h)^2)^{-1}||_{0,0} \\
     &= ||(\eta_l - \eta_m)(1+(h\susp(Q)h)^2)^{-1}||_{0,0} \\
     &\leq \mathrm{low}^{-1}\sqrt{2} \max\{||h\susp(Q)h||_{0,0},||\eta_l - \eta_m||_\infty\} \\
     &\leq \mathrm{low}^{-1}\sqrt{2}\mathrm{C}_a\max\{B,1\} \\
     &\leq \sup\{h(t)^{-1} \, : t \in \supp(\eta_l - \eta_m)\}/c \cdot \sqrt{2}\mathrm{C}_a \max\{B,1\} \xrightarrow{m \to \infty} 0,
 \end{align*}
 because $h$ is proper and $\supp (\eta_l -\eta_m)$ lies in the complement of $(a+m(m+1)/2)I^{n+1}$.
 
 In conclusion, $(1-\eta_0)(1+(h\susp(Q)h)^2)^{-1}$ is compact as it is the norm limit of a sequence of compact operators.
\end{proof}

Again, different admissible choices yield homotopic spectral triples, so we get a well-defined map
\begin{equation*}
    \pi_n(\InvBlockDirac_{g_0},\Dirac_{g_0\oplus \euclmetric}) \rightarrow \Psi(\Cliff[0][d+n+1],\R) \quad \text{via} \quad [P] \mapsto [h\susp(Q)h].
\end{equation*}

\begin{lemma}
 The map just defined composed with the bounded transform $b$ is a group homomorphism provided $n \geq 1$.
\end{lemma}
\begin{proof}
 Consider the diagram
 \begin{equation*}
     \xymatrix{ \pi_n(\InvBlockDirac_{g_0}) \ar[rr]^{[P]\mapsto [hQh]} \ar[rrrd]_{[P] \mapsto \ker(Q)} &&\Psi(\Cliff[0][d+n+1],\R) \ar[r]^-{b} &KK(\Cliff[0][d+n+1],\R) \ar[d]^{\ker(\placeholder)}_\cong \\ &&& KO^{-(d+n+1)}(\mathrm{pt}). }
 \end{equation*}
 The vertical map is the one from Lemma \ref{PictKOHomology - Lemma} and therefore an isomorphism.
 The diagonal map is a group homomorphism by Corollary \ref{BunkesIndexThm - Cor}.
 Thus, the horizontal composition must be a group homomorphism, too.
\end{proof}

\noindent\textbf{Unbounded representatives of} $\alpha_{n+1}$:

\noindent The spectral triple $(L^2(\R^{n+1},\Cliff[n+1][0]), \check{\CliffRight}, \DiracSt_{\R^{n+1}}) \in \Psi(\mathcal{C}_0(\R^{n+1},\Cliff[0][n+1]),\R)$ is well known. 
Its bounded transform represents the fundamental class $\alpha_{n+1} \in KO_{n+1}(\R^{n+1})$ because its complexification represents the fundamental class in $K_{n+1}(\R^{n+1})$, see \cite{higson2000analytic}*{p.314ff}, and complexification induces an isomorphism $KO_{n+1}(\R^{n+1}) \rightarrow K_{n+1}(\R^{n+1})$.

When we calculate the Kasparov product, we will use $h\DiracSt_{\R^{n+1}}h$ instead of $\DiracSt_{\R^{n+1}}$. 
This is unproblematic because the bounded transform of these two operators over $(L^2(\R^{n+1},\Cliff[n+1][0]), \check{\CliffRight})$ represent the same class in $KK(\mathcal{C}_0(\R^{n+1},\Cliff[0][n+1]),\R)$.
Indeed, pick a small open ball $rD^{n+1} \subseteq \R^{n+1}$ such that $h = 1$ on this set.
Extension by zero $j \colon \mathcal{C}_0(rD^{n+1}) \hookrightarrow \mathcal{C}_0(\R^{n+1})$ is a homotopy equivalence (in the sense of $C^\ast$-algebras) because $rD^{n+1}$ and $\R^{n+1}$ are properly homotopy equivalent.
By \cite{higson2000analytic}*{Proposition 10.8.8} and that $h\DiracSt_{\R^{n+1}}h$ and $\DiracSt_{\R^{n+1}}$ agree on $rD^{n+1}$, we have 
\begin{equation*}
    j^\ast \bigl(b[\DiracSt_{\R^{n+1}}]\bigr) = [\![\DiracSt_{\R^{n+1}}|_{rD^{n+1}}]\!] = j^\ast \bigl(b[h\DiracSt_{\R^{n+1}}h]\bigr) \in KK(\mathcal{C}_0(rD^{n+1},\Cliff[0][n+1]),\R).
\end{equation*}
Here, $[\![\DiracSt_{\R^{n+1}}|_{rD^{n+1}}]\!]$ denotes the K-homology class that is defined in \cite{higson2000analytic}*{Definition 10.8.3}.
The claim follows now from homotopy invariance of $KK$-theory, see \cite{blackadar1998k}*{Proposition 17.9.1}.

From now on, we will denote a spectral triple only by its operator if the underlying Hilbert $C^\ast$-module and the representation are clear from the context.
 
\subsubsection*{Calculation of the Kasparov Product}
We have constructed the horizontal homomorphisms in diagram (\ref{eq: DiagramKKcomp}), so it remains to show that it commutes.
If we choose the same interpolation function $\chi$ in the construction of the upper and lower horizontal arrow, we have
\begin{equation*}
    \susp\left(\chi\susp(P) + (1-\chi)\Dirac_{g_0\oplus \euclmetric}\right) = \susp(\chi P + (1-\chi)\Dirac_{g_0}).
\end{equation*}
It thus remains to show $\tau_{\Cliff[0][n+1]}(b([hQh])) \sharp \alpha_{n+1} = b([h\susp(Q)h])$.

The first step is to observe that the underlying Hilbert $C^\ast$-modules are canonically isomorphic
\begin{align*}
    &\quad \  \mathcal{C}_0(\R^{n+1},L^2(\Spinor_{g_0})) \otimes \Cliff[0][n+1] \otimes_{\mathcal{C}_0(\R^{n+1}) \otimes \Cliff[0][n+1],\mathrm{\mathrm{taut}}\otimes \check{\CliffRight}} L^2(\R^{n+1}, \Cliff[n+1][0])  \\
    &\cong \mathcal{C}_0(\R^{n+1},L^2(\Spinor_{g_0})) \otimes_{\mathcal{C}_0(\R^{n+1}),\mathrm{taut}} L^2(\R^{n+1},\Cliff[n+1][0]) \\
    &\cong L^2(M \times \R^{n+1};\Spinor_{g_0}\boxtimes \Cliff[n+1][0]) \cong  L^2(M \times \R^{n+1}, \Spinor_{g_0 \oplus \euclmetric}),
\end{align*}
where $\mathrm{taut}$ is the left action of $\mathcal{C}_0(\R^{n+1},\R)$ on $L^2(\R^{n+1},\Cliff[n+1][0])$ via left multiplication.
The first isomorphism is the universal one, the second isomorphism is induced by 
\begin{align*}
    \mathcal{C}_c(\R^{n+1},L^2(\Spinor_{g_0})) \otimes_{\mathcal{C}_0(\R^{n+1})} \Gamma_c(\R^{n+1},\Cliff[n+1][0]) &\rightarrow \Gamma_c(M \times \R^{n+1}, \Spinor_{g_0} \boxtimes \Cliff[n+1][0]), \\
    \sigma \otimes \varphi &\mapsto \left((m,x) \mapsto \sigma(x)(m) \otimes \varphi(x)\right).
\end{align*}

The next step is to show that $[h\susp(Q)h]$ represents the Kasparov product of $\tau_{\Cliff[0][n+1]}(b([hQh]))$ and $\alpha_{n+1} = b([h\DiracSt_{\R^{n+1}}h])$ by showing that it satisfies the conditions of Theorem \ref{Kucerovski - Theorem}.

\begin{theorem}\label{KKCommutative - Theorem}
 For each $[P] \in \pi_n(\InvPseudDir[\bullet]_{g_0})$, let $Q \deff \chi P + (1-\chi)\Dirac_{g_0}$ for some $\chi \colon \R \rightarrow [0,1]$ with $\chi(t) = 1$ for $t \gg 0$ and $\chi(t) = 0$ for $t \ll 0$.
 Then we have
 \begin{equation*}
     \tau_{\Cliff[0][n+1]}(b([hQh])) \sharp \alpha_{n+1} = b(\tau_{\Cliff[0][n+1]}([hQh]))\sharp \alpha_{n+1} = b([h\susp(Q)h]).
 \end{equation*}
\end{theorem}
\begin{proof}
 To verify the first condition of Theorem \ref{Kucerovski - Theorem}, it suffices to show that 
 \begin{equation*}
     h\susp(Q)h T_\sigma + (-1)^{\mathrm{deg}(\sigma)} T_\sigma h\DiracSt_{\R^{n+1}}h
 \end{equation*}
 extends to a bounded operator on $L^2(\R^{n+1},\Cliff[n+1][0])$ for all $\sigma \in \Gamma_c(M \times \R^{n+1},\Spinor_{g_0})$ of pure degree because, in our case, the second operator in condition $(i)$ of Theorem \ref{Kucerovski - Theorem} is the adjoint of the first operator.
 Recall that $T_\sigma$ was defined as
 \begin{equation*}
    T_\sigma \deff \sigma \otimes (\placeholder) \colon L^2(\R^{n+1},\Cliff[n+1][0]) \rightarrow L^2(M \times \R^{n+1};\Spinor_{g_0 \oplus \euclmetric}). 
 \end{equation*}
 We abuse notation by abbreviating $\evenodd \boxtimes \DiracSt_{\R^{n+1}}$ to $\DiracSt_{\R^{n+1}}$.
 The Leibniz rule implies that the graded commutator $[\DiracSt_{\R^{n+1}},T_\sigma]$ is given by
 \begin{align*}
     [\DiracSt_{\R^{n+1}},T_\sigma] &= \DiracSt_{\R^{n+1}}T_\sigma - (-1)^{\mathrm{deg}(\sigma)} T_\sigma \DiracSt_{\R^{n+1}} \\
     &= (-1)^{\mathrm{deg}(\sigma)}\sum_{j=1}^{n+1}  \partial_{x_j}(\sigma) \otimes \partial_j \cdot \\
     &=: \DiracSt_{\R^{n+1}}(\sigma) \cdot .
 \end{align*}
 We use this identity to compute
 \begin{align*}
    &\quad \ \, h\susp(Q)h T_\sigma - (-1)^{\deg(\sigma)} T_\sigma h\DiracSt_{\R^{n+1}}h \\
    &= h(Q^\extension + \DiracSt_{\R^{n+1}})hT_\sigma - (-1)^{\deg(\sigma)} T_\sigma h\DiracSt_{\R^{n+1}} h \\
    &= hQ^\extension hT_\sigma + h\DiracSt_{\R^{n+1}}(\sigma)\cdot h + (-1)^{\deg(\sigma)} hT_\sigma \DiracSt_{\R^{n+1}}h - (-1)^{\deg(\sigma)}T_\sigma h\DiracSt_{\R^{n+1}} h \\
    &= hQ^\extension h T_\sigma + h\DiracSt_{\R^{n+1}}(\sigma)\cdot h.
 \end{align*}
 Since $\sigma$ is compactly supported, Fubini's theorem and the triangle inequality yield
 \begin{align*}
     || h \susp(Q) h T_\sigma - (-1)^{\deg(\sigma)} T_\sigma h\DiracSt_{\R^{n+1}} h||_{0,0} \leq \left(||Q^\extension\sigma||_0 + ||\DiracSt_{\R^{n+1}}\sigma||_0\right) \cdot ||h^2||_{\infty,\supp \, \sigma},
 \end{align*}
 which proves the first condition of Theorem \ref{Kucerovski - Theorem}.

 The second condition $\mathrm{dom}(h\susp(Q)h) \subseteq \mathrm{dom}(hQ^\extension h)$ is clearly satisfied because, by definition, $\dom(\susp(h\susp(Q)h)) = \dom(hQ^\extension h) \cap \dom(h\DiracSt_{\R^{n+1}}h)$.
 
 To prove the third condition, we need to show that there is a $c \in \R$ such that 
 \begin{equation}\label{eq: Kucerovski3}
     \langle hQ^\extension h u, h \susp(Q) h u\rangle + \langle h \susp(Q) h u, h Q^\extension h u \rangle \geq c ||u||_0^2
 \end{equation}
 for all $u \in \Gamma_c(M \times \R^{n+1}; \Spinor_{g_0 \oplus \euclmetric})$.
 Recall the identity
 \begin{align*}
     \DiracSt_{\R^{n+1}}Q^\extension + Q^\extension\DiracSt_{\R^{n+1}} &= \sum_{j=1}^{n+1} \partial_j \cdot (\partial_jQ)^\extension =: \DiracSt_{\R^{n+1}}(Q^\extension).
 \end{align*}
 from the proof of Lemma \ref{susp quad - Lemma} and that $h$ and $Q^\extension$ commute.
 The calculation
 \begin{align*}
     \langle h \DiracSt_{\R^{n+1}} h u, hQ^\extension h u\rangle &= \langle Q^\extension \DiracSt_{\R^{n+1}} h u, h^3u\rangle \\
     &= \langle -\DiracSt_{\R^{n+1}} Q^\extension h u, h^3 u\rangle + \langle \DiracSt_{\R^{n+1}}(Q^\extension) h u, h^3 u\rangle \\
     \begin{split}
       &= - \langle Q^\extension hu, \grad{h^2} h u\rangle - \langle h Q^\extension h u, h \DiracSt_{\R^{n+1}}h u\rangle \\
       & \qquad \qquad + \langle \DiracSt_{\R^{n+1}}(Q^\extension) h^2 u, h^2 u\rangle
     \end{split}
 \end{align*}
 implies 
 \begin{align*}
     &\quad \ \, \langle hQ^\extension h u, h \susp(Q) h u\rangle + \langle h \susp(Q) h u, h Q^\extension h u \rangle \\
     &= 2||hQ^\extension h u||_0^2 + \langle \grad{h^2}\cdot Q^\extension h u, h u\rangle + \langle \DiracSt_{\R^{n+1}}(Q^\extension) h^2 u, h^2 u\rangle. 
 \end{align*}
 Fubini's theorem gives
 \begin{align*}
   &\quad \ \, ||hQ^\extension h u||_0^2 + \langle \grad{h^2}\cdot Q^\extension h u, h u\rangle + \langle \DiracSt_{\R^{n+1}}(Q^\extension) h^2 u, h^2 u\rangle \\
   &= \int_{\R^{n+1}}h(t)^2\left( 2||Q^\extension u_t||_0^2 + \langle \grad{h^2}Q^\extension u_t, u_t\rangle + \langle \DiracSt_{\R^{n+1}}(Q^\extension) u_t, u_t\rangle \right) \diff t,
 \end{align*}
 where $u_t = u(\placeholder, t) \in \Gamma(M,\Spinor_{g_0} \otimes \Cliff[n+1][0])$ is the restriction to the slice $M = M\times \{t\}$.
 
 We will show the required inequality (\ref{eq: Kucerovski3}) on each slice.
 The difference $B - \Dirac_{g_0}$ is a block map of pseudo differential operators of order zero because $Q \colon \R^{n+1} \rightarrow \InvPseudDir_{g_0}$ is a block map of pseudo Dirac operators over a fixed metric $g_0$.
Hence, $\DiracSt_{\R^{n+1}}(Q^\extension) = \DiracSt_{\R^{n+1}}(B^\extension)$
restricts to a bounded operator on $\Spinor_{g_0} \otimes \Cliff[n+1][0] \rightarrow M \times \{t\}$
for all $t \in \R^{n+1}$.
Their operator norms are uniformly bounded from above.
We conclude
\begin{equation*}
    \langle \DiracSt_{\R^{n+1}}(Q^\extension) u_t, u_t\rangle \geq c_1 ||u_t||_0^2
\end{equation*}
for all $c_1 \geq - \mathrm{sup} \{||\DiracSt_{\R^{n+1}}(B^\extension)\restrict_{M\times \{t\}}||_{0,0} \, : \, t \in \R^{n+1}\} > -\infty$.

On the other hand,  $Q^\extension_t \deff Q^\extension \restrict_{M \times\{t\}} = Q_t \otimes \id_{\Cliff[n+1][0]}$ is an elliptic operator of order $1$ over the closed manifold $M$, so the Hilbert space $L^2(M;\Spinor_{g_0} \otimes \Cliff[n+1][0])$ decomposes into eigenspaces $\mathrm{Eig}(Q_t^\extension,\lambda_i(t))$, where $(\lambda_i(t))_{i\in \Z}$ is the sequence of eigenvalues with $\lambda_i \leq \lambda_{i+1}$.
This sequence of eigenvalues is unbounded in both directions, see \cite{LawsonMichelsonSpin}*{Theorem II.5.8}.

Denote the orthogonal projection onto $\mathrm{Eig}(Q_t^\extension,\lambda_i(t))$ by $\Pi_{\lambda_i(t)}$ and the spectrum of $Q_t^\extension$ by $\Sigma_t$.
The operators $Q^\extension$ and $\grad{h^2} \cdot $ anti-commute.
Thus, if $\grad{h^2}(t) \neq 0$, it restricts to an isomorphism
\begin{equation*}
    \grad{h^2}\cdot \colon \mathrm{Eig}(Q^\extension_t, \lambda_i(t)) \rightarrow \mathrm{Eig}(Q^\extension_t,-\lambda_i(t)).
\end{equation*}
This implies
\begin{align*}
    &\quad \ \, 2h^2||Q^\extension u_t||_0^2 + \langle \grad{h^2} \cdot Q^\extension u_t, u_t\rangle \\
    &= \sum_{\lambda,\mu \in \Sigma_t} 2h^2 \langle (Q^\extension)^2 \Pi_{\lambda} u_t, \Pi_{\mu}u_t \rangle + \langle \grad{h^2}\cdot Q^\extension \Pi_{\lambda}u_t, \Pi_{\mu} u_t\rangle \\
    &= \sum_{\lambda,\mu \in \Sigma_t} 2h^2 \lambda^2(t) \langle  \Pi_{\lambda} u_t, \Pi_{\mu}u_t \rangle + \lambda\langle \grad{h^2}\cdot  \Pi_{\lambda}u_t, \Pi_{\mu} u_t\rangle \\
    &= \sum_{\lambda \in \Sigma_t} 2h^2 \lambda^2(t) \langle  \Pi_{\lambda} u_t, \Pi_{\lambda}u_t \rangle + \lambda\langle \grad{h^2}\cdot  \Pi_{\lambda}u_t, \Pi_{-\lambda} u_t\rangle.
\end{align*}
The spectrum of $Q_t^\extension$ is symmetric and $h \geq 1$, so we can estimate
\begin{align*}
     &\quad \ \, 2h^2||Q^\extension u_t||_0^2 + \langle \grad{h^2} \cdot Q^\extension u_t, u_t\rangle \\
     &\geq \sum_{\lambda \in \Sigma_t} \lambda^2 ||\Pi_{\lambda}u_t||_0^2 - ||\grad{h^2}||_\infty \cdot |\lambda| \cdot ||\Pi_{\lambda}u_t||_0 ||\Pi_{-\lambda}u_t||_0 \\
     &= \sum_{\lambda \in \Sigma_t \cap \R_{\geq 0}} \lambda^2\left(||\Pi_{\lambda}u_t||_0^2 + ||\Pi_{-\lambda}u_t||^2_0\right) - 2||\grad{h^2}||_\infty \cdot |\lambda|\cdot ||\Pi_{\lambda}u_t||_0 ||\Pi_{-\lambda}u_t||_0 \\
     &\geq \sum_{\lambda \in \Sigma_t \cap \R_{\geq 0}} \left(\lambda^2 - ||\grad{h^2}||_\infty\cdot |\lambda|\right)\left(||\Pi_{\lambda}u_t||_0^2 + ||\Pi_{-\lambda}u_t||_0^2\right).
\end{align*}
Since $ ||\grad{h^2}||_\infty \leq B$ for some $B \in \R$, since $Q \colon \R^{n+1} \rightarrow \InvPseudDir_{g_0}$ is a block map, and since eigenvalues of self-adjoint elliptic operators depend continuously on the elliptic operator \cite{booss1993elliptic}*{Section 17A}, the set 
\begin{equation*}
    \{i \in \Z \, : \, \text{there is a }t\in \R^{n+1} \text{ such that } \lambda_i^2(t) < ||\grad{h^2}|| \cdot |\lambda_i(t)|\}
\end{equation*}
is finite.
Thus, there is a $c_2 \in \R$ such that 
\begin{equation*}
    \lambda^2(t) - ||\grad{h^2}||_\infty \cdot |\lambda(t)| \geq c_2
\end{equation*}
for all $\lambda(t) \in \Sigma_t$ and all $t \in \R^{n+1}$.
We conclude
\begin{equation*}
    h^2 ||Q^\extension u_t||^2_0 + \langle \grad{h^2}Q^\extension u_t, u_t\rangle \geq c_2 ||u_t||_0^2.
\end{equation*}
Note that this inequality is trivially satisfied if $\grad{h^2}(t) = 0$.

Integrating these inequalities over $t$ gives
\begin{equation*}
    \langle hQ^\extension h u, h \susp(Q) h u\rangle + \langle h \susp(Q) h u, h Q^\extension h u \rangle \geq \mathrm{min}\{c_1,c_2\} ||u||_0^2,
\end{equation*}
so the third and final condition of Theorem \ref{Kucerovski - Theorem} is satisfied.
We apply it to deduce that 
\begin{equation*}
    \tau_{\Cliff[0][n+1]}(b([hQh])) \sharp \alpha_{n+1} = b([h\susp(Q)h])
\end{equation*}
and the theorem is proven.
\end{proof}
\begin{cor}\label{SuspensionInjection - Cor}
 The operator suspension $\susp_\bullet \colon \InvPseudDir[\bullet]_{g_0} \rightarrow \InvBlockDirac_{g_0}$ induces an injective map on all homotopy groups.
\end{cor}
\begin{proof}
 Recall the split diagram
 \begin{equation*}
    \xymatrix{ \pi_n(\InvPseudDir[\bullet]_{g_0}) \ar[r] \ar[d]_{\pi_n(\susp_\bullet)} & KK(Cl_{0,d}, \mathcal{C}_0(\R^{n+1})) \ar[d]^{\tau_{Cl_{0,n+1}}(\placeholder)\sharp \alpha_{n+1}} \\ 
    \pi_n(\InvBlockDirac_{g_0}) \ar[r] & KK(Cl_{0,d+n+1},\R).}
\end{equation*}
 The upper horizontal map is a bijection by Lemma \ref{UpperArrowIso - Lemma} and the right vertical map is an isomorphism by Kasparov theory.
 The diagram commutes by the previous theorem, so the composition of the left vertical arrow and the lower horizontal arrow is a bijection.
 In particular, $\pi_n(\susp_\bullet)$ is split injective on the set of path components and on all homotopy groups. 
\end{proof}

\subsubsection*{The Determination of $\pi_n(\InvBlockDirac)$}

Recall diagram (\ref{eq: DiagramKKcomp}) from Corollary \ref{SuspensionInjection - Cor}.
We know that three out of four groups are isomorphic to $KO^{-(d+n+1)}(\mathrm{pt})$, so we are left with the following situation
\begin{equation*}
    \xymatrix{KO^{-(d+n+1)}(\mathrm{pt}) \ar[rd] \ar@{^{(}->}[d]&  \\ \pi_n(\InvBlockDirac_{g_0}) \ar[r] & KO^{-(d+n+1)}(\mathrm{pt}).}
\end{equation*}
The group $KO^{-(d+n+1)}(\mathrm{pt})$ is cyclic, so if we could show that there is a surjective map $KO^{-(d+n+1)}(\mathrm{pt}) \twoheadrightarrow \pi_n(\InvBlockDirac)$, then $\pi_n(\InvBlockDirac)$ would be cyclic, too, and an easy diagram chase would imply that $\pi_n(\susp_\bullet)$ is also surjective.

In conclusion, we are left to prove the next theorem.
\begin{theorem}\label{KInvBlockDirSurj - Theorem}
 If $M$ is a closed, $d$-dimensional, spin manifold with Riemannian metric $g_0$ such that $\Dirac_{g_0}$ is invertible, then there is a surjective map
 \begin{equation*}
     KO^{-(d+n+1)}(\mathrm{pt}) \twoheadrightarrow \pi_n(\InvBlockDirac_{g_0}, \Dirac_{g_0})  
 \end{equation*}
 for all $n \in \N_0$.
\end{theorem}

The proof idea is the same as the one of \cite{ebert2017indexdiff}*{Prop 4.3}.
We will compare the set $\Omega_{\Dirac_{g_0}}^n\InvBlockDirac_{g_0}$, suitably topologised, to $(B + \mathbf{Kom})^\times$, the space of invertible, compact pertubations of an invertible operator $B$ and show that these two spaces are weakly homotopy equivalent. 
The latter space is homotopy equivalent to $\Omega \mathrm{Fred}^{d+n}$, see \cite{ebert2017indexdiff}*{Proof of Prop 4.3 and Lemma 4.9}, which is a classifying space for $KO^{-(d+n+1)}(\mathrm{pt})$.

Note that the cited arguments work only for a \emph{fixed} (but arbitrary) Hilbert space with an ample $\Cliff[d+n+1][0]$-right action, $\mathrm{i.e.}$, a right action that contains infinitely many copies of all irreducible $\Cliff[d+n+1][0]$-representations.
The requirement to fix the Hilbert space is the reason why we work with the fibre $\InvBlockDirac_{g_0}$ instead of the total space $\InvBlockDirac$.

Let us begin with fixing notation.
Let $\mathcal{B}$ be the set of all compactly supported, bounded operators on $L^2(M \times \R^{n}, \Spinor_{g_0})$ and let $\mathbf{B}$ be subset of $\mathcal{B}$ consisting of all self-adjoint, odd, and Clifford linear operators.
We endow $\mathbf{B}$ with the operator norm topology so that it becomes a locally convex vector space.
Furthermore, we define $\mathbf{B}_{act} = \mathbf{B} \cap \PseudOp^0_c(\Spinor_{g_0})$ as the set of all compactly supported pseudo differential operators of order zero that are self-adjoint, odd, and Clifford-linear. 

Let further $\mathcal{K}$ be the set of all compact operators on $L^2(M \times \R^{n}, \Spinor_{g_0})$ and let 
\begin{equation*}
    \mathbf{Kom} \deff \{A \in \mathcal{K} \, : \, A \text{ is s.a., odd, Clifford-linear}\}.
\end{equation*}
Finally, we denote the spinor connection on $\Spinor_{g_0}$ with $\nabla$ and define the operator $T \deff (1+\nabla^\ast \nabla)^{-1/4}$.

Recall that we have defined the Sobolev spaces $H^s$ as the closure of all compactly supported sections $\Gamma_c(M \times \R^{n},\Spinor_{g_0})$ under the norm
\begin{equation*}
    ||u||_s^2 \deff \int_{M \times \R^{n}} \langle (1+\nabla^\ast \nabla)^s u, u\rangle_{\Spinor_{g_0}} \diff \vol(g_0)
\end{equation*}
so that $T \colon H^s \rightarrow H^{s+1/2}$ is an isometric isomorphism.

\begin{lemma}\label{BoundedTtrafoIsDense - Lemma}
 The assignment $B \mapsto TBT$ gives rise to well-defined continuous maps
 \begin{equation*}
     T(\placeholder)T \colon \mathcal{B} \rightarrow \mathcal{K} \quad \text{and} \quad T(\placeholder)T \colon \mathbf{B}_{(act)} \rightarrow \mathbf{Kom}
 \end{equation*}
 whose images are dense. 
\end{lemma}
\begin{proof}
 The map $T(\placeholder)T$ is continuous because $T$ is a bounded operator.
 
 The spinor connection $\nabla$ is even and Clifford-right linear, so the operator $TBT$ is self-adjoint, odd, and Clifford-linear if $B$ is self-adjoint, odd, and Clifford-right linear.
 Thus, it remains to show that $TBT$ is compact.
 
 Since $B$ is supported within $M \times RI^n$ for some $R > 0$, every compactly supported function with $\chi \equiv 1$ on $M \times RI^n$ satisfies $B= B \chi \cdot $.
 Thus, the operator $TBT \colon H^0 \rightarrow H^0$ factors as follows
 \begin{equation*}
     \xymatrix{H^0 \ar[r]^-T & H^{1/2} \ar[r]^-{\chi\cdot} & H^{1/2}_{\supp(\chi)} \ar@{^{(}->}[r]& H^0_{\supp(\chi)} \ar[r]^-B & H^0 \ar[r]^-T & H^{1/2} \ar@{^{(}->}[r] & H^0.}
 \end{equation*}
 The inclusion $H^{1/2}_{\supp(\chi)} \hookrightarrow H^0_{\supp(\chi)}$ is compact by Rellich's lemma, so $TBT$ is compact.
 
 We first prove that $T(\placeholder)T \colon \mathcal{B} \rightarrow \mathcal{K}$ has a dense image.
 Since $T \colon H^s \rightarrow H^{s+1/2}$ is an isometry for all $s \in \R$, it induces an homeomorphism $T \colon H^\infty \rightarrow H^\infty$, where $H^\infty = \bigcap_{s\in \R}H^s$ is endowed with the inverse limit topology.
 Thus, $T$ maps dense subsets $W \subseteq H^\infty$ to dense subsets.
 
 For each dense subset $W \subseteq H^\infty$, the linear hull of the set of rank 1 projections 
 $\{\langle \placeholder , w_1\rangle w_2 \, : \, w_j \in W \}$ arising from $W$ is a dense subset of $\mathcal{K}$.
 The identity 
 \begin{equation*}
     T \circ \langle \placeholder ,w_1\rangle w_2 \circ T = \langle T(\placeholder),w_1\rangle T(w_2) = \langle \placeholder, T(w_1) \rangle T(w_2)
 \end{equation*}
 implies that $T(\placeholder)T$ maps bijectively rank 1 projections arising from $W$ to rank $1$ projections arising from $T(W)$.
 The tensor product $\Gamma(M,\Spinor_{g_0}) \otimes \Gamma_c(\R^n, \Cliff[n][0])$ is a dense subspace of $\Gamma_c(M \times \R^n,\Spinor_{g_0\oplus \euclmetric})$, which is itself a dense subset of $H^\infty$.
 Clearly,
 \begin{equation*}
     \{\langle \placeholder, \sigma_1 \otimes \varphi_1\rangle \sigma_2\otimes \varphi_2 \, : \, \sigma_j \in \Gamma(M,\Spinor_{g_0}), \varphi \in \Gamma_c(\R^n,\Cliff[n][0])\} \subseteq \PseudOp_c^0(\Spinor_{g_0}) \subseteq \mathcal{B}
 \end{equation*}
 because the sections are compactly supported.
 It follows that $T\PseudOp_c^0(\Spinor_{g_0})T$ and $T\mathcal{B}T$ are dense in $\mathcal{K}$.
 
 In order to show that the restriction $T(\placeholder)T \colon \mathbf{B} \rightarrow \mathbf{Kom}$ has dense image, we make use of the projection (4.7) of \cite{ebert2017indexdiff}.
 If $\Gamma$ denotes the multiplicative group generated by the standard basis $e_j \in \R^{d+n} \subseteq \Cliff[d+n][0]$, then the projection $\Pi_\mathcal{K} \colon \mathcal{K} \rightarrow \mathbf{Kom}$ is defined by
 \begin{align*}
     K &\mapsto \frac{1}{4|\Gamma|} \left\lbrack\left(\sum_{\gamma \in \Gamma} \gamma K \gamma^{-1} - \evenodd \gamma K \gamma^{-1} \evenodd\right) - \left(\sum_{\gamma \in \Gamma} \gamma K \gamma^{-1} - \evenodd \gamma K \gamma^{-1} \evenodd\right)^\ast\right\rbrack.
 \end{align*}
 If we can show that this assignment gives rise to a projection $\Pi_\mathcal{B} \colon \mathcal{B} \rightarrow \mathbf{B}$, then the commutative diagram
 \begin{equation*}
     \xymatrix{\PseudOp_c^0(\Spinor_{g_0}) \ar@{^{(}->}[r] \ar@{->>}[d]_{\Pi_{\mathcal{B}}} & \mathcal{B} \ar[rr]^{T(\placeholder)T} \ar@{->>}[d]_{\Pi_{\mathcal{B}}} && \mathcal{K} \ar@{->>}[d]^{\Pi_\mathcal{K}} \\ \mathbf{B}_{act} \ar@{^{(}->}[r] & \mathbf{B} \ar[rr]^{T(\placeholder)T} && \mathbf{Kom}}
 \end{equation*}
 implies that the restriction of $T(\placeholder)T$ to $\mathbf{B}$ and $\mathbf{B}_{act}$ has dense image, too.
 
 Note that if $A$ and $B$ are supported within $M \times RI^n$ and $M \times SI^n$, respectively, then $A+B$ is supported within $M \times (R+S)I^n$.
 If $f$ is a vector bundle map over the identity, then $fA$ and $Af$ are still bounded operators and supported with in $M \times RI^{n}$.
 It is easily veryfied that if $A$ is supported within $M\times RI^n$, then also $A^\ast$ is supported within $M \times RI^n$.
 
 In conclusion, $\prod_{\mathcal{B}} \colon \mathcal{B} \rightarrow \mathbf{B}$ is well-defined, which finishes the proof.
\end{proof}

Since we are working over non-compact manifolds the next result does not follow from general theory.

\begin{lemma}\label{TDTBnd - Lemma}
  The operator $T \Dirac_{g_0\oplus \euclmetric} T \colon H^0 \rightarrow H^0$ is bounded.
\end{lemma}
The proof relies on the following abstract operator theoretic result.
\begin{lemma}\label{UnboundedMurphy - Lemma}
 Let $P$ be a possibly unbounded, self-adjoint, non-negative operator and $B$ be a bounded self-adjoint operator with $B \geq 1$.
 Then we have
 \begin{equation*}
     (B + P)^{-1/4} \leq (1 + P)^{-1/4}. 
 \end{equation*}
\end{lemma}
 \begin{proof}
  We clearly have $B + P \geq 1 + P$, so both operators are invertible.
  Since the inverse of an unbounded operator is bounded, we can apply the proof of \cite{murphy2014c}*{Theorem 2.2.5 (4)}, to deduce $(1+P)^{-1} \geq (B + P)^{-1}$. Applying \cite{murphy2014c}*{Theorem 2.2.6} twice, we deduce $(B + P)^{-1/4} \leq (1 + P)^{-1/4}$. 
 \end{proof}
\begin{proof}[Proof of Lemma \ref{TDTBnd - Lemma}]
 For $c\geq 1$, define $T_c \deff (c + \nabla^\ast \nabla)^{-1/4}$. 
 Because multiplication with $c$ and $\nabla^\ast\nabla$ commute, unbounded functional calculus implies that $T_c \circ T^{-1}$ is bounded with operator norm $||T_c \circ T^{-1}||_{0,0} \leq c$.
 Thus, $T \Dirac_{g_0\oplus \euclmetric} T$ is bounded if and only if $T_c \Dirac_{g_0\oplus \euclmetric} T_c$ is bounded for some $c\geq 1$.
 
 Since $M$ is closed, there is a $c>0$ such that $c - \scal(g_0)/4 > 1$.
 Using the Lichnerowicz formula and Lemma \ref{UnboundedMurphy - Lemma} we deduce
 \begin{align*}
     T_c\Dirac_{g_0\oplus \euclmetric}T_c &=(c - \scal(g_0)/4 + \Dirac_{g_0\oplus \euclmetric}^2)^{-1/4}\Dirac_{g_0\oplus \euclmetric}(c - \scal(g_0)/4 + \Dirac_{g_0\oplus \euclmetric}^2)^{-1/4} \\
     &\leq (1 + \Dirac_{g_0\oplus \euclmetric}^2)^{-1/4}\Dirac_{g_0\oplus \euclmetric}(1 + \Dirac_{g_0\oplus \euclmetric}^2)^{-1/4}\\
     &=\Dirac_{g_0\oplus \euclmetric}(1+\Dirac_{g_0\oplus \euclmetric}^2)^{-1/2},
 \end{align*}
 so $T_c\Dirac_{g_0\oplus \euclmetric}T_c$ is bounded.
\end{proof}
\begin{proof}[Proof of Theorem \ref{KInvBlockDirSurj - Theorem}]
 Set $B \deff T \Dirac_{g_0\oplus \euclmetric} T$.
 Recall from Corollary \ref{SymbolExtensionAs - Corollary} that ${\symb}^1_K \colon \PseudOpCl^1_K(\Spinor_{g_0}) \rightarrow \Symb^0_K(\Spinor_{g_0})$ is continuous for all compact sets $K \subseteq M \times \R^n$ and that its kernel is contained in $\mathcal{K}(H^1,H^0)$ by Lemma \ref{ASExactSequence - Lemma}.
 Thus 
 \begin{equation*}
    \mathfrak{B} \deff \{P - \Dirac_{g_0} \, : \, P \in \Omega_{g_0}^n \BlockDirac_{g_0}\} 
 \end{equation*}
 injects into $\mathcal{K}(H^1,H^0)$.
 In fact, it injects even into the smaller space $\mathbf{K}_c(H^1,H^0)$ of all compact operators from $H^1(\Spinor_{g_0})$ to $H^0(\Spinor_{g_0}) = L^2(\Spinor_{g_0})$ that are odd, self-adjoint (considered as unbounded operators) and Clifford linear. 
 Thus, we equip $\mathfrak{B}$ with the subspace topology of the operator norm topology of $\mathcal{K}(H^1,H^0)$. 
 We equip the spaces in the lower row with the topologies such that the vertical arrows in the following diagram become homeomorphisms:
 \begin{equation*}
     \xymatrix{ \mathfrak{B} \ar[d]^\cong_{\Dirac_{g_0\oplus \euclmetric}+} & \ar[l] \mathbf{B}_{act}  \ar[d]^\cong \ar[r]^{T(\placeholder)T} \ar[d]^\cong & \mathbf{Kom} \ar[d]^{B + }_{\cong}\\
     \Omega_{\Dirac_{g_0}}^n\BlockDirac_{g_0}& \ar[l] \Dirac_{g_0} + \mathbf{B}_{act} \ar[r]^{T(\placeholder)T} &  B + \mathbf{Kom}.  }
 \end{equation*}
 The diagram restricts to the subspaces of invertible operators: 
 \begin{equation*}
     \xymatrix{ U_0  \ar[d]^\cong_{\Dirac_{g_0\oplus \euclmetric}+} & \ar[l] U_{1/2} \ar[d]^\cong \ar[r]^{T(\placeholder)T} & U_1 \ar[d]^{B + }_{\cong}\\
     \Omega_{\Dirac_{g_0}}^n\InvBlockDirac_{g_0}& \ar[l] (\Dirac_{g_0} + \mathbf{B}_{act})^\times \ar[r]^{T(\placeholder)T} &  (B + \mathbf{Kom})^\times.}
 \end{equation*}
 Note that $U_{1/2}$ is the preimage of $U_1$ under $T(\placeholder)T$ in $\mathbf{B}_{act}$.
 Furthermore $\mathbf{B}_{act}$ is a locally convex topological vector space and its inclusions into $\mathfrak{B}$ is continuous and has dense image (because infinite smoothing operators are dense in the compact operators as explained in the proof of Lemma \ref{BoundedTtrafoIsDense - Lemma}). 
 By Palais' theorem \cite{palais1966homotopy}*{Corollary of Theorem 12}, the two maps, $T(\placeholder)T \colon U_{1/2} \rightarrow U_1$ and the inclusion $U_{1/2} \hookrightarrow U_0$, are weak homotopy equivalences.
 This implies, in particular, that we have a zig-zag of bijections 
 \begin{equation*}
     \pi_0(\Omega^n_{\Dirac_{g_0}}\InvBlockDirac_{g_0}) \xleftarrow{\cong} \pi_0((\Dirac_{g_0} + \mathbf{B}_{act})^\times)\xrightarrow{\cong} \pi_0((B + \mathbf{Kom})^\times).
 \end{equation*}

 We will prove in Lemma \ref{FactorisationIsGrpHom - Lemma} below that the surjection $\pi_0(\Omega^n_{g_0}\InvBlockDirac_{g_0}) \twoheadrightarrow \pi_n(\InvBlockDirac_{g_0})$ is a group homomorphism so that we have the following zig-zag
 \begin{equation*}
     \pi_n(\InvBlockDirac_{g_0}) \twoheadleftarrow \pi_0(\Omega^n_{\Dirac_{g_0}}\InvBlockDirac_{g_0}) \xrightarrow{\cong} \pi_0((B + \mathbf{Kom})^\times).
 \end{equation*}
 Combined with Ebert's result that $(B + \mathbf{Kom})^\times$ is homotopy equivalent to $\Omega\mathrm{Fred}^{d+n}$, see \cite{ebert2017indexdiff}*{Lemma 4.9 or Proof of Proposition 4.3}, the zig-zag produces, by inverting the isomorphism, the following surjection:
 \begin{equation*}
     KO^{-(d+n+1)}(\mathrm{pt}) \twoheadrightarrow \pi_n(\InvBlockDirac_{g_0}).\qedhere
 \end{equation*}
\end{proof}

\begin{lemma}\label{FactorisationIsGrpHom - Lemma}
 The map $\pi_n(\susp_\bullet) \colon \pi_n(\InvPseudDir[\bullet]_{g_0}) \rightarrow \pi_n(\InvBlockDirac_{g_0})$ factors through $\pi_0(\Omega^n_{\Dirac_{g_0}}\InvBlockDirac_{g_0})$.
\end{lemma}
\begin{proof}
 Clearly, the suspension also induces a map
 \begin{equation*}
     \xymatrix{\pi_n(\InvPseudDir[\bullet]_{g_0}) && \ar[ll]_-{\cong}  \pi_n( B_{\bullet,g_0},\Dirac_{g_0}) \ar[rr]^-{\pi_n(\susp_\bullet)} && \pi_0(\Omega_{\Dirac_{g_0}}^n\InvBlockDirac_{g_0})}
 \end{equation*}
 because a combinatorial homotopy in $B_{n+1,g_0}$ between $P_{-1}$ and $P_1$ suspends to a path between $\susp(P_{-1})$ and $\susp(P_1)$ by suspending along the last $n$ coordinates.
 
 It remains to show that the map $\pi_0(\Omega^n_{\Dirac_{g_0}}\InvBlockDirac_{g_0}) \rightarrow \pi_n(\InvBlockDirac)$ that sends a representative of a path component to its concordance is well-defined.
 Denote by $\mathbf{B}_{act}$ the locally convex vector space of all compactly supported pseudo differential operators of order zero on $\Spinor_{g_0} \rightarrow M \times \R^n$ that are self-adjoint, odd and Clifford linear equipped with the operator norm topology.
 The canonical map $\iota \colon \mathbf{B}_{act} \rightarrow \mathfrak{B}$ is a continuous, injective, linear map with dense image, hence, by Palais' result \cite{palais1966homotopy}*{Corollary of Theorem 12}, it restricts to a weak equivalence $\iota^{-1}(U_0) \rightarrow U_0$.
 
 If $P_{-1}, P_1 \in \Omega^n_{\Dirac_{g_0}}\InvBlockDirac$ represent the same path component, then slight perturbations $\tilde{P}_j$  of $P_j$ such that $\tilde{P}_j - \Dirac_{g_0} \in \iota^{-1}(U_0)$ still represent the same path component.
 We find therefore a sufficiently slowly parametrised path $Q \colon \R \rightarrow \Dirac_{g_0} + \iota^{-1}(U_0)$ that connects $\tilde{P}_{-1}$ with $\tilde{P}_1$ and such that $\susp(Q)$ is invertible.
 
 However, $\susp(Q)$ might not yet be a block operator because the amplitude topology does not control the diameter of a core.
 In order to fix this problem, consider the smooth block map
 \begin{equation*}
     t \mapsto B_t \deff Q_t - \Dirac_{g_0}
 \end{equation*}
 of compactly supported pseudo differential operators of order zero (hence bounded operators).
 
 For the sake of readability, we will drop $\iota$ from the notation.
 For each given $\eps > 0$, we find a sufficiently large $N \in \N$, such that if $t \in [k/N,(k+1)/N]$, then $||B_t - B_{k/N}||_{0,0} < \eps/3$ for all $k \in \Z$.
 
 Assume that $B_{k/N}$ is supported within $M \times R_kI^n$.
 As a block map, $B$ is locally constant near infinity, so we may (and do) assume the sequence $(R_k)_{k \in \Z}$ to be bounded.
 Let $R$ be the supremum of this sequence.
 Fix a function $\chi \colon \R^n \rightarrow [0,1]$ that is identically $1$ on $RI^n$ and identically zero on the complement of $2RI^n$.
 The path $t \mapsto \chi B \chi$ then satisfies
 \begin{align*}
     ||B_t - \chi B_t \chi||_{0,0} & \leq ||B_t - B_{k/N}||_{0,0} + ||B_{k/N} - \chi B_t \chi||_{0,0} \\
                                   & = ||B_t - B_{k/N}||_{0,0} + ||\chi(B_{k/N} -  B_t) \chi||_{0,0} \\
                                   &\leq  \eps/3 + \eps/3 < \eps. 
 \end{align*}
 Since invertibility is an open condition, we may choose $\varepsilon>0$ sufficiently small such that $t \mapsto \Dirac_{g_0} + \chi B_t\chi$ is a block map of invertible operators that connect 
 $\tilde{P}_{-1}$ and $\tilde{P}_1$.
 
 However, we cannot form the suspension yet, because we do not know whether the path $t \mapsto \chi B_t \chi$ is continuous with respect to the \emph{amplitude topology}. 
 Luckily, the path only takes values in operators that are supported within $M \times 2RI^n$, in other words, we have a continuous path $\R \rightarrow \PseudOpCl^0_{M\times 2RI^n}(\Spinor_{g_0})$ (the target is equipped with the operator norm topology and has is locally constant with values $\tilde{P}_{\pm1}$).
 Since $\PseudOp^0_{M\times 2RI^n}(\Spinor_{g_0}) \hookrightarrow \PseudOpCl^0_{M\times 2RI^n}(\Spinor_{g_0})$ is continuous and has dense image,
 we can now use Palais theorem (and smoothing theory) to replace $t \mapsto \chi B_t \chi$ by a smooth path $ \hat{B}_t \colon \R \rightarrow \PseudOp_{M \times 2RI^n}(\Spinor_{g_0})$ of operators such that $\Dirac_{g_0} + \hat{B_t}$ is invertible, that connects $\tilde{P}_{- 1}$ and $\tilde{P}_1$, and whose suspension is an invertible block Dirac operator.

 Different interpolations of $P_{-1}$ and $P_1$ can be connected by convex-combination, which is a suspensionable path.
 Hence, sending an element of a path component to its concordance class gives indeed a well-defined map.
 
 To see that the composition 
 \begin{equation*}
     \xymatrix{\pi_n(B_\bullet,\Dirac_{g_0}) \ar[rr]^-{\pi_n(\susp_\bullet)} && \pi_0(\Omega_{\Dirac_{g_0}}^n\InvBlockDirac_{g_0}) \ar@{->>}[rr] && \pi_n(\InvBlockDirac_{g_0})}
 \end{equation*}
 agrees with $\pi_n(\susp_\bullet) \colon \pi_n(B_{\bullet,g_0},\Dirac_{g_0}) \rightarrow \pi_n(\InvBlockDirac_{g_0})$, just observe that the suspension of the convex combination between $\susp(P)$ and its slight perturbation $\widetilde{\susp(P)}$ give rise to a concordance connecting these two elements in $\pi_n(\InvBlockDirac_{g_0})$.
\end{proof}

An immediate consequence is the proof of Theorem \ref{OperatorSuspWeakEquiv - Theorem} that the operator suspension map is a weak homotopy equivalence.
\begin{proof}[Proof of Theorem \ref{OperatorSuspWeakEquiv - Theorem}]
 As discussed at the beginning of this chapter, it suffices to prove this statement for the restriction to the fibre $B_{\bullet,g_0} \hookrightarrow \InvBlockDirac_{g_0}$.
 By Corollary \ref{SuspensionInjection - Cor} we already know that the suspension induces an injective map on all homotopy groups.
 If $n\geq 1$, then Theorem \ref{KInvBlockDirSurj - Theorem} and Diagram (\ref{eq: DiagramKKcomp}) yield a sequence of surjective maps
 \begin{equation*}
     \xymatrix@C-0.3em{KO^{-(d+n+1)}(\mathrm{pt}) \ar@{->>}[r] & \pi_n(\InvBlockDirac_{g_0}) \ar@{->>}[r] & KK(\Cliff[0][d+n+1],\R) \ar@{=}[r] & KO^{-(d+n+1)}(\mathrm{pt}). }
 \end{equation*}
 Since surjective endomorphisms between cyclic groups are bijective, three out of four maps in Diagram (\ref{eq: DiagramKKcomp}) are isomorphisms, so the remaining homomorphism $\pi_n(\susp_\bullet)$ must be an isomorphism, too.
 
 If $n=0$, we use that $\PseudOp^1(\Spinor_{g_0}) \hookrightarrow \PseudOpCl^1(\Spinor_{g_0})$ is dense to deduce that $\Omega^0_{\Dirac_{g_0}}\InvPseudDir[\bullet]_{g_0} \hookrightarrow \Omega^0_{\Dirac_{g_0}}\InvBlockDirac_{g_0}$ is a weak homotopy equivalence.
 Surjectivity of $\pi_0(\susp_\bullet)$ follows now from Lemma \ref{FactorisationIsGrpHom - Lemma}.
\end{proof}



\chapter{The Factorisation of the Index Difference}\label{Factorisation Indexdiff - Chapter}

We have put much effort into the construction of $\InvBlockDirac$ and into the proof that it is a classifying space for real $K$-theory.
It remains to prove Theorem \ref{Factorisation - MainThm}, namely that the index difference factors through the concordance space.
We will give a precise formulation of this Theorem together with its proof in Section \ref{Section - The Factorisation Theorem}. 
Applications regarding the non-triviality of $\ConcSet$ for special $M$ using the existing literature, in particular Theorem \ref{ConterExampleStrongConcVsIsotopy - MainThm}, will be given in Section \ref{Section - Applications}.

\section{The Factorisation Theorem}\label{Section - The Factorisation Theorem}

\begin{theorem}[The Factorisation Theorem]\label{IndDifFactor - Theorem}
 Let $\AuxInvDir$ be the cubical subset from Definition \ref{AuxInvDir - Definition}.
 There is a weakly homotopy equivalent cubical subset $\AuxSingMet \subseteq \SingMet$ such that the following diagram commutes up to homotopy:
 \begin{equation*}
  \xymatrix{ \SingMet  \ar@{^{(}->}[d]^{\Dirac_\bullet} & \ar@{_{(}->}[l] \AuxSingMet \ar@{^{(}->}[rr]^{\susp} \ar@{^{(}->}[d]^{\Dirac_\bullet} & 
  & \ConcSet \ar@{^{(}->}[d]^{\Dirac_\bullet} 
  \\ \InvPseudDir[\bullet] & \ar@{_{(}->}[l] \AuxInvDir \ar@{^{(}->}[rr]^\susp  &
  & \InvBlockDirac.  } 
 \end{equation*}  
\end{theorem}

The two maps $\Dirac_\bullet \circ \susp$ and $\susp \circ \Dirac_\bullet$ do not agree. 
We will establish this theorem by finding suitable conditions on the block maps of psc metrics so that the operator norm of the \emph{error term} $\susp \circ \Dirac_\bullet - \Dirac_\bullet \circ \susp$ becomes sufficiently small.
Since invertibility is an open condition, convex interpolation serves then as a homotopy between the two maps. 

\begin{definition}\label{AuxSingMet - Definition}
 Let $\AuxSingMet$ be the following cubical subset of $\SingMet$:
 A block map $g \in \SingMet[M][n]$ is contained in $\AuxSingMet[n]$ if and only if it satisfies the following conditions:
 \begin{gather}
  \sum_{j=1}^n |\trace( (\partial_j g^\op)^2)| + |\trace(\partial^2_j g)| +| \trace(\partial_j g)|^2 < \frac{1}{8} \scal(g(t)), \label{Eq: AuxSing1}\\
  \sum_{j=1}^n ||\partial_j g^\op||_{\op, g(t)} < \min \{d^{-1}, 2^{-(d+4)} \lowBnd(t)\}, \label{Eq: AuxSing2}\\
  \sum_{j=1}^n || [\Dirac_{g(t)}^\extension, \nabla_{\partial_j}^{\Spinor_{\susp g}}] ||_{1,0;M \times t} < 1/32 \cdot \lowBnd(t)^2, \label{Eq: AuxSing3}
 \end{gather}
 where $\lowBnd(t) \deff \inf \{||\Dirac_{g(t)}v||_{0,g(t)} \, : \, ||v||_{1,g(t)} = 1\} = \lowBnd_{\Dirac_g}(t)$.
\end{definition}
\begin{lemma}
 The sequence of sets form indeed a cubical subset of $\SingMet$.
\end{lemma}
\begin{proof}
 The proof is similar to the proof that $\AuxInvDir$ is a cubical subset.
 Note that the three defining conditions are local and stable, so this lemma follows from Lemma \ref{local stable criterium - Lemma}.
\end{proof}

Next, we wish to show that the inclusion is a weak homotopy equivalence. We employ the same strategy as before: We first show that $\AuxSingMet$ is a Kan set and then use the usual rescaling trick on elements of the combinatorial homotopy groups.
\begin{proposition}
 The cubical subset $\AuxSingMet$ is a Kan set.
\end{proposition}
\begin{proof}[Proof of Theorem \ref{IndDifFactor - Theorem}]
 The proof is similar to the proof that $\AuxInvDir$ is Kan. 
 In fact, we have written that proof in a manner such that it can be easily complemented to prove the this proposition.
 We are going to employ the same construction, but, in contrast to the proof of Proposition \ref{AuxPseduDir Is Kan - Prop}, we need to choose the constant $R_1>1$ more carefully.
 
 We first recall the construction.
 For a given cubical $n$-horn
 \begin{equation*}
  \CubeBox = \{g_{(j,\omega)} \in \AuxSingMet[n-1] \, : \, \face[\omega]{j}g_{(k,\eta)} = \face[\eta]{k-1}g_{(j,\omega)}; \, j<k; (j,\omega), (k,\eta) \neq (i,\eps)\},
 \end{equation*}
 there is a sufficient large $\rho>0$ such that all block maps $g_{(j,\omega)}$ decompose appropriately outside of $\rho I^{n-1}$.
 We restrict $\degen{j}g_{(j,\omega)}$ to a map $\{\omega x_j > \rho \} \rightarrow \Riem^+(M)$.
 By the compatibility requirement, these restriction agree on the intersection of there domains and we can therefore glue them together to a map
 \begin{equation*}
   g_0 \colon \left(\quader[n][1,1][\rho]\right)^c\rightarrow \Riem^+(M),
 \end{equation*}
 where $\quader[n][i,\eps][\rho] \deff \{x \in \R^n \, : \, x_j \in [-\rho,\rho] \text{ if } j \neq i,\,  \eps x_i \geq -\rho\}$ is the cuboid with diameter $\rho$. 
 Let $K \subseteq \R^{n-1}$ be a compact, convex subset that contains $\rho I^{n-1}$, is point symmetric at the origin, and has a smooth boundary. 
 Let $Q_K \deff \{x \in \R^n \, : \, \CubeProj{i}(x) \in K, \eps x_i \geq -\rho \}$ be the cuboid with base $K$. 
 Restricting $g_0$ yields 
 \begin{equation*}
  g_0 \colon Q_K^c \cap \{\eps x_i < \rho + R_1\} \rightarrow \Riem^+(M),
 \end{equation*}
 for a fixed choice of $R_1>1$. 
 
 For the chosen $K$, there is a unique norm $||\placeholder ||_K$ whose unit ball is $K$. 
 Since $\partial K$ is smooth, the norm is smooth on $\R^{n-1}\setminus \{0\}$.
 For a sufficient large $R_2 = R_2(\rho ,R_1)$ to be determined later, pick $\chi \colon \R \rightarrow \lbrack 0,1\rbrack$ that is identically $1$ on $\R_{\leq 1.2}$, that vanishes on $\R_{\geq R_2-1}$, and whose first and second derivative satisfy $\chi' \leq 1.3/R_2$ and $|\chi''| \leq 10/R_2^2$.
 Furthermore, pick a monotonically increasing function $q \colon \R \rightarrow \R$ that is the identity on $\R_{\leq \rho}$, the constant map with value $\rho + R_1 - 1$ on $\R_{\rho + R_1}$, and whose derivatives satisfy $q' \leq 1$ and $|q''| \leq 1.2/R_1$.
 The map $q$ defines a map 
 \begin{align*}
  q &\colon \R^n \rightarrow \R^n \\
  (x_1, \dots, x_n) &\mapsto (x_1, \dots, x_{i-1}, q(x_i), x_{i+1}, \dots, x_n).
 \end{align*}
 Finally, we define
 \begin{align*}
  H \colon \left(Q_K^c  \cap \{\eps x_i < \rho + R_1\}\right) \times [0,1] &\rightarrow  Q_K^c \cap \{\eps x_i < \rho + R_1\}, \\
  (x,t) &\mapsto x - 2\eps(\rho +R_1)t e_i 
 \end{align*}
 and set 
 \begin{align*}
  g \colon \R^n &\rightarrow \Riem^+(M), \\
  x &\mapsto \begin{cases}
   g_0 \circ H\bigl(q(x), \chi(||\CubeProj{i}(x)||_K)\bigr), & \text{if } x \in Q_K^c, \\
   g_{(i,-\eps)}(\CubeProj{i}(x)), & \text{if } x \in Q_K.
  \end{cases}
 \end{align*}
 Note that $g$ is well defined and smooth because $g_0$ agrees with $g_{(i,-\eps)}\circ \CubeProj{i}$ on the set $\{-\eps x_i > \rho\}$ and $\partial Q_K$ lies in the preimage of $\{-\eps x_i > \rho\}$ under the map $H \circ (q(\placeholder ), \chi(||\CubeProj{i}(\placeholder )||_K))$. 
 
 The same argument as in the proof for Proposition \ref{AuxPseduDir Is Kan - Prop} shows that $g$ is a smooth block map and a filler of $\CubeBox$ in $\SingMet[M][n]$, from which we already know that there must exists one.
 It remains to show that $g$ lies in $\AuxSingMet[n]$, provided $R_1$ and $R_2$ were chosen large enough.
 
 We need to find a sufficient large lower bound for $R_1$ and $R_2$ such that $g$ satisfies the first condition. 
 This is slightly more tricky, because it involves second derivatives. 
 We abbreviate $H(q(\placeholder ), \chi(||\CubeProj{i}(\placeholder )||_K))$ to $\varphi$. 
 Recall the chain rule for the (total) second derivative:
 \begin{equation*}
   D^2(g_0 \circ \varphi) = D^2_{\varphi(\placeholder)} g_0 (D \varphi \cdot, D \varphi \cdot) + D_{\varphi(\placeholder)} g_0 \cdot D^2 \varphi(\cdot, \cdot).
\end{equation*}  
In matrix form, the first and second derivative of $\varphi$ are given by
\begin{align*}
 D \varphi &= DH \cdot (Dq, D_{||\placeholder||_K \circ \CubeProj{i}}\chi \cdot D_{\CubeProj{i}}||\placeholder||_K \cdot p_i(\placeholder )) \\
 &= Dq - \left(2 \eps (\rho + R_1)\cdot D_{||\placeholder||_K \circ p_i}\chi \cdot D_{\CubeProj{i}}||\placeholder||_K \cdot p_i(\placeholder )\right) e_i 
\end{align*}
and 
\begin{align*}
 \begin{split}
  D^2 \varphi &= D^2q - 2\eps(\rho+R_1) \cdot \bigl\lbrack (\chi'(||\CubeProj{i}(\placeholder )||_K))' \cdot D_{\CubeProj{i}}||\placeholder||_K \cdot \CubeProj{i}  \\
  & \quad \quad \quad \quad \quad + \chi'(||\CubeProj{i}(\placeholder )||_K) D(D_{\CubeProj{i}}||\placeholder||_K \cdot \CubeProj{i})\bigr\rbrack e_i 
 \end{split}\\
 \begin{split}
  &= D^2q - 2 \eps(\rho+R_1)\bigl\lbrack  \chi''(||\CubeProj{i}(\placeholder )||_K) \cdot (D_{\CubeProj{i}}||\placeholder||_K \cdot p_i \otimes D_{\CubeProj{i}}||\placeholder||_K \cdot p_i ) \\
  & \qquad \qquad \qquad \qquad  \qquad  + \chi'(||p_i(\placeholder )||_K) \cdot D^2_{\CubeProj{i}} ||\placeholder||_K(\CubeProj{i}(\placeholder ), \CubeProj{i}(\placeholder ))\bigr\rbrack e_i.
 \end{split}
\end{align*}
This yields, for the first and second derivatives of the block map $g$ on $Q_K^c$, the following expressions: 
For $j< i$, we have
\begin{align*}
 \partial_jg &= \partial_j g_0 \circ \varphi - \partial_i g_0 \circ \varphi \cdot 2\eps(\rho+R_1)\chi' \partial_j||\text{-}||_K, \\
 \begin{split}
 \partial_j^2 g &= D^2g(e_j,e_j) = (\partial_j^2g_0)\circ \varphi + (\partial_i^2g_0)\circ \varphi \cdot \bigl(2\eps(\rho+R_1)\chi' \partial_j ||\text{-}||_K\bigr)^2 \\
 \begin{split}
   &\qquad \qquad \qquad \qquad -4 (\partial_i\partial_jg_0) \circ \varphi \cdot \eps(\rho + R_1)\chi'\partial_j ||\placeholder||_K\\
   &\qquad \qquad \qquad \qquad + (\partial_i g_0)\circ \varphi \cdot 2\eps(\rho+R_1)\bigl(\chi'' (\partial_j||\text{-}||_K)^2 + \chi' \partial_j^2||\text{-}||_K\bigr), 
 \end{split}
 \end{split}
\end{align*}
for $j>i$, we have
\begin{align*}
 \partial_jg &= \partial_j g_0 \circ \varphi - \partial_i g_0 \circ \varphi \cdot 2\eps(\rho+R_1)\chi' \partial_j||\text{-}||_K, \\
 \begin{split}
 \partial_j^2 g &= D^2g(e_j,e_j) = (\partial_j^2g_0)\circ \varphi + (\partial_i^2g_0)\circ \varphi \cdot 2\eps(\rho+R_1)\chi' \partial_{j-1} ||\text{-}||_K \\
 \begin{split}
     &\qquad \qquad \qquad \qquad -4(\partial_i\partial_j g_0\circ \varphi) \cdot \eps(\rho + R_1)\chi'\partial_{j-1} ||\placeholder||_K\\
     &\qquad \qquad \qquad \qquad + (\partial_i g_0)\circ \varphi \cdot 2\eps(R+R_1)[\chi'' (\partial_{j-1}||\text{-}||_K)^2 + \chi' \partial_{j-1}^2||\text{-}||_K],
 \end{split}
 \end{split}
\end{align*}
and, for $j = i$, we have
\begin{align*}
 \partial_i g &= (\partial_i g_0) \circ \varphi \cdot q', \\ 
 \partial_i^2g &= (\partial_i^2 g_0)\circ \varphi \cdot (q')^2 + (\partial_i g_0)\circ \varphi q''.
\end{align*}
By definition, $|q''| \leq 1.2/R_1$, so we can choose $R_1$ sufficiently large such that
 \begin{align*}
  &\quad \sup\biggl\{\trace(\partial_i^2 g) - \trace(\partial_i^2 g_0)\circ \varphi \cdot (q')^2 \, : \, t \in Q_K^c \biggr\} \\
  &< \frac{1}{3}\min \biggl\{\frac{1}{8}\scal(g_0) - \sum_{j=1}^n |\trace( (\partial_j g_0^\op)^2)| + |\trace(\partial^2_j g^\op_0)| + \trace(\partial_j g^\op_0)^2 \biggr\}.
 \end{align*}
 
Recall that $\partial_j ||\text{-}||_K$ and $\partial_j^2||\text{-}||_K$ are bounded on every closed subset of $\R^n$ that does not contain the origin. 
Also recall that $g_0$ has bounded first and second derivatives on $Q_K^c$. 
Since $|\chi'|$ or $\chi''$ appear in all partial derivatives $\partial_jg$ and $\partial_j^2 g$ except for $j=i$ and since $|\chi'|\leq 1.3/R_2$ and $|\chi''|\leq 10/R_2^2$, we can choose $R_2$ to be sufficiently large  (in dependence of $R$ and $R_1$) such that 
\begin{align*}
 \begin{split}
  &\quad \sup_{t \in \overline{Q_K^c}}\biggl\{\sum_{j=1}^n |\trace(( (\partial_j g)^\op)^2) - \trace( (\partial_j g_0^\op)^2) \circ \varphi| \\
  &\quad \qquad + \sum_{j\neq i}^n |\trace(\partial_j^2 g) - \trace(\partial_j^2 g_0)\circ \varphi| \\
  &\quad \qquad + \sum_{j=1}^n |\trace(\partial_j g)|^2 - |\trace(\partial_j g_0)|^2 \circ \varphi \biggr\} \\
  &\quad < \frac{1}{3}\min \biggl\{\frac{1}{8}\scal(g_0) - \sum_{j=1}^n |\trace( (\partial_j g^\op_0)^2)| + |\trace(\partial^2_j g_0)| + |\trace(\partial_j g_0)|^2 \biggr\}.
 \end{split}
\end{align*}
Since $g_0$ satisfies the Condition \ref{Eq: AuxSing1}, also $g$ satisfies it by a simple triangle inequality argument.

Thus, if we increase $R_2$ such that it additionally satisfies $R_2 > 2(\rho + R_1)B_1B_2\mathcal{I}^{-1}$, where $\mathcal{I} = \min\{\mathcal{I}_1,\mathcal{I}_2\}$ is the minimum of
\begin{equation*}
 \mathcal{I}_1 \deff \mathrm{inf}\left\{\frac{1}{32} \, \lowBnd_{\Dirac_{g_0}}(t)^2 - \frac{1}{2}\sum_{j=1}^n ||[\Dirac_{g_0}^\extension, \nabla_{j}^{\Spinor_{\susp g_0}}] || \, : \, t \in Q_K^c\right\}
\end{equation*}
and
\begin{equation*}
 \mathcal{I}_2 \deff \mathrm{inf}\left\{ \mathrm{min}\{d^{-1}, 2^{-(d+4)} \lowBnd(t)\} - \sum_{k=1}^n ||\partial_kg_0^\op||_{\op,g_0(t)} \, : \, t\in Q_K^c\right\},
\end{equation*}
then $g$ also satisfies the remaining defining conditions of $\AuxSingMet$. 
Note that $\mathcal{I}>0$ because $g_0|_{Q_K^c}$ is a union of restrictions of degenerate elements of $\AuxSingMet[n]$ and the defining conditions of $\AuxSingMet[n]$ are local.

Thus, $g \in \AuxSingMet[n]$ is a filler for the given cubical $n$-horn and $\AuxSingMet$ is therefore a Kan set.
\end{proof}
\begin{cor}
 The inclusion $\AuxSingMet \hookrightarrow \SingMet$ is a weak homotopy equivalence. 
\end{cor}
\begin{proof}
 The proof is analogous to the proof of Corollary \ref{AuxDir is weak equivalent -  Cor} that $\AuxInvDir \hookrightarrow \InvPseudDir[\bullet]$ is a weak homotopy equivalence.
\end{proof}
\begin{lemma}\label{SuspOnAux - Lemma}
 The suspension restricts to a cubical map $\susp \colon \AuxSingMet \hookrightarrow \ConcSet$.
\end{lemma}
\begin{proof}
 By an iterative application of Lemma \ref{Scalar Curvature Formuala - Prop} we have
 \begin{align*}
  |\scal(\susp(g)) - \scal(g)| &= \left|\sum_{j=1}^n \frac{3}{4}(\trace((\partial_jg)^\op)^2) - \trace(\partial_j^2g) - \frac{1}{4}\trace(\partial_j g)^2 \right| \\
  &< \sum_{j=1}^n |\trace((\partial_jg)^\op)^2)| + |\trace(\partial_j^2g)| + |\trace(\partial_j g)^2| \\
  &< \frac{1}{8} \scal(g). 
\end{align*}  
Thus, $\scal(\susp(g))(m,t) > \frac{7}{8}\scal(g(t))(m) > 0$, which proves the claim.
\end{proof}
The Lichnerowicz formula implies that $\Dirac_\bullet$ restricts to a cubical map $\ConcSet \rightarrow \InvBlockDirac$. 
By the second and third defining condition of  $\AuxSingMet$ the Dirac Operator restricts to a cubical map $\AuxSingMet \hookrightarrow \AuxInvDir$. Thus, we have two cubical maps
\begin{equation*}
 \Dirac_\bullet \circ \susp, \, \susp \circ \Dirac_\bullet \colon \AuxSingMet \rightarrow \InvBlockDirac.
\end{equation*}
Surprisingly, these two maps are not the same.
We will show that they are homotopic by controlling their difference.
\begin{definition}
 Define the \emph{error term} to be
 \begin{equation*}
   \Err(g) \deff \Dirac_{\susp(g)} - \susp(\Dirac_g).
 \end{equation*}
\end{definition}
\begin{lemma}\label{Error Formula - Lemma}
 For each block map $g \in \SingMetNon[M][n]$, the error term $\Err(g)$ is a differential operator of order zero that has block form. It can be expressed by the formula
 \begin{equation*}
   \Err(g) = \sum_{j=1}^n \frac{1}{4}\left(\trace{ \partial_j g^\op}\right)\partial_j \cdot.
\end{equation*}  
\end{lemma}
\begin{proof}
 As a difference of two differential operators of order 1 with the same principal symbol, $\Err(g)$ is a differential operator of order zero. 
 The formula follows from the calculation of the Christoffel symbols in Lemma \ref{CalcChristoffelSymb - Lemma}.
 Since $g$ is a block metric, $\Err(g)$ decomposes accordingly, so it is a block operator of order zero.
\end{proof}

\begin{lemma}\label{Error Estimation Lemma}
 For each $g \in A_n$ and all $u \in H^1(\Spinor_{\susp(g)})$ the following estimates hold true
 \begin{align*}
     ||\Err||^2_{\op,t} &< \frac{1}{64} \scal(g(t)), \\
     ||\Dirac_{\susp (g)}u||_0 - ||\Err(g)u||_0 &\geq \frac{1}{2}\min\{1, 1/8\inf\{\scal(g(t))^{1/2} \, : \, t\in\R^n\}\} \cdot ||u||_1.
 \end{align*}
\end{lemma}
\begin{proof}
 The first estimate follows from Lemma \ref{Error Formula - Lemma} and the fact that the (fibre-wise) operator norm satisfies the $C^\ast$-identity
 \begin{align*}
  ||\Err||^2_{\op,t} &= ||\Err^\ast \Err||_{\op,t} \\
  &= \left|\left| \frac{1}{16} \sum_{k,l} \trace(\partial_k g) \trace(\partial_l g) \partial_k^\ast \partial_l \right|\right|_{\op,t} \\
  &=\left|\left| \frac{1}{16} \sum_{k\leq l} \trace(\partial_k g) \trace(\partial_l g)(\partial_k\partial_l + \partial_l \partial_k) \right|\right|_{\op,t} \\
  &=\left|\left| \frac{2}{16} \sum_{k=1}^n \trace(\partial_kg)^2 \id \right|\right|_{\op,t} \\
  &= \frac{2}{16} \sum_{k=1}^n \trace(\partial_k g(t))^2 < \frac{1}{8}\cdot\frac{1}{8} \scal(g(t)). 
 \end{align*}
 This inequality together with the Lichnerowicz formula implies the second inequality by the following chain of estimations. 
 \begin{align*}
 &\quad \ \, ||\Dirac_{\susp (g)}u||_0 - ||\Err(g)u||_0 \\
  &=\left(||\nabla u||_0^2 + \frac{1}{4}||\scal(g)^{1/2}u||_0^2\right)^{1/2} - ||\Err(g)u||_0  \\
  &\geq \left(||\nabla u||_0^2 + \frac{1}{4}||\scal(g)^{1/2}u||_0^2\right)^{1/2} - \sqrt{\frac{1}{64}} ||\scal(g)^{1/2}u||_0 \\
  &\geq \left(||\nabla u||_0^2 + \frac{1}{4}||\scal(g)^{1/2}u||_0^2\right)^{1/2} - 1/8 ||\scal(g)^{1/2} u||_0 \\
  &\geq \frac{1}{2}\left( ||\nabla u||_0^2 + \frac{1}{8} ||\scal(g)^{1/2}u||_0^2 \right)^{1/2} \\
  &\geq \frac{1}{2}\min\Bigl\{1, 1/8\inf\{\scal(g(t))^{1/2} \, : \, t\in\R^n\}\Bigr\} \cdot ||u||_1. \qedhere
 \end{align*}
\end{proof}

We can use the previous Lemma to re-prove that $\susp(\Dirac_g)$ is invertible, something we already know to be true by Condition \ref{Eq: AuxSing2} and \ref{Eq: AuxSing3}.

\begin{cor}
 For each $g \in \AuxSingMet[n]$ the block Dirac operator $\susp(\Dirac_g)$ is invertible.
\end{cor}
\begin{proof}
 The Lichnerowicz formula implies $||\Dirac_{\susp(g)} u||^2_0 \geq 1/4||\scal(g)^{1/2}u||^2_0$, in particular, this operator is invertible.
 From Lemma \ref{Error Estimation Lemma} and
 \begin{align*}
     ||\susp(\Dirac_g)u||_0 \geq \left| ||\Dirac_{\susp g}u||_0 - ||\Err(g)u||_0\right| 
 \end{align*}
 we conclude that $\susp(\Dirac_g)$ has trivial kernel and closed image. 
 Since $\susp(\Dirac_g)$ is self-adjoint, it is also surjective and hence invertible.
\end{proof}
\begin{definition}
 Let $\chi \colon \R \rightarrow [0,1]$ be a smooth map that is identically zero on $\R_{\leq -1}$ and identically $1$ on $\R_{\geq 1}$ and let 
 \begin{equation*}
  \mathcal{H}_\bullet \colon (\StandCube[ ]{1} \otimes \AuxSingMet[ ])_{\bullet} \rightarrow \InvBlockDirac
\end{equation*}  
 be the following cubical map:
 If $\varphi \in \StandCube[k]{1}$  and $g \in \AuxSingMet[n-k]$ then
 \begin{equation*}
  \mathcal{H}_n(\varphi,g) = \chi \circ \bm{\sigma}(\varphi) \cdot \susp \Dirac_{\degen{1}^k g} + (1 -\chi) \circ \bm{\sigma}(\varphi) \cdot \Dirac_{\susp(\degen{1}^k g)},
\end{equation*}  
where $\bm{\sigma}\varphi = \degen{n}\dots\degen{k+1}\varphi$.
More precisely, $\mathcal{H}_n(\varphi,g) \in \PseudOpCl^1(\Spinor_{\susp(\degen{1}^kg)})$ is the unique operator defined by 
\begin{align*}
 \mathcal{H}_n&(\varphi,g)(u)(m,t) = \\
 &\chi(\bm{\sigma}(\varphi)(t)) \cdot \left(\susp \Dirac_{\degen{1}^kg} u\right)(m,t) + (1 -\chi(\bm{\sigma}(\varphi)(t)) \cdot \left(\Dirac_{\susp(\degen{1}^kg)}u\right)(m,t).
\end{align*}
\end{definition}
\begin{theorem}
 The sequence $\mathcal{H}_\bullet$ consist of well defined maps and assembles to a cubical map. It is a homotopy between $\susp \circ \Dirac_\bullet$ and $\Dirac_\bullet \circ \susp$.
\end{theorem}
\begin{proof}
 For well-definedness, we need to check that $\mathcal{H}_{n+1}(\degen{k+1}\varphi, g) = \mathcal{H}_{n+1}(\varphi,\degen{1}g)$ for all $n \in \N_0$ and that $\mathcal{H}_n(\varphi,g)$ is an invertible block Dirac operator for all $\varphi \in \StandCube[k]{1}$ and all $g \in \AuxSingMet[n-k]$.  
 
 We first show that $\mathcal{H}_n$ preserves the equivalence relation of $(\StandCube[ ]{1} \otimes \AuxSingMet[ ])_\bullet$. 
 Let $\varphi \in \StandCube[k]{1}$ and $g \in \AuxSingMet[n-k]$ be given.
 Since $\degen{k+1}\varphi \in \StandCube[k+1]{1}$, we have
 \begin{align*}
  & \quad \, \mathcal{H}_{n+1}(\degen{k+1}(\varphi),g) \\
  \begin{split}
      &= \chi \circ \degen{n} \circ \dots \circ \degen{k+2}(\degen{k+1}\varphi)  \susp \, \Dirac_{\degen{1}^{k+1}g} \\
      &\qquad + (1 - \chi \circ \degen{n} \circ \dots \circ \degen{k+2}(\degen{k+1}\varphi) )\Dirac_{\susp (\degen{1}^{k+1}g)}
  \end{split}\\
  &= \chi \circ \degen{n} \circ \dots \circ \degen{k+1}\varphi  \, \susp \, \Dirac_{\degen{1}^{k+1}g} + (1 - \chi \circ \degen{n} \circ \dots \circ \degen{k+1}\varphi )\Dirac_{\susp (\degen{1}^{k+1}g)}\\
  &= \mathcal{H}_{n+1}(\varphi,\degen{1}g)
 \end{align*}
 
 During the entire proof, we slightly abuse notation by not indicating the length of composition $\bm{\sigma}$, i.e., we write $\bm{\sigma}\varphi$ instead of $\bm{\sigma}\degen{k+1}\varphi$.
 
The operator $\mathcal{H}_{n}(\varphi,g)$ has block form. 
Indeed, if $\varphi$ is constant then
\begin{align*}
 \mathcal{H}_n(\varphi,g) = \begin{cases}
  \Dirac_{\susp (\degen{1}^kg)}, & \text{if } \varphi \equiv -1,\\
  \susp (\Dirac_{\degen{1}^kg}), & \text{if } \varphi \equiv 1.
 \end{cases}
\end{align*}  
If $\varphi$ is not constant then, by Lemma \ref{Box Morphism Ordering Lemma}, there is a unique coordinate $t_s$ with $1 \leq s \leq k$ such that $\varphi(t) = t_s$. 
It is easy to see that $M \times RI^n$, which is the core of  $\Dirac_{\susp (\degen{1}^kg)}$ and $\susp (\Dirac_{\degen{1}^kg})$, is also the core of $\mathcal{H}_n(g,\varphi)$. 
Since $\mathcal{H}_n(\varphi,g)$ is a differential operator it decreases the support, so that it restricts in particular to each $U_R(\fateps)$.
The $\mathcal{C}^\infty(\R^n)$-linearity of $(\Phi_\fateps)_\ast$ implies 
\begin{align*}
  &(\Phi_\fateps)_\ast\Bigl(\mathcal{H}_n(\varphi,g)|_{U_R(\fateps)}\Bigr) = \\
 &\begin{split}
  \Biggl( \chi(t_{c(\fateps)(s)}) \susp \Dirac_{\degen{1}^kg}(\fateps|_{\{1,\dots,k\}}) &+ (1-\chi({c(\fateps)(s)}))\Dirac_{\susp \sigma g}(\fateps_{\{1,\dots,k\}}) \boxtimes \id   \\ 
   & +\left. \evenodd \boxtimes \sum_{j \in \mathrm{dom}\fateps|_{\{k+1,\dots,n\}}} \partial_j \frac{\partial}{\partial t_j}\right).
 \end{split}
\end{align*}  
Hence, for all sufficient large $R>0$ such that $g|_{U_R(\fateps)}$ does not depend on the variables parametrised by $\supp \, \fateps$, the operator $\mathcal{H}_n(\varphi,g)$ decomposes as required on $U_R(\fateps)$.

The operator $\mathcal{H}_n(\varphi,g)$ is a convex combination of $\susp (\Dirac_{\degen{1}^kg})$ and $\Dirac_{\susp (\degen{1}^kg)}$ using positive real valued function instead of real numbers.
This implies that $\mathcal{H}_n(\varphi,g)$ is a differential operator with the same symbol as $\Dirac_{\susp(\degen{1}^k{\sigma}(g))}$ and the block form implies that it differs from $\Dirac_{\susp(\degen{1}^kg)}$ by a bounded operator.
As a convex combination of symmetric differential operators, $\mathcal{H}_n(\varphi,g)$ is a symmetric differential operator whose principal symbol has finite propagation speed, hence $\mathcal{H}_n(\varphi,g)$ is self-adjoint by \cite{higson2000analytic}*{Prop 10.2.11}

So far, we have shown that $\mathcal{H}_n(\varphi,g)$ is a block Dirac operator. 
Next, we will show that it is invertible by showing that it is bounded from below using Lemma \ref{Error Estimation Lemma}:
\begin{align*}
 ||\mathcal{H}_n(\varphi,g)||_0 &\geq \left| ||\Dirac_{\susp (\degen{1}^kg)} u||_0 - ||(\mathcal{H}_n(\varphi,g) - \Dirac_{\susp (\degen{1}^kg)})u||_0\right| \\
 &= \left| ||\susp\Dirac_{ \degen{1}^kg}u||_0 - ||(\chi \circ \bm{\sigma}(\varphi))\susp \Dirac_{\degen{1}^kg} - (\chi \circ \bm{\sigma}(\varphi)) \Dirac_{\susp (\degen{1}^kg})||_0\right| \\
 &= \left| ||\Dirac_{\susp \degen{1}^kg} u||_0 - ||(\chi \circ \bm{\sigma}(\varphi)) \Err(\degen{1}^kg)u||_0\right| \\
 &\geq ||\Dirac_{\susp (\degen{1}^kg)}u||_0 - ||\Err(\degen{1}^kg) u||_0 \\
 &\geq \frac{1}{2} \min\Bigl\{1, 1/8 \inf\{\scal(\degen{1}^kg)(t)^{1/2} \, : \, t \in \R^n \}\Bigr\} ||u||_1.
\end{align*} 
We verify now that $\mathcal{H}_\bullet$ is a cubical map.
From the proof that $\mathcal{H}_n(\varphi,g)$ is a block map and that $\susp$ and $\Dirac$ are cubical maps we deduce 
\begin{align*}
 \face{i} \mathcal{H}_n(\varphi,g) &= \lim_{R \to \infty} \CubeIncl[R\eps]{i}^\ast H_n(\varphi,g)(i,\eps) \\
            \begin{split}&= \lim_{R \to \infty} \CubeIncl[R\eps]{i}^\ast\left( \chi \circ \bm{\sigma}(\varphi)|_{M \times \R^n(i,\eps)} \cdot \susp \, \Dirac_{\degen{1}^kg}(i,\eps) \right. \\
            & \qquad \left. \left( 1 -\chi \circ \bm{\sigma}(\varphi)|_{M \times \R^n(i,\eps)} \right)\cdot \Dirac_{\susp (\degen{1}^kg)}(i,\eps)\right) \end{split}\\
            \begin{split}
            &= \lim_{R \to \infty} \CubeIncl[R\eps]{i}^\ast(\chi \circ \bm{\sigma}(\varphi)|_{M \times \R^n(i,\eps)}) \cdot \susp \Dirac_{\CubeIncl[R\eps]{i}^\ast \degen{1}^kg} \\
            & \qquad 1-\CubeIncl[R\eps]{i}^\ast(\chi \circ \bm{\sigma}(\varphi)|_{M \times \R^n(i,\eps)})  
            \Dirac_{\susp (\CubeIncl[R\eps]{i}^\ast \degen{1}^kg)}
            \end{split} \\
            &= \begin{cases}
             \mathcal{H}_{n-1}(\face{i}\varphi,g), & \text{if } i \leq k, \\
             \mathcal{H}_{n-1}(\varphi,\face{i-k}g), & \text{if } i > k.
            \end{cases}
\end{align*}
For the degeneracies, we use the simple observation $(f\cdot P)\boxtimes \id = \degen{n+1}(f)\cdot (P\boxtimes \id)$ in the following calculation:
\begin{align*}
 &\degen{i} \mathcal{H}_{n}(\varphi,g) = \Cycl(i,n+1)_\ast\left(\mathcal{H}_n \boxtimes \id + \evenodd \boxtimes \DiracSt_{\R} \right) \\
 \begin{split}
  &= \Cycl(i,n+1)_\ast\left(\chi \circ \degen{n+1}\bm{\sigma}(\varphi) (\susp \Dirac_{\degen{1}^kg}\boxtimes \id) \right. \\
  & \left. \qquad \qquad \qquad \qquad \,  + (1- \chi) \circ \degen{n+1} \bm{\sigma}(\varphi) (\Dirac_{\susp \, \degen{1}^kg}\boxtimes \id) + \evenodd \boxtimes \DiracSt_{\R} \right) 
  \end{split}\\
  \begin{split}
    &= \Cycl(i,n+1)_\ast\left( \chi \circ \degen{n+1}\bm{\sigma}(\varphi) \susp \Dirac_{\degen{n+1} \degen{1}^kg} \right. \\
    & \left. \qquad \qquad \qquad \qquad \,+ (1- \chi \circ \degen{n+1}\bm{\sigma}(\varphi))  \Dirac_{\susp \degen{n+1}\degen{1}^kg}\right)  
  \end{split}\\
  &= \chi \circ \degen{i}\bm{\sigma}(\varphi) \susp(\Dirac_{\degen{i}\degen{1}^kg}) + (1 - \chi)\circ \degen{i}\bm{\sigma}(\varphi) \Dirac_{\susp \degen{i}\degen{1}^kg} \\
 &= \begin{cases}
  \mathcal{H}_{n+1}(\degen{i}\varphi,g), & \text{if } i \leq k, \\
  \mathcal{H}_{n+1}(\varphi,\degen{i-k}g), & \text{if } i > k.
\end{cases}
\end{align*}

Lastly, we need to show that $\mathcal{H}_\bullet$ is a homotopy between $\susp \circ \Dirac_\bullet$ and $\Dirac_\bullet \circ \susp$.
If $\varphi = (-1,\dots,-1) \in \StandCube[k]{1}$ then 
\begin{equation*}
 \mathcal{H}_n(\varphi,g) = \susp (\Dirac_{\degen{1}^kg}) = \degen{1}^k( \susp \Dirac_g)
\end{equation*}
and if $\varphi = (1,\dots, 1) \in \StandCube[k]{1}$ then
\begin{equation*}
 \mathcal{H}_n(\varphi,g) = \Dirac_{\susp (\degen{1}^kg)} = \degen{1}^k( \Dirac_{\susp (g)}),
\end{equation*}
so $\mathcal{H}_{n}$ is indeed a homotopy between $\susp \circ \Dirac_\bullet$ and $\Dirac_\bullet \circ \susp$.
\end{proof}

\section{Applications}\label{Section - Applications}

We sketch how our work relates to the existing literature and give some application of our main theorem.
The presentation here is less detailed as in the rest of this thesis as the missing details can be found in the existing literature.
We also give a counterexample to the stronger concordance-implies-isotopy conjecture, namely that the suspension $\susp \colon \SingMet \dashrightarrow \ConcSet$, defined on a weakly equivalent subset, is a weak equivalence.

\subsubsection*{Relations to the work of Botvinnik, Ebert, Randal-Williams:}

The strongest detection result for non-trivial elements in the homotopy groups of $\Riem^+(M)$ that does not take the fundamental group of the manifold into account was achieved by Botvinnik, Ebert and Randal-Williams in \cite{botvinnik2017infinite}.
Using sophisticated tools from homotopy theory, they proved that the index difference is rationally surjective.
\begin{theorem}[\cite{botvinnik2017infinite}*{Theorem A}]
  Let $M^d$ be a closed spin manifold of dimension $d\geq 6$.
  Then $$\pi_n(\mathrm{inddif})\otimes \mathbb{Q} \colon \pi_n(\Riem^+(M)) \otimes \Q \rightarrow KO^{-(d+n+1)}(\mathrm{pt}) \otimes \Q$$ is surjective.
  If $d+n+1 \equiv 1,2 \mod 8$ so that $KO^{-(d+n+1)}(\mathrm{pt}) \cong \Z_2$, then $\pi_n(\mathrm{inddif})$ is surjective.
\end{theorem}
 By Theorem \ref{IndDifFactor - Theorem}, index difference factors over $\ConcSet$, so we deduce:
 \begin{cor}
  Let $M^d$ be a closed spin manifold of dimension $d\geq 6$.
  Then 
  \begin{equation*}
      \pi_n(\Dirac_{\bullet})\otimes \Q \colon \pi_n(\ConcSet)\otimes \Q \rightarrow \pi_n(\InvBlockDirac)\otimes \Q = KO^{-(d+n+1)}(\mathrm{pt}) \otimes \mathbb{Q}
  \end{equation*}
  is surjective. 
  If $d+n+1 \equiv 1,2 \mod 8$ so that $KO^{-(d+n+1)}(\mathrm{pt}) \cong \Z_2$, then $\pi_n(\Dirac_{\bullet})$ is surjective. 
 \end{cor}

\subsubsection*{Relations to the work of Hitchin:}
Hitchin was the first to discover that $\Riem^+(M)$ can be non-contractible if $M^d$ is any closed spin manifold of dimension $8k$ that carries at least one psc metric. 
Despite being a consequence of the previously discussed by results of Botvinnik, Ebert and Randal-Williams, it is worthwhile to recall Hitchin's results in detail as we will use them below.
Hitchin's argument in \cite{hitchin1974harmonic} is roughly the following one:
The diffeomorphism group $\mathrm{Diff}(M)$ acts on $\Riem^+(M)$ via pull back.
For each choice of a base point $g_0 \in \Riem^+(M)$, the pull back produces elements $[\varphi^\ast g_0]\in \pi_n(\Riem^+(M))$ for each element $[\varphi]\in \pi_n(\mathrm{Diff}(M))$.
The key observation is now that the composition of the (Hitchin) index difference and the pull back action
\begin{equation*}
    \xymatrix{\pi_n(\mathrm{Diff}(M)) \ar[rr]^-{\text{pull back}} && \pi_n(\Riem^+(M)) \ar[rr]^-{\mathrm{inddif}} && KO^{-(d+n+1)}(\mathrm{pt}) }
\end{equation*}
agrees with the topological index of the total space of the clutching construction $M \times D^{n+1} \cup_\varphi M \times D^{n+1} \rightarrow S^{n+1}$, provided the total space is a spin manifold.
More precisely, under the aforementioned conditions, we have
\begin{equation*}
    \mathrm{inddif}(\varphi^\ast g_0,g_0) = \alpha(M\times D^{n+1} \cup_\varphi M \times D^{n+1}).
\end{equation*}

In the case $n=0,1$, there are diffeomorphisms $\varphi \in \pi_n(\mathrm{Diff}(M))$ such that
\begin{equation*}
    M \times D^{n+1} \cup_\varphi M \times D^{n+1} \cong (M \times S^{n+1}) \# \Sigma^{d+n+1},
\end{equation*}
where $\Sigma^{d+n+1}$ is an exotic sphere with $\alpha(\Sigma^{d+n+1}) \neq 0$.
Since the topological index $\alpha \colon \Omega_\ast^{\Spin}(\mathrm{pt}) \rightarrow KO^{-\ast}(\mathrm{pt})$ is a ring-homomorphism, we deduce
\begin{align*}
    \alpha((M\times S^{n+1})\#\Sigma^{d+n+1}) = \alpha(M)\underbrace{\alpha(S^{n+1})}_{=0} + \alpha(\Sigma^{d+n+1}) \neq 0.
\end{align*}
This implies that $\pi_0(\Riem^+(M))$ and $\pi_1(\Riem^+(M))$ are non-trivial.

For later purposes we remark that the special diffeomorphisms are in the image of $\pi_n(\mathrm{Diff}(D^d,S^{d-1}))$, the group of diffeomorphisms on an embedded disc $D^d \subseteq M$ that are the identity near $S^{d-1}$.

The connection to our work is given by Theorem \ref{IndDifFactor - Theorem}. 
The left vertical map is the a cubical map that models the space level description of the (Hitchin) index difference. 
Applying homotopy groups to the diagram in Theorem \ref{IndDifFactor - Theorem} together with Hichtin's result implies the following corollary.
\begin{cor}
 Let $M^d$ be a closed, spin, psc manifold of dimension $8k$.
 Then $\pi_0(\ConcSet)$ and $\pi_1(\ConcSet)$ are non-trivial.
\end{cor}

\subsubsection*{Relations to the work of Ruberman:}

We present an example where $\susp_\bullet \colon \SingMet[M] \dashrightarrow \ConcSet$ is not a weak homotopy equivalence.
The following proposition is just a reformulation of the results in \cite{ruberman2002positive}.
\begin{proposition}\cite{ruberman2002positive}*{Corollary 5.2, Theorem 5.5}
 For $n \geq 2$ and $k > 10n$, define $M \deff \#_{2n}\mathbb{CP}^2 \#_k \overline{\mathbb{CP}^2}$. 
 Then $\susp_\bullet \colon \SingMet[M] \dashrightarrow \ConcSet$ is not injective on the level of connected components.
 In fact, every concordance class contains infinitely many different isotopy classes.
\end{proposition}

\subsubsection*{Non Surjectivity Results:}

We can also apply Hitchin's methods directly to $\ConcSet$ to produce non-trivial elements in its homotopy groups. 
Fix a base point $g_0 \in \Riem^+(M)$. 
The corresponding base point in $\pi_k(\ConcSet)$ is then $g_0 \oplus \euclmetric_{\R^k}$.
Fix further an embedded disc $D^{d+k} \subseteq M \times \R^k$.

Let $\varphi$ be the diffeomorphism used by Hitchin in the subsection above (for $n=0$).
Since this diffeomorphism is the identity near the boundary of the embedded disc and since the disc lies in the interior of some $M \times RI^k$, the pull back $g \deff \varphi^\ast (g_0\oplus \euclmetric_{\R^k})$ is a block metric.
We claim that $[g] \neq 0 \in \pi_k(\ConcSet[\bullet])$ provided $d+k \equiv 0 \mod 8$.

Assume this were not the case, then there would be a block metric $G$ on $M \times \R^{k+1}$ that satisfies $ \face[-1]{1}G = g$ and $\face[\eps]{i}G = g_0 \oplus \euclmetric_{\R^k}$ if $(i,\eps) \neq (1,-1)$.
Pick a sufficiently large $R'>R$ such that $G$ decomposes outside of $M \times R'I^{k+1}$. 
Then $G$ induces a psc metric $\bar{G}$ on 
\begin{equation*}
    M \times \R \times \mathbb{T}^k = M \times \R \times (\R^k/2R'\Z^k),
\end{equation*}
which serves as a concordance between the psc metric $\bar{g} \in \Riem^+(M \times \mathbb{T}^k)$ and $g_0 \oplus \euclmetric_{\mathbb{T}^k}$.
Due to the locality of this construction, $\bar{g}$ is the psc metric we would obtain if we would apply Hitchin's construction to $g_0 \oplus \euclmetric_{\mathbb{T}^k}$ on $M \times \mathbb{T}^k$.
Thus, the Gromov-Lawson index difference of these two metrics vanishes.
On the other hand Theorem \ref{IndDifFactor - Theorem}, or alternatively Theorem A of \cite{ebert2017indexdiff}, together with Hitchin's result now imply that
\begin{align*}
    \mathrm{inddif}_{GL}\bigl(\bar{g},g_0\oplus \euclmetric_{\mathbb{T}^k}\bigr) = \mathrm{inddif}_{H}\bigl(\bar{g},g_0\oplus \euclmetric_{\mathbb{T}^k}\bigr) = \alpha(\Sigma^{d+k+1}) \neq 0,
\end{align*}
which is a contradiction.

If we now assume that $\dim M=2$ or $3$, then $M$ is spin if and only if it is orientable, see \cite{LawsonMichelsonSpin}*{p.86}.
The previous argument now implies that $\pi_6(\ConcSet[\bullet][S^2])$ and $\pi_5(\ConcSet[\bullet][M^3])$ are non-trivial.
However, we know from \cite{rosenberg149metrics}*{Theorem 3.4} and \cite{bamler2019ricci}*{Theorem 1.1} that $\Riem^+(M)$ is contractible for all orientable psc manifolds of dimension 2 or 3.
This argument proves Theorem \ref{ConterExampleStrongConcVsIsotopy - MainThm} of the introduction:
\begin{theorem}\label{CounterHigherConcVsIso - Theorem}
  The suspension $\susp_\bullet \colon \SingMet \dashrightarrow \ConcSet$ is not a weak equivalence if $\dim(M) = 2$ or $3$.
\end{theorem}


 \chapter{The PSC Hatcher Spectral Sequence}\label{The PSC Hatcher Spectral Sequence - Chapter}

A consequence of the Factorisation Theorem \ref{IndDifFactor - Theorem} is that the index difference cannot distinguish different isotopy classes in the same conocrdance class.
This question has to be addressed with tools and techniques beyond index theory.

As a first attempt, we will make full use of the ``spacification'' of the isotopy versus concordance problem we have carried out in Chapter \ref{Chapter - Cubical Versions of psc} and Chapter \ref{Factorisation Indexdiff - Chapter}:
The comparison map $\susp_\bullet \colon \SingMet \dashrightarrow \ConcSet$, defined on a weakly equivalent cubical subset $A_\bullet \subseteq \SingMet$, allows us to the define the homotopy fibre $\mathrm{hofib}\left(A_\bullet \rightarrow \ConcSet\right)$, which we abbreviate to $\hofib$.
Since $A_\bullet$ and $\ConcSet$ are cubical Kan sets, $\hofib$ is also a cubical Kan set by formal reasons. 
The combinatorial description of its homotopy groups allows us to put a filtration on them, the so called \emph{isotopy filtration}.
This filtration, roughly speaking, measures how far away from an isotopy a general concordance is.
The \emph{psc Hatcher spectral sequence}, the main result of this section stated in Theorem \ref{psc Hatcher SS}, approximates the difference between a filtration subgroup and its successor.

This chapter is structures as follows.
First, we will develop the geometric foundation to construct the psc Hatcher spectral sequence in Section \ref{Section - Geometric Foundations}.
This includes the discussion of the isotopy filtration in Definition \ref{Isotopy Filtration - Def}, the defect fibrations in Definition \ref{Restriction is Smooth Fib}, and the connection between those two, which is explained by Lemma \ref{Boundary Comparison - Lemma}.

Secondly, we present in Section \ref{Section - Abstract Construction of the Spectral Sequence} an abstract construction of a spectral sequence that generalises the construction of a spectral sequence out of an exact couple of abelian groups. 
In our setup, we allow some abelian groups to be replaced by monoids or even mere pointed sets.

Finally, in Section \ref{Section - Application to the Defect Fibration}, we apply this abstract spectral sequence construction of Theorem \ref{ExactCoupleYieldSS - Theorem} to the geometric foundation established before to derive the psc Hatcher spectral sequence.


\section{Geometric Foundations}\label{Section - Geometric Foundations}

Recall that we have defined in Chapter \ref{Factorisation Indexdiff - Chapter} a weakly equivalent cubical subset $A_\bullet \hookrightarrow \SingMet$, which is Kan and on which the suspension map takes values in $\ConcSet$.
In this section, it will be more convenient to work with a slightly bigger model given by
\begin{equation*}
    A_n \deff \left\{ g \in \SingMet[][n] \, : \, (1+d)\sum_{j=1}^n |\trace((\partial_j g^{\op})^2)| + |\trace(\partial_j g)|^2 + |\trace(\partial_j^2 g)| < \scal(g)\right\}.
\end{equation*}
Observe that it shares the just mentioned properties of its analog in Chapter \ref{Factorisation Indexdiff - Chapter}.

During the entire section, we fix a base point $g_0 \in \SingMet[M][0] = \ConcSet[0]$. 
We remind the reader that we do not notionally distinguish $g_0$ from its image under degeneracy maps. 
We furthermore abbreviate $\hofib(\susp_\bullet \colon A_\bullet \rightarrow \ConcSet;g_0)$ to $\hofib$ and the dimension of $M$ with $d$.

Let us start with the definition of the isotopy filtration.
For that recall the set of loops of a cubical set from Definition \ref{combitarial loop space - Def}.
\begin{definition}\label{Isotopy Filtration - Def}
 For $0\leq p \leq n$, we set the $p$-th term of the \emph{isotopy filtration} to be
 \begin{align*}
     \begin{split}
         F_p\Omega^n_{g_0}\hofib \deff \Bigl\{&g \in \Omega^n_{g_0} \hofib \, : \, g_{(m,x)} = g^{(p+1)}_{(m,x)} + \sum_{j=p+2}^{n+1} \diff x_j^2 \text{ and} \\
         &(1+d)\sum_{j=p+2}^{n+1}  |\trace((\partial_j g^\op)^2)| + |\trace(\partial_j g)^2| + |\trace(\partial_j^2 g)| < \scal(g)\Bigr\}.
     \end{split}
 \end{align*}
 Here, $g^{(p+1)}_{(m,x)}$ denotes the restriction of $g$ to the subspace $T_m M \oplus \mathrm{span}(e_1,\dots,e_{p+1})$ in $T_mM \oplus T_x\R^{n+1}$.
\end{definition}
The filtration of $\Omega^n_{g_0}\hofib$ induces a filtration of $\pi_n(\hofib)$, the set of concordance classes of $\Omega^n_{g_0}\hofib$.
Informally speaking, the set $F_p\Omega^n_{g_0}\hofib$ consists of block metrics that are isotopies in the last $(n-p)$ coordinates.
This observation inspires Definition \ref{Filtrationspace - Def} and Lemma \ref{SuspensionFiltration - Lemma}.

\begin{definition}\label{Filtrationspace - Def}
 Let $\ConcSet[ ][M \times RI^n]$ be the topology space of all Riemannian metrics with positive scalar curvature that have product structure near each face equipped with the smooth Whitney topology.
 For $p \geq 0$, let $F_p$ be the topological space whose underlying set is
 \begin{equation*}
  F_p \deff \left\{g \in \ConcSet[p+1] \, : \, \face[1]{1}g \in \SingMet[][p], \, \face[\omega]{j}g = g_0 \text{ for } (j,\omega) \neq (1,1)\right\}
 \end{equation*}
 and whose topology is the (subspace topology of the) colimit topology coming from the embeddings $\ConcSet[ ][M \times RI^{p+1}] \hookrightarrow \ConcSet[ ][M \times SI^{p+1}]$ given by elongations.
\end{definition}

We remind the reader that $F_p$ carries the colimit topology induced by embeddings between Hausdorff spaces.
By Lemma \ref{Relative T1} and Lemma \ref{EmbeddingsReloativeT1 - Lemma}, every continuous map from a compact space into $F_p$ factors through some $\ConcSet[ ][M \times RI^{n+1}]$.
This observation will be important in the forthcoming proofs.

\begin{lemma}\label{SuspensionFiltration - Lemma}
 The suspension induces a map $\pi_{n-p}(F_p) \rightarrow \pi_n(\hofib)$ whose image is the set of concordance classes of $F_p\Omega^n_{g_0}\hofib$.
 Furthermore, the suspensions decompose into a sequence of maps
 \begin{equation*}
  \xymatrix@C+1em{\pi_n(F_0) \ar[r]^{\iota_{0,n}} \ar[drrr] & \pi_{n-1}(F_1) \ar[r]^-{\iota_{1,n-1}} \ar[drr]& \dots \ar[r]^-{\iota_{n-1,1}} &\pi_0(F_n) \ar@{->>}[d]^{\iota_{n,0}} \\ &&& \pi_n(\hofib), }
 \end{equation*}
 where all $i_{p,n-p}$ are induced by ``suspending along the first parameter variable'' and $\iota_{n,0}$ sends an isotopy class to its concordance class.
\end{lemma}
\begin{proof}
 For later purposes, we define the morphisms $\iota_{p,q}$ as induced morphisms from continuous maps.
 Set
 \begin{equation*}
     F_p^{slow} \deff \{g \in F_p \, : \, \face[1]{1}g \in  A_p(M)\}
 \end{equation*}
 and
 \begin{align*}
    \begin{split}
      \Omega^{q,\infty} F_p^{(slow)} \deff \{ & g \colon \R^q \rightarrow F_p^{(slow)} \, : \, \text{smooth block map}, \ g = g_0 \text{ near } \infty, \\
    &\  (1 + d) \cdot \sum_{j=1}^q |\trace((\partial_j g)^2)| + |\trace(\partial_j g)^2| + |\trace(\partial_j^2 g)| < \scal(g)\}.
    \end{split}
 \end{align*}
 The latter space includes into the loop space modelled by continuous block maps, which is, of course, weakly homotopy equivalent to the classical loop space.
 Approximating continuous maps by smooth maps and a re-scaling argument show that the inclusions $\Omega^{q,\infty}F_p^{slow} \hookrightarrow \Omega^q F_p^{slow}$ and $F_p^{slow} \hookrightarrow F_p$ are weak homotopy equivalences, see Lemma \ref{FiltrationModels - Lemma} below.
 
 Define 
 \begin{equation*}
     \susp \colon \Omega^{q,\infty} F_p^{(slow)} \rightarrow \Omega^{q-1} F_{p+1}^{(slow)}
 \end{equation*}
 via
 \begin{equation*}
     \susp(g)(t_1,\dots,t_{q-1})_{(m,x_1,\dots,x_{p+2})} \deff g(x_{p+2},t_1,\dots,t_{q-1})_{(m,x_1,\dots,x_{p+1})} + \diff x_{p+2}^2.
 \end{equation*}
 These maps are continuous by the following argument: 
 Recall the adjunction between functions and sections
 \begin{equation*}
     \mathcal{C}^\infty(N;\Riem^+(M \times I^p)) \cong \Gamma^\infty(M \times I^n\times N; \mathrm{Pos}(T(M\times I^p)) \times N),
 \end{equation*}
 where $N$ is a smooth compact manifold with faces and $\mathrm{Pos}(T(M\times I^p))$ is the bundle of all positive definite, symmetric bilinear forms on $T(M\times I^p)$ over $M \times I^p$.
 Under this adjunction, the suspension is induced by the following smooth map of fibre bundles
 \begin{equation}\label{Eq: SuspInducingBundleMap}
     \xymatrix{\mathrm{Pos}(TM \times I^p) \times I \times I^{q-1} \ar[rrr]^-{\pr_1^\ast(\placeholder) + \pr_2^\ast\euclmetric} 
    &&& \mathrm{Pos}(T(M \times I^p \times I)) \times I^{q-1},  }
 \end{equation}
 via post-composition. 
 Thus, it induces a continuous map $\Omega^{q,\infty} F_p \rightarrow \Omega^{q-1} F_{p+1}$.
 In fact, this map factors through the space of smooth loops.
 
 The map $\iota_{p,q}$ is now defined as
 \begin{equation*}
     \xymatrix{\iota_{p,q}\colon \pi_q(F_p) = \pi_0(\Omega^{q,\infty}F_p^{(slow)}) \ar[rr]^-{\pi_0(\susp)} && \pi_0(\Omega^{q-1}F_{p+1}^{(slow)}) = \pi_{q-1}(F_{p+1}),  }
 \end{equation*}
 and $\iota_{n,0} \colon \pi_0(F_{n}^{slow}) \rightarrow \pi_n(\hofib)$ is the surjective map that sends the isotopy class of a block metric to its concordance class.
 
 Clearly, 
 \begin{align*}
     \begin{split}
         \iota_{p+k,q-k} \circ \dots \circ \iota_{p,q} &(g)(t_1,\dots,t_{q-k-1})_{(m,x_1,\dots,x_{p+k+2})} \\
         &= g(x_{p+2},\dots, x_{p+k+2},t_1,\dots,t_{q-k-1})_{(m,x_1,\dots,x_{p+1})} + \sum_{j=p+2}^{p+k+2} \diff x_j^2
     \end{split}
 \end{align*}
 for $k<q$ and all metrics for which the composition is defined.
 Thus, the composition $\pi_q(F_p) \rightarrow \pi_{p+q}(\hofib)$ sends a block map to the concordance class of its suspension.
 
 In remains to prove that the image of the suspension $\pi_{n-p}(F_p) \rightarrow \pi_n(\hofib)$ is the set of concordance classes of $F_p \Omega_{g_0}^n\hofib$.
 Consider the map
 \begin{align*}
     \mathfrak{pr} \colon F_p\Omega_{g_0}^n\hofib &\rightarrow \Omega^{n-p,\infty}F_p^{slow}, \\
     g & \mapsto \left(y \mapsto g_{(\cdot,\cdot,y)}\restrict_{T(M\times \R^{p+1})} \right).
 \end{align*}
 If we turn the set $F_p\Omega^n_{g_0}\hofib$ into a topological space by endowing it with the colimit topology, then $\mathfrak{pr}$ is continuous and the inverse of $\susp\colon \Omega^{n-p,\infty}F_p^{slow} \rightarrow F_p\Omega^{n}_{g_0}\hofib$.
 
 In particular, the suspension maps the group $\pi_{n-p}(F_p)$ onto the isotopy classes of $F_p\Omega^n_{g_0}\hofib$ and thus onto the concordance classes of $F_p\Omega^n_{g_0}\hofib$.
\end{proof}

\begin{lemma}\label{FiltrationModels - Lemma}
 Let $\Omega^{q,\infty}F_p$ and $F_p^{slow}$ be as in the proof of Lemma \ref{SuspensionFiltration - Lemma}. Then the inclusions
 \begin{align*}
    F_p^{slow} \hookrightarrow F_p  \qquad \text{ and } \qquad \Omega^{q,\infty}F_p^{(slow)} \hookrightarrow \Omega^q F_p^{(slow)}
 \end{align*}
 are weak homotopy equivalences.
\end{lemma}
\begin{proof}
 To prove the first statement, we will show the equivalent statement that every map $(D^n,S^{n-1}) \rightarrow (F_p,F_p^{slow})$ can be homotopied through such maps to a map $D^n\rightarrow F_p^{slow}$.
 
 Given a such map $g \colon (D^n,S^{n-1}) \rightarrow (F_p, F^{slow}_p)$. 
 As $F_p$ is a sequential colimit of embeddings between Hausdorff spaces, we can find a sufficiently large $\rho$ such that each $g(t)$ decomposes on $M \times \rho I^{p+1}$.
    
 By assumption, each $g(t)$ decomposes 
 \begin{equation*}
     g(t)|_{\{x_1 > \rho \}} = g(t)^{TM} + \sum_{j=1}^{p+1} \diff x_j^2,
 \end{equation*}
 where $g^{TM}$ is the induced metric on $TM \subseteq T(M\times \R^{p+1})$.
 The block form of $g(t)$ implies its independence of the coordinate $x_1$.
 Thus, there is a unique map $\gamma \colon D^n \times \rho I^p \rightarrow \Riem^+(M)$ such that 
 $g(t)= \susp(\gamma(t,\placeholder)) $ on $\{x_1 > \rho\}$.
 The map $\gamma$ is continuous in the $D^n$-variable and smooth in the $\rho I^p$-variable.
 
 Recall that the defining conditions of $A_p$ are algebraic expressions of the $2$-jet of $g$.
 By compactness, we find therefore a sufficiently large $R>1$ such that the map $x \mapsto \gamma(t,x/R)$ is an element of $A_p$ for all $t \in D^n$.
 
 Pick a sufficiently large $L = L(R) >0$ and a smooth function $\Phi \colon [\rho, \rho + L] \rightarrow [1,R]$ that is monotonically increasing and satisfies 
    \begin{align*}
        & \Phi \equiv 1 \text{ near } \rho, &\Phi \equiv R \text{ near } \rho + L, \\
        &0 \leq \Phi' \ll 1, &|\Phi''| \ll 1.
    \end{align*}
    For each $s \in [0,1]$, set $\Phi_s \deff (1-s)\cdot 1 + s\Phi$.
    Since $\gamma|_{\partial \rho I^p} = \const_{g_0}$, we can extend it via $g_0$ to a block map $\gamma \colon \R^n \rightarrow \Riem^+(M)$.
    The map 
    \begin{align*}
        \gamma_{\aux}^{s} \colon D^n \times [\rho,\rho + L] \times \rho I^p &\rightarrow \Riem^+(M), \\ 
        (t,x_1,\dots,x_{p+1}) &\mapsto \gamma\bigl(t, 1/\Phi_{s}(x_1) \cdot (x_2,\dots,x_{p+1})\bigr)
    \end{align*}
    can be extended to a family of block maps 
    \begin{equation*}
        \gamma_{\aux}^s \colon D^n \times \R^{p+1} \rightarrow \Riem^+(M).
    \end{equation*}
    
    Consider the map
    \begin{align*}
        \mathcal{H}\colon D^n \times [0,1] \rightarrow F_p \quad \text{given by} \quad (t,s) \mapsto g(t)|_{\{x_1 \leq \rho\}} \cup \susp(\gamma_{\aux}^s)|_{\{x_1\geq \rho\}}.
    \end{align*}
    This map is well-defined because $g(t) = \susp(\gamma_{\aux}^s)$ near $\{x_1 = \rho\}$ and $\scal(\gamma_{\aux}^s)>0$ provided $|\Phi'|$ and $|\Phi''|$ are sufficiently small, see Lemma \ref{susp(gammaaux)pos - Lemma} below. 
    It satisfies $\mathcal{H}(\placeholder,0) = \gamma$ and $\mathcal{H}(\placeholder,1)$ takes values in $F^{slow}_p$.
    Since $|(1/\Phi)'|\leq 1$ and $|(1/\Phi)''|\leq 1$, the map $\mathcal{H}$ restricts to a map $S^n \times [0,1] \rightarrow F_p^{slow}$.
    
    Thus, $\mathcal{H}$ is a homotopy of pairs between $g$ and a map with values in $F_p^{slow}$, so the first statement is proven. \\ $ $
    
     In order to show that $\Omega^{q,\infty}F_p^{(slow)} \hookrightarrow \Omega^q F_p^{(slow)}$ is a weak homotopy equivalence we consider 
 \begin{equation*}
    \Bar{\Omega}^{q,\infty}F_p^{(slow)} \deff \{g \colon \R^n \rightarrow F_p^{(slow)} \, : \, g \text{ smooth block map}, g = g_0 \text{ near }\infty\},
 \end{equation*}
 the space of all smooth block maps into $F_p^{(slow)}$ without ``speed constrains''.
 Approximation theory, see Section \ref{Subsection - The Topological Space of Block Maps} for details, yields that the inclusion $\Bar{\Omega}^{q,\infty}F_p^{(slow)} \hookrightarrow \Omega^q F_p^{(slow)}$ is a weak homotopy equivalence.

 We will show that the inclusion $\Omega^{q,\infty}F_p^{(slow)} \hookrightarrow \Bar{\Omega}^{q,\infty}F_p^{(slow)}$ is a homotopy equivalence.

 Consider the map $R \colon \bar{\Omega}^{q,\infty}F_p^{(slow)} \rightarrow \R_{\geq 1} $ given by
 \begin{equation*}
     R(g) \deff \mathrm{\sup}_{M \times \R^{p+1} \times \R^q}\left\{1 + (1+d)\sum_{j=1}^q |\trace((\partial_j g)^2)| + |\trace(\partial_j g)^2| + |\trace(\partial_j^2 g)|\right\}.
 \end{equation*}
 The values of $R$ are always finite because we only plug in block maps that take values in block metrics.
 The map $R$ is also continuous because it is an algebraic expression of the $2$-jet of $g$ considered now as section over the space $M\times \R^{p+1} \times \R^q$.
 
 The same holds for the map $\mathtt{m} \colon \bar{\Omega}^{q,\infty}F_p^{(slow)} \rightarrow (0,1]$ given by
 \begin{equation*}
    g \mapsto \frac{1}{2d} \inf\left(\bigl\{\left(\scal(g(t))(m,x)\right)^{1/2} \, : \, (m,x,t) \in M \times \R^{p+1} \times \R^q\bigr\}\right). 
 \end{equation*}
 
 The inclusion $\varphi \colon \Omega^{q,\infty}F_p^{(slow)} \hookrightarrow \bar{\Omega}^{q,\infty}F_p^{(slow)}$ is a homotopy equivalence.
 Its homotopy inverse is given by $\psi \colon g \mapsto g \circ \left(\mathtt{m}/(1+\mathtt{m}^2)^{1/2} \cdot R^{-1}\right)$.
 The homotopy between the identity and $\varphi \circ \psi$, respectively $\psi \circ \varphi$, is given by 
 \begin{equation*}
     (g,s) \mapsto g \circ \left(s\cdot \mathtt{m}/(1+\mathtt{m}^2)^{1/2} \cdot R^{-1} + (1-s)\cdot 1\right). 
 \end{equation*}
 because if $g \in A_n$, then also all re-parametrisation of $g$ that are ``slower'' than $g$ are also elements of $A_n$. 
\end{proof}

\begin{lemma}\label{susp(gammaaux)pos - Lemma}
 The suspension of $\gamma_{\aux}^s$ defined in the proof of Lemma \ref{FiltrationModels - Lemma} has positive scalar curvature.
\end{lemma}
\begin{proof}
  We know from 
  \begin{align*}
      0 < \scal(\susp(\gamma)) &= \scal(\gamma) + \sum_{j=2}^{p+1} \frac{3}{4} \trace((\partial_j\gamma^{op})^2) - \frac{1}{4} \trace(\partial_j \gamma)^2 - \trace(\partial^2_j \gamma) \\
      &=: \scal(g) + \mathfrak{err}(g)
  \end{align*}
  that $\mathfrak{err}(\gamma) > - \scal(\gamma)$.
  The block form of the map $g$ actually implies $\mathfrak{err}(\gamma) > - \scal(\gamma) + \eps$ for some $\eps > 0$.
  The relevant differentials of $\gamma_{\aux}^s$ are given by
  \begin{equation*}
      \partial_j \gamma_{\aux}^s = \begin{cases}
        \frac{1}{\Phi_s(x_1)} \cdot (\partial_j \gamma)(1/\Phi_s(x_1)\cdot (x_2,\dots,x_{p+1}), & \text{ if } j > 1, \\
        \frac{-\Phi_s'(x_1)}{\Phi_s(x_1)} \cdot \sum_{k=2}^{p+1} (\partial_k \gamma)(1/\Phi_s(x_1) \cdot (x_2,\dots,x_{p+1}))x_k, & \text{ if } j=1,
      \end{cases}
  \end{equation*}
  and
  \begin{align*}
      \partial_j^2 \gamma_{\aux}^s = \frac{1}{\Phi_s(x_1)^2} \cdot (\partial_j^2 \gamma)(1/\Phi_s(x_1)\cdot (x_2,\dots,x_{p+1})), \qquad \text{ if } j>1.
  \end{align*}
  The last partial differential is given by
  \begin{align*}
       \begin{split}
         \partial_1^2 \gamma_{\aux}^s = & \frac{(\Phi_s')^2}{\Phi_s^4} \cdot \sum_{k,l=2}^{p+1} (\partial_k \partial_l \gamma)(1/\Phi_s(x_1)\cdot(x_2,\dots, x_{p+1})) x_k x_l \\
          & \quad + \frac{-\Phi_s'' \Phi_s + 2(\Phi_s')^2}{\Phi_s^4} \cdot \sum_{k=2}^{p+1} (\partial_k \gamma)(1/\Phi_s(x_1)\cdot(x_2,\dots, x_{p+1})) x_k.
      \end{split}
  \end{align*}
  Since $\gamma(t,\placeholder)$ is constant outside of $\rho I^p$ for all $t \in D^n$, the map $\gamma_{\aux}^s$ is constant outside of $\rho R\cdot I^p$.
  Thus, by increasing $L(R)$, we can make $|\partial_1 \gamma_{\aux}^s|$ and $|\partial_1^2\gamma_{\aux}^s|$ arbitrary small (by choosing $\Phi$ with sufficiently small derivatives).
  
  This, in particular, implies we can choose $L$ and $\Phi$ such that
  \begin{align*}
      \mathfrak{err}(\gamma_{\aux}^s) &= \Phi_s(x_1)^{-2}\cdot\mathfrak{err}(\gamma) \circ 1/\Phi_s + \frac{3}{4} \trace((\partial_1\gamma_{\aux}^{s,op})^2) - \frac{1}{4} \trace(\partial_1 \gamma_{\aux}^s)^2 - \trace(\partial^2_1 \gamma_{\aux}^s)^2 \\
      \begin{split}
         &> - \Phi_s(x_1)^{-2}\cdot\scal(\gamma) \circ 1/\Phi_s + \eps  \\
         & \quad + \, \frac{3}{4} \trace((\partial_1\gamma_{\aux}^{s,op})^2) - \frac{1}{4} \trace(\partial_1 \gamma_{\aux}^s)^2 - \trace(\partial^2_1 \gamma_{\aux}^s)^2 
      \end{split}\\
      &> -\Phi_s(x_1)^{-2}\cdot\scal(\gamma_{\aux}^s) \circ 1/\Phi_s = \scal(\gamma_{\aux}^s),
  \end{align*}
  and hence so that $\scal(\susp(\gamma_{\aux}^s))$ is positive.
\end{proof}

The images of $F_p\Omega^n_{g_0}\hofib$ in $\pi_n(\hofib)$ form a filtration of subsets.
We will see later that this filtration is a filtration of subgroups and that the quotients of the filtration groups and its pre-successors agree with the infinity page of a spectral sequence.
To set up the spectral sequence, we first need to realise the maps $\iota_{p,q}$ as the boundary operator of a long exact sequence associated to a fibration.
The next definition introduces the total space of this fibration.
\begin{definition}
 For $p\geq 0$, let $G_p$ be the topological space whose underlying set is given by
 \begin{align*}
     \begin{split}G_p \deff \Bigl\{ g \in \ConcSet[ ][M \times \R^{p+1}] \, : \, &\face[1]{1}g \in \SingMet[M][p], \\ 
     &\face[\omega]{j}g = g_0 \text{ for } (j,\omega) \neq (1,1), \,(p+1,1) \Bigr\}
     \end{split}
 \end{align*}
 and whose topology is the (subspace topology of the) colimit topology coming from the inclusions $\ConcSet[ ][M \times RI^{p+1}] \hookrightarrow \ConcSet[ ][M \times SI^{p+1}]$.
 We further define $G_p^{slow} \deff \{g \in G_p \, : \, \face[1]{1}g \in A_p\}$.
\end{definition}
Note that the proof of Lemma \ref{FiltrationModels - Lemma} also implies that the inclusion $G_p^{slow} \hookrightarrow G_p$ is a weak homotopy equivalence.

There is a canonical map $\restrict \colon G_{p+1}^{(slow)} \rightarrow F_p^{(slow)}$ given by $g \mapsto \face[1]{p+1}g$. 
Identifying this map as a Serre fibration is probably difficult as we would have to deal with arbitrary continuous maps.
Luckily, the following weakening is easier to verify in our setup and is still good enough for our purpose.
\begin{definition}
 A smooth map between two smooth manifolds or  between two open subsets of affine Fréchet spaces $p \colon E \rightarrow B$ is a \emph{fibration in the smooth category}, if for all given solid outer squares
 \begin{equation*}
     \xymatrix{D^n \times \{0\} \ar[rr]^-f \ar@{^{(}->}[d] && E \ar[d]^p \\ 
     D^n \times [0,1] \ar[rr]^-h \ar@{.>}[rru]^{\mathcal{H}} && B}
 \end{equation*}
 with $f$ and $h$ smooth and $h$ stationary near the boundary of $[0,1]$, there exists a smooth map $\mathcal{H}$ that is stationary near the boundary of $[0,1]$ and makes the diagram commutative.
\end{definition}
The next lemma is no surprise as continuous maps can be approximated by smooth ones relatively to closed subsets on which they are already smooth. 
With the help of such approximations, the next lemma will be proved as the analogous statement for Serre fibrations.
\begin{lemma}\label{SmoothFib are HomotopyFib}
 Let $p \colon E \rightarrow B$ be a fibration in the smooth category.
 Then $p$ is a homotopy fibration, that is, the comparison map from each fibre into the homotopy fibre is a weak homotopy equivalence for all path components of $B$.
\end{lemma}
\begin{proof}
    Pick a point $b_0 \in B$ and denote its fibre by $F \deff p^{-1}(\{b_0\})$.
    The fibre comparison map 
    \begin{align*}
        j \colon F \hookrightarrow \mathrm{hofib}(f;b_0) \qquad \text{ given by } \qquad f \mapsto (f, \const_{b_0})
    \end{align*}
    is a weak homotopy equivalence if we find, for each $\varphi \colon (D^n, S^{n-1}) \rightarrow (\mathrm{hofib}(f;b_0),F)$, a homotopy of pairs $H\colon [0,1] \times (D^n, S^{n-1}) \rightarrow (\mathrm{hofib}(f;b_0),F)$ that deforms $\varphi$ into a map whose target is $F$.
    
    Each map $\varphi\colon (D^n,S^{n-1}) \rightarrow (\mathrm{hofib}(p;b_0),F)$ is given by a commutative diagram
    \begin{equation*}
        \xymatrix{ D^n \times \{1\} \ar[rr]^-{\varphi^{(1)}} \ar@{^{(}->}[d] && E \ar[d]^p \\
                   D^n \times [0,1] \ar[rr]^-{\varphi^{(2)}} && B}
    \end{equation*}
    so that $\varphi^{(2)}$ restricts on $S^{n-1}\times [0,1] \cup D^n\times\{0\}$ to the constant map with value $b_0$.
    By approximation theory, there is a homotopy $\varphi^{(1)}_t$ such that $\varphi^{(1)}_0 = \varphi^{(1)}$ and $\varphi^{(1)}_1$ is a smooth map.
    The homotopy can be chosen to be stationary near the boundary.
    We modify this homotopy in three steps to a homotopy of diagrams that is stationary near $D^{n} \times \{0,1\}$.
    
    First, we re-parameterise $\varphi^{(2)} = \varphi_{-1}^{(2)}$ such that the result $\varphi_{0}^{(2)}$ is stationary near the boundary.
    We can choose the interpolating homotopy $\varphi^{(2)}_t$ with $t \in [-1,0]$ to be relative to the boundary.
    Next, we define, for all $t \in [0,1]$, the homotopy
    \begin{align*}
        \varphi^{(2)}_t(x,s) = \begin{cases}
            \varphi^{(2)}(x,(1+t)s), & s \in [0,(1+t)^{-1}], \\
            p(\varphi^{(1)}_{(1+t)s-1}(x)), & s \in [(1+t)^{-1},1].
          \end{cases}
    \end{align*}
    The map is $\varphi^{(2)}_1$ is not smooth yet, but it smooth near $D^n \times \{0,1\}$ and constant on the boundary $D^n\times \{0\} \cup S^{n-1} \times [0,1]$.
    Thus, as the last step, we can homotopy it relative to the boundary to a smooth map $\varphi^{(2)}_2$ that is stationary near $D^n\times \{0,1\}$ using approximation theory.
    
    This diagram homotopy restricts to the diagram homotopy
    \begin{equation*}
        \xymatrix{ S^{n-1} \times \{1\} \ar[rr]^{\varphi^{(1)}_t} \ar@{^{(}->}[d] && F \ar[d]^p \\
                   S^{n-1} \times [0,1] \ar[rr]^{\varphi^{(2)}_t} && b_0}
    \end{equation*}
    
    In conclusion, each map $\varphi\colon (D^n,S^{n-1}) \rightarrow (\mathrm{hofib}(p;b_0),F)$ can be homotopied through map of pairs to a smooth map $\psi$, so that $\psi^{(2)}$ are stationary near $D^n \times \partial [0,1]$.
    
    To ease the notation, assume that $\varphi$ has this property in the first place.
    As $p \colon E \rightarrow B$ is a fibration in the smooth category, we find a lift $\Phi$ with restriction
   \begin{equation*}
        \xymatrix{ D^n \times \{1\} \ar[rr]^{\varphi^{(1)}} \ar@{^{(}->}[d] && E \ar[d]^p && S^{n-1}\times \{1\} \ar@{^{(}->}[d] \ar[rr]^-{\varphi^{(1)}|_{S^{n-1}}} && F \ar[d]  \\
                   D^n \times [0,1] \ar[rr]^{\varphi^{(2)}} \ar@{-->}[rru]^{\Phi} && B, && S^{n-1} \times [0,1] \ar[rr]^-{\varphi^{(2)}|_{S^{n-1}}} \ar@{-->}[rru]^{\Phi|_{S^{n-1}}} && b_0.}
    \end{equation*}
    We use $\Phi$ to define $\mathcal{H} \colon [0,1] \times D^n \times [0,1] \rightarrow E$ as $\mathcal{H} \colon (t,x,s) \mapsto \Phi(x,ts).$
        
    The desired homotopy $H \colon [0,1] \times (D^n,S^{n-1}) \rightarrow (\mathrm{hofib}(f,b_0),F)$ is the continuous map that corresponds to
    \begin{equation*}
        \xymatrix{[0,1]\times D^n \times \{1\} \ar@{^{(}->}[d] \ar[rr]^-{\mathcal{H}} && E \ar[d]^-p \\
        [0,1] \times D^n \times [0,1] \ar[rr]^-{p \circ \mathcal{H}} && B,}
    \end{equation*}
    so the lemma is proven.
\end{proof}
\begin{proposition}\label{Restriction is Smooth Fib}
 For all $p\geq 0$, the restriction to the last front face $\restrict \colon G_{p+1} \rightarrow F_p$ is well-defined and continuous.
 Furthermore, the sequence
 \begin{equation*}
     \xymatrix@C+1em{F_{p+1} \ar@{^{(}->}[r]^{\mathrm{incl}} & G_{p+1} \ar[r]^{\restrict} & F_p }
 \end{equation*}
 is a fibration in the smooth category.
\end{proposition}
\begin{proof}
  The cubical identities imply that $\restrict = \face[1]{p+2}$ maps $G_{p+1}$ to $F_p$.
  
  The restriction to the last front face is a continuous map $\ConcSet[ ][M \times RI^{p+2}] \rightarrow \ConcSet[ ][M \times RI^{p+1}]$ as it is induced by a pull back of a vector bundle map.
  Since these restriction maps are compatible with the inclusions $\ConcSet[ ][M \times RI^q] \hookrightarrow \ConcSet[ ][M \times SI^q]$, they induce a map $G_{p+1} \rightarrow F_p$, which agrees with $\restrict$.
  
  Given a commutative diagram of smooth maps
   \begin{equation*}
     \xymatrix{D^n \times \{0\} \ar[rr]^-f \ar@{^{(}->}[d] && G_{p+1} \ar[d]^\restrict \\ 
     D^n \times [0,1] \ar[rr]^-h  && F_p}
 \end{equation*}
 such that $h$ is stationary near $D^n \times \partial [0,1]$.
  By Lemma \ref{Relative T1}, the maps $f$ and $h$ factor through $\ConcSet[ ][M \times \rho I^{p+2}]$ and $\ConcSet[ ][M \times \rho I^{p+1}]$, respectively, as $G_{p+1}$ and $F_p$ are colimits of relative $T_1$ inclusions.
 Let $\Phi \colon D^n \times \partial [0,1]^2 \rightarrow F_p$ be given by
 \begin{align*}
     & h(\placeholder,0), \text{ on } D^n \times 0 \times [0,1], & h, \text{ on } D^n \times 1 \times [0,1], \\
     & h(\placeholder,0), \text{ on } D^n \times [0,1] \times 0, & h, \text{ on } D^n \times [0,1] \times 1 .
 \end{align*}
 Since $h$ is stationary near $\partial [0,1]$, we can extend $\Phi$ to a smooth map on $D^n \times I^2$ that is independent of the normal variable of each face of $I^2$ near $\partial I^2$. 
 Let $\Psi(x,t)$ be the block map extending $[\rho, \rho + 1] \ni s \mapsto \Phi(x,s-\rho,t)$ by elongation and let $\mathcal{H} \colon D^n \times [0,1] \rightarrow G_{p+1}$ be given by 
 \begin{equation*}
     \mathcal{H}(x,t) \deff f \cup \susp \, \Phi\left(x,\frac{(\placeholder) -\rho}{R},t\right) =: f|_{\{x_{p+2} \leq \rho\}} \cup \susp \left({}_R\Psi(x,t)\right)|_{\{x_{p+2} \geq \rho\}}.
 \end{equation*}
 If $R$ is so large that
 \begin{equation*}
     |\mathrm{tr}((\partial_s{}_R\Psi)^2)| + |\mathrm{tr}(\partial_s{}_R\Psi)^2| + |\mathrm{tr}(\partial_s^2{}_R\Psi)| < \scal({}_R\Psi(x,t))(s),
 \end{equation*}
 then $\susp({}_R\Psi)(x,t)$ is a psc metric on $M \times \R^{p+1} \times \R$, which also agrees with $f(x)$ near $M \times \R^{p+1} \times \{\rho\}$.
 Thus, $\mathcal{H} \colon D^n \times [0,1] \rightarrow G_{p+1}$ is well-defined. 
 The map $\mathcal{H}$ is smooth because $\Phi$ smooth and $\susp$ is a smooth map between smooth Fréchet spaces as $\susp$ is induced by vector bundle maps, see Equation \ref{Eq: SuspInducingBundleMap}.
 From ${}_R\Psi(x,t)(\rho+R) = h(x,t)$ follows that $\restrict(\mathcal{H}(x,t)) = h(x,t)$, so $\mathcal{H}$ lifts the homotopy $h$.
 The lift $\mathcal{H}$ is also stationary near $\partial[0,1]$ because $\Phi$ is stationary near the corner $(1,1) \in [0,1]^2$.
 Thus, the restriction to the last front face $\restrict$ is a fibration in the smooth category.
\end{proof}
Lemma \ref{SmoothFib are HomotopyFib} and Proposition \ref{Restriction is Smooth Fib} yield a long exact sequence of homotopy groups.
We would like to relate the connecting homomorphisms $\partial_q \colon \pi_q(F_p) \rightarrow \pi_{q-1}(F_{p+1})$ to $\iota_{p,q}$.
Before we formulate the next lemma, let us recall that every fibration $F \hookrightarrow E \rightarrow B$ induces a long exact sequence on homotopy groups via the Puppe sequence, see for example \cite{davis2001lecture}*{Theorem 6.42}.
The boundary map $\partial_1 \colon \pi_1(B) \rightarrow \pi_0(F)$, in this picture, is given by evaluating the fibre transport at the base point of the fibre.
With this in mind, we denote by $(-1)\partial_1$ the map that evaluates the inverse fibre transport at the base point.
\begin{lemma}\label{Boundary Comparison - Lemma}
 For all $p\geq 0$ and $q \geq 1$, we have
 \begin{equation*}
    \iota_{p,q} = (-1)^{q} \partial_q \colon \pi_q(F_p) \rightarrow \pi_{q-1}(F_{p+1}). 
 \end{equation*}
\end{lemma}
\begin{proof}
  We consider the case $q=1$ first.
  Recall that $\partial_1 \colon \pi_1(F_p) \rightarrow \pi_0(F_{p+1})$ comes from the Puppe sequence.
  Each map $d$, possibly only defined on a weakly equivalent subspace, that makes the diagram
  \begin{equation*}
      \xymatrix{ \Omega F_p \ar[d]_{\simeq_w} \ar@{-->}[r]^d & F_{p+1} \ar[d]^{\simeq_w} \\ 
      \hofib[ ] \ar[r] & \hofib[ ](\restrict)  }
  \end{equation*}
  homotopy commutative gives $\partial_1 = \pi_0(d)$.
  Here, $\hofib[ ]$ is the homotopy fibre of the canonical map $\hofib[ ](\restrict) \rightarrow G_{p+1}$.
  The homotopy fibres are defined using block maps, so for example,
  \begin{equation*}
      \begin{split}\hofib[ ](\restrict) = \{(g,\gamma) \, | \, &g \in G_{p+1}, \  \gamma \colon \R \rightarrow F_p \text{ block map}, \\ 
      & \lim_{t\to -\infty}\gamma(t) = g_0, \lim_{t \to \infty} \gamma(t) = \restrict(g)\}.
      \end{split}
  \end{equation*}
  
  It is straighfoward to check that the composition
  \begin{equation*}
     \xymatrix{ \delta \colon \Omega^{1,\infty}F_p  \ar@{^{(}->}[r]^-\simeq & \Omega F_p \ar[r] & \hofib[ ] \ar[r] & \hofib[ ](\restrict) }
  \end{equation*}
  is given by
  \begin{equation*}
      (t \mapsto g(t)) \mapsto (g_0,t \mapsto g(t))
  \end{equation*}
  and that the composition 
  \begin{equation*}
     \xymatrix{ \susp \colon \Omega^{1,\infty}F_p  \ar@{^{(}->}[rr]^-\susp && F_{p+1} \ar[r]^-{\simeq} & \hofib[ ](\restrict) }
  \end{equation*}
  is given by
  \begin{equation*}
      (t \mapsto g(t)) \mapsto (\susp(g),\const_{g_0}).
  \end{equation*}
  
  We will construct a path between $\delta(\Bar{g})$ and $\susp(g)$ for all $g \in \Omega^{1,\infty}F_p$, where $\Bar{\cdot}$ is the standard involution on the loop space.
  This will then imply the claim for $q=1$.
  
  Choose a $\rho >0$ such that the smooth block map $g \colon \R \rightarrow F_p$ is constant (with value $g_0$) outside of $M \times [-\rho,\rho]$.
  Let $\Phi \colon [-\rho, \rho]\times [-1,1] \rightarrow [-\rho,\rho]$ be a smooth map that satisfies
  \begin{align*}
      \Phi(\placeholder,-1) = -\rho \quad \text{ and } \quad \Phi(\placeholder, s) = \begin{cases}
         -\rho, & \text{ near } \R_{\leq -\rho}, \\
         s\cdot \rho, & \text{ near } \R_{\geq \rho}.
      \end{cases}
  \end{align*}
  Fix $R\geq 1$ such that $R^{-1}\left( |\partial_t\Phi(t,s)| + |\partial_t^2\Phi(t,s)|\right) < 1$ for all $(t,s) \in [-\rho,\rho]\times[-1,1]$ and define the path
  \begin{gather*}
      H_g \colon [-1,1] \rightarrow \hofib[ ](\restrict), \\
      s \mapsto [\susp(t \mapsto g(\Phi(t/R,s))), \, t \mapsto g(-\Phi(t,-s)) ].
  \end{gather*}
  
  We need to show that this path is a well-defined map.
  
  First, we show that 
  the scalar curvature of $\susp(t \mapsto g(\Phi(t/R,s)))$ is positive for all $s \in [-\rho,\rho]$.
   We abbreviate $\Phi(\placeholder/R,s)$ to ${}_R\Phi$ as $s$ is not important in the calculation.
 The chain rule gives
 \begin{align*}
     & \quad \, \left| \scal(\susp(g \circ {}_R\Phi)) - \scal(g\circ {}_R\Phi) \right| \\
     \begin{split}
         &= \biggl| \frac{3}{4}\trace((\partial_t g^{op})^2)\cdot R^{-2}{}_R(\partial_t\Phi)^2 
     -\frac{1}{4} \trace(\partial_t g)^2 \cdot R^{-2}{}_R(\partial_t\Phi)^2 \\
     & \qquad - \trace(\partial_t^2 g)\cdot R^{-2}{}_R(\partial_t\Phi)^2 
     - \trace(\partial_t g)\cdot R^{-2}{}_R(\partial_t^2\Phi)\biggr|
     \end{split} \\
     &< (1 + d) \cdot \sum_{j=1}^q |\trace((\partial_j g)^2)| + |\trace(\partial_j g)^2| + |\trace(\partial_j^2 g)|,
 \end{align*}
 provided $R$ is sufficiently large.
 This implies $\scal(\susp(g\circ {}_R\Phi)) > 0$.
  
  From $\susp(t \mapsto g(\Phi(t/R,s))) = g(s\rho) + \diff t^2$ on $\{t \geq R\}$ we conclude $\restrict \susp(t \mapsto g(\Phi(t/R,s))) = g(s\rho)$, which agrees with
  \begin{align*}
      \lim_{t \to \infty} g(-\Phi(t/R,-s)) = g(-\Phi(\rho R/R,-s)) = g(s\rho).
  \end{align*}
  Thus, $H_g$ takes values in $\hofib[ ](\restrict)$. It is also smooth, as $\Phi$ is a smooth map.
  
  The elements $H_g(1)$ and $\susp(g)$ are connected through the two paths 
  \begin{equation*}
      [0,1] \ni \lambda \mapsto [\susp(g(\lambda \cdot(\placeholder/R) + (1-\lambda)\Phi(\placeholder /R,1))), \const_{g_0}],
  \end{equation*}
  and
  \begin{equation*}
      [1,R] \ni r \mapsto [\susp(g(\placeholder /r)),\const_{g_0}],
  \end{equation*}
  while $\delta(\overline{g})$ and $H_g(-1)$ are connected through the path 
  \begin{equation*}
      \lambda \mapsto [ g_0, g(-\lambda \id - (1-\lambda)\Phi(\placeholder/R,1)) ].
  \end{equation*}
  Thus, the elements $\susp(g)$ and $\delta(\Bar{g})$ lie in the same path-component, which implies the claim for the case $q=1$.
  
  The general case $q\geq 1$ can be deduced from the special one as follows.
  By definition, $\partial_q = \pi_0(\Omega^{q-1}\delta)$, so it is enough to compare the two maps
  \begin{equation*}
    \Omega^{q-1}\susp_1, \, \susp_q \colon \Omega^{q,\infty}F_p \rightarrow \Omega^{q-1}F_{p+1}.  
  \end{equation*}
  The first map is defined as
  \begin{equation*}
    \bigl((t_1,\dots,t_q) \mapsto g(t_1,\dots,t_q)\bigr) \mapsto (t_1,\dots,t_{q-1}) \mapsto \Omega^{q-1}\susp(g)(t_1,\dots,t_{q-1}),  
  \end{equation*}
  where 
  \begin{equation*}
      \Omega^{q-1}\susp(g)(t_1,\dots,t_{q-1})_{(m,x_1,\dots,x_{p+2})} = g(t_1,\dots,t_{q-1},x_{p+2})_{(m,x_1,\dots,x_{p+1})} + \diff x_{p+2}^2
  \end{equation*}
  while the second map $\susp_q$ was defined as
  \begin{equation*}
      \susp_q(g)(t_1,\dots,t_{q-1})_{(m,x_1,\dots,x_{p+2})} = g(x_{p+2},t_1,\dots, t_{q-1})_{(m,x_1,\dots,x_{p+1})} + \diff x_{p+2}^2.
  \end{equation*}
  Thus, 
  \begin{equation*}
      \susp_q(g) = \Omega^{q-1} \susp_1 \circ (1 \ 2 \ \dots \ q),
  \end{equation*}
  where $(1 \ 2 \ \dots \ q)(g)(t_1,\dots,t_q) = g(t_q,t_1,\dots, t_{q-1})$.
  
  This identity now implies the claim via 
  \begin{align*}
      \iota_{p,q} &= \pi_0(\susp_q) = \pi_0(\Omega^{q-1}\susp_1 \circ (1 \ 2 \ \dots \ q)) \\
      &=(-1)^{q-1}\pi_0(\Omega^{q-1}\susp_1) = (-1)^{q-1} \pi_0(\delta \circ \Bar{\cdot}) \\
      &= (-1)^q \pi_0(\Omega^{q-1}\delta) = (-1)^q \partial_q.
  \end{align*}
\end{proof}

If $p \geq 1$, then every element $g \in G_{p+1}$ or $g\in F_p$ restricts to the base point $g_0$ on at least one pair of opposite faces.
This allows us to glue two block metrics together along their common faces. 
Sadly, the gluing construction does not give an $H$-space structure because of the following reason:
The cubes outside of those the metrics decompose are of unknown length, and the function that assigns to a block metric the diameter of its smallest cube outside it decomposes is not continuous.
However, different choices yield elements that can be canonically deformed into each other, so we can still introduce  ``$H$-space like group structures'' on the homotopy groups.

Recall that $\Shift_{j,R}\colon x \mapsto x + Re_j$ is a diffeomorphism and the abbreviation $\Shift_{j,R}(g) \deff \Shift_{j,R\ast}(g) \deff \Shift_{j,-R}^\ast(g)$.  

\begin{lemma}\label{MonoidStructuresDifectFib - Lemma}
 For $p \geq 1$ and $q\geq 0$, there are $p$-many monoid structures on $\pi_q(G_{p+1})$ that are defined as follows:
 
 Let $g,h \colon D^n \rightarrow G_{p+1}$  represent $[g]$, $[h] \in \pi_q(G_{p+1})$ and assume that they factor through $\ConcSet[ ][M \times RI^{p+2}]$.
 Then, for all $2\leq j\leq p+1$, we define
 \begin{equation*}
     (g +_{j,R} h)(t) \deff \Shift_{j,-2R}(g(t))|_{\{x_j \leq 0\}} \cup \Shift_{j,2R}(h(t))|_{\{x_j \geq 0\}}
 \end{equation*}
 and 
 \begin{equation*}
     [g] +_j [h] \deff [g +_{j,R} h] \in \pi_q(G_{p+1}).
 \end{equation*}
 
 If $p\geq 2$, then all of these monoid structures are Eckmann-Hilton related to each other.
 If $q \geq 1$, then these monoids structures are Eckmann-Hilton related to the group structure of the homotopy group.
 
 The same statement is true for $\pi_q(F_p)$.
\end{lemma}
\begin{proof}
  Since $\Shift_{j,-2R}(g(t))$ and $\Shift_{j,2R}(h(t))$ agree with the base point near $\{x_j = 0\}$, the metric $(g+_{j,R} h)$ is well-defined and has positive scalar curvature. 
  The shift operators also preserves the defining conditions of $G_{p+1}$, so $g +_{j,R} h \in G_{p+1}$.
  
  For $R_0 \leq R_1$ and all $s \in [0,1]$, we set $R_s \deff (1-s)R_0 + sR_1$.
  A homotopy relative to $\partial D^n$ between $g +_{j,R_0} h$ and $g +_{j,R_1}h$ is given by $s \mapsto g +_{j,R_s}h$. 
  
  Let $g,h \colon (D^n \times I^n,\partial D^n \times I) \rightarrow (G_{p+1},g_0)$ homotopies relative to the boundary.
  For a sufficiently large $R$, the map $g +_{j,R} h$ is a homotopy relative to the boundary between $g(\placeholder,0) +_{j,R} h(\placeholder,0)$ and $g(\placeholder, 1) +_{j,R} h(\placeholder,1)$.
  
  It is straightforward to check that these structures are associative and that the neutral element is given by the base point $g_0$.
  
  The proof that these structures are Eckmann-Hilton related to each other if $p \geq 2$ is formally identical with the proof that higher homotopy groups are abelian.
  
  The proof that these monoid structures are Eckmann-Hilton related to the group structure from the homotopy groups is formally the same as the proof for $H$-spaces.
\end{proof}

\begin{cor}
 For $q \geq 1$, the groups $\pi_q(G_{p+1})$ and $\pi_q(F_p)$ are abelian if $p + q \geq 2$.
\end{cor}
 
 The reader certainly has notice the similarities in notation and formulation to Theorem \ref{Geometric Addition - Thm}.
 This is no coincidence as it carries over to the homotopy fibre.
 \begin{proposition}\label{Geometric Addition hofib - Prop}
  There are $n$-many group structures on $\pi_n(\hofib,g_0)$, for all $n\geq 1$, defined as follows:
  If $g, h$ represent $[g], [h] \in \pi_n(\hofib)$ and decompose outside $M \times R I^{n+1}$, then 
  \begin{equation*}
      [g] +_j [h] \deff [g +_{j,R} h]
  \end{equation*}
  for all $2 \leq j \leq n+1$.
  
  If $r_j$ denotes the reflection at the hyperplane $\{x_j = 0\}$, then the inverse element of $[g]$ for $+_j$ is given by $[r_j^\ast g]$.
  The structures are Eckmann-Hilton related to each other.
  Furthermore, $+_2$ agrees with the group structure provided by cubical set theory.
 \end{proposition}
 \begin{proof}
  As the construction $g \mapsto g^\angle$ preserves suspensions, see Proposition \ref{Construction-AngleRotation - Prop}, all arguments in the proof of Theorem \ref{Geometric Addition - Thm} carry over immediately to the homotopy fibre.
 \end{proof}
\begin{cor}\label{suspension monoidho - Cor}
 The map $\iota_{p,q} \colon \pi_{q}(F_p) \rightarrow \pi_{q-1}(F_{p+1})$ is a monoid homomorphism for all $p \in \N_0$ and all $q \geq 1$.
 Furthermore, $\iota_{p,0} \colon \pi_0(F_{p}) \rightarrow \pi_0(\hofib[p])$ is a monoid homomorphism for $p \geq 1$.
\end{cor}
\begin{proof}
 For $q > 1$, the claim follows from Lemma \ref{Boundary Comparison - Lemma}, so it remains to check the case $q=1$.
 
 If we realise elements of $\pi_1(F_p)$ as homotopy classes of block maps $\gamma \colon \R \rightarrow F_p$, then the group structure is given $[\gamma_1] \cdot [\gamma_2] = [\gamma_1(\cdot - 2\rho) \cup \gamma_2(\cdot + 2\rho)]$ for all sufficiently large $\rho > 0$ such that $\gamma_j$ is constant (with value $g_0$) outside of $[-\rho,\rho]$.
 
 The claim now follows from
 \begin{align*}
     \iota_{p,1}([\gamma_1] \cdot [\gamma_2]) &= [\susp(\gamma_1(\cdot + 2\rho) \cup \gamma_2(\cdot - 2\rho))] \\
     &= [{\Shift_{-2\rho}}(\susp(\gamma_1))|_{\{x_{p+2} \leq 0\}} \cup {\Shift_{2\rho}} (\susp(\gamma_2))|_{\{x_{p+2} \geq 0\}}] \\
     &= \iota_{p,1}([\gamma_1]) +_{p+2} \iota_{p,1}([\gamma_2]).
 \end{align*}
 
 To see that $\iota_{p+1,0}$ is a monoid homomorphism, one uses the geometric addition on $\pi_{p+1}(\hofib)$, which agrees with the usual one by Proposition \ref{Geometric Addition hofib - Prop}.  
\end{proof}

\begin{definition}\label{psc-pseudoisotopies - def}
  Define the space of \emph{psc pseudo isotopies} to be the following subspace
  \begin{equation*}
      \pscPseudoIso{p} \deff \{g \in G_p \, : \, \face[1]{1}g = g_0\}
  \end{equation*}
  of $G_p$.
\end{definition}
\begin{lemma}\label{PSCPseudoInclusion - Lemma}
 The inclusion $\pscPseudoIso{p} \hookrightarrow G_p^{slow}$ is a weak homotopy equivalence.
\end{lemma}
\begin{proof}
 We will actually show the following slightly stronger statement:
 Every block map $g \colon \R^n \rightarrow G_p^{slow}$ can be homotopied through block maps into a block map $\R^n \rightarrow \pscPseudoIso{p}$ such that if $g(t) \in \pscPseudoIso{p}$, then the homotopy at $t \in \R^n$ stays in $\pscPseudoIso{p}$.
 
 We find a sufficiently large $R>0$ such that $g(t)$ decomposes away from $M \times RI^{p+1}$ for all $t \in \R^n$.
 Since each block map is uniquely determined by its restriction to a sufficiently large cube, we can find a sufficiently large $\rho > R$ such that $(\face[1]{1}g(t))^\angle_\rho$ is a psc metric for all $t \in \R^n$, where $(\face[1]{1}g(t))^\angle_\rho$ refers to the construction of Proposition \ref{Construction-AngleRotation - Prop} in the $x_1$ - $x_{p+1}$ plane in $\R^{p+1}$ around the point $(\rho, 0, \dots, 0, \rho)$, see Figure \ref{fig:HomotopyPscConcordance} for a visualisation.
 \begin{figure}[htpb]
     \centering
     \includegraphics[width = \textwidth]{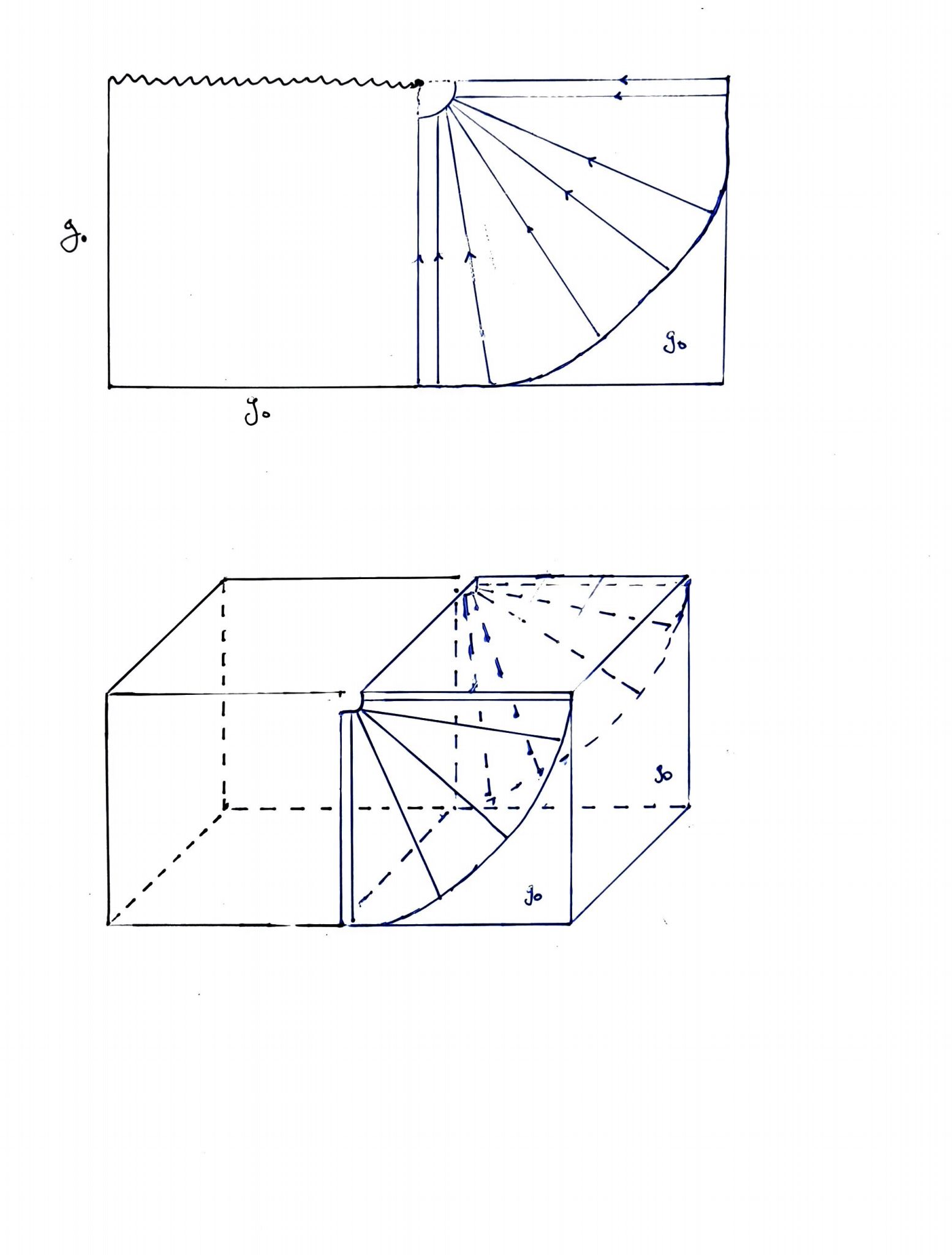}
     \caption{The construction for $p = 1$ and $p=2$.}
     \label{fig:HomotopyPscConcordance}
 \end{figure}
 So far, $(\face[1]{1}g(t))^\angle_\rho$ is only a metric on $\{x_1 \geq \rho, x_{p+1} \leq \rho\}$. 
 We use its product form near $\{x_{p+1} = \rho\}$ to extend it to $\{x_{1} \geq \rho\}$. 
 Since $g(t)$ and $(\face[1]{1}g(t))^\angle_\rho$ agree near $\{x_1 = \rho\}$, the map
 \begin{equation*}
     t \mapsto g(t)|_{\{x_1 \leq \rho\}} \cup (\face[1]{1}g(t))^\angle_\rho |_{\{x_1 \geq \rho\}}
 \end{equation*}
 is a smooth block map $\R^n \rightarrow G_p$, which we denote by $g \cup (\face[1]{1} g)^\angle_\rho$.
 If $g(t)$ lies in  $\pscPseudoIso{p}$, then $g(t) \cup (\face[1]{1}g(t))_\rho^\angle$ also lies in there.
 
 We claim that this map is homotopic to $g$ via a homotopy with the desired properties.
 Since $\face[1]{1}g(t) = \susp(h(t,\placeholder))$, where $h(t,\placeholder) \colon \R^{p} \rightarrow \Riem^+(M)$ is a block map that lies in the subset $A_p$, the extension of $(\face[1]{1} g)^\angle_\rho$ to the half plane satisfies
 \begin{equation*}
     (\face[1]{1} g(t))^\angle_\rho = \susp(\mathfrak{h}_\rho(t,\placeholder)),
 \end{equation*}
 where $\mathfrak{h}_\rho(t,\placeholder) \colon \{x_1 \geq \rho\} \subseteq \R^{p+1} \rightarrow \Riem^+(M)$ is independent of $x_1$ near $\{x_1 = \rho\}$ and $\{x_1 \geq 2 \rho\}$, see Proposition \ref{Construction-AngleRotation - Prop}.
 We may assume, after rescaling $h(t,\placeholder)$ in the second component (with a constant independent of $t$)\footnote{Note that this does not change the homotopy type of the block map $h$}, that $h(t,\placeholder)$ is sufficiently slowly parameterised such that $\mathfrak{h}_\rho(t,\placeholder)$ satisfies the same speed-constrains as an element of $A_{p+1}$.
 
 If $\varphi \colon [\rho , 2\rho] \rightarrow [\rho,2\rho]$ is a smooth, monotonically increasing, surjective function that is constant near boundary, and such that $\supp(\varphi' - 1)$ is contained in a sufficiently small neighbourhood of $\{x_1 = \rho\}$ and $\{x_1 = 2\rho\}$ on which $\mathfrak{h}_\rho(t;\placeholder)$ is still independent of $x_1$ for all $t \in \R^n$, then $\susp(\mathfrak{h}_\rho(t;\varphi(\placeholder),\placeholder))$ is still a psc metric and convex combinations yield the following chain of homotopies
 \begin{align*}
     g(t) \cup (\face[1]{1} g(t))_\rho^\angle &= g(t) \cup \susp(\mathfrak{h}_\rho(t,\placeholder)) \\
     &\simeq g(t) \cup  \susp(\mathfrak{h}_\rho(t,\varphi(\placeholder),\placeholder)), &\text{via } s \mapsto (1-s)\id + s \varphi,  \\
     &\simeq g(t) \cup  \susp(\mathfrak{h}_\rho(t,\rho,\placeholder)), &\text{via } s \mapsto (1-s)\varphi + s \rho. \,  
 \end{align*}
 Note that the last homotopy is a homotopy of positive scalar curvature metrics because $\mathfrak{h}_\rho(t,(1-s)\rho+s\rho,\placeholder)$ satisfies the speed-constrains of $A_{p+1}$ since $\mathfrak{h}_\rho(t,\varphi(\placeholder))$ does it.
 
 If $g \in \pscPseudoIso{p}$ so that $\face[1]{1}g(t) = g_0$, then $(\face[1]{1}g(t))_\rho^\angle = g_0$ and the homotopies stay in $\pscPseudoIso{p}$. 
\end{proof}

\begin{cor}\label{MonoidPSCPseudoIso - Cor}
  The set $\pi_q(G_{p+1})$ has $(p+1)$-many Eckmann-Hilton related monoid structures that extends the ones from Lemma \ref{MonoidStructuresDifectFib - Lemma}. 
\end{cor}
\begin{proof}
 The defined additions in Lemma \ref{MonoidStructuresDifectFib - Lemma} restrict to $\pi_q(\pscPseudoIso{p})$.
 The additional monoid structure is given by $+_1$.
\end{proof}

\section{Abstract Construction of the Spectral Sequence}\label{Section - Abstract Construction of the Spectral Sequence}

We describe the algebraic machinery to construct the psc Hatcher spectral sequence.
This forces us to leave the realm of homological algebra as not all sets we are considering are abelian groups.

The starting point is the following exact couple of $\Z$-bigraded, pointed sets
\begin{equation}\label{eq - AbstractExactCouple}
    \xymatrix{D_{p-1,q+1} \ar[rr]_{i_{p-1,q+1}}^{(1,-1)} && D_{p,q} \ar@/^1pc/[ld]^{(0,0)}_{j_{p,q}} \\
              & E_{p,q} \ar@{-->}@/^1pc/[lu]^{(-1,0)}_{k_{p,q}}.&}
\end{equation}
Here, the tuples in the diagram refer to the bigrading of the corresponding map.
We further require that the sets $D_{p,q}$ satisfy the following properties:
\begin{itemize}
    \item $D_{p,q} = 0$ if $p<0$.
    \item $D_{p,0}$ is a monoid if $p \geq 1$, it is abelian if $p\geq 2$.
    \item $D_{p,q}$ is a group if $q \geq 1$, it is abelian if $p+q\geq 2$.
\end{itemize}
In particular, $D_{p,q}$ is not a monoid only if $p=q=0$, and it is not an abelian monoid only if $p+q \leq 1$.
We require the sets $E_{p,q}$ to satisfy
\begin{itemize}
    \item $E_{p,q}=0$ if $p<0$ or $q<0$.
    \item $E_{p,0}$ is a monoid if $p\geq 2$, it is abelian if $p\geq 3$.
    \item $E_{p,q}$ is a group if $q \geq 1$, it is abelian if additionally $p+q \geq 2$.
\end{itemize}
In particular, $E_{p,q}$ is not a monoid only if $q=0$ and $p \leq 1$, it is not an abelian monoid only if $p+q \leq 1$ or $(p,q) = (2,0)$.
We require further that
\begin{itemize}
    \item all maps in diagram (\ref{eq - AbstractExactCouple}) are homomorphisms of monoids if domain and target are monoids.
\end{itemize}
This implies, in particular, that the kernel of a map, $\mathrm{i.e.}$, the pre-image of the base point is a submonoid or subgroup whenever the domain of the map is a monoid or a group.
Exact sequences of (abelian) monoids are in general badly behaved; monoid homomorphisms are not necessarily injective if the kernel is trivial, and the homomorphism-theorem for monoids fails to hold\footnote{A counterexample is $(\N_0,+) \rightarrow (\Z_2,\times)$,  $n \mapsto [2^n]$. This monoid homomorphisms is not injective but has zero kernel and the homomorphism theorem fails to hold because $\N_0/\{0\} = \N_0 \neq (\Z_2,\times)$}.

If the target is an abelian group, then the homomorphism theorem still holds true.
\begin{lemma}\label{MonoidHomTheorem - Lemma}
 Given an exact sequence of abelian monoids 
 \begin{equation*}
     \xymatrix{0 \ar[r] & L \ar@{^{(}->}[r]^f & M \ar@{->>}[r]^g & A \ar[r] & 0.}
 \end{equation*}
 If $A$ is an abelian group, then $M/f(L) = \mathrm{coker}(f)$ is also an abelian group and $g$ induces an isomorphism $\bar{g} \colon M/f(L) \xrightarrow{\cong} A$.
\end{lemma}
\begin{proof}
 Recall that $M/f(L) = M/\sim$ is the set of equivalence classes with respect to the equivalence relation 
 \begin{equation*}
     m_1 \sim m_2 :\Longleftrightarrow  \exists l_1, l_2 \in L \text{ s.t. } m_1 + f(l_1) = m_2 + f(l_2).
 \end{equation*}
 It inherits the abelian monoid structure from $M$.
 By exactness, the canonical map $\bar{g}$ is well-defined and surjective.
 
 To prove that $M/f(L)$ is an abelian group we argue as follows:
 For all $x \in M$, there is, since $g$ is surjective, a $y\in M$ such that $g(y) = - g(x)$.
 Hence, by exactness, $[y] \in M/f(L)$ is an additive inverse of $[x] \in M/f(L)$, so $M/f(L)$ is an abelian group.
 
 We already know that $\bar{g}$ is surjective.
 Since $M/f(L)$ is an abelian group, it remains to show that $\ker(\bar{g})$ is trivial.
 But this follows immediately from the fact that $g^{-1}(0) = f(L)$ and $M \rightarrow M/f(L)$ that is surjective.
\end{proof}
 Since $E_{p,q} = 0$ if $q<0$, the previous lemma implies the following consequence.
\begin{cor}\label{Dpq are stable - Cor}
 The maps $\iota_{p,q} \colon D_{p,q} \rightarrow D_{p+1,q-1}$ are isomorphisms of monoids if $p\geq 1$ and $q<0$.
\end{cor}

Although we cannot apply the machinery of homological algebra to turn this exact couple into a spectral sequence, we can still construct the spectral sequence in an ad-hoc manner by mimicking the construction of a spectral sequence out of an exact couple.
For $r \geq 1$, we define
\begin{align*}
    Z_{p,q}^r &\deff k_{p,q}^{-1}\left(\mathrm{im}[i^{r-1}\colon D_{p-r,q+r-1} \rightarrow D_{p-1,q}]\right),\\
    B_{p,q}^r &\deff j_{p,q}\left(\ker [i^{r-1}\colon D_{p,q} \rightarrow D_{p+r-1,q-r+1}]\right).
\end{align*}
These sets fit into a chain of inclusions
\begin{equation*}
    0 = B_{p,q}^1 \subseteq B^2_{p,q} \subseteq \cdots \subseteq Z^2_{p,q} \subseteq Z^1_{p,q} = E_{p,q}
\end{equation*}
and we use this chain to define
\begin{equation*}
    Z^\infty_{p,q} \deff \bigcap_{r\geq 1} Z^r_{p,q} = Z^{p+1}_{p,q} \quad \text{and} \quad B^\infty_{p,q} \deff \bigcup_{r \geq 1} B^r_{p,q} = B^{p+2}_{p,q},
\end{equation*}
where we used Corollary \ref{Dpq are stable - Cor} to deduce $B^\infty_{p,q} = B^{q+2}_{p,q}$.

For  $q\geq 1$, we define
\begin{equation*}
    E^r_{p,q} \deff Z^r_{p,q}/B^r_{p,q}.
\end{equation*}
Note that the groups $Z^r_{0,1}$ and $B^r_{0,1}$ do not need to be abelian, so $E^r_{0,1}$ might be just a pointed set of left cosets.

We now define the differentials.
For $p,q\geq 0$, we set 
\begin{equation*}
    d^r_{p,q} \deff j_{p-r,q+r-1} \circ (i^{r-1})^{-1} \circ k_{p,q} \colon Z^r_{p,q} \rightarrow E^r_{p-r,q+r-1}.
\end{equation*}
For all other sets, we define the differentials to be zero, $\mathrm{i.e.}$, the constant map whose value is the base point.

\begin{lemma}\label{HSS-Differentials - Lemma}$ $
  \begin{itemize}
      \item[(i)] The differentials are well-defined.
      \item[(ii)] $\mathrm{im} \, d^{r}_{p,q} = B^{r+1}_{p-r,q+r-1}/B^r_{p-r,q+r-1}$ for $r \geq 2$ and $\mathrm{im}\, d^1_{p,q} = B^2_{p-1,q}$.
      \item[(iii)] $\ker d^r_{p,q} = Z^{r+1}_{p,q}$ for $r \geq 1$.
  \end{itemize}
\end{lemma}
The proof of this lemma is inspired by the discussion of \cite{boardman1999conditionally}. 
However, we cannot rely on isomorphism-theorems to identify various quotients, so we have to calculate everything ``by hand''.
\begin{proof}
  We start with the proof of $(i)$.
  There is nothing to prove for $d^1_{p,q}$, so we assume $r \geq 2$.
  
  By definition, $k_{p,q} \colon Z^r_{p,q} \rightarrow \mathrm{im}[ i^{r-1} \colon D_{p-r,q+r-1} \rightarrow D_{p-1,q}]$.
  As $r \geq 2$, the domain of $i^{r-1}$ is a group.
  We may also exclude the case $(p,q) = (1,0)$ as the differential $d^r_{1,0}$ has zero as target, so that the target of $i^{r-1}$ is a monoid, too.
  
  Thus, each pair of elements $x,y \in D_{p-r,q+r-1}$ with $i^{r-1}(x) = i^{r-1}(y) = k_{p,q}(a)$ are automatically in $\ker\lbrack i^{r} \colon D_{p-r,q+r-1} \rightarrow D_{p,q-1} \rbrack$.
  By exactness, there is an element $e \in \ker \iota^{r-1}$ such that $x\cdot e = y$. 
  Clearly,
  \begin{equation*}
      j_{p-r,q+r-1}(e) \in j_{p-r,q+r-1}(\ker i^{r-1}) = B^r_{p-r,q+r-1},
  \end{equation*}
  so the map $d^r_{p,q}\colon a \mapsto j_{p-r}(x) \in B^{r+1}_{p-r,q+r-1}/B^r_{p-r,q+r-1} \subseteq E^r_{p,q}$ is well-defined. 
  
  For $(ii)$, we consider the case $r=1$ first.
  By exactness, $k_{p,q} \colon Z^1_{p,q} \rightarrow \mathrm{im}\, k_{p,q}$, so that 
  \begin{align*}
      \mathrm{im} \, d_{p,q}^1 &= j_{p-1,q}(\mathrm{im} \, k_{p,q}) = j_{p-1,q}(\ker i_{p-1,q}) = B^2_{p-1,q}.
  \end{align*}
  
  For $r\geq 2$, the argument is similar. By exactness, $k_{p,q}\colon Z^r_{p,q} \rightarrow \mathrm{im} \, i^{r-1} \cap \ker i_{p-1,q}$.
  Since $(i^{r-1})^{-1}(\im i^{r-1} \cap \ker i_{p-1,q}) = \ker i^r$, we find for all $x \in \ker i^r$ an element $a \in Z^r_{p,q}$ with $i^{r-1}(x) = k_{p,q}(a)$.
  Thus, $d^r_{p,q}$ is surjective onto $B^{r+1}_{p-r,q+r-1}/B^r_{p-r,q+r-1}$.
  
  We now prove $(iii)$: If $x \in Z^r_{p,q}$ lies in the kernel of $d^{r}_{p,q}$, then there is a $\xi \in D_{p-r,q+r-1}$ such that $k_{p,q}(x) = i^{r-1}(\xi)$ and
  \begin{align*}
      j_{p-r,q+r-1}(\xi) &\in j_{p-r,q+r-1}\left( \ker i^{r-1} \colon D_{p-r,q+r-1} \rightarrow D_{p-1,q} \right) = B^r_{p-r,q+r-1}.
  \end{align*}
  Hence, there is an $\tilde{\eta} \in \ker j_{p-r,q+r-1} \subseteq D_{p-r-1,q+r}$ (since this is a group if $q\geq 0$) such that $\xi \cdot i_{p-r,q+r}(\tilde{\eta}) \in \ker i^{r-1}$.
  By exactness, there is an $\eta \in D_{p-r-1,q+r}$ such that $\tilde{\eta} = \iota_{p-r,q+r}(\eta)$ so that $\xi \cdot \iota_{p-r,q+r}(\eta) \in \ker \iota^{r-1}$.
  If $(p,q) \neq (1,0)$, then domain and image of $i^{r-1}$ are monoids and $i^{r-1}$ is a map of monoids, so 
  \begin{align*}
      0 &= i^{r-1}(\xi \cdot i_{p-r,q+r}(\eta)) = i^{r-1}(\xi)i^{r-1}(i_{p-r,q+r}(\eta)) \\
        &= k_{p,q}(x) \cdot i^{r}(\eta),
  \end{align*}
  or equivalently, $k_{p,q}(x) = i^{r}(\eta^{-1})$.
  Thus, $x \in k_{p,q}^{-1}(\mathrm{im} \, i^{r}) = Z^r_{p,q}$.
  
  If $(p,q) = (1,0)$, then $d^r_{1,0} = 0$ for $r \geq 2$ as the target is zero.
  The claim for $d^1_{1,0}$ follows immediately from exactness of pointed sets.
\end{proof}

\begin{cor}
  $E_{p,q}^{r+1} = H^r(E^r,d^r)_{p,q}$ for all $q \geq 1$.
\end{cor}

\begin{rem}
   The reader may wonder how we define $E^r_{p,0}$ for $r\geq 2$ as $E^1_{p,0}$ is only a monoid.
   The answer is that we don't! 
   Since $B^2_{p,0} = B^\infty_{p,0}$ there is no gain in forming the quotient $E^2_{p,0} = Z^2_{p,0}/B^2_{p,0}$, even if these monoids were abelian groups. 
   
   If all pointed sets are abelian groups, and all maps are homomorphisms, then the here constructed spectral sequence agrees with the usual one derived from exact couple as described in \cite{boardman1999conditionally} or \cite{weibel1996homological}. 
\end{rem}

\subsubsection*{Convergence of the spectral sequence}

Let us now discuss covergence of this spectral sequence.
We will show, in analogy to the spectral sequence derived from an exact couple of abelian groups, that the constructed spectral sequence converges to 
\begin{equation*}
    \colim{} [D_{p-1,q+1} \rightarrow D_{p,q}] 
\end{equation*}
under the assumption that a certain comparison map exists.
We use Corollary \ref{Dpq are stable - Cor} to identify this colimit with $D_{p+q+1,-1}$ and define
\begin{equation*}
    F_{p,q} \deff \mathrm{im}\, \mathrm\lbrack D_{p,q} \xrightarrow{i^{q+1}} D_{p+q+1,-1}\rbrack.
\end{equation*}
Note that $F_{p,q}$ is an abelian group if $p+q\geq 2$ because $D_{p+q,0}$ is abelian in this case and that $\iota_{p+q,0}\colon D_{p+q,0}\twoheadrightarrow D_{p+q+1,-1}$ is surjective.

\begin{lemma}
 Let $p,q \geq 0$ and $p+q\geq 2$.
 If the map
 \begin{equation*}
     \xymatrix{\mathfrak{I}_{p,q}\colon Z_{p,q}^\infty \ar[rr]^-{j_{p,q}^{-1}} && D_{p,q} \ar[rr]^-{i^{q+1}} && F_{p,q}/F_{p-1,q+1}}
 \end{equation*}
 is well-defined, then it is surjective and has $B^\infty_{p,q}$ as kernel.
 In particular, it induces an isomorphism
 \begin{equation*}
     \xymatrix{Z^\infty_{p,q}/B^\infty_{p,q} \ar[rr]^-\cong &&} F_{p,q}/F_{p-1,q+1}. 
 \end{equation*}
\end{lemma}

\begin{rem}
   A sufficient condition for $\mathfrak{I}_{p,q}$ to be well-defined is the following one: 
   For all $x_1,x_2 \in D_{p,q}$ with $j_{p,q}(x_1)=j_{p,q}(x_2)$, there are $y_1,y_2 \in \ker j_{p,q}$ such that $x_1 +y_1 = x_2 + y_2$. 
   This condition is satisfied, for example, if $D_{p,q}$ is a (not necessarily abelian) group.
\end{rem}
\begin{proof}
 By exactness, $Z^\infty_{p,q} = \ker k_{p,q} = j_{p,q}(D_{p,q})$, so if $\mathfrak{I}_{p,q}$ exists, then it is clearly surjective.
 
 To see that $\ker \mathfrak{I}_{p,q} = B^\infty_{p,q}$, we argue as follows:
 Let $x \in D_{p,q}$ such that $\mathfrak{I}_{p,q}(x) = 0$, or equivalently, that $i^{q+1}(x) \in F_{p-1,q+1} \subseteq F_{p,q}$.
 By definition, there is a $y \in D_{p-1,q+1}$ such that $i^{q+2}(y) = \iota^{q+1}(x)$.
   Then $x + \iota_{p-1,q+1}(-y) \in \ker i^{q+1} \subseteq D_{p,q}$.
   Exactness implies
   \begin{equation*}
       j_{p,q}(\ker i^{q+1}) \ni j_{p,q}(x - i_{p-1,q+1}(y)) = j_{p,q}(x).
   \end{equation*}
   Thus, the pre-image of the base point under $\mathfrak{I}_{p,q} \colon Z^\infty_{p,q} \rightarrow F_{p,q}/F_{p-1,q+1}$ agrees with $j_{p,q}(\ker i^{q+1}) = B^\infty_{p,q}$.
   
   We apply Lemma \ref{MonoidHomTheorem - Lemma} to the short exact sequence
   \begin{equation*}
       \xymatrix{0 \ar[r] & B^\infty_{p,q} \ar@{^{(}->}[r] & Z^\infty_{p,q} \ar[r]^-{\mathfrak{I}_{p,q}} \ar@{->>}[r] & F_{p,q}/F_{p-1,q+1} \ar[r] & 0 }
   \end{equation*}
   to deduce that $\mathfrak{I}_{p,q} \colon Z^\infty_{p,q}/B^\infty_{p,q} \rightarrow F_{p,q}/F_{p-1,q+1}$ induces an isomorphism. 
\end{proof}

We summarise the achivements of this subsection in the following theorem.
\begin{theorem}\label{ExactCoupleYieldSS - Theorem}
 For each exact couple of pointed sets satisfying the conditions below equation (\ref{eq - AbstractExactCouple}), there is a first quadrant spectral sequence $(E^r_{p,q},d^r_{p,q})$ that starts with $E^1_{p,q} = E_{p,q}$, and whose differentials on the first page are given by $d^1_{p,q} = j_{p-1,q} \circ k_{p,q}$.
 If the comparision maps $\mathfrak{I}_{p,q} \colon Z_{p,q}^\infty \rightarrow F_{p,q}/F_{p-1,q+1}$ are well-defined, then the spectral sequence converges to $D_{p+q+1,-1}$ for $p+q\geq 2$.
\end{theorem}

\section{Application to the Defect Fibration}\label{Section - Application to the Defect Fibration}

We are now ready to construct the psc Hatcher spectral sequence.
Recall from Proposition \ref{Restriction is Smooth Fib} the fibrations in the smooth category
\begin{equation*}
     \xymatrix@C+1em{F_{p+1} \ar@{^{(}->}[r]^{\mathrm{incl}} & G_{p+1} \ar[r]^{\restrict} & F_p }
 \end{equation*}
for $p \geq 1$. 
Since $F_0 = G_0$, we can extend the fibrations to $p \leq -1$ by setting $F_p = G_p = \{\ast\}$ in this case.
By Lemma \ref{Boundary Comparison - Lemma}, the fibrations induce the following long exact sequences on homotopy groups
\begin{equation*}
    \xymatrix{ \bigoplus_{\substack{p \in \Z \\ q \geq -1}} \pi_{q+1}(F_p) \ar[rr]^{\iota_{p,q}} && \bigoplus_{\substack{p \in \Z \\ q \geq 0}} \pi_{q} (F_{p+1}) \ar[ld]^{\pi_q(\mathrm{incl})} \\
    & \bigoplus_{\substack{p \in \Z \\ q \geq 0}} \pi_q(G_{p+1}) \ar@{-->}[lu]^{\pi_q(\restrict)}. &}
\end{equation*}
We will abbreviate $\pi_q(\restrict)$ to $\restrict_{p+1,q}$ and $\pi_q(\mathrm{incl})$ to $\incl_{p+1,q}$.
We are not yet in the setup of the previous section because the sequences end at $\pi_0(F_p)$.
The next lemma allows us to extend the couple to the negative half-plane $q<0$.
\begin{lemma}
  The sequence of pointed sets
  \begin{equation*}
      \xymatrix{\pi_0(G_{p+1}) \ar[r]^-{\restrict_{p+1,0}} & \pi_0(F_p) \ar[r]^-{\iota_{p,0}} & \pi_p(\hofib) \ar@{->>}[r] & 0}
  \end{equation*}
  is exact.
\end{lemma}
\begin{proof}
 For $p < 0$, there is nothing to prove.
 For all $p \geq 0$, we have the short exact sequence
 \begin{equation*}
     \xymatrix{0 \ar[r] & \ker \iota_{p,0} \ar@{^{(}->}[r] & \pi_0(F_{p}) \ar@{->>}[r]^-{\iota_{p,0}} & \pi_{p}(\hofib) \ar[r] & 0. }
 \end{equation*}
 An element $\sigma \in \Omega^{p}_{g_0}\hofib$ represents the zero element if and only if there is an $h \in \hofib[p+1] \subseteq \ConcSet[p+2]$ that satisfies
 \begin{align*}
     &{}^{\mathrm{hofib}} \face[1]{1}h = \face[1]{2}h = \sigma, & {}^{\mathrm{hofib}} \face[-1]{1}h = \face[-1]{2}h = g_0, \\
     &{}^{\mathrm{hofib}} \face[\omega]{j}h = \face[\omega]{j+1}h = g_0 \text{ for }j\geq 2. & { }
 \end{align*}
 If $\cycl(2:p+2)$ denotes the cyclic permutation $(2 \ 3 \ \dots \ p+2)$, then one calculates
 \begin{align*}
     \cycl(2:p+2) \circ \CubeIncl{2}(x_1,\dots,x_{p+1}) &= \cycl(2:p+2)(x_1,\eps,x_2,\dots,x_{p+1})\\
     &= (x_1,x_2,\dots,x_{p+1},\eps) = \CubeIncl{p+2}.
 \end{align*}
 
 This implies that $h$ is a homotopy in $\hofib[p+1]$ between $g_0$ and $\sigma$ if and only if $\cycl(2,p+2)_\ast h \in G_p$ with $\restrict(h) = \face[1]{p+2}(h) = \sigma$.
 
 In particular, $\ker \iota_{p,0} = \mathrm{im}\, \restrict_{p+1,0}$.
\end{proof}

If we define
\begin{align*}
    D_{p,q} &\deff \pi_q(F_p) \text{ if } q\geq 0 \text{ and } D_{p,q} = \pi_{p+q}(\hofib) \text{ if } q<0,\\
    E_{p,q} &\deff \pi_q(G_p) \text{ if } q\geq 0 \text{ and } E_{p,q} = 0 \text{ if } q<0,
\end{align*}
then we are in the setup of the previous section. 
Indeed, the required monoid and group structures follow from Lemma \ref{MonoidStructuresDifectFib - Lemma}.
The maps $j_{p,q}$ and $k_{p,q}$ are group homomorphisms if $q\geq 1$ because they are induced by continuous maps. 
For $q=0$ and appropriate $p$, this follows because the geometric addition is compatible with inclusions and restriction to faces.
The maps $\iota_{p,q}$, for $p+q\geq 1$, are monoid homomorphisms by Corollary \ref{suspension monoidho - Cor}.

By Theorem \ref{ExactCoupleYieldSS - Theorem} this leads to a first quadrant spectral sequence starting with $E_{p,q} = \pi_q(G_p)$ and whose differentials on the first page are given by $d^1_{p,q} = \incl_{p-1,q}\circ \restrict_{p,q}$.
By Lemma \ref{PSCPseudoInclusion - Lemma}, we can replace it by the homotopy groups of psc pseudo isotopies $E_{p,q}^1 = \pi_q(C^{\mathrm{psc}}_p(M))$.
Convergence of the spectral sequence follows from the following lemma that will be proved geometrically.

\begin{lemma}
 The map
 \begin{equation*}
      \mathfrak{I}_{p,q} \colon Z^\infty_{p,q} \xrightarrow{\incl_{p,q}^{-1}} \pi_q(F_{p}) \xrightarrow{\iota^{q+1}} F_{p,q}/F_{p-1,q+1}.
  \end{equation*}
  is well-defined.
\end{lemma}
\begin{proof}
   It suffices to show that each pair of elements $[g_{-1}],[g_1] \in \pi_q(F_{p})$ related by $\incl_{p,q}([g_{-1}]) = \incl_{p,q}([g_{1}])$ satisfy $\iota^{q+1}([g_{-1}]) - \iota^{q+1}([g_{1}]) \in F_{p-1,q+1}$.
   Let $h \colon \R^{q} \times \R \rightarrow G_p$ be a block map that serves as a homotopy between $g_{-1}$ and $g_1$.
   We may assume $h$ to be smooth, that $g_j$ are elements of $\Omega^{q+1,\infty}F_{p}$, and, after a possible reparametrisation, that $\susp(h)$ has positive scalar curvature.
   
   Since $\face[1]{p+q+1}\susp(h) = \susp(\face[1]{q}h)$, it follows that $[\face[1]{p+q+1}\susp(h)] \in F_{p-1,q+1}$.
   Assume that $\susp(h)$ decomposes outside of $M \times \rho I^{p+q+2}$ for a sufficiently large $\rho$.
   Applying Proposition \ref{Construction-AngleRotation - Prop} to $\face[1]{p+q+1}\susp(h)$ in the $x_{p+q+1}$-$x_{p+q+2}$ direction and glue it to $\susp(h)|_{\{x_{p+q+1}\leq \rho\}}$, we obtain a block metric $\susp(h) \cup (\face[1]{p+q+1}\susp(h))^\angle$ that serves at concordance between $\iota^{q+1}(g_1) +_{p+q+1} \face[1]{p+q+1}(\susp(h))$.  
   Using that the geometric addition agrees with the addition in $\pi_{p+q+1}(\hofib)$ from cubical set theory, we conclude $[\iota^{q+1}(g_1)] + [\face[1]{p+q+1}\susp(h)] = [\iota^{q+1}(g_{-1})]$, so that we have $[\iota^{q+1}(g_1)] - [\iota^{q+1}(g_{-1})] \in F_{p-1,q+1}$.
   We conclude that $\mathfrak{I}_{p,q}$ is well-defined.
\end{proof}

Thus, we have proved the main result of this chapter.
\begin{theorem}[psc Hatcher Spectral Sequence]\label{psc Hatcher SS}
  There is a first-quadrant spectral sequence that starts with $E^1_{p,q} = \pi_q(\pscPseudoIso{p})$, whose differentials on the first page are given by $d^1_{p,q} = \pi_q(\restrict \circ \incl)\colon \pi_q(\pscPseudoIso{p}) \rightarrow \pi_q(\pscPseudoIso{p-1})$, and that converges to $\pi_{p+q}(\hofib)$ for $p+q \geq 2$.
 \end{theorem}

 \chapter{The Remains of the Day}\label{Chapter - The Remains of the Day}

Our overall goal was to develop the foundations to carry out a program for psc metrics similar to the study of $\Diff_\bullet(M)$ via $\BlockDiff(M)$.
This goal was achieved.
We constructed the main player $\ConcSet$, the cubical set of all psc block metrics on $M$ (Section \ref{Section - Foundations on the Concordance Set}), showed that it is a Kan set (Theorem \ref{Concordance Set is Kan - Theorem}), related it to the cubical model $\SingMet$ of the space of psc metric $\Riem^+(M)$ via the suspension map $\susp_\bullet \colon \SingMet \dashrightarrow \ConcSet$ (Lemma \ref{SuspOnAux - Lemma}), and established the first tool to study the difference between those spaces beyond index theory, namely the psc version of the Hatcher spectral sequence (Theorem \ref{psc Hatcher SS}).

This is important because we showed that the index difference factors through $\ConcSet$ (Theorem \ref{IndDifFactor - Theorem}), so the index difference is as concordance invariant as possible - it cannot detect any higher isotopy informations.
To achieve the factorisation result, we constructed a new model for real $K$-theory based on the notion of invertible block Dirac operators.
The main feature of this new model is that the Gromov-Lawson index difference can be modeled by the cubical map that assigns a psc block metric its Dirac operator.

We now describe how one can proceed from here.
To get information on $\BlockDiff(M)$, one uses the fibration
\begin{equation*}
    \xymatrix{{\hAut_\bullet(M)}/{\BlockDiff(M)} \ar[rr] && B\BlockDiff(M) \ar[rr] && B\hAut_\bullet(M)   }
\end{equation*}
because the fibre can be approached using surgery theory and the base can be approached using obstruction theory.
There is no immediate counterpart for homotopy automorphism in the world of psc metrics, but guided by the principle that generalisations may lead to easier computable objects, we could consider the bordism version of $\ConcSet$.
This is the cubical set $\BordPscSet$ whose $0$-cubes are pairs of a closed, spin manifold and a positive scalar curvature on it, whose $1$-cubes are bordism of those, whose $2$-cubes are bordisms of bordisms and so on.
In general, an $n$-cube of $\BordPscSet$ is a $(d+n)$-dimensional compact spin manifold with $2n$-many (possibly empty) faces with a positive scalar curvature metric on it that has product structure near each face\footnote{Actually, it should be the elongation of those so that we end up with a ``block manifold''; but for the sake of simplicity we ignore this here.}.
The corresponding cubical set $\BordSet$ without the data of a psc metric should be a cubical model for $\Omega^{\infty + d}MSpin$ because it is the cubical analog of Quinn's bordism space.

Both of these sets are Kan, in fact, the proof that $\BordPscSet$ is Kan is even easier than the proof that $\ConcSet$ is Kan.

There is an obvious map from $\ConcSet[\bullet][M^d]$ into the homotopy fibre of the forgetful map
\begin{equation*}
    \ConcSet[\bullet][M^d] \rightarrow \hofib[\bullet]\bigl(\BordPscSet \rightarrow \BordSet,M^d\bigr).
\end{equation*}
The results of Stolz \cite{stolz1998concordance} and straightforward generalisations thereof should imply that this map is a weak homotopy equivalence if $M$ is simply connected.
In fact, the homotopy fibration
\begin{equation*}
    \xymatrix{\ConcSet[\bullet][M^d] \ar[r] & \BordPscSet \ar[r] & \BordSet }
\end{equation*}
should be a space-level description of the exact sequence in \cite{stolz1998concordance} (in the case of simply connected spin manifolds).
If $M$ is not simply connected, then one has to enrich the model of $\BordPscSet$ and $\BordSet$ with reference maps to $B\pi_1(M)$.
\\

In our index theoretic consideration, we have completely ignored the fundamental group.
In future work, one could construct a version of $\InvPseudDir[\bullet]$ and $\InvBlockDirac$ that is suitable for higher index theory.
This approach would require to extend the theory developed in this thesis to pseudo differential operators on Hilbert $C^\ast$-modules. 
Furthermore, we would need a new proof of Theorem \ref{OperatorSuspWeakEquiv - Theorem} as the $K$-theory groups of a general $C^\ast$-algebra may not be cyclic (this assumption was used in the proof of Theorem \ref{OperatorSuspWeakEquiv - Theorem} to deduce bijectivity from surjectivity) and the theorem of Kucerovsky \cite{kucerovsky1997kk} is not directly applicable (because the positivity condition is not preserved by small perturbations).
However, all these obstacles are likely to be overcome.

The advantage of the models $\InvPseudDir[\bullet]$ and $\InvBlockDirac$ for the classifying space of real $K$-theory is that they emphasise the role of the Dirac operator.
The analysis developed in this thesis allows us to define the cubical set $\mathcal{R}^{\mathrm{inv}}_\bullet(M)$ of all block metrics whose Dirac operator is invertible.
The proof that $\InvBlockDirac$ is Kan also shows that $\mathcal{R}^{\mathrm{inv}}_\bullet(M)$ is Kan because a careful analysis of this proof shows that all applied constructions preserve the Dirac operator; in other words, if we carry out such a construction on a Dirac operator, we produce a Dirac operator.
\\

We have constructed a psc version for the Hatcher spectral sequence, but, so far, this spectral sequence cannot be used for explicit calculations.
The differential $d^1_{p,q}$ in the Hatcher spectral sequence for diffeomophism groups can be computed in the range $q\ll p$ because \emph{Igusa's stabilisation map}
\begin{equation*}
    \Sigma \colon C(M \times I^p) \rightarrow C(M\times I^{p+1})
\end{equation*}
induces an isomorphism on all homotopy groups of degree $\leq q$.
With the methods we have developed here, it is easy to define an analogous stabilisation map for the psc pseudo isotopies.
A psc analog of Igusa's stabilisation result would be of great and also independent interest.

 
 \appendix

 \chapter{Sobolev Spaces}\label{Sobolev Spaces - Chapter}

We will recall the theory of Sobolev spaces for the convenience of the reader. 
This also serves as the opportunity to fix notations and conventions. 
As a disclaimer: Nothing presented in this section is original work from the author but a collection from the following sources: \cite{atiyah1968indexI}, \cite{bleecker2013index}, \cite{LawsonMichelsonSpin} and \cite{palais1965seminar}.

Throughput this section $E \rightarrow M$ denotes a complex vector bundle\footnote{This discussion also works for real and Real bundles, of course.
One simply has to replace the Hermitian metric by a Riemannian metric. 
Then the Sobolev spaces are real and Real vector spaces.} over a not necessarily compact manifold $M^d$ (without boundary), $h$ is a Hermitian metric $E$. 
A Riemannian metric $g$ on $M$ defines a canonical Borel regular measure $\diff \vol_g$ on $M$. 
If $M$ is oriented the measure is given by its volume form, if $M$ is not orientable, it is given by half of the volume form on its oriented cover.
In this section we call the tupel $(E,h,\nabla) \rightarrow (M,g)$ a \emph{geometric bundle}.

\begin{definition}\label{L2spaces - Definition}
 Let $(E,h) \rightarrow (M,g)$ be a Hermitian vector bundle over a Riemannian manifold.
 Its $L^2$-\emph{inner product} on $\Gamma_c(E)$ is given by
 \begin{equation*}
  \langle u_1,u_2 \rangle \deff \int_M h_x(u_1(x),u_2(x)) \diff \mathrm{vol}_g(x).
 \end{equation*}
 The Hilbert space of \emph{square integrable sections} $L^2(M;E)$ (or simply $L^2(E)$) 
 \nomenclature{$L^2(E)$}{Hilbert space of square integrable sections of $E$}
 is the completion of $\Gamma_c(E)$ with respect to the norm induced by the $L^2$-inner product.
\end{definition}\label{SobolevNorms - Definition}
 A connection $\nabla \colon \Gamma(E) \rightarrow \Gamma(T^\vee M \otimes E)$ is a linear map that satisfies the Leibniz rule: $\nabla(f u) = \diff f \otimes u + f\nabla u$.
 One can show without using Sobolev theory that every connection has a formal adjoint, see \cite{bleecker2013index}*{Section 6.6} that is an operator $\nabla^\ast$ that is uniquely satisfied by
 \begin{equation*}
  \langle \nabla^\ast u,v\rangle = \langle u,\nabla v\rangle. 
 \end{equation*}  
 Unbounded operator theory \cite{grubb2008distributions}*{Theorem 12.11} implies that $\nabla^\ast \nabla$ is an (unbounded) self adjoint operator, so we can apply unbounded functional calculus to it. 
\begin{definition}\label{SobolevSpaces - Definition}
 Let $(E,h,\nabla) \rightarrow (M,g)$ be a Hermitian vector bundle with connection over a Riemannian manifold.
 For each $s \in \R$ define the \emph{Sobolev} $s$\emph{-norm} on $\Gamma_c(E)$ via 
 \begin{equation*}
  ||u||_s^2 \deff \int_M h_x((\id + \nabla^\ast\nabla)^s u(x),u(x)) \diff \vol_g(x).
 \end{equation*}
 The \emph{Sobolev} $s$\emph{-space} $H^s(E)$
 \nomenclature{$H^s(E)$}{Sobolev $s$-space of the bundle $E$}
 is defined as the completion of $\Gamma_c(E)$ with respect to the Sobolev $s$-norm.
\end{definition}
In fact, Sobolev spaces are Hilbert spaces and $H^0(E) = L^2(E)$. We will use both notations interchangeably throughout this text.

\textbf{Warning:} The definition of the Sobolev norms and the Sobolev spaces depends dramatically on the choices of Hermitian metric $h$, the connection $\nabla$, and the measure $\vol_g$. Different choices may yield to non-isomorphic \emph{topological} vector spaces. In other words, for different choices, the identity on $\Gamma_c(E)$ may not extend to a continuous linear map.

Let $K$ be a compact set that is the closure of its interior. 
There are two inequivalent ways to define Sobolev sections on $K$.
The first way is the completion of $\Gamma_c(K^\circ,E|_{K^\circ})$ with respect to the restricted Sobolev norms. 
We denote the result with $H^s(K,E)$, $H^s(E|_K)$ or $H_K^s(E)$ 
\nomenclature{$H_K^s(E)$}{Sobolev sections with support in $K$}
depending on the interpretation, if we want to consider the space in its own right or as a subspace\footnote{usually, the space $H^s_K(E)$ denotes the set of all sections $u$ supported within $K$, which means that $\langle u,v\rangle_{L^2} = 0$ for all $v \in \Gamma_c(E)$. 
The differences between these two space are minimal if $K$ is the closure of its interior: With this interpretation, $H_K^s(E)$ equals the intersection of $H^s(E|_{L^\circ})$ of all relative compact open subsets $L^\circ$ of $K$.} of $H^s(M;E)$

The second way only works if $s\geq 0$.
We define $H^s(K \subseteq M,E) \deff \{u|_K \, : \, u \in H^s(E)\}$. 
There are two norms $H^s(K \subseteq M,E)$. 
The first one is just the usual Sobolev norm where the integration domain is just changed:
\begin{equation*}
 ||u||_{s,K}^2 \deff \int_K h_x((\id + \nabla^\ast \nabla)^s u(x),u(x)) \diff \vol_g(x).
\end{equation*} 
The second norm is given by the infimum of the Sobolev norm of all extensions: 
\begin{equation*}
||u||_{s,K\subseteq M} \deff \{||\tilde{u}||_s \, : \, \tilde{u}|_K = u\}.
\end{equation*}
In general $||u||_{s,K} \leq ||u||_{s,K\subseteq M}$ but they agree in general only if $s=0$.

Also the two variants of the Sobolev spaces are not the same. 
In general $H(E|_K) \subseteq H(K \subseteq M,E)$ but, for example, every section $u \in H^1(E|_K)$ vanishes at the boundary, \cite{grubb2008distributions}*{Remark 4.26}.

When we discuss pseudo differential operators on non-compact manifolds, we use other versions of Sobolev spaces which are not normed spaces. 
The next definition is slightly modified from the one in \cite{atiyah1968indexI}.
\begin{definition}
 For every $s \in \R$ and each pair $(L,K)$ of relative compact open subsets with $\overline{K} \subseteq L$, extension by zero yields an isometric embedding $H^s(K,E|_K) \hookrightarrow H^s(L,E|_L)$.
 Define $H^s_{cpt}(E)$ as the colimit over all these inclusions. 
\end{definition}

\begin{definition}
 For each $s\geq 0$, define $H^s_{loc}(E)$ 
 \nomenclature{$H^s_{loc}(E)$}{Space of local Sobolev $s$-sections}
 be the sub-vector space of all locally square integrable function $L^2_{loc}(E)$ such that $x \mapsto h_x((1 + \nabla^\ast\nabla)^s u(x),u(x))$ is locally square integrable.
 For negative $s$ we define $H^s_{loc}(E)$ as the dual space of $H^{-s}_{cpt}(E^\vee)$.
\end{definition}
Note that $\Gamma(E) \subseteq H^s_{loc}(E)$ for all $s \geq 0$. 
It is easy to see if $s$ is an integers; the general case follows from $||(1+\nabla^\ast\nabla)^tu||_0 \leq ||(1+\nabla^\ast\nabla)^su||_0$ if $0 \leq t \leq s$ (where we allow the norms to be infinity).

The vector space $H^s_{loc}(E)$ carries a natural Fréchet topology induced by the semi-norms $||\placeholder||_{s,K \subseteq M}$ if $s\geq 0$. 
Although $H_{cpt}^s(E)$ includes into $H^s_{loc}(E)$,  the colimit topology on $H^s_{cpt}(E)$ is not the subspace topology of $H^s_{loc}(E)$. 
Also note that both topologies are independent of the choices of the geometric data $h$, $\nabla$, and $\vol_g$.

\begin{example}
 If $\C^N \rightarrow \R^n$ is the trivial bundle endowed with the standard metrics and the trivial connection.
 Then $\nabla^\ast\nabla = - \sum_{j=1}^n \partial_{x_j}^2$ is the standard Laplacian.
 Fourier transformation allows us to give a different description of the (standard) Sobolev $s$-norm
 \begin{equation*}
  ||u||_s^2 = \int_{\R^n} (1+|\xi|^2)^s |\hat{u}(\xi)|^2\diff \xi.
 \end{equation*}
 There are positive constants $c_s\leq C_s$, depending only on $n \in \N$ and $s \in \N$ such that 
 \begin{equation*}
  c_s(1+|\xi|^2)^s \leq \sum_{|\alpha|\leq n} |\xi|^{2\alpha} \leq C_s (1 + |\xi|^2)^s,
 \end{equation*}
 so, for the integral $s$, the (standard) Sobolev $s$-norm is equivalent to
 \begin{equation*}
  |||u|||_s^2 \deff \sum_{k=0}^s ||\nabla^k u||_0^2. 
 \end{equation*}  
\end{example}
For $s=1$, the obvious generalisation of the previous example $|||\placeholder|||_1$ to  a geometric bundle agrees with $||\placeholder||_1$.
This provides us the opportunity to do inductions proofs, when we are discussing block pseudo operators.

The standard Sobolev norms on $\R^n$ gives another possibility to construct Sobolev norms on vector bundles by choosing local trivialisations of $E$ and underlying charts of $M$, defining the Sobolev norms locally in the charts and trivialisations, and then gluing them together with a partition of unity.
The Sobolev norms then drastically depend on the choices of the involved charts, trivialisations and the partition of unity. In practise, it requires tedious work to verify that two different choices yield the same topological space $H^s$, not to mention equivalent Sobolev norms.
However, one can show, for example see \cite{bleecker2013index}*{p.197 ff}, that $u \in H^s_{loc}$ if and only if every restriction to some relative compact trivialisations domain $U$ lies in $H^s(U \subseteq \R^n;\C^N)$.
This observation turns out to be quite useful then we discuss pseudo differential operators for their defining properties are of local nature.  

The proof sketch presented in \cite{LawsonMichelsonSpin}*{p.176} carries over to non-compact manifolds in the following manner. 
\begin{theorem}[Sobolev Regularity]
 For all compact neighbourhoods $K \subseteq M$ and each $s > k + n/2$ there is a constant $C(K,s)$ such that all smooth sections $u \in \Gamma(E)$ supported within $K$ satisfy
 \begin{equation*}
  ||u||_{C^k} \leq C(K,s) ||u||_{s}.
 \end{equation*}
 In particular, $H^{s}(E) \subseteq \Gamma^k(M;E)$.
\end{theorem}

\begin{theorem}[Rellich Lemma]
 For each relative compact open subset $K \subseteq M$ and every pair of real numbers $s<t$, the inclusion $H^t(K;E|_K) \hookrightarrow H^s(E)$ is a compact operator.
\end{theorem}

The next result is often stated for compact manifolds only although the proof for $\R^n$ can be easily adapted to the geometric set up.
\begin{proposition}[perfect pairing]\label{Perfect Pairing - Proposition}
 For each $(E,h,\nabla) \rightarrow (M,g)$ and all section $u,v \in \Gamma_c(E)$ the inner product
 \begin{equation*}
  \langle u , v \rangle \deff \int_M h_x(u(x), v(x)) \diff \vol_g(x)
 \end{equation*}
 extends to a unique non-degenerate bilinear form
 \begin{equation*}
  \langle \placeholder , \placeholder \rangle \colon H^s(E) \times \overline{H^{-s}(E)} \rightarrow \C.
 \end{equation*}
 for all $s \in \R$.
\end{proposition}
\begin{proof}
 We need to show that $\euclmetric^\flat \colon \overline{H^{-s}(E)} \rightarrow H^s(E)^\vee$ is an isomorphism; in fact, it will be an isometric isomorphism.
 The Cauchy-Schwarz inequality $\langle u, v \rangle \leq ||u||_s ||v||_{-s}$ implies that $\euclmetric$ extends to a continuous bilinear form.
 Furthermore, it shows $||\langle \placeholder, v\rangle||_{\mathrm{dual}} \leq ||v||_{-s}$.
 
 On the other hand, for all $v \in \Gamma_c(E)$, we have
 \begin{align*}
     ||v||_{-s}^2 &=  \langle (1+\nabla^\ast\nabla)^{-s} v,v\rangle \leq ||(1+\nabla^\ast \nabla)^{-s}v||_s \cdot ||\langle \placeholder, v\rangle||_{\mathrm{dual}} \\
     &= ||v||_{-s} \cdot ||\langle \placeholder, v\rangle||_{\mathrm{dual}},
 \end{align*}
 which implies $||v||_{-s} \leq ||\langle \placeholder, v\rangle||_{\mathrm{dual}}$ and hence equality.
 Thus, $\euclmetric^\flat$ is isometric, in particular, injective with closed image.
 
 If $\euclmetric^\flat$ were not surjective, then we would find by the Hahn-Banach theorem a non-zero functional $\theta \in H^s(E)^{\vee\vee}$ that vanishes identically on the image of $\overline{H^{-s}(E)}$.
 Since $H^s(E)$ is a Hilbert space, it is reflexive, so there is a unique $u_\theta \in H^s(E)$ with $\mathrm{ev}_{u_\theta} = \theta$.
 This implies 
 \begin{equation*}
     0 = \theta(\langle\placeholder,v\rangle) = \mathrm{ev}_{u_\theta}(\langle\placeholder,v\rangle) = \langle u_\theta, v\rangle,
 \end{equation*}
 for all $v \in H^{-s}(E)$.
 Since $\Gamma_c(E)$ lies dense in $H^s(E)$ and $H^{-s}(E)$, this implies $u_\theta = 0$, contradicting our initial assumption $\theta \neq 0$.
\end{proof}

If $H^\infty \deff \mathrm{invlim}_k H^k(E) = \bigcap_k H^k(E)$ is topologised with the inverse limit topology, then $H^\infty(E)$ is a subspace of $\Gamma(M;E)$ with respect to the smooth (weak) topology. 
Similarly, $H^\infty_{cpt} = \Gamma_c(E)$ and $H^\infty_{loc} = \Gamma(E)$, where the smooth sections are equipped with the smooth (weak) topology \cite{atiyah1968indexI}.
If we define $H^{-\infty} = \mathrm{colim}_k H^k = \bigcup H^k$ endowed with the colimit topology, then the perfect pairing yields $H^{-\infty} = (H^\infty)^\vee$.

\textbf{Important remark:} While the Hilbert space structure $L^2(E) = H^0(E)$ is geometrically motivated, the inner products on the other Sobolev spaces are not. 
Their purpose is to provide regularity for their elements.
Thus, the precise Sobolev $s$-norm is of minor importance (except for $s=0$).
Only the underlying topological vector space is of interest.
This leads to the notion of an \emph{hilbertisable} vector space and chains of hilbertisable vector spaces \cite{palais1965seminar}.
In fact, Proposition \ref{Perfect Pairing - Proposition}, identifies $\{H^k(E)\}_{k \in \Z}$ as a discrete chain of hilbertian spaces in the sense of \cite{palais1965seminar}*{Definition VIII.1.1}.

\begin{definition}
  A \emph{discrete chain of hilbertian} spaces is a sequence of hilbertian spaces $\{H^k\}_{k \in \N_0}$ such that
  \begin{itemize}
      \item[(i)] if $k \leq l$, then $H^l$ is a linear subspace of $H^k$ and the inclusion is continuous,
      \item[(ii)] the intersection $H^\infty \deff \bigcap_k H^k$ is dense in $H^0$ (and hence in all $H^k$), and
      \item[(iii)] $H^0$ is a Hilbert space.
  \end{itemize}
\end{definition}
For negative $k$, we set $H^k \deff (H^{-k})^\vee$.
One can show that that there is a inclusion $H^l \hookrightarrow H^k$ for all $k \leq l \in \Z$, see \cite{palais1965seminar}*{p.126}.

We would like to know how the construction of Sobolev spaces behaves under exterior tensor products.
\begin{lemma}[\cite{palais1965seminar}*{Theorem XIV.2}]
For each pair of geometric bundles $(E_j,h_j,\nabla_j) \rightarrow (M_j,g_j)$ and for all $s \geq 0$, the Sobolev spaces of the exterior tensor product $(E_1 \boxtimes E_2, h_1 \boxtimes h_2, \nabla_1 \boxtimes \id + \id \boxtimes \nabla_2 \rightarrow M_1 \times M_2, g_1 \times g_2)$ is given by
\begin{align*}
 &H^s_{cpt}(E_1 \boxtimes E_2) = \bigcap_{0\leq s\leq t} H^t_{cpt}(E_1) \otimes  H^{s-t}_{cpt}(E_2) 
\end{align*}
and
\begin{align*}
    H^s_{loc}(E_1 \boxtimes E_2) = \bigcap_{0\leq t\leq s} H^t_{loc}(E_1) \otimes H^{s-t}_{loc}(E_2)
\end{align*}
\end{lemma}
In \cite{palais1965seminar} the lemma is proven only for Hermitian vector bundles over closed manifolds. 
However, the proof carries over, because we only consider (semi-)norms on compact subsets.
Alternatively, one could prove the lemma for the trivial bundles $\C^N\rightarrow \R^n$ equipped with the standard Hermitian and Riemannian metric first and then use local patching arguments and the fact that two norms obtained in this manner are equivalent on compact subsets.

The next definition is a slight modification of \cite{atiyah1968indexI}, which generalises the analogous definition of \cite{palais1965seminar} to non-compact manifolds.
\begin{definition}
Let $E,F \rightarrow X$ be two Hermitian vector bundles over $(M,g)$ and let $r \in \Z$.
Denote by $Op^r_k(E,F)$ the set of all continuous linear maps $H^s_{cpt}(E) \rightarrow H^{s-r}_{loc}(F)$ and by $Op^r(E,F)$ the set of all linear maps  $H^\infty_{cpt}(E) \rightarrow H^\infty_{loc}(F)$ that extend to continuous linear map $H^k_{cpt}(E) \rightarrow H^{k-r}_{loc}(F)$ for all $k \in \Z$.
\end{definition}
Note that, in the definition of $Op^r$, 
\nomenclature{$Op^r$}{Space of operators of order $r$}
we do not need to demand continuity of the linear maps $H^\infty \rightarrow H^\infty$ because $H^\infty$ carries the inverse limit topology, as continuity follows from the continuity of the extensions. 
The reader may wonder, why we restrict in the definition of $Op^r$ to Sobolev spaces of integral order.
The reason is that it is good enough for our purpose and in many cases easier to check that an operator is in $Op^r$ because we can avoid functional calculus arguments for $(1+\nabla^\ast\nabla)^s$ if $s$ is integral\footnote{if $s$ is negative, we can use the perfect pairing to reduce it to non-negative $s$.} 

We can  turn $Op^r_s$ it into a Fréchet space using the topology of bounded convergence.
This means that a sequence of continuous linear maps $T_\alpha \colon H^s_{cpt} \rightarrow H^{s-r_{loc}}$ converges to $T$ if and only if $||T_\alpha - T||_{s,s-r;K} \rightarrow 0$ for all compact neighbourhoods $K$, where 
\begin{equation*}
    ||S||_{s,s-r;K} \deff \sup\left\{||Sf||_{s-r;K} \, : \, f \in H^{s}(K,E|_K) \text{ with } ||f||_{s} = 1\right\}.
\end{equation*}
The canonical inclusion 
\begin{equation*}
    Op^k \hookrightarrow \prod_{m \in \Z} Op_m^k \quad \text{ given by } \quad A \mapsto (\dots, A,A,A,\dots)
\end{equation*}
has closed image, so that $Op^k$ is also a Fréchet space.
The space $Op^k$ is local in the sense that $A \in Op^k$ if and only if $fAg \in Op^k$ for all smooth and compactly supported real values functions $f$ and $g$.

 \chapter{Foundation of Pseudo Differential Operators}\label{Appendix Pseudos - Chapter}

We will recall the foundations of pseudo differential operators to set up notation and for the convenience of the reader because the general literature often develops the theory for compact manifolds only, while we mostly work with non-compact manifolds.
Nothing presented here is really new, and the author does not claim any originality for these passages.
In this thesis, we need to deal with several topologies on the space on pseudo differential operators. 
Many existing proofs actually show continuity statement if one carefully keeps track of the constant, and this is why we wrote this appendix.    
Those who are familiar with pseudo differential operators, should pay attention to the construction of the \emph{Atiyah-Singer closure} $\PseudOpCl^k(E,F)$ in Definition \ref{ASClosure - Definition} originally defined in \cite{atiyah1968indexI}, the generalised symbol exact sequence Lemma \ref{ASExactSequence - Lemma}, and, the reason why Atiyah and Singer considered the broader class, Theorem \ref{ASTensor - Theorem} and \ref{FamilyASTensor - Theorem}.

We will define pseudo differential operators on vector bundles as operators that look locally like pseudo differential operators on trivial bundles over open subsets of the euclidean space. 
Therefore, we will consider this model case first.

Since the definition of pseudo differential operators requires Fourier transformations, we cannot define them on real vector bundles but we have to define them on Real vector bundles. 
Recall that a Real vector bundle is complex vector bundles $E$ that is equipped with a $\C$-antilinear involution $\overline{\cdot} \colon E \rightarrow E$. 
In this thesis, we will only deal with Real vector bundles that arise from real vector bundles by complexification, meaning, every Real vector bundle is of the form $E = E' \otimes_{\R} \C$ and
the involution is just the complex conjugation on the second factor.
Moreover, every construction we will pursue in this thesis is a construction on the underlying real vector bundle.
For this reason, only in this section we will distinguish real and Real vector bundles.

\begin{definition}
 Let $U \subseteq \R^n$ be open and let $k$ be a real number. 
 A smooth map $p \colon U \times \R^n \rightarrow \End[\C]{\C^{N_1}}{\C^{N_2}}$ is an \emph{amplitude of order k}, if on each compact subset $K \subseteq U$,
 it has the following asymptotic behaviour
 \begin{equation}\label{eq: AmplitudeEstimation}
  \forall \alpha,\beta \in \N_0^n \, \exists C_{\alpha, \beta, K} : ||\partial^\beta_x \partial^\alpha_\xi p(x,\xi)|| \leq C_{\alpha, \beta, K}(1+||\xi||)^{k - |\alpha|}.
 \end{equation}
 An amplitude is \emph{Real}, if it additionally satisfies
 \begin{equation}
  \overline{p(x,\xi)} = p(x, -\xi).
 \end{equation}
 We denote the set of all Real amplitudes by $A_k(U \times \R^n; \C^{N_1}, \C^{N_2})$.
\end{definition} 

The next theorem, which is taken from \cite{bleecker2013index}, says that amplitudes induce linear operators on smooth functions. We follow the authors using the convention $\diffj\xi = (2\pi)^{-n/2} \diff \xi$.
\begin{theorem}[\cite{bleecker2013index}*{Theorem 8.3}]
Every amplitude $p$ on $U \times \R^n$ defines a linear map
 \begin{equation*}
  \mathrm{Op}(p) \colon \mathcal{C}_c^\infty(U, \C^{N_1}) \rightarrow \mathcal{C}^\infty(U,\C^{N_2})
 \end{equation*}
 via the formula
 \begin{equation}
  \left(\mathrm{Op}(p)\right)(x) \deff \int_{\R^n}\exp(\I\langle x,\xi\rangle)p(x,\xi)\hat{u}(\xi) \diffj\xi.
 \end{equation}
 Even better, if $U = \R^n$, then this formula induces a map between the Schwartz spaces
 \begin{equation*}
  \mathrm{Op}(p) \colon \mathcal{S}(\C^{N_1}) \rightarrow \mathcal{S}(\C^{N_2}).
 \end{equation*}
\end{theorem}
We follow \cite{bleecker2013index} by calling the assignment $p \mapsto \mathrm{Op}(p)$ the \emph{quantisation of} $p$.
It is natural to ask whether an amplitude $p$ is uniquely defined through it quantisation $\mathrm{Op}(p)$ and under which condition $\mathrm{Op}(p)$ restricts to a map between real valued function. 
The answer is given in the next Proposition.
\begin{proposition}\label{AmplitudeVsOperator - Prop}
 Let $p$ be an amplitude on  $U \times \R^n$. Then 
 \begin{enumerate}
  \item[(i)] $p = 0$ if and only if $\mathrm{Op}(p)|_{\mathcal{C}_c^\infty(U,\R^{N_1})} = 0$.
  \item[(ii)] $p$ is Real if and only if $\mathrm{Op}(p)$ descends to $\mathcal{C}^\infty_c(U,\R^{N_1})\rightarrow \mathcal{C}^\infty(U,\R^{N_2})$.
 \end{enumerate}
\end{proposition}
\begin{proof}
 For (i), we prove the non-trivial direction by contra-position. 
 Assume that there is a point $(x_0,\xi_0)$ for which $p(x_0,\xi_0) \neq 0$. 
 Since the Fourier transformation is a bijection on the set of Schwartz functions \cite{konigsberger2004analysisII}*{p.331} the decay conditions for amplitudes imply that $\xi \mapsto p(x,\xi)\hat{u}(\xi)$
 is a Schwartz function. 
 Indeed, for each $\alpha \in \N_0^n$ and $l \geq 0$, we have
 \begin{align*}
  |\partial^\alpha_\xi p(x,\xi)\hat{u}(\xi)| &\leq \sum_{\beta \leq \alpha} ||\partial_\xi^\beta p(x,\xi)||\cdot ||\partial_\xi^{\alpha - \beta}\hat{u}(\xi)|| \\
  &\leq \sum C_{\beta,0,K}(1 + ||\xi||)^{k-|\beta|}\cdot ||\partial_\xi^{\alpha - \beta} \hat{u}(\xi)|| \\
  &\leq \sum C_{\beta, 0 , K} D_{\alpha - \beta, l + k}(1 + ||\xi||)^{k-|\beta|}(1+ ||\xi||)^{-(k-|\beta|+l)} \\
  &\leq \underbrace{\left(\sum C_{\beta, 0 , K} D_{\alpha - \beta, l + k}\right)}_{=: D'_{\alpha, l}}  (1+||\xi||)^{-l}.
\end{align*}
Thus, if we can find a real valued function $u$ with compact support in $U$ and $p(x_0,\xi_0)\hat{u}(\xi_0) \neq 0$, we are done, because the inverse Fourier transformation is a unitary operator on all square-integrable functions \cite{konigsberger2004analysisII}*{p.333}.

Note that $A|_{\R^{N_1}} = 0$ implies $A = 0$ for any $A \in \End[\C]{\C^{N_1}}{\C^{N_2}}$, so we can pick a real vector $v_0$ with $p(x_0,\xi_0)v_0 \neq 0$. 
By Lemma \ref{FouriertrafoAuxilarryFunction - Lemma} below, we can find an arbitrary large $R \in (\R_{>0})^n$ such that the Fourier transformation of the extension of $\prod_{j=1}^n \exp(-1/(1-R_j^2x_j^2))$ to $\R^n$ by zero evaluated at $\xi_0$ is not zero.
Thus,
\begin{equation*}
 \varphi(x) \deff \prod_{j=1}^n \exp\left(\frac{-1}{1-R^2(x_j - (x_0)_j)^2}\right)
\end{equation*}
has compact support within $U$ (since $R$ was chosen sufficiently large) and its Fourier transformation does not vanish at $\xi_0$.

For (ii), let $p$ be a Real amplitude.  If $u$ is an eigenfunction of $\overline{\cdot}$, so is $\mathrm{Op}(p)(u)$ as the following computation shows
\begin{align*}
 \overline{\mathrm{Op}(p)(x)} &= \overline{\int_{\R^n} \exp(\I\langle x,\xi\rangle)p(x,\xi)\hat{u}(\xi) \diffj\xi} = \int_{\R^n} \overline{\exp(\I\langle x,\xi\rangle)p(x,\xi)\hat{u}(\xi)} \diffj\xi \\
 &= \int_{\R^n} \exp(\I\langle x,\xi\rangle)\overline{p(x,\xi)}\overline{\hat{u}(\xi)} \diffj\xi \\
 &= \int_{\R^n} \exp(\I\langle x,-\xi\rangle)p(x,-\xi)\hat{u}(-\xi) \diffj\xi \\
 &= \int_{\R^n} |\det(-\id)| \exp(\I\langle x,\xi\rangle)p(x,\xi)\hat{u}(\xi) \diffj\xi = \mathrm{Op}(p)(x).
\end{align*}   
For the other direction, note that $p(x,\xi) -\overline{p(x,-\xi)}$ is also an amplitude.
If $\mathrm{Op}(p)$ restricts to a linear map $\mathcal{C}_c^\infty(U,\R^{N_1}) \rightarrow
\mathcal{C}^\infty(U, \R^{N_2})$, then the previous calculation shows that
\begin{equation*}
 \mathrm{Op}(p - \overline{p}(\placeholder, -\id(\placeholder)) = 0.
\end{equation*}
By part ($i$), this equality implies $p - \overline{p}(\placeholder, -\id(\placeholder)) = 0$, so $p$ is Real. 
\end{proof}
\begin{lemma}\label{FouriertrafoAuxilarryFunction - Lemma}
 For every $\xi_0 \in \R$ and there is an $R_0 > 0$ such that for all $R > R_0$ 
 the Fourier transformation of $x \mapsto \exp(-1/(1 - R^2x^2))$ does not vanish at $\xi_0$.
\end{lemma}
\begin{proof}
 Substitution of variables show
 \begin{align*}
  \int_{\R} \exp(-\I x\xi_0) \exp\left(\frac{-1}{1-R^2x^2}\right) \diffj x &= \int_{\R} \exp\left(-\I Rx\frac{\xi_0}{R}\right) \exp\left(\frac{-1}{1-R^2x^2}\right) \diffj x \\
  &= R^{-1} \int_{\R} \exp\left(-\I y\frac{\xi_0}{R}\right) \exp\left(\frac{-1}{1-y^2}\right) \diffj y \\
  &= R^{-1} \int_{-1}^1 \exp\left(-\I y\frac{\xi_0}{R}\right)\exp\left(\frac{-1}{1-y^2}\right) \diffj y\\
  &=: h(R).
 \end{align*}
 Lebesgue's dominated convergence theorem yields 
 \begin{equation*}
  \lim_{R \rightarrow \infty} R\cdot h(R) = \int_{-1}^1\exp(-\I y \cdot 0)\exp\left(\frac{-1}{1-y^2}\right) \diffj y > 0.
 \end{equation*}
 Since $h$ depends continuously on $R \in \R_{>0}$, the claim follows.
\end{proof}

\begin{definition}
 A linear map $P \colon \mathcal{C}_c^\infty(U,\C^{N_1}) \rightarrow \mathcal{C}^\infty(U,\C^{N_2})$ is a \emph{canonical pseudo differential operator of order} $k$ if there is a (necessarily unique) amplitude $p$ of order $k$ such that $P = \mathrm{Op}(p)$.
 
 In the real case, a linear map $P \colon \mathcal{C}_c^\infty(U,\R^{N_1}) \rightarrow \mathcal{C}^\infty(U,\R^{N_2})$ \emph{canonical pseudo differential operator of order} $k$ if its complex linear extension is a pseudo differential operator of order $k$.
\end{definition}
Note that the amplitude of the operator in the real case is necessarily Real by Proposition \ref{AmplitudeVsOperator - Prop}.
However, canonical pseudo differential operators are not good enough for our purposes. 
We need to restrict to a smaller subclass of operators that allow the definition of a principal symbol.
\begin{definition}
 A canonical pseudo differential operator of order $k$ is called \emph{principal classical} if its 
 amplitude $p$ satisfies the following assumptions:
 \begin{itemize}
  \item[(i)] For each $x \in U$ and each $\xi \in \R^n \setminus \{0\}$ the limit
  \begin{equation*}
   \symb_k(p)(x,\xi) \deff \lim_{\lambda \to \infty} \lambda^{-k} \cdot p(x,\lambda\xi)
  \end{equation*}
  exists. It is called the \emph{principal symbol} of $p$. 
  \item[(ii)] For some cut-off function (and therefore every) $\chi$ with
  \begin{align*}
   \chi(\xi) = \begin{cases}
                0, & \text{if }|\xi| \text{ is small,} \\
                1, & \text{if }|\xi| \geq 1, 
   \end{cases}
  \end{align*}
  the function $p(x,\xi) - \chi(\xi)\cdot\symb_k(p)(x,\xi)$ is a (canonical) amplitude of order $k-1$.
 \end{itemize}
\end{definition}

\begin{example}\label{ClassicalPseudos - Example}
 The following list of examples are taken form \cite{bleecker2013index}.
 \begin{enumerate} 
  \item Every differential operator $P = \sum_{|\alpha| \leq k} a_\alpha(x) \partial^\alpha$ of order $k$ is a principal classical pseudo differential operator of order $k$. 
  Indeed, using the Fourier inversion formula, it is not hard to see that $P$ has the amplitude
  $p(x,\xi) = \sum_{|\alpha| \leq k} i^{|\alpha|}a_\alpha(x)\xi^\alpha$. 
  
  The principal symbol is given by $\symb_k(P)(x,\xi) = \sum_{|\alpha| = k} i^k a_\alpha(x)\xi^\alpha$.
  This observation was actually the motivation for the defining formula of $\mathrm{Op}(p)$.
  
  \item As a special case of $1$, the identity $\id = \mathrm{Op}(\id)$ is a principal classical pseudo differential operator of order zero.
  More generally, the multiplication with a matrix valued function $f \in \mathcal{C}^\infty(U;\R^{N_1 \times N_2})$ is a pseudo differential operator of order zero. 
  Of course, this operator is Real if and only if $f$ is real a valued function.
  The principal symbol is given by $\symb_0(f \cdot)(x,\xi) = f(x)\cdot$.
  
  \item Every smooth function $K \colon U \times U \rightarrow \End{\R^{N_1}}{\R^{N_2}}$ such that
  $K(x,\cdot)$ has compact support for all $x \in U$ defines a principal classical
  pseudo differential of order $-\infty$, which means that it is of order $k$ for all $k \in \Z$. The principal symbol is therefore $0$.
 \end{enumerate}
\end{example}
\begin{proof}[Proof of Example 3]
 We will calculate the amplitude of $K$.
 Since $K(x,\cdot)$ has compact support and $\eta \mapsto \hat{u}(\eta)$ is a Schwartz function,
 the map $(y,\eta) \mapsto K(x,y)\hat{u}(\eta)$ is integrable, so all operations in the following chain of equations are valid
 \begin{align*}
  (Ku)(x) &= \int_{\R^n}K(x,y)u(y) dy \\
  &= \int_{\R^n}K(x,y)\int_{\R^n}\exp(\I\langle y,\eta\rangle)\hat{u}(\eta) \diffj\eta dy \\
  &= \int_{\R^n\times \R^n} \exp(\I\langle y,\eta\rangle)K(x,y)\hat{u}(\eta) \diffj\eta dy \\
  &= \int_{\R^n\times \R^n} \exp(\I\langle x,\eta\rangle)\exp(\I\langle y-x,\eta\rangle)K(x,y)\hat{u}(\eta) \diffj\eta dy \\
  &= \int_{\R^n\times \R^n} \exp(\I\langle x,\eta\rangle)\exp(\I\langle z,\eta\rangle)(2\pi)^{n/2}K(x,z-x)\hat{u}(\eta) \diffj z \diffj\eta \\
  &= \int_{\R^n}\exp(\I\langle x,\eta\rangle)k(x,\eta) \diffj\eta. 
 \end{align*}
 Here, $k(x,-\eta)$ is the Fourier transformation of $ z \mapsto (2\pi)^{n/2}K(x,z-x)$.
 Since $z \mapsto (\partial^\alpha_xK)(x,z-x)$ is compactly supported for all $\alpha \in \N_0^n$, its Fourier transformation is a Schwartz function (in $\eta$-direction). Thus, $K = \mathrm{Op}(k)$ is a pseudo differential operator of order $s$ for all $s \in \Z$.
\end{proof}
\begin{rem}
 An pseudo differential of order $-\infty$ is called \emph{smoothing operator}. 
 The name comes from the fact that such operators applied to arbitrary square integrable functions yield smooth functions, c.f. Theorem \ref{PseudosInduceBndMapsSobolev - Theorem} below.
\end{rem}
It is a deep computational result, for example presented in \cite{bleecker2013index}*{Theorem 8.19}, that canonical and principal classical pseudo differential operators are invariant under diffeomorphism in the following sense:

Suppose that $\kappa \colon U \rightarrow V$ is a diffeomorphism between relatively compact open subsets of $\R^n$ and that $P = \mathrm{Op}(p)$ is a canonical or principal classical pseudo differential operator of order $k$ those amplitude has compact support in the first variable, which means that there is a compact subset such that $p(x,\placeholder) = 0$ for all $x \in K^c$.
Then the pull-back $\kappa^\ast(P) \deff P(\cdot \circ \kappa^{-1})\circ \kappa$ is again a canonical or principal classical pseudo differential operator of order $k$ whose amplitude has compact support in the first variable.
Furthermore, if $q$ denotes the amplitude of $\kappa^\ast P$, then 
\begin{equation*}
 \symb_k(q)(x,\xi) = \symb_k(p)(\kappa(x),{}^T(D_x\kappa)^{-1}(\xi)).
\end{equation*}
In other words, $\symb_k(P)$ transforms like a section of $\mathrm{Hom}(\pi^\ast\C^{N_1},\pi^\ast\C^{N_2}) \rightarrow T^\vee U$.

We are now in the position to define pseudo differential operators on arbitrary manifolds.
\begin{definition}
 Let $E,F \rightarrow M$ be Real bundles over the manifold $M$. 
 A linear map $P \colon \Gamma_c(E) \rightarrow \Gamma(F)$ is a \emph{principal classical (canonically) pseudo differential operator of order k} 
 if for all $x_0 \in M$ there exists a local chart $(U, \kappa)$ around $x_0$,
 local trivialisations
 $\Phi\colon E_U \rightarrow U \times \C^{N_1}$, $\Psi\colon F_U \rightarrow U \times \C^{N_2}$,
 as well as a smooth function $f \colon M \rightarrow \lbrack 0, 1\rbrack$ that is identically $1$
 in a neighbourhood of $x_0$ and has compact support within $U$ such that, under the obvious identification of sections with functions, the linear map
 \begin{equation*}
  \kappa_\ast \left(\left(\Psi \circ f\cdot P \circ \Phi^{-1}\right)\right) \colon \mathcal{C}^\infty_c(\kappa(U),\R^{N_1}) \rightarrow \mathcal{C}^\infty(\kappa(U),\R^{N_2})
 \end{equation*}
 is a principal classical (canonical) pseudo differential operator of order $k$.
 If $p_f$ is the representing amplitude, we define
 \begin{equation*}
  \symb_k(P)(x_0,\xi) \deff \lim_{\lambda \to \infty} \frac{p_f(x_0,\lambda \xi)}{\lambda^k}.
 \end{equation*}
 We call such a choice $(U,\kappa,\Phi,\Psi,f)$ an \emph{amplitude datum}.
 The set of all principal classical pseudo differential operators is denoted by $\Pseud[k]$,
 \nomenclature{$\Pseud[k]$}{The space of (principal classical) pseudo differential operators of order $k$ between the vector bundles $E$ and $F$}
 the set of all canonical pseudo differential operators is denoted by $C\Pseud[k]$.
\end{definition}
\begin{rem}{$ $}
 The \emph{(principal) symbol} transforms accordingly under a change of coordinates and therefore gives rise a smooth homomorphism
  \begin{equation*}
   \symb_k(P) \colon \pi^\ast E \rightarrow \pi^\ast F,
  \end{equation*}
  where $\pi \colon T^\vee M\setminus\{0\} \rightarrow M$. The symbol is positive homogeneous of degree $k$,
  meaning $\symb_k(P)(x,\lambda \xi) = \lambda^k\symb(P)(x,\xi)$ for all $\lambda > 0$. 
  \nomenclature{$\symb_k(P)$}{Principal symbol of order $k$ of $P$}
  For a proof see \cite{bleecker2013index}*{Theorem 8.19}.
\end{rem}
Note that $Pu$ may not be compactly supported. 
However, when applying analytical tools, it will be convenient to control the increase of the support. Following \cite{LawsonMichelsonSpin}*{p.183} we make the following definition.
\begin{definition}
 A pseudo differential operator $P \in \PseudOp^k(E,F)$ has \emph{support in K} if $\supp Pu \subseteq K$ for all $u \in \Gamma_c(E)$ and $Pu=0$ if $\supp u \cap K = \emptyset$.  
\end{definition}
It follows from Proposition \ref{AmplitudeVsOperator - Prop} that $p_f \equiv 0$ for all $f \in \mathcal{C}_c^\infty(M)$ whose support does not intersect $K$.
In particular, $\symb_k(P)(x,\placeholder) = 0$ for all $x \in K^c$.
\begin{definition}\label{SetOfSymbols - Def}
 For each $k \in \R$ and each compact subset $K \subseteq M^n$, let 
 \begin{equation*}
  \Symb^k(E,F) \deff \{\sigma \in \Gamma(T^\vee M\setminus \{0\}, \mathrm{Hom}(\pi^\ast E, \pi^\ast F)) \, | \, \sigma(x,\lambda\xi) = \lambda^k\sigma(x,\xi)\} \nomenclature{$\Symb^k(E,F)$}{Space of Symbols}
 \end{equation*}  
 the vector space of all (smooth) symbols of order $k$
 and 
 \begin{equation*}
  \Symb^k_K \deff \{\sigma \in \Symb^k(E,F) \, : \, \sigma(x,\placeholder) = 0 \text{ if } x\in K^c\}
 \end{equation*}
 be its subspace of all symbols supported in $K$.
 The map that sends a pseudo differential operator to its principal symbol is, of course, denoted by
 \begin{equation*}
  \symb_k \colon \Pseud[k] \rightarrow \Symb^k(E,F).
 \end{equation*}
\end{definition}
If we choose a Riemannian metric, then we can can identify the set of all symbols with
$\Gamma(S(T^\vee M), \mathrm{Hom}(\pi^\ast E, \pi^\ast F))$, where $S(T^\vee M)$ denotes the unit-sphere bundle of the cotangent bundle $T^\vee M$ defined by the chosen Riemannian metric.

The defining conditions of pseudo differential operators are of local nature.
If one wants to verify a property of pseudo differential operators, for example, that the composition of two pseudo differential operators is again a pseudo differential operator, then one might be tempted to verify this property in local coordinates and then use a partition of unity argument, to reduce the global statement to a local one.
However, an arbitrary partition $\{\chi_j\}$ of unity is not good enough to deduce global properties of pseudo differential operators from local ones, even if each $\chi_j$ is supported in a coordinate chart. 
This comes from the fact that we need two different charts to detect $\chi_j P\chi_i u$ and we have no model for this situation to compare it against.

The following lemma, which is a generalisation of \cite{grubb2008distributions}*{Lemma 8.4} to non-compact manifolds, allows us to overcome this difficulty.

\begin{lemma}[Covering Lemma]\label{CoveringLemma - Lemma}
 Let $E,F \rightarrow M$ be vector bundles over a smooth manifold $M$. Then there exists an open covering of charts $\{\kappa_i \colon U_i \rightarrow V_i\}_{i\in I}$ over which the bundles trivialise and for which there is a subordinate partition of unity $\{\chi_i\}_{i \in I}$ such that each quadruple $(\chi_j,\chi_k,\chi_l,\chi_m)$ is supported in some $U_i$.
\end{lemma}
\begin{proof}
 By \cite{greene1978complete}*{Theorem 2'} there is a complete Riemannian metric $g$ on $M$ that has a positive injectivity radius.
 Using this information, it is easy to construct a countable, locally finite collection of charts $\{\kappa_j \colon U_j \rightarrow V_j\}_{j \in \N}$ whose targets $V_j$ are pair-wise disjoint and whose domains are contractible.
 
 We do this as follows:
 Let $1/3>\delta>0$ be smaller than the injectivity radius of $g$ and let $d = \mathrm{dist}(\placeholder , x_0)$ be distance function centered at $x_0 \in M$.
 Since $g$ is complete, $d$ is proper.
 Cover the compact sets $K_n \deff d^{-1}([n-1,n])$ by finitely many geodesic balls of radius $\delta$ if they are non empty.
 These geodesic balls together form a countable, locally finite open covering $\{U_j\}_{j \in \N}$ of $M$.
 If $\log$ denotes the inverse of the geodesic exponential map $\exp \colon B_\delta(0) \rightarrow U_j$ then $\{j + \log \colon U_j \rightarrow B_\delta(j)\}_{j \in \N}$ is a locally finite, countable set of charts with disjoint targets. 
 In particular, the vector bundles trivialise over the chart-domains as every bundle over contractible paracompact spaces are trivial, \cite{tomDieck2008algebraic}*{Theorem 14.3.1}.
 A close inspection of the proof of Lemma 8.4 in \cite{grubb2008distributions} shows that it is still applicable to the just constructed cover and that it yields a cover and a partition of unity with the stated properties. 
\end{proof}  

All examples presented in Example \ref{ClassicalPseudos - Example} carry over to vector bundles over manifolds. 
These examples are, informally speaking, quite local in nature.
However, there are globally defined linear operators on sections and it might be hard to guess a local their local amplitude to identify them as pseudo differential operators.
The first application of Lemma \ref{CoveringLemma - Lemma} is the following useful characterisation of pseudo differential operators that puts the operator in the centre instead of the amplitude.
Before we state it, recall from Appendix \ref{Sobolev Spaces - Chapter} that $\Gamma_c(M,E) = H^\infty_{cpt}$ and $\Gamma(M,E) = H^\infty_{loc}$, which allows us to topologise these spaces (with the inverse limit topology of the Sobolev spaces).
\begin{theorem}[\cite{atiyah1968indexI}*{p.509}]
  Let $E,F \rightarrow M$ be Real vector bundles over $M$. A linear map $P \colon \Gamma_c(M,E) \rightarrow \Gamma(M,F)$ is a pseudo differential operator of order $k$ if and only if $P$ is continuous and for each point $x_0 \in M$ there is a chart $(U,\kappa)$ around $x_0$, trivialisations $\Phi\colon E_U \rightarrow \C^{N_1}$, $\Psi\colon F_U \rightarrow U \times \C^{N_2}$, and a smooth function $f \colon M \rightarrow \lbrack 0, 1\rbrack$ that is identically $1$ near $x_0$ and compactly supported within $U$ such that, under the identification of sections and functions, the function
   \begin{equation*}
    p_f(y,\xi) \deff \exp(\I\langle y,\xi\rangle)\Psi \circ P( f \cdot \Phi^{-1}(\cdot, \exp(-\I \langle \kappa(\cdot), \xi\rangle)) 
   \end{equation*} 
   is a Real amplitude of order $k$.   
\end{theorem}
\begin{proof}[Sketch of Proof]
 Since the partition of unity from the Covering Lemma is locally finite, $P$ is continuous if and only if each summand $\chi_j P \chi_k$ is continuous.
 
 Since $\chi_j$ and $\chi_k$ are supported in the contractible chart $U_i$, it is enough to verify the claimed equivalence for all summands $\chi_jP\chi_k$ in local coordinates. 
 But this is done in \cite{grubb2008distributions}*{p.173ff}.
\end{proof} 
The Covering Lemma can also be used to reduce the following proposition on composition properties, which are stated in \cite{atiyah1968indexI}*{p.510} for the complex case, to the local case.
\begin{proposition}\label{App: Composition PSiDOs - Prop}
 Let $E,F,G$ be Real vector bundles over $M$.
 \begin{itemize}
  \item[(i)] If $P \in \Pseud[k]$, $Q \in  \Pseud[l][F][G]$, and $f \in \mathcal{C}_c^\infty(M)$, then $PfQ \in \Pseud[{k+l}][E][G]$ and $\symb(PfQ) = f\symb(P)\symb(Q)$.
  In particular, if $P \in \PseudOp^k_K(E,F)$, $Q \in \PseudOp_L^l(F,G)$, then $P Q \in \PseudOp_{K \cup L}^{k+l}(E,G)$.
  \item[(ii)] $P \in \PseudOp^k_{(K)}(E,F)$ implies $P^\vee \in \PseudOp^k_{(K)}(F^\vee,E^\vee)$ and $P^\ast = ({h_E}_\# \circ P^\vee \circ h_F^\flat) \in \PseudOp^k_{(K)}(\bar{F},\bar{E})$, if $E$ and $F$ are equipped with Riemannian metrics $h_E$ and $h_F$ respectively.  
  Furthermore, $\symb(P^\vee) = \symb(P)^\vee$ and $\symb(P^\ast) = \symb(P)^\ast$.
 \end{itemize}
\end{proposition}
\begin{proof}[Sketch of Proof]
 Use the Covering Lemma to reduce the statements to the case where $M$ is an open neighbourhood of $\R^d$, and $E$ and $F$ are trivial bundles.
 The statement now follows from the local version, as presented in \cite{bleecker2013index}*{Theorem 8.27}. 
\end{proof}

The proposition implies that $\PseudOp_K^0(E,E)$ is a $\ast$-algebra and that $\symb$ is a $\ast$-homomorphism. 
If $M$ is closed, then $\PseudOp_M^0(E,E) = \PseudOp^0(E,E)$ is unital.
Furthermore, in this case, $\Gamma_c(M,E) = \Gamma(M,E)$, so asking whether a pseudo differential operator is invertible actually makes sense.
Although the following result is a simple consequence of Lemma \ref{SymbolExactSequence - Lemma} below, we already state it here.
\begin{proposition}
 Let $E ,F\rightarrow M$ be a Real bundles over a closed manifold $M$. If $P \in \Pseud[k]$ is invertible, then $P^{-1} \in \Pseud[-k][F][E]$.
\end{proposition}
We use Proposition \ref{AmplitudeVsOperator - Prop} to endow $\Pseud[k]$ with a structure of a 
Fréchet space, even a Fréchet algebra, if we restrict our consideration to $\PseudOp_K^k(E,E)$. 

\begin{theorem}\label{PseudosFrechetNorm - Theorem}
Let $E, F, G$ be Real bundles over $M$.
 There is a Fréchet structure on $\Pseud[k]$ induced by the following semi-norms:
 For $A=(U,\kappa,\Phi,\Psi,f,\alpha,\beta)$ consisting of a chart $\kappa \colon U \rightarrow V$, trivialisations $\Phi$, $\Psi$ of $E_U$ and $F_U$, respectively, a smooth function $f$ with support within $U$, and $\alpha, \beta \in \N_0^n$, we set 
 \begin{equation}
   ||P||_A \deff \sup\left\{\frac{\partial_x^\alpha \partial_\xi^\beta p_f(x,\xi)}{(1 + |\xi|)^k}\right\}.
 \end{equation}
 Moreover, the compositions from Proposition \ref{App: Composition PSiDOs - Prop} are (jointly) continuous.
\end{theorem}
\begin{proof}[Sketch of Proof]
 The result about the Fréchet vector space structure is classical, mentioned, for example, in \cite{atiyah1971indexIV}. 
 The proof that the composition of two pseudo differential operator is a pseudo differential operator also shows that the composition is continuous with respect to this Fréchet structure by carefully tracking the constants in the involved estimations. 
\end{proof}

There are other topologies on $\Pseud[k]$ coming from operator theory. 
They arise from the fact that pseudo differential operators extend to Sobolev spaces.

\begin{theorem}[\cite{atiyah1968indexI}*{p.512}]\label{PseudosInduceBndMapsSobolev - Theorem}
 Every pseudo differential operator of order $k$ 
 \begin{equation*}
  P \colon \Gamma_c(M,E) \rightarrow \Gamma(M,F)
 \end{equation*}
 extends, for each $s \in \R$, to a continuous linear operator
 \begin{equation}
  P_s \colon \Sobolev[cpt]{s}{M,E} \rightarrow \Sobolev[loc]{s-k}{M,F}. 
 \end{equation}
 Moreover, these extensions induce a continuous map
 \begin{equation*}
     \PseudOp^k(E,F) \rightarrow Op^k(E,F).
 \end{equation*}
\end{theorem}
The theorem follows immediately from the Covering Lemma and the following local version that can be found in \cite{wells2008differential}*{Theorem 4.3.4} if one simply keeps track of the constants of the presented estimations.
\begin{lemma}\label{OperatorVsAmplitudeLocal - Lemma}
 Let $p \colon U \times \R^n \rightarrow \mathrm{Hom}(\C^{N_1},\C^{N_2})$ be an amplitude of order $k$ such that the support of the first variable lies in the compact set $K$.  
 
 Then $\mathrm{Op}(p) \colon \Gamma_c(U,\C^{N_1}) \rightarrow \Gamma_K(U,\C^{N_2})$ extends to a bounded map between Sobolev spaces $H^s(U,\C^{N_1}) \rightarrow H^{s-k}_K(U;\C^{N_2})$ with operator norm 
 \begin{align*}
  ||\mathrm{Op}(p)||_{s,s-k} &\leq  C \mathrm{vol}(K) \cdot \sum_{|\beta| \leq |s-k| + n + 1} ||p||_{\beta,0}.  
 \end{align*}
 Here, $||p||_{\beta,0}$ denotes the infimum of all valid bounds $C_{\beta,0,K}$ in the amplitude estimation \ref{eq: AmplitudeEstimation}.
\end{lemma}

Using the previous theorem, we can make the following definition.
\begin{definition}\label{ASClosure - Definition}
 Denote by $\PseudCl[k]$ the closure of $\Pseud[k]$ 
 \nomenclature{$\PseudOpCl^k_s(E,F)$}{Space of pseudo differential operators of order $k$ in the Atiyah-Singer sense}
 in $Op^k$ using the canonical inclusion. 
  We call this Fréchet space the \emph{Atiyah-Singer closure} of $\Pseud[k]$.
\end{definition}
We will also consider a third completion.
\begin{definition}
 For each $s \geq k$, we denote by $\PseudOpCl^k_s(E,F)$ the closure of $\Pseud[k]$ in $Op^k_s(E,F)$. 
 We will call this set the \emph{Sobolev closure} of $\Pseud[k]$.
\end{definition}
We can give $\Pseud[k]$ three different topologies, Amp, the Fréchet structure described via all representing amplitudes, AS, the subspace topology from the Atiyah-Singer closure, and $\mathrm{S}_s$, the subspace topology coming from the Sobolev closure.
For these three topologies, the operations of Proposition \ref{App: Composition PSiDOs - Prop} are continuous and extend to continuously to the closures.
They are related as follows:
\begin{theorem}
 For every two Real vector bundles $E, F$ over the compact manifold $M$, the identity give rise to continuous maps
 \begin{equation*}
  \xymatrix@C+0em{ (\Pseud[k],\mathrm{Amp}) \ar[r] & (\Pseud[k],\mathrm{AS})\ar[r] & (\Pseud[k],\mathrm{S}_s) }
 \end{equation*}
\end{theorem}
\begin{proof}
 The continuity of the first inclusion follows from Theorem \ref{PseudosInduceBndMapsSobolev - Theorem}.
 The second map is obviously continuous.
\end{proof}
Besides the smooth Fréchet topology, there is another reasonable topology on $\Symb_K^k(E,F)$ that turns $\Symb_K^k(E,F)$ into a Banach space. 
The norm is given by the supremum of the operator norm taken over all points of the sphere bundle.
Its completion will be denoted by $\overline{\Symb}_K^k(E,F)$. It can be identified with the set
of continuous section from the sphere bundle $S(T^\vee M)$ to $\mathrm{Hom}(\pi^\ast E, \pi^\ast F)$ that are supported within $K$.  

\begin{lemma}\label{SymbolExactSequence - Lemma}
 If $\PseudOp^k(E,F)$ carries the amplitude topology and $\mathrm{Symb}^k(E,F)$ carries the smooth topology, then the symbol map is continuous and has a continuous right-inverse $\rho$.
 The two maps restrict to the subspaces of operators and symbols with support in $K$, so we have the following two split exact sequences
 \begin{equation*}
  \xymatrix{0 \ar[r] & C\PseudOp^{k-1}_{(K)}(E,F) \ar[r] & \PseudOp^k_{(K)}(E,F) \ar[r]^{\mathrm{symb}^k} & \Symb_{(K)}^k(E,F) \ar[r] \ar@{.>}@/^1pc/[l]^{\rho} & 0} 
 \end{equation*}  
\end{lemma}
\begin{proof}
 We first construct the right-inverse $\rho$ closely following \cite{wells2008differential}. 
 
 We choose an open cover $\{U_\mu\}$ with a subordinate partition of unity $\{\phi_\mu\}$ as in the Covering Lemma.
 For each $\phi_\mu$, we pick a smooth function supported in $U_\mu$ that is identically $1$ on $\supp \phi_\mu$.
 Finally, we fix a cut-off function $\chi\colon \R \rightarrow \lbrack 0, 1\rbrack$ with $\chi \equiv 0$ on $\R_{\leq 0}$ and $\chi \equiv 1$ on $\R_{\geq 1}$.
 
 Let $\sigma \in \Gamma(T^\vee M \setminus \{0\}, \mathrm{Hom}(\pi^\ast E, \pi^\ast F))$ be a given symbol.
 Under the identifications $T^\vee U_\mu = U_\mu  \times \R^n$ and $E_{U_\mu} = U_\mu \times \C^{N_1}$, $F_{U_\mu} = U_\mu \times \C^{N_2}$, the symbol $\sigma$ corresponds to a smooth, matrix-valued map
 \begin{equation*}
  s_\mu \colon U_\mu \times \R^n \setminus \{0\} \rightarrow \mathrm{Hom}(\C^{N_1}, \C^{N_2}).
\end{equation*}
 Set $p_\mu(x,\xi) = \chi(||\xi||)\sigma_\mu(x,\xi)$ and extending by zero to $U_\mu \times \R^n$.
 It is an amplitude of order $k$ as $s_\mu$ is homogeneous of order $k$ in the second component.
 Define $P_\mu$ is defined via
 \begin{equation*}
  \xymatrix@C+1em{
  \mathcal{C}^\infty(U_\mu,\C^{N_1}) \ar[r]^{\psi_\mu \mathrm{Op}(p_\mu)} &
  \mathcal{C}_c^\infty(U_\mu,\C^{N_2}) \ar[r]^\cong &
  \Gamma_c(U_\mu, F) \ar@{^{(}->}[d] \\
  \Gamma_c(U_\mu, E) \ar[u]^\cong \ar[rr]_{P_\mu}& & \Gamma(M,F)}
\end{equation*}
and set 
\begin{equation*}
    \rho(\sigma) \deff P \deff \sum_\mu P_\mu.
\end{equation*}
Of course, the identifications are given by the charts and the local trivialisations.
It is clear that $P \in C\Pseud[k]$ because it is locally the quantisation of an amplitude.

In the chosen trivialisations, $\symb_k(P_\mu)$ is expressed by
\begin{align*}
 \symb_k(P_\mu)(x,\xi) &= \symb_k\left(\psi_\mu \cdot \xi(|| \placeholder||)s_\mu(\placeholder\,, \,\placeholder)\right)(x,\xi) \\
                       &= \lim_{\lambda \rightarrow \infty} \frac{\psi_\mu(x)\chi(||\lambda \xi||)s_\mu(x,\lambda \xi)}{\lambda^k} \\
                       &= \psi_\mu(x) \chi(||\lambda\xi||)s_\mu(x,\xi) = \psi_\mu(x)s_\mu(x,\xi).
\end{align*}     
With abuse of notation we have
\begin{equation*}
 \symb_k(P) = \sum_\mu \symb_k(P_\mu)\phi_\mu = \sum_\mu \psi_\mu s_\mu \phi_\mu = \sum_\mu \phi_\mu s_\mu = \sigma.
\end{equation*}

Clearly, $p_\mu(x,\xi) - \chi(||\xi||)\symb(p_\mu)(x,\xi) = 0$ is an amplitude of order $\leq k-1$.
Thus, each $P\cdot \phi_\mu$ is a principal classical pseudo differential operator of order $k$ and so is $P$. 
Consequently, $\rho$ is indeed a right-inverse of $\symb^k$.
Note that $\rho$ is $\mathcal{C}^\infty(M,\R)$-linear and that $\rho(\sigma)$ is Real, if $\sigma$ is Real.

Since $\rho$ was constructed by patching local right-inverse together it restricts to a right-inverse $\rho \colon \Symb_K^k(E,F) \rightarrow \PseudOp^k_K(E,F)$ of the symbol map.\\
$ $

The two maps, $\symb_k$ and $\rho$, are linear, so it suffices to check that they are continuous at $0$.

In the construction of $\rho$, we have chosen a partition of unity as in the Covering Lemma, so it is enough to verify its continuity in local coordinates.
We need to control the amplitude semi-norms of $\psi_\mu p_\mu$ in terms of the semi norms of $s_\mu$.
To this end, we calculate
\begin{align*}
 \partial_x^\alpha \partial_\xi^\beta \psi_\mu(x) \chi(\xi) s_\mu(x,\xi) = \sum_{\substack{\alpha' + \alpha'' = \alpha \\ \beta' + \beta'' = \beta}} \partial_x^{\alpha'}\psi_\mu(x) \partial_{\xi}^{\beta'}\chi(\xi) \partial_x^{\alpha''}\partial_\xi^{\beta''} s_\mu(x,\xi).
\end{align*}
Each summand $\partial_x^{\alpha''} s_\mu(x,\xi)$ is homogeneous of order $k$, so the chain rule implies
\begin{equation*}
 ||\partial_x^{\alpha''} \partial_\xi^{\beta''} s_\mu(x,\xi/||\xi||)|| \leq C_{\alpha'',\beta''} \sum_{l \leq |\beta''|} ||D^{\otimes l}|_{TS^{n-1}} s_\mu(x,\xi)|| \cdot (1 + ||\xi||)^{k - |\beta''|}
\end{equation*}
on $\mathrm{supp}\, \chi$, where $D^{\otimes l}|_{TS^{n-1}} s_\mu(x,\xi/||\xi||)$ denotes the $l$-th iterated (total) differential of $s_\mu(x,\placeholder)$ restricted to $(\R\xi^\perp)^{\otimes l}$.
The function $\chi$ vanishes near the origin and is constantly $1$ outside some compact set.
Thus,
\begin{equation*}
 \partial_\xi^{\beta'} \chi \leq C_{\beta'} (1+||\xi||)^{k-|\beta'|}
\end{equation*} 
for all $\beta' \neq 0$.
Putting these inequalities together, we conclude
\begin{equation*}
 ||\psi_\mu \chi s_m||_{\alpha,\beta;U_\mu} \leq C \sum_{\substack{l \leq |\alpha| + |\beta|}} ||D_{(x,\xi)}^{\otimes l} s_\mu||,
\end{equation*}
where $D_{(x,\xi)}^{\otimes l}s_\mu$ denotes the $l$-iterated total differential of $s_\mu$.
This shows that $\rho$ is continuous.

We deduce the continuity of the symbol map from the following inequality of amplitudes that are compactly supported within a chart domain
\begin{equation*}
 ||\symb_k(p)||_{C^\alpha;U_\mu} \leq \sum_{\alpha' + \alpha'' = \alpha} ||p||_{\alpha',\alpha'';U_\mu}, 
\end{equation*}
where $||p||_{\alpha',\alpha'';U_\mu}$ denotes the smallest constant in the amplitude defining inequality \ref{eq: AmplitudeEstimation} for some compact set that contains the support of $p$.
\end{proof}

\begin{lemma}\label{ASExactSequence - Lemma}
 The symbol map is also continuous if $\PseudOp^k_K(E,F)$ carries the $s$-Sobolev topology $S_s$ and $\Symb^k_K(E,F)$ carries the $C^0$-topology.
 The continuous extension to the completions yield a short exact sequences of Banach spaces
 \begin{equation*}
  \xymatrix{0 \ar[r] & \mathcal{K}(H^s,H^{s-k}) \ar[r] & \PseudOpCl^k_{s,K}(E,F) \ar[r]^{\overline{\symb}^k} & \overline{\Symb}^k_K(E,F) \ar[r] & 0.}
 \end{equation*}
\end{lemma}
\begin{proof}
 We will prove that the symbol map is continuous by showing
 \begin{equation*}
  ||\symb_k(P)||_{\infty} \leq C_{s,k}||P||_{s,s-k}.
 \end{equation*}
 We start with the special case $s = k = 0$.
 The support of $\symb_0(P)$ is compact, so we find a $(x_0,\xi_0) \in S(T^\vee M)$ with $||\symb_0(P)(x_0,\xi_0)||_{\mathrm{op}} = ||\symb_0(P)||_\infty$. 
 By a result attributed to Gohberg, see \cite{seeley1965integro}*{Theorem 2.2}, or the original \cite{gohberg1960Multidim}, applied in local coordinates, we find a sequence of smooth section $(u_n)_{n \in \N}$ such that
 \begin{itemize}
  \item[(i)] $u_n(x) = 0$ if $|x - x_0| > 1/n$,
  \item[(ii)] $||u_n||_{0} = 1$,
  \item[(iii)] $||\chi P u_n - \symb_0(P)(x_0,\xi_0) u_n||_0 \xrightarrow{n \to \infty} 0$ for some (and hence all) smooth function $\chi \leq 1$ that vanish on $|x - x_0| > 1$ and are identically $1$ on $|x-x_0|<1/2$.
  \end{itemize} 
  We may further assume that 
 \begin{itemize}
  \item[(iv)] all $u_n(x_0,\xi_0)$ are colinear, non-zero, and satisfy 
  \begin{equation*}
    ||\symb(P)(x_0,\xi_0)u_n(x_0,\xi)|| = ||\symb(P)(x_0,\xi_0)||_{\mathrm{op}}||u_n(x_0,\xi_0)||,
  \end{equation*}    
\end{itemize}   
  otherwise, we may replace $u_n$ by $A_n\cdot u_n$, where $A_n$ is a local endomorphism on $E$ that maps $u_n(x_0,\xi_0)$ to a vector realises the operator norm of $\symb_k(P)(x_0,\xi_0)$. 
  
  The claim for the special case follows now from property (iii) and 
  \begin{equation*}
   \lim_{n \to \infty} ||\mathrm{symb}_0(P)(x_0,\xi_0)u_n||_0 = ||\mathrm{symb}_0(P)(x_0,\xi_0)||_{\mathrm{op}} = ||\symb_0(P)||_\infty.
  \end{equation*}
  
  We reduce the general case to the special one as follows: 
  For some compact neighbourhood $L$ that contains $K$ in its interior, we pick  $\Lambda_{-s} \in \PseudOp^{-s}_L(E,E)$ and $\Lambda_{s-k} \in \PseudOp^{s-k}_L(F,F)$ that satisfy $\symb_j(\Lambda_j) = \id$ on $S(T^\vee M)|_K$.
  Then, $\symb_0(\Lambda_{s-k}P\Lambda_s) = \symb_k(P)$ on $S(T^\vee M)$ and 
  \begin{equation*}
   ||\Lambda_{s-k}P\Lambda_{-s}||_{0,0} \leq ||\Lambda_{s-k}||_{s-k,0} ||P||_{s,s-k} ||\Lambda_{-s}||_{0,s},
  \end{equation*}
  which implies the general statement.
  
  It remains to show that the continuous extension yields the claimed short exact sequence.
  The compact operators lie in the kernel of $\overline{\symb}_k$ because, with respect to the operator norm, every finite rank operator can be approximated by an infinitely smoothing operator and the latter lie in the kernel of $\symb_k$.
  
  For the converse, assume we are given $P \in \PseudOpCl_K^k(E,F)$ that does not lie in the kernel.
  First, we consider the special case $s=k=0$.
  By a diagonal sequence argument, we can extend Gohberg's result to the completion.
  Assume that $\overline{\symb}(P)(x_0,\xi_0) \neq 0$ and that $P$ were compact.
  If $u_n$ is a sequence with property (i) - (iv), then $Pu_n$ would have a convergent subsequence.
  To ease the notation, we assume the whole sequence converges.
  By property (iii), the sequence $\overline{\symb}_0(P)(x_0,\xi_0)u_n$ must be convergent.
  Condition (i) implies that the limit is zero.
  But condition (ii) and (iv) imply
  \begin{equation*}
   \lim_{n \to \infty}||\overline{\symb}_0(P)(x_0,\xi_0)u_n||_0 = ||\overline{\symb}_0(P)(x_0,\xi_0)||_{\mathrm{op}} \neq 0,
  \end{equation*}
  which is a contradiction.
  
  Assume, for general $s$ and $k$, there is a compact $P \in \PseudOpCl^k_K(E,F)$ such that $\overline{\symb}_k(P) \neq 0$.
  Then $\Lambda_{s-k}P\Lambda_{-s} \in \PseudOpCl^0(E,F)$ would be a compact operator with non-zero symbol contradicting the special case.
  
  We first prove that $\overline{\symb}_k$ is surjective for the special case $s = k = 0$ and $E = F$.
  In this case, $\symb_0$ is a $\ast$-homomorphism between $C^\ast$-algebras with dense image, so is must be surjective. 
  
  In the case $E \neq F$, we identify $\PseudOpCl^0_K(E,F)$ with the closed subspace of odd operators in $\PseudOpCl^0_K(E\oplus F,E \oplus F)$ and $\overline{\Symb}^0_K(E,F)$ with the closed subspace of $\overline{\Symb}^0_K(E \oplus F,E \oplus F)$ whose elements take values in the odd endomorphisms.
  These inclusions are isometric and commute with the symbol maps, so the image of the symbol map must be closed. 
  Since we already know from Lemma \ref{SymbolExactSequence - Lemma} that the image of $\overline{\symb}_0$  is dense, the symbol map must be surjective.
  
  For general $s$ and $k$, we argue as before: For every given symbol 
  \begin{equation*}
     \sigma \in \Gamma_K\bigl(S(T^\vee M),\mathrm{Hom}(\pi^\ast E,\pi^\ast F)\bigr)  
  \end{equation*}
   pick $Q \in \PseudOpCl^0(E,F)$ with $\overline{\symb}_0(Q) = \sigma$ and a pseudo differential $\Lambda_{k} \in \PseudOp_K^{s}(F)$  whose principal symbols is the identity.
  Then $\Lambda_{k}Q \in \PseudOpCl_K^k(E,F)$ and its principal symbol satisfies 
\begin{equation*}
 \symb_k(P) = \symb_0(Q) = \sigma
\end{equation*}   
on $S(T^\vee M)$.  
\end{proof}

\begin{cor}\label{SymbolExtensionAs - Corollary} 
 The symbol map extends to the Atiyah-Singer closure
 \begin{equation*}
  \symb_k \colon \PseudOpCl_K^k(E,F) \rightarrow \overline{\Symb}^k_K(E,F).
 \end{equation*}
\end{cor}

It is worth pointing out that the right inverse $\rho$ \emph{does not} extend in general to a continuous section $\overline{\rho} \colon \overline{\Symb}^k_K(E,F) \rightarrow \PseudOpCl^k_K(E,F)$. 
However, the right-inverse extends to continuous section if we control higher jets of the symbols.

\begin{lemma}
 The right inverse $\rho$ is continuous if $\PseudOpCl^k_{s,K}(E,F)$ carries the $s$-Sobolev topology $S_s$ and $\Symb_K^k(E,F)$ carries the $C^{|s-k|+\mathrm{dim}\, M +1}$-topology.
 Moreover, if $E$ and $F$ are Riemannian vector bundles with connections and if $\nabla$ denotes the induced connection on $\mathrm{Hom}(E,F)$, then 
 \begin{align*}
  ||\rho(\sigma)||_{s,s-k} \leq  C \sum_{r \leq |s-k| + \mathrm{dim}\, M + 1} ||\nabla^{\otimes r} \sigma||_{\infty; K}.
 \end{align*}
\end{lemma} 
\begin{proof}
 By construction, $\rho(\sigma) = \sum_\mu \mathrm{Op}(\psi_\mu p_\mu) \cdot \phi_\mu$, where $\psi_\mu$ and $\phi_\mu$ are supported in a domain of an amplitude datum and $\psi_\mu \phi_\mu = \phi_\mu$.
 This implies using Lemma \ref{OperatorVsAmplitudeLocal - Lemma}
 \begin{align*}
  ||P||_{s,s-k} &\leq \sum_\mu ||\mathrm{Op}(\psi_\mu p_\mu) \phi_\mu||_{s,s-k} \\
  & \leq \sum_\mu C(\phi_\mu) ||\mathrm{Op}(\psi_\mu p_\mu)||_{s,s-k;U_\mu} \\
  &\leq \sum_\mu C(U_\mu) \sum_{r \leq |s-k| + n +1} ||\psi_\mu p_\mu||_{r,0;U_\mu}.
 \end{align*}
 The calculations we made when we constructed $\rho$ imply in this special case, where we do not differentiate in $\xi$-direction (the fibre direction), that 
 \begin{align*}
  ||P||_{s,s-k} \leq \sum_\mu \tilde{C}(U_\mu) \sum_{|\alpha| \leq |s-k| + n +1} ||\partial^\alpha_x s_\mu||_{\infty,U_\mu}.
 \end{align*}
 Since $U_\mu$ is a chart domain over which the chart map extends to a slightly bigger open subset and in which $U_\mu$ is relative compact, all connections on $\mathrm{Hom}(E,F)$ restricted to $U_\mu$ yield equivalent $C^{|s-k|+\mathrm{dim}\, M + 1}$-norms.
 Thus, by increasing the constants, we deduce
 \begin{equation*}
     ||P||_{s,s-k} \leq \sum_\mu \tilde{C}(U_\mu) \sum_{r \leq |s-k| + n +1} ||\nabla^{\otimes r} s_\mu||_{\infty;U_\mu} \leq \tilde{C} \sum_{r \leq |s-k| + n +1} ||\nabla^{\otimes r} \sigma||_{\infty;K}. \qedhere
 \end{equation*}
\end{proof}

The reason why we mostly to work with the Atiyah-Singer closure instead of the Sobolev closure is that operators in the Atiyah-Singer closure share many properties of (actual) pseudo differential operators.
\begin{proposition}\label{PropertiesASPseudos - Proposition}
 {$ $}
 
 \begin{enumerate}
  \item Each $P \in \PseudCl[k]$ induces a continuous map $\Gamma_c(E) \rightarrow \Gamma(F)$.
  \item The space $\PseudCl[k]$ is a local space, which means that $P \in \PseudCl[k]$ if and only if $P|_U \in \PseudOpCl^k(E|_U,F|_U)$, where $P|_U$ denotes the restriction of $P$ to $\Gamma_c(U,E|_U) \rightarrow \Gamma(U,F)$.
 \end{enumerate}
\end{proposition}
\begin{proof}
 For the first statement see \cite{atiyah1968indexI}*{p.512}, for the second statement, see \cite{atiyah1968indexI}*{p.513}.
\end{proof}
The most important property is, in contrast to actual pseudo differential operators, that $\PseudCl[k]$ is stable under exterior tensor products.
\begin{theorem}[\cite{atiyah1968indexI}*{Theorem 5.3 and Theorem 5.4}]\label{ASTensor - Theorem}
 Let $E, F$ be Real bundles over $M$ and $G$ a Real bundle over $N$. Then for all $k>0$ there is a continuous map
 \begin{equation*}
   - \boxtimes \id \colon \PseudCl[k][M;E][F] \rightarrow \PseudCl[k][M\times N; E\boxtimes G][F \boxtimes G]
 \end{equation*}
 defined by sending $P$ to the continuous linear extension of 
 \begin{equation*}
  u \otimes v \mapsto P(u) \otimes v, 
 \end{equation*}
 for all $u \in \Gamma_c(M,E)$, $v \in \Gamma_c(N,G)$.
 Moreover, 
 \begin{equation*}
  \symb_k(P \boxtimes \id)(m,n;\xi,\eta) = \symb_k(P)(m,\xi)\otimes \id
 \end{equation*}
\end{theorem}
The proof of this theorem generalises to smooth families, provided we have the right notion of
smoothness.
\begin{theorem}\label{FamilyASTensor - Theorem}
 Let $P \colon \R^n \rightarrow (\Pseud[k],\mathrm{Ampl})$ be a smooth family and $k >0$. 
 Then there is a unique operator $P \boxtimes \id \in \PseudCl[k][E\boxtimes G][F\boxtimes G]$ that satisfies 
 \begin{equation*}
  (P \boxtimes \id)(u \otimes v)(m,t) = P_t(u)(m) \otimes v(t)  
 \end{equation*}
 and this assignment yields a continuous map
 \begin{equation*}
   \placeholder \boxtimes \id \colon \mathcal{C}^\infty(\R^n,\Pseud[k][M;E][F]) \rightarrow \PseudCl[k][M\times \R^n; E\boxtimes G][F \boxtimes G].
 \end{equation*}
\end{theorem}
 The proof of Theorem 5.4 in \cite{atiyah1968indexI} carries over without essential change. 
 The reason why we have to take the amplitude Fréchet structure is that we have to guarantee the 
 amplitude estimates also for the differentials with respect to the parameters $t \in \R^n$. 
 If we would chosen the topology coming from the Atiyah-Singer closure, then we could only conclude
 that $P \boxtimes \id$ is in $Op^k(M\times \R^n,E\boxtimes G, F \boxtimes G)$. With this choice of topology, $P$ would be simply a special smooth $Op^k$-valued map.
 
 Unfortunately, the proof of these theorems breaks down in the case $k=0$.
 Fortunately, for this thesis, the following weaker statement, which is a family version of \cite{atiyah1968indexI}*{Theorem 5.3}, is good enough. 
 Again, we will not prove this theorem as its proof essentially agrees with the one in \cite{atiyah1968indexI}. 
\begin{theorem}\label{App: ExternalTensor is Continuous - Theorem}
 Let $P \colon \R^n \rightarrow Op^k(M;E,F)$ be a smooth map of operators and $k\geq 0$.
 Let $G \rightarrow \R^n$ be a bundle over $\R^n$.
 Then there is a unique operator
 \begin{equation*}
  P \tilde{\boxtimes} \id \in Op^k(M \times \R^n; E \boxtimes G, F \boxtimes G)
 \end{equation*}
 that satisfies 
 \begin{equation*}
  (P \tilde{\boxtimes} \id)(u \otimes v)(m, t) = P_{t}(m) \otimes v(t),
 \end{equation*}
 and this assignment defines a continuous map
 \begin{equation*}
  \placeholder \, \tilde{\boxtimes} \id \colon \mathcal{C}^\infty(\R^n,Op^k(M;E,F)) \rightarrow Op^k(M \times \R^n; E \boxtimes G, F \boxtimes G).
 \end{equation*}
\end{theorem}

 \chapter{Formulas in Riemannian Geometry}

We prove formulas from differential geometry and geometric analysis that are need in the main text but may interrupt the reading flow of the text.

\section{The scalar curvature for a foliation of level sets}

If $(M,g)$ is a Riemannian manifold of dimension $d$ and $\mathfrak{d} \colon M \rightarrow \R$ is a function without critical values, then a standard result from differential topology, see \cite{hirsch1997differential}*{p.153ff}, says that the (shifted) gradient flow provides an isomorphism 
\begin{equation*}
    \xymatrix{\Phi \colon \mathfrak{d}^{-1}(\{a\}) \times [a,b] \ar[rd]_{\mathrm{pr}_2} \ar[rr] && \mathfrak{d}^{-1}([a,b]) \ar[ld]^{\mathfrak{d}} \\
    & [a,b] & }
\end{equation*}
Abbreviate $\mathfrak{d}^{-1}(\{a\})$ to $H$.
Let $v = [\gamma] \in T_mH$ be a tangent vector represented by a curve $\gamma \colon (-\eps,\eps) \rightarrow H$.
The calculation
\begin{align*}
    (\Phi^\ast g)_{(m,t_0)}(\partial_t,v) &= g\left(D_{(m,t_0)}\Phi(\partial_t),D_{(m,t_0)}\phi(v)\right) = g\left((\partial_t\Phi)_{(m,t_0)},D_{(m,t_0)}\phi(v)\right) \\
    &= g(||\grad{\mathfrak{d}}||_g^{-1/2}\grad{\mathfrak{d}}, [\Phi_{t_0}\circ \gamma]) \\
    &= ||\grad{\mathfrak{d}}||_g^{-1/2} \cdot \diff \, \mathfrak{d}([\Phi_{t_0}\circ \gamma]) \\
    &= ||\grad{\mathfrak{d}}||_g^{-1/2} \cdot [s \mapsto t_0] = 0
\end{align*}
implies the decomposition 
\begin{align*}
    (\Phi^\ast g)_{(m,t_0)} &= \bigl((\Phi_{t_0})^\ast g\bigr)_{m} + ||\grad{\mathfrak{d}}||_{g_{m,t_0}}^{-2} \diff t^2 \\
                            &=: h(t_0)_m + f^2(m,t_0) \diff t^2,
\end{align*}
where $h \colon [a,b] \rightarrow \Riem(H)$ is a curve of Riemannian metrics.

$ $

This example motivates the following setup.
Given a smooth curve of Riemannian metrics $h \colon (a,b) \rightarrow \Riem(M)$ and a smooth function $f \colon M \times (a,b) \rightarrow \R_{>0}$.
This data defines a Riemannina metric on $M \times (a,b)$ via 
\begin{equation*}
    g_{m,t} \deff h(t)_m + f(m,t)^2 \diff t^2.
\end{equation*}
We are interested in the relation between the curvatures of $g$ and $h$.
The following Proposition is a generalisation of \cite{bar2005generalized}*{Proposition 4.1}.
The proof presented there straighfowardly generalises to our set up.
 Only in the derivation of equation (\ref{Eq: LevelRiemanCurvIII}) below we need to be careful because the Riccati equation as stated in loc. cit. is not valid in our set up and needs to be replaced.
 Nonetheless, we present all calculations as a service to the reader.
 We also adopt the notation from their proof and write $\langle \placeholder , \placeholder \rangle$ instead of $g$ and 
 also denote the first and second derivative of $h$ with $\Dot{h}$ and $\Ddot{h}$, respectively. 
 They satisfy the relation
 \begin{align*}
     \Dot{h}(t_0)(v,w) &= \left. \frac{d}{dt}\Bigl( h(t)(v,w) \Bigr) \right|_{t=t_0},\\
     \Ddot{h}(t_0)(v,w) &= \left. \frac{d^2}{dt^2} \Bigl( h(t)(v,w) \Bigr) \right|_{t=t_0}
 \end{align*}
 for all tangent vectors $v,w \in TM_{t_0} = T(M \times \{t_0\}) \subseteq T(M \times (a,b))$.

 Instead of secretly identifying tangent vectors with local integrable vector fields, we will work with \emph{time-independent} vector fields that are tangent to the level sets.
 Because of the canonical split $T_{m,t}M \times (a,b) = T_mM \oplus T_t(a,b)$, each vector field $X$ on $T(M \times (a,b))$ can be written as $X(m,t) = [\Bar{X}(m,t),\lambda(m,t)]$, where $\Bar{X}(m,t)$ takes values in $TM_m$ and $\lambda(m,t)$ in $\R = T_{t}(a,b)$.
 A vector field $X$ on $T(M \times (a,b))$ is \emph{time-independent} if $\Bar{X}$ and $\lambda$ are independent of $t$.
 A vector field is \emph{tangent to the level sets} if $\lambda(m,t)=0$ for all $m \in M$ and $t \in (a,b)$.
 
 \begin{proposition}\label{Scalar Curvature Formuala - Prop}
  Let $h \colon (a,b) \rightarrow \Riem(M)$ be a smooth curve of Riemannian metrics and let $f \colon M \times (a,b) \rightarrow \R_{>0}$ be a smooth function. 
  If $g$ is the Riemannian metric on $M \times (a,b)$ given by $g_{(x,t)} = h(t)_x + f^2(x,t) \diff t^2$, if $X,Y,U,V$ are time independent vector fields that are tangent to the level sets, and $\nu = f^{-1}\partial_t$ is the unit normal vector field, then we have the following identities:
  \begin{align}
  \langle W(X),Y \rangle &= -\frac{1}{2} f^{-1} \Dot{h}(X,Y), \qquad \text{where } W(X) \deff - \nabla_X\nu \label{Eq: level-Weingarten}\\
  \langle R(U,V)X,Y\rangle &= \langle R^{M_t}(U,V)X,Y\rangle +\frac{1}{4}f^{-2}\bigl\lbrack \Dot{h}(U,X)\Dot{h}(V,Y) - \Dot{h}(U,Y)\Dot{h}(V,X)\bigr\rbrack \label{Eq: LevelRiemanCurvI} \\
  \begin{split}
    \langle R(X,Y)U,\nu \rangle &= -\frac{1}{2}f^{-1} \bigl\lbrack \langle (\nabla_X^{M_t}\Dot{h})(Y),U\rangle\rbrack - \langle(\nabla_Y^{M_t}\Dot{h})(X),U\rangle  \bigr\rbrack \\
    &\quad \ +\frac{1}{2}f^{-2} \bigl\lbrack X(f) \Dot{h}(Y,U) - Y(f)\Dot{h}(X,U) \bigr\rbrack
  \end{split}\label{Eq: LevelRiemanCurvII} \\
  \begin{split}
    \langle R(X,\nu)\nu,Y\rangle &= \frac{1}{2} f^{-2}\bigl\lbrack f^{-1}\Dot{f}\Dot{h}(X,Y) - \Ddot{h}(X,Y) + \frac{1}{2}\dot{h}({\dot{h}}^\mathrm{op}(X),Y) \bigr\rbrack \\
    & \quad \ +f^{-2}X(f)Y(f) - f^{-1}\mathrm{Hess}_{M_t}(f)(X,Y)
  \end{split}\label{Eq: LevelRiemanCurvIII} \\
  \begin{split}\label{Eq: Level-ScalarCurv}
      \scal(g) &= \scal(h) + f^{-2}\left(\frac{3}{4} \trace((\Dot{h}^{\mathrm{op}})^2) - \frac{1}{4} \trace(\Dot{h})^2  - \trace(\Ddot{h}) +f^{-1}\Dot{f}\trace(\Dot{h}) \right) \\
      &\qquad - 4f^{-2}|\diff f|^2_{M_t} + 2 f^{-1}\Delta_{M_t}(f).
   \end{split}
  \end{align}
 Here, $b^{\mathrm{op}} = b^{op,h}$ is the unique endomorphism that satisfies $b(v,w) = h(b^{\mathrm{op}}(v),w)$, where $b$ is an arbitrary bilinear form, $|\diff f|_{M_t}$ is the norm of the differential of $f|_{M_t}$ with respect to $h_t$, and $\Delta_{M_t}$ is the Laplace-Beltrami operator\footnote{the trace of the Hessian $\mathrm{Hess}_{M_t}$ with respect  to $h_t$.} of $M_t$ with respect to $h_t$.
 \end{proposition}
\begin{proof}
  If $X$ is tangent to the level sets, then $\langle X, \nu\rangle = 0$. 
  If $X$ is additionally time independent, then 
 \begin{equation}\label{Eq: Level-Commutator}
     [X,\nu] = [X,f^{-1} \partial_ t] = X(f^{-1})\partial_t + f^{-1}[X,\partial_t] = -f^{-2} X(f) \partial_t = -f^{-1}X(f)\nu.
 \end{equation}
 From 
 \begin{align*}
     0 &= \nu \langle X, \nu\rangle = \langle \nabla_\nu X, \nu\rangle + \langle X, \nabla_\nu \nu \rangle \\
       &= \langle W(X) , \nu \rangle - \langle [X,\nu], \nu \rangle + \langle X, \nabla_\nu \nu \rangle \\
       &= - \langle [X,\nu], \nu \rangle + \langle X, \nabla_\nu \nu \rangle \\
       &= -f^{-1}\cdot X(f) + \langle X, \nabla_\nu \nu \rangle
 \end{align*}
 and $0 = \nu 1 = \nu \langle \nu,\nu\rangle = 2 \langle \nabla_\nu \nu, \nu \rangle$ we conclude
 \begin{equation}\label{eq: Level-NormalDerivateunitNormal}
     \nabla_\nu \nu = -f^{-1} \grad{f_t}, 
 \end{equation}
 where $f_t \deff f(\placeholder,t)$.
 In particular, the vector field $\nabla_\nu \nu$ is always tangent to the level set.
 
 Since the Weingarten operator $W$ is self-adjoint, Equation (\ref{Eq: level-Weingarten}) follows from
 \begin{align*}
     \langle W(X),Y \rangle &= \langle -\nabla_X \nu, Y \rangle = - \langle \nabla_\nu X,Y\rangle \\
     &= - \nu \left( \langle X,Y\right) + \langle X, \nabla_\nu Y \rangle \\
     &= - f^{-1} \Dot{h}(X,Y) - \langle X, W(Y) \rangle.
 \end{align*}
 
 Equation (\ref{Eq: LevelRiemanCurvI}) then follows by plugging equation (\ref{Eq: level-Weingarten}) into the Gauß equation, see \cite{petersen2006riemannian}*{p.44} 
 \begin{align*}
     \begin{split}
         \langle R(U,V)X,Y\rangle &= \langle R^{M_t}(U,V)X,Y\rangle + \langle W(U),X\rangle \langle W(V),Y\rangle \\
         &\quad -\langle W(U),Y\rangle \langle W(V),X\rangle .
     \end{split}
 \end{align*}
 
 We derive equation (\ref{Eq: LevelRiemanCurvII}) by plugging equation (\ref{Eq: level-Weingarten}) into the Codazzi-Mainardi equation  \cite{petersen2006riemannian}*{p.44}\footnote{We use the opposite sign convention for the Weingarten operator} 
 \begin{equation*}
    \langle R(X,Y)U,\nu \rangle = \Bigl \langle (\nabla_X^{M_t}W)(Y),U  \Bigr \rangle - \Bigl \langle (\nabla_Y^{M_t}W)(X),U \Bigr \rangle 
 \end{equation*}
 and applying the product formula.
 
 To prove equation (\ref{Eq: LevelRiemanCurvIII}) we need to generalise the Riccati equation to our set up.
 The following calculation, as presented in \cite{petersen2006riemannian}*{p.44} with opposite sign convention, shows
 \begin{align*}
     (\nabla_\nu W)(X) - W^2(X) &= \nabla_\nu(W(X)) - W(\nabla_\nu X) - W(W(X)) \\
     &= - \nabla_{\nu}(\nabla_X \nu) + \nabla_{\nabla_\nu X}\nu - \nabla_{\nabla_X \nu} \nu \\
     &= - \nabla_{\nu}(\nabla_X \nu) - \nabla_{[X,\nu]} \nu \\
     &= R(X,\nu)\nu - \nabla_{X}\nabla_\nu \nu.
 \end{align*}
 Together with equation (\ref{eq: Level-NormalDerivateunitNormal}) we get the generalised Riccati equation
 \begin{equation}\label{Eq: Level-Riccati}
     R(X,\nu)\nu = (\nabla_\nu W)(X) - W^2(X) -  f^{-2}X(f) \cdot \grad{f_t} + f^{-1}\nabla_X \grad{f_t} 
 \end{equation}
 From the generalised Riccati equation, we derive equation (\ref{Eq: LevelRiemanCurvIII}) in two steps.
 We first consider the summands involving the Weingarten maps.
 Equation (\ref{Eq: Level-Commutator}) and equation (\ref{eq: Level-NormalDerivateunitNormal}) yield
 \begin{equation*}
     \langle W(X), [\nu,Y]\rangle = 0 \qquad \text{ and } \qquad \langle W([\nu,X]),Y\rangle = 0
 \end{equation*}
 because the Weingarten operator restricts to an endomorphism on $TM_{t}$.
 This implies
 \begin{align*}
    &\quad  \langle (\nabla_\nu W)(X),Y\rangle - \langle W^2(X), Y \rangle \\
    &= \nu \langle W(X),Y\rangle - \langle W(\nabla_\nu X), Y \rangle - \langle W(X), \nabla_\nu Y \rangle - \langle W^2(X),Y\rangle \\
    \begin{split}
        &= \nu\langle W(X),Y \rangle - \langle W(X),\nabla_\nu Y\rangle - \langle W(\nabla_X \nu + [\nu,X]),Y\rangle \\
        &\quad \ -\langle W^2(X),Y \rangle
    \end{split} \\
    &=\nu \langle W(X),Y\rangle + \langle W^2(X),Y\rangle \\
    &= f^{-1} \partial_t\bigl\lbrack \frac{-1}{2}f^{-1}\dot{h}(X,Y)\bigr\rbrack - \frac{1}{2}f^{-1}\dot{h}(W(X),Y) \\
    &= \frac{1}{2} f^{-2}\Bigl\lbrack f^{-1}\dot{f} \dot{h}(X,Y) - \Ddot{h}(X,Y) + \frac{1}{2} \dot{h}({\dot{h}}^{\mathrm{op}}(X),Y)\Bigr\rbrack
 \end{align*}
 
 For the other two summand of (\ref{Eq: Level-Riccati}), we use that the Levi-Cevita connection of the induced metric on a level set agrees with the tangential part of the Levi-Cevita of the ambient manifold to we deduce $\mathrm{Hess}_{M_t}(X,Y) = \langle \nabla_X \grad{f_t}, Y \rangle$. \
 Therefore, the remaining summand of (\ref{Eq: Level-Riccati}) can be  written as 
 \begin{align*}
    \langle f^{-2}X(f) \grad{f_t} - f^{-1}\nabla_X \grad{f_t}, Y \rangle = f^{-2}X(f)Y(f) - f^{-1}\mathrm{Hess}_{M_t}(f)(X,Y). 
 \end{align*}
 
 Splicing the results together, the generalised Riccati equation gives  equation (\ref{Eq: LevelRiemanCurvIII}).
 
 
 The equation for the scalar curvature now follows from plugging the previous results into the formula
 \begin{align*}
     \scal(g) &= \sum_{i=1}^d \mathrm{Ric}(e_i,e_i) + \mathrm{Ric}(\nu,\nu) = \sum_{i,j=1}^d \langle R(e_j,e_i)e_i,e_j\rangle + 2 \sum_{j=1}^d \langle R(e_j,\nu)\nu,e_j\rangle.
 \end{align*}
 Indeed, plugging equation (\ref{Eq: LevelRiemanCurvI}) into the first summand gives
 \begin{align*}
     \begin{split}
     \sum_{i,j=1}^d \langle R(e_j,e_i)e_i,e_j\rangle &= \sum_{i,j=1}^d \langle R^{M_t}(e_j,e_i)e_i,e_j\rangle + \frac{1}{4} f^{-2}\Bigl( \Dot{h}(e_j,e_i) \Dot{h}(e_i,e_j)  \\ &\qquad \qquad \qquad \qquad \qquad \qquad \qquad -\Dot{h}(e_j,e_j)\Dot{h}(e_i,e_i) \Bigr) 
     \end{split}\\
     &= \scal(h) + \frac{1}{4}f^{-2} (\trace((\Dot{h}^{\mathrm{op}})^2) - \trace(\Dot{h})^2),
 \end{align*}
 while plugging equation (\ref{Eq: LevelRiemanCurvIII}) into the second summand gives
 \begin{align*}
     \begin{split}
         2 \sum_{j=1}^d \langle R(e_j,\nu)\nu,e_j\rangle &= \sum_{j=1}^d f^{-2}\left(f^{-1} \Dot{f} \Dot{h}(e_j,e_j) - \Ddot{h}(e_j,e_j) + \frac{1}{2}\Dot{h}(\Dot{h}^{\mathrm{op}}(e_j),e_j)\right)\\ 
         & \qquad +2 f^{-2}e_j(f)e_j(f) - 2f^{-1} \mathrm{Hess}_{M_t}(f_t)(e_j,e_j) 
     \end{split}\\
     \begin{split}
         &= f^{-2}\left(\frac{1}{2}\trace((\dot{h}^{\mathrm{op}})^2) - \trace(\Ddot{h}) + f^{-1}\dot{f}\trace(\dot{h})\right) \\
         & \qquad +2f^{-2}|\diff f|^2_{M_t} - 2f^{-1}\Delta_{M_t}(f).
     \end{split}
 \end{align*}
 Their sum yields equation (\ref{Eq: Level-ScalarCurv}) so the proposition is proven.
\end{proof}

More importantly than the precise formula is how the scalar curvature behaves under rescaling.

\begin{definition}
 Let $U$ be a subset of $\R^n$ and $X$, $Y$ be sets. For every map $f \colon X \times U \rightarrow Y$ we denote with ${}_R f \colon X \times R\cdot U \rightarrow Y$ the map
 \begin{align*}
      {}_Rf(x,t) \deff f(x,R^{-1}\cdot t).
 \end{align*}
\end{definition}

\begin{cor}\label{Rescaling Behaviour - Cor}
 Let $h$ and $f$ be as in the previous proposition but assume further that their $2$-jet is uniformly bounded.
 Then 
 \begin{equation*}
  \scal(R^2 {}_Rh + {}_Rf^2 \diff t^2 )  \xrightarrow{R \to \infty} 0.
 \end{equation*}
 If, additionally, $f$ only depends on $t \in (a,b)$, then 
 \begin{equation*}
     \scal({}_Rh + {}_R f^2 \diff t^2 ) - {}_R\scal(h) \xrightarrow{R \to \infty} 0.
 \end{equation*}
 In particular, two isotopic psc metrics are concordant.
\end{cor}
\begin{proof}
 We first recall the rescaling behaviour for various objects of Riemannian geometry.
 In the following, let $g$ be an arbitrary Riemannian metric on $M$.
 
 For each smooth function $f\colon M \rightarrow \R$ the gradient satisfies
 \begin{equation*}
     \mathrm{grad}_{R^2g}(f) = R^{-2}\mathrm{grad}_g(f),
 \end{equation*}
 which implies
 \begin{align*}
     |\diff f|_{R^2g}^2 = R^2g(\mathrm{grad}_{R^2g}(f),\mathrm{grad}_{R^2g}(f)) = R^{-2}g(\mathrm{grad}_g(f),\mathrm{grad}_g(f)) = R^{-2}|\diff f|_g^2.
 \end{align*}
 
 If $(e_i)_i$ is a local orthonormal frame for $g$, then $(\Bar{e}_i = R^{-1}e_i)_i$ is a local orthonormal frame for $R^2g$.
 For each bilinear form $b$, this gives
 \begin{align*}
    \trace_{R^2g}(b) = \sum_{i=1}^n b(\Bar{e}_i, \Bar{e_i}) = R^{-2} \sum_{i=1}^n b(e_i,e_i) = R^{-2}\trace_{g}(b). 
 \end{align*}
 This result carries over to $\trace(b^{\mathrm{op}})$, because of the relation $b^{op,R^2g} = R^{-2}b^{op,g}$.
 
 The metrics $g$ and $R^2g$ have the same Levi-Cevita connection $\nabla$ because $\nabla$ is a metric connection with respect to both metrics and being torsion-free is a metric independent condition.
 This implies
 \begin{align*}
     \mathrm{Hess}_{R^2g}(f)(\Bar{e_i},\Bar{e_i}) &= R^2g(\nabla_{\Bar{e_i}} \mathrm{grad}_{R^2g}(f),\Bar{e_i}) \\
     &= g(\nabla_{{e_i}} \mathrm{grad}_{R^2g}(f),e_i) \\
     &= R^{-2}g(\nabla_{{e_i}} \mathrm{grad}_{g}(f),e_i) = R^{-2}\mathrm{Hess}_g(f)(e_i,e_i).
 \end{align*}
 In particular, 
 \begin{equation*}
     \Delta_{R^2g}(f) = R^{-2}\Delta_g(f).
 \end{equation*}
 
 Finally, we observe $\partial_t({}_R h) = R^{-1} \cdot {}_R (\partial_t h)$ and recall that $\scal(R^2h) = R^{-2}\scal(h)$.
 
 $ $
 
 With these information, we derive  
 \begin{equation*}
   \scal(R^2{}_Rh + {}_R f^2 \diff t^2) = R^{-2}\cdot {}_R\scal(h + f^2 \diff t^2).  
 \end{equation*}
 The assumption on the $2$-jet implies that $\scal(h + f^2\diff t^2)$ is uniformly bounded, so the first statement follows.
 
 The second statement follows from the rescaling behaviour and the previous lemma because the differential of $f$ ``vanishes in $M$-direction'':
 \begin{align*}
     \scal({}_Rh &+ {}_R f^2 \diff t^2) - {}_R\scal(h) \\
     &= R^{-2}\cdot {}_R\left(f^{-2}\left(\frac{3}{4} \trace((\Dot{h}^{\mathrm{op}})^2) - \frac{1}{4} \trace(\Dot{h})^2  - \trace(\Ddot{h}) +f^{-1}\Dot{f}\trace(\Dot{h}) \right)\right).  
 \end{align*}
 Since the $2$-jets of $f$ and $h$ are uniformly bounded, the right hand side converges to zero.
 
 The final statement follows the second statement applied to the special case $f = 1$.
\end{proof}
 
 \printnomenclature

 \bibliographystyle{alpha}
 \bibliography{Literatur}

\begin{bibdiv}
\begin{biblist}

\bib{antolini_cubical_2000}{article}{
      author={Antolini, Rosa},
       title={Cubical structures, homotopy theory},
        date={2000},
        ISSN={0003-4622},
     journal={Ann. Mat. Pura Appl. (4)},
      volume={178},
       pages={317\ndash 324},
         url={https://doi.org/10.1007/BF02505901},
      review={\MR{1849392}},
}

\bib{atiyah1968indexI}{article}{
      author={Atiyah, M.~F.},
      author={Singer, I.~M.},
       title={The index of elliptic operators. {I}},
        date={1968},
        ISSN={0003-486X},
     journal={Ann. of Math. (2)},
      volume={87},
       pages={484\ndash 530},
         url={https://doi.org/10.2307/1970715},
      review={\MR{236950}},
}

\bib{atiyah1969indexskew}{article}{
      author={Atiyah, M.~F.},
      author={Singer, I.~M.},
       title={Index theory for skew-adjoint {F}redholm operators},
        date={1969},
        ISSN={0073-8301},
     journal={Inst. Hautes \'{E}tudes Sci. Publ. Math.},
      number={37},
       pages={5\ndash 26},
         url={http://www.numdam.org/item?id=PMIHES_1969__37__5_0},
      review={\MR{285033}},
}

\bib{atiyah1971indexIV}{article}{
      author={Atiyah, M.~F.},
      author={Singer, I.~M.},
       title={The index of elliptic operators. {IV}},
        date={1971},
        ISSN={0003-486X},
     journal={Ann. of Math. (2)},
      volume={93},
       pages={119\ndash 138},
         url={https://doi.org/10.2307/1970756},
      review={\MR{279833}},
}

\bib{baum2009eichfeldtheorie}{book}{
      author={Baum, Helga},
       title={Eichfeldtheorie. {Eine} {Einf{\"u}hrung} in die
  {Differentialgeometrie} auf {Faserb{\"u}ndeln}},
    language={German},
   publisher={Berlin: Springer},
        date={2009},
        ISBN={978-3-540-38292-8; 978-3-540-38293-5},
}

\bib{bleecker2013index}{book}{
      author={Bleecker, David~D.},
      author={Boo\ss~Bavnbek, Bernhelm},
       title={Index theory---with applications to mathematics and physics},
   publisher={International Press, Somerville, MA},
        date={2013},
        ISBN={978-1-57146-264-0},
      review={\MR{3113540}},
}

\bib{booss1993elliptic}{book}{
      author={Boo\ss~Bavnbek, Bernhelm},
      author={Wojciechowski, Krzysztof~P.},
       title={Elliptic boundary problems for {D}irac operators},
      series={Mathematics: Theory \& Applications},
   publisher={Birkh\"{a}user Boston, Inc., Boston, MA},
        date={1993},
        ISBN={0-8176-3681-1},
         url={https://doi.org/10.1007/978-1-4612-0337-7},
      review={\MR{1233386}},
}

\bib{botvinnik2017infinite}{article}{
      author={Botvinnik, Boris},
      author={Ebert, Johannes},
      author={Randal-Williams, Oscar},
       title={Infinite loop spaces and positive scalar curvature},
        date={2017},
        ISSN={0020-9910},
     journal={Invent. Math.},
      volume={209},
      number={3},
       pages={749\ndash 835},
         url={https://doi.org/10.1007/s00222-017-0719-3},
      review={\MR{3681394}},
}

\bib{bar2005generalized}{article}{
      author={B\"{a}r, Christian},
      author={Gauduchon, Paul},
      author={Moroianu, Andrei},
       title={Generalized cylinders in semi-{R}iemannian and {S}pin geometry},
        date={2005},
        ISSN={0025-5874},
     journal={Math. Z.},
      volume={249},
      number={3},
       pages={545\ndash 580},
         url={https://doi.org/10.1007/s00209-004-0718-0},
      review={\MR{2121740}},
}

\bib{baaj1983theorie}{article}{
      author={Baaj, Saad},
      author={Julg, Pierre},
       title={Th\'{e}orie bivariante de {K}asparov et op\'{e}rateurs non
  born\'{e}s dans les {$C^{\ast} $}-modules hilbertiens},
        date={1983},
        ISSN={0249-6291},
     journal={C. R. Acad. Sci. Paris S\'{e}r. I Math.},
      volume={296},
      number={21},
       pages={875\ndash 878},
      review={\MR{715325}},
}

\bib{bamler2019ricci}{article}{
      author={Bamler, Richard~H},
      author={Kleiner, Bruce},
       title={Ricci flow and contractibility of spaces of metrics},
        date={2019},
     journal={arXiv preprint arXiv:1909.08710},
}

\bib{blackadar1998k}{book}{
      author={Blackadar, Bruce},
       title={{$K$}-theory for operator algebras},
     edition={Second},
      series={Mathematical Sciences Research Institute Publications},
   publisher={Cambridge University Press, Cambridge},
        date={1998},
      volume={5},
        ISBN={0-521-63532-2},
      review={\MR{1656031}},
}

\bib{berglund2012homological}{article}{
      author={Berglund, Alexander},
      author={Madsen, Ib},
       title={Homological stability of diffeomorphism groups},
        date={2012},
     journal={arXiv preprint arXiv:1203.4161},
}

\bib{boardman1999conditionally}{article}{
      author={Boardman, J.~Michael},
       title={Conditionally convergent spectral sequences},
        date={1999},
      volume={239},
       pages={49\ndash 84},
         url={https://doi.org/10.1090/conm/239/03597},
      review={\MR{1718076}},
}

\bib{salamon2018functional}{book}{
      author={B\"{u}hler, Theo},
      author={Salamon, Dietmar~A.},
       title={Functional analysis},
      series={Graduate Studies in Mathematics},
   publisher={American Mathematical Society, Providence, RI},
        date={2018},
      volume={191},
        ISBN={978-1-4704-4190-6},
         url={https://doi.org/10.1090/gsm/191},
      review={\MR{3823238}},
}

\bib{buggisch2019spectral}{book}{
      author={Buggisch, Lukas~Werner},
       title={The spectral flow theorem for families of twisted dirac
  operators},
     edition={PhD-Thesis},
   publisher={Westfaelische Wilhelms-Universitaet Muenster (Germany)},
        date={2019},
}

\bib{bunke2009index}{book}{
      author={Bunke, Ulrich},
       title={Index theory, eta forms, and {D}eligne cohomology},
        date={2009},
      volume={198},
      number={928},
        ISBN={978-0-8218-4284-3},
         url={https://doi.org/10.1090/memo/0928},
      review={\MR{2191484}},
}

\bib{bunke1995Callias}{article}{
      author={Bunke, Ulrich},
       title={A {$K$}-theoretic relative index theorem and {C}allias-type
  {D}irac operators},
        date={1995},
        ISSN={0025-5831},
     journal={Math. Ann.},
      volume={303},
      number={2},
       pages={241\ndash 279},
         url={https://doi.org/10.1007/BF01460989},
      review={\MR{1348799}},
}

\bib{carr1988construction}{article}{
      author={Carr, Rodney},
       title={Construction of manifolds of positive scalar curvature},
        date={1988},
        ISSN={0002-9947},
     journal={Trans. Amer. Math. Soc.},
      volume={307},
      number={1},
       pages={63\ndash 74},
         url={https://doi.org/10.2307/2000751},
      review={\MR{936805}},
}

\bib{cuntz2017ktheory}{book}{
      author={Cuntz, Joachim},
      author={Echterhoff, Siegfried},
      author={Li, Xin},
      author={Yu, Guoliang},
       title={{$K$}-theory for group {$C^*$}-algebras and semigroup
  {$C^*$}-algebras},
      series={Oberwolfach Seminars},
   publisher={Birkh\"{a}user/Springer, Cham},
        date={2017},
      volume={47},
        ISBN={978-3-319-59914-4; 978-3-319-59915-1},
      review={\MR{3618901}},
}

\bib{cisinski2006prefaisceaux}{book}{
      author={Cisinski, Denis-Charles},
       title={Les prefaisceaux comme modeles des types d'homotopie},
   publisher={These de doctorat de l’Universite Paris VII},
        date={2002},
}

\bib{cuntz2007topological}{book}{
      author={Cuntz, Joachim},
      author={Meyer, Ralf},
      author={Rosenberg, Jonathan~M.},
       title={Topological and bivariant {$K$}-theory},
      series={Oberwolfach Seminars},
   publisher={Birkh\"{a}user Verlag, Basel},
        date={2007},
      volume={36},
        ISBN={978-3-7643-8398-5},
      review={\MR{2340673}},
}

\bib{curtis1971simplicial}{article}{
      author={Curtis, Edward~B.},
       title={Simplicial homotopy theory},
        date={1971},
        ISSN={0001-8708},
     journal={Advances in Math.},
      volume={6},
       pages={107\ndash 209 (1971)},
         url={https://doi.org/10.1016/0001-8708(71)90015-6},
      review={\MR{279808}},
}

\bib{dubrovin1990modern}{book}{
      author={Dubrovin, BA},
      author={Fomenko, AT},
      author={Novikov, SP},
       title={Modern geometry -- methods and applications, part iii:
  Introduction to homology theory},
      series={Graduate texts in mathematics},
   publisher={Springer},
        date={1990},
      volume={124},
}

\bib{dugger2004topological}{article}{
      author={Dugger, Daniel},
      author={Isaksen, Daniel~C.},
       title={Topological hypercovers and {$\Bbb A^1$}-realizations},
        date={2004},
        ISSN={0025-5874},
     journal={Math. Z.},
      volume={246},
      number={4},
       pages={667\ndash 689},
         url={https://doi.org/10.1007/s00209-003-0607-y},
      review={\MR{2045835}},
}

\bib{davis2001lecture}{book}{
      author={Davis, James~F.},
      author={Kirk, Paul},
       title={Lecture notes in algebraic topology},
      series={Graduate Studies in Mathematics},
   publisher={American Mathematical Society, Providence, RI},
        date={2001},
      volume={35},
        ISBN={0-8218-2160-1},
         url={https://doi.org/10.1090/gsm/035},
      review={\MR{1841974}},
}

\bib{dold1963partitions}{article}{
      author={Dold, Albrecht},
       title={Partitions of unity in the theory of fibrations},
        date={1963},
        ISSN={0003-486X},
     journal={Ann. of Math. (2)},
      volume={78},
       pages={223\ndash 255},
         url={https://doi.org/10.2307/1970341},
      review={\MR{155330}},
}

\bib{ebert2016elliptic}{article}{
      author={Ebert, Johannes},
       title={Elliptic regularity for dirac operators on families of noncompact
  manifolds},
        date={2016},
     journal={arXiv preprint arXiv:1608.01699},
}

\bib{ebert2017indexdiff}{article}{
      author={Ebert, Johannes},
       title={The two definitions of the index difference},
        date={2017},
        ISSN={0002-9947},
     journal={Trans. Amer. Math. Soc.},
      volume={369},
      number={10},
       pages={7469\ndash 7507},
         url={https://doi.org/10.1090/tran/7133},
      review={\MR{3683115}},
}

\bib{eschenburg2021bott}{article}{
      author={Eschenburg, Jost-Hinrich},
      author={Hanke, Bernhard},
       title={Bott-{T}hom isomorphism, {H}opf bundles and {M}orse theory},
        date={2021},
        ISSN={1982-6907},
     journal={S\~{a}o Paulo J. Math. Sci.},
      volume={15},
      number={1},
       pages={127\ndash 174},
         url={https://doi.org/10.1007/s40863-021-00215-6},
      review={\MR{4258892}},
}

\bib{farrell1978rational}{inproceedings}{
      author={Farrell, F.~T.},
      author={Hsiang, W.~C.},
       title={On the rational homotopy groups of the diffeomorphism groups of
  discs, spheres and aspherical manifolds},
        date={1978},
   booktitle={Algebraic and geometric topology ({P}roc. {S}ympos. {P}ure
  {M}ath., {S}tanford {U}niv., {S}tanford, {C}alif., 1976), {P}art 1},
      series={Proc. Sympos. Pure Math., XXXII},
   publisher={Amer. Math. Soc., Providence, R.I.},
       pages={325\ndash 337},
      review={\MR{520509}},
}

\bib{frenck2021SphericalKappa}{article}{
      author={Frenck, Georg},
       title={Sphericity of kappa-classes and positive curvature via block
  bundles},
        date={2021},
}

\bib{gloeckle2019master}{article}{
      author={Gl\"ockle, Jonathan},
       title={Initial value spaces in general relativity},
        date={2019},
     journal={Master's Thesis (unpublished)},
}

\bib{gromov1980classification}{article}{
      author={Gromov, Mikhael},
      author={Lawson, H.~Blaine, Jr.},
       title={The classification of simply connected manifolds of positive
  scalar curvature},
        date={1980},
        ISSN={0003-486X},
     journal={Ann. of Math. (2)},
      volume={111},
      number={3},
       pages={423\ndash 434},
         url={https://doi.org/10.2307/1971103},
      review={\MR{577131}},
}

\bib{gromov1983positive}{article}{
      author={Gromov, Mikhael},
      author={Lawson, H.~Blaine, Jr.},
       title={Positive scalar curvature and the {D}irac operator on complete
  {R}iemannian manifolds},
        date={1983},
        ISSN={0073-8301},
     journal={Inst. Hautes \'{E}tudes Sci. Publ. Math.},
      number={58},
       pages={83\ndash 196 (1984)},
         url={http://www.numdam.org/item?id=PMIHES_1983__58__83_0},
      review={\MR{720933}},
}

\bib{gohberg1960Multidim}{article}{
      author={Gohberg, I.~C.},
       title={On the theory of multidimensional singular integral equations},
        date={1960},
        ISSN={0197-6788},
     journal={Soviet Math. Dokl.},
      volume={1},
       pages={960\ndash 963},
      review={\MR{0124704}},
}

\bib{gray1974volume}{article}{
      author={Gray, Alfred},
       title={The volume of a small geodesic ball of a {R}iemannian manifold},
        date={1973},
        ISSN={0026-2285},
     journal={Michigan Math. J.},
      volume={20},
       pages={329\ndash 344 (1974)},
         url={http://projecteuclid.org/euclid.mmj/1029001150},
      review={\MR{339002}},
}

\bib{greene1978complete}{article}{
      author={Greene, R.~E.},
       title={Complete metrics of bounded curvature on noncompact manifolds},
        date={1978/79},
        ISSN={0003-889X},
     journal={Arch. Math. (Basel)},
      volume={31},
      number={1},
       pages={89\ndash 95},
         url={https://doi.org/10.1007/BF01226419},
      review={\MR{510080}},
}

\bib{grubb2008distributions}{book}{
      author={Grubb, Gerd},
       title={Distributions and operators},
      series={Graduate Texts in Mathematics},
   publisher={Springer, New York},
        date={2009},
      volume={252},
        ISBN={978-0-387-84894-5},
      review={\MR{2453959}},
}

\bib{hatcher2001algebraic}{book}{
      author={Hatcher, Allen},
       title={Algebraic topology},
   publisher={Cambridge University Press, Cambridge},
        date={2002},
        ISBN={0-521-79160-X; 0-521-79540-0},
      review={\MR{1867354}},
}

\bib{hatcher1978concordance}{inproceedings}{
      author={Hatcher, A.~E.},
       title={Concordance spaces, higher simple-homotopy theory, and
  applications},
        date={1978},
   booktitle={Algebraic and geometric topology ({P}roc. {S}ympos. {P}ure
  {M}ath., {S}tanford {U}niv., {S}tanford, {C}alif., 1976), {P}art 1},
      series={Proc. Sympos. Pure Math., XXXII},
   publisher={Amer. Math. Soc., Providence, R.I.},
       pages={3\ndash 21},
      review={\MR{520490}},
}

\bib{hijazi1986conformal}{article}{
      author={Hijazi, Oussama},
       title={A conformal lower bound for the smallest eigenvalue of the
  {D}irac operator and {K}illing spinors},
        date={1986},
        ISSN={0010-3616},
     journal={Comm. Math. Phys.},
      volume={104},
      number={1},
       pages={151\ndash 162},
         url={http://projecteuclid.org/euclid.cmp/1104114937},
      review={\MR{834486}},
}

\bib{hirsch1997differential}{book}{
      author={Hirsch, Morris~W.},
       title={Differential topology},
      series={Graduate Texts in Mathematics},
   publisher={Springer-Verlag, New York},
        date={1994},
      volume={33},
        ISBN={0-387-90148-5},
        note={Corrected reprint of the 1976 original},
      review={\MR{1336822}},
}

\bib{hitchin1974harmonic}{article}{
      author={Hitchin, Nigel},
       title={Harmonic spinors},
        date={1974},
        ISSN={0001-8708},
     journal={Advances in Math.},
      volume={14},
       pages={1\ndash 55},
         url={https://doi.org/10.1016/0001-8708(74)90021-8},
      review={\MR{358873}},
}

\bib{higson2000analytic}{book}{
      author={Higson, Nigel},
      author={Roe, John},
       title={Analytic {$K$}-homology},
      series={Oxford Mathematical Monographs},
   publisher={Oxford University Press, Oxford},
        date={2000},
        ISBN={0-19-851176-0},
        note={Oxford Science Publications},
      review={\MR{1817560}},
}

\bib{hanke2014space}{article}{
      author={Hanke, Bernhard},
      author={Schick, Thomas},
      author={Steimle, Wolfgang},
       title={The space of metrics of positive scalar curvature},
        date={2014},
        ISSN={0073-8301},
     journal={Publ. Math. Inst. Hautes \'{E}tudes Sci.},
      volume={120},
       pages={335\ndash 367},
         url={https://doi.org/10.1007/s10240-014-0062-9},
      review={\MR{3270591}},
}

\bib{jardine2002cubical}{article}{
      author={Jardine, J.~F.},
       title={Cubical homotopy theory: a beginning},
        date={2002},
     journal={preprint},
       pages={39},
}

\bib{jardine2006categorical}{article}{
      author={Jardine, J.~F.},
       title={Categorical homotopy theory},
        date={2006},
        ISSN={1532-0073},
     journal={Homology Homotopy Appl.},
      volume={8},
      number={1},
       pages={71\ndash 144},
         url={http://projecteuclid.org/euclid.hha/1140012467},
      review={\MR{2205215}},
}

\bib{KanAbstractHomotopyI}{article}{
      author={Kan, Daniel~M},
       title={{Abstract} {Homotopy} {I}},
        date={1955},
     journal={Proceedings of the National Academy of Sciences of the United
  States of America},
      volume={41},
      number={12},
       pages={1092},
}

\bib{kasparov1980operator}{article}{
      author={Kasparov, G.~G.},
       title={The operator {$K$}-functor and extensions of {$C^{\ast}
  $}-algebras},
        date={1980},
        ISSN={0373-2436},
     journal={Izv. Akad. Nauk SSSR Ser. Mat.},
      volume={44},
      number={3},
       pages={571\ndash 636, 719},
      review={\MR{582160}},
}

\bib{kato1966perturbation}{book}{
      author={Kato, Tosio},
       title={Perturbation theory for linear operators},
      series={Die Grundlehren der mathematischen Wissenschaften, Band 132},
   publisher={Springer-Verlag New York, Inc., New York},
        date={1966},
      review={\MR{0203473}},
}

\bib{kriegl1997convenient}{book}{
      author={Kriegl, Andreas},
      author={Michor, Peter~W.},
       title={The convenient setting of global analysis},
      series={Mathematical Surveys and Monographs},
   publisher={American Mathematical Society, Providence, RI},
        date={1997},
      volume={53},
        ISBN={0-8218-0780-3},
         url={https://doi.org/10.1090/surv/053},
      review={\MR{1471480}},
}

\bib{konigsberger2004analysisII}{book}{
      author={K{\"o}nigsberger, Konrad},
       title={Analysis 2},
     edition={5},
   publisher={Springer-Verlag},
        date={2004},
}

\bib{kucerovsky1997kk}{article}{
      author={Kucerovsky, Dan},
       title={The {$KK$}-product of unbounded modules},
        date={1997},
        ISSN={0920-3036},
     journal={$K$-Theory},
      volume={11},
      number={1},
       pages={17\ndash 34},
         url={https://doi.org/10.1023/A:1007751017966},
      review={\MR{1435704}},
}

\bib{lee2013smooth}{book}{
      author={Lee, John~M.},
       title={Introduction to smooth manifolds},
     edition={Second},
      series={Graduate Texts in Mathematics},
   publisher={Springer, New York},
        date={2013},
      volume={218},
        ISBN={978-1-4419-9981-8},
      review={\MR{2954043}},
}

\bib{lichnerowicz1963spineurs}{article}{
      author={Lichnerowicz, Andr\'{e}},
       title={Spineurs harmoniques},
        date={1963},
        ISSN={0001-4036},
     journal={C. R. Acad. Sci. Paris},
      volume={257},
       pages={7\ndash 9},
      review={\MR{156292}},
}

\bib{Lueck2020surgery}{misc}{
      author={L{\"u}ck, Wolfgang},
      author={Macko, Tibor},
       title={Surgery theory: foundations},
   publisher={Ongoing Book Project},
        date={2020},
}

\bib{LawsonMichelsonSpin}{book}{
      author={Lawson, H.~Blaine, Jr.},
      author={Michelsohn, Marie-Louise},
       title={Spin geometry},
      series={Princeton Mathematical Series},
   publisher={Princeton University Press, Princeton, NJ},
        date={1989},
      volume={38},
        ISBN={0-691-08542-0},
      review={\MR{1031992}},
}

\bib{may1992simplicial}{book}{
      author={May, J.~Peter},
       title={Simplicial objects in algebraic topology},
      series={Chicago Lectures in Mathematics},
   publisher={University of Chicago Press, Chicago, IL},
        date={1992},
        ISBN={0-226-51181-2},
        note={Reprint of the 1967 original},
      review={\MR{1206474}},
}

\bib{milnor1957geometric}{article}{
      author={Milnor, John},
       title={The geometric realization of a semi-simplicial complex},
        date={1957},
        ISSN={0003-486X},
     journal={Ann. of Math. (2)},
      volume={65},
       pages={357\ndash 362},
         url={https://doi.org/10.2307/1969967},
      review={\MR{84138}},
}

\bib{murphy2014c}{book}{
      author={Murphy, Gerard~J.},
       title={{$C^*$}-algebras and operator theory},
   publisher={Academic Press, Inc., Boston, MA},
        date={1990},
        ISBN={0-12-511360-9},
      review={\MR{1074574}},
}

\bib{palais1965seminar}{article}{
      author={Palais, Richard~S.},
       title={Seminar on the {A}tiyah-{S}inger index theorem},
        date={1965},
       pages={x+366},
        note={With contributions by M. F. Atiyah, A. Borel, E. E. Floyd, R. T.
  Seeley, W. Shih and R. Solovay},
      review={\MR{0198494}},
}

\bib{palais1966homotopy}{article}{
      author={Palais, Richard~S.},
       title={Homotopy theory of infinite dimensional manifolds},
        date={1966},
        ISSN={0040-9383},
     journal={Topology},
      volume={5},
       pages={1\ndash 16},
         url={https://doi.org/10.1016/0040-9383(66)90002-4},
      review={\MR{189028}},
}

\bib{petersen2006riemannian}{book}{
      author={Petersen, Peter},
       title={Riemannian geometry},
     edition={Second},
      series={Graduate Texts in Mathematics},
   publisher={Springer, New York},
        date={2006},
      volume={171},
        ISBN={978-0387-29246-5; 0-387-29246-2},
      review={\MR{2243772}},
}

\bib{rosenberg2006manifolds}{article}{
      author={Rosenberg, Jonathan},
       title={Manifolds of {Positive} {Scalar} {Curvature}: {A} {Progress}
  {Report}},
    language={en},
        date={2006},
        ISSN={10529233, 21644713},
     journal={Surveys in Differential Geometry},
      volume={11},
      number={1},
       pages={259\ndash 294},
  url={http://www.intlpress.com/site/pub/pages/journals/items/sdg/content/vols/0011/0001/a009/},
}

\bib{rosenberg149metrics}{article}{
      author={Rosenberg, Jonathan},
      author={Stolz, Stephan},
       title={Metrics of positive scalar curvature and connections with
  surgery},
        date={2001},
      volume={149},
       pages={353\ndash 386},
      review={\MR{1818778}},
}

\bib{ruberman2002positive}{article}{
      author={Ruberman, Daniel},
       title={Positive scalar curvature, diffeomorphisms and the
  {S}eiberg-{W}itten invariants},
        date={2001},
        ISSN={1465-3060},
     journal={Geom. Topol.},
      volume={5},
       pages={895\ndash 924},
         url={https://doi.org/10.2140/gt.2001.5.895},
      review={\MR{1874146}},
}

\bib{rutter1997spaces}{book}{
      author={Rutter, John~W.},
       title={Spaces of homotopy self-equivalences},
      series={Lecture Notes in Mathematics},
   publisher={Springer-Verlag, Berlin},
        date={1997},
      volume={1662},
        ISBN={3-540-63103-8},
         url={https://doi.org/10.1007/BFb0093736},
        note={A survey},
      review={\MR{1474967}},
}

\bib{schick1998counterexample}{article}{
      author={Schick, Thomas},
       title={A counterexample to the (unstable) {G}romov-{L}awson-{R}osenberg
  conjecture},
        date={1998},
        ISSN={0040-9383},
     journal={Topology},
      volume={37},
      number={6},
       pages={1165\ndash 1168},
         url={https://doi.org/10.1016/S0040-9383(97)00082-7},
      review={\MR{1632971}},
}

\bib{seeley1965integro}{article}{
      author={Seeley, R.~T.},
       title={Integro-differential operators on vector bundles},
        date={1965},
        ISSN={0002-9947},
     journal={Trans. Amer. Math. Soc.},
      volume={117},
       pages={167\ndash 204},
         url={https://doi.org/10.2307/1994203},
      review={\MR{173174}},
}

\bib{stolz1992simply}{article}{
      author={Stolz, Stephan},
       title={Simply connected manifolds of positive scalar curvature},
        date={1992},
        ISSN={0003-486X},
     journal={Ann. of Math. (2)},
      volume={136},
      number={3},
       pages={511\ndash 540},
         url={https://doi.org/10.2307/2946598},
      review={\MR{1189863}},
}

\bib{stolz1998concordance}{article}{
      author={Stolz, Stephan},
       title={Concordance classes of positive scalar curvature metrics},
        date={1998},
     journal={Preprint},
}

\bib{stocker2013algebraische}{book}{
      author={St{\"o}cker, Ralph},
      author={Zieschang, Heiner},
       title={Algebraische topologie: eine einf{\"u}hrung},
     edition={2},
   publisher={Springer-Verlag},
        date={1994},
}

\bib{tomDieck2008algebraic}{book}{
      author={tom Dieck, Tammo},
       title={Algebraic topology},
      series={EMS Textbooks in Mathematics},
   publisher={European Mathematical Society (EMS), Z\"{u}rich},
        date={2008},
        ISBN={978-3-03719-048-7},
         url={https://doi.org/10.4171/048},
      review={\MR{2456045}},
}

\bib{weibel1996homological}{book}{
      author={Weibel, Charles~A.},
       title={An introduction to homological algebra},
      series={Cambridge Studies in Advanced Mathematics},
   publisher={Cambridge University Press, Cambridge},
        date={1994},
      volume={38},
        ISBN={0-521-43500-5; 0-521-55987-1},
         url={https://doi.org/10.1017/CBO9781139644136},
      review={\MR{1269324}},
}

\bib{wells2008differential}{book}{
      author={Wells, Raymond~O., Jr.},
       title={Differential analysis on complex manifolds},
     edition={Third},
      series={Graduate Texts in Mathematics},
   publisher={Springer, New York},
        date={2008},
      volume={65},
        ISBN={978-0-387-73891-8},
         url={https://doi.org/10.1007/978-0-387-73892-5},
        note={With a new appendix by Oscar Garcia-Prada},
      review={\MR{2359489}},
}

\bib{werner2006funktionalanalysis}{book}{
      author={Werner, Dirk},
       title={Funktionalanalysis},
    language={German},
     edition={6th corrected ed.},
      series={Springer-Lehrb.},
   publisher={Berlin: Springer},
        date={2007},
        ISBN={978-3-540-72533-6},
}

\end{biblist}
\end{bibdiv}
 
\end{document}